\documentclass[a4paper, reqno, 14pt]{amsart}

\usepackage[usenames,dvipsnames]{color}

\usepackage{amsthm,amsfonts,amssymb,amsmath,amsxtra}
\usepackage[all]{xy}
\SelectTips{cm}{}
\usepackage{xr-hyper}
\usepackage[colorlinks=
   citecolor=Black,
   linkcolor=Red,
   urlcolor=Blue,
   backref=page]{hyperref}
\usepackage{verbatim}

\usepackage[margin=1.25in]{geometry}

\usepackage[T1]{fontenc}
\usepackage{imakeidx}
\makeindex
\makeindex[name=NO,title=Index of Notation]

\usepackage{mathrsfs}

\RequirePackage{xspace}
\RequirePackage{etoolbox}
\RequirePackage{varwidth}
\RequirePackage{enumitem}
\RequirePackage{tensor}
\RequirePackage{mathtools}
\RequirePackage{longtable}
\RequirePackage{multirow}

\newcommand{\sB}{\ensuremath{\mathscr{B}}\xspace}

\newcommand{\BA}{\ensuremath{\mathbb {A}}\xspace}

\newcommand{\BC}{\ensuremath{\mathbb {C}}\xspace}

\newcommand{\BF}{\ensuremath{\mathbb {F}}\xspace}
\newcommand{\BG}{\ensuremath{\mathbb {G}}\xspace}

\newcommand{\BP}{\ensuremath{\mathbb {P}}\xspace}
\newcommand{\BQ}{\ensuremath{\mathbb {Q}}\xspace}
\newcommand{\BR}{\ensuremath{\mathbb {R}}\xspace}

\newcommand{\BX}{\ensuremath{\mathbb {X}}\xspace}
\newcommand{\BY}{\ensuremath{\mathbb {Y}}\xspace}
\newcommand{\BZ}{\ensuremath{\mathbb {Z}}\xspace}

\newcommand{\CA}{\ensuremath{\mathcal {A}}\xspace}

\newcommand{\CF}{\ensuremath{\mathcal {F}}\xspace}

\newcommand{\CJ}{\ensuremath{\mathcal {J}}\xspace}

\newcommand{\CL}{\ensuremath{\mathcal {L}}\xspace}
\newcommand{\CM}{\ensuremath{\mathcal {M}}\xspace}
\newcommand{\CN}{\ensuremath{\mathcal {N}}\xspace}
\newcommand{\CO}{\ensuremath{\mathcal {O}}\xspace}
\newcommand{\CP}{\ensuremath{\mathcal {P}}\xspace}

\newcommand{\CW}{\ensuremath{\mathcal {W}}\xspace}

\newcommand{\bK}{\ensuremath{\mathbf K}\xspace}

\newcommand{\ad}{{\mathrm{ad}}}

\DeclareMathOperator{\Aut}{Aut}

\DeclareMathOperator{\Coker}{Coker}

\newcommand{\cl}{{\mathrm{cl}}}

\DeclareMathOperator{\End}{End}

\DeclareMathOperator{\Gal}{Gal}
\newcommand{\GL}{\mathrm{GL}}

\newcommand{\GU}{\mathrm{GU}}

\DeclareMathOperator{\Hom}{Hom}

\newcommand{\id}{\ensuremath{\mathrm{id}}\xspace}

\newcommand{\Ind}{{\mathrm{Ind}}}

\newcommand{\inv}{{\mathrm{inv}}}
\DeclareMathOperator{\Isom}{Isom}

\DeclareMathOperator{\Ker}{Ker}

\DeclareMathOperator{\Lie}{Lie}

\DeclareMathOperator{\Nm}{Nm}

\DeclareMathOperator{\ord}{ord}

\DeclareMathOperator{\rank}{rank}

\newcommand{\PGL}{{\mathrm{PGL}}}

\newcommand{\red}{\ensuremath{\mathrm{red}}\xspace}

\DeclareMathOperator{\Res}{Res}

\newcommand{\SL}{{\mathrm{SL}}}
\DeclareMathOperator{\Spec}{Spec}
\DeclareMathOperator{\Spf}{Spf}

\newcommand{\Sp}{{\mathrm{Sp}}}
\newcommand{\SU}{{\mathrm{SU}}}

\DeclareMathOperator{\sgn}{sgn}

\DeclareMathOperator{\tr}{tr}

\newcommand{\U}{\mathrm{U}}

\newcommand{\wt}{\widetilde}
\newcommand{\wh}{\widehat}

\newcommand{\ov}{\overline}

\newcommand{\lra}{\longrightarrow}

\newcommand{\bs}{\backslash}

\newtheorem{theorem}{Theorem}
\newtheorem{proposition}[theorem]{Proposition}
\newtheorem{lemma}[theorem]{Lemma}

\newtheorem{corollary}[theorem]{Corollary}

\theoremstyle{definition}
\newtheorem{definition}[theorem]{Definition}
\newtheorem{example}[theorem]{Example}

\newtheorem{remark}[theorem]{Remark}
\newtheorem{remarks}[theorem]{Remarks}

\newenvironment{altenumerate}
   {\begin{list}
      {\textup{(\theenumi)} }
      {\usecounter{enumi}
       \setlength{\labelwidth}{0pt}
       \setlength{\labelsep}{0pt}
       \setlength{\leftmargin}{0pt}
       \setlength{\itemsep}{\the\smallskipamount}
       \renewcommand{\theenumi}{\roman{enumi}}
      }}
   {\end{list}}
\newenvironment{altitemize}
   {\begin{list}
      {$\bullet$}
      {\setlength{\labelwidth}{0pt}
	   \setlength{\itemindent}{5pt}
       \setlength{\labelsep}{5pt}
       \setlength{\leftmargin}{0pt}
       \setlength{\itemsep}{\the\smallskipamount}
      }}
   {\end{list}}

\numberwithin{equation}{subsection}
\numberwithin{theorem}{subsection}

\setitemize[0]{leftmargin=*,itemsep=\the\smallskipamount}
\setenumerate[0]{leftmargin=*,itemsep=\the\smallskipamount}

\renewcommand{\to}{%
   \ifbool{@display}{\longrightarrow}{\longrightarrow}%
   }
\let\shortmapsto\mapsto
\renewcommand{\mapsto}{%
   \ifbool{@display}{\longmapsto}{\shortmapsto}%
   }
\newlength{\olen}
\newlength{\ulen}
\newlength{\xlen}
\newcommand{\xra}[2][]{%
   \ifbool{@display}%
      {\settowidth{\olen}{$\overset{#2}{\longrightarrow}$}%
       \settowidth{\ulen}{$\underset{#1}{\longrightarrow}$}%
       \settowidth{\xlen}{$\xrightarrow[#1]{#2}$}%
       \ifdimgreater{\olen}{\xlen}%
          {\underset{#1}{\overset{#2}{\longrightarrow}}}%
          {\ifdimgreater{\ulen}{\xlen}%
             {\underset{#1}{\overset{#2}{\longrightarrow}}}
             {\xrightarrow[#1]{#2}}}}%
      {\xrightarrow[#1]{#2}}
   }
\makeatother
\newcommand{\xyra}[2][]{%
   \settowidth{\xlen}{$\xrightarrow[#1]{#2}$}%
   \ifbool{@display}%
      {\settowidth{\olen}{$\overset{#2}{\longrightarrow}$}%
       \settowidth{\ulen}{$\underset{#1}{\longrightarrow}$}%
       \ifdimgreater{\olen}{\xlen}%
          {\mathrel{\xymatrix@M=.12ex@C=3.2ex{\ar[r]^-{#2}_-{#1} &}}}%
          {\ifdimgreater{\ulen}{\xlen}%
             {\mathrel{\xymatrix@M=.12ex@C=3.2ex{\ar[r]^-{#2}_-{#1} &}}}
             {\mathrel{\xymatrix@M=.12ex@C=\the\xlen{\ar[r]^-{#2}_-{#1} &}}}}}%
      {\mathrel{\xymatrix@M=.12ex@C=\the\xlen{\ar[r]^-{#2}_-{#1} &}}}%
   }
\makeatletter
\newcommand{\xla}[2][]{%
   \ifbool{@display}%
      {\settowidth{\olen}{$\overset{#2}{\longleftarrow}$}%
       \settowidth{\ulen}{$\underset{#1}{\longleftarrow}$}%
       \settowidth{\xlen}{$\xleftarrow[#1]{#2}$}%
       \ifdimgreater{\olen}{\xlen}%
          {\underset{#1}{\overset{#2}{\longleftarrow}}}%
          {\ifdimgreater{\ulen}{\xlen}%
             {\underset{#1}{\overset{#2}{\longleftarrow}}}
             {\xleftarrow[#1]{#2}}}}%
      {\xleftarrow[#1]{#2}}
   }
\newcommand{\isoarrow}{%
   \ifbool{@display}{\overset{\sim}{\longrightarrow}}{\xrightarrow\sim}%
   }
   
\DeclareMathOperator{\Bil}{Bil}
\DeclareMathOperator{\Trace}{Tr}
\DeclareMathOperator{\Nilp}{Nilp}
\DeclareMathOperator{\height}{height}

\subjclass[2010]{14G20,  14L05, 14G35}
\keywords{Shimura curves, $p$-adic uniformization,  displays of formal $p$-divisible groups.}

\begin{document}


\title{On the $p$-adic uniformization of unitary Shimura curves}
\author{Stephen Kudla}
\address{Department of Mathematics,
University of Toronto,
40 St. George St., BA6290,
Toronto, ON M5S 2E4, Canada.}
\email{skudla@math.toronto.edu}
\author{Michael Rapoport}
\address{Mathematisches Institut der Universit\"at Bonn, Endenicher Allee 60, 53115 Bonn, Germany, and University of Maryland, Department of Mathematics, College Park, MD 20742, USA}
\email{rapoport@math.uni-bonn.de}
\author{Thomas Zink}
\address{Fakult\"at f\"ur Mathematik,
Universit\"at Bielefeld,
Postfach 100131,
33501 Bielefeld, Germany}
\email{zink@math.uni-bielefeld.de}

\date{\today}
\begin{abstract}

We prove $p$-adic uniformization for Shimura curves attached to the group of unitary similitudes of certain binary skew hermitian spaces $V$ with respect to an 
arbitrary CM field $K$ with maximal totally real subfield $F$. For a place $v|p$ of $F$ that is not split in $K$ and for which $V_v$ is anisotropic, 
let $\nu$ be an extension of $v$  to the reflex field $E$. We define an integral model of the corresponding Shimura curve over $\Spec O_{E, (\nu)}$ by means of a moduli problem for abelian schemes with suitable 
polarization and level structure prime to $p$. The formulation of the moduli problem 
involves a \emph{Kottwitz condition}, an \emph{Eisenstein condition}, and an \emph{adjusted invariant}. 
The uniformization of the formal completion of this model along its special fiber is given 
in terms of the formal Drinfeld upper half plane $\widehat\Omega_{F_v}$ for $F_v$. The proof  relies on the construction of the \emph{contracting functor} which relates a relative Rapoport-Zink space 
for strict formal $O_{F_v}$-modules with a Rapoport-Zink space of $p$-divisible groups which arise from the  moduli problem, where the $O_{F_v}$-action 
is usually not strict when $F_v\ne \BQ_p$. Our main tool is the theory of displays, in particular   the \emph{Ahsendorf functor}.

\end{abstract}
\maketitle
\tableofcontents

\section{Introduction} 
\subsection{History of uniformization}
One of the major results of the Mathematics of the 19th century is the \emph{uniformization theorem}. It states that any non-singular projective algebraic curve $X$ of genus $g(X)\geq 2$ can be uniformized, i.e., can be written as
\begin{equation}
X\simeq \Gamma\backslash\Omega_{\BR},  
\end{equation}
where $\Omega_\BR=\BP^1(\BC)\setminus\BP^1(\BR)$\index[NO]{ZZWA@$\Omega_\BR$} is the union of the upper and the lower half plane and  $\Gamma$ denotes a discrete cocompact subgroup of $\PGL_2(\BR)$. This notation reinforces the analogy with the $p$-adic uniformization discussed below.
The history of this theorem is very complicated, and involves the names of many mathematicians, among them Poincar\'e, Hilbert and Koebe, comp. \cite{Gray}. Inspired by the uniformization theorem, Poincar\'e gave a systematic construction of cocompact discrete subgroups of $\PGL_2(\BR)$. For this he used the exceptional isomorphism between inner forms of $\PGL_2$ and special orthogonal groups of  ternary quadratic forms. In fact, for his construction, he used arithmetic subgroups of the special orthogonal group of an indefinite anisotropic ternary quadratic form over $\BQ$, cf. \cite{Gray}. 

\medskip

Now let $p$ be a prime number. The history of the $p$-adic uniformization of algebraic curves starts with \index{Tate uniformization} Tate's uniformization theory of elliptic curves. It turns out that not all elliptic curves over $p$-adic fields admit a $p$-adic uniformization, but only those  with (split) multiplicative reduction \cite[\S 6]{Ta}.

The next step was Mumford's $p$-adic uniformization theory of algebraic curves of higher genus, \cite{Mum}\index{Mumford}. Again, it turns out that not all such algebraic curves  over $p$-adic fields admit a $p$-adic uniformization, but only those with totally degenerate reduction \cite{Mum}. In view of Mumford's results, it becomes interesting to single out  classes of algebraic curves with totally degenerate reduction.  Such classes are exhibited by Cherednik \cite{Ch}.

 Cherednik's discovery is that certain \emph{quaternionic Shimura curves}\index{quaternionic Shimura curve}, i.e., Shimura curves  associated to 
quaternion algebras over a totally real field $F$, admit $p$-adic uniformization. The quaternion algebra has to satisfy the following conditions. It  is required to be split at precisely one archimedean place $w$ of $F$ (and ramified at all other archimedean places), and to be  ramified at a non-archimedean place $v$ of residue characteristic $p$. In this case, the reflex field can be identified with $F$. Then one obtains $p$-adic uniformization by the Drinfeld halfplane \index{Drinfeld halfplane} associated to $F_v$, provided that the
level structure is prime to $v$.   It follows  that if $X$ is a connected component of the Shimura tower for such a level, considered as an algebraic curve over $\bar F$, then there is an isomorphism of algebraic curves over $\bar F_v$, 
\begin{equation}\label{Chered}
X\otimes_{\bar F} \bar F_v\simeq(\bar\Gamma\backslash\Omega_{F_v})\otimes_{F_v}\bar F_v .
\end{equation}
Here $\Omega_{F_v}=\BP^1_{F_v}\setminus \BP^1(F_v)$
\index[NO]{ZZWB@$\Omega_{F_v}$} denotes the Drinfeld  halfplane \index{Drinfeld halfplane} for the local field $F_v$, and $\bar\Gamma$ denotes a discrete cocompact subgroup of $\PGL_2(F_v)$. Recall that $\Omega_{F_v}$ is a rigid-analytic space over $F_v$. The isomorphism \eqref{Chered} is to be interpreted as follows: the rigid-analytic space $\bar\Gamma\backslash\Omega_{F_v}$ is (uniquely) algebraizable by a  projective algebraic curve over $F_v$. After extension of scalars $F_v\to\bar F_v$, there exists an isomorphism as in \eqref{Chered}. We thus see that \eqref{Chered} allows us to pass from the original complex uniformization $X\otimes_{\bar F}\BC\simeq 
\Gamma\backslash\Omega_{\BR}$, where $\Gamma$ is a congruence subgroup maximal at $v$, to $p$-adic uniformization.

 Let us comment on the proof of Cherednik's theorem. When $F=\BQ$, these quaternionic  Shimura curves are moduli spaces of abelian varieties with additional structure, and  Drinfeld \cite{Dr} gives a moduli-theoretic proof of Cherednik's theorem\index{Cherednik's theorem} in this special case.  Furthermore, he proves an `integral version' of this theorem (which has the original version as a corollary).  For this, Drinfeld extends the moduli problem integrally  and then relates the integral version to a theorem on formal moduli spaces of $p$-divisible groups, which is in fact the deepest part of Drinfeld's paper. When $F\neq \BQ$, Cherednik's quaternionic Shimura curves do not represent a moduli problem of abelian varieties, and  Drinfeld's approach runs into problems. Cherednik's  approach \cite{Ch} seems to only use arguments involving the generic fiber.

There are also higher-dimensional versions of $p$-adic uniformization. Drinfeld's method  has been generalized by Rapoport and Zink \cite{RZ} to Shimura varieties associated to certain \emph{ fake unitary groups}. These are associated to central division algebras over a CM-field equipped with an involution of the second kind; for Rapoport-Zink uniformization, one has to assume that the  $p$-adic place of the totally real subfield splits in the CM-field. This higher-dimensional generalization also includes integral uniformization theorems. In \cite{RZ}, these integral uniformization theorems appear as a special instance of a general \emph{non-archimedean uniformization theorem}, which describes the formal completion of PEL-type Shimura varieties along a fixed isogeny class. In the case of $p$-adic uniformization, the whole special fiber forms a single isogeny class.

The method of \cite {RZ} has been applied by Boutot and Zink \cite{BZ} to prove  Cherednik's original theorem and an integral variant of it  by embedding Cherednik's quaternionic Shimura curves into Shimura curves obtained by the Rapoport-Zink method; in an update \cite{BZ1}, some  gaps in \cite{BZ} are filled. The integral uniformization theorems in \cite{BZ1} have the draw-back  that they only show that there exists some integral model of the Shimura curve for which one has integral uniformization. There is a characterization of this integral model as the unique stable model in the sense of Deligne-Mumford \cite{DM} but this characterization is of  a rather abstract nature since there is no moduli-theoretic description of it. 

 A variant of Cherednik's method  has been developed by Varshavsky\index{Varshavsky} \cite{Va1, Va2} to obtain $p$-adic uniformization of certain higher-dimensional Shimura varieties associated  to fake unitary groups, again at a split place. We refer to Boutot's Bourbaki talk \cite{B} for an account of all these developments.

In the present paper, we deal with Shimura curves attached to unitary similitude groups associated to anti-hermitian\footnote{It turns out to be more natural to consider anti-hermitian forms, rather than hermitian forms, cf. below.} vector spaces $V$ of dimension  $2$ over a CM-field  $K$ with totally real subfield $F$ of arbitrary degree. Our results generalize those in \cite{KRnew}, where the case $F=\BQ$  is considered. Like Cherednik, we assume that  $V$ is split at precisely one archimedean
place $w$ of $F$ (and ramified at all other archimedean places). We also assume that $V$ is ramified at a non-archimedean place $v$ of residue characteristic $p$ of $F$. However, in contrast to the cases of $p$-adic uniformization mentioned above, we assume that $v$ does \emph{not} split in $K$.  Of course, these Shimura curves are closely related to the Shimura curves considered by Cherednik (we refer to \cite{KRnew} for a general discussion of the relation between quaternion algebras and two-dimensional hermitian vector spaces). However, they are different. In particular, they have the enormous advantage that they always represent a moduli problem of abelian varieties. Our uniformization theorem is optimal when the level structure imposed is prime to $p$,  in the sense that it extends to an integral uniformization that allows an explicit moduli interpretation of the points in the reduction modulo $p$.  

As in Drinfeld's approach, our uniformization theorem relies on a theorem on formal moduli spaces of $p$-divisible groups. In fact, the main work in proving our theorems is to establish an isomorphism of our formal moduli spaces with the moduli space of  Drinfeld. Such an isomorphism is also constructed by Scholze and Weinstein \cite{Sch}. Their construction relies on Scholze's theory of \emph{local Shimura varieties} and his \emph{integral $p$-adic Hodge theory}, as well as on results in a preliminary (unpublished) version of the present paper on  \emph{local models}. Our construction here is more direct and more elementary; it relies on the \emph{theory of displays}, cf. \cite{Zi}. We do not see any direct connection between the isomorphism in \cite{Sch} and the one constructed here. 

 Drinfeld's version of Cherednik's theorem for $F=\BQ$ has found numerous arithmetic applications, to \emph{level raising, level lowering} and \emph{bounding Selmer groups}, at the hands of Ribet, Bertolini, Darmon, Nekovar and many others, comp. also the references in the introduction of \cite{KRnew}. It is to be hoped that our direct construction can be the basis of similar such applications for general totally real fields $F$. 
 
 Our results are an expression of  the exceptional isomorphism between an inner twist of the adjoint group of $\GL_2$ and an inner twist of the adjoint group of ${\U}_2$. Just as for Poincar\'e's exceptional isomorphism of inner forms of $\PGL_2$  and special orthogonal groups of ternary quadratic forms, there is no higher rank analogue.

\subsection{Global results}\label{ss:globalres}Now let us state our global results. Let $K$ be a CM-field,  with  totally real subfield $F$. We denote the non-trivial $F$-automorphism of $K$ by $a\mapsto \bar a$. Let  $V$ be a two-dimensional $K$-vector space, equipped with an alternating $\BQ$-bilinear form $\psi\colon V\times V\to \BQ$ such that
\begin{equation}\label{introalt}
\psi(a x, y)=\psi(x, \bar a y), \quad x, y\in V, \, a\in K . 
\end{equation}
 There is a unique anti-hermitian form $\varkappa$ on $V$ such that 
 \begin{equation}\label{introback}
  \Trace_{K/\mathbb{Q}_p} a \varkappa(x, y) = \psi(ax, y), \quad x, y\in V, \, a\in K. 
\end{equation}
Conversely, the anti-hermitian form $\varkappa$ determines the alternating bilinear form $\psi$ with \eqref{introalt}. We  say that $\varkappa$ arises from $\psi$ by contraction. \index{contraction}
Recall that anti-hermitian spaces $V$ are determined up to isomorphism by their signature at the archimedean places of $F$ and their local invariants $\inv_v(V)$ at the non-archimedean places $v$ of $F$  (see \S \ref{ss:binform} for  the definition of $\inv_v(V)$). Let $w$ be  an archimedean place such that $V_w$ has signature $(1, 1)$ and such that $V_{w'}$ is definite for all archimedean places $w'\neq w$. Let us be more precise. Let
\index[NO]{ZZSA@$\Phi$} $\Phi = \Hom_{\text{$\mathbb{Q}$-Alg}}(K, \mathbb{C})$.
Let $r$ be a \emph{generalized special CM-type of rank 2, special w.r.t. $w$},\index{generalized special CM-type} i.e., a function \index[NO]{RAA@$r$}
\begin{equation}
r: \Phi\lra \BZ_{\geqslant 0}, \qquad \,\, \varphi \mapsto r_\varphi,
\end{equation}
such that $r_{\varphi}+r_{\bar{\varphi}} =2$ for all $\varphi\in \Phi$, and such that  for  the extensions $\{\varphi_0, \bar{\varphi}_0\}$ of $w$ we have  $r_{\varphi_0} = r_{\bar{\varphi}_0} = 1$ and  with  $r_\varphi\in \{0, 2\}$ for $\varphi\notin \{\varphi_0, \bar\varphi_0\}$, comp. \cite{KRnew}. Then we demand that the signature of $V_\varphi=V\otimes_{K, \varphi} \BC$ be equal to  $(r_\varphi, 2-r_\varphi)$. 

We denote the reflex field \index{reflex field} of $r$  by $E=E_r$. It is a subfield of $\ov\BQ$, the algebraic closure of $\BQ$ in $\BC$. Note  that $F$ embeds via $\varphi_0$ into $E$, and that the  archimedean place of $F$ induced by
\begin{equation}\label{introFintoE}
  F \overset{\varphi_0}{\longrightarrow} E \longrightarrow \BC
  \end{equation}
is equal to $w$.  If $F=\BQ$, then $E=F$.

Associated to these data, there is a \emph{Shimura pair} $(G, \{h\})$. \index{Shimura pair}Here    $G $ denotes the group of unitary 
similitudes of $V$, with similitude factor in $\BG_m$, an algebraic subgroup of ${\rm GSp}(V, \psi)$ over $\BQ$. For an open compact subgroup $\bK \subset G(\BA_f)$, there is a Shimura variety ${\rm Sh}_{\bK}$, with canonical model over the reflex field $E$,
whose complex points are given by 
$${\rm Sh}_{\bK}(\BC) \simeq G(\BQ) \bs [\Omega_\BR \times G(\BA_f)/{\bK}] .$$
Here  $\Omega_\BR$ is acted on by
$G(\BR)$ via the projection to ${\rm GU}(V_w)_{\rm ad}$ and a fixed
isomorphism  
${\rm GU}(V_w)_{\rm ad} \simeq \PGL_2(\BR)$.  

Consider the following moduli problem on $({\rm Sch}/E)$.  It associates to an $E$-scheme $S$ the set of isomorphism classes of tuples $(A, \iota, \bar\lambda, \bar\eta)$. Here
\begin{itemize}
\item $A$ is an abelian scheme up to isogeny of dimension $2[F:\BQ]$ over $S$.
\item  $\iota\colon O_K\to \End(A)$ is an action of $O_K$ on $A$ such that 
$$
{\rm Tr}(\iota(a)|\Lie A)=\sum\nolimits_{\varphi\in\Phi}r_\varphi\varphi(a), \quad \text{ for all } a\in O_K .
$$
\item $\bar\lambda$ is a $\BQ$-homogeneous polarization of $A$ such that its Rosati involution induces the conjugation on $K/F$. 
\item a  $\bK$-orbit   of $K$-linear similitudes $\bar\eta\colon V\otimes\BA_f \isoarrow \wh V(A)$. 
\end{itemize} 
Here the rational Tate module is equipped with its natural anti-hermitian form arising by contraction from its polarization form. 

This moduli problem is represented by a quasi-projective scheme $\CA_{\bK, E}$ which is the \emph{canonical model} of ${\rm Sh}_{\bK}$ over $E$. It is a projective scheme when the existence of $v$ as below is imposed. 

Let $p$ be a prime number and let $v$ be  a $p$-adic place of $F$ which is non-split in $K$ and such that $V_v$ is a non-split $K_v/F_v$-anti-hermitian space, i.e., $\inv_v(V)=-1$. We take  the open compact subgroup of the form  
${\bK}^\star = {\bK}^\star_p \cdot {\bK}^p$, where ${\bK}^p$ is an arbitrary open compact subgroup of  $G(\BA_f^{ p})$,  and where ${\bK}^\star_p$ 
\index[NO]{KAA@${\bK}^\star$}\index[NO]{KAB@${\bK}^\star_p$} has the following shape. Let 
\begin{equation*}
 V \otimes \BQ_p = \bigoplus\nolimits_{{\mathfrak p} \mid p} V_{\mathfrak p} \, 
\end{equation*}
be the orthogonal decomposition according to the prime ideals of $F$ over $p$. Note that the prime ideal ${\mathfrak p}_v$ corresponding to $v$ occurs as an index here. Then 
$$
G(\BQ_p)\subset \prod\nolimits_{\mathfrak p} G_{\mathfrak p}(\BQ_p)\, ,
$$
where $G_{\mathfrak p}$ denotes the group of unitary similitudes of $V_{\mathfrak p}$ with similitude factor in $\BG_{m}$. We take ${\bK}^\star_p$
  of the form 
\begin{equation}\label{formKstar}
{\bK}^\star_p=G(\BQ_p)\cap {\bK}_v{\bK}_p^{\star, v}\, ,
\end{equation}
where ${\bK}_v$ is the unique maximal compact subgroup of $G_{{\mathfrak p}_v}(\BQ_p)$, and where ${\bK}_p^{\star,v}\subset \prod\nolimits_{{\mathfrak p}\neq {\mathfrak p}_v} G_{\mathfrak p}(\BQ_p)$ is an arbitrary open compact subgroup.

Let $J$ \index[NO]{JAA@$J$} be the inner form of $G$  which is anisotropic at $w$ and quasi-split  at $v$, and which locally coincides 
with $G$ at all places $\neq v, w$ of $F$.  Then there exists an identification of the adjoint group $J_{v , \rm ad}(\BQ_p)$ with $\PGL_2(F_v)$ and an action of 
$J(\BQ)$ on $G(\BA_f)/{\bK}^\star$ (which is, however, not induced by an action of $J(\BQ)$ on $G(\BA_f)$). 

We now formulate our main theorem in the version over a $p$-adic field, cf. Corollary \ref{rigmain7}. Recall the embedding \eqref{introFintoE} of $F$ into $E$.  We choose a place $\nu$ \index[NO]{ZZMA@$\nu$}  of $E$ over $v$. Throughout the paper, we always assume\footnote{In the light of the results of Kirch \cite{Ki}, it should be possible to remove this blanket assumption.} $p\neq 2$ if $v$ is ramified in $K$. 
\begin{theorem}\label{ThmA} Let $\bK^\star=\bK_p^\star\bK^p$, where $\bK_p^\star$ is of the form \eqref{formKstar}. Let $\breve{{E}}_\nu$ be  the completion of the maximal unramified extension  of
$E_\nu$ in $\ov\BQ_p$. 
 There is  a finite abelian extension $\breve {E}_\nu^\star$ of  $\breve{{E}}_\nu$ and an
isomorphism of algebraic curves over $\breve {E}_\nu^\star$,
\begin{equation*}
{\CA}_{\bK^\star, E} \times_{\Spec\,E}\Spec\, \breve{{E}}_\nu^\star
\simeq \big(J(\BQ)\bs [{\Omega}_{F_v} \times  G(\BA_f)/ {\bK}^\star]\big)\times_{\Spec\,F_v}\Spec\,\breve{{E}}_\nu^\star
\end{equation*}

\end{theorem}
Here, as before,  ${\Omega}_{F_v} $ denotes the Drinfeld halfplane relative to the local field $F_v$, and the interpretation is as before that the scheme on the LHS is the algebraization of the rigid-analytic variety on the RHS.    If  $\bK_p$ is of the form \eqref{introplev} below, then we may take $\breve E_\nu^\star=\breve E_\nu$, cf. Theorem \ref{ThmB} below; but in general, one needs a non-trivial  extension, comp. Theorem \ref{ThmC}.  

From this theorem we deduce an analogue of Cherednik's isomorphism \eqref{Chered}, noting that any geometric connected component $X$ of ${\rm Sh}_{\bK^\star}$ is defined over the maximal abelian extension $E^{\rm ab}$ of $E$,
\begin{equation}\label{modChered}
X\otimes_{E^{\rm ab}} \breve E^{\rm ab}_\nu\simeq(\bar\Gamma\backslash\Omega_{F_v})\otimes_{F_v} \breve E^{\rm ab}_\nu.
\end{equation}
Here $ \breve E^{\rm ab}_\nu$ denotes the completion of the maximal abelian extension of $E_\nu$, and $\bar\Gamma$ is a cocompact discrete subgroup of $\PGL_2(F_v)$. Since the Cherednik Shimura datum is a central twist of $(G, \{h\})$, the geometric connected components of ${\rm Sh}_{\bK}$ can be identified with those appearing in Cherednik's theorem, so that in fact \eqref{Chered} follows from Theorem \ref{ThmA}.

By extending  the moduli problem for ${\rm Sh}_{\bK}$ integrally over $\Spec O_{E, (\mathfrak p_{\nu})}$, we obtain semi-global integral models \index{semi-global integral model} of these Shimura varieties. This gives us the possibility of formulating an `integral' version of this theorem.  Let us explain the moduli problem in question.  

We first explain the level structure. For every ${\mathfrak p}|p$, we fix a lattice $\Lambda_{\mathfrak p}$ in $V_{\mathfrak p}$. We assume that $\Lambda_{\mathfrak p}$ is a self-dual lattice  (for the alternating form $\psi$) when ${\mathfrak p}$ is either split in $K$ or ramified.  When ${\mathfrak p}$ is unramified in $K$, we assume that $\Lambda_{{\mathfrak p}}$ is selfdual when $\inv_{{\mathfrak p}}(V)=1$, and almost selfdual 
 when $\inv_{{\mathfrak p}}(V)=-1$. Let 
 \begin{equation}\label{introplev}
 \bK_p=\{g\in G(\BQ_p)\mid g\Lambda_{\mathfrak p}= \Lambda_{\mathfrak p}, \text{ for all } \mathfrak p|p \} .
 \end{equation} 
 We also fix an open compact subgroup $\bK^p\subset G(\BA_f^p)$ and set $\bK=\bK_p\bK^p$. We continue to assume that for the distinguished $p$-adic place $v$ we have $\inv_{{\mathfrak p}_v}(V)=-1$.

 We now define a functor $\CA_\bK$\index[NO]{ABA@$\mathcal{A}_{\bK}$} 
 on the category of $O_{E, ({\mathfrak p}_\nu)}$-schemes. Let $S\in({\rm Sch}/ O_{E, ({\mathfrak p}_\nu)})$. Then a point of $\CA_\bK(S)$ consists of an equivalence class of quadruples $(A, \iota, \bar\lambda, \bar\eta^p)$. Here
\begin{itemize}
\item $A$ is an abelian scheme over $S$ and $\iota\colon O_K\to \End(A)\otimes \BZ_{(p)}$ is an action of $O_K$ on $A$.
\item $\bar\lambda$ is a $\BQ$-homogeneous polarization of $A$ such that its Rosati involution induces the conjugation on $K/F$. 
\item $\bar\eta^p\colon  V\otimes\BA_f^p \isoarrow \wh V^p(A)$ is  $\bK^p$-class of $K$-linear similitudes. 
\end{itemize} 
Here the  prime-to-$p$-rational Tate module $\wh V^p(A)$
\index[NO]{VAA@$\wh V^p(A)$}  of $A$ is equipped with its natural anti-hermitian form arising by contraction from its polarization form.  Two quadruples $(A, \iota, \bar\lambda, \bar\eta^p)$ and $(A', \iota', \bar\lambda', \bar\eta^{p,\prime})$ are equivalent, if there exists an isogeny $A\to A'$ of degree prime to $p$ compatible with the remaining data. 

We impose the following conditions on  the quadruples $(A, \iota, \bar\lambda, \bar\eta^p)$. First, for   the action of $O_K$ on $\Lie A$ induced by $\iota$, we impose  the \emph{Kottwitz condition} $({\rm KC}_r)$ relative to $r$, see \eqref{signature.condition}, comp. \cite{KRnew}. In addition, we demand that this action also satisfies the \emph{Eisenstein condition} $({\rm EC}_r)$ relative to $r$. This condition is defined in section \ref{s:almostprinc}, and is a  key novelty of this paper. The condition $({\rm EC}_r)$ follows from the Kottwitz condition $({\rm KC}_r)$ when $S$ is an $E$-scheme but is quite subtle when $p$ is nilpotent in $\CO_S$. Imposing this condition ensures the flatness of the moduli scheme. 

Secondly, we demand that there exists a polarization $\lambda\in\bar\lambda$  such that, for every ${\mathfrak p}| p$,  the localization of the kernel of the polarization $\lambda$ at the place ${\mathfrak p}$  satisfies 
\begin{equation}
 \vert (\Ker\,\lambda)_{\mathfrak p}\vert= [\Lambda_{\mathfrak p}^\vee:\Lambda_{\mathfrak p}] .
 \end{equation}

Thirdly, we impose that for each $\mathfrak{p}|p$,  the {\it $r$-adjusted invariant} $\inv_{\mathfrak p}^r(A,\iota,\lambda)$  
coincides with the invariant $\inv_{\mathfrak p}(V)$ of the anti-hermitian space $V$. Here the $r$-adjusted invariant of the triple $(A,\iota,\lambda)$,  defined in \S \ref{ss:radj},  is another key novelty of this paper.  This condition is automatically satisfied when $\mathfrak{p}_v$ is the only prime ideal of $F$ over $p$. In general, this condition cuts out the open and closed part  of the moduli scheme defined by the Shimura variety. The reason for the name \emph{$r$-adjusted} is that this adjusts the definition of the invariant in \cite{KRnew}, where it was erroneously asserted that the invariant is  locally constant in families. We prove here that this local constancy indeed holds for the $r$-adjusted invariant, cf. Proposition \ref{contcorr}. 
     \begin{proposition}Let $r$ be a generalized CM-type of \emph{even} rank $n$, with associated reflex field $E$,  cf. \cite[\S 2]{KRnew}.
Let $S$ be an $O_E$-scheme. Let $(A, \iota, \lambda)$ be a CM-triple over $S$ 
which satisfies the Kottwitz condition $({\rm KC}_r)$, cf. \S \ref{ss:radj}. Let $c \in \{\pm1\}$. Then for every non-archimedean place
$v$ of $F$, the set of points $s \in S$ such that
\begin{displaymath}
\inv^r_v(A_s, \iota_s, \lambda_s) = c 
  \end{displaymath}
is open and closed in $S$. 
  \end{proposition}

Here, now, is our  main theorem in the context  of schemes over $p$-adic integer rings, cf. Theorem \ref{maint} and Corollary \ref{cortomain}.
\begin{theorem}\label{ThmB} Let $\bK=\bK_p\bK^p$, where $\bK_p$ is of the form \eqref{introplev},  and where $\bK^p$ is sufficiently small. 

\noindent (i) The functor $\CA_\bK$ is representable by a projective flat $O_{E, ({\mathfrak p}_\nu)}$-scheme  of relative dimension one,  which is the unique stable model in the sense of Deligne-Mumford \cite{DM} of its generic fiber. Its generic fiber $\CA_\bK\otimes_{O_{E, ({\mathfrak p}_\nu)}} E$ is identified with $\CA_{\bK, E}$ and its  complex fiber $\CA_\bK\otimes_{O_{E, ({\mathfrak p}_\nu)}} \BC$  with ${\rm Sh}_{\bK}$.

\smallskip

\noindent (ii) Let $\hat\CA_\bK$ be the formal completion of $\CA_\bK$ along its special fiber, which is a formal scheme over $\Spf O_{E_\nu}$. Then there exists an isomorphism of formal schemes over $\Spf O_{\breve E_\nu}$,
$$
\hat\CA_\bK\times_{{\rm Spf}\,O_{E_{\nu}}}{\rm Spf}\, O_{\breve E_\nu}  \simeq   J(\BQ)\bs\big[\big(\widehat{\Omega}_{F_v} \times_{{\rm Spf}\,O_{F_v}}{\rm Spf}\, O_{\breve E_\nu}\big) \times
 G(\BA_f)/{\bK}\big]  \, .
$$
For varying $\bK^p$, this isomorphism is compatible with the action of $G(\BA_f^p)$ through Hecke correspondences on both sides. 
The natural Weil descent datum on the LHS is given on the RHS by $(\xi, g)\mapsto (\tau_{E_\nu}(\xi), w'_r (g))$,
where $\tau_{E_\nu}$ is the natural Weil descent datum down to $O_{E_\nu}$ on $\wh\Omega_{F_v}\times_{\Spf O_{F_v}}\Spf O_{\breve{E}_\nu}$, and where $w'_r$ is a certain automorphism of  $G(\BQ_p)/\bK_p$ commuting with the action of $G(\BQ_p)$ by left translations. If the inertia index $f_{E_\nu}$ is even and if furthermore there are no prime ideals $\mathfrak p|p$  which split in $K$, then $w'_r$ is given by multiplication by $p$. 

\end{theorem} 

Here $\widehat{\Omega}_{F_v}$\index[NO]{ZZWC@$\widehat{\Omega}_{F_v}$}  denotes the formal scheme version  of ${\Omega}_{F_v}$ over ${\rm Spf}\, O_{F_v}$ due to Deligne, Drinfeld and Mumford, cf. \cite{Dr}. In section \ref{s:attpu} we give a variant of the RHS, which allows us  to express the automorphism $w'_r$ explicitly.  Theorem \ref{ThmB} is optimal in the sense that it describes explicitly the scheme $\CA_\bK$ over $O_{E_\nu}$ and its $p$-adic uniformization. 

 If we assume that there are prime ideals $\mathfrak p|p$ different from $\mathfrak p_v$, we may  pass to deeper level structures and still prove an integral version of $p$-adic uniformization. 
 Let $\bK_p^\star\subset G(\BQ_p)$ be of the form 
 \begin{equation}\label{intKPp}
 \bK_p^\star=G(\BQ_p)\cap \bK_v\bK_p^{\star, v} ,  
 \end{equation}
 where  $\bK_v$ is the stabilizer of $\Lambda_{\mathfrak{p}_v}$, and where $\bK_p^{\star, v}$ is an arbitrary open compact subgroup of $G^{v}(\BQ_p)=\prod\nolimits_{\mathfrak{p}\neq{\mathfrak{p}_v}}G_{\mathfrak{p}}(\BQ_p)$.    The system of such subgroups is stable under conjugation with elements of $G(\BQ_p)$.  For such subgroups, we have the following version of our main theorem, cf. Corollary \ref{genKPschemes}. 
\begin{theorem}\label{ThmC}
Let $\bK^\star=\bK_p^\star\bK^p$,  for a choice of  sufficiently small  ${\bf K}^p\subset G(\BA_f^p)$, where $\bK_p^\star$ is of the form \eqref{intKPp}.  There 
  exists a normal scheme $\mathcal{A}^{\star}_{\mathbf{K}^{\star}}$ over $\Spec O_{\breve E_\nu}$ 
  such that for the $p$-adic completion of this scheme there is an isomorphism 
  \begin{displaymath}
\hat{\mathcal{A}}^{\star}_{\mathbf{K}^{\star}} \simeq
J(\mathbb{Q}) \backslash [(\hat{\Omega}_{F_{v}} \times_{\Spf O_{F_{v}}}
\Spf O_{\breve{E}_{\nu}})\times  G^{{v}}(\mathbb{Q}_p)/\mathbf{K}_{p}^{\star, {v}}
\times G(\mathbb{A}_f^p)/\mathbf{K}^p].
  \end{displaymath}
 For varying $\mathbf{K}^{\star}$, these schemes form a tower with an action of the group
$G(\mathbb{Q}_p) \times G(\mathbb{A}_f^{p})$, where the action of $G(\BQ_p)$ factors through  $G(\mathbb{Q}_p) \rightarrow G^{{v}}(\mathbb{Q}_p)$. The isomorphism of formal schemes 
is compatible with these actions. 
  
The general fiber of $\mathcal{A}^{\star}_{\mathbf{K}^{\star}}$ is a Galois twist of
 $\mathcal{A}_{ \mathbf{K}^{\star}, E}\times_{\Spec E}\Spec \breve E_\nu$ by an abelian character
$\chi_0^{\rm h}$. The Galois twist
respects the Hecke operators.  
\end{theorem}
The scheme $\mathcal{A}^{\star}_{\mathbf{K}^{\star}}$ represents a moduli problem of abelian varieties with additional structure over $O_{\breve E_\nu}$, cf.  section \ref{ss:deeper}.  We refer to section \ref{ss:TTm} for the explicit determination of $\chi_0^{\rm h}$. 

It should be pointed out that Theorem \ref{ThmC} is not optimal since we cannot describe  the descent to $E_\nu$. Also, when $v$ is ramified in $K$, we can only give the character $\chi_0^{\rm h}$ explicitly after restricting to a subgroup of index $2$.  This is in contrast to Theorem \ref{ThmB}. The deeper reason for this deficiency lies in the fact that the natural context for Theorem \ref{ThmC} is the class of Shimura varieties appearing in \cite{RSZ}. Let $\Psi\subset \Phi$ be a CM-type for $K/F$ such that 
\begin{equation}
\Psi\cap (\Phi\setminus \{\varphi_0,\bar\varphi_0\})=\{\varphi\in\Phi \setminus \{\varphi_0,\bar\varphi_0\}\mid r_\varphi=2\}.
\end{equation}
There are two possibilities for $\Psi$. Let $E_\Psi$ be the reflex field of $\Psi$ and let $\tilde E$ be the composite of $E_\Psi$ and $E=E_r$. Then $\tilde E$ is an extension of degree one or two of $E_r$. Associated to $(V, \psi, \Psi)$, there is a finite number of Shimura varieties ${\rm Sh}_{\tilde\bK}(\tilde G, \{\tilde h\})$ with reflex field $\tilde E$, cf. \cite[\S 3]{RSZ}. Here $\tilde G$ maps surjectively to $G$ with kernel a central torus, hence the Shimura varieties ${\rm Sh}_{\tilde\bK}(\tilde G, \{\tilde h\})$ are central twists of the Shimura variety ${\rm Sh}_{\bK}(G, \{ h\})$.  Each one represents a moduli problem on $({\rm Sch}/\tilde E)$. In a sequel to this paper, we will construct  semi-global integral models of these Shimura varieties over $\Spec O_{\tilde E, (\tilde\nu)}$. These are described by  moduli problems of abelian varieties on $({\rm Sch}/O_{\tilde E, (\tilde\nu)})$ and admit $p$-adic uniformization in the strong sense of Theorem \ref{ThmB}, when the congruence condition on the open compact subgroup $\tilde\bK$ is prime to the chosen place $v$.  The trade-off in comparison with our Shimura variety is that the corresponding reflex field $\tilde E$ is larger than the reflex field $E$ of our Shimura variety (which, in turn, is larger than the reflex field $F$ of Cherednik's Shimura variety).

 Both Theorems \ref{ThmA} and \ref{ThmB}  are proved in \cite{KRnew} when $F_v=\BQ_p$. Most of the work in \cite{KRnew} was local, and an essential ingredient was the {\it alternative moduli interpretation} of the Drinfeld halfplane in \cite{KRalt}.  Once this is accomplished, the proof of the global theorems follows in a relatively straightforward way from the general non-archimedean uniformization theory of \cite[Chap. 6]{RZ}. The same is true here. 
In \cite{KRnew},  we expressed the hope that it might be possible to eliminate the strong limitation $F_v=\BQ_p$ made there, and this hope 
is achieved in the present paper.  As explained in \cite{KRnew}, the main issue is the contrast between the condition 
on the action of $O_{F_v}$ on the Lie algebras 
of the $p$-divisible groups in the local moduli problems.  On the one hand, for the moduli problem represented by the Drinfeld half-plane 
$\wh \Omega_{F_v}$, the action of $O_{F_v}$ on the Lie algebra is required to be strict, i.e., to factor through the structure morphism of the base scheme $S$. 
On the other hand, in the global moduli problem, the Lie algebras of the relevant abelian schemes are often free $O_F\otimes_\BZ \CO_S$-modules. 
The main results of the present paper, and in particular the contracting functor defined in section 4, provide the  bridge between the two types 
of moduli problems. 
 
 \subsection{Local results} Let us now formulate our local results, referring to section \ref{s:almostprinc} for  more details and more explanations of some terms used here. Let $p$ be a prime number, and let $F$ be a finite extension of degree $d=[F:\BQ_p]$ of $\BQ_p$ and let $K/F$ be a quadratic extension. Let $\Phi = \Hom_{\text{$\mathbb{Q}_p$-Alg}}(K, \bar\BQ_p)$, and fix a pair $\{\varphi_0, \bar\varphi_0\}$ of conjugate elements in $\Phi$. Here $\bar\varphi_0(a)=\varphi_0(\bar a)$. 
Let $r$ be a \emph{local  CM-type of rank 2} which is special w.r.t $\{\varphi_0, \bar\varphi_0\}$, i.e., a function 
\begin{equation}
r: \Phi\lra \BZ_{\geqslant 0}, \qquad \,\, \varphi \mapsto r_\varphi,
\end{equation}
such that $r_{\varphi}+r_{\bar{\varphi}} =2$ for all $\varphi\in \Phi$, and   $r_{\varphi_0} = r_{\bar{\varphi}_0} = 1$ and $r_\varphi\in \{0, 2\}$ for $\varphi\notin \{\varphi_0, \bar\varphi_0\}$, comp. \cite{KRnew}. We denote the reflex field of $r$  by $E$. It is a subfield of $\bar\BQ_p$.

For an $O_E$-scheme $S$, we consider triples $(X, \iota, \lambda)$, where $X$ is a $p$-divisible group of height $4d$ and dimension $2d$ over $S$, where $\iota\colon O_K\to\End(X)$ is an action of $O_K$ on $X$, and where $\lambda\colon X\to X^\vee$ is a polarization of $X$ such that its Rosati involution induces on $O_K$ the conjugation involution over $O_F$. We impose the Kottwitz condition $({\rm KC}_r)$ and the Eisenstein conditions $({\rm EC}_r)$ on the action of $O_K$ on $\Lie X$. Furthermore, we assume that $\lambda$ is a \emph{principal} polarization if $K/F$ is ramified, and that $\lambda$ is an \emph{almost principal} polarization if $K/F$ is unramified. 

We fix such a triple $(\BX, \iota_\BX, \lambda_\BX)$ over the algebraic closure $\bar k$ of the residue field $\kappa_E$ of $E$, and refer to it as a \emph{framing object}. When $K/F$ is unramified, then any two such  triples are isogenous by  an $O_K$-isogeny of height zero which preserves the polarizations.  The same is true when $K/F$ is ramified, provided we impose that the $r$-adjusted invariant $\inv^r(\BX, \iota_\BX, \lambda_\BX)$ is $-1$ (this last condition is automatic when $K/F$ is unramified). In either case, the group $J(\BQ_p)$ of $O_K$-self-quasi-isogenies of $(\BX, \iota_\BX)$ preserving the polarization $\lambda_\BX$ up to a factor in $\BQ_p^\times$ can be identified with the group of unitary similitudes, with similitude factor in $\BQ_p^\times$, of a \emph{split}  $K/F$-anti-hermitian space of dimension $2$. Let $J^1(\BQ_p)$ denote the special unitary group.

We consider the Rapoport-Zink space  $\CM_r$ over $\Spf\, O_{\breve E}$ representing the functor on   $ ({\rm Sch}/\Spf O_{\breve E})$ which associates to $S\in  ({\rm Sch}/\Spf O_{\breve E})$  the set of isomorphism classes of $4$-tuples $(X, \iota, \lambda, \rho)$, where $(X, \iota, \lambda)$ is as above, and where $\rho$ is a \emph{framing} of height zero, with framing object $(\BX, \iota_\BX, \lambda_\BX)$. Our main local result may now be formulated as follows. We fix an isomorphism $J^1(\BQ_p)\simeq \SL_2(F)$. 
\begin{theorem}\label{mainlocalintro}  
  The RZ-space  $ \CM_r$ is isomorphic to  $\wh{\Omega}_F \hat{\otimes}_{O_F, \varphi_0} O_{\breve E}$. More precisely, there exists a unique isomorphism of formal schemes 
  $$
  \CM_r\simeq \wh{\Omega}_F {\times}_{\Spf O_F} \Spf O_{\breve E} ,
  $$
  which is equivariant with respect to the fixed identification $J^1(\BQ_p)\simeq \SL_2(F)$. In particular, $\CM_r$ is flat over $\Spf O_{\breve E}$ with semi-stable reduction.
\end{theorem}
 
 It is more honest to formulate this theorem as follows.  Let $\CM$ be the \emph{relative} RZ-space over $\Spf O_{\breve F}$ from \cite{KRalt}. It parametrizes tuples $(X', \iota', \lambda', \rho')$, where $X'$ is a \emph{strict} formal $O_F$-module of \emph{relative} height $4$ and dimension $2$, and where $\iota'$ is an $O_K$-action on $X$ which is of signature $(1, 1)$, and where $\lambda'$ is a \emph{relative} polarization compatible with $\iota'$, which is principal if $K/F$ is ramified and almost principal if $K/F$ is unramified. Also, $\rho'$ is a framing of height zero with a suitable framing object $(\BX', \iota'_{\BX'}, \lambda'_{\BX'})$. We fix an extension $\breve \varphi_0\colon O_{\breve F}\to O_{\breve E}$ of $\varphi_0\colon O_F\to O_E$. Then our main local result is the construction of a \emph{contracting functor}  
 \begin{equation}\label{contrintro}
 \CM_r\to \CM\times_{\Spf O_{\breve F}}\Spf O_{\breve E}
 \end{equation}
 and the proof that it induces an isomorphism of formal schemes over $\Spf O_{\breve E}$. The construction of the contracting functor is another key novelty of this paper. Theorem \ref{mainlocalintro} then follows  by combining \eqref{contrintro} with the  alternative interpretation of the Drinfeld halfplane of  \cite{KRalt}, which yields an isomorphism $\CM\simeq \wh{\Omega}_F \hat{\otimes}_{O_F} O_{\breve F}$. 
 
 Our construction of the contracting functor is based on the theory of displays. Let $R$ be a $p$-adic ring, and let $W(R)$ be its ring of Witt vectors. Displays over $R$ are certain modules over $W(R)$ with additional structures. Under suitable hypotheses, the category of $p$-divisible groups over $R$ is equivalent to the category of displays over $R$. The contracting functor is the composition of two functors.  The first functor  associates to a tuple $(X, \iota, \lambda)$ as above a new tuple $(\tilde X, \tilde\iota, \tilde\lambda)$, where $\tilde X$ is a $p$-divisible group of height $4d$ and dimension $2$, where  $\tilde\iota$ is an $O_K$-action such that its restriction to $O_F$ is strict and which is of signature $(1, 1)$, and where $\tilde\lambda$ is a polarization 
 \emph{with values in the Lubin-Tate group}   compatible with $\tilde\iota$. We call this a polarization in the sense of Faltings. Note that because of the values for the height and the dimension of $\tilde X$, there cannot exist a polarization in the usual sense on $\tilde X$. The  second functor is 
  the {\it Ahsendorf functor} from \cite{ACZ}. It   associates to  $(\tilde X, \tilde\iota, \tilde\lambda)$ a relative tuple $(X', \iota', \lambda')$ as above. The Ahsendorf functor is the analogue for displays of the Drinfeld functor which associates to a Cartier module of a $p$-divisible group with strict $O_F$-action its \emph{relative} Cartier module, cf. \cite[\S 2]{Dr}.  We also use the theory of displays to give a new proof of (a slight refinement of) the alternative interpretation of the Drinfeld halfplane, which is the third proof after the original proof in \cite{KRalt} and the proof of Kirch \cite{Ki}. 
  
 Let us formulate the main contribution of this paper to the theory of displays, cf. Theorem \ref{Adorf1t} and Theorem \ref{LTD10p}. It compares displays for the \emph{Witt frame} $\CW(R)$  with displays for the \emph{relative Witt frame} $\CW_{O_F}(R)$, comp. Definition \ref{Rah2d}. We give a simpler construction of the Ahsendorf functor and use this to prove the following theorem. 
\begin{theorem}\label{maindisplayintroduction}
  Let $R$ be an $O_F$-algebra such that $p$ is nilpotent in $R$. 
  The Ahsendorf functor is a functor 
  \begin{displaymath}
  \mathfrak{A}_{O_F/\BZ_p, R}:
\left(
  \begin{array}{l}
  \text{  $\mathcal{W}_{}(R)$-displays}\\
    \text{with strict $O_F$-action} 
    \end{array}
\right)
  \; \longrightarrow \; \Big(\text{ $ \mathcal{W}_{O_F}(R)$-displays}\Big) .  
  \end{displaymath}
  It induces an equivalence of categories
  \begin{displaymath}
  \mathfrak{A}_{O_F/\BZ_p, R}:
\left(
  \begin{array}{l}
  \text{nilpotent   $\mathcal{W}_{}(R)$-displays}\\
    \text{with strict $O_F$-action} 
    \end{array}
\right)
  \; \longrightarrow \; \Big(\text{nilpotent $ \mathcal{W}_{O_F}(R)$-displays}\Big) .  
  \end{displaymath}
 Let $\mathcal{P}_1$ and $\mathcal{P}_2$ be $\CW(R)$-displays over $R$ with a strict
  $O_F$-action. We denote by $\mathcal{P}_{1, {\rm a}}$ and $\mathcal{P}_{2, {\rm a}}$ their images
by the Ahsendorf functor $\mathfrak{A}_{O_F/\mathbb{Z}_p, R}$. Then there is a natural homomorphism between groups of bilinear forms of displays, 
\begin{equation*}
  \Bil_{\text{$O_F${\rm-displays}}}(\mathcal{P}_1 \times \mathcal{P}_2, \mathcal{L}_R)
\longrightarrow  
  \Bil_{\text{$\CW_{O_F}(R)${\rm-displays}}}(\mathcal{P}_{1, {\rm a}} \times \mathcal{P}_{2, {\rm a}},
  \mathcal{P}_{m, \mathcal{W}_{O_F}(R)}(\pi^{ef}/p^f)) .
  \end{equation*}
  If  the dual $(\mathcal{P}_{1, {\rm a}})^{\vee}$ of  $\mathcal{P}_{1, {\rm a}}$ and
  $\mathcal{P}_{2, {\rm a}}$ are nilpotent $\mathcal{W}_{O_F}(R)$-displays, then this homomorphism is an isomorphism. 
  
  Here $\CL_R$ denotes the Lubin-Tate $\CW(R)$-display, and $  \mathcal{P}_{m, \mathcal{W}_{O_F}(R)}(\pi^{ef}/p^f)$  the twist of the multiplicative $\mathcal{W}_{O_F}(R)$-display   by the unit $\xi=\pi_F^{ef}/p^f$. 

\end{theorem}  
 Let us now put the local results of this paper in perspective. We  address in our special case the general problem of identifying a \emph{basic} Rapoport-Zink space associated to the pair $(G, \{\mu\})$ with a twist of the   basic Rapoport-Zink space associated to the pair $(G, \{\mu'\})$, where $\mu'$ differs from $\mu$ by a central character, cf. the Introduction of \cite{RZdrin}. This problem is also addressed by Scholze in \cite[Chap. 23]{Sch}, in both the case considered here and in the \emph{fake Drinfeld case} of \cite{RZdrin}. As mentioned above,  Scholze's proof uses in an essential way our formulation of the local moduli problem, via  the theory of {\it local models} (and hence implicitly  the linear algebra lemma \cite[Lem. 4.9]{RZdrin}).  One of the main reasons that we are successful in constructing the contracting functor  in the case treated here is that here we are able to develop a good understanding of the Kottwitz condition $({\rm KC}_r)$, even in unequal characteristic. Our failure to do the same in the fake Drinfeld case is the essential reason that in \cite{RZdrin} we only succeeded in defining the contracting functor in the special fiber.  
 The contracting functor is an expression of the exceptional isomorphism between the quasi-split special unitary group in two variables and the special linear group in two variables. We restricted ourselves here  to the case of curves; it would have been possible to  prove a higher-dimensional version where the uniformizing space ${\Omega}_{F_v} $  is replaced by a product of such spaces, comp. \cite{KRnew} and \cite[\S 6]{RZ}.

\subsection{Layout of the paper}We now explain the layout of the paper. The whole paper, with the exception of section \ref{s:attpu}, is devoted to the local theory.  In section \ref{s:almostprinc} we explain in detail the definition of the formal moduli spaces of (polarized) $p$-divisible formal groups, including the Kottwitz conditions relevant here and the Eisenstein conditions; in particular, Subsection \ref{ss:defofform} contains the detailed statements of our main local results. Section \ref{s:relativeDT} summarizes the relevant facts on relative Dieudonn\'e theory and relative display theory. The most important fact proved in this section is the relation established by the Ahsendorf functor between the \emph{Lubin-Tate display} and the \emph{relative multiplicative display}.  In section \ref{s:tcf} we first consider the relation between the Kottwitz condition and the Eisenstein condition; this is used in the rest of the section to construct the contracting functor. More precisely, we first consider the first step in its construction which we call the \emph{pre-contracting functor}, cf. above. After this, we complete the second step in the case of a \emph{special} generalized CM-type. In the final subsection of section \ref{s:tcf}, we consider the second step in the case of a \emph{banal} generalized CM-type. Section \ref{s:altmodprob} is devoted to an alternative proof of the main result of \cite{KRalt}, based on the theory of displays. In section \ref{s:modformal} we prove the main local results, namely Theorem \ref{mainlocalintro} and its banal counterpart. In the appendix, section \ref{s:appA}, we give the  correct version of the sign factor  of \cite{KRnew}  by defining the {\it adjusted invariant} of a CM-triple of generalized CM-type $r$ of even rank $n$, and  investigate its behaviour under the contracting functor. Section \ref{s:attpu} deduces the global results from the local theory.

\subsection{Acknowledgements}  
We are grateful to P.~Scholze for very helpful discussions. In particular, he helped us locate the mistake in the definition of the invariant of a CM-triple (we had discovered a discrepancy caused by this mistake in 2014, but it took us 10 months to locate the error (e-mail from P.~Scholze of 5 Oct. 2015)). We also thank the anonymous referee of \cite{RSZ} whose report gave us the idea of defining the $p$-adic \'etale sheaf in Corollary \ref{padicetale}.  

We also acknowledge the hospitality of the MSRI during the fall of 2014 when the work on this paper was begun.  We also thank the referee for helpful remarks.

\subsection{ Notation}  
\noindent $\bullet$ If $R$ and $R'$ are $\BZ_p$-algebras, we often write $R\otimes R'$ for $R\otimes_{\BZ_p} R'$. Also, we often write $X\otimes_A B$ for $X\times_{\Spec A}\Spec B$.

\noindent $\bullet$  If $F$ is a finite extension of $\BQ_p$, we write $\breve F$ \index[NO]{FAA@$\breve F$}  for the completion of a maximal unramified extension, and $F^t$ for the maximal subfield unramified over $\BQ_p$. We write $d=e f$, where $d=[F:\BQ_p]$ and $f=[F^t:\BQ_p]$ and $e=[F:F^t]$. We denote by $O_F$, resp. $O_{F^t}$, resp. $O_{\breve F}$ the rings of integers. 

\noindent $\bullet$  Let $V$ be an $\BC/\BR$-anti-hermitian vector space.  The signature of $V$ is $(a, b)$ if the anti-hermitian form is equivalent to ${\rm diag}({\bf i}^{(a)}, {\bf i}^{(b)})$, where $\bf i$ is the imaginary unit. 

\noindent $\bullet$  Let $F$ be a finite extension of $\BQ_p$ and let $K/F$ be a quadratic extension.  Let  $V$ be a  $K$-vector space, equipped with an alternating $\BQ_p$-bilinear form $\psi\colon V\times V\to \BQ_p$ satisfying \eqref{introalt}. Let $\Lambda$ be a $O_K$-lattice in $V$. Then the dual $O_K$-lattice is $\Lambda^\vee=\{x\in V\mid \psi(x, y)\in \BZ_p \text{ for all } y\in \Lambda\}$. \index[NO]{ZZKA@$\Lambda^\vee$} The lattice $\Lambda$ is called \index{self-dual} self-dual if $\Lambda=\Lambda^\vee$; it is called almost self-dual  if $\Lambda$ is contained in $\Lambda^\vee$ with colength one. \index{almost self-dual}

\noindent $\bullet$  If $O$ is a discrete valuation ring with uniformizer $\pi$, we write ${\rm Nilp}_{O}$\index[NO]{NAA@${\rm Nilp}_{O}$}  for the category of $O$-algebras $R$ such that $\pi$ is  locally on $\Spec R$ nilpotent. Similarly, we denote by $({\rm Sch}/\Spf O)$ the category of $O$-schemes such that $\pi\CO_S$ is a locally nilpotent ideal sheaf. 

\noindent $\bullet$  Given modules $M$ and $N$ over a ring $R$, we write $M\subset^r N$ to indicate that $M$ is an $R$-submodule of $N$ of finite colength $r$. 

\medskip

{\bf Warning.} It is customary to denote a finite extension of $\BQ_p$ and the Frobenius by the same symbol $F$. This should not lead to confusions.

\section{Main local statements}\label{s:almostprinc}

In this section we formulate our main results in the local theory. We
fix a prime number $p$ and an algebraic closure $\ov \BQ_p$  
of $\BQ_p$.   Here, as in the rest of the paper, we assume that $p\neq 2$. Let $F$ be a finite field extension of $\BQ_p$,  with residue class field $\kappa_F$. 
 We set $d=[F:\BQ_p]$, 
 $f = [\kappa_F : \mathbb{F}_p]$ and define $e$ through $d = ef$.
 \index[NO]{DAA@$d = ef$} \index[NO]{EAA@$e$} \index[NO]{FAB@$f$} 
We  let $K/F$ be  an \'etale algebra of degree $2$. We
denote the non-trivial automorphism of $\Gal(K/F)$ by $a \mapsto
\bar{a}$.\index[NO]{AAD@$a\mapsto\bar{a}$} 

In the case where $K/F$ is a ramified extension of local fields
(ramified case) we choose a prime element $\Pi \in O_K$ \index[NO]{ZZOA@$\Pi$} 
such that
$\bar{\Pi} = - \Pi$. Then $\pi = -\Pi^2$ is a prime element 
of $F$. In the case where $K/F$ is unramified extension of local
fields (unramified case) or $K \cong F \times F$ (split case) we
choose a prime element $\pi \in F$ and we set $\Pi = \pi$. 

Let $\Phi=\Phi_K=\Hom_{\text{$\BQ_p$-Alg}}(K, \ov\BQ_p)$ be the set of
algebra homomorphisms.   

\subsection{Special and banal local CM-types}\label{ss:specialbanal}  
Let $r$  be  a generalized local CM-type of rank $2$ (relative to
$K/F$) in the sense of \cite[section 5]{KRnew}, i.e.,  a function 
\begin{equation}
r: \Phi\lra \BZ_{\geqslant 0}, \qquad \,\, \varphi \mapsto r_\varphi,
\end{equation}
such that $r_{\varphi}+r_{\ov\varphi} =2$ for all $\varphi\in \Phi$.  Here $\ov\varphi(a) = \varphi(\ov a)$, 
where $a\mapsto \ov a$ is the non-trivial automorphism of $K$ over $F$. 
The corresponding reflex field $E=E(r)$ is the subfield of $\ov\BQ_p$ fixed by 
$$\Gal(\ov\BQ_p/E):=\{ \tau\in \Gal(\ov\BQ_p/\BQ_p)\mid r_{\tau\varphi} = r_\varphi, \ \forall \varphi\}.$$
Let  $O_E$ be the ring of integers of $E$.

When we fix an embedding $\varphi_0:F\to \ov \BQ_p$, we denote by $\varphi_0, \ov\varphi_0$ the two extensions of $\varphi_0$ to $K$ (by abuse of notation).

\begin{definition}\label{def:special} 
A local CM-type\index{special local CM-type} $r$ of rank $2$  is called {\it special} relative to the choice of embedding $\varphi_0: F\to \ov\BQ_p$ if $K/F$ is a field extension and 
$$
r_{\varphi_0}=r_{\ov\varphi_0}=1, \text{ and } r_\varphi\in \{ 0, 2 \}, \quad \text { for all } \varphi\in \Phi\setminus \{\varphi_0, \ov\varphi_0 \} . 
$$
It is called {\it banal non-split}\index{generalized  local CM-type, special or banal (split or non-split)} if $K/F$ is a field extension and $ r_\varphi\in \{ 0, 2 \}$,  for all $\varphi\in \Phi$. It is called {\it banal split} if $K\simeq F\oplus F$ and $ r_\varphi\in \{ 0, 2 \}$,  for all $\varphi\in \Phi$.
\end{definition}

From now on, we will assume  $r$ to be  either special (relative to a fixed choice of $\varphi_0$) or banal (non-split or split). We will  consider $p$-divisible groups $X$ with an action of $O_K$ over $O_E$-schemes $S$. We will want to impose certain conditions on the induced action of $O_K$ on $\Lie X$.  

\subsection{The Kottwitz and the Eisenstein conditions}\label{ss:kotteis}
Let $S$ be an $O_E$-scheme, and let $\CL$ be a locally free
$\CO_S$-module, equipped with an action
\begin{displaymath}
\iota: O_K \longrightarrow \End_{\mathcal{O}_S} \mathcal{L}.  
  \end{displaymath}
of $O_K$.

We say that $(\mathcal{L}, \iota)$ satisfies the Kottwitz condition\index{Kottwitz condition}
$({\rm KC}_r)$\index[NO]{KAC@$({\rm KC}_r)$} relative to $r$ if the identity of polynomials with
coefficients in $\CO_S$ holds    
\begin{equation}\label{signature.condition}
  {\rm char} (T, \iota (a) \vert \CL) = i\big(\underset{\varphi \in
    \Phi} \prod (T-\varphi(a))^{r_\varphi}\big),\quad \text{ for all
    $a\in O_K$ }, 
\end{equation}
where $i:O_E\longrightarrow \mathcal O_S$ is the structure homomorphism
(compare \cite{RZdrin}).  

We denote by $F^t \subset F$ \index[NO]{FAC@$F^t$} the maximal subextension which is
unramified over $\mathbb{Q}_p$.  We similarly define $K^t\subset K$ when $K$ is a field;  in the split case $K \cong F \times F$ we set
$K^t = F^t \times F^t$. We set 
$\Psi =\Psi_K= \Hom_{\text{$\mathbb{Q}_p$-Alg}}(K^{t},\ov{\mathbb{Q}}_{p})$.
\index[NO]{ZZUA@$\Psi_K$} 
We call $\psi \in \Psi$ \emph{banal} if $r_{\varphi} \in \{ 0, 2 \}$ for each
$\varphi \in \Phi$ such that $\varphi \mid \psi$. If this is not the
case we call $\psi$ \emph{special}. We use the notation
\index[NO]{ZZSB@$\Phi_\psi$} 
\begin{equation}\label{notPhipsi}
\Phi_\psi=\{\varphi\in\Phi\mid \varphi_{|K^t}=\psi \}, \quad \psi\in\Psi .
\end{equation}
 The subfield $\psi(K^t)$ of $\bar\BQ_p$ is unramified over $\BQ_p$, and hence is normal and independent of $\psi\in\Psi$. 
We denote by $E'$  the compositum of $E$ with $\psi(K^t)$. Note that  $E'/E$ is an unramified extension
of local fields. 

Let $S$ be an $O_E$-scheme. Let $\alpha: S \longrightarrow \Spec O_{E'}$ 
be a morphism of $O_E$-schemes. Then  $\alpha$ gives rise to an
isomorphism of $O_{K^t}\otimes_{\BZ_p}\CO_S$ algebras 
\begin{equation}\label{unrdecomp}
  O_{K^t}\otimes_{\BZ_p}\CO_S=\bigoplus\nolimits_{\psi\in\Psi}\CO_S\, ,
\end{equation}
where the action of $O_{K^t}$ on the $\psi$-th factor is via
$\psi$. Hence for a locally free $\CO_S$-module $\CL$ with action by
$O_{K}$, we obtain a decomposition 
into locally free $\CO_S$-modules,  
\begin{equation}\label{eigensp}
  \CL=\bigoplus\nolimits_{\psi\in\Psi} \CL\!_\psi \,.
\end{equation}
If $(\mathcal{L}, \iota)$ satisfies the Kottwitz condition we obtain
from (\ref{signature.condition}) applied to $a \in O_{K^t}$ that
\begin{equation}\label{rank.condition}
  \rank\, \CL\!_\psi = \sum_{\varphi \in\Phi_\psi} r_{\varphi} 
\end{equation}
We say that $(\mathcal{L}, \iota)$ satisfies the rank condition\index{rank condition} 
$({\rm RC}_r)$, if (\ref{rank.condition}) is satisfied for all $\psi$.
\index[NO]{RAB@$({\rm RC}_r)$}The rank
condition does not depend on the $\alpha$ chosen above because a second
$\alpha'$ differs from $\alpha$ by an automorphism of $E'$ over $E$
if $S$ is connected. If there is no $\alpha$ we use base change
$\Spec O_{E'} \times_{\Spec O_E} S \longrightarrow S$ to define the
condition $({\rm RC}_r)$. This agrees with the old definition if $\alpha$
exists.

We consider a pair $(\mathcal{L}, \iota)$ 
that  satisfies $({\rm RC}_r)$. Then we will define the  Eisenstein
condition $({\rm EC}_r)$\index[NO]{EAB@$({\rm EC}_r)$} (this definition is analogous to
\cite[section 2]{RZdrin}, but different). 
We introduce the notation
\begin{equation}\label{abpsi1e}
  \begin{aligned}
    A_\psi &= \{\varphi: K \longrightarrow \ov{\mathbb{Q}}_p \;|\;
    \varphi_{|K^{t}} = \psi, \; \text{and} \; r_{\varphi} = 2  \}\\
    B_{\psi} &= \{\varphi: K \longrightarrow \ov{\mathbb{Q}}_p \;|\;
    \varphi_{|K^{t}} = \psi, \; \text{and} \; r_{\varphi} = 0  \}.
  \end{aligned}
\end{equation} 
\index[NO]{AAA@$A_\psi, B_\psi$}\index[NO]{AAB@$a_\psi, b_\psi $}
\index[NO]{ZZEA@$\epsilon_\psi$} We note that under the action of the non-trivial automorphism of $K/F$, 
\begin{equation}\label{underconj}
   \bar{A}_\psi=B_{\ov\psi} . 
\end{equation}
Also, let $a_\psi=|A_\psi|$ and $b_\psi=|B_\psi|$.

With this notation we may rewrite the rank condition $({\rm RC}_r)$
\begin{equation}
  \rank\, \CL\!_\psi =2a_\psi \,  +\epsilon_{ \psi} \,,
\end{equation}
where 
\begin{displaymath}
  \epsilon_\psi=
  \begin{cases}
    \begin{aligned}
      0,  &\text{ if $\psi$ is banal}\\
      1  , &\text{ if $\psi$ is special and $K/F$ is unramified}\\
      2,  &\text{ if $\psi$ is special and $K/F$ is ramified.}
    \end{aligned}
  \end{cases}
\end{displaymath}
In the case where $K/F$ is ramified we have $K^t=F^t$, $[K:K^t]=2e$, and for each $\psi \in \Psi$
\begin{displaymath} 
  \varphi \mid \psi \quad \Rightarrow \bar{\varphi} \mid \psi. 
\end{displaymath}
Therefore $a_\psi=b_\psi$, and the rank condition reads, in the ramified case,
\begin{displaymath}
  \rank \mathcal{L}_{\psi} = 2e, 
\end{displaymath}
regardless of whether $r$ is banal or not.

Consider the  Eisenstein
polynomial $\mathbf{E}(T)$ of $\Pi$ in $O_{K^t}[T]$. We consider the image
$\mathbf{E}_\psi(T)$ of $\mathbf{E}(T)$ in $\ov\BQ_p[T]$ under $\psi$, for
$\psi\in\Psi$. In 
$\ov\BQ_p[T]$ this has a decomposition into linear factors,
\index[NO]{EAC@$\mathbf{E}_\psi(T)$}
\begin{equation}\label{decoveralgcl}
  \mathbf{E}_\psi(T) =
  \prod_{\varphi\in\Phi_\psi} (T-\varphi(\Pi)).   
\end{equation}
We define
\begin{equation}\label{Eisen1e}
  \mathbf{E}_{A_\psi}(T)=\prod_{\varphi\in
    A_\psi}(T-\varphi(\Pi)), \quad \mathbf{E}_{B_\psi}(T)=\prod_{\varphi\in
    B_\psi}(T-\varphi(\Pi)) .
\end{equation}\index[NO]{EAD@$\mathbf{E}_{A_\psi}$, $\mathbf{E}_{B_\psi}$}
The action of $\Gal(\ov\BQ_p/E')$ stabilizes the
corresponding subsets in the index set on the right hand sides of
\eqref{decoveralgcl} and \eqref{Eisen1e}. Therefore all three
polynomials lie in $O_{E'}[T]$.  

If $r$ is special we fix an embedding
$\varphi_0: K \longrightarrow \bar{\mathbb{Q}}_p$ such that
$r_{\varphi_0} = 1$. We denote by $\psi_0$ the restriction of
$\varphi_0$ to $K^{t}$. In the ramified case we have
$\psi_0 = \bar{\psi_0}$ and in the unramified case
$\psi_0 \neq \bar{\psi_0}$. 

We define $\mathbf{S}_{\psi}$\index[NO]{SAA@$\mathbf{S}_{\psi}$}  by the following factorization in
$O_{E'}[T]$,
\begin{equation}\label{Qpsi1e}
  \mathbf{E}_\psi(T) =
  \mathbf{S}_\psi(T) \cdot \mathbf{E}_{A_\psi}(T)\cdot
  \mathbf{E}_{B_{\psi}}(T). 
\end{equation}
Hence 
\begin{displaymath}
  \mathbf{S}_\psi(T)= 
  \begin{cases}
    \begin{array}{ll}
      1 , &\text{if $\psi$ is banal}\\
      (T-\varphi_0(\Pi)) (T-\ov\varphi_0(\Pi)),
      &\text{if $\psi=\psi_0$ and $K/F$ is ramified}\\ 
      T-\varphi_0(\Pi) ,  &\text{if $\psi=\psi_0$  and $K/F$ is unramified}\\
      T-\ov\varphi_0(\Pi) ,  &\text{if $\psi=\ov\psi_0$  and $K/F$ is
        unramified}.\\ 
    \end{array}
  \end{cases}
\end{displaymath}

Now using the structure morphism $O_{E'}\to \CO_S$, each of the three factors in \eqref{Qpsi1e}, when evaluated on $\Pi$, 
defines an endomorphism   of the
$\CO_S$-module $\CL_\psi$. These endomorphisms are denoted by  $ \mathbf{E}_{A_\psi}(\iota(\Pi)
\vert\CL_\psi )$, resp. $ \mathbf{E}_{B_\psi}(\iota(\Pi) \vert\CL_\psi )$, resp. 
$\mathbf{S}_\psi(\iota(\Pi)\vert \CL_{\psi})$.  

We say that $(\mathcal{L}, \iota)$ satisfies the {\it Eisenstein conditions}\index{Eisenstein
condition}
$({\rm EC}_r)$\index[NO]{EAC@$({\rm EC}_r)$} if $({\rm RC}_r)$ is fulfilled and if for each $\psi$  
\begin{equation}\label{Eisen}
  \begin{aligned}
    \big( \mathbf{S}_\psi\cdot
    \mathbf{E}_{A_{\psi}}\big)(\iota(\Pi) \mid \CL_{\psi} ) & = 0, \\
    \bigwedge^{4-[K^t:F^t]} \big( \mathbf{E}_{A_{\psi}}(\iota(\Pi)
    \mid \mathcal{L}_{\psi}) \big)& =  0.
  \end{aligned}
\end{equation}
In the case where $\psi$ is banal the first condition says 
\begin{equation}\label{Eisenbanal}
  \begin{aligned}
    \mathbf{E}_{A_{\psi}}(\iota(\Pi) \vert\CL_{\psi} ) & = 0 ,
    \text{ for all } \psi\in\Psi. 
  \end{aligned}
\end{equation}
and the second condition follows from the first.

The Eisenstein conditions do not depend on the $O_E$-morphism
$\alpha: S \longrightarrow  \Spec O_{E'}$. Indeed, if $S$ is connected, any
other choice of $\alpha$ differs by an automorphism $\rho \in \Gal(E'/E)$. 
In the decomposition (\ref{eigensp}) $\mathcal{L}_{\psi}$
is then replaced by $\mathcal{L}_{\rho \psi}$ and $\mathbf{E}_{\psi}$ is replaced
by $\rho (\mathbf{E}_{\psi}) = \mathbf{E}_{\rho \psi}$. Here the last
identity holds by the definition of the reflex field $E$. Therefore changing $\alpha$
does not change the Eisenstein conditions $({\rm EC}_r)$. If there exists
no $\alpha$, we use base change 
$\Spec O_{E'} \times_{\Spec O_E} S \longrightarrow S$ to define the
condition $({\rm EC}_r)$. The same arguments apply to the condition
$({\rm KC}_r)$.

We first note the following statement.
\begin{proposition} \label{unram} Let $S$ be an $O_E$-scheme and $\CL$ a locally free $\CO_S$-module with an   $O_K$-action $\iota\colon O_K\to\End_{\CO_S}(\CL)$. 
\begin{altenumerate}
\item The Eisenstein conditions $({\rm EC}_r)$ are independent of the uniformizer $\Pi$.

  \item When $K/\BQ_p$ is unramified, the Eisenstein conditions $({\rm EC}_r)$
  are implied by the Kottwitz condition $({\rm KC}_r)$. The same conclusion
 holds if $F=\BQ_p$ and $K/F$ is ramified.
 
 \item When $S$ is an $E$-scheme,  the  Eisenstein conditions $({\rm EC}_r)$ hold automatically.

 \end{altenumerate}
\end{proposition}
\begin{proof}
  Let us prove (i). 
  Let $\Pi'$ be another uniformizer. It is enough to show that the elements of $ O_K \otimes_{O_{K^t}, \psi} O_{E'}$,
  \begin{displaymath}
    \begin{array}{l} 
    \mathbf{E}_{A_{\psi}}(\Pi \otimes 1), \; \mathbf{E}_{A_{\psi}}(\Pi' \otimes 1) , \quad  \text{resp.}\\[2mm]   
    \mathbf{S}_{\psi}(\Pi \otimes 1) \mathbf{E}_{A_{\psi}}(\Pi \otimes 1), \; 
    \mathbf{S}_{\psi}(\Pi' \otimes 1)\mathbf{E}_{A_{\psi}}(\Pi' \otimes 1),
      \end{array}
    \end{displaymath}
  differ by a unit in $O_K \otimes_{O_{K^t}, \psi} O_{E'}$. Indeed, let $\tilde{E}'$ be the
  normalization of $E'$ in $\bar{\mathbb{Q}}_p$. Since
  $O_K \otimes_{O_{K^t},\psi} O_{E'}\rightarrow O_K\otimes_{O_{K^t},\psi} O_{\tilde{E}'}$
  is a flat extension of local rings, we can replace $E'$ by $\tilde{E}'$.
  By the definitions (\ref{Eisen1e}) and (\ref{Qpsi1e}), it suffices to show
  that the elements $\Pi \otimes 1 - 1 \otimes \varphi(\Pi)$ and
  $\Pi' \otimes 1 - 1 \otimes \varphi(\Pi')$ differ by a unit in
  $O_K \otimes_{O_{K^t}, \psi} O_{E'}$. But by \cite[Lem. 6.11]{unitshim} the   
  elements $\Pi\otimes1-1\otimes \Pi$ and $\Pi'\otimes1-1\otimes \Pi'$ of
  $O_K\otimes_{O_{K^t}}O_K$ differ by a unit, whence the assertion.

 Now  we prove (ii).  Let us only treat the case where $r$ is special; the banal case
  is similar. When $K=K^t$ is unramified over $\BQ_p$, then
  $\mathbf{E}(T)=T-\pi$ is a linear polynomial. Furthermore,
  $A_\psi$ has at most one element for $\psi\notin \{ \psi_0, \ov\psi_0\}$, and
  $A_{\psi_0}=A_{\ov\psi_0}=\emptyset$. Let $\psi\notin \{ \psi_0, \ov\psi_0\}$. If
  $A_\psi=\emptyset$, then $\CL_\psi =(0)$ and the Eisenstein
  condition relative to the index $\psi$ is empty; if $A_\psi$ has one
  element, the Eisenstein condition relative to the index $\psi$ is
  just equivalent to the definition of the $\psi$-th eigenspace in the
  decomposition \eqref{eigensp}. Something analogous applies to the
  indices $\psi_0, \ov\psi_0$. The case when $F=\BQ_p$ is handled in the same way. 
  
  Finally we prove (iii). Let $\tilde K$ be the normal closure of $K$ in $\ov\BQ_p$. It suffices to prove the assertion after replacing $S$ by its base change $S\times_{\Spec E} \Spec \tilde K$.  Then we have a decomposition
$$
O_K\otimes_{\BZ_p}\CO_S=\bigoplus\nolimits_{\varphi\in\Phi}\CO_S .
$$
Correspondingly, we have $\CL=\oplus \CL_\varphi$, and the endomorphism $\iota(\Pi)$ is diagonal with respect to this decomposition, with entries $\varphi(\Pi)\id_{\CL_\varphi}$. It is easy to see  that  $({\rm KC}_r)$ is equivalent to the condition
\begin{equation}
{\rm rank}\, \CL_\varphi=r_\varphi, \quad \forall \varphi\in\Phi .
\end{equation}
The Eisenstein conditions $({\rm EC}_r)$ involve endomorphisms of $\CL$ which are products of endomorphisms of the form $(\iota(\Pi)_{|\CL_\varphi}-\varphi(\Pi)\id_{\CL_\varphi})\oplus_{\varphi'\neq\varphi}\id_{\CL_{\varphi'}}$. From this, the conditions follow trivially. 
\end{proof}

Let us make the Eisenstein conditions more explicit in the case where $r$
is special. For this, we distinguish between the case when $K/F$ is
ramified and the case when $K/F$ is unramified.
Let $S$ be an $O_{E'}$-scheme and $\CL$ be a locally free $\CO_S$-module
satisfying $({\rm RC}_r)$. 
\begin{altitemize}
\item {\it $K/F$ ramified.}  In this case, we have $K^t=F^t$, and
  $\psi=\ov\psi$ for all $\psi\in\Psi$.
  Hence, \eqref{underconj} implies in this case
  \begin{equation}
    a_\psi=b_\psi= 
    \begin{cases}
      \begin{array}{ll}
        e, &\text{ if $\psi\neq\psi_0$}\\
        e-1,  &\text{ if $\psi=\psi_0 .$}
      \end{array}
    \end{cases}
  \end{equation} 
  We have   
\begin{displaymath}
  \mathbf{S}_{\psi_0}(T)=(T-\varphi_0(\Pi))(T-\ov\varphi_0(\Pi)).
  \end{displaymath}
The Eisenstein conditions become in this case 
  \begin{equation}\label{furtherDr2ram}
    \begin{aligned}
      \big( \mathbf{S}_{\psi_0}\cdot
      \mathbf{E}_{A_{\psi_0}}\big)(\iota(\Pi) \vert\CL_{\psi_0} ) & = 0, \\
      \bigwedge^{3} \big( \mathbf{E}_{A_{\psi_0}}(\iota(\Pi) \vert\CL_{\psi_0}
      ) \big)& =  0, \\ 
      \mathbf{E}_{A_\psi}(\iota(\Pi) \vert\CL_\psi ) &   = 0 ,
      \, \text{ for all }\psi\neq \psi_0 .
    \end{aligned}
  \end{equation}
\item {\it $K/F$ unramified.} In this case, $[K^t:F^t]=2$, and $\psi\neq\ov\psi$ for all $\psi\in\Psi$. Furthermore, 
  $a_\psi=b_{\ov\psi}$ and
  \begin{equation}\label{furtherDr2unram1e}
    \sum_{\psi\in\Psi\setminus\{\psi_0, \ov\psi_0\}} a_\psi= e(f-1),
    \quad a_{\psi_0}+a_{\ov\psi_0}=e-1 . 
  \end{equation}
\end{altitemize}
In this case, the Eisenstein conditions become 

\begin{equation}\label{furtherDr2unram}
  \begin{aligned}
    \big( \mathbf{S}_{\psi_0}\cdot
    \mathbf{E}_{A_{\psi_0}}\big)(\iota(\pi) \vert\CL_{\psi_0} ) & = 0, \\
    \bigwedge^{2} \big( \mathbf{E}_{A_{\psi_0}}(\iota(\pi) \vert\CL_{\psi_0}
    ) \big)& =  0, \\ 
    \big( \mathbf{S}_{\ov\psi_0}\cdot
    \mathbf{E}_{A_{\ov\psi_0}}\big)(\iota(\pi) \vert\CL_{\ov\psi_0} ) & = 0, \\
    \bigwedge^{2} \big( \mathbf{E}_{A_{\ov\psi_0}}(\iota(\pi) \vert\CL_{{\ov\psi_0}}
    ) \big)& =  0, \\ 
    \mathbf{E}_{A_\psi}(\iota(\pi) \vert\CL_\psi ) &   = 0 , \, \text{ for all }\psi\neq \psi_0, \ov\psi_0 .
  \end{aligned}
\end{equation}

\subsection{Local CM-pairs and CM-triples}\label{ss:loctriples}
Let $S$ be an $O_E$-scheme such that $p$ is locally nilpotent i.e.
a scheme over $\Spf O_E$. 
A {\it local CM-pair of type} $r$ \index{local CM-pair} is a pair $(X, \iota)$ such that
$X$ is a $p$-divisible group of height $4d$ and dimension $2d$ and  
$\iota$ is an $\mathbb{Z}_p$-algebra homomorphism 
\begin{displaymath}
  \iota: O_K \longrightarrow \End X 
\end{displaymath}
such that the rank condition $({\rm RC}_r)$ is satisfied for the
induced action of $O_K$ on $\Lie X$. In the split case
$O_K = O_F \times O_F$ we require moreover that in the induced
decomposition $X = X_1 \times X_2$ each factor is a $p$-divisible
group of height $2d$. 

Later we will introduce displays $\CP$, and these have a Lie algebra $\Lie\CP$, cf. Definition \ref{Rah2d}. Therefore, we can also speak of local CM-pairs $(\mathcal{P}, \iota)$ of type $r$,  where $\mathcal{P}$ is a display over
$S$, cf. section \ref{s:relativeDT}.

Let $S = \Spec k$ be a perfect field of characteristic $p$ which is endowed with an
$O_{E'}$-algebra structure. In this case,  a display in the same thing
as a Dieudonn\'e module $\mathcal{P} = (P,F,V)$, where $P$ is a
finitely generated free module over the ring of Witt vectors
$W(k)$. If $\mathcal{P}$ is the Dieudonn\'e module of $X$, there is
a canonical isomorphism of $k$-vector spaces $\Lie X \cong P/VP$.

Via $\iota$ we regard $P$ as a $O_K \otimes_{\mathbb{Z}_p} W(k)$-module. 
The homomorphisms
\begin{equation}\label{localpairs2e}
  \psi: O_{K^t} \longrightarrow O_{E'} \longrightarrow k,
\end{equation}
$\psi \in \Psi$ lift uniquely to homomorphisms
\begin{equation}\label{localpairs3e}
  \tilde{\psi}: O_{K^t} \longrightarrow W(k). 
\end{equation}
We obtain a ring isomorphism
\begin{displaymath}
  O_K \otimes_{\mathbb{Z}_p} W(k) =
  \prod_{\psi \in \Psi} O_K \otimes_{O_{K^t}, \tilde{\psi}} W(k). 
\end{displaymath}
This induces a decomposition
\begin{equation}\label{localpairs4e}
  P = \oplus_{\psi \in \Psi} P_{\psi}. 
\end{equation}
More explicitly
\begin{displaymath}
  P_{\psi} = \{x \in P \; | \;
  \iota(a)x = \tilde{\psi}(a) x, \; \text{for} \; a \in O_{K^t} \}.  
\end{displaymath}  

Let us denote by $\sigma$
the Frobenius automorphism of $W(k)$. The
operators $F$ and $V$ on $P$ induce $\sigma$-linear maps
\begin{equation}\label{localpairs1e} 
  F: P_{\psi} \longrightarrow P_{\sigma \psi}, \quad V: P_{\sigma \psi}
  \longrightarrow P_{\psi}.   
\end{equation}
Here $\sigma \psi$ denotes the composite of (\ref{localpairs2e}) with
the absolute Frobenius of $k$.  
\begin{lemma}\label{localpair1l}
  Let $(P, \iota)$ be local CM-pair of type $r$ over a perfect field
  $k$. Then $P$ is a free $O_K \otimes_{\mathbb{Z}_p} W(k)$-module  of rank 2. 
\end{lemma}
\begin{proof}
  Since $F V = p$ it follows that
  \begin{equation}\label{localpairs3e4} 
    \rank_{O_K \otimes_{O_{K^t}, \tilde{\psi}} W(k)} P_{\psi} = \rank_{O_K
      \otimes_{O_{K^t}, \sigma \tilde{\psi}} W(k)} P_{\sigma \psi}. 
  \end{equation}
  Since $\rank_{W(k)} P = 4d$, and by the extra condition in the split
  case, this implies that the common rank of (\ref{localpairs3e4}) is
  $2$. This proves the Lemma. 
\end{proof}

To each local CM-pair $(X, \iota)$ we define the
{\it conjugate dual} $(X^{\vee}, \iota^{\wedge})$. \index{conjugate dual local CM-pair}\index[NO]{XAA@$(X^{\vee}, \iota^{\wedge})$}Here $X^{\vee}$ is
the dual $p$-divisible group of $X$ but we change the action dual to $\iota$
by the conjugation of $K/F$, i.e.,  $\iota^{\wedge}(a)=\iota^\vee(\bar a)$. We
will denote the conjugate dual simply by $X^{\wedge}$. 
\begin{lemma}\label{conjdual}
The  conjugate dual of a CM-pair $(X, \iota)$ of type $r$ is again a local CM-pair of type $r$. If $(X, \iota)$ satisfies the Kottwitz condition $({\rm KC}_r)$, resp., the Eisenstein conditions $({\rm EC}_r)$, then so does its conjugate dual. 
\end{lemma}
\begin{proof} For the first assertion, we may assume that we are over an algebraically closed field. We use the Dieudonn\'e module $\mathcal{P}$. We set
\begin{displaymath}
  P^{\vee} = \Hom_{W(k)} (P, W(k)).
\end{displaymath}
We use the canonical pairing
\begin{equation}\label{localpairs5e}
  \langle\, \; , \; \rangle: P \times P^{\vee} \longrightarrow W(k).  
\end{equation}
The operators $F$ and $V$ on the dual Dieudonn\'e module
$\mathcal{P}^{\vee}$ are defined by the equations
\begin{displaymath}
  \begin{aligned}
    \langle\, Vx , Vx^{\vee} \rangle & =  p \sigma^{-1}
    (\langle\, x , x^{\vee} \rangle), & \quad x \in P, \; x^{\vee} \in
    P^{\vee}\\
    \sigma (\langle\, Vx , x^{\vee} \rangle) & = 
    \langle\, x , Fx^{\vee} \rangle. &  
  \end{aligned}
\end{displaymath}
One of these equations implies the other. It follows that
$VP/pP \subset P/pP$ and
$VP^{\vee}/pP^{\vee} \subset P^{\vee}/pP^{\vee}$ are orthogonal
complements with respect to the non-degenerate pairing of $k$-vector
spaces,
\begin{displaymath}
  P/pP \times P^{\vee}/pP^{\vee} \longrightarrow k.  
\end{displaymath}
If we use the action $\iota^{\wedge}$,  we write for the decomposition
(\ref{localpairs4e})
\begin{displaymath}
  P^{\vee} = \oplus_{\psi \in \Psi} P^{\wedge}_{\psi} = P^{\wedge}. 
\end{displaymath}
Then $P_{\psi_1}$ and $P_{\psi_2}^{\wedge}$ are for
$\psi_1 \neq \bar{\psi}_2$ orthogonal with respect to
(\ref{localpairs5e}) and  
\begin{equation}\label{localpairs6e}
  \langle\, \; , \; \rangle: P_{\psi} \times P^{\wedge}_{\bar{\psi}}
  \longrightarrow W(k).  
\end{equation}
is a perfect pairing. The $k$-vector spaces $VP_{\sigma \psi}/pP_{\psi}$ and
$VP^{\wedge}_{\sigma \bar{\psi}}/pP^{\wedge}_{\bar{\psi}}$ are orthogonal
complements with respect to the induced non-degenerate $k$-bilinear form
\begin{displaymath}
  P_{\psi}/pP_{\psi} \times
  P^{\wedge}_{\bar{\psi}}/pP^{\wedge}_{\bar{\psi}} \longrightarrow k. 
\end{displaymath}
Let us assume that $K/F$ is unramified or split. In this case
Lemma \ref{localpair1l} implies $\rank_{W(k)} P_{\psi} = 2e$ and by 
(\ref{localpairs6e}) $\rank_{W(k)} P^{\wedge}_{\psi} = 2e$.  
Since $P$ satisfies $({\rm RC}_r)$ we find by the orthogonality
above
\begin{displaymath}
  \rank_k P^{\wedge}_{\bar{\psi}}/VP^{\wedge}_{\sigma\bar{\psi}} = 2e -
  \rank_k P_{\psi}/VP_{\sigma \psi} = 2e - \sum_{\varphi | \psi} r_{\varphi} =
  \sum_{\varphi | \psi} (2 - r_{\varphi}) = \sum_{\varphi | \psi} r_{\bar{\varphi}}.
\end{displaymath}
This shows that the conjugate dual satisfies $({\rm RC}_r)$. The case
$K/F$ ramified is similiar.

For the proof of the assertion concerning  $({\rm KC}_r)$, we refer to Proposition \ref{dualKC}.  For the proof of the assertion concerning  $({\rm EC}_r)$, we refer to  Corollary \ref{dualEisen2c}  in the case when $K/F$ is unramified or split, resp., Corollary \ref{dualECram} when $K/F$ is ramified.
\end{proof}

The notion of a   {\it local CM-triple of type} \index{local CM-triple, polarized local CM-pair} $r$ over $S$ was
introduced in \cite{KRnew}.  This is a triple $(X, \iota, \lambda)$,
where $(X, \iota)$ is a local CM-pair of type $r$ and
$\lambda: X\longrightarrow X^\vee$ is an anti-symmetric isogeny (also
called a {\it polarization})  such that the corresponding Rosati involution
induces the non-trivial automorphism on $K/F$. In particular $\lambda$
induces a morphism of local CM-pairs
\begin{displaymath}
  \lambda: (X,\iota) \longrightarrow (X^{\wedge}, \iota^{\wedge}).  
\end{displaymath}
In the present paper, we will also say that $(X, \iota, \lambda)$ is a \emph{polarized local CM-pair} (of type $r$). We call the polarized local CM-pair $(X, \iota, \lambda)$ {\it principal} if $\Ker \lambda=0$; we call it   {\it almost principal} \index{almost principal polarization}when $\Ker \lambda\subset X[\iota(\pi)]$
and $\Ker \lambda$ has order $p^{2f}$. We will distinguish  principal polarized local CM-pairs from almost principal  ones  by attaching the integer ${\tt h}=0$ to the principal case, and
${\tt h}=1$ to the almost principal case, i.e., $\height \lambda=2f{\tt h}$.  We will see in \S \ref{ss:uniquefram} below that, when $K/F$ is a ramified field extension, then the almost principal case does not occur. Compare Lemmas \ref{Kneun0l} and \ref{Kneun01l} for the analogue in linear algebra. 

\subsection{The invariant of a local CM-triple}\label{ss:inv}  

Let $K/F$ be a field extension.  We recall from \cite{KRnew} the definition of the invariant of a CM-triple, in a slightly more general context.  

Let $k$ be an algebraically closed field of characteristic
$p$ and let $W(k)$ be the ring of Witt vectors. We set
$W(k)_{\mathbb{Q}} = W(k) \otimes \mathbb{Q}$. Let $(M, F, V)$ be a
Dieudonn\'e module of height $4d$ and dimension 
$2d$ which is endowed with a $\mathbb{Z}_p$-algebra homomorphism
\begin{displaymath}
  \iota: K \longrightarrow \End^{0} (M, F, V).  
\end{displaymath}
We set $M_{\mathbb{Q}} = M \otimes \mathbb{Q}$. 
We assume that $M$ is endowed with a non degenerate alternating bilinear form
\begin{displaymath}
  \beta: M_{\mathbb{Q}} \times M_{\mathbb{Q}} \longrightarrow W(k)_{\mathbb{Q}}
\end{displaymath}
Let us denote by $\sigma$ the Frobenius automorphism of $W(k)$. We
require the following properties:
\begin{displaymath}
  \begin{aligned}
    \beta (Vx, Vy) &= p \sigma^{-1}(\beta(x,y)), & x, y \in
    M_{\mathbb{Q}}\\ 
    \beta (\iota(a) x, y) &= \beta(x, \iota(\bar{a})y), & a \in K.
  \end{aligned}
\end{displaymath}
We will associate to such a set of data $(M, \iota, \beta)$ an invariant
$\inv(M, \iota, \beta) \in \{\pm 1\}$. 
We set $\Psi = \Hom_{\text{$\BQ_p$-Alg}}(K^t, W(k)_{\mathbb{Q}})$.

The ring $O_K \otimes_{\mathbb{Z}_p} W(k)$ decomposes 
\begin{equation}\label{Kneunpsi4e}
  K \otimes_{\mathbb{Z}_p} W(k) =
  \prod_{\psi} K \otimes_{O_{K^t}, \tilde{\psi}} W(k).
\end{equation}
If $\xi = (\xi_{\psi})$ is an element of (\ref{Kneunpsi4e}). Then we
set
\begin{equation}\label{Kneunord1e} 
  \ord_{K\otimes W(k)} \xi = \ord_p \Nm_{K/\mathbb{Q}_p} \xi 
  = \sum_{\psi} \ord_{\Pi} \xi_{\psi} \in \mathbb{Z}.
\end{equation}
The Frobenius homomorphism $\sigma$ acts via the second factor on
$K \otimes_{\mathbb{Z}_p} W(k)$. The $\sigma$-conjugacy class of an
element $\xi \in (K \otimes W(k))^{\times}$ is uniquely determined by
$\ord_{K\otimes W(k)} \xi$.  

We view $M_{\mathbb{Q}}$ as a $K \otimes_{\mathbb{Z}_p} W(k)$-module
and suppress the notation $\iota$. This is a free module of rank $2$.
We define an anti-hermitian form $\varkappa=\varkappa_\beta$,
\begin{displaymath} 
 \varkappa: M_{\mathbb{Q}} \times M_{\mathbb{Q}} \longrightarrow
  K \otimes_{\mathbb{Z}_p} W(k), 
\end{displaymath}
on the $K \otimes_{\mathbb{Z}_p} W(k)$-module $M_{\mathbb{Q}}$ 
by the formula
\begin{equation}\label{Kneun03e}
  \Trace_{K/\mathbb{Q}_p} (a\varkappa(x, y)) = \beta(ax, y), \quad
  x, y \in M_{\mathbb{Q}}, \; a \in K \otimes W(k).   
\end{equation}
Then $\varkappa$ satisfies
\begin{equation}\label{Kneun1e}
 \varkappa(Vx, Vy) = p \sigma^{-1}( \varkappa(x, y)). 
\end{equation}

We write
$\wedge^2 M_{\mathbb{Q}} := \bigwedge^{2}_{K \otimes W(k)} M_{\mathbb{Q}}$
for the exterior product as a $K \otimes_{\mathbb{Z}_p}
W(k)$-module. This is a free 
$K \otimes_{\mathbb{Z}_p} W(k)$-module of rank $1$. According to
(\ref{Kneunpsi4e}) we have decompositions
\begin{displaymath}
  \begin{aligned}
    M_{\mathbb{Q}} &= \bigoplus\nolimits_{\psi} M_{\mathbb{Q},\psi}\\
    \wedge^{2} M_{\mathbb{Q}} &= \bigoplus\nolimits_{\psi} \big(
    \bigwedge^2_{K \otimes_{K^t,\psi} W_{\mathbb{Q}}(k)} M_{\mathbb{Q}, \psi} \big).
  \end{aligned}
\end{displaymath}
We choose an isomorphism
$\wedge^2 M_{\mathbb{Q}} \cong K \otimes_{\mathbb{Z}_p} W(k)$. Then we can write
\begin{displaymath}
\wedge^2 V (z) = \gamma \sigma^{-1}(z).   
  \end{displaymath}
We have 
\begin{displaymath}
  \begin{aligned} 
    \ord_{K\otimes W(k)} \wedge^2 V &= \ord_p \Nm_{K/\mathbb{Q}_p}
    \det\textstyle_{K \otimes_{\mathbb{Z}_p} W(k)} (V |M_{\mathbb{Q}})\\
    &=\ord_p {\rm det}_{W(k)} (V |M_{\mathbb{Q}}) = \dim M = 2d.
  \end{aligned}
\end{displaymath}
Therefore we find $\ord_{K\otimes W(k)} \gamma = 2d$. Since
$\ord_{K\otimes W(k)} p = 2d$, the elements
$p, \gamma \in K \otimes_{\mathbb{Z}_p} W(k)$ are in the same $\sigma$-conjugacy
class by the remark after (\ref{Kneunord1e}). We conclude that there is a 
generator $x \in \wedge^2 M_{\mathbb{Q}}$  such that
\begin{equation}\label{Kneun02e}
  \wedge^2 V ( x) = px.
\end{equation}
Note that the last equation is equivalent to $\wedge^2 F (x) = px$.
Any other generator with this property has the form $ux$, where
$u \in K^{\times}$. We consider the hermitian form
\begin{displaymath}
  h = \wedge^2\varkappa: \wedge^2 M_{\mathbb{Q}} \times
  \wedge^2 M_{\mathbb{Q}} \longrightarrow K \otimes W(k). 
\end{displaymath}
One deduces that $h (x, x)$ is an element of
$F \subset K \otimes W(k)$ which is $\neq 0$ because $\beta$ is non-degenerate
by assumption. We denote by 
\begin{equation}\label{Kneun2e} 
  \mathfrak{d}_{K/F} (M, \iota, \beta) \in F^{\times}/\Nm_{K/F} K^{\times} 
\end{equation}
the class of $h(x,x)$. This class is called the \emph{discriminant} \index{discriminant}of $(M, \iota, \beta)$ and is independent of the choice of $x$. 

Let $a \in F$. Then $a  \varkappa$ is again an anti-hermitian form which
satisfies (\ref{Kneun1e}). We can replace $ \varkappa$ by $a \varkappa$
in the definition of (\ref{Kneun2e}) without changing the discriminant. We denote by
\begin{equation}\label{Kneun3e} 
  \inv (M, \iota, \beta) \in \{\pm 1\} 
\end{equation}\index[NO]{IAA@$\inv (M, \iota, \beta)$} 
the image of $\mathfrak{d}_{K/F} (M, \iota, \beta)$ by the canonical
isomorphism $F^{\times}/\Nm_{K/F} K^{\times} \simeq \{\pm 1\}$.

Let $r$ be a local CM-type of rank $2$. Let $E$ be the reflex
field. Let $O_E \longrightarrow k$ an algebra structure of the
algebraically closed field $k$.  
Let $(X, \iota, \lambda)$ be a local triple of CM-type $r$ over $k$. Let
$(M, \iota, \beta)$ be the associated Dieudonn\'e module with its
polarization $\beta$. Then we set\index[NO]{IAB@$\inv (X, \iota, \lambda)$}
\begin{displaymath}
  \inv (X, \iota, \lambda) := \inv (M, \iota, \beta).
\end{displaymath}
For CM-triples of CM-type $r$, we use also the adjusted invariant  $\inv^{r} (X, \iota, \lambda)=\inv^{r} (M, \iota, \beta)$, cf. section \ref{ss:radj}. In the case at hand we have 
\begin{equation}\label{Kneun4e}  
  \inv^{r} (M, \iota, \beta) =
  \begin{cases}
    (-1)^{d-1}\inv (M, \iota, \beta), \; \text{for} \; r \; \text{special},\\[2mm]
    (-1)^{d}\inv (M, \iota, \beta), \; \text{for} \; r \; \text{banal}.
  \end{cases}
\end{equation}

\subsection{Uniqueness of framing objects}\label{ss:uniquefram}

In this subsection, we discuss the existence and uniqueness of \emph{framing objects} that are used in the formulation of the formal moduli problems. The proofs of these statements are given later in the paper.  

  Let $r$ be a generalized local CM-type of rank $2$ for $K/F$. Let $\ov k$ be an algebraic closure  of the residue field $\kappa_E$ of
  $O_E$. Consider CM-triples $(X, \iota, \lambda)$ over $\ov k$ which satisfy $({\rm KC}_r)$ and $({\rm EC}_r)$.
  \smallskip

(i) \label{uniqueframe,i} 
 {\it  Assume that  $r$ is special.  If $K/F$ is ramified, then a local CM-triple of type $r$ over $\ov k$ as above such that the polarization is principal and with $r$-adjusted invariant $-1$ is isoclinic. When  $K/F$ is unramified, then a local CM-triple of type $r$ over $\ov k$ as above such that the polarization is almost principal has    $r$-adjusted invariant $-1$ and is isoclinic. In either case, any two such CM-triples are isogenous by a $O_K$-linear quasi-isogeny of height zero that preserves the polarizations. 

  Furthermore, the group of $O_K$-linear self-isogenies of such a local CM-triple, preserving the polarization,  can be identified with the unitary group of the  split $K/F$-anti-hermitian space $C$ of dimension $2$. 
  }

  \smallskip
 The assertions concerning slopes follow from Corollary \ref{C'onslopes1c}. The uniqueness assertion is in the ramified case the content of Proposition \ref{ramFrame2p}, and in the unramified case of Proposition \ref{unrFrame2p}.  The last part of the assertion follows from the fact that the contraction functor is an equivalence of categories. 
 
 \smallskip
 
  (ii), a)\label{uniqueframe,iia} 
  {\it Let $r$ be banal non-split. Any  local CM-triple of type $r$ over $\ov k$ as above is isoclinic.  The group of $O_K$-linear self-isogenies of such a local CM-triple, preserving the polarization,  can be identified  with the unitary group of a  $K/F$-anti-hermitian space $C$ of dimension $2$. When $K/F$ is unramified and the polarization is principal, then the  anti-hermitian space $C$ is split and the $r$-adjusted invariant is $1$;  when   $K/F$ is ramified and the polarization is principal, then the  anti-hermitian space $C$ is non-split and the $r$-adjusted invariant is $-1$; when $K/F$ is unramified and the polarization is almost principal, then the anti-hermitian space $C$ is non-split and the $r$-adjusted invariant is $-1$. The case $K/F$ ramified and almost principal polarization does not occur. Any two  CM-triples with the same  $r$-adjusted invariant are  isogenous by a $O_K$-linear
   quasi-isogeny of height zero that preserves the polarizations.
   }
   
   \smallskip
   
   The assertions concerning slopes follow from Corollary \ref{C'onslopes1c}. The uniqueness assertion is  the content of Proposition \ref{CFbanal5p}. Similar arguments apply to the banal split case. 

\smallskip
   
 (ii), b)\label{uniqueframe,iib} 
  {\it  Let  $r$ be banal split. Then the $p$-divisible group underlying a  local CM-triple of type $r$ is the direct product of two isoclinic $p$-divisible groups of slope $\lambda$, resp. $1-\lambda$, where $\lambda$ depends only on $r$. Any two  local CM-triples of type $r$ over $\ov k$ are  isogenous by a $O_K$-linear
   quasi-isogeny of height zero that preserves the polarizations. The group of $O_K$-linear self-isogenies of such a local CM-triple, preserving the polarization,  can be identified with $\Res_{F/\BQ_p}(\GL_2)$. 
}

  \bigskip

 The  chart below summarizes the discussion above. The last column lists the group of $O_K$-linear self-isogenies preserving the polarization. In all cases listed below, the framing object is unique up to quasi-isogeny of height zero. 

\begin{table}[h]
\begin{center}
\begin{tabular}{l c r c c} 
\rule[-9pt]{0pt}{14pt}
type $r$ &\, $K/F$&\,$\inv^r$&\, polarization&\quad  group of self-isogenies\\ \hline\smallskip
special&\quad  ramified& $-1$&\quad  principal&\quad  quasi-split unitary group\\ \hline\smallskip
special &\quad  unramified&$-1$&\quad almost principal&\quad  quasi-split unitary group\\ \hline\smallskip 
banal&\quad  ramified&$-1$&\quad principal&\quad non-quasi-split unitary group\\ \hline\smallskip
banal&\quad  unramified&$1$&\quad principal&\quad  quasi-split unitary group\\ \hline\smallskip
banal&\quad unramified&$-1$&\quad almost principal&\quad  non-quasi-split unitary group\\ \hline\smallskip
banal&\quad split&1&\quad principal&\quad  $\GL_2/F$\\ \hline
{}&{}&{}&{}&
\end{tabular}
\caption{Framing objects}
\end{center}
\end{table}

\begin{remark} The statement (i) above  generalizes  \cite[Prop. 5.4]{KRnew}. However, the proof of the uniqueness assertion given there is incomplete. Note that for a local CM-type of the first kind in the sense of loc.~cit. we have  imposed $F=\BQ_p$; therefore the condition in loc.~cit. that $\varepsilon=\inv(X, \iota, \lambda)=-1$ implies that the associated  anti-hermitian space $(C, h)$ is split (in this case the $r$-adjusted invariant coincides with the invariant).
\end{remark}
\begin{remark}
  The statement (i)  is closely related to  the fact that $B(G, \{\mu\})$ has only one element, cf. \cite{K2}, \S 6.  Here $G=\Res_{F/\BQ_p}(\GU)$ is the linear algebraic group over $\BQ_p$ associated to the group of unitary similitudes of the non-split  anti-hermitian space of dimension $2$ over $K$, and $\{\mu\}$ is the conjugacy class of cocharacters with component $(1, 0)$ for $\varphi_0$ and central component for 
  $\varphi\neq\varphi_0$. In fact, it seems that the essential contents of  the calculations in section \ref{ss:invcont} is to show that the Frobenius element of a local CM-triple of type $r$ with $r$-adjusted invariant $-1$  over $\ov k$ defines an element in $B(G, \mu)$.  Recall from the introduction of \cite{K2} that $p$-adic uniformization can only be expected when the pair $(G, \mu)$ is \emph{uniform}, i.e.,  $B(G, \mu)$ consists of a single element. 
\end{remark}

\subsection{Formal moduli spaces}\label{ss:defofform}
In this subsection we are going to define RZ-spaces of formal local CM-triples, and formulate our main results about them. We fix $K/F$ as before. 

First let $r$ be special, so that  $K/F$ is a field extension.
We fix  a local CM-triple $(\BX, \iota_{\BX}, \lambda_\BX)$ of type $r$ over
$\bar{\kappa}_E$ as in (i) of subsection \ref{ss:uniquefram} (a {\it framing object}). \index{framing object} We assume that if $K/F$ is ramified, then $\lambda_\BX$ is principal and that, if $K/F$ is unramified, then $\lambda_\BX$ is almost principal. Then, in either case,  the $r$-adjusted invariant equals  $-1$.  We identify $\bar{\kappa}_E$ with the residue
class field of $O_{\breve E}$, the ring of integers in the completion of
the maximal unramified extension of $E$. Let $({\rm Sch}/{O_{\breve E}})$ be the category of $O_{\breve E}$-schemes $S$ such that the ideal sheaf $\pi\CO_S$ is locally nilpotent.

\begin{definition}\label{CRFd2} We set ${\tt h}=0$ \index[NO]{HDA@${\tt h}$} if $K/F$ is ramified,
  and ${\tt h}=1$ if $K/F$ is unramified. 

  We define a functor $\mathcal{M}_{K/F, r}$\index[NO]{MBA@$\mathcal{M}_{K/F, r}$} on   $ ({\rm Sch}/{O_{\breve E}})$. A point of
  $\mathcal{M}_{K/F, r}(S)$ consists of an isomorphism class of the following
  data: 
  \begin{altenumerate}
  \item[(1) ]\, Two local CM-pairs $(X_0, \iota_0)$, $(X_1, \iota_1)$ of
    CM-type $r$ over $S$ which satisfy the Eisenstein conditions
    $({\rm EC}_r)$ relative to a fixed uniformizer $\pi$ of $F$ and
    the Kottwitz condition $({\rm KC}_r)$.
  \item[(2) ] Two isogenies of $p$-divisible $O_K$-modules
    \begin{displaymath}
      X_0 \longrightarrow X_1 \longrightarrow X_0  ,
    \end{displaymath}
    which have both height $2f{\tt h}$ and such that the composite is
    $\iota_0(\pi)^{{\tt h}}\,\id_{X_0}$.
  \item[(3) ] An isomorphism of $p$-divisible $O_K$-modules
    \begin{displaymath}
      \varkappa: X_1 \isoarrow X^{\wedge}_0. 
    \end{displaymath}
    We require that the composite
    $\lambda\colon X_0 \longrightarrow X_1  \isoarrow X^{\wedge}_0$ is a polarization of
    $X_0$, i.e., this map is anti-symmetric, and that  this polarization  is principal when $K/F$ is 
    ramified, and almost principal when $K/F$ is unramified.\footnote{It can be proved that, when $K/F$ is unramified, the fact that $\Ker \lambda\subset X_0[\pi]$ follows automatically from the assumption that $\deg\lambda=p^{2f}$, cf. Proposition  \ref{SFB701p}. We impose the condition  $\Ker \lambda\subset X_0[\pi]$ in order to make transparent that the moduli problem $\mathcal{M}_{K/F, r}$ is of the kind considered in \cite{RZ}.  }
  \item[(4) ] A quasi-isogeny of height zero
    of $p$-divisible $O_K$-modules
    \begin{displaymath}
      \rho_X: X_0 \times_{S}  \bar{S} \longrightarrow \mathbb{X}
      \times_{\Spec \bar{\kappa}_E}  \bar{S}, 
    \end{displaymath}
    such that the pullback quasi-isogeny $\rho^*(\lambda_\BX)$ differs from
    $\lambda_{|X_0 \times_S \ov{S}}$ by a scalar in $F^{\times}$, locally on $\ov S$.
    Here $\ov{S} = S \otimes_{O_{\breve E}} \bar{\kappa}_E$. We call $\rho_X$
    a {\it framing}. 
  \end{altenumerate}
\end{definition}
We denote the data above simply by $(X, \iota, \lambda, \rho)$. Another datum
$(X', \lambda', \rho')$ defines the same point of $\mathcal{M}_{K/F, r}(S)$
iff there are $O_K$-isomorphisms $X_0 \isoarrow X'_0$ and $X_1 \isoarrow X'_1$
which commute with the data (2) and (4) above. This implies that the
isomorphism $X_0 \isoarrow X'_0$ respects the polarizations up to a
factor in $O_F^{\times}$. 

To ease the notation we write $\mathcal{M}_r = \mathcal{M}_{K/F, r}$.
\index[NO]{MBB@$\mathcal{M}_r=\mathcal{M}_{K/F, r}$} If $R$ is a $p$-adic
$O_{\breve{E}}$-algebra we set
$\mathcal{M}_r(R) = \varprojlim_{n} \; \mathcal{M}_r(\Spec R/p^n R)$. 

It follows by the methods of \cite{RZ} that $\mathcal{M}_r$ is
representable by a formal scheme which is locally formally of finite type
over $\Spf O_{\breve{E}}$. 
Let $J$\index[NO]{JAA@$J$}  be the algebraic group over $\BQ_p$ of unitary $K$-linear quasi-automorphisms of $(\BX, \iota_{\BX}, \lambda_\BX)$ which preserve the polarization up to a scalar in $\BQ_p^\times$. Let $J^1$
\index[NO]{JAB@$J^1$}  denote the derived group of $J$. Then $J^1(\BQ_p)$ acts on the functor $\CM_r$ by changing the framing. It follows from (i) in subsection \ref{ss:uniquefram} that $J^1$ can be identified with $\Res_{F/\BQ_p}(\SU)$, where $\SU$ denotes the quasi-split special unitary group in two variables over $F$. Note that $\SU$ is isomorphic to  $\SL_2/F$. 

The first main result in the local case can now be stated as follows. 
\begin{theorem}\label{mainUR} Let $r$ be special. 
  Then the  functor $ \CM_{K/F, r}$ is represented by $\wh{\Omega}_F \hat{\otimes}_{O_F, \varphi_0} O_{\breve E}$. More precisely, there exists a unique isomorphism of formal schemes 
  $$
  \CM_{K/F, r}\simeq \wh{\Omega}_F \hat{\otimes}_{O_F, \varphi_0} O_{\breve E} ,
  $$
  which is equivariant with respect to a fixed identification $J^1(\BQ_p)\simeq \SL_2(F)$. In particular, $\CM_{K/F, r}$ is flat over $\Spf O_{\breve E}$ with semi-stable reduction.
\end{theorem}

Now let $r$ be banal. Fix  a local CM-triple $(\BX, \iota_{\BX}, \lambda_\BX)$ over $\ov k$ as in (ii) a) or (ii) b) in subsection \ref{ss:uniquefram}.
  We write the height of $\lambda_\BX$ as $2f{\tt h}$, where ${\tt h}\in\{ 0,1\}$. We assume that ${\tt h}=0$ when $r$ is banal split, or when   $r$ is banal  and $K/F$ is a ramified field extension.  Recall from  from  (ii) a) in subsection \ref{ss:uniquefram}  that, when  $r$ is non-split,   there is a anti-hermitian space $C=C(\BX, \iota_{\BX}, \lambda_\BX)$ attached to $(\BX, \iota_{\BX}, \lambda_\BX)$.  By Proposition \ref{CFbanal5p}, the framing object $(\BX, \iota_{\BX}, \lambda_\BX)$ is uniquely defined up to isogeny by the $r$-adjusted invariant $\inv^r (\BX, \iota_{\BX}, \lambda_\BX)=\inv (C)\in \{\pm1\}$ (see Proposition \ref{KneunBa1l} for this last identity). To make our statements uniform, we set $\inv^r (\BX, \iota_{\BX}, \lambda_\BX)=1$ in the banal split case. 

We may now define a variant for banal $r$ of the  functor $\mathcal{M}_{K/F, r}$ of Definition \ref{CRFd2}. Since the functor depends not only on $K/F$ but also on $\inv^r (\BX, \iota_{\BX}, \lambda_\BX)$, we denote this functor by $\mathcal{M}_{K/F, r, \varepsilon}$, \index[NO]{MBC@$\mathcal{M}_{K/F, r, \varepsilon}$} where $\inv^r (\BX, \iota_{\BX}, \lambda_\BX)=\varepsilon$. When $r$ is banal split, we have $\varepsilon=1$; when $r$ is banal non-split and $K/F$ is unramified, then $\varepsilon=(-1)^{\tt h}$, cf. Proposition \ref{CFbanal5p}.  

Let  $S\in({\rm Sch}/{O_{\breve E}})$.  A point of
 $\mathcal{M}_{K/F, r, \varepsilon}(S)$ consists of an isomorphism class of exactly 
 the same data as in Definition \ref{CRFd2}. 
\begin{theorem}\label{main-banal} 
Let $r$ be banal, and let $\varepsilon\in\{\pm1\}$. 
 The formal scheme $\CM_{K/F, r, \varepsilon}$ is isomorphic to $(\Spf\, O_{\breve E})\times (J(\BQ_p)^o/C_{\ov M})$, where $J(\BQ_p)^o$ denotes the subgroup of elements of $J(\BQ_p)$ which preserve the polarization up to a scalar in $\BZ_p^\times$. More precisely, there exists a unique isomorphism of formal schemes 
  $$
  \CM_{K/F, r, \varepsilon}\simeq (\Spf\, O_{\breve E})\times (J(\BQ_p)^o/C_{\ov M}) ,
  $$
  which is equivariant for the action of $J(\BQ_p)^o$. In particular, $\CM_{K/F, r, \varepsilon}$ is formally \'etale over $\Spf\, O_{\breve E}$. 
  
  Here, when $r$ is non-split,  $J(\BQ_p)^o$ can be identified with the group of $K$-linear automorphisms of $C=C(\BX, \iota_{\BX}, \lambda_\BX)$  preserving the anti-hermitian form  up to a factor in $\BZ_p^\times$, and $C_{\ov M}$ is the stabilizer in $J(\BQ_p)^o$ of a lattice $\ov M$ in $C$ which is self-dual when  ${\tt h}=0$ and almost self-dual when  ${\tt h}=1$.  When $r$ is split, then  $J(\BQ_p)^o/C_{\ov M}$ can be identified with  the set of lattices in  the two-dimensional standard $F$-vector space of dimension $2$.
\end{theorem}
In the later part of the paper, we write simply $J(\BQ_p)^o/C_{\ov M}$ for the formal scheme $(\Spf\, O_{\breve E})\times (J(\BQ_p)^o/C_{\ov M})$ over $\Spf O_{\breve E}$.

\section{Background on Display Theory}\label{s:relativeDT}
In this section, $K/\mathbb{Q}_p$ is an arbitrary finite field extension with
ring of integers $O = O_K$, and  $\Nilp_O$ will denote the category of  $O$-algebras $R$ such that $p$ is nilpotent
in $R$. We recall the classification of strict formal
$p$-divisible $O$-modules over $R\in\Nilp_O$ proved in \cite{ACZ}.
A main ingredient is the Ahsendorf functor, which we present in a new form
which is better suited for our applications.

\subsection{Displays}\label{ss:reldisp}

We fix a prime element $\pi \in O$. We
denote by $q$ the number of elements in the residue class field
$\kappa$ of $O$.  

\begin{definition}(\cite[Def. 3.1]{ACZ}, \cite{L2},
  \cite{Zi2}\label{defOframe}) 
  Let $R$ be an $O$-algebra. A frame\index{frame} $\mathcal{F}$
  \index[NO]{FBA@$\mathcal{F}$} for $R$ consists of the
  following data:
  \begin{enumerate}
  \item[(1)] An $O$-algebra $S$ and a surjective $O$-algebra homomorphism
    $S \rightarrow R$. We denote the kernel by $I$.  
  \item[(2)] An $O$-algebra endomorphism $\sigma: S \rightarrow S$.  
  \item[(3)] A $\sigma$-linear map of $S$-modules
    $\dot{\sigma}: I \rightarrow S$. 
  \end{enumerate}
The following conditions are required. 
  \begin{enumerate}
  \item[(i)]  $I + pS$ is contained in the radical of $S$. 
  \item[(ii)] $\sigma(s) \equiv s^{q} \mod \pi S$ for all $s \in S$. 
  \item[(iii)] $\dot{\sigma}(I)$ generates $S$ as an $S$-module. 
  \end{enumerate}
\end{definition}

We will denote a frame by $\mathcal{F} = (S, I, R,\sigma, \dot{\sigma})$ and
we will sometimes make the identification $S/I = R$.

A {\it morphism of $O$-frames}\index{morphism of $O$-frames}
$\alpha\colon \mathcal{F} = (S, I,R, \sigma, \dot{\sigma}) \longrightarrow
\mathcal{F'} = (S', I',R', \sigma', \dot{\sigma}')$
is an $O$-algebra homomorphism $\alpha\colon S \longrightarrow S'$
such that $\alpha(I) \subset I'$ and such that 
\begin{displaymath}
  \dot{\sigma}'(\alpha(a)) = \alpha(\dot{\sigma}(a)) , \quad 
  \; a \in I. 
\end{displaymath} 
The last equation implies that 
\begin{displaymath}
  \sigma'(\alpha(s)) = \alpha(\sigma(s)), \quad s \in S.
\end{displaymath}

Let $\mathcal{F} = (S, I,R, \sigma, \dot{\sigma})$ be an
$O$-frame. Then  there exists a unique  element $\theta \in S$ in the
radical of $S$ such that 
\begin{equation}\label{deftheta}
\sigma (a) = \theta \dot{\sigma}(a), \quad \text{  for all
$a\in I$,}
\end{equation}
 cf. \cite[Lem.~3.2]{ACZ}. In the frames 
below we have $\theta = \pi$. 
\begin{example}\label{Wittdisplay} Let $R$ be a $p$-adic $O$-algebra. Then the
  {\it Witt ring $W_{O}(R)$ relative to $O$}\index{Witt ring relative to $O$} 
  with respect to the chosen uniformizer $\pi \in O$ is a $p$-adic
  $O$-algebra. The {\it Witt polynomials relative to $O$}\index{Witt polynomials relative to $O$} 
  \begin{displaymath}
    \mathbf{w}_{O,n} = X_0^{q^n} + \pi X_1^{q^{n-1}} + \pi^{2}X_2^{q^{n-2}} +
    \ldots + \pi^{n-1}X_{n-1}^{q} + \pi^{n}X_n, 
  \end{displaymath}
  define $O$-algebra homomorphisms
  $\mathbf{w}_{O,n}: W_{O}(R) \longrightarrow R$.
  We denote by $F$ and $V$ the Frobenius and the Verschiebung acting
  on $W_{O}(R)$, cf. \cite{Dr}. In the case where $k = R$ is a perfect field, 
  the ring $W_{O}(k)$ is the complete discrete valuation with residue class
  field $k$ which is unramified over $O$.
  
  The {\it Witt frame relative to $O$}\index{Witt frame relative to $O$} for $R$ is the $O$-frame defined as \index[NO]{WBA@$\mathcal{W}_{O}(R)$}
  \begin{equation}
    \mathcal{W}_{O}(R) = (W_{O}(R), VW_{O}(R), R, \sigma, \dot{\sigma}).
  \end{equation}
  Here $\sigma = F: W_O(R) \longrightarrow W_O(R)$ is the Frobenius endomorphism
  written as $\sigma(\xi) = ~^{F}\xi$, and $\dot{\sigma}(^{V}\!\xi) = \xi$, for
  $\xi \in W_O(R)$.  We use also the notation $\dot{F} := \dot{\sigma}$ and
  $I_O(R) = VW_O(R)$. If $K = \mathbb{Q}_p$ and $\pi = p$, we obtain the
  classical ring of Witt vectors $W(R) = W_{\mathbb{Z}_p}(R)$. We write
  $\mathcal{W}(R)$ for the $\mathbb{Z}_p$-frame $\mathcal{W}_{\mathbb{Z}_p}(R)$. 
\end{example}

\begin{example}\label{ex:crystdisp} Let $S \longrightarrow R$ be a surjective
  homomorphism of $p$-adic $O$-algebras. We assume that the kernel
  $\mathfrak{a}$
  is endowed with divided powers relative to $O$ (\cite{ACZ}, 1.2.2). They
  make sense out of the expression $''a^{q}/\pi''$. We also call this an $O$-$pd$-thickening. 
We denote by
  $\mathfrak{a}_{[F^n]}$ the ideal $\mathfrak{a}$ considered as an
  $W_{O}(S)$-module via restriction of scalars relative to
  $\mathbf{w}_{O,n}: W_{O}(S) \longrightarrow S$. The divided powers give rise
  to {\it divided Witt polynomials}\index{divided Witt polynomials} $\dot{\mathbf{w}}_{O,n}$. They are
  homomorphisms of $W_{O}(S)$-modules 
  $\dot{\mathbf{w}}_{O,n}: W_{O}(\mathfrak{a}) \longrightarrow \mathfrak{a}_{[F^n]}$
  such that $\pi^n \dot{\mathbf{w}}_{O,n} = \mathbf{w}_{O,n}$.   
  They give rise to an isomorphism of $W_{O}(S)$-modules 
  \begin{displaymath}
    \prod_{n \geq 0} \dot{\mathbf{w}}_{O,n}:  W_{O}(\mathfrak{a})
    \overset{\sim}{\longrightarrow} \prod_{n \geq 0} \mathfrak{a}_{[F^n]} ,
  \end{displaymath}
  cf. \cite{ACZ}, 1.2.2. 
 The inverse image in $W_{O}(\mathfrak{a})$  of an element
 $[a, 0, 0, \ldots]$ from the right hand side is called the 
 {\it logarithmic Teichm\"uller representative} \index{logarithmic Teichm\"uller representative} of $a\in\mathfrak a$.
 The logarithmic
  Teichm\"uller representatives of elements of $\mathfrak{a}$ form an
  ideal $\tilde{\mathfrak{a}} \subset W_{O}(S)$. 
  The ideal $\mathcal{J} = \tilde{\mathfrak{a}} \oplus I_O(S)$ is the
  kernel of the composition  
  \begin{displaymath}
    W_{O}(S) \overset{\mathbf{w}_{O,0}}{\longrightarrow} S \longrightarrow R.
  \end{displaymath} 
  Then $\dot{F}: I_{O}(S) \longrightarrow W_{O}(S)$ extends uniquely to a
  $F$-linear homomorphism $\dot{F}: \mathcal{J} \longrightarrow W_{O}(S)$ such that
  $\dot{F}(\tilde{\mathfrak{a}}) = 0$.  We define the relative {\it Witt
    frame} \index{Witt frame for $S/R$ relative to $O$} for $S \longrightarrow R$ as  \index[NO]{WBB@$\mathcal{W}_{O}(S/R)$}
  \begin{equation}\label{relWittFrame1e} 
    \mathcal{W}_{O}(S/R) = (W_{O}(S), \mathcal{J}, R,  F, \dot{F}).
  \end{equation}
  This is an $O$-frame. Later we use the more precise notation
  \begin{displaymath}
I_O(S/R) = \mathcal{J} = W_O(\mathfrak{a}) + I_O(S).  
    \end{displaymath}
  
\end{example}
\begin{definition}[\cite{ACZ}, Def. 3.3]\label{Rah2d}
  Let $\mathcal{F}= (S, I, R,\sigma, \dot{\sigma})$ be an $O$-frame. An {\it $\mathcal{F}$-display}\index{display over a frame}
  $\mathcal{P}=(P, Q, F, \dot{F})$ consists of 
  the following data: a finitely generated projective $S$-module $P$,  
  a submodule $Q\subset P$, and two $\sigma$-linear maps
  \begin{displaymath}    
    F\colon P\longrightarrow P, \quad  \dot{F}\colon Q\longrightarrow P.
  \end{displaymath} 

  \noindent
  The following conditions are required.

  \begin{enumerate}
  \item[(i)]  $IP \subset Q.$
  \item[(ii)] The factor module $P/Q$ is a finitely generated projective
    $R$-module. 
  \item[(iii)]  The following relation holds for  $a \in I$ and $x\in P$,
    \begin{displaymath}    
      \dot{F}(a x)= \dot{\sigma}(a) F(x).
    \end{displaymath} 
  \item[(iv)] $\dot{F}(Q)$ generates $P$ as an $S$-module.
  \item[(v)] The projective $R$-module $\Lie \mathcal{P}=P/Q$ lifts to a finitely
    generated projective $S$-module. It is called the \emph{Lie algebra} \index{Lie algebra of a display}of $\CP$. 
  \end{enumerate}
\end{definition}
If the rank of $\Lie \mathcal{P}$ is
constant, we call it the \emph{dimension of $\mathcal{P}$}. If the $S$-module 
$P$ is of constant rank, we call it the \emph{height of $\mathcal{P}$}. If we
want to be precise, we say $\mathcal{F}$-height. \index{dimension of a display}\index{height of a display}

$\CF$-displays form a category in the obvious way. In the case
$O = \mathbb{Z}_p$ and $\mathcal{F} = \mathcal{W}(R)$ for a $p$-adic ring $R$,
we speak simply of a {\it display} over $R$. Displays for general frames $\CF$
were originally called $\CF$-windows, cf. \cite[Def. 3.3]{ACZ}. We note that
for the $O$-frames $\mathcal{W}_O(R)$, the condition (v) of Definition
\ref{Rah2d} is automatically satisfied, cf. \cite[Lem. 2]{Zi}.  

\begin{example}\label{def:multdisplay}
  For each $O$-frame $\mathcal{F} = (S, I, R, \sigma, \dot{\sigma})$ we
  \index[NO]{PBA@$\mathcal{P}_{m, \mathcal{F}}$}
  have the {\it multiplicative $\mathcal{F}$-display}\index{multiplicative $\mathcal{F}$-display}  
  \begin{displaymath}
    \mathcal{P}_m = \mathcal{P}_{m, \mathcal{F}} = (S, I, \sigma, \dot{\sigma}) .
  \end{displaymath}
\end{example}

\begin{example}\label{def:twistdisplay} 
  Let $\mathcal{P}$ be an $\mathcal{F}$-display. Let $\varepsilon \in S$ be a
  unit. The display \index[NO]{PBB@$\mathcal{P}(\varepsilon)$}
  \begin{displaymath}
\mathcal{P}(\varepsilon) = (P, Q, \varepsilon F, \varepsilon \dot{F}).
  \end{displaymath}
is called the \emph{twist} of $\mathcal{P}$ by $\varepsilon$. \index{twist of a display}
\end{example}

Recall the element $\theta$ from \eqref{deftheta}. The conditions in Definition \ref{Rah2d} imply that 
\begin{equation}\label{thetadot}
  F(y) = \theta \dot{F}(y), \quad y \in Q .
\end{equation}

We can always find a direct sum decomposition $P = T \oplus L$ such that
$Q = I T \oplus L$.  Such a decomposition  is called a {\it normal decomposition} \index{normal decomposition of a display} of $P$. The
$\sigma$-linear homomorphism
\begin{equation}\label{Adorf1e} 
  \Phi := F_{|T} \oplus \dot{F}_{|L} \colon T \oplus L \longrightarrow P 
\end{equation}
is a $\sigma$-linear isomorphism, i.e., 
\eqref{Adorf1e} corresponds to the  {\it linearization  isomorphism}, 
\begin{equation}\label{lineariziso}
  F^{\sharp} \oplus \dot{F}^{\sharp}\colon (S \otimes_{\sigma, S} T)
  \oplus (S \otimes_{\sigma, S} L) \isoarrow P.
\end{equation}
Conversely, an arbitrary
$\sigma$-linear isomorphism (\ref{Adorf1e}) defines an
$\mathcal{F}$-display in the obvious way.

For each display $\mathcal{P}$ there is a homomorphism of $S$-modules
(\cite[Def. 3.3]{ACZ})  
\begin{equation}\label{Vsharp1e}
  V^{\sharp}: P \longrightarrow S \otimes_{\sigma, S} P 
\end{equation}
which is uniquely determined by
\begin{displaymath}
  V^{\sharp}(s\dot{F}y) = s \otimes y, \quad V^{\sharp}(Fx) = \theta
  \otimes x, \quad x \in P, \; y \in Q, \; s \in S. 
\end{displaymath}
We have
\begin{displaymath}
  V^{\sharp} \circ F^{\sharp} = \theta\, \id_{S \otimes_{\sigma, S} P},
  \quad F^{\sharp} \circ V^{\sharp} = \theta\, \id_P. 
\end{displaymath}

Any morphism of $O$-frames $\alpha\colon \mathcal{F} \longrightarrow \mathcal{F}'$
defines a base change functor $\alpha_{\ast}$ \index{base change of a display} from the category of
$\mathcal{F}$-displays to the category of $\mathcal{F}'$-displays as
follows, cf. \cite[Def. 3.8]{ACZ}.
Let $\mathcal{P}$ be an $\mathcal{F}$-display. Then we define
$\alpha_{\ast}( \mathcal{P}) = \mathcal{P}' = (P', Q', F', \dot{F}')$ as
follows:
\begin{equation}\label{framechange1e}
  \begin{aligned}
    &  P' = S' \otimes_{S} P, & Q' = \Ker \big(S' \otimes_{S} P \longrightarrow R'
    \otimes_{R} (P/Q)\big), & \quad F' = \sigma' \otimes F\colon
    P' \longrightarrow P'.    
  \end{aligned}
\end{equation}
 Here $Q'$ is the image of
$I' \otimes_S P \oplus S' \otimes_S Q \rightarrow S' \otimes_S P$. 
The map $\dot{F}': Q' \rightarrow P'$ is uniquely determined by
\begin{displaymath}
  \dot{F}'(\xi \otimes x) = \dot{\sigma}'(\xi) \otimes F(x), \; 
  \dot{F}'(\eta \otimes y) = \sigma'(\eta) \otimes \dot{F}(y), \quad
  \text{for} \; \xi \in I', \eta \in S', x \in P, y \in Q. 
  \end{displaymath} 
 If $\mathcal{P}$ is given
in terms of a normal decomposition (\ref{Adorf1e}), we obtain $\mathcal{P}'$ 
from the $\sigma'$-linear extension of $\Phi$,  
\begin{displaymath} 
  \Phi'\colon (S' \otimes_{S} T) \oplus (S' \otimes_{S} L) \longrightarrow P'.  
\end{displaymath}

\begin{example}\label{basechangemult}
The base change of the multiplicative display for the frame $\CF$ under $\alpha\colon\CF\to\CF'$ is the multiplicative display for $\CF'$. 
\end{example}
 If $R$ is a perfect ring of
  characteristic $p$ with an $O$-algebra structure,  the category of
  $\mathcal{W}_O(R)$-displays is
  equivalent to the more classical category of Dieudonn\'e modules.
We describe this equivalence in its natural generality.
\begin{definition}\label{Adorf1ex} 
 \noindent (a) A \emph{perfect $O$-frame} \index{perfect $O$-frame} is an $O$-frame
  $\mathcal{F} = (S, I, R, \sigma, \dot{\sigma})$  such that
  $\sigma: S \longrightarrow S$ is bijective.

 Let $\mathcal{F} = (S, I, R, \sigma, \dot{\sigma})$ be a perfect $O$-frame.
It follows from (\ref{deftheta}) and Definition \ref{defOframe} that
$\dot{\sigma}: I \rightarrow S$ is bijective. Let $u \in I$ such that  
$\dot{\sigma}(u) = 1$. Again by (\ref{deftheta}) we obtain
$\sigma(u) = \theta$. One can see that the elements $u$ and $\theta$ are non
zero divisors in $S$.

  \smallskip
  
  \noindent (b)  A \emph{Dieudonn\'e module} $(M, F, V)$ for the perfect\index{Dieudonn\'e module of a perfect display}
  $O$-frame $\mathcal{F}$ consists of a finitely generated 
  projective $S$-module $M$ and two additive maps $F: M \longrightarrow M$,
  $V: M \longrightarrow M$ such that the following conditions are satisfied.
  \begin{equation*}
  \begin{aligned}
 (i)&\,\,F(sx) = \sigma(s) F(x), \quad V(sx) = \sigma^{-1}(s) V(x), \quad
      x \in P, \; s \in S.\\
 (ii)&
\,\,F \circ V = \theta\, \id_{M}, \quad V \circ F = u\, \id_{M}. \\
(iii)&\,\, \text{The $R $-module $M/VM$ is projective and lifts to a
    finitely generated projective $S$-module. }
  \end{aligned}
  \end{equation*}
  \end{definition}

  If $R$ is a perfect $O$-algebra, then  $\mathcal{F} = \mathcal{W}_O(R)$ is a perfect
  $O$-frame and we have $u = ~^{V}\! 1 = \pi = \theta$.
  \begin{proposition}\label{perfFrame1p}
    Let $\mathcal{F} = (S, I, R, \sigma, \dot{\sigma})$ be a perfect
    $O$-frame. Let $u, \theta \in S$ as defined above. Then the category of Dieudonn\'e modules for $\mathcal{F}$ is
    equivalent to the category of $\mathcal{F}$-displays.
  \end{proposition}
  \begin{proof}
    Let $(M,F,V)$ be a Dieudonn\'e module. Since
    $u$ and $\theta$ are not zero divisors, the maps $F:M \longrightarrow M$
    and $V: M \longrightarrow M$ are injective. Therefore we can define a
    display $(P,Q,F, \dot{F})$ by setting
    \begin{displaymath}
      P = M, \; Q = VM, \; F= F, \; \dot{F} = V^{-1}. 
    \end{displaymath}
    Conversely, let $(P,Q, F, \dot{F})$ be a display. We set
    $(M, F) := (P, F)$. We have the bijective map
    \begin{displaymath}
      \nu: S \otimes_{\sigma, S} P \longrightarrow P, \quad \nu(s \otimes x) :=
      \sigma^{-1}(s) x.  
    \end{displaymath}
    Then we define $V = \nu \circ V^{\sharp}$. More explicitly, we have
    \begin{displaymath}
      V(s \dot{F}(y)) = \sigma^{-1}(s) y, \quad  y \in Q, \; s
      \in S. 
    \end{displaymath}
    This implies that $V(P) = Q$. Moreover, we obtain
    \begin{displaymath}
      \begin{aligned}
        FV(s \dot{F}y) &= F(\sigma^{-1}(s) y) = s F(y) = \theta s
        \dot{F}y\\ 
        VF(x) &= V(\dot{F}(ux)) = ux.   
      \end{aligned}
    \end{displaymath}
    Therefore $(M,F,V)$ is a Dieudonn\'e module.
  \end{proof}
  In our basic example $\mathcal{F} = \mathcal{W}_{O}(R)$ for a perfect
  $O$-algebra $R$, we can replace the condition (iii) above by the weaker condition that $M/VM$
  is a projective $R$-module. We note that for this frame
  $F \circ V = \pi\, \id_{M},  V \circ F = \pi\, \id_{M}$.

  We refer to \cite[Def. 3.3]{ACZ} or \cite{Zi3} for the definition of a 
  \emph{nilpotent} $\CF$-display. \index{nilpotent display}If $R$ is a perfect $O$-algebra, a
  $\mathcal{W}_O(R)$-display is nilpotent iff for the corresponding
  Dieudonn\'e module $(M, F, V)$ the endomorphism $V$ of $M/\pi M$ is
  nilpotent. For an arbitrary $O$-algebra $R$ such that $\pi$ is nilpotent
  in $R$, a $\mathcal{W}_O(R)$-display $\mathcal{P}$ is nilpotent iff for any
  homomorphism of $O$-algebras to a perfect field $R \longrightarrow k$, the base change of $\mathcal{P}$ by the morphism of frames
  $\mathcal{W}_O(R) \longrightarrow \mathcal{W}_O(k)$ is nilpotent.  

\begin{definition}\label{def:str}
  Let $R$ be an $O$-algebra. Let $X$ be a $p$-divisible group over $R$
  endowed with a $\mathbb{Z}_p$-algebra homomorphism
  $\iota: O \longrightarrow \End X$. We call the action $\iota$ \emph{strict}\index{strict action}
  if the induced action on $\Lie X$ coincides with the $O$-action on this
  $R$-module given by restriction of scalars $O \longrightarrow R$. We say that
  $(X, \iota)$ is a strict $p$-divisible $O$-module.  
\end{definition}

The following main result of \cite{ACZ} was known before for
$O = \mathbb{Z}_p$ \cite{Zi}, \cite{L3}.  
\begin{theorem}[\cite{ACZ}, Thm. 1.1]\label{maindisp}
  Let $R \in \Nilp_O$. There is an equivalence of categories 
  \begin{equation*}\label{DspFunktor1e}
     \text{\big(nilpotent $\mathcal{W}_O(R)$-displays\big)} \longrightarrow  
 \text{\big(strict formal $p$-divisible $O$-modules over $R$\big)} 
  \end{equation*}
  which is functorial in $R$.
\end{theorem} 
The theorem
extends to $p$-adic $R$ if we require the properties ''nilpotent'' and
''formal'' only after base change to $R/pR$. 

A nilpotent $\mathcal{W}_O(R)$-display gives rise to a crystal, as
follows. Let $S \longrightarrow  R$ be a  $O$-$pd$-thickening, cf. 
Example \ref{ex:crystdisp}. We assume that $p$ is nilpotent in $S$. 
The ring homomorphism $W_O(S) \longrightarrow W_O(R)$
defines morphisms of $O$-frames,

\begin{displaymath}
  \mathcal{W}_O(S) \longrightarrow \mathcal{W}_O(S/R) \longrightarrow
  \mathcal{W}_O(R). 
\end{displaymath}
The base change of displays with respect to the first arrow goes as follows.
Let $\mathcal{P} = (P, Q, F, \dot{F})$ be a $\mathcal{W}_O(S)$-display and
let $\mathcal{P}' = (P', Q', F',\dot{F}')$ be the $\mathcal{W}_O(S/R)$-display
obtained by the base change functor defined before Example
\ref{basechangemult}. Explicitly, it is given as follows: $P' = P$, and 
$Q' = Q + \mathcal{J} P = Q \oplus \tilde{\mathfrak{a}}$, and $F' = F$, and
$\dot{F}'_{\mid Q} = \dot{F}$, and $\dot{F}'(ax) = 0$ for 
$a \in \bar{\mathfrak{a}}$. Note that the last equation is necessary because
$\dot{F}'(ax) = \dot{F}(a) F'(x) = 0$, since $\dot{F}(a) = 0$ by the definition
of (\ref{relWittFrame1e}). 
\begin{theorem}[\cite{L2}, \cite{Zi}]\label{DspKristall1t}
  Let $S \longrightarrow R$ be an $O$-$pd$-thickening such that $p$ is
  nilpotent in $S$. The base change functor
  \begin{displaymath}
    \big(\text{nilpotent $\mathcal{W}_O(S/R)$-displays}\big) \longrightarrow
    \big(\text{nilpotent $ \mathcal{W}_O(R)$-displays}\big)  
  \end{displaymath}
  is an equivalence of categories. \qed
\end{theorem} 
\begin{remark}\label{LaufunctorRm}
In the case $O = \mathbb{Z}_p$, Lau \cite{L} has defined a functor
\begin{equation}\label{Laufunctor1e} 
  \big(p\text{\it-divisible groups over } R\big) \rightarrow \big(\mathcal{W}(S/R)\text{\it-displays}\big)
\end{equation}
which gives a quasi-inverse of the functor in  Theorem \ref{DspKristall1t}, when restricted
to formal $p$-divisible groups. In particular this functor associates to
an arbitrary $p$-divisible group over $R$ a display.  
\end{remark}
Let $\mathcal{P}$ be a nilpotent $\mathcal{W}_O(R)$-display. Let
$\tilde{\mathcal{P}}$ be the unique $\mathcal{W}_O(S/R)$-display
associated to $\mathcal{P}$ by Theorem \ref{DspKristall1t}. Then we set
\index[NO]{DDA@$\mathbb{D}_{\mathcal{P}}$}
\begin{equation}\label{displaycrystal1e} 
  \mathbb{D}_{\mathcal{P}}(S) = \tilde{P}/I_O(S)\tilde{P}. 
\end{equation}
This is a finitely generated projective $S$-module. It is a crystal in the
following sense. If $S' \longrightarrow R$ is a another $O$-$pd$-thickening such
that $p$ in nilpotent in $S'$ and $S' \longrightarrow S$ is a morphism of
$O$-$pd$-thickenings, then there is a canonical isomorphism
\begin{displaymath}
  S \otimes_{S'} \mathbb{D}_{\mathcal{P}}(S') \cong
  \mathbb{D}_{\mathcal{P}}(S). 
\end{displaymath}
This crystal corresponds to the Grothendieck-Messing crystal of a
$p$-divisible group via Theorem \ref{maindisp}. From Theorem
\ref{DspKristall1t} one obtains the Grothendieck-Messing criterion for
displays in the following formulation. 
\begin{corollary}\label{GM1c}
  Let $\mathcal{P}$ be a nilpotent $\mathcal{W}_O(R)$-display. Let
  $S \longrightarrow R$ be an $O$-$pd$-thickening. Each
  $\mathcal{W}_O(S)$-display $\tilde{\mathcal{P}}$ which lifts $\mathcal{P}$
  defines a lifting
  $\widetilde{\mathrm{Fil}} := \tilde{Q}/I_O(S)\tilde{P} \subset
  \mathbb{D}_{\mathcal{P}}(S)$
  of the  Hodge filtration 
  $\mathrm{Fil} := Q/I_O(R)P \subset \mathbb{D}_{\mathcal{P}}(R)$. 

For a fixed   $O$-$pd$-thickening $S \longrightarrow R$, consider the category of pairs $(\mathcal{P}, \widetilde{\mathrm{Fil}})$,
where $\mathcal{P}$ is a nilpotent display and
$\widetilde{\mathrm{Fil}} \subset \mathbb{D}_{\mathcal{P}}(S)$ is a lifting of
the Hodge filtration associated to $\mathcal{P}$. The functor which maps a pair $(\mathcal{P}, \tilde{\mathcal{P}})$  
to the pair $(\mathcal{P}, \widetilde{\mathrm{Fil}})$ is an equivalence 
of categories. \qed
\end{corollary}
\begin{proof}
  Let $\mathcal{P}'$ the unique $\mathcal{W}_O(S/R)$-display which corresponds
  to $\mathcal{P}$ by Theorem \ref{DspKristall1t}. We note that
  $P = W_O(R) \otimes_{W_O(S)} P'$ by definition of the base
  change. A lifting of the Hodge-filtration
  $P \rightarrow P/Q = P'/Q'$ corresponds to a $\mathcal{W}_O(S)$-module
  $\tilde{Q}$ such that 
  $P' \supset \tilde{Q} \supset I_O(S) P'$. Since $\tilde{Q} \subset Q'$,
  we obtain a $\mathcal{W}_O(S)$-display $\tilde{\mathcal{P}}$ by restricting 
  $\dot{F}'$ to $\tilde{Q}$, i.e.,
  $\tilde{\mathcal{P}} = (P', \tilde{Q}, F', \dot{F}')$. From this the
  equivalence of categories follows. 
  \end{proof} 
The following fact is well-known, but we give a proof.

\begin{lemma}\label{Displaykristall1l}
  Let $R$ be a $p$-adic $O$-algebra. 
  Let $\mathcal{P}$ be a $\mathcal{W}_O(R)$-display. Let $\tilde{O}$ be
  a discrete valuation ring which is a finite extension of $O$. Let
  \begin{displaymath}
    \tilde{O} \longrightarrow \End \mathcal{P}, 
  \end{displaymath}
  be an $O$-algebra homomorphism. Then $P$ is a locally on $\Spec R$ a
  free $\tilde{O} \otimes_{O} W_O(R)$-module. 

  Let $S \longrightarrow R$ be an $O$-$pd$-thickening such that $p$ is
  nilpotent in $S$. We assume that $\mathcal{P}$ is nilpotent.
  Then $\mathbb{D}_{\mathcal{P}}(S)$ is locally on $\Spec S$ a free
  $\tilde{O} \otimes_{\mathbb{Z}_p} S$-module.   
\end{lemma}

\begin{proof}
  We start with the case where $S = R = k$ is a perfect field which
  contains the residue class field of $\tilde{O}$. Let $\tilde{O}^t$ be
  the maximal unramified extension of $O$ contained in $\tilde{O}$. Let
  $\sigma$ be the Frobenius automorphism of $\tilde{O}^t$ relative to
  $O$. To each $O$-algebra homomorphism
  $\psi: \tilde{O}^t \longrightarrow k$ there is a unique Frobenius
  equivariant $O$-algebra homomorphism 
  \begin{displaymath}
    \tilde{\psi}: \tilde{O}^t \longrightarrow W_O(k) 
  \end{displaymath}
  which induces $\psi$ when composed with
  $\mathbf{w}_{O,0}: W_O(k) \longrightarrow k$. This follows from the remark after
  the definition of $W_O(R)$, cf. Example \ref{Wittdisplay}.   The decomposition 
  \begin{displaymath}
    \tilde{O} \otimes_{O} W_{O}(k) = \prod_{\psi} \tilde{O}
    \otimes_{\tilde{O}^t, \tilde{\psi}} W_O(k)  
  \end{displaymath}
  induces a decomposition
  \begin{displaymath}
    P = \oplus_{\psi} P_{\psi}. 
  \end{displaymath}
  Each $P_{\psi}$ is a free module over the discrete valuation ring
  $\tilde{O} \otimes_{\tilde{O}^t, \tilde{\psi}} W_O(k)$. The Frobenius
  $F\colon P\longrightarrow P$ induces maps 
  $P_{\psi} \longrightarrow P_{\psi \sigma}$. This shows that all $P_{\psi}$
  have the same rank as $W_O(k)$-modules. This proves the case where $R = k$
  is a perfect field containing the residue class field of $\tilde{O}$. 
  
  Now let $k$ be an arbitrary field of characteristic $p$. It suffices to
  show that $P \otimes_{W_O(k)} k$ is a free $\tilde{O} \otimes_{O} k$-module.
  By base change
  $\alpha\colon k \longrightarrow \bar{k}$ this follows from the previous case because
  $P \otimes_{W_O(k)} \bar{k}=\alpha_*(P)\otimes \bar{k}$ is a free
  $\tilde{O} \otimes_{O} \bar{k}$-module.

  If $R$ is a local ring with residue class field of characteristic $p$
  we conclude by Nakayama's lemma that $P \otimes_{W_O(R)} R$ is a free
  $\tilde{O} \otimes_{O} R$-module. The generalization
  to arbitrary $R$ is immediate. This proves the first assertion of the Lemma.

  If $\mathcal{P}$ is nilpotent, the crystal 
  $\mathbb{D}_{\mathcal{P}}$ is defined. The case
  $\mathbb{D}_{\mathcal{P}}(R) = P/I_O(R)P$ was proved above. For
   arbitrary $S \longrightarrow R$ we can apply  again Nakayama's
  lemma. 
\end{proof}
\begin{remark}\label{Displaykristall1Rm}
  Let $R$ be a ring such that $p$ is nilpotent in $R$. Let $X$ be a $p$-divisible group over $R$ with a ring homomorphism
  \begin{displaymath}
    O \longrightarrow \End X. 
  \end{displaymath}
  Let $S \longrightarrow R$ be a nilpotent $pd$-thickening. Then the value
  of the Grothendieck-Messing crystal $\mathbb{D}_{X}(S)$ is locally on
  $\Spec S$ a free $O \otimes_{\mathbb{Z}_p} S$-module. This can be
  shown by the same arguments as above. 
\end{remark}

Finally we discuss isogenies of $\mathcal{W}_O(R)$-displays, where $R$
is an $O$-algebra such that $p$ is nilpotent in $R$.  We assume
moreover that $\Spec R$ is connected. 
Let
$\alpha: \mathcal{P}_1 \longrightarrow \mathcal{P}_2$ be a morphism of
displays of the same height and dimension, cf. the remark after Definition
\ref{Rah2d}. 
Locally on $\Spec R$ the
$W_O(R)$-modules $P_1$ and $P_2$ are free of the same rank. We may
choose a basis in each of these modules and write
$\det \alpha \in {W}_O(R)$. This is locally defined up to a
unit in ${W}_O(R)$. More invariantly one can write exterior
powers.
\begin{definition}
A morphism of $\CW_O(R)$-displays of the same height and dimension 
$\alpha: \mathcal{P}_1 \longrightarrow \mathcal{P}_2 $ is called an isogeny\index{isogeny of displays}
  if $\det \alpha \neq 0$.  
\end{definition}
\begin{proposition}[\cite{Zi3}, Prop. 17.6.2.]\label{Isg1p}
  Let $R$ be an $O$-algebra such that $p$ is nilpotent in $R$ and such that
  $\Spec R$ is connected. 
  Let $\alpha: \mathcal{P}_1 \longrightarrow \mathcal{P}_2$ be an isogeny of
  $\mathcal{W}_{O}(R)$-displays.   Then there exists a natural number ${\tt h} \in \mathbb{Z}_{\geq 0}$ such that locally
  on $\Spec R$\index{height of an isogeny of displays} 
  \begin{displaymath}
    \det \alpha = \pi^{\tt h} \epsilon, \quad \epsilon \in
    {W}_O(R)^{\times}.
  \end{displaymath}\qed
\end{proposition}
We call ${\tt h}$ the \emph{$O$-height} of $\alpha$, and write ${\tt h}=\height_O \alpha$.
If $O = \mathbb{Z}_p$, we write simply $\height \alpha$. An abbreviation for the
Proposition is:
\begin{displaymath}
  \height_O \alpha = \ord_{\pi} \det \alpha. 
\end{displaymath}
\begin{proposition}[\cite{Zi3}, Prop. 17.6.4.] 
  Assume that the ideal of nilpotent elements in $R$ is nilpotent and
  that $\Spec R$ is connected.   
  Let $\alpha: \mathcal{P}_1 \longrightarrow \mathcal{P}_2$ be an isogeny of
  $O$-height ${\tt h}$. Then there exists locally on $\Spec R$ a morphism of
  $\mathcal{W}_O(R)$-displays 
  $\beta: \mathcal{P}_2 \longrightarrow \mathcal{P}_1$ such that
  \begin{displaymath}
    \beta \circ \alpha = \pi^{\tt h}\,\id_{\mathcal{P}_1}, \quad
    \alpha \circ \beta = \pi^{\tt h}\,\id_{\mathcal{P}_2}. 
  \end{displaymath}\qed
\end{proposition}

\begin{proposition}
  With the assumptions of Proposition \ref{Isg1p}, let $a: X_1 \longrightarrow X_2$
  be a morphism of strict formal $p$-divisible $O$-modules over $R$. Let $\alpha: \mathcal{P}_1 \longrightarrow \mathcal{P}_2$ be the induced morphism of the associated $\CW_O(R)$-displays, cf. Theorem
  \ref{maindisp}. The morphism $a$ is an isogeny of height ${\tt h}$ if and only if $\alpha$
is an isogeny of height ${\tt h}$. 
\end{proposition}
\begin{proof}
  This can be reduced to the case of a perfect field $R = k$ where it is
  well-known by Dieudonn\'e theory. 
  \end{proof}
Let $R$ be an $O$-algebra and
let $a\colon X_1\to X_2$ be a morphism of strict formal $p$-divisible $O$-modules. By Theorem
\ref{maindisp}, there is an associated morphism $\alpha: \mathcal{P}_1
\longrightarrow \mathcal{P}_2$ of $\mathcal{W}_{O}(R)$-displays. We set
\begin{equation}\label{def:height}
  \height_{O} a = \height_{O} \alpha, \quad \height_{O} X_1 =
  \height_{\mathcal{W}_O(R)} \mathcal{P}_1. 
\end{equation}
The last height was defined after Definition \ref{Rah2d}. It is equal to the
$O$-height of the endomorphism of $\mathcal{P}_1$ given by multiplication by
$\pi$. We also write
\begin{displaymath}
  \height_{O} \mathcal{P}_1 = \height_O (\pi | \mathcal{P}_1) = 
  \height_{\mathcal{W}_O(R)} \mathcal{P}_1. 
  \end{displaymath}

\subsection{Bilinear forms of displays}\label{ss:Bfod} 
Let $\mathcal{F} = (S,I,R,\sigma, \dot{\sigma})$ be an $O$-frame and
let $\theta \in S$ be the element from \eqref{deftheta}.  

\begin{definition}\label{BF1d} 
      Let\index{bilinear form of $\mathcal{F}$-displays}
      $\mathcal{P}_1, \mathcal{P}_2, \mathcal{P}$ be $\mathcal{F}$-displays.
      A {\it bilinear form of $\mathcal{F}$-displays}  
      \begin{equation*}
        \beta: \mathcal{P}_1 \times \mathcal{P}_2 \longrightarrow \mathcal{P} 
      \end{equation*}
      is a bilinear form of $S$-modules
      \begin{equation}\label{bilin1e}
        \beta: P_1 \times P_2 \longrightarrow P
      \end{equation}
      with the following properties: 
      \begin{enumerate}
      \item[(i)] The restriction of $\beta$ to $Q_1 \times Q_2$ takes values
        in $Q$.
      \item[(ii)] For $y_1 \in Q_1$ and $y_2 \in Q_2$,
        \begin{displaymath}
          \dot{F} \beta (y_1, y_2) = \beta (\dot{F}_1y_1, \dot{F}_2 y_2).
        \end{displaymath} 
      \end{enumerate}
      We will denote the $O$-module of all bilinear forms by
      \begin{displaymath}
\mathrm{Bil}(\mathcal{P}_1 \times \mathcal{P}_2, \mathcal{P}). 
        \end{displaymath}
\end{definition} 
      \begin{lemma}
        The following equations hold 
        \begin{displaymath}
          \begin{aligned}
          F \beta (x_1, y_2) &= \beta ( F_1x_1, \dot{F}_2 y_2), &\quad  
          x_1 \in P_1, \; y_2 \in Q_2,\\
          F \beta (y_1, x_2) &= \beta (\dot{F}_1y_1, F_2 x_2),  &\quad 
          y_1 \in Q_1, \; x_2 \in P_2,\\[2mm] 
          \theta F \beta (x_1, x_2) &= \beta (F_1x_1, F_2 x_2), &\quad  
          x_1 \in P_1, \; x_2 \in P_2.
          \end{aligned}
          \end{displaymath}
      \end{lemma}
      \begin{proof}
        We omit the verification which is, for classical displays, contained in
        \cite{Zi}. 
        \end{proof} 
Let $R$ be a perfect $O$-algebra and let $\mathcal{F} = \mathcal{W}_O(R)$.
Then we may equivalently consider Dieudonn\'e modules $(P, F, V)$ and
$(P_i, F_i, V_i)$ for $i = 1,2$, cf. Proposition \ref{perfFrame1p}. 
We can reformulate the Definition \ref{BF1d} as follows: A bilinear form
of Dieudonn\'e modules is a bilinear form of $W_O(R)$-modules
$\beta: P_1 \times P_2 \longrightarrow P$ such that
\begin{equation}\label{bilin2e}
\beta(V_1 x_1, V_2 x_2) = V \beta(x_1, x_2). 
  \end{equation}

\begin{proposition}\label{basebil}
  Let $\beta: \mathcal{P}_1 \times \mathcal{P}_2 \longrightarrow \mathcal{P}$
  be a bilinear form of $\mathcal{F}$-displays. Let
  $\alpha: \mathcal{F} \longrightarrow \mathcal{F}'$ be a morphism of frames. Denote by $\mathcal{P}'_1$, $\mathcal{P}'_2$, and $\mathcal{P}'$ the
  displays obtained by base change with respect to $\alpha$. Because
  $P'_i = S' \otimes_{S} P_i$, and $P' = S' \otimes_{S} P$, there is an  induced 
   $S'$-bilinear form $\beta': P'_1 \times P'_2 \longrightarrow P'$. This is  a bilinear form of $\mathcal{F}'$-displays
  \begin{displaymath}
\mathcal{P}'_1 \times \mathcal{P}'_2 \longrightarrow \mathcal{P}'.
    \end{displaymath}
\end{proposition}
\begin{proof}
We omit the straightforward verification.
  \end{proof}

Let $\mathcal{P} = (P, Q, F, \dot{F})$ be an $\mathcal{F}$-display. We 
are going to define the {\it dual $\CF$-display}
$\mathcal{P}^{\vee} = (P^{\vee}, Q^{\vee}, F^{\vee}, \dot{F}^{\vee})$.
For an $S$-module $M$, we define $M^{\ast} = \Hom_{S}(M, S)$. We set
$P^{\vee} :=  P^{\ast}$, and 
\begin{displaymath}
  Q^{\vee} = \{ \psi \in P^{\vee} \; | \; \psi (Q) \subset I \} .
\end{displaymath}
We note that we have a natural perfect pairing
\begin{displaymath}
  P/IP \times P^{\vee}/I P^{\vee} \longrightarrow R.
\end{displaymath}
We deduce that $Q^{\vee}/IP^{\vee}$ is the orthogonal complement
of $Q/IP$ and is therefore a direct summand of $P/IP$. We claim that there are $\sigma$-linear maps 
\begin{displaymath}
  F^{\vee}\colon  P^{\vee} \longrightarrow P^{\vee}, \quad \dot{F}^{\vee}\colon
  Q^{\vee} \longrightarrow P^{\vee}   
\end{displaymath}
which are uniquely determined by the
following conditions. We denote by $<\; , \; > \colon P \times P^{\vee}
\longrightarrow S$ the natural perfect pairing.  Then we require for
$x \in P, \; y \in Q, \phi \in P^{\vee}, \;  \psi \in Q^{\vee}$: 
\begin{equation}\label{dual1e}
  \begin{aligned}
    < \dot{F}(y), F^{\vee}(\phi)> &= \sigma (< y, \phi>), &
    <F(x), F^{\vee}(\phi)> &= \theta \sigma (< x, \phi >),  \\
    <\dot{F}(y), \dot{F}^{\vee}(\psi)> &= \dot{\sigma}(< y, \psi >), &
    <F(x), \dot{F}^{\vee}(\psi)> &= \sigma (<x, \psi >).
  \end{aligned}
\end{equation}
Since $\dot{F}$ is a $\sigma$-linear surjection,  the maps $F^{\vee}$ and 
$\dot{F}^{\vee}$ are uniquely determined by these identities.
To verify the existence of these maps,  we consider a normal
decomposition,  
\begin{displaymath}
  P = T \oplus L. 
\end{displaymath}
Let $L^{\vee} \subset P^{\vee}$ be the orthogonal complement of
$L$ and  $T^{\vee} \subset P^{\vee}$ the orthogonal complement of
$T$. Hence  there are canonical isomorphisms
\begin{displaymath}
  L^{\vee} \cong T^{\ast}, \quad T^{\vee} \cong L^{\ast}.  
\end{displaymath} 
We obtain the normal decomposition
\begin{displaymath}
  P^{\vee} = T^{\vee} \oplus L^{\vee}, \quad Q^{\vee} =
  IT^{\vee} \oplus L^{\vee}.
\end{displaymath}
For  $\psi \in L^{\vee} =T^{\ast}$, we set 
\begin{displaymath}
  \dot{F}^{\vee}(\psi)(\dot{F}(\ell)) = 0, \quad
  \dot{F}^{\vee}(\psi)(\dot{F}(t)) = \sigma(\psi(t)), \quad
  \ell \in L, \; t \in T.  
\end{displaymath}
This definition makes sense because of the linearization  isomorphism
\eqref{lineariziso}. 
Finally, we define
$F^{\vee}(\phi)$ for $\phi \in T^{\vee} = L^{\ast}$ by the equations 
\begin{displaymath}
  F^{\vee}(\phi)(\dot{F}(\ell)) = \sigma(\phi(\ell)), \quad
  F^{\vee}(\phi)(F(t)) = 0. 
\end{displaymath}
One verifies that, with these definitions, the identities 
(\ref{dual1e}) are satisfied. It follows from the symmetry of the equations (\ref{dual1e})  
that we have a natural isomorphism
\begin{displaymath}
\mathcal{P} \cong (\mathcal{P}^{\vee})^{\vee}. 
\end{displaymath}
By the equations (\ref{dual1e}) we have a natural bilinear form of displays
\begin{equation}\label{dual3e}
\mathcal{P} \times \mathcal{P}^{\vee} \longrightarrow \mathcal{P}_m  
\end{equation}
with values in the multiplicative display
$\mathcal{P}_m = \mathcal{P}_{m,\mathcal{F}}$. If $\mathcal{P}'$ is another
$\mathcal{F}$-display, the bilinear form (\ref{dual3e}) induces an isomorphism
\begin{equation}\label{dual4e}
  \Hom_{\text{${\mathcal{F}}$-displays}}(\mathcal{P}', \mathcal{P}^{\vee}) \isoarrow
  \mathrm{Bil}(\mathcal{P}' \times \mathcal{P}, \mathcal{P}_m). 
  \end{equation}

We deduce a variant of the Grothendieck-Messing
criterion. Let $\mathcal{P}$ and $\mathcal{P}'$ be
$\mathcal{W}_{O}(R)$-displays such that $\mathcal{P}^{\vee}$ and $\mathcal{P}'$ 
are nilpotent. Let  $S \longrightarrow R$ be a $O$-$pd$-thickening in $\Nilp_O$,
cf. Example \ref{ex:crystdisp}. 
We denote by $\mathcal{P}_{\rm rel}^{\vee}$ and $\mathcal{P}'_{\rm rel}$
the associated $\mathcal{W}_O(S/R)$-displays, which exist by Theorem
\ref{DspKristall1t}. We define
$\mathcal{P}_{\rm rel} = (\mathcal{P}_{\rm rel}^{\vee})^{\vee}$, where the last $~^{\vee}$
denotes the dual in the category of $\mathcal{W}_O(S/R)$-displays.
We set $\mathbb{D}_{\mathcal{P}}(S) = P_{\rm rel}/I(S) P_{\rm rel}$. Then we obtain a
crystal which is dual to the crystal $\mathbb{D}_{\mathcal{P}^{\vee}}(S)$, cf.
(\ref{displaycrystal1e}). This crystal agrees with $\mathbb{D}_{\mathcal{P}}(S)$
defined earlier, if $\mathcal{P}$ is nilpotent. 
It follows form (\ref{dual4e}) that each bilinear
from
\begin{equation}
  \beta: \mathcal{P}' \times \mathcal{P} \longrightarrow
  \mathcal{P}_{m,\mathcal{W}_{O}(R)}
  \end{equation}
induces a bilinear form
\begin{displaymath}
  \beta_{\rm rel}: \mathcal{P}'_{\rm rel} \times \mathcal{P}_{\rm rel} \longrightarrow 
  \mathcal{P}_{m, \mathcal{W}_{O}(S/R)}
  \end{displaymath}
and, in particular, a $S$-bilinear form 
\begin{equation}
  \beta_{\rm crys}: \mathbb{D}_{\mathcal{P}'}(S) \times \mathbb{D}_{\mathcal{P}}(S)
  \longrightarrow S.
\end{equation}

\begin{proposition}\label{GM2p}
  Let $R \in \Nilp_O$ and let $S \longrightarrow R$ be an $O$-$pd$-thickening in
  $\Nilp_O$. Let $\mathcal{P}$ and $\mathcal{P}'$ be
  $\mathcal{W}_{O}(R)$-displays and
  assume that $\mathcal{P}^{\vee}$ and $\mathcal{P}'$ are nilpotent. Let
  $\tilde{\mathcal{P}}$ and $\tilde{\mathcal{P}}'$
  be liftings which correspond to liftings of the Hodge filtrations
  $\widetilde{\mathrm{Fil}} \subset \mathbb{D}_{\mathcal{P}}(S)$ and
  $\widetilde{\mathrm{Fil}}' \subset \mathbb{D}_{\mathcal{P}'}(S)$, cf. Corollary \ref{GM1c}. Then a bilinear form
  $\beta: \mathcal{P}' \times \mathcal{P} \longrightarrow
  \mathcal{P}_{m,\mathcal{W}_{O}(R)}$
  lifts to a bilinear form
  $\tilde{\beta}: \tilde{\mathcal{P}}' \times \tilde{\mathcal{P}} \longrightarrow
  \mathcal{P}_{m,\mathcal{W}_{O}(S)}$ iff
  \begin{displaymath}
\beta_{\rm crys}(\widetilde{\mathrm{Fil}}', \widetilde{\mathrm{Fil}}) = 0.
    \end{displaymath}\qed
\end{proposition}
\begin{proof}
This is a consequence of Corollary \ref{GM1c} and (\ref{dual4e}). 
  \end{proof}

We go back to an arbitrary $O$-frame $\mathcal{F}$ and add a remark on the map
$V^{\sharp}$, cf.  (\ref{Vsharp1e}).  
If $P$ is an $S$-module we set
$P^{(\sigma)} = S \otimes_{\sigma, S} P$. If $P$ is projective and
finitely generated,  the perfect pairing  $<\; , \; >$ induces a perfect
pairing 
\begin{equation*}
  \begin{aligned}
    <\; , \; >_{(\sigma)}: \;  & P^{(\sigma)} \times (P^{\ast})^{(\sigma)} &
    \longrightarrow &  S\\ 
  &  (s_1 \otimes x,   s_2 \otimes \phi) & \mapsto & s_1s_2 \sigma(\phi(x)).
  \end{aligned}
\end{equation*}
Let $\mathcal{P}$ be an $\mathcal{F}$-display and let $\mathcal{P}^{\vee}$ be
the dual display. 
The maps $(F^\vee)^{\sharp}$ and $V^{\sharp}$ are dual in the following sense 
\begin{displaymath}
  <V^{\sharp}x, s \otimes \phi >_{(\sigma)} = <x, (F^\vee)^{\sharp}(s \otimes \phi)>.
\end{displaymath}

\index{polarization of an $\mathcal{F}$-display}\index{height of a polarization of an $\mathcal{F}$-display}\index{principal polarization of an $\mathcal{F}$-display}
    \begin{definition}\label{def:pol}
      A {\it polarization} of an $\mathcal{F}$-display $\mathcal{P}$ is a
      bilinear form
      \begin{equation*}
        \beta: \mathcal{P} \times \mathcal{P} \longrightarrow
        \mathcal{P}_{m,\mathcal{F}}
      \end{equation*} 
      such that the underlying bilinear form $P \times P \longrightarrow S$ is
      alternating and its determinant is non-zero.

      If $\mathcal{F} = \mathcal{W}_{O}(R)$,  the \emph{height} of $\beta$ is the
      height of the associated isogeny
      $\mathcal{P} \longrightarrow \mathcal{P}^{\vee}$ (cf. Proposition \ref{Isg1p}).
      We write $\height_O \beta$ for the height of $\beta$. Then we have
      $\height_O \beta = \ord_{\pi} \det \beta$ in the notation of 
      Proposition \ref{Isg1p}.  The polarization is called \emph{principal} if $\height_O(\beta)=0$. 
    \end{definition}

    \begin{remark}\label{OFdual} 
      Let $X$ be a  strict formal  $p$-divisible $O$-module over \index{dual of a strict formal  $p$-divisible $O$-module}
      $R \in \Nilp_O$. Let $\mathcal{P}$ be the $\mathcal{W}_{O}(R)$-display
      of $X$ in the sense of Theorem \ref{maindisp}.
      If the dual display $\mathcal{P}^{\vee}$ is nilpotent, it corresponds
      to a strict formal $p$-divisible $O$-module $X^{\vee}$, called
       the {\it $O$-dual} of $X$. In this case,  a polarization is given
      by an anti-symmetric $O$-module homomorphism $X \longrightarrow X^{\vee}$.   
    \end{remark}

    \subsection{The Ahsendorf functor}\label{ss:ahsfunctor}
    We will give here an alternative definition of the Ahsendorf functor of 
    \cite{ACZ} which is better suited to our purposes. One step of this
    definition is contained in the Appendix of \cite{M}. We use
    a Lubin-Tate frame introduced by Mihatsch in loc.~cit., but for us it will be important
    to make a specific choice, cf. Definition \ref{LTF3p}. 

Let $\mathbb{Q}_p \subset \mathfrak{k} \subset K$ be a subfield. We denote
by $\mathfrak{o}$ the ring of integers in $\mathfrak{k}$. We fix a prime
element $\varpi \in \mathfrak{o}$. If $R$ is an $\mathfrak{o}$-algebra we
denote by $W_{\mathfrak{o}}(R)$ the Witt vectors relative to $\mathfrak{o}$ and
$\varpi$. The Frobenius and the Verschiebung will be denoted by $\mathfrak{f}$
and $\mathfrak{v}$.  We set $[K:\mathfrak{k}] = ef$ where $e$ is the
ramification index and $f$ is the inertia index. Beginning with 
section \ref{s:tcf} we will only consider  the case where
$\mathfrak{k} = \mathbb{Q}_p$.

    For an $O$-algebra $R$ we have the Drinfeld homomorphism \index{Drinfeld homomorphism}
    \begin{equation}\label{Drinfmap2e}
\mu: W_{\mathfrak{o}}(R) \longrightarrow W_{O}(R) ,  
\end{equation}
cf. \cite[Prop. 1.2]{Dr}. It is functorial in $R$ and satisfies
$\mathbf{w}_{O,n}(\mu(\xi)) = \mathbf{w}_{\mathfrak{o},fn}(\xi)$, for
$\xi \in W_{\mathfrak{o}}(R)$. This implies the following properties: 

\begin{equation}\label{Drinfmap1e}
  \mu(~^{\mathfrak{f}^{f}}\xi) = ~^{F}\mu(\xi), \quad \mu(~^{\mathfrak{v}}\xi) =
  \frac{\varpi}{\pi} ~^{V}(\mu(~^{\mathfrak{f}^{f-1}}\xi)), \quad
  \mu([u]) = [u],  
\end{equation}
for $\xi \in W_{\mathfrak{o}}(R)$, $u \in R$. The last equation says that 
the Teichm\"uller representative $[u] \in W_{\mathfrak{o}}(R)$, is
mapped by $\mu$ to the Teichm\"uller representative $[u] \in W_{O}(R)$.

We have $\mu(I_{\mathfrak{o}}(R)) \subset I_{O}(R)$. Therefore we may
rewrite the second equation of (\ref{Drinfmap1e}) as 
\begin{equation}\label{Adorf16e}
  ~^{\dot{F}}(\mu(\eta)) = \frac{\varpi}{\pi}
  \mu(~^{\mathfrak{f}^{f-1} \dot{\mathfrak{f}}}\eta), \quad
  \eta \in I_{\mathfrak{o}}(R). 
\end{equation}
The following definition extends Definition \ref{def:str} to the relative case. 

\begin{definition}\label{Adorf1d}
Let $\mathcal{F} = (S,I,R,\sigma, \dot{\sigma})$ be an $\mathfrak{o}$-frame, where $R$ is a $p$-adic $O$-algebra. Let
$\mathcal{P} = (P,Q,F, \dot{F})$ be an $\CF$-display.
\emph{A strict $O$-action on $\mathcal{P}$} is a homomorphism of
$\mathfrak{o}$-algebras $O \longrightarrow \End \mathcal{P}$
such that the induced action on the $R$-module $P/Q$ coincides with
the $O$-module structure on $P/Q$ obtained by restriction of
scalars $O \longrightarrow R$.
\end{definition} 
For a $p$-adic $O$-algebra $R$ we will define a functor
\begin{equation}\label{Adorf11e}
  \mathfrak{A}_{O/\mathfrak{o}, R}:
\left(
  \begin{array}{l}
    \text{$\mathcal{W}_{\mathfrak{o}}(R)$-displays}\\
    \text{with strict $O$-action} 
    \end{array}
\right)
  \; \longrightarrow \; \Big(\text{$\mathcal{W}_{O}(R)$-displays}\Big).  
\end{equation}
We call this functor \emph{the Ahsendorf functor}. \index{Ahsendorf functor}
The image of a $\mathcal{W}_{\mathfrak{o}}(R)$-display $\mathcal{P}$ as in
Definition \ref{Adorf1d} will be denoted by $\mathcal{P}_{\rm a}= \mathfrak{A}_{O/\mathfrak{o}, R}(\CP)$. 
The main theorem on the Ahsendorf functor is: 
\begin{theorem}\label{Adorf1t}
  Let $R$ be an $O$-algebra such that $p$ is nilpotent in $R$. 
  The Ahsendorf functor induces an equivalence of categories
  \index[NO]{ACA@$\mathfrak{A}_{O/\mathfrak{o}, R}$}
  \begin{displaymath}
  \mathfrak{A}_{O/\mathfrak{o}, R}:
\left(
  \begin{array}{l}
  \text{nilpotent   $\mathcal{W}_{\mathfrak{o}}(R)$-displays}\\
    \text{with strict $O$-action} 
    \end{array}
\right)
  \; \longrightarrow \; \Big(\text{nilpotent $ \mathcal{W}_{O}(R)$-displays}\Big) .  
  \end{displaymath}
 Furthermore, the Ahsendorf functor canonically associates  to a bilinear form
\begin{equation}\label{Adorf7e}
\beta: \mathcal{P}' \times \mathcal{P}'' \longrightarrow \mathcal{P} 
\end{equation}
of $\mathcal{W}_{\mathfrak{o}}(R)$-displays with a strict $O$-actions such
that $\beta$ is also $O$-bilinear, a bilinear form of $\mathcal{W}_{O}(R)$-displays 
\begin{displaymath}
\mathcal{P}'_{\rm a} \times \mathcal{P}''_{\rm a} \longrightarrow \mathcal{P}_{\rm a}.
\end{displaymath}
 \end{theorem}
 \begin{proof}
  The first statement is the main result of \cite{ACZ}. The second statement is shown  in Proposition \ref{Adorf4c}.
 \end{proof}
\begin{remark}\label{Ahsrem}
By Theorem \ref{Adorf1t} we obtain a functor
\begin{equation}\label{Adorf26e}
 \Big(\text{strict formal $p$-divisible $O$-modules over $R$}\Big) \longrightarrow  
 \Big( \text{$\mathcal{W}_O(R)$-displays}\Big), 
  \end{equation}
which is defined as follows. By \cite{L}, \cite{Zi3} there is a functor
from the first category to the category of $\mathcal{W}(R)$-displays with a
strict $O$-action. Composing this with $\mathfrak{A}_{O/\mathbb{Z}_p, R}$ we
obtain (\ref{Adorf26e}). In particular this gives a quasi-inverse functor
to the functor of Theorem \ref{maindisp}. 
\end{remark}

We will now define the Ahsendorf functor. 
We denote by $K^t \subset K$ the maximal subextension which is unramified
over $\mathfrak{k}$. Let $O^t $ be the ring of integers of $K^t$.
We
consider the Witt vectors $W_{O^t}(R)$ with respect to the prime element
$\varpi \in O^t$. The Frobenius resp. the Verschiebung acting on $W_{O^t}(R)$
will be denoted by $F'$ and $V'$. We will define
$\mathfrak{A}_{O/\mathfrak{o}, R}$ as the composite of two functors
\begin{equation}\label{Adorf8e}
  \begin{array}{l} 
    \mathfrak{A}_{O^t/\mathfrak{o}, R}:
    \left(
  \begin{array}{l}
    \mathcal{W}_{\mathfrak{o}}(R)-\text{displays}\\
    \text{with strict $O$-action} 
  \end{array}
  \right)
  \; \longrightarrow \;
\left(
  \begin{array}{l}
    \mathcal{W}_{O^t}(R)-\text{displays}\\
    \text{with strict $O$-action} 
  \end{array}
  \right),\\[4mm] 
  \mathfrak{A}_{O/O^t, R}:
  \left(
  \begin{array}{l}
    \mathcal{W}_{O^t}(R)-\text{displays}\\
    \text{with strict $O$-action} 
  \end{array}
  \right)
  \; \longrightarrow \; \Big( \mathcal{W}_{O}(R)-\text{displays} \Big).  
    \end{array}
  \end{equation}

We begin with the definition of $\mathfrak{A}_{O^t/\mathfrak{o}, R}$. 
\begin{lemma}\label{Adorf1l}
Let $S$ be an $O$-algebra which has no $\pi$-torsion. Let 
$\tau: S \longrightarrow S$ be a $O$-algebra homomorphism such that 
\begin{displaymath} 
\tau(s) \equiv s^q \; (\hspace{-3mm} \mod \pi). 
\end{displaymath}
Let $u_0, u_1, \ldots, u_n, \ldots \in S$. Then there exists  
$\xi \in W_{O}(S)$ such that $\mathbf{w}_{O,n}(\xi) = u_n$ for all $n$ iff
  \begin{displaymath}
\tau(u_{n-1}) \equiv u_n \; \mod \pi^n S, \quad \text{for} \quad n \geq 1. 
    \end{displaymath}
The element $\xi$ is uniquely determined. 
\end{lemma}
\begin{proof}
  The proof is up to obvious changes identical with the proof for the
  classical case $O = \mathbb{Z}_p$, cf. \cite[ IX, \S1, 2, Lemme 2]{Bour}.
  \end{proof}
We denote by $\sigma \in \Gal(K^t/\mathfrak{k})$ the Frobenius automorphism.
By Lemma \ref{Adorf1l}, there is a homomorphism 
$\lambda: O^t \longrightarrow W_{\mathfrak{o}}(O^t)$, defined by 
$\mathbf{w}_{\mathfrak{o},n}(\lambda(a)) = \sigma^n(a)$ for $a \in O^t$ and all $n$. 
We obtain a ring homomorphism  
\begin{equation}\label{maplambda} 
\varkappa: O^t \overset{\lambda}{\longrightarrow}
W_{\mathfrak{o}}(O^t) \longrightarrow W_{\mathfrak{o}}(R).   
  \end{equation}
We introduce the \emph{ Ahsendorf frame} \index{Ahsendorf frame} with respect to the
unramified extension $O^t/\mathfrak{o}$ for a $p$-adic $O^t$-algebra $R$,   
\begin{equation}\label{Aframe1e}
\mathcal{A}_{\mathfrak{o}}(R) = (W_{\mathfrak{o}}(R), I_{\mathfrak{o}}(R),
R, \mathfrak{f}^f, \mathfrak{f}^{f-1}\dot{\mathfrak{f}}).   
\end{equation}
This is an $O^t$-frame via $\varkappa$.

Let $\mathcal{P} = (P, Q, F, \dot{F})$ be a
$\mathcal{W}_{\mathfrak{o}}(R)$-display with a strict $O$-action. We set
\begin{displaymath}
P_m = \{x \in P \;|\; \iota(a)x = \varkappa(\sigma^{m}(a))x, \;
\text{for}\; a \in O^t\}, \quad m \in
\mathbb{Z}/f\mathbb{Z}. 
\end{displaymath} 
The $W_{\mathfrak{o}}(R)$-module $P$ decomposes as 
\begin{equation}\label{Aunverzweigt1e}
P = \oplus_{m\in \mathbb{Z}/f\mathbb{Z}} P_m.
\end{equation}
There is a similiar decomposition for $Q$. 
The maps $F$ and $\dot{F}$ of $\mathcal{P}$ are graded of degree one,
\begin{displaymath}
F: P_m \longrightarrow P_{m+1}, \quad  \dot{F}: Q_m \longrightarrow P_{m+1}.  
\end{displaymath}
If the action $\iota$ is strict, we have $Q_m = P_m$ for $m \neq 0$. Then we define 
the $\mathcal{A}_{\mathfrak{o}}(R)$-display $\mathcal{P}_{\rm ua}$: 
\begin{equation}\label{Adorf14e}
  P_{\rm u a} = P_0,\quad Q_{\rm u a} = Q_0,\quad F_{\rm u a} = \dot{F}^{f-1}F,\quad
  \dot{F}_{\rm u a} = \dot{F}^{f}.
\end{equation}
It is clear that $O$ acts strictly on $\mathcal{P}_{\rm u a}$.

It follows from (\ref{Drinfmap1e}) that
$\mu: W_{\mathfrak{o}}(R) \longrightarrow W_{O^t}(R)$ induces a morphism of
$O^t$-frames
\begin{equation}\label{Adorf20e} 
\mu: \mathcal{A}_{\mathfrak{o}}(R) \longrightarrow \mathcal{W}_{O^t}(R).
\end{equation}
By base change we obtain from $\mathcal{P}_{\rm u a}$ a
$\mathcal{W}_{O^t}(R)$-display  $\mathcal{P}_{\text{t}}=\mu_*(\CP_{\rm u a})$. The strict action of $O$
on $\mathcal{P}_{\rm u a}$ induces a strict action of $O$ on $\mathcal{P}_{\rm t}$
because the tangent space remains unchanged by this base change. 
\begin{definition}\label{DefAdunr}
  The Ahsendorf functor $\mathfrak{A}_{O^t/\mathfrak{o}, R}$ is the functor which associates to a
  $\mathcal{W}_{\mathfrak{o}}(R)$-display $\mathcal{P}$ with a strict $O$-action
  the $\mathcal{W}_{O}(R)$-display $\mathcal{P}_{\rm t}$ defined above.    
\end{definition}

The Ahsendorf functor is compatible with bilinear forms as follows. Let $\beta: \mathcal{P}' \times \mathcal{P}'' \longrightarrow \mathcal{P}$
as in (\ref{Adorf7e}). Because $\beta$ is $O^t$-bilinear, $\beta$ induces for each $m\in\BZ/f\BZ$ a pairing 
\begin{equation*}
  \beta: P'_m \times P''_m \longrightarrow P_m,
  \end{equation*}
and $P'_i$ and $P''_j$ are orthogonal for $i \neq j$. For $y' \in Q'_0$
and $y'' \in Q''_0$ we find
\begin{displaymath}
\beta((\dot{F}')^f y', (\dot{F}'')^f y'') = \dot{F}^f \beta(y', y''). 
  \end{displaymath}
Therefore the restriction of $\beta$,
\begin{displaymath}
\beta_{\rm u a}: P'_0 \times P''_0 \longrightarrow P_0 
\end{displaymath}
induces a bilinear form of $\mathcal{A}_{\mathfrak{o}}(R)$-displays
\begin{displaymath}
\mathcal{P}'_{\rm u a} \times \mathcal{P}''_{\rm u a} \longrightarrow \mathcal{P}_{\rm u a}. 
\end{displaymath}
Applying Proposition \ref{basebil}, we obtain 
a bilinear form in the category of $\mathcal{W}_{O^t}(R)$-displays,
\begin{equation}\label{bilt}
\beta_{\rm t}: \mathcal{P}'_{\rm t} \times \mathcal{P}''_{\rm t} \longrightarrow \mathcal{P}_{\rm t}. 
  \end{equation}

Now we define  $\mathfrak{A}_{O/O^t, R}$. First we introduce the
Lubin-Tate frame. We choose a
finite normal extension $L$ of $K^t$ which contains $K$. We set
$\Phi = \Hom_{\text{$K^t$-Alg}}(K, L)$. Let $\varphi_0: K \longrightarrow L$ be 
the identity embedding. 

Let $\mathbf{E}_K \in O^t[T]$ be the Eisenstein polynomial of
$\pi \in O$ over $O^t$. In $O_L[T]$ it decomposes as
\index[NO]{EAE@$\mathbf{E}_K$, $\tilde{\mathbf{E}}_K$} 
\begin{displaymath}
\mathbf{E}_K(T) = \prod_{\varphi \in \Phi} (T - \varphi(\pi)). 
  \end{displaymath}
We set \index[NO]{EAF@$\mathbf{E}_{K,0}$, $\tilde{\mathbf{E}}_{K,0}$}
\begin{displaymath}
  \mathbf{E}_{K,0}(T) = \prod_{\varphi \in \Phi, \varphi \neq \varphi_0}
  (T - \varphi(\pi)) \in O_L[T] 
\end{displaymath}
One sees easily that $\mathbf{E}_{K,0} \in O[T]$. 
We lift these polynomials via $\mathbf{w}_{O^t,0}$ to the ring of Witt vectors,
\begin{equation}\label{Adorf21e}
  \begin{array}{ll}
    \tilde{\mathbf{E}}_K(T) = \prod_{\varphi \in \Phi} (T - [\varphi(\pi)]) \in
    W_{O^t}(O^t)[T], \\[3mm]  
    \tilde{\mathbf{E}}_{K,0}(T) = \prod_{\varphi \in \Phi, \varphi \neq \varphi_0}
    (T - [\varphi(\pi)]) \in  W_{O^t}(O)[T].  
    \end{array}
\end{equation}
The Frobenius $F'$ and the Verschiebung $V'$ act
via the second factor on $O \otimes_{O^t} W_{O^t}(R)$. We set
\begin{displaymath}
  \dot{F}' = (V')^{-1}: O \otimes_{O^t} I_{O^t}(R) \longrightarrow
  O \otimes_{O^t} W_{O^t}(R). 
\end{displaymath}

\begin{proposition}\label{NK1p}
  The element
  \begin{displaymath}
    ~^{\dot{F}'}(\tilde{\mathbf{E}}_K(\pi \otimes 1)) \in 
    O \otimes_{O^t} W_{O^t}(O^t) 
  \end{displaymath}
  is a unit of the form
  \begin{equation}\label{NK2e}
    (\frac{\pi^e}{\varpi} \otimes 1)\delta, \quad \delta \in
    O \otimes_{O^t} W_{O^t}(O^t) 
    \end{equation}
  such that $\delta - (1 \otimes 1)$ lies in the kernel of
  \begin{displaymath}
    O \otimes_{O^t} W_{O^t}(O^t) \longrightarrow 
    O \otimes_{O^t} W_{O^t}(\kappa) ,
    \end{displaymath}
  and hence in the Jacobson radical of $O \otimes_{O^t} W_{O^t}(O^t)$. 
\end{proposition}
\begin{proof} The element is defined because
  $\id \otimes \mathbf{w}_{O^t,0}: O \otimes_{O^t} W_{O^t}(O^t) \longrightarrow O$
  maps $\tilde{\mathbf{E}}_K(\pi \otimes 1)$ to $0$. 
  In the following computation we pass to
  $O \otimes_{O^t} W_{O^t}(O_L)$. We find from the  definitions:
    \begin{equation}\label{NK3e}
    \begin{aligned} 
    ~^{\dot{F}'}(\tilde{\mathbf{E}}_K(\pi \otimes 1)) &= 
      ~^{\dot{F}'}(\prod_{\varphi \in \Phi}(\pi \otimes 1 - 1 \otimes [\varphi(\pi)])
      \\
       &=    \frac{1}{\varpi} \prod_{\varphi \in \Phi}
    (\pi \otimes 1 - 1 \otimes [\varphi(\pi)^q]) = 
    \frac{1}{\varpi} \sum_{i=0}^{e} \pi^{e-i} \otimes ~^{F'}s_i.
    \end{aligned}
  \end{equation}

  Here we denote by $s_i$ the elementary symmetric polynomial of degree $i$ evaluated at
  the $e$ arguments $[\varphi(\pi)]$. By definition $s_0 =1$. We claim that
  for $i >0$ the elements $~^{F'}s_i \in W_{O^t}(O^t)$ are
  divisible by $\varpi$. Clearly $\mathbf{w}_{O^t,0}(s_i)$ is divisible by $\pi$.
  On the other hand, $\mathbf{w}_{O^t,0}(s_i) \in O^{t}$ and therefore
  is divisible by $\varpi$. We find expressions in $W_{O^t}(O^t)$,  
  \begin{displaymath}
    s_i = [\varpi c_i] + ~^{V'}\xi_i, \quad c_i \in O^t, \; \xi_i \in
    W_{O^t}(O^t). 
    \end{displaymath}
  Therefore $~^{F'}s_i = [\varpi^q][c_i^q] + \varpi \xi^i$ is divisible
  by $\varpi$. Indeed, using Lemma \ref{Adorf1l},  one  shows 
  as in the proof of \cite[Lem. 28]{Zi} that  $\varpi$ divides  $[\varpi^q]$.

  Now we may write the last term of (\ref{NK3e}) as 
  \begin{displaymath}
    \frac{\pi^e}{\varpi} \otimes 1  
    + \sum_{i=1}^{e} \pi^{e-i} \otimes \frac{~^{F'}s_i}{\varpi}.
  \end{displaymath}
  Finally, $\frac{~^{F'}\!s_i}{\varpi}$ lies for $i > 0$ in the kernel of
  $W_{O^t}(O^t) \longrightarrow W_{O^t}(\kappa)$. Indeed, the
  elements $[\varphi(\pi)] \in W_{O^t}(O_L)$ are mapped to zero in
  $W_{O^t}(\kappa_L)$ and therefore a fortiori the symmetric functions 
  $s_i$. We conclude that $s_i$ and then $~^{F'}s_i$ become zero in
  $W_{O^t}(\kappa)$ for $i > 0$. Because $\varpi$ is not a zero divisor in
  $W_{O^t}(\kappa)$ the elements $\frac{~^{F'}s_i}{\varpi}$ are then also in the
  kernel. 
\end{proof}

We write in the ring $W_{O}(O)$, 
\begin{equation}\label{NK16e}
\pi - [\pi] = ~^{V}\varepsilon.
\end{equation}
One checks that $\varepsilon \in W_{O}(O)$ is a unit. If we apply $F$
to the last equation, we obtain
\begin{equation}\label{NK15e}
\pi - [\pi^q] = \pi \varepsilon. 
\end{equation}
In particular $[\pi^q]$ is divisible by $\pi$. 

\begin{lemma}\label{NK5l}
  The image of the element
  $\big(~^{\dot{F}'}\tilde{\mathbf{E}}_K(\pi \otimes 1)\big)^{-1} \cdot
  ~^{F'}\tilde{\mathbf{E}}_{K,0}(\pi \otimes 1)$ 
  under the Drinfeld homomorphism
  \begin{displaymath}
\mu: O \otimes_{O^t} W_{O^t}(O) \longrightarrow W_{O}(O) 
  \end{displaymath}
  equals $\varepsilon^{-1}(\varpi/\pi)$. 
\end{lemma}
\begin{proof}
  It is enough to show the same assertion for
  $O \otimes_{O^t} W_{O^t}(O_L) \longrightarrow W_{O}(O_L)$. The image of
  $~^{\dot{F}'}\tilde{\mathbf{E}}_K(\pi \otimes 1)$ by the last map is
  $\varpi^{-1} \prod_{\varphi}(\pi - [\varphi(\pi)^q])$. Here we used that $\varpi$
  is not a zero divisor in the participating rings. Our assertion is
  equivalent to the equation
  \begin{displaymath}
    \varpi^{-1} \prod_{\varphi}(\pi - [\varphi(\pi)^q]) \varepsilon^{-1}
    \frac{\varpi}{\pi} = \prod_{\varphi \neq \varphi_0}(\pi - [\varphi(\pi)^q]). 
  \end{displaymath}
  But this is a consequence of (\ref{NK15e}). 
  \end{proof}

The free $W_{O^t}(R)$-module $O \otimes_{O^t} W_{O^t}(R)$ has the basis
\begin{equation}\label{NK4e}
1 \otimes 1, \; \pi \otimes 1 - 1 \otimes [\pi], \ldots, 
\pi^{m} \otimes 1 - 1 \otimes [\pi]^{m}, \ldots, 
\pi^{e-1} \otimes 1 - 1 \otimes [\pi]^{e-1}. 
\end{equation}
To ease the notation, here $[\pi]$ denotes the Teichm\"uller representative of
the image of $\pi$ by the morphism $O \longrightarrow R$.  
Let
\begin{displaymath}
\CJ=\Ker\big(O \otimes_{O^t} W_{O^t}(R) \longrightarrow R\big)  ,
  \end{displaymath}
  where the map is induced by $\mathbf{w}_0: W_{O^t}(R) \longrightarrow R$. The ideal $\mathcal{J}$ is contained in the radical of
$O \otimes_{O^t} W_{O^t}(R)$.  As a
$W_{O^t}(R)$-module, $\mathcal{J}$ is the direct sum of
$O \otimes_{O} I_{O}(R)$ and the direct
summand generated by the last $e-1$ elements of (\ref{NK4e}). 
In particular we obtain 
\begin{equation}\label{NK5e}
  \mathcal{J} = O \otimes_{O^t} I_{O^t}(R) + (\pi \otimes 1 - 1 \otimes [\pi])
  (O \otimes_{O^t} W_{O^t}(R)). 
\end{equation}

We define maps
$\sigma_{\rm{lt}}: O \otimes_{O^t} W_{O^t}(R) \longrightarrow
O \otimes_{O^t} W_{O^t}(R)$,
$\dot{\sigma}_{\rm{lt}}: \mathcal{J} \longrightarrow O \otimes_{O^t} W_{O^t}(R)$
by
\begin{displaymath}
  ~^{\sigma_{\rm{lt}}} \xi = ~^{F'} \xi, \quad \quad 
  ~^{\dot{\sigma}_{\rm{lt}}} \eta =
(~^{\dot{F}'}\tilde{\mathbf{E}}_K)^{-1} ~^{\dot{F}'}(\tilde{\mathbf{E}}_{K,0} \eta),
  \quad  \xi \in O \otimes_{O^t} W_{O^t}(R), \;
  \eta \in \mathcal{J}. 
\end{displaymath}
The map
$\dot{\sigma}_{\rm{lt}}: \mathcal{J} \longrightarrow O \otimes_{O^t} W_{O^t}(R)$
is $\sigma_{\rm{lt}}$-linear. 
Then we obtain
\begin{equation}\label{LTF1e}
  ~^{\dot{\sigma}_{\rm{lt}}}(\pi \otimes 1 - 1 \otimes [\pi]) =
  (~^{\dot{F}'}\tilde{\mathbf{E}}_K)^{-1}
  ~^{\dot{F}'}(\tilde{\mathbf{E}}_K) = 1. 
\end{equation}
\begin{definition}\label{LTF3p}{\rm (comp. \cite[Def. 2.7]{M})}
 The \emph{Lubin-Tate frame}  \index{Lubin-Tate frame} for $O$ is the $O$-frame 
$$\mathcal{F}_{\rm{lt}}(R) := (O \otimes_{O^t} W_{O^t}(R), \mathcal{J}, R,
  \sigma_{\rm{lt}}, \dot{\sigma}_{\rm{lt}}).$$
 \end{definition}
This is indeed an $O$-frame: the only thing we need to check is
  \begin{displaymath}
    ~^{\sigma_{\rm{lt}}} \xi \equiv \xi^q \;
    (\!\!\!\!\mod (\pi \otimes 1) O \otimes_{O^t} W_{O^t}(R)) ,  
  \end{displaymath}
  and this follows because $a \equiv a^q \!\!\mod \pi$ for $a \in O$ and
  $^{F'}\!\eta \equiv \eta^q  \!\!\mod \varpi$ for $\eta \in W_{O^t}(R)$. 

We remark that by (\ref{LTF1e}) 
\begin{displaymath}
  (\pi \otimes 1 - 1 \otimes [\pi^q]) ~^{\dot{\sigma}_{\rm{lt}}} \eta =
  ~^{\sigma_{\rm{lt}}} \eta, \quad  \eta \in \mathcal{J}. 
  \end{displaymath}

Now we start with a $\mathcal{W}_{O^t}(R)$-display  $\mathcal{P}$ with a strict 
$O$-action.  The last condition can be reformulated as 
\begin{equation}\label{NKstrict1e}
\mathcal{J}P \subset Q.
  \end{equation}
  We refer to \cite[Prop. 2.26]{ACZ}  for the  proof of the following lemma.  
\begin{lemma}
Let $P$ be a $W_{O^t}(R)$-module with an action of $O$, i.e., a
homomorphism of $O$-algebras
\begin{displaymath}
O \longrightarrow \End_{W_{O^t}(R)} P.
\end{displaymath}
Assume that locally on $\Spec R$ the $W_{O^t}(R)$-module $P$ is free. 
Then $P$ is locally on $\Spec R$ a 
finitely generated free $O \otimes_{O^t} W_{O^t}(R)$-module. \qed
\end{lemma}

\begin{lemma}\label{NK2l}
Let $\mathcal{P}$ be a $\mathcal{W}_{O^t}(R)$-display with a strict 
$O$-action. Let $x \in P$. By (\ref{NKstrict1e})
$(\pi \otimes 1 - 1 \otimes [\pi])x \in Q$. The following equation holds,
\begin{displaymath}
  Fx = \big(~^{\dot{F}'}\tilde{\mathbf{E}}_K(\pi \otimes 1)\big)^{-1} \cdot
  ~^{F'}\tilde{\mathbf{E}}_{K,0}(\pi \otimes 1) \cdot
  \dot{F}\big((\pi \otimes 1 - 1 \otimes [\pi])x\big). 
  \end{displaymath} 
\end{lemma}
\begin{proof}
  From the definition of the polynomials $\tilde{\mathbf{E}}_K$ and
  $\tilde{\mathbf{E}}_{K,0}$, we find since $\varphi_0(\pi) = \pi$, 
  \begin{displaymath}
    \tilde{\mathbf{E}}_K(\pi \otimes 1)  = 
    \tilde{\mathbf{E}}_{K,0}(\pi \otimes 1)
    \cdot (\pi \otimes 1 - 1 \otimes [\pi]). 
  \end{displaymath}
  Therefore 
   \begin{equation*}
    \begin{aligned} 
\dot{F}(\tilde{\mathbf{E}}_K(\pi \otimes 1) x) & =
\dot{F}(\tilde{\mathbf{E}}_{K,0}(\pi \otimes 1))
\cdot (\pi \otimes 1 - 1 \otimes [\pi]) x)\\ &  =
~^{F'}(\tilde{\mathbf{E}}_{K,0}(\pi \otimes 1))
\dot{F}((\pi \otimes 1 - 1 \otimes [\pi])x).
      \end{aligned} 
  \end{equation*}
  Because $\tilde{\mathbf{E}}_K(\pi \otimes 1) \in O \otimes_{O^t} I_{O^t}(R)$,
  we obtain
  \begin{displaymath}
    \dot{F}(\tilde{\mathbf{E}}_K(\pi \otimes 1) x) =
    ~^{\dot{F}'}(\tilde{\mathbf{E}}_K(\pi \otimes 1)) Fx.
  \end{displaymath}
  We conclude by Lemma \ref{NK1p}. 
\end{proof}

We now associate to the $\mathcal{W}_{O^t}(R)$-display  $\mathcal{P} = (P, Q, F, \dot{F})$ with
a strict $O$-action a
$\mathcal{F}_{\rm{lt}}$-display\index[NO]{PBC@$\mathcal{P}_{\rm{lt}}$}
$\mathcal{P}_{\rm{lt}} = (P_{\rm{lt}}, Q_{\rm{lt}}, F_{\rm{lt}},
\dot{F}_{\rm{lt}})$. We set $P_{\rm{lt}} = P$, $Q_{\rm{lt}} = Q$,
$\dot{F}_{\rm{lt}} = \dot{F}$, and
\begin{equation}\label{LTF7e} 
  F_{\rm{lt}}(x) = \dot{F}((\pi \otimes 1 - 1 \otimes [\pi])x), \quad
   x \in P. 
\end{equation}
\begin{proposition}\label{FLT4p}
  $\mathcal{P}_{\rm{lt}}$ is an $\mathcal{F}_{\rm{lt}}(R)$-display. 
\end{proposition}
\begin{proof}
  The only thing we have not checked is the equation
  \begin{equation}\label{LTF3e}
    \dot{F}_{\rm{lt}}(\eta x) = ~^{\dot{\sigma}_{\rm{lt}}}\eta F_{\rm{lt}} x,
    \quad  \eta \in \mathcal{J}.
  \end{equation}
  We begin with the case $\eta = ~^{V'}\xi$. We apply Lemma \ref{NK2l}
  \begin{equation}\label{LTF4e}
    \dot{F}_{\rm{lt}}(\eta x) = \dot{F}(~^{V'}\xi x) = \xi F(x) =
    \xi (~^{\dot{F}}\tilde{\mathbf{E}}_K)^{-1} \cdot
  ~^{F}\tilde{\mathbf{E}}_{K,0} \cdot
  \dot{F}((\pi \otimes 1 - 1 \otimes [\pi])x). 
  \end{equation}
  By definition
  \begin{displaymath}
  ~^{\dot{\sigma}_{\rm{lt}}} \eta =
    (~^{\dot{F}'}\tilde{\mathbf{E}}_K)^{-1}
    ~^{\dot{F}'}(\tilde{\mathbf{E}}_{K,0} ~^{V}\xi) = 
(~^{\dot{F}'}\tilde{\mathbf{E}}_K)^{-1} ~^{F'}(\tilde{\mathbf{E}}_{K,0}) \xi.  
  \end{displaymath}
  Using the definition (\ref{LTF7e}), we can write the right hand side of
  (\ref{LTF4e}) as $~^{\dot{\sigma}_{\rm{lt}}} \eta F_{\rm{lt}}(x)$, hence we are done in this case. 

  Next we consider the case where
  $\eta = (\pi \otimes 1 - 1 \otimes [\pi])\xi$. Then we find
  \begin{displaymath}
    \dot{F}_{\rm{lt}}(\eta x) =
    \dot{F}((\pi \otimes 1 - 1 \otimes [\pi])\xi x) =
    ~^{F}\xi F_{\rm{lt}}(x). 
    \end{displaymath}
  But by (\ref{LTF1e}) we have 
  $~^{\dot{\sigma}_{\rm{lt}}}((\pi \otimes 1 - 1 \otimes [\pi])\xi) = ~^{F'}\xi$. 
  Therefore in this case (\ref{LTF3e}) is true as well.
 \end{proof}

We use the same symbol $\varepsilon \in W_O(R)^\times$ for the image of the element
$\varepsilon \in W_O(O)^\times$ defined by (\ref{NK16e}).  We define the frame
\index[NO]{WBC@$\mathcal{W}^{\varepsilon}_{O}(R)$} 
\begin{equation}\label{Adorf10e} 
\mathcal{W}^{\varepsilon}_{O}(R) =
(W_{O}(R), I_{O}(R), F, \varepsilon^{-1} \dot{F}). 
  \end{equation} 
We note that the categories of displays over $\mathcal{W}_{O}(R)$ and
$\mathcal{W}^{\varepsilon}_{O}(R)$ are canonically isomorphic. Indeed,
if $\mathcal{P} = (P, Q, F, \dot{F})$ is a $\mathcal{W}_{O}(R)$-display,
then $\mathcal{P} = (P, Q, \varepsilon F, \dot{F})$ is a
$\mathcal{W}^{\varepsilon}_{O}(R)$-display.

Recall that we denote the Frobenius and the
Verschiebung acting on $W_{O^t}(R)$ by $F'$ and $V'$. 
We consider the Drinfeld homomorphism $\mu: W_{O^t}(R) \longrightarrow W_{O}(R)$,
cf. (\ref{Drinfmap2e}).   This is a functorial ring homomorphism such that
$\mathbf{w}'_n(\mu(\xi)) = \mathbf{w}_n(\xi)$ which has the following
properties
\begin{equation}\label{LTF6e}
\mu(~^{F'}\xi) = ~^{F}\mu(\xi), \quad  \mu(~^{V'}\xi) = \frac{\varpi}{\pi}
~^{V}\mu(\xi), \quad \mu([a]) = [a], \; \text{for} \; a \in R.
\end{equation}
The Drinfeld homomorphism extends to a ring homomorphism
\begin{equation}\label{LTF5e}
\mu: O \otimes_{O^t} W_{O^t}(R) \longrightarrow W_{O}(R) 
  \end{equation}
which we denote by the same letter.
\begin{proposition}\label{LTF2p}
  The Drinfeld homomorphism induces a morphism of $O$-frames
  \begin{displaymath}
  \mu:  \mathcal{F}_{\rm{lt}}(R) \longrightarrow  \mathcal{W}^{\varepsilon}_{O}(R).
  \end{displaymath} 
  \end{proposition} 
\begin{proof}
  We have to check that the image of 
  $\mathcal{J} \subset O \otimes_{O^t} W_{O^t}(R)$ by $\mu$ is contained in
  $I_{O}(R)$. This is immediate because
  $\mu(\pi \otimes 1 - 1 \otimes [\pi]) = \pi - [\pi] = ~^{V}\varepsilon$.
  It remains to prove the equations for $\xi \in O \otimes_{O^t} W_{O^t}(R)$
  and $\eta \in \mathcal{J}$,
  \begin{equation}
    \mu(~^{\sigma_{\rm{lt}}}\xi) = ~^{F}\mu(\xi), \quad
    \mu(~^{\dot{\sigma}_{\rm{lt}}}\eta) = \varepsilon^{-1} ~^{\dot{F}}\mu(\eta). 
  \end{equation}
  The first equation follows from (\ref{LTF6e}). To prove the second equation,
  it is enough to consider the following two cases separately:
  $\eta = ~^{V'}\xi$ and $\eta = (\pi \otimes 1 - 1 \otimes [\pi])\xi$.
  In the first case we have
  \begin{displaymath}
    ~^{\dot{\sigma}_{\rm{lt}}}\eta =
    (~^{\dot{F}'}\tilde{\mathbf{E}}_K)^{-1}
    ~^{\dot{F}'}(\tilde{\mathbf{E}}_{K,0} ~^{V'}\xi) =
    (~^{\dot{F}'}\tilde{\mathbf{E}}_K)^{-1}
    (~^{F'}\tilde{\mathbf{E}}_{K,0}) \xi.     
  \end{displaymath}
  Applying Lemma \ref{NK5l}, we obtain
  \begin{displaymath}
\mu(~^{\dot{\sigma}_{\rm{lt}}}\eta) = \varepsilon^{-1}(\varpi/\pi)\mu(\xi).
  \end{displaymath}
  On the other hand, we have by (\ref{LTF6e})
  \begin{displaymath}
    \varepsilon^{-1} ~^{\dot{F}}\mu(^{V'}\xi) = \varepsilon^{-1}
    ~^{\dot{F}}((\varpi/\pi)~^{V}\mu(\xi)) = \varepsilon^{-1}
    ((\varpi/\pi)\mu(\xi) , 
  \end{displaymath}
  as desired.

  Now we consider the case $\eta = (\pi \otimes 1 - 1 \otimes [\pi])\xi$.
  We have
  \begin{displaymath}
    \mu((\pi \otimes 1 - 1 \otimes [\pi])\xi) = (\pi - [\pi])\mu(\xi)
    = ~^{V}\varepsilon \mu(\xi) = ~^{V}(\varepsilon ~^{F}\mu(\xi)). 
  \end{displaymath}
  We obtain
  \begin{displaymath}
\varepsilon^{-1} ~^{\dot{F}}\mu(\eta) = ~^{F}\mu(\xi). 
  \end{displaymath}
  On the other hand, we find by (\ref{LTF1e}) and (\ref{LTF6e}) 
  \begin{displaymath}
\mu(~^{\dot{\sigma}_{\rm{lt}}} \eta) = \mu(~^{F'}\xi) = ~^{F}\mu(\xi) , 
  \end{displaymath}
  as desired.
\end{proof}
Starting now with a $\mathcal{W}_{O^t}(R)$-display
 $\mathcal{P} = (P, Q, F, \dot{F})$ with a strict $O$-action, we  
 have the associated  $\mathcal{F}_{\rm{lt}}(R)$-display
$\mathcal{P}_{\rm{lt}} = (P_{\rm{lt}}, Q_{\rm{lt}}, F_{\rm{lt}},
\dot{F}_{\rm{lt}})$, cf. Proposition \ref{FLT4p}. After taking the base change by the morphism of frames of Proposition
\ref{LTF2p}, we obtain a $\mathcal{W}^{\varepsilon}_{O}(R)$-display
$\mathcal{P}^\varepsilon_{\rm a} = (P^\varepsilon_{\rm a}, Q^\varepsilon_{\rm a},  F^\varepsilon_{\rm a}, \dot{F}^\varepsilon_{\rm a})$. Then
$\mathcal{P}_{\rm a} = (P^\varepsilon_{\rm a}, Q^\varepsilon_{\rm a}, \varepsilon^{-1}F^\varepsilon_{\rm a}, \dot{F}^\varepsilon_{\rm a})$ is an $\mathcal{W}_{O}(R)$-display.
\begin{definition}\label{Adorf2d}
  The Ahsendorf functor $\mathfrak{A}_{O/O^t, R}$ is the functor which associates to the 
  $\mathcal{W}_{O^t}(R)$-display $\mathcal{P}$ with a strict $O$-action
  the $\mathcal{W}_{O}(R)$-display $\mathcal{P}_{\rm a}$ defined above.    
\end{definition}
From the construction we obtain that 
\begin{equation}\label{Adorf9e} 
P_{{\rm{a}}} = W_{O}(R) \otimes_{O \otimes_{O^t} W_{O^t}(R)} P,
\end{equation}
and that $Q_{{\rm{a}}}$ is the kernel of the natural map
\begin{displaymath} 
W_{O}(R) 
\otimes_{O \otimes_{O^t} W_{O^t}(R)} P \longrightarrow P/Q.
\end{displaymath}
We note that the canonical map $P \longrightarrow P_{{\rm{a}}}$ induces a
map $Q \longrightarrow Q_{{\rm{a}}}$. 

\begin{proposition}
  Let $\mathcal{P}$ be a $\mathcal{W}_{O^t}(R)$-display with a strict action
  of $O$. Let $\mathcal{P}_{\rm a}=\mathfrak{A}_{O/O^t, R}(\CP)$ be its image by the Ahsendorf functor. The following diagram is commutative
  \begin{displaymath}
   \xymatrix{
          Q \ar[r]^{\dot{F}} \ar[d] & P \ar[d]\\
          Q_{{\rm{a}}} \ar[r]_{\dot{F}_{\rm a}} & P_{{\rm{a}}}\\
        }
    \end{displaymath}
\end{proposition}
\begin{proof}
  This follows from the definition of $\mathcal{P}_{\rm{lt}}$ before
  Proposition \ref{FLT4p} and the definition of base change (via the morphism
  of Proposition \ref{LTF2p}). 
\end{proof}
We note that this diagram determines the map $\dot{F}_{\rm a}$ uniquely. Indeed, consider the  
 following equation in $P_{\rm a}$ under the identification (\ref{Adorf9e}),
\begin{displaymath}
^{V}\varepsilon \otimes x = 1 \otimes (\pi \otimes 1 - 1 \otimes [\pi])x .
  \end{displaymath}
Applying $\dot{F}_{\rm a}$, we obtain from the diagram that
\begin{displaymath}
  \varepsilon F_{\rm a}(1 \otimes x) = 1 \otimes
  \dot{F}((\pi \otimes 1 - 1 \otimes [\pi])x).
\end{displaymath}
This shows that $F_{\rm a}$ is uniquely determined. Because the image of $Q$ and
$I_{O}(R)P_{\rm a}$ generate $Q_{\rm a}$ as a $W_{O}(R)$-module, the map $\dot{F}_{\rm a}$ is
then also uniquely determined. 

We return to the notation that $\mathcal{P}$ is an $\mathfrak{o}$-display
with a strict $O$-action. Applying the functors $\mathfrak{A}_{O^t/\mathfrak{o}, R}$ and $\mathfrak{A}_{O/O^t, R}$, we
obtain first $\mathcal{P}_{\rm t}$ and then $\mathcal{P}_{\rm a}$. 
We find by our definitions that, with the notation of (\ref{Aunverzweigt1e}),
\begin{equation}\label{Adorf12e}
  P_{\rm a} = W_{O}(R) \otimes_{\big(O \otimes_{O^t, \varkappa} W_{\mathfrak{o}}(R)\big)} P_0 =
  W_{O}(R) \otimes_{\big(O \otimes_{\mathfrak{o}} W_{\mathfrak{o}}(R)\big)} P. 
  \end{equation}
 We note that
$P_0 = (O \otimes_{O^t, \varkappa} W_{\mathfrak{o}}(R))
\otimes_{O \otimes_{\mathfrak{o}} W_{\mathfrak{o}}(R)} P$.

We already noted that the Ahsendorf functor $\mathfrak{A}_{O^t/\mathfrak{o}, R}$  is
compatible with bilinear forms. Similar remarks are also valid for the Ahsendorf functor $\mathfrak{A}_{O/O^t, R}$: first one checks that the functor
$\mathcal{P} \mapsto \mathcal{P}_{\rm{lt}}$ is compatible with bilinear
forms of displays, and then applies Proposition \ref{basebil} for the compatibility of base change with bilinear forms.  Taken together with \eqref{bilt}, we obtain the following property of the Ahsendorf functor $\mathfrak{A}_{O/\mathfrak{o}, R}$. 
\begin{proposition}\label{Adorf4c} 
Consider a bilinear form of $\mathcal{W}_{\mathfrak{o}}(R)$-displays,
\begin{displaymath}
\beta: \mathcal{P}' \times \mathcal{P}'' \longrightarrow \mathcal{P} ,
  \end{displaymath}
which is also $O$-bilinear.  Then the bilinear form $\beta: P' \times P'' \longrightarrow P$ 
induces by (\ref{Adorf12e}) 
a $W_O(R)$-bilinear form $\beta_{\rm a}: P'_{\rm a} \times P''_{\rm a} \longrightarrow P_{\rm a}$.
The bilinear form $\beta_{\rm a}$ is a bilinear form of $\mathcal{W}_O(R)$-displays,
\begin{displaymath}
\beta_{\rm a}: \mathcal{P}'_{\rm a} \times \mathcal{P}''_{\rm a} \longrightarrow \mathcal{P}_{\rm a}.
  \end{displaymath}\qed
\end{proposition}

\begin{remark}\label{Adorf1r}
In the case where $R = k$ is a perfect field, the description of the Ahsendorf
functor is very simple. We consider the functor $\mathfrak{A}_{O/\mathbb{Z}_p, k}$
which is relevant for us. As a prime element of $\mathbb{Z}_p$ we choose $p$.
The element $\varepsilon \in W_O(k)$ is $1$. As above, we denote by
$\mathfrak{f}$, resp. $\mathfrak{v}$, the Frobenius, resp. the Verschiebung,
of the ring of Witt vectors $W(k)$. In this case the morphism (\ref{Adorf20e}) of frames
$\mu: \mathfrak{A}_{O/\mathbb{Z}_p, k} \longrightarrow \mathcal{W}_{O^t}(k)$, 
 and the morphism of frames
$\mu: \mathcal{F}_{\rm{lt}}(k) \longrightarrow \mathcal{W}_{O}(k)$ 
of Proposition \ref{LTF2p} are isomorphisms. Therefore we identify
$\mathcal{W}_O(k)$ with the frame
\begin{equation}\label{Adorf22e}
  (O \otimes_{O^t} W(k), \pi O \otimes_{O^t} W(k), k, \mathfrak{f}^f,
  \mathfrak{f}^f \pi^{-1}). 
\end{equation}
This is a perfect frame with $u = \theta = \pi$, cf. Definition \ref{Adorf1ex}.

Let $(P,F,V)$ be a $\mathcal{W}(k)$-Dieudonn\'e module with a strict $O$-action.
We have the decomposition $P = \oplus_{m} P_m$ cf. (\ref{Aunverzweigt1e}).
The summand $P_0$ is an $O \otimes_{O^t} W(k)$-module. Since the action of
$O$ is strict, we find
\begin{displaymath}
\pi P_0 \subset Q_0 = V^fP_0. 
\end{displaymath}
Therefore we can define
\begin{equation}\label{Adorf23e}
V_{\rm a} = V^f, \; F_{\rm a} = V^{-f}\pi: P_0 \longrightarrow P_0.  
\end{equation}
Then $(P_0, F_{\rm a}, V_{\rm a})$ is a Dieudonn\'e module for the frame (\ref{Adorf22e}).
It is the image of $(P, F, V)$ by the Ahsendorf functor
$\mathfrak{A}_{O/\mathbb{Z}_p, k}$. 
  \end{remark}

\index{heights and the Ahsendorf functor}
\begin{proposition}\label{Adorfslopes1p}
  Let $R \in \Nilp_O$. We assume that $\Spec R$ is connected. Let
  $\mathcal{P}$ be a $\mathcal{W}(R)$-display with a strict $O$-action, and
  let $\mathcal{P}_{\rm a}$ be the image by the Ahsendorf functor
  $\mathfrak{A}_{O/\mathbb{Z}_p, k}$. Then 
  \begin{displaymath}
\height \mathcal{P} = [O:\mathbb{Z}_p] \height_O \mathcal{P}_{\rm a} .
    \end{displaymath}
  The right hand side denotes the height of the $\mathcal{W}_O(R)$-display
  $\mathcal{P}_{\rm a}$ in the sense of Definition \ref{Rah2d}. 

Let $\alpha: \mathcal{P}_1 \longrightarrow \mathcal{P}_2$ be an isogeny of
$\mathcal{W}(R)$-displays with strict $O$-action, and let
$\alpha_{\rm a}: \mathcal{P}_{1, a} \longrightarrow \mathcal{P}_{2,a}$ be the image
by the Ahsendorf functor. Then 
\begin{displaymath}
\height \alpha = [O^t:\mathbb{Z}_p] \height_O \alpha_{\rm a} .
\end{displaymath} 

Let $R = k$ be a perfect field. Let $\lambda_1 < \ldots < \lambda_m$ be the
slopes of $\mathcal{P}$. Then the slopes of $\mathcal{P}_{\rm a}$ are
$[O:\mathbb{Z}_p]\lambda_1, \ldots, [O:\mathbb{Z}_p]\lambda_m$. The display 
$\mathcal{P}$ with its strict $O$-action is isogenous to a direct sum of
displays with a strict $O$-action $\oplus_{i=1}^m \mathcal{P}(\lambda_i)$
such that $\mathcal{P}(\lambda_i)$ is isoclinic of slope $\lambda_i$. 
\end{proposition}
\begin{proof}
  It suffices to consider the case where $R = k$ is a perfect field.
  Then it is a consequence of the description of the Ahsendorf functor given
  above, cf.  (\ref{Adorf23e}). 
  \end{proof}

In the end of this subsection, we relate explicitly the deformation
theory of a display with
a strict $O$-action and its image by the Ahsendorf functor. 
Let $S \longrightarrow R$ be an epimorphism of $O$-algebras which are $p$-adic.
We assume that the kernel $\mathfrak{a}$ of this epimorphism is endowed
with divided powers $\gamma$ relative to $\mathfrak{o}$.
Then $\gamma$ induces also divided powers $\gamma_{\rm t}$ on $\mathfrak{a}$ relative
to $O^t$. Indeed, let $q_{\mathfrak{o}}$ be the number of elements in the residue
class field of $\mathfrak{o}$. Then we set 
\begin{equation}\label{Adorf19e} 
  \gamma_{\rm t}(a) = \gamma(a) a^{q - q_{\mathfrak{o}}} =  {''a^q/\varpi}'', \quad
  a \in \mathfrak{a}  
  \end{equation}
 By setting $\gamma_{\rm a}(a) = \gamma_{\rm t}(a) (\varpi/\pi)$, we  obtain divided powers $\gamma_{\rm a}$ relative to $O$ on $\mathfrak{a}$.

Let $\mathcal{P} = (P, Q, F, \dot{F})$ be a
$\mathcal{W}_{\mathfrak{o}}(S/R)$-display with a strict action
\begin{displaymath}
\iota: O \longrightarrow \End \mathcal{P}. 
\end{displaymath}
The definition of strictness is literally the same as Definition \ref{Adorf1d}. 
Since $(S \longrightarrow R, \gamma_{\rm a})$ is an $O$-$pd$-thickening, the $O$-frame
$\mathcal{W}_{O}(S/R)$ is
defined, cf. Example \ref{ex:crystdisp}. The Ahsendorf functor generalizes to a \emph{Ahsendorf functor for $S/R$}\index{Ahsendorf functor for $S/R$}
\index[NO]{ACB@$\mathfrak{A}_{O/\mathfrak{o}, S/R}$}
\begin{equation}\label{Adorf13e}
  \mathfrak{A}_{O/\mathfrak{o}, S/R}:
\left(
  \begin{array}{l}
    \mathcal{W}_{\mathfrak{o}}(S/R)-\text{displays}\\
    \text{with strict $O$-action} 
    \end{array}
\right)
  \; \longrightarrow \; \Big(\mathcal{W}_{O}(S/R)-\text{displays}\Big).  
\end{equation}
The construction is the same but uses some additional arguments, which we
will indicate now.

 We define the \emph{Ahsendorf frame}\index{Ahsendorf frame for $S/R$}
for $S \longrightarrow R$ relative to $\mathfrak{o}$,   
\begin{equation}\label{Aframe2e}
\mathcal{A}_{\mathfrak{o}}(S/R) = (W_{\mathfrak{o}}(S), I_{\mathfrak{o}}(S/R),
R, \mathfrak{f}^f, \mathfrak{f}^{f-1}\dot{\mathfrak{f}}),    
\end{equation}
where $\dot{\mathfrak{f}}: I_{\mathfrak{o}}(S/R) \longrightarrow W_{\mathfrak{o}}(S)$
is defined as in Example \ref{ex:crystdisp}.  This is an $O^t$-frame by the homomorphism
$\varkappa: O^{t} \longrightarrow W_{\mathfrak{o}}(S)$. 

From the $O^t$-action on the $W_{\mathfrak{o}}(S)$-module $P$ we obtain a
decomposition, comp (\ref{Aunverzweigt1e}),
\begin{equation}\label{Adorf25e}
  P = \oplus_{m \in \mathbb{Z}/f\mathbb{Z}} P_m, \quad  Q =
  \oplus_{m \in \mathbb{Z}/f\mathbb{Z}} Q_m.
  \end{equation}
We obtain an $\mathcal{A}_{\mathfrak{o}}(S/R)$-display
$\mathcal{P}_{\rm u a} = (P_{\rm u a}, Q_{\rm u a}, F_{\rm u a}, \dot{F}_{\rm u a})$ by the
formulas (\ref{Adorf14e}).

\begin{lemma}
The Drinfeld homomorphism
$\mu: W_{\mathfrak{o}}(S) \longrightarrow W_{O^t}(S)$ induces a morphism of frames
\begin{equation}\label{Adorf18e}
\mathcal{A}_{\mathfrak{o}}(S/R) \longrightarrow \mathcal{W}_{O^t}(S/R). 
  \end{equation}
\end{lemma}
\begin{proof}
We have to prove the formula
\begin{displaymath}
  ~^{\dot{F}'}(\mu(\eta)) = 
  \mu(~^{\mathfrak{f}^{f-1} \dot{\mathfrak{f}}}\eta), \quad
  \eta \in I_{\mathfrak{o}}(S/R). 
\end{displaymath}
For $\eta \in I_{\mathfrak{o}}(S)$, this is (\ref{Adorf16e}). Therefore the
formula follows if we show that $\mu: W_{\mathfrak{o}}(S) \longrightarrow W_{O^t}(S)$
maps logarithmic Teichm\"uller representatives of elements in 
$\mathfrak{a}$ to logarithmic Teichm\"uller representatives. 
Let $\dot{\mathbf{w}}_{\mathfrak{o},n}$ be the divided Witt polynomials defined 
by $\gamma$ and let $\dot{\mathbf{w}}_{O^t,n}$ be the divided Witt polynomials
defined by $\gamma_{\rm t}$. It follows from the definition of the Drinfeld
homomorphism (\ref{Drinfmap2e}) that
\begin{equation}\label{Adorf17e}
  \dot{\mathbf{w}}_{O^t,n}(\mu(\xi)) = \varpi^{(f-1)n}
  \dot{\mathbf{w}}_{\mathfrak{o},fn} (\xi), \quad \xi \in
  W_{\mathfrak{o}}(\mathfrak{a}). 
\end{equation}
This is verified by reducing to a universal case where $\mathfrak{a}$ is
without $p$-torsion. If now $\xi = \tilde{a} \in W_{\mathfrak{o}}(\mathfrak{a})$
is a logarithmic Teichm\"uller representative, the right hand side of
(\ref{Adorf17e}) is $0$ for $n \neq 0$, and is $a$ for $n = 0$. This shows that
$\mu(\mathfrak{a})$ is the logarithmic Teichm\"uller representative of $a$ in
$W_{O^t}(\mathfrak{a})$.
  \end{proof}

Applying now base change to $\mathcal{P}_{\rm u a}$ relative to (\ref{Adorf18e}),
we obtain a $\mathcal{W}_{O^t}(S/R)$-display $\mathcal{P}_{\rm t}$ with a strict
$O$-action. The assignment $\mathcal{P} \mapsto \mathcal{P}_{\rm t}$ defines
the functor
\begin{displaymath}
  \mathfrak{A}_{O^t/\mathfrak{o}, S/R}:
\left(
  \begin{array}{l}
    \mathcal{W}_{\mathfrak{o}}(S/R)-\text{displays}\\
    \text{with strict $O$-action} 
    \end{array}
\right)
  \; \longrightarrow \; 
\left(
  \begin{array}{l}
    \mathcal{W}_{O^t}(S/R)-\text{displays}\\
    \text{with strict $O$-action} 
    \end{array}
\right) .
\end{displaymath}

Next we  define the functor
\begin{equation}\label{Adorf24e}
  \mathfrak{A}_{O/O^t, S/R}:
\left(
  \begin{array}{l}
    \mathcal{W}_{O^t}(S/R)-\text{displays}\\
    \text{with strict $O$-action} 
    \end{array}
  \right)
  \; \longrightarrow \;
   \Big( \mathcal{W}_{O}(S/R)-\text{displays}\Big). 
\end{equation}
We begin with the definition of the Lubin-Tate frame
$\mathcal{F}_{\rm{lt}}(S/R)$ for $S \rightarrow R$ relative to $O$. We start
with the frame
$\mathcal{W}_{O^t}(S/R) = (W_{O^t}(S), I_{O^t}(S/R), F', \dot{F}')$. Recall
that $I_{O^t}(S/R) = \tilde{\mathfrak{a}} \oplus I_{O^t}(S)$, where the ideal
$\tilde{\mathfrak{a}}$ consists of the logarithmic Teichm\"uller
representatives $\tilde{a}$ of elements $a \in \mathfrak{a}$, with respect to
the divided powers $\gamma_t$. We have by definition $\dot{F}'(\tilde{a}) = 0$.
Tensoring with $O \otimes_{O^t}$, we obtain
\begin{displaymath}
  \begin{aligned}
   F'\colon & O \otimes_{O^t} W_{O^t}(S)& \longrightarrow  O \otimes_{O^t} W_{O^t}(S), \\ 
  \dot{F}'\colon & O \otimes_{O^t} I_{O^t}(S/R)& \longrightarrow O \otimes_{O^t} W_{O^t}(S).
    \end{aligned}
\end{displaymath}

We define an ideal in $O \otimes_{O^t} W_{O^t}(S)$,
\begin{equation}\label{NK7e}
  \mathcal{J}(S/R) = O \otimes_{O^t} I_{O^t}(S/R) +
  (\pi \otimes 1 - 1 \otimes [\pi])
  (O \otimes_{O^t} W_{O^t}(S)). 
\end{equation}
For an element $\eta \in \mathcal{J}(S/R)$ we find
\begin{displaymath}
\tilde{\mathbf{E}}_{K,0}(\pi \otimes 1) \eta \in O \otimes_{O^t} I_{O^t}(S/R).
  \end{displaymath}
Indeed, the factor ring
$O \otimes_{O^t} W_{O^t}(S)/O \otimes_{O^t} I_{O^t}(S/R) = O \otimes_{O^t} R$
is annihilated by $\mathbf{E}_K(\pi \otimes 1)$.
As before in the definition of the Lubin-Tate frame, we define maps
$\sigma_{\rm{lt}}: O \otimes_{O^t} W_{O^t}(S) \longrightarrow
O \otimes_{O^t} W_{O^t}(S)$ and 
$\dot{\sigma}_{\rm{lt}}: \mathcal{J}(S/R) \longrightarrow O \otimes_{O^t}
W_{O^t}(S)$
by
\begin{displaymath}
  ~^{\sigma_{\rm{lt}}} \xi = ~^{F'} \xi, \; \; 
  ~^{\dot{\sigma}_{\rm{lt}}} \eta =
  ~^{\dot{F}'}\tilde{\mathbf{E}}_K(\pi \otimes 1)^{-1}
  ~^{\dot{F}'}(\tilde{\mathbf{E}}_{K,0}(\pi \otimes 1) \eta),
\end{displaymath}
with $\xi \in O \otimes_{O^t} W_{O^t}(S)$, $\eta \in \mathcal{J}(S/R)$.  
The justification of the following definition is analogous to the justification of Definition \ref{LTF3p} of the Lubin-Tate frame. 
\begin{definition}
  The $O$-frame 
  \begin{displaymath}
    \mathcal{F}_{\rm{lt}}(S/R) :=
    (O \otimes_{O^t} W_{O^t}(S), \mathcal{J}(S/R),
    {\sigma_{\rm{lt}}}, {\dot{\sigma}_{\rm{lt}}})
  \end{displaymath}
  is called the \emph{Lubin-Tate frame} \index{Lubin-Tate frame for $S/R$}
  for the epimorphism of $O$-algebras $S \longrightarrow R$ and the divided
  powers $\gamma_t$ relative to $O^t$ on the kernel $\mathfrak{a}$.
\end{definition}
\begin{lemma}
 The Drinfeld homomorphism
  \begin{displaymath}
\mu: O \otimes_{O^t} W_{O^t}(S) \longrightarrow W_{O}(S) 
  \end{displaymath}
  defines a morphism of $O$-frames
  \begin{equation}\label{FLT8e} 
\mathcal{F}_{\rm{lt}}(S/R) \longrightarrow \mathcal{W}^{\varepsilon}_{O}(S/R). 
  \end{equation}
The last frame is defined by (\ref{Adorf10e}). 
\end{lemma}
\begin{proof}
  One can argue exactly as in the proof of Proposition \ref{LTF2p},
  but we need that
  \begin{displaymath}
\mu: W_{O^t}(\mathfrak{a}) \longrightarrow W_O(\mathfrak{a}) 
  \end{displaymath}
  maps logarithmic Teichm\"uller representatives
  $\tilde{a} \in W_{O^t}(\mathfrak{a})$ with respect to the $O^t$-divided powers
  $\gamma_t$ to logarithmic Teichm\"uller representatives
  $\tilde{a} \in W_{O}(\mathfrak{a})$ with respect to the $O$-divided powers
  $\gamma_{\rm a}$.
  This is a consequence of the following relation of  divided Witt
  polynomials,
  \begin{displaymath}
    \dot{\mathbf{w}}_{O,n}(\mu(\alpha)) = \left(\frac{\varpi}{\pi}\right)^n
    \dot{\mathbf{w}}_{O^t,n}(\alpha), \quad \alpha \in W_{O^t}(\mathfrak{a}). 
  \end{displaymath}
  Again we may restrict to the $p$-torsionfree case, where this formula follows
  immediately from the definition of $\mu$, cf.  (\ref{Drinfmap2e}). 
\end{proof}

Let $\mathcal{P}_{\rm t}$ be a $\mathcal{W}_{O^t}(S/R)$-display with a strict
$O$-action. Then the $R$-module $P_{\rm t}/Q_{\rm t}$ is annihilated by
$\pi \otimes 1 - 1 \otimes [\pi]$. As in Proposition \ref{FLT4p}, we define
$F_{\rm{lt}}: P \longrightarrow P$ by
\begin{displaymath}
 F_{\rm{lt}}(x) = \dot{F}_{\rm t}((\pi \otimes 1 - 1 \otimes [\pi])x), \quad
  x \in P. 
  \end{displaymath}
We set $P_{\rm{lt}} = P_{\rm t}$, $Q_{\rm{lt}} = Q_{\rm t}$,
$\dot{F}_{\rm{lt}} = \dot{F}_{\rm t}$. Then we obtain a
$\mathcal{F}_{\rm{lt}}(S/R)$-display
$\mathcal{P}_{\rm{lt}} = (P_{\rm{lt}}, Q_{\rm{lt}}, F_{\rm{lt}},
\dot{F}_{\rm{lt}})$.
If we apply the base change by (\ref{FLT8e}), we obtain a
$\mathcal{W}^{\varepsilon}_{O}(S/R)$-display
$\mathcal{P}_{\rm a}^{\varepsilon} = (P^\varepsilon_{\rm a}, Q^\varepsilon_{\rm a},  F^\varepsilon_{\rm a}, \dot{F}^\varepsilon_{\rm a})$.
The assignement
$\mathcal{P}_{\rm t} \mapsto \mathcal{P}_{\rm a} = (P^\varepsilon_{\rm a}, Q^\varepsilon_{\rm a}, \varepsilon^{-1}F^\varepsilon_{\rm a}, \dot{F}^\varepsilon_{\rm a})$ is the
desired relative Ahsendorf functor $\mathfrak{A}_{O/O^t, S/R}$. 

\begin{proposition}\label{Adorf5p}
  Let $\mathcal{P}$ be a $\mathcal{W}_{\mathfrak{o}}(S/R)$-display with a strict
  $O$-action. Let $\mathcal{P}_{\rm a}$ \index[NO]{PBD@$\mathcal{P}_{\rm a}$}
  be the $\mathcal{W}_O(S/R)$-display associated
  to it by the relative Ahsendorf functor $\mathfrak{A}_{O/\mathfrak{o}, S/R}$. Then there is a canonical isomorphism
  \begin{equation}\label{Adorf15e}
P_{\rm a}/I_{O}(S)P_{\rm a} \cong S \otimes_{O \otimes_{\mathfrak{o}} S} (P/I_{\mathfrak{o}}(S) P).  
    \end{equation}
\end{proposition}
\begin{proof}
This is an immediate consequence of (\ref{Adorf12e})
  \end{proof}
With the notation $P = P_{\rm u a}$ of \eqref{Adorf25e}, we may write
\begin{equation}\label{Adorf27e}
  P_{\rm a}/I_{O}(S)P_{\rm a} = P_{\rm u a}/I_{\mathfrak{o}}(S)P_{\rm u a} +
  (\pi \otimes 1 - 1 \otimes [\pi])P_{\rm u a}. 
\end{equation}
To see this, one uses that $\pi \otimes 1 - 1 \otimes \pi$ generates the
kernel of the canonical map $O \otimes_{O^t} O \longrightarrow O$ as an ideal.  
We see that $P_{\rm a}/I_{O}(S)P_{\rm a}$ is the biggest quotient of
$P_{\rm u a}/I_{\mathfrak{o}}(S)P_{\rm u a}$ such that the action via $\iota$ and via the
structure homomorphism $O \longrightarrow S$ agree.

Let $R$ be an $O$-algebra $R$ such that $\pi$ is nilpotent in $R$. 
Let  $\mathcal{P}$ be a $\mathcal{W}_{\mathfrak{o}}(R)$-display with strict $O$-action.
We assume that $\mathcal{P}$ is nilpotent. Then $\mathcal{P}_{\rm a}$ is also
nilpotent. Then there is a crystal $\mathbb{D}_{\mathcal{P}}$ on the category
of $\mathfrak{o}$-$pd$-thickenings and a crystal $\mathbb{D}_{\mathcal{P}_{\rm a}}$ on
the category of ${O}$-$pd$-thickenings associated to these displays. 

\begin{corollary}\label{Adorf5c}
  Let $\mathcal{P}$ be a nilpotent  $\mathcal{W}_{\mathfrak{o}}(R)$-display with a
  strict $O$-action. Then the image $\mathcal{P}_{\rm a}$ by the
  Ahsendorf functor is a nilpotent $\mathcal{W}_O(R)$-display.
  Let $S \longrightarrow R$ be a surjective map of $O$-algebras which are $p$-adic.
Assume that the kernel $\mathfrak{a}$ of this epimorphism is endowed
with divided powers $\gamma$ relative to $\mathfrak{o}$. Let   $\gamma_{\rm a}$ be the corresponding $O$-divided powers on $\mathfrak{a}$. There is a canonical isomorphism
  \begin{displaymath}
    \mathbb{D}_{\mathcal{P}_{\rm a}}(S, \gamma_{\rm a}) = S \otimes_{(O \otimes_{\mathfrak{o}} S)}
    \mathbb{D}_{\mathcal{P}}(S, \gamma).  
    \end{displaymath}
\end{corollary}
\begin{proof}
Indeed, $\mathbb{D}_{\mathcal{P}}(S)$ is computed from a
$\mathcal{W}_{\mathfrak{o}}(S/R)$-display $\tilde{\mathcal{P}}$ which lifts
$\mathcal{P}$ and which is unique up to isomorphism. But then the relative 
Ahsendorf functor applied to $\tilde{\mathcal{P}}$ gives a
$\mathcal{W}_{O}(S/R)$-display $\tilde{\mathcal{P}}_{\rm a}$ which lifts
$\mathcal{P}_{\rm a}$. We conclude by Proposition \ref{Adorf5p}.   
  \end{proof}
\begin{corollary}
  With the notation of Corollary \ref{Adorf5c}, the Ahsendorf functor
  $\mathfrak{A}_{O/\mathfrak{o}, S}$ defines a bijection between the
  liftings of $\mathcal{P}$ to a $\mathcal{W}_{\mathfrak{o}}(S)$-display
  with a strict $O$-action and the liftings of $\mathcal{P}_{\rm a}$ to a
  $\mathcal{W}_{O}(S)$-display. 
\end{corollary}
\begin{proof}
  We show that each lifting of $\mathcal{P}_{\rm a}$ is in the essential image of
  $\mathfrak{A}_{O/\mathfrak{o}, S}$. A lifting $\tilde{\mathcal{P}}_{\rm a}$ 
  of $\mathcal{P}_{\rm a}$ corresponds, by
  Grothendieck-Messing for  nilpotent displays, to a direct summand
  $U_{\rm a} \subset \mathbb{D}_{\mathcal{P}_{\rm a}}(S)$. Let
  $U \subset \mathbb{D}_{\mathcal{P}}$ be the preimage of $U_{\rm a}$ by the natural
  epimorphism $\mathbb{D}_{\mathcal{P}} \longrightarrow \mathbb{D}_{\mathcal{P}_{\rm a}} $.
  Then $U$
  defines a lifting of $\mathcal{P}$ which is mapped by the Ahsendorf functor
  to $\tilde{\mathcal{P}}_{\rm a}$. 
\end{proof}
It is straightforward to deduce from the last Corollary 
Ahsendorf's Theorem \ref{Adorf1t} for an artinian local ring with perfect
residue class field, i.e., we reproved a special case of \cite{ACZ}.

    \subsection{The Lubin-Tate display} \label{ss:tltd} 
Let $K$ be a finite extension of $\mathbb{Q}_p$.  
Let $K^t \subset K$ be the maximal subextension which is unramified over
$\mathbb{Q}_p$. 
We denote by $O^t \subset O$ the rings of integers and by
$\kappa$ the common residue class field. 
We fix a prime element $\pi \in O$. 
Let $L$ be a normal extension of $\mathbb{Q}_p$ which contains $K$. We set
$L^t = K^t$. Let $\Phi = \Hom_{\text{$\mathbb{Q}_p$-Alg}}(K, L)$ and
$\Psi = \Hom_{\text{$\mathbb{Q}_p$-Alg}}(K^t, L)$. We denote by $\varphi_0 \in \Phi$
and $\psi_0 \in \Psi$ the identity embeddings. We denote by $\Phi_{\psi}$ the
preimage of $\psi$ by the restriction map $\Phi \longrightarrow \Psi$. We define
\begin{displaymath}
  \mathbf{E}_{\psi}(T) = \prod_{\varphi \in \Phi_{\psi}} (T - \varphi(\pi)) \in
  O_{L}[T]. 
  \end{displaymath}
Clearly this polynomial has coefficients in $O_{L^t} \subset O_{L}$.
Let $\mathbf{E}\in O^t[T]$
be the Eisenstein polynomial of $\pi$ in the extension $K/K^t$. Then
$\mathbf{E}_{\psi}$ is the image of $\mathbf{E}$ by $\psi$ in $O_{L^t}[T]$.
We consider the surjective $O_L$-algebra homomorphismus
\begin{displaymath}
O_L[T] \longrightarrow O \otimes_{O^t, \psi} O_L, 
  \end{displaymath}
which maps $T$ to $\pi \otimes 1$. Then 
$\mathbf{E}_{\psi}(\pi \otimes 1) = 0$.

We lift the polynomials $\mathbf{E}_{\psi}$ to the Witt ring
\begin{displaymath}
  \tilde{\mathbf{E}}_{\psi}(T) =
  \prod_{\varphi \in \Phi_{\psi}} (T - [\varphi(\pi)]) \in
    W(O_{L^t})[T]. 
\end{displaymath}
We consider the decomposition
\begin{displaymath}
O \otimes_{\mathbb{Z}_p} O_{L^t} = \prod_{\psi \in \Psi} O \otimes_{O^t, \psi} O_{L^t}. 
\end{displaymath}
Let $\sigma \in \Gal(K^t/\mathbb{Q}_p)$ be the Frobenius automorphism. 
We have the morphism $\lambda: O^t \longrightarrow W(O^t)$ from \eqref{maplambda}. We define
$\tilde{\psi}$ as the composite
\begin{displaymath}
  \tilde{\psi}: O^t \overset{\lambda}{\longrightarrow} W(O^t)
  \overset{W(\psi)}{\longrightarrow} W(O_{L^t}). 
  \end{displaymath}
Then we obtain the decomposition
\begin{equation}\label{LTD1e}
  O \otimes_{\mathbb{Z}_p} W(O_{L^t}) = \prod_{\psi \in \Psi}
  O \otimes_{O^t, \tilde{\psi}} W(O_{L^t}). 
\end{equation} 
Let $\tilde{\mathbf{E}}_{\psi}(\pi \otimes 1)$ be the image of
$\tilde{\mathbf{E}}_{\psi}$ by the homomorphism
$W(O_{L^t})[T] \longrightarrow O \otimes_{O^t, \tilde{\psi}} W(O_{L^t})$ 
which maps $T$ to $\pi \otimes 1$. Since
$\mathbf{E}_{\psi}(\pi \otimes 1) = 0$, we conclude that
\begin{equation}\label{LTD6e}
  \tilde{\mathbf{E}}_{\psi}(\pi \otimes 1) \in O \otimes_{O^t, \tilde{\psi}}
  I(O_{L^t}). 
  \end{equation}
For an arbitrary $O_{L^t}$-algebra $R$, the decomposition (\ref{LTD1e}) induces  
\begin{equation}\label{LTD11e} 
  O \otimes_{\mathbb{Z}_p} W(R) = \prod_{\psi \in \Psi}
  O \otimes_{O^t, \tilde{\psi}} W(R). 
\end{equation}
The Frobenius and the Verschiebung act on the left hand side via the second
factor, and this induces on the right hand side the maps
\begin{equation}\label{LTD5e} 
  \begin{aligned}
    F\colon & O \otimes_{O^t, \tilde{\psi}} W(R)&  \longrightarrow&
    \quad O \otimes_{O^t, \widetilde{\psi \sigma}} W(R)\\ 
    & a \otimes \xi&  \mapsto&   \quad a \otimes ~^{F}\xi\\
    V\colon & O \otimes_{O^t, \widetilde{\psi\sigma}} W(R)&  \longrightarrow&
    \quad O \otimes_{O^t, \tilde{\psi}} W(R)\\
    & a \otimes \xi&  \mapsto&   \quad a \otimes ~^{V}\xi
    \end{aligned}
\end{equation}
We note that $\widetilde{\psi \sigma} = \tilde{\psi} \circ \sigma$. 
We will write
$\dot{F} = V^{-1}: O \otimes_{O^t, \tilde{\psi}} I(R)  \longrightarrow  
    O \otimes_{O^t, \widetilde{\psi \sigma}} W(R)$.

\begin{proposition}\label{Adorf2p}
    The element
$~^{\dot{F}}\tilde{\mathbf{E}}_{\psi}(\pi \otimes 1) \in
    O \otimes_{O^t, \widetilde{\psi \sigma}} W(O_{L^t})$
   is a unit of the form 
  \begin{equation}\label{Adorf5e}
    ~^{\dot{F}}\tilde{\mathbf{E}}_{\psi}(\pi \otimes 1) = (\frac{\pi^e}{p}
    \otimes 1)\delta, \quad \delta \in
    O \otimes_{O^t, \widetilde{\psi \sigma}} W(O_{L^t}) , 
    \end{equation}
  where $\delta - 1\otimes 1$ is in the kernel of
 $O \otimes_{O^t, \widetilde{\psi \sigma}} W(O_{L^t}) \longrightarrow  
    O \otimes_{O^t, \widetilde{\psi \sigma}} W(\kappa_{L^t})$. 
    In particular, $\delta - 1\otimes 1$ lies in the Jacobson radical
    of $O \otimes_{O^t, \widetilde{\psi \sigma}} W(O_{L^t})$.   
    \end{proposition}
\begin{proof}
  The proof is the same as that of Proposition \ref{NK1p}.
  \end{proof} 

The polynomial $\mathbf{E}_{\psi_0} $ has the decomposition
\begin{displaymath}
  \mathbf{E}_{\psi_0}(T) = (T - \pi) \cdot \mathbf{E}_{0}(T).
\end{displaymath}
These polynomials have coefficients in $K \subset L$. Recall that
$\varphi_0(\pi) = \pi$. We set 
\begin{displaymath}
 \tilde{\mathbf{E}}_{0}(T) =
 \prod_{\varphi \in \Phi_{\psi_0}, \varphi \neq \varphi_0} (T - [\varphi(\pi)]) \in
 W(O_L)[T]. 
  \end{displaymath}
This polynomial lies in $W(O)[T] \subset W(O_L)[T]$. 
We set 
\begin{equation}\label{Adorf2e} 
P_{\mathcal{L}} =   O \otimes_{\mathbb{Z}_p} W(O) = \oplus_{\psi \in \Psi}
  O \otimes_{O^t, \tilde{\psi}} W(O).
\end{equation}

We denote by $Q_{\psi_0} \subset O \otimes_{O^t, \tilde{\psi_0}} W(O)$ the kernel
of the map
\begin{equation}\label{Adorf3e}
  O \otimes_{O^t, \tilde{\psi_0}} W(O) 
  \xra{\id \otimes \mathbf{w}_0}O \otimes_{O^t, \psi_0} O \xra{\rm mult.}O. 
\end{equation}
\begin{lemma}\label{LTD5l}
  \begin{displaymath}
  Q_{\psi_0} = \{x \in O \otimes_{O^t, \tilde{\psi_0}} W(O) \; | \;
  \tilde{\mathbf{E}}_{\psi_0, 0}(\pi \otimes 1) x \in
  O \otimes_{O^t, \tilde{\psi_0}} I(O) \}
    \end{displaymath}
\end{lemma} 
\begin{proof}
  Let $\bar{Q}_{\psi_0}$ be the kernel of the second map of (\ref{Adorf3e}). Then
  we can reformulate our assertion as 
  \begin{displaymath}
    \bar{Q}_{\psi_0} = \{x \in O \otimes_{O^t, \psi_0} O \; | \;
    \mathbf{E}_{\psi_0, 0}(\pi \otimes 1) x = 0 \}. 
  \end{displaymath}
We write
  \begin{equation*}
    O \otimes_{O^t, \psi_0} O \simeq \big(O^t[T]/\mathbf{E}_{\psi_0}(T) O^t[T]\big)
    \otimes_{O^t, \psi_0} O \simeq  O[T]/\mathbf{E}_{\psi_0}(T) O[T]. 
  \end{equation*}
  We see that
  $\bar{Q}_{\psi_0} = (T- \pi)O[T]/\mathbf{E}_{\psi_0}(T)O[T]$ which coincides 
with $\mathbf{E}_{0}(T)O[T]/\mathbf{E}_{\psi_0}(T)O[T]$. 
  \end{proof}
We let $F$ and $V$ act via the second factor on
$O \otimes_{\mathbb{Z}_p} W(O)$ and therefore on the right hand side of
(\ref{Adorf2e}) by the formulas (\ref{LTD5e}). 
\begin{definition}\label{LTD1d}
  The \emph{Lubin-Tate display}
  \index{Lubin-Tate display}\index[NO]{LBA@$\mathcal{L}$} is the
  $\CW(O)$-display
  $\mathcal{L} =(P_{\mathcal{L}}, Q_{\mathcal{L}}, F_{\mathcal{L}},\dot{F}_{\mathcal{L}})$,
  defined as follows. 
  Let $P_{\mathcal{L}} = O \otimes_{\mathbb{Z}_p} W(O)$. Then
  $P_{\mathcal{L}} = \oplus_{\psi} P_{\psi, \mathcal{L}}$ with
  $P_{\psi, \mathcal{L}} = O \otimes_{O^t, \tilde{\psi}} W(O)$. Set 
  $Q_{\psi} = P_{\psi, \mathcal{L}}$ for $\psi \neq \psi_0$ and
   $Q_{\psi_0} \subset P_{\psi_0, \mathcal{L}}$ as above, and  define 
  \begin{displaymath}
Q_{\mathcal{L}} = \oplus_{\psi \in \Psi} Q_{\psi} \subset P_{\mathcal{L}}  .
  \end{displaymath}
The maps 
  $F_{\mathcal{L}}$ and $\dot{F}_{\mathcal{L}}$ are defined as the direct sum of the following
  maps for all $\psi$. For $\psi \neq \psi_0$ we define 
  \begin{displaymath}
    \begin{array}{rcccrccc} 
      F_{\mathcal{L}}: &  O \otimes_{O^t, \tilde{\psi}} W(O) & \longrightarrow & 
      O \otimes_{O^t, \widetilde{\psi \sigma}} W(O), &
      \dot{F}_{\mathcal{L}}: & O \otimes_{O^t, \tilde{\psi}} W(O) & \longrightarrow & 
O \otimes_{O^t, \widetilde{\psi \sigma}} W(O).\\[2mm] 
& x & \mapsto & ~^{F}(\tilde{\mathbf{E}}_{\psi}(\pi \otimes 1)x), &
& y & \mapsto & ~^{\dot{F}}(\tilde{\mathbf{E}}_{\psi}(\pi \otimes 1)y) 
  \end{array}
  \end{displaymath}

\noindent For $\psi_0$ we define 
  \begin{displaymath}
    \begin{array}{rcccrccc} 
F_{\mathcal{L}}: & O \otimes_{O^t, \tilde{\psi}_0} W(O) & \longrightarrow & 
O \otimes_{O^t, \widetilde{\psi_0 \sigma}} W(O), &
\dot{F}_{\mathcal{L}}: & Q_{\psi_0} & \longrightarrow & 
O \otimes_{O^t, \widetilde{\psi_0 \sigma}} W(O).\\[2mm]
& x & \mapsto & ~^{F}(\tilde{\mathbf{E}}_{0}(\pi \otimes 1)x), &
& y & \mapsto & ~^{\dot{F}}(\tilde{\mathbf{E}}_{0}(\pi \otimes 1)y) 
      \end{array}
    \end{displaymath}
The action of $O$ by
  multiplication via the first factor on $O \otimes_{\mathbb{Z}_p} W(O)$
  defines a strict $O$-action on $\mathcal{L}$. If $R$ is a $p$-adic
  $O$-algebra we denote by $\mathcal{L}_R$ the base change of $\mathcal{L}$
  via the morphism of frames $\CW(O) \longrightarrow \CW(R)$.  
\end{definition}
The tuple $ (P_{\mathcal{L}}, Q_{\mathcal{L}}, F_{\mathcal{L}}, \dot{F}_{\mathcal{L}})$ is indeed a $\CW(O)$-display. The only non-trivial point is that $\dot{F}_{\mathcal{L}}$ is surjective. But this follows from Proposition
\ref{Adorf2p}. 

We will now apply the Ahsendorf functor $\mathfrak{A}_{O/\mathbb{Z}_p, O}$ to the Lubin-Tate display
$\mathcal{L}$. We use the previous section in the case
$\mathfrak{o} = \mathbb{Z}_p$ and $\varpi = p$.  
We set $\psi_i = \psi_0 \sigma^i: O^t \longrightarrow O$, where $\psi_0$ is the
identity. We have the decomposition
\begin{displaymath}
  P_{\mathcal{L}} = \oplus_{i=0}^{f-1} P_{\mathcal{L}, \psi_i}, \quad
  P_{\mathcal{L}, \psi_i} = O \otimes_{O^t, \tilde{\psi_i}} W(O). 
  \end{displaymath}
We denote the Frobenius on $W(O)$ by $F$. In the last section we have
associated to $\mathcal{L}$ an $\mathcal{A}(O)$-display $\mathcal{P}_{\rm u a}$
for the $O^t$-frame
\begin{displaymath}
\mathcal{A}_{\BZ_p}(O)=\mathcal{A}(O) = (W(O), I(O), F^f, F^{f-1}\dot{F}). 
  \end{displaymath}
By definition $P_{\rm u a} = P_{\mathcal{L}, \psi_0}$, $Q_{\rm u a} = Q_{\mathcal{L}, \psi_0}$,  
and $\dot{F}_{\rm u a} = \dot{F}_{\mathcal{L}}^f$. More explicitely, for
$y \in Q_{\rm u a}$,
\begin{displaymath}
  \dot{F}_{\rm u a}(y) = ~^{\dot{F}}\!\tilde{\mathbf{E}}_{\psi_{f-1}}(\pi \otimes 1)
  ~^{F\dot{F}}\!\tilde{\mathbf{E}}_{\psi_{f-2}}(\pi \otimes 1) \cdot \ldots \cdot
  ~^{F^{f-2}\dot{F}}\!\tilde{\mathbf{E}}_{\psi_1}(\pi \otimes 1)
  ~^{F^{f-1}\dot{F}}\!(\tilde{\mathbf{E}}_{0}(\pi \otimes 1)y). 
  \end{displaymath}
We set
\begin{displaymath}
\mathbf{n} = ~^{\dot{F}}\!\tilde{\mathbf{E}}_{\psi_{f-1}}(\pi \otimes 1)
  ~^{F\dot{F}}\!\tilde{\mathbf{E}}_{\psi_{f-2}}(\pi \otimes 1) \cdot \ldots \cdot
  ~^{F^{f-2}\dot{F}}\!\tilde{\mathbf{E}}_{\psi_1}(\pi \otimes 1)
  ~^{F^{f-1}\dot{F}}\!\tilde{\mathbf{E}}_{\psi_0}(\pi \otimes 1). 
  \end{displaymath}
Then we may write
\begin{equation}
  \dot{F}_{\rm u a}(y) = \mathbf{n}
  \big(~^{F^{f-1}\dot{F}}\!\tilde{\mathbf{E}}_{\psi_0}(\pi \otimes 1)\big)^{-1}
  ~^{F^{f-1}\dot{F}}\!(\tilde{\mathbf{E}}_{0}(\pi \otimes 1)y). 
  \end{equation}
To $\mathcal{P}_{\rm u a}$ we apply base change with respect to the Drinfeld
morphism $\mu: \mathcal{A}(O) \longrightarrow \mathcal{W}_{O^t}(O)$, cf.
(\ref{Adorf20e}). We obtain the $\mathcal{W}_{O^t}(O)$-display
$\mathcal{P}_{\rm t}$, where
\begin{displaymath}
P_{\rm t} = O \otimes_{O^t} W_{O^t}(O)
\end{displaymath}
and where  $Q_{\rm t}$ is the kernel of the homomorphism
$O \otimes_{O^t} W_{O^t}(O) \longrightarrow O$ induced by $\mathbf{w}_{O^t,0}$, cf.
(\ref{Adorf3e}). The polynomials $\tilde{\mathbf{E}}_{\psi_0}$ and
$\tilde{\mathbf{E}}_{0}$ are mapped by $\mu$ to the polynomials
$\tilde{\mathbf{E}}_{K}$ and $\tilde{\mathbf{E}}_{K,0}$ of (\ref{Adorf21e}).
We denote by $\mathbf{n}_{\rm{lt}}$ the image of $\mathbf{n}$ by $\mu$.
Therefore we obtain
\begin{displaymath}
  \dot{F}_{\rm t}(y) = \mathbf{n}_{\rm{lt}}
  \big(~^{\dot{F}'}\tilde{\mathbf{E}}_{K}(\pi \otimes 1)\big)^{-1}
  ~^{\dot{F}'}(\tilde{\mathbf{E}}_{K,0} (\pi \otimes 1)y), \quad y \in Q_{\rm t}.  
\end{displaymath}
By Proposition \ref{FLT4p}, we associate to $\mathcal{P}_{\rm t}$ a
$\mathcal{F}_{\mathfrak{lt} }(O)$-display $\mathcal{P}_{\rm{lt}}$.  
In terms of the Lubin-Tate frame (comp. Definition \ref{LTF3p}),  we may rewrite the last equation as 
\begin{displaymath}
  \dot{F}_{\rm{lt}}(y) = \mathbf{n}_{\rm{lt}}
  \dot{\sigma_{\rm{lt}}}(y), \quad
  y \in Q_{\rm t} = Q_{\rm{lt}} = \mathcal{J}. 
\end{displaymath}

\begin{proposition}\label{LTD8p}
  The Ahsendorf functor $\mathfrak{A}_{O/\mathbb{Z}_p, O}$  maps the Lubin-Tate
  display $\mathcal{L}$ to a $\CW_O(O)$-display which is canonically isomorphic to
  $\mathcal{P}_{m, \mathcal{W}_{O} (O)}(\pi^{ef}/p^f)$, i.e., the twist of the
  multiplicative display by $\pi^{ef}/p^f\in O \subset W_O(O)$,  
  cf. Example \ref{def:twistdisplay}.  
\end{proposition}
\begin{proof}
 Let $\mathcal{P}_{m, \mathfrak{lt}}$ be the multiplicative
  $\mathcal{F}_{\rm{lt}}(O)$-display. 
  The above identities show that the display  $\mathcal{L}_{\mathfrak t}$ is equal to the $\mathcal{F}_{\rm{lt}}(O)$-display
  $\mathcal{P}_{m, \mathfrak{lt}}(\mathbf{n}_{\rm{lt}})$. Applying to this display the base change by the Drinfeld morphism of frames, cf. Proposition \ref{LTF2p}, we obtain
  essentially (i.e., neglecting $\varepsilon$) the image of $\mathcal{L}$ by the Ahsendorf functor. 
  By Proposition \ref{Adorf2p} the image of the element $\mathbf{n}$ by the
  map
  $O \otimes_{O^t, \tilde{\psi}_0} W(O) \longrightarrow O \otimes_{O^t, \psi_0} W(\kappa) =
  O$
  is $\pi^{ef}/p^f$. This implies that the image of 
  $\mathbf{n}_{\rm{lt}} \in O \otimes_{O^t} W_{O^t}(O)$ in
  $O \otimes_{O^t} W_{O^t}(\kappa)$ is $\pi^{ef}/p^f$. 
  By the following Lemma there is a uniquely determined unit
  $\xi \in O \otimes_{O^t} W_{O^t}(O)$
  such that $\xi \mathbf{n}_{\rm{lt}} = (\pi^{ef}/p^f \otimes 1) ~^{F^f} \xi$.
  This gives a canonical isomorphism
  $\mathcal{P}_{m, \mathfrak{lt}}(\mathbf{n}_{\rm{lt}}) \isoarrow \mathcal{P}_{m, \mathfrak{lt}}(\pi^{ef}/p^f \otimes 1)$. 
  The base change with respect to a morphism of frames maps the multiplicative
  display to the multiplicative display, cf. Example \ref{basechangemult}. Therefore the last display is mapped
  by the base change of Proposition \ref{LTF2p} to the
  $\mathcal{W}^{\varepsilon}_{O}(O)$-display
  \begin{displaymath}
\big(W_O(O), I_O(O), F, \varepsilon^{-1} \dot{F})(\pi^{ef}/p^f \otimes 1\big). 
  \end{displaymath}
  Here we denote by $F$ the Frobenius acting on $W_O(O)$ and by $\dot{F}$ the
  inverse of the Verschiebung. The element $\varepsilon$ is defined by
(\ref{NK16e}). 

  Therefore the definition of the Ahsendorf functor  gives
  $\mathcal{P}_{m,\mathcal{W}_{O}(O)}(\varepsilon^{-1} (\pi^{ef}/p^f \otimes 1))$
  as the image of $\mathcal{L}$. The image of $\varepsilon$ by the homomorphism
  $W_O(O) \longrightarrow W_O(\kappa)$ is $1$. A variant of the next Lemma
  shows that there is a unique $\xi \in W_O(O)$ such that
  $~^{F}\xi \xi = \varepsilon$. This shows that the last display is
  canonically isomorphic to
  $\mathcal{P}_{m,\mathcal{W}_{O}(O)}(\pi^{ef}/p^f \otimes 1)$.  
\end{proof}
\begin{lemma}
  Let $\alpha \in O \otimes_{O^t} W_{O^t}(O)$ be a unit whose image in
  $O \otimes_{O^t} W_{O^t}(\kappa)$ is $1$. Then there exists a unique unit 
  $\xi \in O \otimes_{O^t} W_{O^t}(O)$ whose image in
  $O \otimes_{O^t} W_{O^t}(\kappa)$ is $1$ and such that
  \begin{displaymath}
~^{F^f}\!\xi \cdot \xi^{-1} = \alpha. 
    \end{displaymath}
\end{lemma}
\begin{proof}
  One proves this by induction on $n$ for $O \otimes_{O^t} W_{O^t}(O/\pi^n O)$.
  Alternatively one can use Grothendieck-Messing for frames due to Lau
  and show that the multiplicative display of
  $\mathcal{F}_{\rm{lt}}(\kappa)$ has no nontrivial deformation
  with respect to
  $\mathcal{F}_{\rm{lt}}(O) \longrightarrow
  \mathcal{F}_{\rm{lt}}(\kappa)$.  
\end{proof}

\begin{remark}\label{LTD1r}
  Let $k$ a perfect field which is an extension of $\kappa$. We regard it as
  an $O$-algebra via the residue class map $O \longrightarrow \kappa$. Then we can
  describe the Dieudonn\'e module
  $(P_{\mathcal{L}_k}, F_{\mathcal{L}_k}, V_{\mathcal{L}_k})$ of $\mathcal{L}_k$
  as follows.
  Let $\tilde{\psi}_0: O^t = W(\kappa) \longrightarrow W(k)$ be the canonical map.
  The set $\Psi$ consists of the maps $\tilde{\psi}_0 \sigma^i$ where $\sigma$
  is the Frobenius of the extension $O^t/\mathbb{Z}_p$. We have
  \begin{displaymath}
    P_{\mathcal{L}_k} = O \otimes_{\mathbb{Z}_p} W(k) = \prod_{\tilde{\psi} \in \Psi}
    O \otimes_{O^t, \tilde{\psi}} W(k). 
  \end{displaymath}
  The Frobenius and the Verschiebung
  \begin{equation}\label{LTD15e}
    O \otimes_{O^t, \tilde{\psi}} W(k)
    \begin{array}{c}
      \overset{F_{\mathcal{L}_k}}{\longrightarrow}\\[-2mm]
      \underset{V_{\mathcal{L}_k}}{\longleftarrow}
    \end{array}
    O \otimes_{O^t, \tilde{\psi} \sigma} W(k)
  \end{equation}
  are defined as follows:
   \begin{equation*}
    \begin{aligned}
      F_{\mathcal{L}_k} (x_{\psi}) &= \pi^e ~^{F}x_{\psi}, 
     \quad &V_{\mathcal{L}_k} (x_{\psi \sigma}) &=
      (p/\pi^e)\, ^{F^{-1}}\!x_{\psi \sigma} , & \text{for} \; \psi \neq \psi_0,\\
      F_{\mathcal{L}_k} (x_{\psi_0}) &= \pi^{e-1} ~^{F}x_{\psi_0}, \quad   &
      V_{\mathcal{L}_k} (x_{\psi_0 \sigma}) &= (p/\pi^{e-1})\, ^{F^{-1}}\!x_{\psi_0 \sigma} .& 
    \end{aligned}
  \end{equation*}
  The upper left index $F$ denotes the action of the Frobenius via the second
  factor on $O \otimes_{\mathbb{Z}_p} W(k)$.
  This description follows easily because $\tilde{\mathbf{E}}_{\psi}(T) = T^e$,
  and therefore
 $\tilde{\mathbf{E}}_{\psi}(\pi \otimes 1) = \pi^e \otimes 1 \in O \otimes_{O^t, \tilde{\psi}} W(k)$
  and $\tilde{\mathbf{E}}_{0}(\pi \otimes 1) = \pi^{e-1}\otimes 1$. Moreover,
  we have $Q_{\psi_0} = \pi O \otimes_{O^t, \tilde{\psi}} W(k)$. 

  We have identified the frame $\mathcal{W}_O(k)$ with
  \begin{displaymath}
    \big(O \otimes_{O^t, \tilde{\psi}_0} W(k), \pi O \otimes_{O^t, \tilde{\psi}_0}
    W(k), k, F^f, F^f \pi^{-1}\big).
  \end{displaymath}
  By  Remark \ref{Adorf1r}, the Ahsendorf functor associates to
  the Dieudonn\'e module
  $(P_{\mathcal{L}_k}, F_{\mathcal{L}_k}, V_{\mathcal{L}_k})$ of $\mathcal{L}_k$
  the $\mathcal{W}_O(k)$-Dieudonn\'e module
  \begin{equation}\label{LTD19e}
    \big(O \otimes_{O^t, \tilde{\psi}_0} W(k), \frac{\pi^{ef}}{p^f} F^f,
    \frac{p^f}{\pi^{ef-1}} F^{-f}\big). 
  \end{equation}
  This is equal to the Dieudonn\'e module of the twisted multiplicative
  display $\mathcal{P}_{m, \mathcal{W}_{O}(k)}(\pi^{ef}/p^f)$, in agreement with  Proposition \ref{LTD8p}. 
  \end{remark}

In the end of this subsection we discuss  the Faltings dual of a display
$\mathcal{P} = (P, Q, F, \dot{F})$ with a strict
$O$-action over an $O$-algebra $R$. We begin with a recipe how to construct such $\mathcal{P}$.
We consider the decomposition induced by (\ref{LTD11e}), 
\begin{displaymath}
  P = \oplus_{\psi} P_{\psi}, \quad Q = Q_{\psi_0} \oplus (\oplus_{\psi \neq \psi_0}
  P_{\psi}).
  \end{displaymath}
Let 
\begin{equation}\label{LTD14e}
  J_{\psi_0}=\Ker\big(O \otimes_{O^t, \tilde{\psi}_0} W(R) \longrightarrow R\big) ,
  \end{equation}
where the map is induced by the structure homomorphism $O \longrightarrow R$ and the homomorphism
$\mathbf{w}_0$. Then
\begin{displaymath}
  J_{\psi_0} = O \otimes_{O^t, \tilde{\psi}_0} I(R) +
  (\pi \otimes 1 - 1 \otimes [\pi]) (O \otimes_{O^t, \tilde{\psi}_0} W(R)) ,
\end{displaymath}
cf. (\ref{NK5e}). To find a normal decomposition $P=T\oplus L$, we start with 
\begin{displaymath}
  P_{\psi_0} = T_{\psi_0} \oplus L_{\psi_0}, \quad Q_{\psi_0} =
  J_{\psi_0} T_{\psi_0} \oplus L_{\psi_0}. 
\end{displaymath}
We define $\phi_{T_{\psi_0}}:T_{\psi_0} \longrightarrow P_{\psi_0 \sigma}$,
\begin{displaymath}
\phi_{T_{\psi_0}}(t) = \dot{F}\big((\pi \otimes 1 - 1 \otimes [\pi])t\big). 
\end{displaymath}
Then we set $T = T_{\psi_0} \subset P$ and
$L = L_{\psi_0} \oplus (\oplus_{\psi \neq \psi_0} P_{\psi})$. Let
$\phi_L: L \longrightarrow P$ be the restriction of $\dot{F}$ to $L$ and
let $\phi_T = \phi_{T_{\psi_0}}$. 
The restriction to $L_{\psi_0}$ is denoted by $\phi_{L_{\psi_0}}$.
\begin{lemma}\label{LTD3l}
  The map
  \begin{displaymath}
\phi_T \oplus \phi_L: T \oplus L \longrightarrow P 
  \end{displaymath}
  is an $F$-linear $O \otimes_{\mathbb{Z}_p} W(R)$-module isomorphism, where
  the Frobenius $F$ acts on $O \otimes_{\mathbb{Z}_p} W(R)$ via the second
  factor. 
  \end{lemma}
\begin{proof}
  For $\psi \neq \psi_0$ the map
  $\dot{F}_{\psi}: P_{\psi} \longrightarrow P_{\psi \sigma}$ is an $F$-linear
  isomorphism. Therefore it suffices to show that the map
  \begin{equation}\label{LTD7e}
    \phi_T \oplus \phi_{L_{\psi_0}}: T_{\psi_0} \oplus L_{\psi_0} \longrightarrow
    P_{\psi_0 \sigma} 
  \end{equation}
  is an $F$-linear isomorphism or, equivalently, an $F$-linear epimorphism.
  Since $P_{\psi_0 \sigma}$ is generated by $\dot{F} (Q_{\psi_0})$ it suffices
  to show that the following elements are in the image of the linearization
  of (\ref{LTD7e}):
  \begin{displaymath}
    \dot{F}(\ell),\; \,\dot{F}((\pi \otimes 1 - 1 \otimes [\pi])t), \; \,
    \dot{F}(~^{V}\eta t), \quad \, \ell \in L_{\psi_0}, t \in T, \eta \in
    O \otimes_{\mathbb{Z}_p} W(R).
  \end{displaymath}
  For the first two elements this is obvious. For the last element this
  follows from the formula
  \begin{equation}\label{LTD17e} 
  Fx = \big(~^{\dot{F}}\tilde{\mathbf{E}}_{\psi_0}(\pi \otimes 1)\big)^{-1} \cdot
  ~^{F}\tilde{\mathbf{E}}_{0}(\pi \otimes 1) \cdot
  \dot{F}\big((\pi \otimes 1 - 1 \otimes [\pi])x\big), \quad x \in P_{\psi_0} ,
  \end{equation}
  which is proved in the same way as Lemma \ref{NK2l}. 
\end{proof}
We omit the proof of the following proposition.
\begin{proposition}\label{LTD5p}
  Let $R$ be a $p$-adic $O$-algebra. 
  Let $T_{\psi_0}$, $L_{\psi_0}$ be $O \otimes_{O^t, \tilde{\psi}_0} W(R)$-modules
  and let $P_{\psi}$ for $\psi \neq \psi_0$ be
  $O \otimes_{O^t, \tilde{\psi}} W(R)$-modules which are free locally on $\Spec R$.
  Set $P_{\psi_0} = T_{\psi_0} \oplus L_{\psi_0}$. Assume  given
  $F$-linear isomorphisms
  \begin{displaymath}
    \phi_{T_{\psi_0}} \oplus \phi_{L_{\psi_0}}: T_{\psi_0} \oplus L_{\psi_0} \longrightarrow
    P_{\psi_0 \sigma}, \quad 
    \phi_{\psi}: P_{\psi} \longrightarrow P_{\psi \sigma}, \; \text{for}\;
    \psi \neq \psi_0, 
  \end{displaymath}
  cf. (\ref{LTD5e}) for the meaning of $F$-linear.
  
  Then there exists a unique display $\mathcal{P} = (P, Q, F, \dot{F})$  with a strict action of $O$ over
  $R$
  such that $P = \oplus_{\psi} P_{\psi}$ and
  $Q = J_{\psi_0} T_{\psi_0}\oplus (\oplus_{\psi \neq \psi_0} P_{\psi})$ with $T = T_{\psi_0}$ and
  $L = L_{\psi_0} \oplus (\oplus_{\psi \neq \psi_0} P_{\psi})$ and such that
  $\phi_T = \phi_{T_{\psi_0}}\colon T\to P$ and
  $\phi_L = \phi_{\psi_0} \oplus (\oplus_{\psi \neq \psi_0} \phi_{\psi}):
  L \longrightarrow P$
  are given by the display structure of $\mathcal{P}$ as in Lemma \ref{LTD3l}. \qed
  \end{proposition}

Let $\mathcal{P}$ be a display  with a strict action by $O$ over $R$, as above.
Then we set
\begin{equation*}
  \begin{aligned} 
  P^{\nabla} &= \Hom_{O \otimes_{\mathbb{Z}_p} W(R)} (P, O \otimes_{\mathbb{Z}_p} W(R)) =
  \oplus_{\psi} P^{\nabla}_{\psi}, \quad \text{where}\\
  P^{\nabla}_{\psi} &= 
  \Hom_{O \otimes_{O^t, \tilde{\psi}} W(R)} (P_{\psi}, O \otimes_{O^t, \tilde{\psi}} W(R)).
    \end{aligned}  
\end{equation*}
We define
\begin{displaymath}
  Q^{\nabla}_{\psi_0} =
  \{ \hat{x} \in P^{\nabla}_{\psi_0}\; | \; \hat{x}(Q_{\psi_0}) \subset J_{\psi_0}\}
  \subset P^{\nabla}_{\psi_0}, \quad 
  Q^{\nabla} = Q^{\nabla}_{\psi_0} \oplus (\oplus_{\psi \neq \psi_0} P^{\nabla}_{\psi}). 
\end{displaymath}
Let
\begin{equation}\label{LTD8e}
\langle \; , \; \rangle_{\rm can}: P \times P^{\nabla} \longrightarrow O \otimes_{\mathbb{Z}_p} W(R) 
\end{equation}
be the canonical perfect $O \otimes_{\mathbb{Z}_p} W(R)$-bilinear form. It
induces pairings
$P_{\psi} \times P^{\nabla}_{\psi} \longrightarrow O \otimes_{O^t ,\tilde{\psi}} W(R)$.
Under the perfect $R$-bilinear form
\begin{displaymath}
P_{\psi_0}/J_{\psi_0} P_{\psi_0} \times P^{\nabla}_{\psi_0}/J_{\psi_0} P^{\nabla}_{\psi_0}
\longrightarrow (O \otimes_{O^t, \tilde{\psi}} W(R))/J_{\psi_0} \simeq R ,
\end{displaymath}
the $R$-submodules $Q_{\psi_0}/J_{\psi_0} P_{\psi_0}$ and
$Q^{\nabla}_{\psi_0}/J_{\psi_0} P^{\nabla}_{\psi_0}$ are orthogonal complements.
\begin{proposition}\label{LTD7p}
  Let $\mathcal{P}$ be a display with a strict $O$-action over $R$. Let
  $P^{\nabla}$ and $Q^{\nabla}$ as above. Then there are unique $F$-linear
  homomorphisms of $O \otimes_{\mathbb{Z}_p} W(R)$-modules
  \begin{displaymath}
    F^{\nabla}: P^{\nabla} \longrightarrow P^{\nabla}, \quad
    \dot{F}^{\nabla}: Q^{\nabla} \longrightarrow P^{\nabla} 
    \end{displaymath}
  such that
  $\mathcal{P}^{\nabla} = (P^{\nabla}, Q^{\nabla}, F^{\nabla}, \dot{F}^{\nabla})$ 
  becomes a display and such that the bilinear form (\ref{LTD8e}) defines
  a bilinear form of displays with a strict $O$-action,
  \begin{equation}\label{LTD9e}
    \langle \; , \; \rangle_{\rm can}: \mathcal{P} \times \mathcal{P}^{\nabla} \longrightarrow
    \mathcal{L} .
    \end{equation}
  We call $\mathcal{P}^{\nabla}$\index[NO]{PBE@$\mathcal{P}^{\nabla}$}  the
  \emph{Faltings dual}\index{Faltings dual} of $\mathcal{P}$. 
\end{proposition}
\begin{proof}
  It follows from the definition that
  \begin{displaymath}
\langle Q  , Q^{\nabla} \rangle_{\rm can} \subset Q_{\mathcal{L}} \subset O \otimes_{\mathbb{Z}_p} W(R). 
    \end{displaymath}
  We will define a display
  $\mathcal{P}^{\nabla} = (P^{\nabla}, Q^{\nabla}, F^{\nabla}, \dot{F}^{\nabla})$.
  Then we will show that
  \begin{equation}\label{LTD10e} 
    \langle\dot{F} y , \dot{F}^{\nabla} \hat{y}\rangle_{\rm can} = \dot{F}_{\mathcal{L}}
    \langle y , \hat{y} \rangle_{\rm can} \quad y \in Q, \; \hat{y} \in Q^{\nabla}. 
  \end{equation}
  This will show that  (\ref{LTD9e}) is a bilinear form of displays with an
  $O$-action. Since the pairing (\ref{LTD8e}) is perfect, the map
  $\dot{F}^{\nabla}$ is uniquely determined by (\ref{LTD10e}).

  We begin  by defining $\mathcal{P}^{\nabla}$. We chose a normal decomposition
  \begin{equation}\label{LTD12e}
    P_{\psi_0} = T_{\psi_0} \oplus L_{\psi_0}, \quad
    Q_{\psi_0} = J_{\psi_0} T_{\psi_0} \oplus L_{\psi_0}
    \end{equation}
  Let $T^{\nabla}_{\psi_0}$ be the orthogonal complement of $T_{\psi_0}$ and let
  $L^{\nabla}_{\psi_0}$ be the orthogonal complement of $L_{\psi_0}$ with respect
  to the perfect bilinear form
  \begin{displaymath}
P_{\psi_0} \times P^{\nabla}_{\psi_0} \longrightarrow O \otimes_{O^t \tilde{\psi}_0} W(R).  
    \end{displaymath}
  We consider the maps $\phi_{T_{\psi_0}}$ and $\phi_{L_{\psi_0}}$. We define maps
  \begin{displaymath}
    \phi_{T^{\nabla}_{\psi_0}}: T^{\nabla}_{\psi_0} \longrightarrow P^{\nabla}_{\psi_0 \sigma},
    \quad
    \phi_{L^{\nabla}_{\psi_0}}: L^{\nabla}_{\psi_0} \longrightarrow P^{\nabla}_{\psi_0 \sigma} 
  \end{displaymath}
  by the equations
  \begin{equation}\label{LTD16e}
    \begin{aligned}
    \langle \phi_{T_{\psi_0}}(t) + \phi_{L_{\psi_0}}(\ell), \phi_{T^{\nabla}_{\psi_0}}(\hat{t})\rangle_{\rm can}
    &= ~^{\dot{F}}\tilde{\mathbf{E}}_{\psi_0}(\pi \otimes 1)
    ~^{F} \langle\ell , \hat{t} \rangle_{\rm can}, \quad &t \in T_{\psi_0}, \ell \in L_{\psi_0}\\
    \langle\phi_{T_{\psi_0}}(t) + \phi_{L_{\psi_0}}(\ell),
    \phi_{L^{\nabla}_{\psi_0}}(\hat{\ell})\rangle_{\rm can}
    &= ~^{\dot{F}}\tilde{\mathbf{E}}_{\psi_0}(\pi \otimes 1)
    ~^{F} \langle t , \hat{\ell} \rangle_{\rm can},  \quad &\hat{t} \in T^{\nabla}_{\psi_0}, \;
    \hat{\ell} \in L^{\nabla}_{\psi_0}. 
      \end{aligned}
  \end{equation}
  This definition makes sense because
  $\phi_{T_{\psi_0}} \oplus \phi_{L_{\psi_0}}: T_{\psi_0} \oplus L_{\psi_0} \longrightarrow
  P_{\psi_0 \sigma}$
  is an $F$-linear isomorphism.
  For $\psi \neq \psi_0$ the map $\dot{F}: P_{\psi} \longrightarrow P_{\psi \sigma}$
  is an $F$-linear isomorphism. Therefore we can define
  $\dot{F}^{\nabla}: P^{\nabla}_{\psi} \longrightarrow P^{\nabla}_{\psi \sigma}$ by the
  equation
  \begin{displaymath}
    \langle\dot{F} z, \dot{F}^{\nabla} \hat{z}\rangle_{\rm can} =
    ~^{\dot{F}}\tilde{\mathbf{E}}_{\psi}(\pi \otimes 1)
    ~^{F}\langle z, \hat{z}\rangle_{\rm can}, \quad
    z \in P_{\psi}, \; \hat{z} \in P^{\nabla}_{\psi}. 
  \end{displaymath}
  We now apply Proposition \ref{LTD5p} to the modules
  $T^{\nabla}_{\psi_0}, L^{\nabla}_{\psi_0}$, and $P^{\nabla}_{\psi}$ for $\psi \neq \psi_0$,
  and to the maps $\phi_{T^{\nabla}_{\psi_0}}, \phi_{L^{\nabla}_{\psi_0}}$, and
  $\phi_{\psi} = \dot{F}^{\nabla}_{\psi}$ for $\psi \neq \psi_0$. This concludes
  the definition of the display
  $\mathcal{P}^{\nabla} = (P^{\nabla}, Q^{\nabla}, F^{\nabla}, \dot{F}^{\nabla})$.

  Now we verify  (\ref{LTD10e}). 
  Let  $\psi \neq \psi_0$.  If $y \in P_{\psi}$ and $\hat{y} \in P^{\nabla}_{\psi}$, the right hand side is by definition
  \begin{displaymath}
    ~^{\dot{F}}\tilde{\mathbf{E}}_{\psi}(\pi \otimes 1)
    \cdot \!\!\! ~^{F}\!(y, \hat{y})_{\rm can}. 
  \end{displaymath}
  Therefore (\ref{LTD10e}) holds in this case by the definition of
  $\dot{F}^{\nabla}_{\psi}$.
  For $\psi_0$ we use the normal decomposition (\ref{LTD12e}) and the induced
  normal decomposition
  $Q^{\nabla}_{\psi_0} = J_{\psi_0} T^{\nabla}_{\psi_0} \oplus L^{\nabla}_{\psi_0}$.
  Using the definition \eqref{LTD14e} of $J_{\psi_0}$, the identity (\ref{LTD10e}) becomes for the
  $\psi_0$-components a series of equations:
    \begin{enumerate}
    \item[(1)] $\langle\dot{F}\big((\pi \otimes 1 - 1 \otimes [\pi]) t\big),
      \dot{F}^{\nabla} \hat{\ell}\rangle_{\rm can} = \dot{F}_{\mathcal{L}}
      \langle(\pi \otimes 1 - 1 \otimes [\pi])t, \hat{\ell}\rangle_{\rm can},$ 
    \item[(2)] $\langle\dot{F}(~^{V}\eta t), \dot{F}^{\nabla} \hat{\ell}\rangle_{\rm can} =
      \dot{F}_{\mathcal{L}} \langle~^{V}\eta t, \hat{\ell}\rangle_{\rm can}, \quad \eta
      \in O \otimes_{O^t, \tilde{\psi}_0 \sigma} W(R),$
    \item[(3)] $\langle\dot{F} \ell,
      \dot{F}^{\nabla}\big((\pi \otimes 1 - 1 \otimes [\pi])\hat{t}\big)\rangle_{\rm can} =
      \dot{F}_{\mathcal{L}} \langle\ell, (\pi \otimes 1 - 1 \otimes [\pi])
      \hat{t}\rangle_{\rm can},$
\item[(4)] $\langle \dot{F} \ell, \dot{F}^{\nabla}(~^{V}\eta \hat{t})\rangle_{\rm can} =
  \dot{F}_{\mathcal{L}} \langle\ell, ~^{V}\eta \hat{t}\rangle_{\rm can},$
\item[(5)] $\langle\dot{F} \ell, \dot{F}^{\nabla} \hat{\ell}\rangle_{\rm can} = 0,$
\item[(6)] $\langle\dot{F} (J_{\psi_0} T_{\psi_0}),
  \dot{F}^{\nabla} (J_{\psi_0} T^{\nabla}_{\psi_0})\rangle_{\rm can} = 0.$
    \end{enumerate}
    We compute the right hand side of equation (1):
    \begin{displaymath}
      {\rm RHS}(1)= ^{\dot{F}}\!\tilde{\mathbf{E}}_{0}(\pi \otimes 1)
      (\pi \otimes 1 - 1 \otimes [\pi]) ~^{F}\langle t, \hat{\ell}\rangle_{\rm can} =
      ~^{\dot{F}}\tilde{\mathbf{E}}_{\psi_0}(\pi \otimes 1)
      ~^{F}\langle t, \hat{\ell}\rangle_{\rm can}.
    \end{displaymath}
    Therefore equation (1) is exactly the second equation of (\ref{LTD16e})
    for $\ell = 0$. The equation (3) follows in the same way. 

    We prove now the equation (2). For the right hand side we find:
    \begin{displaymath}
       {\rm RHS}(2)= \eta F_{\mathcal{L}} \langle t, \hat{\ell}\rangle_{\rm can} = \eta
      ~^{F}\tilde{\mathbf{E}}_{0}(\pi \otimes 1) ~^{F}\langle t,\hat{\ell}\rangle_{\rm can}. 
      \end{displaymath}
    We compute the left hand side of (2) by applying (\ref{LTD17e}) to
    $\dot{F}(~^{V}\eta t) = \eta Ft$:
      \begin{equation*}
      \begin{aligned}
      \langle \dot{F}(~^{V}\eta t), \dot{F}^{\nabla} \hat{\ell}\rangle_{\rm can} &=
      \eta ~^{\dot{F}}\tilde{\mathbf{E}}_{\psi_0}(\pi \otimes 1)^{-1}
      ~^{F}\tilde{\mathbf{E}}_{0}(\pi \otimes 1)
      \langle\dot{F}((\pi \otimes 1 - 1 \otimes [\pi])t),
      \dot{F}^{\nabla} \hat{\ell}\rangle_{\rm can} \\
&=\eta ~^{\dot{F}}\tilde{\mathbf{E}}_{\psi_0}(\pi \otimes 1)^{-1}
~^{F}\tilde{\mathbf{E}}_{0}(\pi \otimes 1)
~^{\dot{F}}\tilde{\mathbf{E}}_{\psi_0}(\pi \otimes 1) ~^{F}\langle t,\hat{\ell}\rangle_{\rm can}. 
      \end{aligned}
    \end{equation*}
    Here the last equation follows from (1). This proves (2). In the same
    way we obtain (4) from equation (3). The equation (5) follows from the second equation of
    (\ref{LTD16e}) for $t = 0$.

    Finally we prove equation (6). The special case
    \begin{equation}\label{LTD18e}
      \langle \dot{F} ((\pi \otimes 1 - 1 \otimes [\pi])t),
      \dot{F}^{\nabla} ((\pi \otimes 1 - 1 \otimes [\pi])\hat{t})\rangle_{\rm can} = 0
    \end{equation}
    is exactly the first equation of (\ref{LTD16e}) for $\ell = 0$.
    We have to show that the same holds if we replace the first argument
    of the bilinear form in (\ref{LTD18e}) by $\dot{F}(~^{V}\eta t)$ or
    the second argument by $\dot{F}^{\nabla}(~^{V}\eta \hat{t})$. But this
    may be reduced to (\ref{LTD18e}) in the same way as in the proof of
    equation (2). This finishes the proof of (\ref{LTD10e}). 
\end{proof}

\begin{proposition}\label{LTD9p}
  Let $\mathcal{P}_1$ and $\mathcal{P}_2$ be displays over $R$ with a strict
  $O$-action. Then the natural map
  \begin{displaymath}
    \Hom_{\text{$O${\rm-displays}}}(\mathcal{P}_2, \mathcal{P}^{\nabla}_{1}) \longrightarrow
    \Bil_{\text{$O${\rm-displays}}}(\mathcal{P}_1 \times \mathcal{P}_2, \mathcal{L}_R)
  \end{displaymath}
  is an isomorphism. Here we consider bilinear forms of displays with a strict
  $O$-action which are also $O$-bilinear. 
  \end{proposition}
\begin{proof}
  We define the inverse map. Let $\beta$ be an element from the right hand
  side. This is in particular a $O \otimes W(R)$-bilinear form
  \begin{displaymath}
\beta: P_1 \times P_2 \longrightarrow O \otimes W(R). 
  \end{displaymath}
  On the other hand, we have the canonical perfect $O \otimes W(R)$-bilinear
  form
$\langle\; , \;\rangle_{\rm can}: P_1 \times {P}^{\nabla}_{1} \longrightarrow O \otimes W(R)$. 
  Since this is perfect, we can define a $O \otimes W(R)$-module homomorphism
  $\alpha: P_2 \longrightarrow {P}^{\nabla}_{1}$ by
  \begin{displaymath}
\beta (x_1, x_2) = \langle x_1, \alpha(x_2)\rangle_{\rm can}. 
    \end{displaymath}
  We omit the straightforward verification that $\alpha$ defines a morphism of
  displays. 
\end{proof}

\begin{theorem}\label{LTD10p}
  Let $R$ be an $O$-algebra such that $p$ is nilpotent in $R$.
  Let $\mathcal{P}_1$ and $\mathcal{P}_2$ be displays over $R$ with a strict
  $O$-action. We denote by $\mathcal{P}_{1, {\rm a}}$ and $\mathcal{P}_{2, {\rm a}}$ their images
by the Ahsendorf functor $\mathfrak{A}_{O/\mathbb{Z}_p, R}$. 
 Proposition \ref{Adorf4c} and Proposition \ref{LTD8p} define a
homomorphism 
\begin{equation}\label{LTD20e}
  \Bil_{\text{$O${\rm-displays}}}(\mathcal{P}_1 \times \mathcal{P}_2, \mathcal{L}_R)
\longrightarrow  
  \Bil_{\text{$\CW_O(R)${\rm-displays}}}(\mathcal{P}_{1, {\rm a}} \times \mathcal{P}_{2, {\rm a}},
  \mathcal{P}_{m, \mathcal{W}_{O}(R)}(\pi^{ef}/p^f)) .
  \end{equation}

  If the displays $\mathcal{P}^{\nabla}_{1}$ and $\mathcal{P}_{2}$ are
  nilpotent, the homomorphism (\ref{LTD20e}) is an isomorphism. Equivalently,
  (\ref{LTD20e}) is an isomorphism if $(\mathcal{P}_{1, {\rm a}})^{\vee}$ and
  $\mathcal{P}_{2, {\rm a}}$ are nilpotent $\mathcal{W}_{O_F}(R)$-displays. 
\end{theorem}
\begin{proof}
  We apply (\ref{LTD20e}) to $\mathcal{P}_2 = \mathcal{P}_1^{\nabla}$ and
  the canonical bilinear form.  We obtain a bilinear form
  \begin{displaymath}
    \mathcal{P}_{1, {\rm a}} \otimes (\mathcal{P}_1^{\nabla})_{\rm a} \longrightarrow
    \mathcal{P}_{m, \mathcal{W}_{O} (O)}(\pi^{ef}/p^f) ,
  \end{displaymath}
  which is perfect by Proposition \ref{Adorf4c}. After twisting, we obtain also
  a perfect pairing of $\mathcal{P}_{1, {\rm a}}$ and
  $(\mathcal{P}_1^{\nabla})_{\rm a}((\pi^{ef}/p^f)^{-1})$ with values in
  $\mathcal{P}_{m, \mathcal{W}_{O} (O)}$. Therefore we have an identification with
  the dual display
  \begin{displaymath}
(\mathcal{P}_{1, {\rm a}})^{\vee} \cong (\mathcal{P}_1^{\nabla})_{\rm a}\big((\pi^{ef}/p^f)^{-1}\big). 
  \end{displaymath}
  By Theorem \ref{Adorf1t} we have an isomorphism
  \begin{displaymath}
    \Hom_{\text{$O${\rm-displays}}}(\mathcal{P}_2, \mathcal{P}_1^{\nabla}) \isoarrow
    \Hom_{\text{$\CW_O(R)${\rm-displays}}}(\mathcal{P}_{2, {\rm a}},
    (\mathcal{P}_1^{\nabla})_{\rm a}).  
  \end{displaymath}
  Here the left hand side agrees with the left hand side of (\ref{LTD20e})
  by Proposition \ref{LTD9p}.
  We have seen that the right hand side is the same as
  \begin{displaymath}
    \Hom_{\text{$\CW_O(R)${\rm-displays}}}(\mathcal{P}_{2, {\rm a}},
    \mathcal{P}_{1, {\rm a}}^{\wedge}(\pi^{ef}/p^f)) \cong
    \Bil_{\text{$\CW_O(R)${\rm-displays}}}(\mathcal{P}_{1, {\rm a}} \times \mathcal{P}_{2, {\rm a}},
  \mathcal{P}_{m, \mathcal{W}_{O}(R)}(\pi^{ef}/p^f)). 
  \end{displaymath}
  The last isomorphism follows from (\ref{dual4e}). 
  \end{proof}
\begin{definition}\label{LTD3d} 
  Let $R \in \Nilp_O$ and let $\mathcal{P}$ be a display with a strict
  $O$-action. A \emph{relative polarization} \index{relative polarization} of $\mathcal{P}$ with respect to $O$
  is a polarization of the $\mathcal{W}_O(R)$-display $\mathcal{P}_{\rm a}$ obtained
  from $\mathcal{P}$ by the Ahsendorf functor, cf. Definition \ref{def:pol}.   
\end{definition}

Let $\breve{O}$ be the completion of the maximal unramified extension of $O$. 
We consider  Theorem \ref{LTD10p} in the case of an $\breve{O}$-algebra $R$.
We denote by $\tau \in \Gal(\breve{O}/O)$ the Frobenius automorphism.
Since $\pi^e/p \in \breve{O}$ is a unit we find
$\eta_0 \in \breve{O}^{\times}$ such that
$\tau(\eta_0) \eta_0^{-1} = \pi^{e}/p$. By Lemma \ref{Adorf1l} there is a
$\tau-F$-equivariant homomorphism
$\breve{O} \longrightarrow W_{O}(\breve{O})$ such that the composite with
$\mathbf{w}_0$ is the identity.

Let $R \in \Nilp_{\breve{O}}$.
We denote by $\eta_{0,R}$ the image of $\eta_0$ by the homomorphism
\begin{displaymath}
  \breve{O} \longrightarrow W_{O}(\breve{O}) \longrightarrow W_{O}(R). 
  \end{displaymath}
Then  multiplication by $\eta_{0,R}^{f}$ defines an isomorphism of
$\mathcal{W}_{O}(R)$-displays,
\begin{equation}\label{P'Pol7e}
  \mathcal{P}_{m, \mathcal{W}_{O}(R)}(\pi^{ef}/p^f) \isoarrow
  \mathcal{P}_{m, \mathcal{W}_{O}(R)}. 
\end{equation}
Therefore we may write  Theorem  \ref{LTD10p}  without the twist $(\pi^{ef}/p^f)$.

\begin{corollary}\label{LTD9c}
 Let $R \in \Nilp_{\breve{O}}$. Fix $\eta_0\in \breve O^\times$ with $\tau(\eta_0) \eta_0^{-1} = \pi^{e}/p$, which defines the isomorphism \eqref{P'Pol7e}.  Let $\mathcal{P}$ be a display  with a
  strict $O$-action over $R$ such that $\mathcal{P}$ and $\mathcal{P}^{\nabla}$ are
  nilpotent.  Then a relative polarization on $\mathcal{P}$ with respect to $O$ is the
  same thing as an isogeny of displays with an $O$-action
  $\mathcal{P} \longrightarrow \mathcal{P}^{\nabla}$ such that the induced
  bilinear form 
\begin{displaymath}
\mathcal{P} \times \mathcal{P} \longrightarrow \mathcal{L}_R 
\end{displaymath}
is alternating. \qed
\end{corollary}
In the situation of the corollary,  $\mathcal{P}$ is the display of a formal
$p$-divisible group $X$ with a strict $O$-action and $\mathcal{P}^{\nabla}$
is the display of a formal $p$-divisible group with a strict $O$-action which
we denote by $X^{\nabla}$. We call $X^{\nabla}$ the \emph{Faltings dual}\index{Faltings dual} of $X$. However, we do not relate our definition to that of Faltings in \cite{F}, which operates directly in the realm of $p$-divisible groups. We can
consider a relative polarization as an isogeny of $p$-divisible groups with an $O$-action, 
\begin{equation}\label{LTD24e}
X \longrightarrow X^{\nabla} .
  \end{equation}

    \section{The contracting functor}\label{s:tcf} 

We return to the notation of section \ref{s:almostprinc}.  In particular,
throughout this section, $K/F$ denotes an etale extension of degree two of
a $p$-adic field $F$, and $r$ denotes a generalized CM-type of rank $2$.

Let $E$ be the reflex field of $r$, and let $\tilde{E} \subset \bar{\mathbb{Q}}_p$
be its normal closure. As in subsection \ref{ss:kotteis},  $E' \subset \tilde{E}$ is the composite
of $E$ and the normal closure of $K^t$ in $\bar{\mathbb{Q}}_p$. 
\subsection{The aim of this section}

\begin{definition}\label{KatCMpairs} Let $S$ be a scheme over $\Spec O_{E'}$ such that $p$ is locally nilpotent in $\CO_S$. 
  We denote\footnote{The symbol $\mathfrak P$ is to remind us that this is a category of local CM-pairs.} by $\mathfrak{P}_{r, S}$
  \index[NO]{PCA@$\mathfrak{P}_{r, S}$}  the category of local CM-pairs of type $r$ over $S$ which satisfy
  the Eisenstein conditions $({\rm EC}_r)$.
  If $S = \Spec R$, we will also write $\mathfrak P_{r, R}$ or simply
  $\mathfrak P_R$.
  
   The category of local CM-pairs $(\mathcal{P}, \iota)$
  of type $r$ in the sense of displays which satisfy the Eisenstein conditions will be denoted by
  $\mathfrak{d}\mathfrak P_{r,S}$\index[NO]{PCB@$\mathfrak{d}\mathfrak P_{r,S}$}, resp., $\mathfrak{d}\mathfrak P_{r,R}$. 
  \end{definition}

We will define a functor $\mathfrak{C}'_{r, R}$ that associates to a CM-pair
$(\mathcal{P}, \iota)\in \mathfrak{d}\mathfrak{P}_{r, R}$  a new display $(\mathcal{P}', \iota')$ over
$R$ with an action $\iota': O_K \longrightarrow \End \mathcal{P}'$. When  $r$ is special relative to $\varphi_0\colon F\to\ov\BQ_p$, cf. Definition \ref{def:special}, then 
the restriction of $\iota'$ to $O_F$ is strict with respect to
$O_F \overset{\varphi_0}{\longrightarrow} O_{E'} \longrightarrow R$. If $r$ is banal, then $\mathcal{P}'$ is \'etale. We will call the functor $\mathfrak{C}'_{r, R}$ the \emph{pre-contracting functor}. Under suitable hypotheses,  the pre-contracting functor will be an equivalence of categories. 

We will also  describe what
$\mathfrak{C}'_{r,R}$ does with polarizations, and define a functor $\mathfrak{C}^{\prime, {\rm pol}}_{r, R}$ on the category $\mathfrak P^{\rm pol}_{r,S}$, defined as follows. 
\begin{definition}\label{KatCMtriple1d} 
  We denote by  $\mathfrak{P}^{\rm pol}_{r,S}$
  \index[NO]{PCC@$\mathfrak{P}^{\rm pol}_{r,S}$}  the category of polarized local
  CM-triples $(X, \iota, \lambda)$ of type $r$ over $S$ such that
  $(X, \iota)$ satisfies the Eisenstein conditions $({\rm EC}_r)$.
  If $S = \Spec R$, we  also write $\mathfrak P^{\rm pol}_{r,R}$. We denote by
  $\mathfrak{d}\mathfrak{P}^{\rm pol}_{r,S}$
  \index[NO]{PCD@$\mathfrak{d}\mathfrak{P}^{\rm pol}_{r,S}$} the corresponding category of local
  CM-triples in the sense of displays.
Explicitly, $\mathfrak{d}\mathfrak{P}^{\rm pol}_{r,R}$ denotes the category of
triples $(\mathcal{P}, \iota, \beta)$ where
$(\mathcal{P}, \iota) \in \mathfrak{d}\mathfrak{P}_{r,R}$ and where 
$\beta\colon \mathcal{P} \times \mathcal{P} \longrightarrow \mathcal{P}_{m, \mathcal{W}(R)}$
is a polarization such that
\begin{displaymath}
  \beta(\iota(a) x_1, x_2) = \beta(x_1, \iota(\bar{a})x_2), \quad a \in O_K,
  \; x_1, x_2 \in P . 
\end{displaymath}
In the sequel, we will abbreviate  the last condition into saying that $\beta$ is \emph{anti-linear for the $O_K$-action}. \index{Anti-linear polarization form}
  \end{definition}
In a second step, we will use the functors $\mathfrak{C}^{\prime}_{r, R}$, resp. $\mathfrak{C}^{\prime, {\rm pol}}_{r, R}$, to define  \emph{contracting functors}. Here, we make a distinction between the case when $r$ is special and the case when $r$ is banal. In the case when $r$ is special, the image of the contracting functor is a $\CW_{O_F}(R)$-display $\CP_{\rm{c}}$, endowed with an action $\iota$ of $O_K$ such that $\CP_{\rm{c}}$ is of height $4$ and dimension $2$ (and, in the polarized case, with a polarization). When $r$ is banal, the image is a $p$-adic \'etale sheaf $G$ in $\BZ_p$-flat modules of rank $4d$, endowed with an action of $O_K$ (and with a polarization form in the polarized case). 
\begin{remark}
The pre-contracting functor is analogous to the functor in \cite[Thm. 4.12]{RZdrin}, with two important differences. First, in loc.~cit., there is no polarization in play, which makes the definition simpler. Second, in loc.~cit. the functor is only considered for schemes $S$ with $p\CO_S=0$. This is due to the fact that, in the context of \cite{RZdrin}, we were not able to handle the Kottwitz condition in the general case when $p$ is only locally nilpotent in $\CO_S$. 

A very similar pre-contracting functor also appears in \cite{M}.   Mihatsch considers the case where $K/F$ is an unramified field extension, and assumes that the  generalized local CM-type $r$ is \emph{unramified}. Let us explain this in the case when $r$ is of rank $2$ special relative to $\varphi_0$. Let $K^t$ be the maximal unramified subextension of $K$. With the notation of \S \ref{ss:kotteis}, we fix a disjoint decomposition $\Psi=\Psi_0\amalg\Psi_1$, where the summands are exchanged by the generator $\sigma\in\Gal(K^t/\BQ_p)$. Then $r_{\varphi_0}=r_{\bar\varphi_0}=1$ and for $\varphi\notin\{\varphi_0, \bar\varphi_0\}$
 \begin{equation}
    r_\varphi=
    \begin{cases}
      \begin{array}{ll}
        0, &\text{ if $\varphi_{|K^t}\in \Psi_0$}\\
         2, &\text{ if $\varphi_{|K^t}\in \Psi_1$}   .   \end{array}
    \end{cases}
  \end{equation} 
  In this case, the Kottwitz condition $({\rm KC}_r)$ and the Eisenstein conditions $({\rm EC}_r)$ can be replaced by the simpler conditions \cite[Def. 2.8]{M}, which makes the definition of the pre-contracting functor in \cite{M} easier than in the case of a general $r$. 
\end{remark}

\subsection{The Kottwitz and the Eisenstein condition for CM-pairs}

The Kottwitz condition $({\rm KC}_r)$ \index{Kottwitz condition} can be formulated in terms of
polynomial functions.  Let $\mathcal{L}$ be a locally free $R$-module
equipped with an action of $O_K$. If $S$ is an $R$-algebra, we write
    $\mathcal{L}_S = \mathcal{L} \otimes_R S$.

    Let us assume for a moment that $R$ is endowed with an
    $O_{\tilde{E}}$-algebra structure. For $\varphi \in \Phi =
    \Hom_{\text{$\mathbb{Q}_p$-Alg}}(K, \tilde{E})$ we 
    consider the induced map
    \begin{displaymath}
      \varphi_R: O_K  \overset{\varphi}{\longrightarrow} O_{\tilde{E}}  
      \longrightarrow R. 
    \end{displaymath}
    \begin{definition}\label{defKO}
      We say that $(\mathcal{L}, \iota)$ satisfies the Kottwitz condition
      relative to $r$ if for each $O_{\tilde{E}}$-algebra $S$ endowed with an
      $O_E$-algebra homomorphism $R \longrightarrow S$ 
      \begin{equation}\label{KottwitzC1e} 
 {\det}_S(a \mid \mathcal{L}_S) = \prod_{\varphi \in \Phi} \varphi_S(a)^{r_{\varphi}},
        \quad \text{for all}\; a \in O_K \otimes_{\mathbb{Z}_p} S. 
      \end{equation}
    \end{definition}
    Let $\mathbb{A} = \mathbb{V}(O_K)$ be the affine space over
    $\Spec \mathbb{Z}_p$. The right hand side of this equation may be
    considered as a polynomial function on
    $\mathbb{A}_{O_{\tilde{E}}}$. By base change to $\tilde{E}$, it is
    easily shown that this function is defined on $\mathbb{A}_{O_E}$.
    We note that each factor  of the right
    hand side of (\ref{KottwitzC1e}) is a linear polynomial function such
    that some coefficient is a unit in $O_{\tilde{E}}$. Therefore these
    factors are non-zero divisors in
    $\Gamma(\mathbb{A}_S, \mathcal{O}_{\mathbb{A}_S})$ for each $S$. Here we
    remark that a polynomial in $S[U_1, \ldots, U_r]$ is a non-zero
    divisor if one of its coefficients is a non-zero divisor in $S$.   

    Because the right hand side of (\ref{KottwitzC1e}) is a polynomial function
    on $\mathbb{A}_{O_E}$, the condition does not depend on the $O_{\tilde{E}}$-algebra
    structure on $S$. By a theorem of Amitsur,  condition (\ref{KottwitzC1e})
    is equivalent to the Kottwitz condition $({\rm KC}_r)$ of \eqref{signature.condition}
    (compare \cite{RZdrin}). 

    For a $O_{E'}$-algebra $S$ we have a decomposition 
    \begin{equation}\label{KottwitzC3e} 
      O_K \otimes_{\mathbb{Z}_p} S = \prod_{\psi \in \Psi}
      O_K \otimes_{O_{K^t}, \psi_S} S. 
    \end{equation}
    Here $\psi_S$ denotes the composite 
    $O_{K^t} \overset{\psi}{\longrightarrow} O_{E'} \longrightarrow S$. 
    Let $\mathbf{E}_{\psi_S}$ be the image of the Eisenstein polynomial
    $\mathbf{E} \in O_{K^t}[T]$ by the last homomorphism. We have a natural
    isomorphism
    \begin{equation}\label{KottwitzC2e} 
      S[T]/\mathbf{E}_{\psi_S} S[T] \isoarrow O_K \otimes_{O_{K^t}, \psi_S} S,
      \quad T \mapsto \Pi \otimes 1. 
    \end{equation}
    Therefore we may regard an $O_K \otimes_{O_{K^t}, \psi_S} S$-module $M$ as an
    $S[T]$-module. If $U \in S[T]$ and
    $x \in M$, we write $Ux = U(\Pi \otimes 1) x$.  If $U_0 \in O_{E'}[T]$, with image  $U \in S[T]$, then we write simply $U_0x = Ux$. 

    Returning to the $R$-module $\CL$ with action by $O_K$, the decomposition \eqref{KottwitzC3e} induces a decomposition
    \begin{displaymath}
      \mathcal{L}_S = \oplus \mathcal{L}_{S, \psi}. 
    \end{displaymath}

    By considering, for given $\psi$, an element $a$ of (\ref{KottwitzC3e})
    whose components are zero for $\psi' \neq \psi$, we obtain
    \begin{equation}\label{KottwitzC4e}
      {\det}_S(a \mid \mathcal{L}_{S, \psi}) = \prod_{\varphi| \psi}
      \varphi_S(a)^{r_{\varphi}}, 
      \quad \text{for all}\; a \in O_K \otimes_{O_{K^t}, \psi} S. 
    \end{equation}
    We call this condition $({\rm KC}_{\psi,r})$.
    \index[NO]{KAD@$({\rm KC}_{\psi,r})$} The condition
    $({\rm KC}_{r})$ holds iff the conditions $({\rm KC}_{\psi,r})$
    hold for each $\psi$.   

    We will call $\psi$ \emph{banal}\index{banal $\psi$} with respect to $r$ if $r_{\varphi} \in \{0, 2 \}$ for each
    $\varphi | \psi$. We call $\psi$ \emph{special} \index{special $\psi$} with respect to $r$ if there exists $\varphi | \psi$
    such that $r_{\varphi} = 1$ and if for any other $\varphi'| \psi$ with
    $r_{\varphi'} = 1$ we have $\varphi' = \bar{\varphi}$. We note that
    another $\varphi'$ can only exist if $\bar{\psi} = \psi$. 

    We also use the conditions $({\rm EC}_{\psi,r})$.
    \index[NO]{EAF@$({\rm EC}_{\psi,r})$}This means that we
    consider \eqref{Eisen} for a fixed $\psi$. 
        
    We consider CM-pairs $(X, \iota)$ of type $r$ over $R \in \Nilp_{O_E}$, cf. section \ref{ss:loctriples}.
    Thus $X$ is a $p$-divisible group of height $4d$ and dimension $2d$
    and $\iota$ is a $\mathbb{Z}_p$-algebra homomorphism
    \begin{displaymath}
      \iota: O_K \longrightarrow \End X 
    \end{displaymath}
    such the rank condition $({\rm RC}_r)$ is satisfied. 
    If we speak about the Kottwitz or Eisenstein condition we refer to the
    induced action on $\Lie X$. We use a similar terminology when we consider
    CM-pairs $(\mathcal{P}, \iota)$ in the sense of displays. This means that
    $\mathcal{P}=(P, Q, F, \dot{F})$ is a $\mathcal{W}(R)$-display of height
    $4d$ and dimension $2d$ endowed with a ring homomorphism
    \begin{displaymath}
      \iota: O_K \longrightarrow \End \mathcal{P} ,
    \end{displaymath}
    such that the rank condition $({\rm RC}_r)$ is satisfied for the
    induced action on $P/Q$.

   By Remark \ref{LaufunctorRm}, display theory provides a functor from the
    category of CM-pairs $(X, \iota)$ of type $r$ to the category of CM-pairs
    $(\mathcal{P}, \iota)$ of type $r$. We set
    $\mathbb{D}_{\mathcal{P}} = P/I_RP$ \index[NO]{DDA@$\mathbb{D}_{\mathcal{P}}$}
    and $\mathcal{L}_{\mathcal{P}} = P/Q$.
    \index[NO]{LBB@$\mathcal{L}_{\mathcal{P}}$}If
    $\mathcal{P}$ is the display of $X$, we have the identifications
    \begin{displaymath}
      \mathbb{D}_{\mathcal{P}} = \mathbb{D}_X(R), \quad
      \mathcal{L}_{\mathcal{P}} = \Lie X.
    \end{displaymath}
Here $\mathbb{D}_X(R)$ is the Grothendieck-Messing crystal evaluated
at $R$. For an $R$-algebra $S$, we will write
$\mathbb{D}_{\mathcal{P}, S} := \mathbb{D}_{\mathcal{P}} \otimes_{R} S$ and
$\mathcal{L}_{\mathcal{P}, S} := \mathcal{L}_{\mathcal{P}} \otimes_{R} S$.
If $\mathcal{P}$ is fixed, we write simply $\mathbb{D}_S$ and $\mathcal{L}_S$.
If $S$ is a $O_{E'}$-algebra, (\ref{KottwitzC3e}) gives a decomposition
\begin{displaymath}
 \mathbb{D}_S = \oplus_{\psi \in \Psi} \mathbb{D}_{S, \psi}. 
  \end{displaymath}
    
    \begin{proposition}\label{KottwitzC2p}
      Let $\psi$ be banal with respect to $r$. Let $(\mathcal{P},\iota)$ be a
      CM-pair of type $r$ over an $O_E$-algebra $R$.

      Then the Eisenstein condition  $({\rm EC}_{\psi,r})$ is satisfied iff
      $\mathbf{E}_{A_{\psi}} \mathbb{D}_{\psi}$ is the kernel of the canonical
      map $\mathbb{D}_{\psi} \longrightarrow \mathcal{L}_{\mathcal{P}, \psi}$.
      Moreover $({\rm EC}_{\psi,r})$ implies the Kottwitz condition
      $({\rm KC}_{\psi,r})$. 
      
      {\rm  Here $\mathbf{E}_{A_{\psi}} $ denotes the operator $\mathbf{E}_{A_{\psi}} \big(\iota(\Pi)\big)$ on the module in question, for a fixed choice of $\Pi$, cf. \eqref{Eisen}.}
    \end{proposition}
    \begin{proof} 
      We reduce to the case where $S$ is an $R$-algebra endowed with an $O_{\tilde{E}}$-algebra structure.
      Then $\mathbf{E}_{A_{\psi, S}} \in S[T]$ is defined as the image of
      $\mathbf{E}_{A_{\psi}}$ by $O_{\tilde{E}}[T] \longrightarrow S[T]$. It acts
      on any $O_K \otimes_{O_{K^t}, \psi} S$-module by (\ref{KottwitzC2e}). 
      
      Via  $\iota$, we view $\mathbb{D}_{S}$ and $\mathcal{L}_{S}$ as $O_K
      \otimes_{\mathbb{Z}_p} S$-modules.  We
      consider the canonical surjective map
      \begin{displaymath}
        \mathbb{D}_{S} \longrightarrow \mathcal{L}_{S}.
      \end{displaymath}
      The decomposition (\ref{KottwitzC3e}) induces decompositions,
      \begin{equation}\label{KottwitzC14e}
        \mathbb{D}_S = \oplus_{\psi} \mathbb{D}_{\psi}, \quad
        \mathcal{L}_S = \oplus_{\psi} \mathcal{L}_{\psi}.
      \end{equation} 
      We allowed ourselves to omit the index $S$ on the right hand side of
      these equations.

      The Eisenstein condition $({\rm EC}_{\psi,r})$ for banal $\psi$ says
      that $\mathcal{L}_{\psi}$ is annihilated by
      $\mathbf{E}_{A_{\psi}}$, cf. \eqref{Eisenbanal}. Therefore it is clearly
      implied by the condition of the Proposition. If conversely
      $({\rm EC}_{\psi,r})$ holds, we obtain a surjective map
      \begin{equation}\label{KottwitzC6e}
        \mathbb{D}_{\psi}/\mathbf{E}_{A_{\psi}}\mathbb{D}_{\psi}
        \longrightarrow \mathcal{L}_{\psi}.
      \end{equation} 
By  Lemma \ref{Displaykristall1l},  $\mathbb{D}_{\psi}$ is locally on
$\Spec S$ a free $O_K \otimes_{O_{K^t}, \psi} S$-module. It has rank $2$ because
the height of $\mathcal{P}$ is $4d$. We may assume that
      \begin{displaymath}
        \mathbb{D}_{\psi} \cong (O_K \otimes_{O_{K^t}, \psi} S)^{2} =
        S[T]^2/\mathbf{E}_{\psi_S}S[T]^2.
      \end{displaymath}
We see that both sides of (\ref{KottwitzC6e}) are locally free $S$-modules
of the same rank $r_{\psi} = \sum_{\varphi|\psi} r_{\varphi}$. Therefore this map is
an isomorphism.
    
      The condition $({\rm KC}_{\psi,r})$ would follow from 
      \begin{equation}\label{KottwitzC7e}
        {\det}_S(a \mid S[T]/ \mathbf{E}_{A_{\psi}} S[T]) = \prod_{\varphi \in
          A_{\psi}} \varphi_S(a), \quad a \in S[T]. 
      \end{equation}
      We have
      \begin{displaymath}
        \mathbf{E}_{A_{\psi}}(T) = \prod_{\varphi \in A_{\psi}} (T -
        \varphi_S(\Pi \otimes 1)), \quad \varphi_S(a) = a(\varphi_S(\Pi \otimes 1)). 
      \end{displaymath}
       We obtain (\ref{KottwitzC7e}) from the following Lemma.
    \end{proof}

    \begin{lemma}
      Let $R$ be a ring. Let
      \begin{displaymath}
        \mathbf{E}(T) = \prod_{i=1}^s (T - \Pi_i), \quad \Pi_i \in R
      \end{displaymath}
      be a polynomial. A polynomial $f(T) \in R[T]$ defines by
      multiplication an endomorphism of the free $R$-module
      $R[T]/\mathbf{E}R[T]$. Then
      \begin{displaymath}
        {\det}_{R} (f(T) \mid R[T]/\mathbf{E}R[T]) = \prod_{i=1}^{s} f(\Pi_i). 
      \end{displaymath}
    \end{lemma}
    \begin{proof}
      One can easily reduces the question to the case where $R$ is a field of
      characteristic $0$ such that $\mathbf{E}(T)$ is a product of
      different linear factors. For the reduction one starts with a ring
      homomorphism 
      \begin{displaymath}
        \mathbb{Z}[\underline{X}, \underline{Y}] \longrightarrow R ,
      \end{displaymath}
      where for the first set of variables $\underline{X} = (X_1, \ldots, X_s)$,
      $X_i$ is mapped to $\Pi_i$ and where the second set of variables $\underline{Y}$ is
      mapped to the coefficients of $f$.

      If now $R$ is a field and $\mathbf{E}$ is separable, we have a canonical
      isomorphism of $R$-modules
      \begin{displaymath}
        \begin{aligned}
          R[T]/\mathbf{E}R[T] & \isoarrow & \oplus_{i=1}^s R .\\
          T              & \mapsto     &   (\Pi_i)_i
        \end{aligned}
      \end{displaymath}
      Multiplication by $f(T)$ on the left hand side acts on the right
      hand side on the $i$-th factor by multiplication by $f(\Pi_i)$. This
      implies the assertion. 
    \end{proof}
\begin{corollary}
Let $r$ be banal. Let $(\CP, \iota)$ be a CM-pair of type $r$ over an $O_E$-algebra $R$. Then the Eisenstein condition $({\rm EC}_r)$ implies the Kottwitz condition $({\rm KC}_r)$.\qed
\end{corollary}
    We consider next the case where $\psi$ is special. This means by
    definition that there is exactly one pair
    $\{\varphi, \bar{\varphi}\}$ such that $\varphi | \psi$ and
    $r_{\varphi} = r_{\bar{\varphi}} = 1$. When $K/F$ is ramified,
    we have $E' = E$ and when $K/F$ is  unramified,   we have $[E':E] \leq 2$. 
    \begin{proposition}\label{KottwitzC4p}
      Let $\psi$ be special with respect to $r$. Let $R$ be a $O_{E}$-algebra such that $p$ is nilpotent in $R$.   
      Let $(\mathcal{P},\iota)$ be a CM-pair of type $r$ over $R$ which
      satifies the Eisenstein condition $({\rm EC}_{\psi,r})$. Let $S$ be a
      $O_{E'}$-algebra which is endowed with an $O_E$-algebra homomorphism
      $R \longrightarrow S$. Then, with the notations of (\ref{KottwitzC14e}), the canonical map
      \begin{displaymath}
        \mathbb{D}_{\psi}/\mathbf{E}_{A_{\psi}} \mathbb{D}_{\psi}
        \longrightarrow \mathcal{L}_{\psi}/\mathbf{E}_{A_{\psi}} \mathcal{L}_{\psi} 
      \end{displaymath}
      is an isomorphism.
    \end{proposition}
    \begin{proof}
      Clearly we may assume that $S$ is a local ring with residue class
      field $k$. We postpone the verification that the assertion holds for
      $S = k$ (compare (\ref{KottwitzC22e}) and (\ref{KottwitzC23e}) below).   

      We begin with the case $K/F$ unramified. Then
      $\rank_S \mathcal{L}_{\psi} = r_{\psi} = 2a_{\psi} + 1$. Let
      \begin{equation}\label{KottwitzC8e}
        f: \mathcal{L}_{\psi} \longrightarrow \mathcal{L}_{\psi}
      \end{equation}
      be the $S$-module homomorphism given by multiplication with
      $\mathbf{E}_{A_{\psi}}$. From the case of 
      a field, we deduce that
      $\dim_k \mathcal{L}/f(\mathcal{L}) \otimes_S k = 2a_{\psi}$. By
      $({\rm EC}_{\psi, r})$ we have
      \begin{displaymath}
        \bigwedge^2 f = 0 ,
      \end{displaymath}
      cf. \eqref{furtherDr2unram}. 
      Therefore we can apply Lemma 4.9 of \cite{RZdrin} with $s = 1$.
      This says that $\mathcal{L}_\psi/f(\mathcal{L}_\psi)$ is a free $S$-module of
      rank $2a_{\psi}$. Therefore the canonical map of the proposition is a
      surjection of free $S$-modules of the same rank, and hence an
      isomorphism.

      The argument in the case $K/F$ ramified is similiar. In this case, we
      have $\rank_S \mathcal{L}_{\psi} = r_{\psi} = 2a_{\psi} + 2 = 2e$. The
      dimension $\dim_k \mathcal{L}_\psi/f(\mathcal{L}_\psi) \otimes_S k = 2a_{\psi}$, as before. In this case the condition $({\rm EC}_{\psi, r})$ says
      \begin{displaymath}
        \bigwedge^3 f = 0 ,
      \end{displaymath}
      cf. \eqref{furtherDr2ram}. Therefore we may apply Lemma 4.9 loc.cit. with $s = 2$. We conclude as
      before. 
    \end{proof}

    \begin{proposition}\label{KottwitzC5p}
      Let  $r$ be special and $K/F$  unramified.
      Let $R$ be a $O_{E}$-algebra such that $p$ is nilpotent in $R$. Let
      $(\mathcal{P},\iota)$ be a CM-pair of type $r$ over  $R$ which satifies
      the Eisenstein condition $({\rm EC}_r)$.  Then the condition $({\rm KC}_r)$ is also satisfied. 
    \end{proposition}
    \begin{proof}
      We consider an algebra $S$ as in the last proposition. We keep the
      notation of (\ref{KottwitzC14e}). We need only to
      verify $({\rm KC}_{\psi_0, r})$ since the Kottwitz condition is
      satisfied for $\psi$ banal by Proposition \ref{KottwitzC2p}.  We have to
      verify that
      \begin{equation}\label{KottwitzC12e}
        {\det}_S(a \mid \mathcal{L}_{\psi_0}) = \varphi_{0,S}(a) \cdot
        \prod_{\varphi \in A_{\psi_0}} \varphi_S(a)^{2}, \quad a \in
        O_K \otimes_{O_{K^t}, \psi_0} S. 
      \end{equation}
      Since $\mathbb{D}_{\psi_0}$ is locally on $S$ a free
      $O_K \otimes_{O_{K^t}, \psi_0} S$-module of rank $2$, we obtain from the
      isomorphism of Proposition \ref{KottwitzC4p}
      \begin{displaymath}
        {\det}_S (a \mid \mathcal{L}_{\psi_0}/\mathbf{E}_{A_{\psi_0}}
        \mathcal{L}_{\psi_0}) = \prod_{\varphi \in A_{\psi_0}} \varphi_S(a)^{2}.
      \end{displaymath}
      The proposition also shows that
      $\mathbf{E}_{A_{\psi_0}} \mathcal{L}_{\psi_0}$ is a locally free $S$-module
      of rank $1$. It follows from the Eisenstein condition  \eqref{furtherDr2unram} that this module
      is annihilated by $(T - \varphi_{0, S}(\Pi \otimes 1))$. Hence  an element
      $a \in O_K \otimes_{O_{K^t}, \psi_0} S = 
      S[T]/\mathbf{E}_{\psi_0}S[T]$
      acts on $\mathbf{E}_{A_{\psi_0}} \mathcal{L}_{\psi_0}$ as
      $\varphi_{0, S}(a)$. In particular
      \begin{displaymath}
        {\det}_S (a \mid \mathbf{E}_{A_{\psi_0}} \mathcal{L}_{\psi_0}) =
        \varphi_{0, S}(a). 
      \end{displaymath}
      The formula (\ref{KottwitzC12e}) follows. 
    \end{proof}
    We reformulate the Eisenstein condition in the case where $K/F$ is
    unramified.
    \begin{proposition}\label{KottwitzC8p} 
  Let  $r$ be special and $K/F$  unramified. Let $\varphi_0, \bar{\varphi}_0 \in \Phi$ be the two
  embeddings such that $r_{\varphi_0} = r_{\bar{\varphi}_0} = 1$, and let $\psi_0$,
  resp. $\bar{\psi_0}$, the embeddings induced by $\varphi_0$,
  resp. $\bar{\varphi}_0$. Let $R$ be a $O_{E'}$-algebra, and let $(\mathcal{P}, \iota)$ be a CM-pair of
  type $r$ over $R$. Let $\mathbb{D} = \mathbb{D}_{\mathcal{P}}(R)$. The CM-pair $(\mathcal{P}, \iota)$ satisfies the Eisenstein conditions iff
  the following conditions hold.
  \begin{enumerate}
  \item[(1)] If $\psi \in \Psi$ is banal,  then the canonical map
    \begin{displaymath}
      \mathbb{D}_{\psi}/\mathbf{E}_{A_{\psi}}\mathbb{D}_{\psi} \longrightarrow
      \Lie_{\psi} X
    \end{displaymath}
    is an isomorphism.
  \item[(2)] If $\psi \in \{\psi_0, \bar{\psi_0}\}$, then the canonical map
    \begin{displaymath}
      \mathbb{D}_{\psi}/\mathbf{E}_{A_{\psi}}\mathbb{D}_{\psi} \longrightarrow
      \Lie_{\psi} X/ \mathbf{E}_{A_{\psi}} \Lie_{\psi} X
    \end{displaymath}
    is an isomorphism.
  \item[(3)] The $R$-modules $\mathbf{E}_{A_{\psi_{0}}} \Lie_{\psi_0}$,
    resp. $\mathbf{E}_{A_{\bar{\psi}_{0}}} \Lie_{\bar{\psi}_0}$, are locally
    free of rank $1$ and $O_K$ acts on them via
    \begin{displaymath}
      \varphi_{0,R}: O_K \longrightarrow O_{E'} \longrightarrow R, \quad \text{resp.}
      \quad \bar{\varphi}_{0,R}: O_K \longrightarrow O_{E'} \longrightarrow R. 
      \end{displaymath}
    \end{enumerate}
  \end{proposition} 
\begin{proof}
  This is a consequence of Proposition \ref{KottwitzC4p} and the proof of
  Proposition \ref{KottwitzC5p}.  
\end{proof}

We next consider what happens to the Eisenstein conditions when passing from a display to its  conjugate dual, cf.  Lemma \ref{conjdual}. We note that we already checked in loc.~cit. that the  condition $({\rm RC}_r)$ is preserved. Recall that, if 
$(\mathcal{P}, \iota)$ is a CM-pair,  the conjugate dual
$(\mathcal{P}^{\wedge}, \iota^{\wedge})$ is defined by
\begin{displaymath}
  \big(\mathcal{P}^{\wedge} = \mathcal{P}^{\vee}, \quad \iota^{\wedge}(a) =
  \iota(\bar{a})^{\wedge}\big). 
  \end{displaymath}
\begin{corollary}\label{dualEisen2c} Let $K/F$ be unramified and let $r$ be
  arbitrary or let $K/F$ be split.  
  Let $(\mathcal{P}, \iota)$ be a CM-pair over an $O_E$-algebra $R$ such that
  $p$ is nilpotent in $R$. If $(\mathcal{P}, \iota)$ satisfies the Eisenstein
  condition $({\rm EC}_r)$, then the conjugate dual
  $(\mathcal{P}^{\wedge}, \iota^{\wedge})$ also satisfies $({\rm EC}_r)$. 
    \end{corollary}
 \begin{proof}
  We have a canonical
  isomorphism $\mathbb{D}^{\wedge} = \Hom_R(\mathbb{D}, R)$. The resulting
  perfect pairing
  \begin{equation}\label{dualEisen3e}
\langle\; , \; \rangle : \mathbb{D} \times \mathbb{D}^{\wedge} \longrightarrow R 
  \end{equation}
  satisfies
  \begin{equation}\label{dualEisen2e}
    \langle ax , \hat{x} \rangle = \langle x , \bar{a}\hat{x} \rangle, \quad
    a \in O_K, \; x \in \mathbb{D}, \; \hat{x} \in \mathbb{D}^{\wedge}. 
  \end{equation}
  This implies that for $\psi_1 \neq \bar{\psi}_2$ the modules
  $\mathbb{D}_{\psi_1}$ and $\mathbb{D}^{\wedge}_{\psi_2}$ are orthogonal and
  that for any $\psi$ the induced pairing
  \begin{equation}\label{dualEisen1e}
\mathbb{D}_{\psi} \times \mathbb{D}^{\wedge}_{\bar{\psi}} \longrightarrow R 
  \end{equation}
  is perfect. 
  Let $\mathbb{D}^{1}_{\psi} \subset \mathbb{D}_{\psi}$ be the kernel of the map
  $\mathbb{D}_{\psi} \longrightarrow \mathcal{L}_{\psi} := \Lie_{\psi} \mathcal{P}$
  and let
  $\mathbb{D}^{\wedge,1}_{\bar{\psi}} \subset \mathbb{D}^{\wedge}_{\bar{\psi}}$
  be defined in the same way. By definition of the dual display, 
  $\mathbb{D}^{1}_{\psi}$ and $\mathbb{D}^{\wedge,1}_{\bar{\psi}}$ are orthogonal
  complements of each other with respect to (\ref{dualEisen1e}).
  We consider the case $\psi \in \{ \psi_0, \bar{\psi}_0 \}$. Recall that
  this is possible only in the non-split case. The
  Eisenstein condition for $\mathcal{P}$ says that we have a split filtration
  of direct summands of $\mathbb{D}_{\psi}$ 
  \begin{equation}\label{dualEisen5e}
    \mathbf{S}_{\psi} \mathbf{E}_{A_{\psi}} \mathbb{D}_{\psi}\subset 
    \mathbb{D}^{1}_{\psi} \subset \mathbf{E}_{A_{\psi}} \mathbb{D}_{\psi},
  \end{equation}
  such that the factor modules are locally free of rank $1$. 
   We claim that the orthogonal complement of
  $\mathbf{E}_{A_{\psi}} \mathbb{D}_{\psi}$ is
$\mathbf{S}_{\bar{\psi}} \mathbf{E}_{A_{\bar{\psi}}} \mathbb{D}^{\wedge}_{\bar{\psi}}$. 
  Indeed, by (\ref{dualEisen2e}), we have
  \begin{equation}\label{dualEisen6e}
    \langle\; \mathbf{E}_{A_{\psi}} x , \; \hat{x} \rangle =
    \langle\; x , \; \mathbf{E}_{B_{\bar{\psi}}} \hat{x} \rangle. 
  \end{equation}
  This implies that $\mathbf{E}_{A_{\psi}} \mathbb{D}_{\psi}$ and
  $\mathbf{S}_{\bar{\psi}} \mathbf{E}_{A_{\bar{\psi}}} \mathbb{D}^{\wedge}_{\bar{\psi}}$
  are orthogonal. Because (\ref{dualEisen3e}) is perfect we obtain a surjection
  of $R$-modules
  \begin{equation}\label{dualEisen4e}
    \mathbb{D}_{\psi}/ \mathbf{E}_{A_{\psi}} \mathbb{D}_{\psi} \longrightarrow
    \Hom_R(\mathbf{S}_{\bar{\psi}} \mathbf{E}_{A_{\bar{\psi}_R}}
    \mathbb{D}^{\wedge}_{\bar{\psi}}, R). 
  \end{equation}
  Recall that $\mathbb{D}_{\psi}$ is locally on $\Spec R$ a free
  $O_K \otimes_{O_{K^t}, \psi} R$-module of rank $2$. It follows that both sides
  of (\ref{dualEisen4e}) are locally free $R$-modules of the same rank
  $2a_{\psi} = 2e - 2a_{\bar{\psi}} - 2$. Therefore this map is an isomorphism.
  This proves our claim about the orthogonal complement. By the same argument,
  $\mathbf{E}_{A_{\bar{\psi}_R}} \mathbb{D}^{\wedge}_{\bar{\psi}}$ is the orthogonal
  complement of $\mathbf{S}_{\psi} \mathbf{E}_{A_{\psi}} \mathbb{D}_{\psi}$. 

  We take the orthogonal complement of (\ref{dualEisen5e}) and obtain the
  filtration
  \begin{displaymath}
\mathbf{S}_{\bar{\psi}} \mathbf{E}_{A_{\bar{\psi}}} \mathbb{D}^{\wedge}_{\bar{\psi}}
\subset \mathbb{D}^{\wedge,1}_{\bar{\psi}} \subset
\mathbf{E}_{A_{\bar{\psi}}} \mathbb{D}^{\wedge}_{\bar{\psi}} ,
  \end{displaymath}
  and conclude that the factor modules are locally free of rank $1$.

  Now let $\psi$ be banal. We have to prove that
  $\mathbf{E}_{A_{\psi}} \mathbb{D}^{\wedge}_{\psi}\subset\mathbb{D}^{\wedge,1}_{\psi}$.
  The right hand side is the orthogonal complement of
  $\mathbb{D}^1_{\bar{\psi}} = \mathbf{E}_{A_{\bar{\psi}}} \mathbb{D}_{\bar{\psi}}$.
  Therefore we have to prove that
  \begin{displaymath}
    \langle\; \mathbf{E}_{A_{\bar{\psi}}} x , \; \mathbf{E}_{A_{\psi}}\hat{x} \rangle
    = 0, \quad \text{for}\; x \in \mathbb{D}_{\bar{\psi}}, \;
    \hat{x} \in \mathbb{D}^{\wedge}_{\psi}. 
  \end{displaymath}
  Using (\ref{dualEisen6e}), we find for the right hand side
  \begin{displaymath}
    \langle\; x , \; \mathbf{E}_{B_{\psi}}\mathbf{E}_{A_{\psi}}\hat{x} \rangle = 
         \langle\; x , \; \mathbf{E}_{\psi}\hat{x} \rangle = 0. 
    \end{displaymath}
\end{proof}

Before proving the analogue of Corollary \ref{dualEisen2c} in the case when $K/F$ is ramified, we further analyze in this case the Eisenstein conditions. 
\begin{proposition}\label{KottwitzC6p}
Let $r$ be special and  $K/F$  ramified. Let $R$ be a $O_{E}$-algebra
such that $p$ is nilpotent in $R$. Let $(\mathcal{P},\iota)$ be a CM-pair of
type $r$ over  $R$. Since $E' = E$, the decomposition (\ref{KottwitzC3e}) is defined for $S = R$.
Then the Eisenstein condition $({\rm EC}_{\psi_0,r})$ holds iff the
following conditions are satisfied.
\begin{enumerate}
\item[(1)] The $R$-module
  $\mathbf{E}_{A_{\psi_0}} \mathcal{L}_{\psi_0} \subset \mathcal{L}_{\psi_0}$ 
  is a direct summand which is locally free of rank $2$.
\item[(2)] The action of $\iota(\pi)$ on
  $\mathbf{E}_{A_{\psi_0}} \mathcal{L}_{\psi_0}$ coincides with the action of
  $\varphi_0(\pi) \in R$, i.e., the action of the image of $\pi$ by the homomorphism
  $O_F \overset{\varphi_0}{\longrightarrow} O_E \longrightarrow R$.   
  \end{enumerate}

Furthermore, a CM-pair $(\mathcal{P}, \iota)$ which satisfies $({\rm EC}_{r})$ also 
satisfies the Kottwitz condition $({\rm KC}_{r})$ if and only if 
      \begin{equation}\label{KottwitzC15e} 
      \Trace_R (\iota(\Pi) \mid \mathbf{E}_{A_{\psi_0}} \mathcal{L}_{\psi_0}) = 0. 
      \end{equation}
    \end{proposition}
    \begin{proof} 
      For the proof we may pass from $R$ to an $R$-algebra $S$ which is endowed
      with a $O_{\tilde{E}}$-algebra stucture. We continue with the notations of
      (\ref{KottwitzC14e}). The first assertion of the proposition is then an 
      immediate consequence of Proposition \ref{KottwitzC4p}.

      To verify the last sentence on the Proposition it suffices by
      Proposition \ref{KottwitzC2p} to consider $({\rm KC}_{\psi_0,r})$.
      This condition reads 
      \begin{displaymath}
        {\det}_S(a \mid \mathcal{L}_{\psi_0}) = \varphi_{0,S}(a) \cdot 
        \bar{\varphi}_{0,S}(a) \cdot \prod_{\varphi \in A_{\psi_0}}
        \varphi_S(a)^{2}, \quad a \in O_K \otimes_{O_{K^t}, \psi_0} S.
      \end{displaymath}
      In this case $\mathbf{E}_{A_{\psi_0}} \mathcal{L}_{\psi_0}$ is a locally
      free $S$-module of rank $2$. By Proposition \ref{KottwitzC4p}, it is
      enough to show
      \begin{equation}\label{KottwitzC13e}
        {\det}_S(a \mid \mathbf{E}_{A_{\psi_0}} \mathcal{L}_{\psi_0}) =
        \varphi_{0,S}(a) \cdot \bar{\varphi}_{0,S}(a) .
      \end{equation} 
      In this case, the Eisenstein condition says that
      \begin{equation}\label{KottwitzC11e}
        (T - \varphi_{0, S}(\Pi \otimes 1))(T - \bar{\varphi}_{0, S}(\Pi
        \otimes 1)) = T^2 + \psi_0(\pi) 
      \end{equation}
      annihilates $\mathbf{E}_{A_{\psi_0}} \mathcal{L}_{\psi_0}$, cf \eqref{furtherDr2ram}. Note
      that the action of $T^2$ on
      $\mathbf{E}_{A_{\psi_0}} \mathcal{L}_{\psi_0}$ is by definition
      the action of $\iota(\Pi^2) = -\iota(\pi)$. Therefore the action of
      $O_F$ on $\mathbf{E}_{A_{\psi_0}} \mathcal{L}_{\psi_0}$ via
      $\iota$ coincides with the action via $O_F
      \overset{\varphi_0}{\longrightarrow} O_{\tilde{E}} \longrightarrow  S$.
      The polynomial (\ref{KottwitzC11e}) is the
      characteristic polynomial of $\iota(\Pi)$ acting on the $S$-module 
      $\mathbf{E}_{A_{\psi_0}} \mathcal{L}_{\psi_0}$. This follows from
      the trace condition of the proposition. Therefore the desired equation
      (\ref{KottwitzC13e}) is a consequence of Lemma \ref{lemdet} below.

      Conversely, assume that the Kottwitz condition holds. By Proposition
      \ref{KottwitzC4p}, this  implies
      \begin{displaymath}
        {\det}_S(a \mid \mathbf{E}_{A_{\psi_0}} \mathcal{L}_{\psi_0}) \cdot
        \prod_{\varphi \in A_{\psi_0}} \varphi_S(a)^{2} = 
        \varphi_{0,S}(a) \cdot \bar{\varphi}_{0,S}(a) \cdot
        \prod_{\varphi \in A_{\psi_0}} \varphi_S(a)^{2} , \quad a\in O_K\otimes S . 
      \end{displaymath}
      We already remarked right after Definition \ref{defKO} that $\varphi_S$ is a non-zero divisor in the ring of
      polynomial functions. Therefore we conclude
      \begin{displaymath}
        {\det}_S(a \mid \mathbf{E}_{A_{\psi_0}} \mathcal{L}_{\psi_0}) = 
        \varphi_{0,S}(a) \cdot \bar{\varphi}_{0,S}(a) , \quad \text{ for all } a\in O_K\otimes S . 
      \end{displaymath}
This implies that the characteristic polynomial of $\iota(\Pi)$ acting
on  $\mathbf{E}_{A_{\psi_0}} \mathcal{L}_{\psi_0}$ is the
polynomial (\ref{KottwitzC11e}). Therefore the trace is $0$. 
    \end{proof}
We state the needed Lemma without proof.
      \begin{lemma}\label{lemdet}
      Let $S$ be a ring. Let $L$ be a locally free $S$-module of rank
      $2$. Let $\alpha: L \longrightarrow L$ be an endomorphism with
      characteristic polynomial
      \begin{displaymath}
        {\det}_S(T\id_L - \alpha \mid L) = T^2 - s_1T + s_2.
      \end{displaymath}
Then for all $\lambda, \mu \in S$
      \begin{displaymath}
        {\det}_S(\mu \id_L - \lambda \alpha) = \mu^2 - s_1\mu \lambda +
        s_2\lambda^2. 
      \end{displaymath}
      Assume that the characteristic polynomial splits
      \begin{displaymath}
        T^2 - s_1T + s_2 = (T- \rho_1)(T - \rho_2) .
      \end{displaymath}
      Consider $L$ as $S[T]$-module, and let $\phi_i: S[T] \longrightarrow S$ be the
      $S$-algebra homomorphism such that $\phi_i(T) = \rho_i$. Then for each
      polynomial $a \in S[T]$
      \begin{displaymath}
        {\det}_S(a \mid L) = \phi_1(a) \cdot \phi_2(a). 
      \end{displaymath}\qed
    \end{lemma}

    \begin{remark}
      Let $A \subset B$ be $R$-modules. Then we write $A \overset{c}{\subset} B$
      \index[NO]{AAC@$A \overset{c}{\subset} B$}
if the factor module $B/A$ is a finitely generated projective $R$-module
of rank $c$.

Let $(\mathcal{P}, \iota)$ as in Proposition \ref{KottwitzC6p} such that
$({\rm EC}_{\psi_0,r})$ is satisfied. We write $\mathbb{D} = P/I_RP$. 
Let $\bar{Q}_{\psi_0}$ the kernel of
$\mathbb{D}_{\psi_0} \longrightarrow \mathcal{L}_{\psi_0}$. By Lemma
\ref{Displaykristall1l},  $\mathbb{D}_{\psi_0}$ is a free
$O_K \otimes_{O_{F^t}, \psi_0} R$-module  of rank $2$. We obtain that 
\begin{displaymath}
  \mathbf{E}_{A_{\psi_0}}\mathbb{D}_{\psi_0} \overset{2(e-1)}{\subset}
  \mathbb{D}_{\psi_0} .
\end{displaymath}
On the other hand, the condition $(1)$ of  Proposition \ref{KottwitzC6p} says
\begin{displaymath}
  \bar{Q}_{\psi_0} \overset{2}{\subset} \mathbf{E}_{A_{\psi_0}}\mathbb{D}_{\psi_0}
  + \bar{Q}_{\psi_0} \overset{2(e-1)}{\subset} \mathbb{D}_{\psi_0}.
\end{displaymath}
This implies $\bar{Q}_{\psi_0} \subset \mathbf{E}_{A_{\psi_0}}\mathbb{D}_{\psi_0}$.
Therefore we may reformulate the two conditions in  Proposition \ref{KottwitzC6p} in one line:
\begin{equation}\label{KottwitzC36enew}
  \mathbf{S}_{\psi_0} \mathbf{E}_{A_{\psi_0}} \mathbb{D}_{\psi_0} 
  \overset{2}{\subset} \bar{Q}_{\psi_0} \overset{2}{\subset}
  \mathbf{E}_{A_{\psi_0}} \mathbb{D}_{\psi_0}.
  \end{equation}
      \end{remark}

    \begin{corollary}\label{dualECram}Let $K/F$ be ramified. 
  Let $(\mathcal{P}, \iota)$ be a CM-pair over an $O_E$-algebra $R$ such that
  $p$ is nilpotent in $R$. If $(X, \iota)$ satisfies the Eisenstein condition
  $({\rm EC}_r)$, then the conjugate dual
  $(\mathcal{P}^{\wedge}, \iota^{\wedge})$ also satisfies $({\rm EC}_r)$. 
    \end{corollary}
    \begin{proof}
We use the notation of the last remark. 
The banal $\psi$ are treated as in the unramified case. We need only to
check that the conjugate dual satisfies $({\rm EC}_{\psi_0,r})$. The
orthogonal complement of $\mathbf{E}_{A_{\psi_0}} \mathbb{D}_{\psi_0}$ in
$\Hom_R(\mathbb{D}_{\psi_0}, R)$ is
$\mathbf{S}_{\psi_0}\mathbf{E}_{A_{\psi_0}}\Hom_R(\mathbb{D}_{\psi_0}, R)$,
where in the last formula we use the action via $\iota^{\vee}$. We obtain
the result by taking the orthogonal complement of \eqref{KottwitzC36enew}. 
    \end{proof}
To end this subsection, we check that the Kottwitz condition $({\rm KC}_{r})$ is preserved under passage to the conjugate dual. 
    \begin{proposition}\label{dualKC}
      Let $K/F, r$ be arbitrary. Let $(\mathcal{P}, \iota)$ be a CM-pair which
      satisfies $({\rm KC}_{r})$. Then the conjugate dual $(\mathcal{P}^{\wedge}, \iota^{\wedge})$ satisfies
      $({\rm KC}_{r})$.
      \end{proposition}
    \begin{proof}
      We may assume that $R$ is endowed with the structure of an
      $O_{\tilde{E}}$-algebra. We use the notation of the proof of Corollary
      \ref{dualEisen2c}. In particular $\mathbb{D}_{\mathcal{P},R} = \mathbb{D}$
      and $\mathbb{D}_{\mathcal{P}^{\wedge},R} = \mathbb{D}^{\wedge}$, and we write
      $\mathcal{L}$ and $\mathcal{L}^{\wedge}$ for the Lie algebras of
      $\mathcal{P}$ and $\mathcal{P}^{\wedge}$. 
      We have to show that for each $R$-algebra $S$
      and for each $\psi \in \Psi$ 
      \begin{displaymath}
        {\det}_S(a \mid \mathcal{L}^{\wedge}_{S, \psi}) =
        \prod_{\varphi | \psi} \varphi_S(a)^{r_{\varphi}},
        \quad \text{for all}\; a \in O_K \otimes_{O_{K^t},\psi_R} S. 
      \end{displaymath}
      To show this, we may replace $\mathcal{P}$ by its base change
      $\mathcal{P}_S$. Therefore it is enough to consider the case $S = R$.
      Since $\mathbb{D}_{\psi}$ is locally on $\Spec R$ a free
      $O_K \otimes_{O_{K^t},\psi_R} R$-module of rank $2$, we find
      \begin{displaymath}
        \det(a \mid \mathbb{D}_{\psi}) = \prod_{\varphi\in\Phi_\psi}
        \varphi_R(a)^{2}, \quad \text{for} \; a \in O_K \otimes_{O_{K^t},\psi_R} R. 
      \end{displaymath}
      Since $\mathcal{L}_{\psi} = \mathbb{D}_{\psi}/\mathbb{D}^1_{\psi}$ we find
      \begin{displaymath}
        \det(\bar{a} \mid \mathbb{D}^1_{\bar{\psi}}) =
        \prod_{\bar{\varphi}|\bar{\psi}} \bar{\varphi}_R(\bar{a})^{(2 - r_{\bar{\varphi}})}
        = \prod_{\varphi\in\Phi_\psi} \varphi_R(a)^{r_{\varphi}}.
      \end{displaymath}
      The perfect pairing (\ref{dualEisen2c}) induces a perfect pairing
      \begin{displaymath}
\mathbb{D}^1_{\bar{\psi}} \times \mathcal{L}^{\wedge}_{\psi} \longrightarrow R.  
      \end{displaymath}
      Therefore we obtain
      \begin{displaymath}
      \det(a \mid \mathcal{L}^{\wedge}_{\psi}) =
      \det(\bar{a} \mid \mathbb{D}^1_{\bar{\psi}}) =
      \prod_{\varphi\in\Phi_\psi} \varphi_R(a)^{r_{\varphi}}.
      \end{displaymath}
      Therefore $({\rm KC}_{\bar{\psi},r})$ for $\mathcal{P}$ implies
      $({\rm KC}_{\psi, r})$ for $\mathcal{P}^{\wedge}$. 
      \end{proof}

\subsection{The pre-contracting functor} \label{ss:tcffcm}

Let $(\mathcal{P}, \iota) = (P,Q,F,\dot{F}, \iota)$ be a CM-pair of type $r$
over an $O_{E'}$-algebra $R$ such that $p$ is nilpotent in $R$. We assume that 
$(\mathcal{P}, \iota)$  satisfies $({\rm EC}_{r})$. In other words,
$(\mathcal{P}, \iota) \in \mathfrak{d}\mathfrak{P}_{r,R}$, cf. Definition
\ref{KatCMpairs}. We will define a functor that associates to
$(\mathcal{P}, \iota)$ a new display $\mathcal{P}' = (P',Q',F',\dot{F}')$ with 
an action 
    \begin{displaymath}
      \iota': O_K \longrightarrow \End \mathcal{P}', 
    \end{displaymath}
such that $(P, \iota) = (P', \iota')$ and $Q \subset Q' \subset P' = P$.
In particular we obtain a natural surjection
$\Lie \mathcal{P} \rightarrow \Lie \mathcal{P}'$ of $O_K$-modules.     
In the case where $r$ is banal, we define $Q' = P$. Therefore the display
$\mathcal{P}'$ will be \'etale in this case. In the case where $r$ is special,
the restriction of the  action $\iota'$ to $O_F$ will be strict with respect to
$\varphi_{0, R} : O_F \longrightarrow O_E \longrightarrow R$. 
We will call this functor the \emph{pre-contracting functor}.\index{pre-contracting functor}  

Let us first restrict our attention to the case where $K/F$ is a field
extension. The case $K = F \times F$ will be treated separately  because it
needs  different notations, see p. \pageref{splitpage}, starting before
eq. \eqref{zerlegt1e}. 
Each $\psi: K^t \longrightarrow \bar{\mathbb{Q}}_p$ induces a homomorphism 
\begin{equation}\label{SchlangePsi} 
  \tilde{\psi}: O_{K^t} \longrightarrow  W(O_{K^t}) \xra{W(\psi)} W(O_{E'}). 
\end{equation}
For an $O_{E'}$-algebra $R$ we deduce a homomorphism 
$\tilde{\psi}_R: O_{K^t} \longrightarrow W(R)$ that is equivariant 
with respect to the Frobenius homomorphisms on both sides.
This induces decompositions  
\begin{equation}\label{OWdecomp1e}
  \begin{aligned} 
     O_{K^t} \otimes_{\mathbb{Z}_p} W(R) &\cong
     \prod\nolimits_{\psi\in \Psi} W(R),\\ 
     O_K \otimes_{\mathbb{Z}_p} W(R) &\cong 
\prod\nolimits_{\psi\in \Psi} O_K \otimes_{O_{K^t}, \tilde{\psi}_R} W(R) 
    \end{aligned}
\end{equation}
which lift the decomposition (\ref{KottwitzC3e}).
Let $\sigma \in \Gal(K^t/ \mathbb{Q}_p)$ be the Frobenius automorphism. The
operators $F$ and $V$ act via $W(R)$ on the left hand side of 
(\ref{OWdecomp1e}). On the right hand side this induces maps
\begin{displaymath}
   O_K \otimes_{O_{K^t}, \tilde{\psi}_R} W(R) 
  \begin{array}{c}
    \overset{F}{\longrightarrow}\\[-2mm] 
    \underset{V}{\longleftarrow}
    \end{array}
  O_K \otimes_{O_{K^t}, \tilde{\psi}_R \circ \sigma} W(R),  
  \end{displaymath}
cf. (\ref{LTD5e}). 
Recall that 
    \begin{displaymath}
      \mathbf{E}_{A_{\psi}}(T) = \prod_{\varphi \in A_{\psi}} (T - \varphi(\Pi))
      \in O_{E'}[T]. 
    \end{displaymath}
We lift this to a polynomial with coefficients in $W(O_{E'})$
by taking the Teichm\"uller representatives of the roots,
\begin{equation}\label{Kottwitz36e}
  \tilde{\mathbf{E}}_{A_{\psi}} (T)= \prod_{\varphi \in A_{\psi}}
  (T - [\varphi(\Pi)]) \in W(O_{E'})[T].
\end{equation}
The image of this polynomial by the homomorphism
$W(O_{E'}) \longrightarrow W(R)$ is denoted by
$\tilde{\mathbf{E}}_{A_{\psi},R}(T)$. If we reduce with respect to
$\mathbf{w}_0: W(O_{E'}) \longrightarrow O_{E'}$, we obtain
the polynomial $\mathbf{E}_{A_{\psi}}(T)$. We note that in the case where
$R$ is a $\kappa_{E'}$-algebra, we have
\begin{equation}\label{KottwitzC33e}
  \tilde{\mathbf{E}}_{A_{\psi},R}(T) = T^{a_{\psi}}. 
\end{equation}
We consider the ring homomorphism
\begin{equation}\label{Kottwitz21e}
  \begin{aligned} 
  W(O_{E'})[T] \longrightarrow   W(R)[T] & \longrightarrow 
  O_K \otimes_{O_{K^t}, \tilde{\psi}_R} W(R) .\\
   T & \mapsto  \Pi \otimes 1
    \end{aligned}
  \end{equation}
We denote by $\tilde{\mathbf{E}}_{A_{\psi},R}(\Pi \otimes 1)$ the image of
$\tilde{\mathbf{E}}_{A_{\psi}}(T)$ under (\ref{Kottwitz21e}). 

Let now $(\CP, \iota)\in\mathfrak{d}\mathfrak{P}_{r, R}$. We obtain
decompositions of the $O_K \otimes_{\mathbb{Z}_p} W(R)$-modules $P$ and $Q$,  
\begin{equation}\label{OWdecomp2e}
P = \oplus_{\psi \in \Psi} P_{\psi}, \quad Q = \oplus_{\psi \in \Psi} Q_{\psi} .
  \end{equation}For $x \in P_{\psi}$,
we write
\begin{equation}\label{KottwitzC41e}
  \tilde{\mathbf{E}}_{A_{\psi}} x =
  \tilde{\mathbf{E}}_{A_{\psi},R}(\Pi \otimes 1) x. 
\end{equation}
On the left hand side we consider $P_{\psi}$ as a $W(O_{E'})[T]$-module via
(\ref{Kottwitz21e}). 

We give first the recipe for the construction of $\mathcal{P}'$ for any $R \in \Nilp_{O_{E'}}$. Then we will discuss  the case of
a perfect field. This special case is then used to prove that
 $\mathcal{P}'$ is indeed a display.
 
We begin with the case where $r$ is banal (and $K/F$ is a field extension).
Let $(\mathcal{P}, \iota)$ as above. We define
    \begin{displaymath}
      P' = \oplus_{\psi} P'_{\psi}, \quad Q' = \oplus_{\psi} Q'_{\psi}
    \end{displaymath}
    as follows: 
    for all $\psi$ we set
    \begin{equation}\label{KottwitzC16e}
      P'_{\psi} = Q'_{\psi} = P_{\psi}. 
    \end{equation}
By the Eisenstein condition \eqref{Eisenbanal}, we have
$\tilde{\mathbf{E}}_{A_{\psi}} P_{\psi} \subset Q_{\psi}$.
Then we may define
   \begin{equation}\label{KottwitzC17e}
      \begin{aligned}
        \dot{F}': Q'_{\psi} \longrightarrow P'_{\psi \sigma}, \quad& \dot{F}'(x) =
        \dot{F}(\tilde{\mathbf{E}}_{A_{\psi}} x), & x \in P_{\psi}\\
        F': P'_{\psi} \longrightarrow P_{\psi \sigma},\quad & F'(x) =
        F(\tilde{\mathbf{E}}_{A_{\psi}} x), & x \in P_{\psi} .
      \end{aligned}
    \end{equation}
We define $F': P' \longrightarrow P'$ and $\dot{F}': Q' \longrightarrow P'$ as
the direct sum of the maps above. We have to prove that
$\mathcal{P}' = (P',Q',F', \dot{F}')$ is a display. The only non-trivial
property we have to check is that $\dot{F}': P_{\psi} \longrightarrow P_{\psi \sigma}$
is an $F$-linear isomorphism. We postpone the verification, cf. p. \pageref{page:banal}, below  \eqref{KottwitzC37e}. 

    We now define the pre-contracting functor in the case where $r$ is special and $K/F$ 
    unramified. In this case we have
    $\psi_0 \neq \bar{\psi}_0$.  If $\psi$ is banal, i.e., if
    $\psi \notin \{\psi_0, \bar{\psi}_0\}$, we keep the definitions
    (\ref{KottwitzC16e}) and (\ref{KottwitzC17e}).
    We set $P'_{\psi_0} = P_{\psi_0}$ and we define $Q'_{\psi_0}$
    as the kernel of the following map,
     \begin{equation}\label{KottwitzC19e}
      P'_{\psi_0} = P_{\psi_0} \longrightarrow P_{\psi_0}/Q_{\psi_0}
       \xra{\mathbf{E}_{A_{\psi_0,R}} }
     \mathbf{E}_{A_{\psi_0,R}}(P_{\psi_0}/Q_{\psi_0}) \subset
      P_{\psi_0}/Q_{\psi_0}. 
    \end{equation}
    It follows from Proposition \ref{KottwitzC4p} that
    $\mathbf{E}_{A_{\psi_0,R}}(P_{\psi_0}/Q_{\psi_0})$ is locally free
    of rank $1$ and is a direct summand of
    $P_{\psi_0}/Q_{\psi_0}$. Therefore
    \begin{equation}\label{KottwitzC20e}
      P'_{\psi_0}/Q'_{\psi_0} \cong \mathbf{E}_{A_{\psi_0,R}}(P_{\psi_0}/Q_{\psi_0})
    \end{equation} 
is locally free of rank $1$ and, as remarked at the end of the proof of
Proposition \ref{KottwitzC5p}, an element $a \in O_K \otimes_{O_{K^t}, \psi_{0}} R$
acts on (\ref{KottwitzC20e}) by multiplication with $\varphi_{0,R}(a)$.
This makes sense because $\varphi_0: O_K \longrightarrow O_{\tilde{E}}$ factors
through $O_{E'} \subset O_{\tilde{E}}$. We define 
      \begin{equation}\label{KottwitzC21e}
      \begin{aligned}
        F': P'_{\psi_0} \longrightarrow P'_{\psi_0 \sigma}, \quad &
        F'(x) = F(\tilde{\mathbf{E}}_{A_{\psi_{0,R}}} x),\\ 
        \dot{F}': Q'_{\psi_0} \longrightarrow P'_{\psi_0 \sigma}, \quad &
        \dot{F}'(y) = \dot{F}(\tilde{\mathbf{E}}_{A_{\psi_{0,R}}} y). 
      \end{aligned}
    \end{equation}
    The last equation makes sense because, by definition,
    $\tilde{\mathbf{E}}_{A_{\psi_0,R}} Q'_{\psi_0} \subset Q_{\psi_0}$. 
    The definitions of the modules $P'_{\bar{\psi}_0}, Q'_{\bar{\psi}_0}$
    and the restrictions of $F'$ and $\dot{F}'$ to these modules are
    defined by interchanging the roles of $\psi_0$ and $\bar{\psi}_0$.
    This completes the definition of
    \begin{equation}
      F': \oplus_{\psi} P'_{\psi} \longrightarrow \oplus_{\psi} P'_{\psi}, \quad
      \dot{F}': \oplus_{\psi} Q'_{\psi} \longrightarrow \oplus_{\psi} P'_{\psi}.
    \end{equation}
Again we postpone the verification that $\mathcal{P}'$ is a display, cf. below \eqref{KottwitzC22e}. The tangent space $\mathcal{L}' = P'/Q'$ is a locally free $R$-module
    of rank $2$. It has a decomposition
    \begin{displaymath}
      \mathcal{L}' = \mathcal{L}'_{\psi_0} \oplus \mathcal{L}'_{\bar{\psi}_0} ,
    \end{displaymath}
    where an element $a \in O_K \otimes_{O_{K^t}, \psi_{0}} R$ acts on the
    first summand by multiplication with $\varphi_{0,R}(a)$ and an element
    $a \in O_K \otimes_{O_{K^t}, \bar{\psi}_{0}} R$ acts on the 
    second summand by multiplication with $\bar{\varphi}_{0,R}(a)$.

Next we define the pre-contracting functor in the case where $r$ is special and
$K/F$ ramified. In this case we have
 $\psi_0 = \bar{\psi_0}$. For banal $\psi$, we keep the definitions 
 (\ref{KottwitzC16e}) and (\ref{KottwitzC17e}). The $R$-module
 $\mathbf{E}_{A_{\psi_{0,R}}}(P_{\psi_0}/Q_{\psi_0}) \subset P_{\psi_0}/Q_{\psi_0}$
 is a direct summand which is locally free
 of rank $2$. We set $P'_{\psi_0} = P_{\psi_0}$ and we define
 $Q'_{\psi_0}$ as the kernel of (\ref{KottwitzC19e}). We define $F'$
 and $\dot{F}'$ by (\ref{KottwitzC21e}). Then
 $P'_{\psi_0}/Q'_{\psi_0} \cong \mathbf{E}_{A_{\psi_{0,R}}}(P_{\psi_0}/Q_{\psi_0})$ 
 is locally free and we define as before 
 $\mathcal{P}' = (P',Q',F', \dot{F}')$, with its 
 $O_K$-action $\iota'$. It follows from Proposition \ref{KottwitzC6p} that
 the action of $O_F$ on $\mathbf{E}_{A_{\psi_{0,R}}}(P_{\psi_0}/Q_{\psi_0})$ via
 $\iota$ coincides with the action of via $\varphi_{0}$, i.e., the
 action via $\iota'$ on $\mathcal{P}'$ is strict. That $\CP'$ is a display is proved around \eqref{KottwitzC23e}.

Now we consider the case where $R = k$ is a perfect field in more detail.
We know that $P_{\psi}$ is a free module of rank $2$ over the
discrete valuation ring $O_K \otimes_{O_{K^t}, \tilde{\psi}} W(k)$. Therefore
$P_{\psi}/\Pi P_{\psi}$ is a $k$-vector space of dimension two. In the perfect field case, we have
now also the operator $V$,
\begin{equation}
  F: P_{\psi} \longrightarrow P_{\psi \sigma}, \quad V: P_{\psi \sigma}
  \longrightarrow P_{\psi}, \quad V(P_{\psi \sigma}) = Q_{\psi}.     
\end{equation}
We will see that in all cases the Eisenstein condition implies that
$V(P_{\psi \sigma}) \subset \Pi^{a_{\psi}}P_{\psi}$. Therefore we may define
operators $F'$ and $V'$:
\begin{equation}\label{KottwitzC39e}
  F' = \Pi^{a_{\psi}}F: P_{\psi} \longrightarrow P_{\psi \sigma}, \quad
  V' = \Pi^{-a_{\psi}}V : P_{\psi \sigma} \longrightarrow P_{\psi}. 
\end{equation}
The Dieudonn\'e module of the display $\mathcal{P}'$ in the sense of
Proposition \ref{perfFrame1p} will then be $(P, F',V')$.

\label{page:banal}We begin with the case when $r$ is banal.
Now $\mathbf{E}_{A_{\psi},k} (T)= T^{a_{\psi}}$ acts
on $P_{\psi}$ as multiplication by $\Pi^{a_{\psi}}$.  
By the Eisenstein condition \eqref{Eisenbanal},  $\Pi^{a_{\psi}}$annihilates $P_{\psi}/Q_{\psi}$. This implies
$\Pi^{a_{\psi}}P_{\psi} \subset V(P_{\psi \sigma})$. By the rank condition,
the factor $P_{\psi}/V(P_{\psi \sigma})$ has length $2a_{\psi}$ as
$O_K \otimes_{O_{K^t}, \tilde{\psi}} W(k)$-module. Since the same is true for the
factor module $P_{\psi}/\Pi^{a_{\psi}} P_{\psi}$, we obtain
    \begin{equation}\label{KottwitzC37e}
\Pi^{a_{\psi}} P_{\psi} = V(P_{\psi \sigma}).
      \end{equation}
Therefore
$\dot{F}'=\dot{F} \Pi^{a_{\psi}} = V^{-1} \Pi^{a_{\psi}}: P_{\psi} \longrightarrow P_{\psi \sigma}$
is bijective. This shows that $\mathcal{P}'$ is a display.
We set
\begin{displaymath}
F' = \oplus_{\psi} F \Pi^{a_{\psi}}, \quad V' = \oplus_{\psi} \Pi^{-a_{\psi}} V. 
  \end{displaymath}
Then $(P, F', V')$ is the Dieudonn\'e module associated to $\mathcal{P}'$. 

We obtain from (\ref{KottwitzC37e}) in the ramified case that
$V^{2f} P_{\psi} = \Pi^{2ef} P_{\psi} = p^f P_{\psi}$ for all $\psi$ and
in the unramified case that $V^{2f}P_{\psi} = \pi^{ef} P_{\psi}$. This implies
that in both cases $\mathcal{P}$ is isoclinic of slope $1/2$. 

Now we can verify that $\mathcal{P}'$ is a display for $r$ banal, for an   arbitrary
$O_{E'}$-algebra $R$.
Let $(\mathcal{P}, \iota)\in \mathfrak{d}\mathfrak{P}_{r, R}$. We must show that $\dot{F}: P \longrightarrow P$  is a
Frobenius-linear isomorphism. We may assume
that $P$ is a free $W(R)$-module. Let $\det \dot{F}$ be the determinant of
the matrix of $\dot{F}$ with respect to any given basis of the $W(R)$-module
$P$. We must verify that $\det \dot{F}$ is a unit in $W(R)$. We have shown
that, for each homomorphism $R \longrightarrow k$  to a perfect field $k$, the image
of $\det \dot{F}$ by $W(R) \longrightarrow W(k)$ is a unit in $W(k)$. In
particular $\mathbf{w}_0(\det \dot{F}) \in R$ has a nonzero image under any
homomorphism $R \longrightarrow k$. But then $\mathbf{w}_0(\det \dot{F})$ is
a unit in $R$, and this implies that $\det \dot{F} \in W(R)$ is a unit. 
This finishes  in the banal case the proof that $\mathcal{P}'$ is a display.

Next we consider the  case when $r$ is special and $K/F$ unramified.
By our conventions, $\Pi = \pi$ is the prime element of $F$. Let $R = k$
be a perfect field and let $\mathcal{P} = (P, F, V)$, regarded as a
Dieudonn\'e module. If $\psi \in \Psi$ is banal, we find as above that
$VP_{\psi \sigma} = \pi^{a_{\psi}}P_{\psi}$.
Now let $\psi \in \{\psi_0, \bar{\psi}_0 \}$. Since $\pi^{a_{\psi} + 1}$
annihilates $\Lie_{\psi} \mathcal{P}$ by the Eisenstein condition
\eqref{furtherDr2unram}, we obtain
$\pi^{a_{\psi}+1}P_{\psi} \subset VP_{\psi \sigma}$. We note that $P_{\psi}$ is
a $O_K \otimes_{O_{K^t}, \tilde{\psi}} W(k)$-module of rank $2$. Therefore 
the factor module of the last inclusion is, by the rank condition,
a $O_K \otimes_{O_{K^t}, \tilde{\psi}} W(k)$-module of length $1$ and is
therefore annihilated by $\pi$. This implies 
    \begin{equation*} 
      \pi^{a_{\psi}+1}P_{\psi} \subset VP_{\psi \sigma} \subset
      \pi^{a_{\psi}} P_{\psi}.
    \end{equation*}
    In particular
    \begin{equation}\label{KottwitzC22e}
      P_{\psi}/\pi^{a_{\psi}}P_{\psi} \overset{\sim}{\longrightarrow}
      \Lie_{\psi} \mathcal{P} / \pi^{a_{\psi}}\Lie_{\psi} \mathcal{P} 
    \end{equation}
    is an isomorphism, as claimed in the beginning of the proof of Proposition 
\ref{KottwitzC4p}. 
By definition (\ref{KottwitzC19e}) 
we have $Q'_{\psi} = \pi^{-a_{\psi}} VP_{\psi \sigma}$. The map
$\dot{F}' = \pi^{a_{\psi}} \dot{F}: Q'_{\psi} \longrightarrow P_{\psi \sigma}$
is therefore surjective. Since we know this fact also for banal $\psi$
we conclude 
that $(P, Q', F', \dot{F}')$ is a display. The associated Dieudonn\'e module
is $(P, F', V')$, where 
\begin{equation}\label{moreeq}
\begin{aligned}
  F'_\psi=\pi^{a_\psi}F_\psi&:  P_{\psi } \longrightarrow P_{\psi\sigma}, \\  
  V'_{\psi\sigma} = \pi^{-a_\psi}V_{\psi\sigma} &: P_{\psi \sigma} \longrightarrow P_{\psi} . \\
    \end{aligned}
\end{equation}
Now we return to an arbitrary $O$-algebra $R$ such that $p$ is nilpotent in
$R$. We note that the definition of $(P', Q')$ commutes with arbitrary base
change because $Q'_{\psi_0}/I(R)P'_{\psi_0}$ is defined as the kernel of an
epimorphism of projective $R$-modules, 
\begin{displaymath}
P_{\psi_0}/I(R)P_{\psi_0} \longrightarrow \mathbf{E}_{A_{\psi_0}}(P_{\psi_0}/Q_{\psi_0}).  
\end{displaymath}
We choose a normal decomposition of $(P', Q')$,
\begin{displaymath}
P' = T' \oplus L' ,
  \end{displaymath}
together with the Frobenius-linear endomorphism $\Phi': P' \longrightarrow P'$
of the $W(R)$-module $P'$ such that the restriction of $\Phi'$ to $T'$ is
$F'$ and the restriction to $L'$ is $\dot{F}'$. We have to show that the
determinant of $\Phi$ in a locally chosen basis is a unit. Since we know
that this is true after any base change $R \longrightarrow k$ with $k$ a perfect
field, this follows  as in the banal case.

We can determine the possible slopes of $\mathcal{P}$ when $r$ is special and $K/F$  unramified. Let $\CP=(P, F, V)$ over the perfect field $k$. By (\ref{KottwitzC39e}) we have $(V')^{2f} = \pi^{-(ef-1)}V^{f}$. Let
\begin{displaymath}
P_{\mathbb{Q}} = \oplus_{\lambda} N(\lambda) 
\end{displaymath}
be the decomposition into isoclinic components. Fix $\lambda = r/s$. Then we find
a $W(k)$-lattice $\Lambda \subset N(\lambda)$ such that
$V^s \Lambda = p^r \Lambda$. From $V^{2fs} \Lambda = p^{2fr} \Lambda$ we obtain
\begin{displaymath}
  (\pi^{ef-1}(V')^{2f})^s \Lambda = p^{2rf} \Lambda, \; \text{i.e.,} \; 
  (V')^{2fs} \Lambda = \pi^{-efs}\pi^s p^{2rf} \Lambda.
  \end{displaymath}
We write the right hand side as $p^{-sf}p^{s/e}p^{2rf} \Lambda$. This shows that
$N(\lambda) \subset P'_{\mathbb{Q}}$ is an isoclinic rational Dieudonn\'e
submodule of slope
\begin{displaymath}
\frac{-sf + (s/e) + 2fr}{2fs} = -\frac{1}{2} + \frac{1}{2d} + \lambda. 
\end{displaymath}
Let us apply the Ahsendorf functor to $\mathcal{P}'$. We obtain a
$\mathcal{W}_{O_F}(k)$-Dieudonn\'e module $(P_{\rm{c}}, F_{\rm{c}}, V_{\rm{c}})$ of height $4$ and
dimension $2$. The slopes of $P_{\rm{c}}$ are by Proposition \ref{Adorfslopes1p} 
\begin{equation}\label{CFslopes1e} 
d(\lambda - \frac{1}{2}) + \frac{1}{2}. 
\end{equation}
The action of $O_K \otimes_{O_F} W_{O_F}(k) \cong W_{O_F}(k) \times W_{O_F}(k)$
on $P_{\rm{c}}$ defines a decomposition $P_{\rm{c}} = P_{{\rm c}, 0} \oplus P_{{\rm c}, 1}$ such that 
$V_{\rm{c}} (P_{{\rm{c}}, 0}) \subset P_{{\rm{c}}, 1}$ and $V_{\rm{c}} (P_{{\rm{c}}, 0}) \subset P_{{\rm{c}}, 1}$. The
$W_{O_F}(k)$-module $P_{{\rm{c}}, 0}$ with the semi-linear operator $V_{\rm{c}}^2$ is of height
$2$ and dimension $2$. Therefore the possible slopes of $(P_{{\rm{c}}, 0}, V_{{\rm{c}}}^2)$
are with multiplicities $(1,1)$ or $(0,2)$. We conclude that the slopes of
$(P_{\rm{c}}, V_{\rm{c}})$ are with multiplicities $(1/2, 1/2, 1/2, 1/2)$ or
$(0, 0, 1, 1)$. From (\ref{CFslopes1e}) we find that in the first case
all slopes $\lambda$ of $\mathcal{P}$ are $1/2$, while in the second case
we obtain the two slopes $\lambda = 1/2 - 1/2d$ and
$\lambda = 1/2 + 1/2d$.

Now we consider the case where $r$ is special and $K/F$ is ramified.
As in the last case, it is enough to verify that $\mathcal{P}'$ is a display
 when  $R = k$ is a perfect field. Recall that $a_{\psi} = e$ for $\psi$ banal
and that $a_{\psi_0} = e -1$. As above we find $VP_{\psi \sigma} = \Pi^e P_{\psi}$
for $\psi$ banal. By the Eisenstein condition (\ref{furtherDr2ram}),
$\mathcal{L}_{\psi_0}$ is annihilated by $\Pi^{e+1}$ and the $k$-vector space
$\Pi^{e-1} \mathcal{L}_{\psi_0}$ has dimension at most $2$. 
We consider the following
      filtration by subvector spaces,
      \begin{displaymath}
        \mathcal{L}_{\psi_0} \supset \Pi \mathcal{L}_{\psi_0} \supset
        \Pi^2 \mathcal{L}_{\psi_0} \supset \ldots \supset
        \Pi^e \mathcal{L}_{\psi_0} \supset \Pi^{e+1} \mathcal{L}_{\psi_0} = 0.
      \end{displaymath}
We have $\dim_{k} \Pi^m \mathcal{L}_{\psi_0}/ \Pi^{m+1} \mathcal{L}_{\psi_0} \leq 2$
for all $m \geq 0$ since $\CL_{\psi_0}$ is a quotient of $P_{\psi_0}$, which is a
free $O_K \otimes_{O_{K^t}, \tilde{\psi}_0} W(k)$-module of rank $2$. Therefore we
find 
      \begin{displaymath}
        \dim_{k} \mathcal{L}_{\psi_0} = \dim_{k}
        \mathcal{L}_{\psi_0}/\Pi^{e-1} \mathcal{L}_{\psi_0} + \dim_{k}
        \Pi^{e-1} \CL_{\psi_0} \leq 2(e-1) + 2 = 2e = \dim_{k} \mathcal{L}_{\psi_0}. 
      \end{displaymath}
      We must have the equality
      \begin{displaymath}
        \dim_{k} \mathcal{L}_{\psi_0}/ \Pi^{e-1} \mathcal{L}_{\psi_0} = 2(e-1),
        \quad \dim_{k} \Pi^{e-1}\mathcal{L}_{\psi_0} = 2. 
      \end{displaymath}
      The first equation shows that the natural map
      \begin{equation}\label{KottwitzC23e}
        P_{\psi_0}/ \Pi^{e-1} P_{\psi_0} \longrightarrow
        \mathcal{L}_{\psi_0}/ \Pi^{e-1} \mathcal{L}_{\psi_0} 
      \end{equation}
is an isomorphism of vector spaces, as asserted in the beginning of the
proof of Proposition \ref{KottwitzC4p}.       
Finally we have by definition
$Q'_{\psi_0} = \Pi^{-e+1}Q_{\psi_0} = \Pi^{-e+1}VP_{\psi_0 \sigma}$. Therefore
\begin{displaymath}
\dot{F}' = \Pi^{e-1}V^{-1}: Q'_{\psi_0} \longrightarrow P_{\psi_0 \sigma}  
  \end{displaymath}
is bijective. We conclude that $(P', Q', F', \dot{F}')$ is a display. The
associated Dieudonn\'e module is $(P, F', V')$, where
\begin{equation}\label{defP'1e}
  \begin{aligned}
  V' = \Pi^{-e} V&:  P_{\psi \sigma} \longrightarrow P_{\psi}, &
 \psi\neq \psi_0\\  
  V' = \Pi^{-e+1} V&: P_{\psi_0 \sigma} \longrightarrow P_{\psi_0} . &\\
    \end{aligned}
\end{equation}
As in the unramified case we conclude for an arbitrary $R \in \Nilp_{O_{E'}}$
that our definitions (\ref{KottwitzC16e}), (\ref{KottwitzC17e}),
(\ref{KottwitzC21e}) give a display $\mathcal{P}' = (P', Q', F', \dot{F}')$.

In the case of a perfect field $k$, the slopes of $\mathcal{P}'$ are computed
in the same way as in the unramified case. We have the equation
$(V')^f = \Pi^{-(ef-1)}V^f$. Let $N(\lambda) \subset P_{\mathbb{Q}}$ be an isoclinic
component. We find a lattice $\Lambda \subset N(\lambda)$ such that
$V^s \Lambda = p^r \Lambda$. We obtain that
\begin{displaymath}
(V')^{sf} \Lambda = \Pi^{-efs} \Pi^s p^{rf} \Lambda.
  \end{displaymath}
Since $\Pi^{2e}$ and $p$ differ by a unit, this implies that
$N(\lambda) \subset P'_{\mathbb{Q}}$ is isoclinic of slope
\begin{displaymath}
  \frac{-(fs/2) + (s/2e) + rf}{sf} = -1/2 + 1/2d + \lambda =
  (\lambda - 1/2) + 1/2d. 
\end{displaymath}
If we apply the Ahsendorf functor, we obtain an
$\mathcal{W}_{O_F}(k)$-Dieudonn\'e module $(P_{\rm{c}}, F_{\rm{c}}, V_{\rm{c}})$ with slopes
$d(\lambda - 1/2) + 1/2$. If we consider the
$O_K \otimes_{O_F} W_{O_F}(k)$-module $P_{\rm{c}}$ with the semi-linear operator $V_{\rm{c}}$
the possible slopes with multiplicity are $(1,1)$ or $(0,2)$ because
$(P_{\rm{c}},V_{\rm{c}})$ is of height $2$ and dimension $2$. If we regard $(P_{\rm{c}}, V_{\rm{c}})$ over
$W_{O_F}(k)$, the heights are multiplied with $2$ and then the possible
heights are $(1/2, 1/2, 1/2, 1/2)$ or $(0, 0, 1, 1)$. As in the
unramified case we conclude that $\mathcal{P}$ is either isoclinic of slope
$1/2$ or has exactly two slopes $1/2 - 1/2d$ and $1/2 + 1/2d$.

\label{splitpage}Finally we consider the case where $r$ is banal and $K = F \times F$.
We set $\Theta = \Hom_{\text{$\mathbb{Q}_p$-Alg}} (F^t, \bar{\mathbb{Q}}_p)$.
\index[NO]{ZZHA@$\Theta$} In this case
$\sigma$ will denote the Frobenius automorphism in $\Gal(F^t/F)$. If
we compose $\theta \in \Theta$ with the first, resp. second, projection
$K^t = F^t \times F^t \longrightarrow F^t$ we obtain $\theta_1, \theta_2 \in \Psi$.
Via the first, resp. the second,  projection, we obtain isomorphisms
\begin{displaymath}
  (O_F \times O_F) \otimes_{(O_{F^t} \times O_{F^t}), \theta_i} O_{E'} \cong
  O_F \otimes_{O_{F^t}, \theta} O_{E'}, \quad i = 1, 2. 
  \end{displaymath}
This leads to the decomposition
\begin{equation}\label{zerlegt1e}
  O_K \otimes O_{E'} = \prod_{\psi \in \Psi} O_K \otimes_{O_{K^t}, \psi} O_{E'} = 
  \big(\prod_{\theta \in \Theta} O_F \otimes_{O_{F^t}, \theta} O_{E'}\big) 
  \times \big(\prod_{\theta \in \Theta} O_F \otimes_{O_{F^t}, \theta} O_{E'}\big). 
\end{equation} 
Assume that $\psi\in \Psi$ factors through $\theta \in \Theta$. We define
$\tilde{\psi}$ as the composite 
\begin{equation}\label{zerlegt7e}
  O_{K^t} = O_{F^t} \times O_{F^t} \xra{{\rm proj.}}O_{F^t}
  \longrightarrow W(O_{F^t}) \xra{W(\theta)} W(O_{E'}).
\end{equation}
The first map  is the projection to the first or second factor according
to $\psi$. We denote by $\tilde{\theta}$ the composite of the last two arrows
in \eqref{zerlegt7e}. We obtain the decomposition
\begin{equation}\label{zerlegt1ea}
  \begin{aligned}
    O_K \otimes_{\mathbb{Z}} W(O_{E'}) = & \prod_{\psi \in \Psi}
    O_K  \otimes_{O_{K^t}, \tilde{\psi}} W(O_{E'})\\
  = &  
  \Big(\prod_{\theta \in \Theta} O_F  \otimes_{O_{F^t}, \tilde{\theta}} W(O_{E'})\Big) 
  \times \Big(\prod_{\theta \in \Theta} O_F \otimes_{O_{F^t}, \tilde{\theta}} W(O_{E'})\Big). 
    \end{aligned}
  \end{equation}
On the right hand side, the first set of factors correspond to those $\psi$ which
factor over the first projection and the second set of factors correspond to
those $\psi$ which  factor over the second projection.

We consider a CM-pair $(\mathcal{P}, \iota)$ over
$R \in \Nilp_{O_{E'}}$ which satisfies the Eisenstein condition. By
(\ref{zerlegt1e}) we obtain a decomposition
\begin{equation}\label{zerlegt8e}
  P = P_1 \times P_2 = (\bigoplus_{\theta \in \Theta} P_{1,\theta}) \oplus 
  (\bigoplus_{\theta \in \Theta} P_{2, \theta}) .
  \end{equation}
This decomposition corresponds to the decomposition into displays
$\mathcal{P} = \mathcal{P}_1 \oplus \mathcal{P}_2$ induced by the
$O_F \times O_F$-action on $\mathcal{P}$. By the definition of a CM-pair
(at the beginning of subsection \ref{ss:loctriples}), the displays $\mathcal{P}_1$ and
$\mathcal{P}_2$ have both height $2d$. 

The maps $F$ and $\dot{F}$ of the display $\mathcal{P}$ induce maps
  \begin{equation}
   F: P_{i, \theta} \longrightarrow P_{i, \theta \sigma}, \quad
  \dot{F}: Q_{i, \theta} \longrightarrow P_{i, \theta \sigma}. 
  \end{equation}
The polynomial $\tilde{\mathbf{E}}_{A_{\psi}} \in W(O_{E'})[T]$ is defined
as before, cf. (\ref{Kottwitz36e}).
For $i = 1,2$ we define the displays
$\mathcal{P}'_i = (P'_i, Q'_i, F'_i, \dot{F}'_i)$ as follows 
\begin{displaymath}
  P'_i = Q'_i = P_i,  \quad \dot{F}'_i(x) =
  \dot{F}(\tilde{\mathbf{E}}_{A_{\theta_i}} x), \;
  F'_i(x) = F(\tilde{\mathbf{E}}_{A_{\theta_i}} x) , \quad
  x \in P_{i\theta}.
  \end{displaymath}
Here, by the convention \eqref{KottwitzC41e}, $\tilde{\mathbf{E}}_{A_{\theta_i}}$ acts as the
  multiplication by
  $\tilde{\mathbf{E}}_{A_{\theta_i},R}(\pi \otimes 1) \in
  O_F \otimes_{O_{F^t}, \tilde{\theta_i}} W(R)$. We set $\mathcal{P}' = \mathcal{P}'_1 \oplus \mathcal{P}'_2$. 
  As in the unramified banal case, the  verification that $\mathcal{P}'$ is a
  display reduces to the case of a perfect field. However, when $R$ is a $\kappa_{E'}$-algebra, then 
  $\tilde{\mathbf{E}}_{A_{\theta_i},R}(\pi \otimes 1) = \pi^{a_{\theta_i}} \otimes 1$.
  If $R = k$ is a perfect field, we consider the Dieudonn\'e module
  $(P_i, F_i, V_i)$ of $\mathcal{P}_i$.  We have
  \begin{equation}\label{zerlegt10e}
V_i(P_{i, \theta \sigma}) = \pi^{a_{\theta_i}}P_{i, \theta} .
    \end{equation}
  We define
  \begin{displaymath}
    F'_i = \pi^{a_{\theta_i}} F_i: P_{i, \theta} \longrightarrow P_{i, \theta \sigma}, \quad 
    V'_i = \pi^{-a_{\theta_i}} V_i: P_{i, \theta \sigma} \longrightarrow P_{i, \theta}. 
    \end{displaymath}
Then $(P_i, F'_i, V'_i)$ is the Dieudonn\'e module of $\mathcal{P}'_i$. 
Finally we determine the slopes of $\mathcal{P}$. If we iterate
(\ref{zerlegt10e}) we find
\begin{displaymath}
  V^f P_i = \pi^{\sum_{\theta} a_{\theta_i}} P_i.
  \end{displaymath} 
We\index[NO]{AAD@$a_i$} set $a_i=\sum_{\theta} a_{\theta_i}$. Then $a_1 + a_2 = ef = d$ because
$a_{\theta_{1}} + a_{\theta_{2}} = e$. We obtain easily that 
\begin{equation}\label{sumsplit}
2a_i = \dim \mathcal{P}_i.
  \end{equation}
It follows that $\CP_i$ is isoclinic of slope
$\lambda_i = a_i/d$ and that $\lambda_1 + \lambda_2 = 1$. 

We summarize the properties of our constructions.
\begin{definition}\label{Kontrakt1d}
  Let $R \in \Nilp_{O_{E'}}$. 
  We define  categories
  $\mathfrak{d}\mathfrak{P}'_{r,R}$ and
  \index[NO]{PCF@$\mathfrak{d}\mathfrak{P}'_{r,R}$} $\mathfrak{P}'_{r,R}$
  \index[NO]{PCE@$\mathfrak{P}'_{r,R}$}  as follows.  
  \begin{enumerate}
 \item If $r$ is banal, then  $\mathfrak{d}\mathfrak{P}'_{r,R}$ is 
  the category of pairs $(\mathcal{P}', \iota')$, where $\mathcal{P}'$
  is an \'etale  display (i.e., $P' = Q'$) of height $4d$ and where $\iota'$ is an
  $O_K$-action. In the split case $O_K = O_F \times O_F$, we require in addition that
  in the induced decomposition
  $\mathcal{P}' = \mathcal{P}'_1 \oplus \mathcal{P}'_2$ both factors have
  height $2d$. 
  
  \item If $r$ is special and $K/F$ is unramified, then  the category
  $\mathfrak{d}\mathfrak{P}'_{r,R}$ is 
  the category of pairs $(\mathcal{P}', \iota')$, where $\mathcal{P}'$
  is a display of height $4d$ and dimension $2$ with an action
  $\iota': O_K \longrightarrow \End \mathcal{P}'$ such that the action
  of $\iota'$ restricted to $O_F$ is strict with respect to
  $\varphi_{0,R}: O_F \overset{\varphi_0}{\longrightarrow} O_{E'} \longrightarrow R$  
  and such that $\Lie \mathcal{P}' = P'/Q'$ is locally on $\Spec R$ a free
  $O_K \otimes_{O_F, \varphi_{0,R}} R$-module of rank 1.
  
  \item If $r$ is special and $K/F$ is ramified, then the category
  $\mathfrak{d}\mathfrak{P}'_{r,R}$ is 
  the category of pairs $(\mathcal{P}', \iota')$, where $\mathcal{P}'$
  is a display of height $4d$ and dimension $2$ with an action
  $\iota': O_K \longrightarrow \End \mathcal{P}'$ such that the action
  of $\iota'$ restricted to $O_F$ is strict with respect to
  $\varphi_{0,R}: O_F \overset{\varphi_0}{\longrightarrow} O_{E'} \longrightarrow R$. 
  \end{enumerate}
  The category $\mathfrak{P}'_{r,R}$ is the category of formal $p$-divisible
  groups $X'$ with an $O_K$-action $\iota'$ such that the associated
  display $(\mathcal{P}', \iota')$   is an object of
  $\mathfrak{d}\mathfrak{P}'_{r,R}$.
  
  Let $r$ be special. We call $(\mathcal{P}, \iota)$ \emph{supersingular}\index{supersingular object of $\mathfrak{d}\mathfrak{P}_{r,R}$} if
  $(\mathcal{P}', \iota')$ satisfies the nilpotence condition.
  We denote the full subcategory of supersingular objects of
  $\mathfrak{d}\mathfrak{P}_{r,R}$ by
   $\mathfrak{d}\mathfrak{P}_{r,R}^{\rm ss}$.
\end{definition}
    
    \begin{theorem}\label{KottwitzC7p}  
      Let $R \in \Nilp_{O_{E'}}$ be such  that the ideal of nilpotent elements
      in $R$ is nilpotent.        
      The construction above defines the \emph{pre-contracting functor}\footnote{Later we will also have a \emph{contracting functor} $ \mathfrak{C}_{r,R}$, which explains our notation.}\index{pre-contracting functor}
      \index[NO]{CCA@$\mathfrak{C}'_{r,R}$}
\begin{equation*}\label{Cstrich1e}
  \mathfrak{C}'_{r,R}: \mathfrak{d}\mathfrak{P}_{r,R} \longrightarrow
  \mathfrak{d}\mathfrak{P}'_{r,R} 
  \end{equation*}
which commutes with arbitrary base change with respect to $R$. Furthermore, 
\begin{altenumerate}

\item if $r$ is banal, the functor $ \mathfrak{C}'_{r,R}$ is an equivalence of
categories. 

\item if $r$ is special and the ring $R$ is reduced, the functor $\mathfrak{C}'_{r,R}$
is an equivalence of categories.

\item if $r$ is special and $R$ is arbitrary,  $\mathfrak{C}'_{r,R}$
induces an equivalence of categories
\begin{equation*}\label{KottwitzC29e}
  \mathfrak{C}'_{r,R}: \mathfrak{d}\mathfrak{P}_{r,R}^{\rm ss} \longrightarrow
  \mathfrak{d}\mathfrak{P}_{r,R}^{\prime\,{\rm nilp}},
  \end{equation*}
where the right hand side is the full subcategory of nilpotent displays.
\index[NO]{PCH@$\mathfrak{d}\mathfrak{P}_{r,R}^{\prime\,{\rm nilp}}$}
\end{altenumerate}
Let   $r$ be special and $K/F$ be ramified. Let
$(\mathcal{P}, \iota) \in \mathfrak{d}\mathfrak{P}_{r,R}$ and let
$(\mathcal{P}', \iota')$ be its image by $\mathfrak{C}'_{r,R}$. Then
$(\mathcal{P}, \iota)$ satisfies $({\rm KC}_r)$ if and only if  
\begin{displaymath}
 \Trace_{R} (\iota'(\Pi) \; \mid \; P'/Q') = 0.  
  \end{displaymath}
    \end{theorem}
    Before proving this, we state a Corollary which we already proved in the
    construction of $\mathfrak{C}'_{r,R}$ above. 
    \begin{corollary}\label{C'onslopes1c}
      Let $k \in \Nilp_{O_{E'}}$ be a perfect field. Let
      $\mathcal{P} \in \mathfrak{d}\mathfrak{P}_{r,k}$ and let $\mathcal{P}'$
      be its image by the functor $\mathfrak{C}'_{r,R}$.
      \begin{altenumerate}
      \item[(1) ]  Let $r$ be banal and $K/F$  a field extension. Then  the display
        $\mathcal{P}$ is isoclinic of slope $1/2$ and $\mathcal{P}'$ is
        \'etale.
      \item[(2) ]  Let $K/F$ be split (and then $r$ is banal). Then  $\mathcal{P}$ decomposes into $\CP=\CP_1\oplus\CP_2$, where $\CP_1$ is isoclinic of slope  $\lambda$ and $\CP_2$ is isoclinic of slope $1-\lambda$. The
        number $\lambda$ depends only on $r$. The display $\mathcal{P}'$ is
        \'etale. 
      \item[(3) ]  Let $r$ be special. Then $\mathcal{P}$ is either isoclinic
        of slope $1/2$ (supersingular case) or it has the two slopes
        $1/2 - 1/2d$ and $1/2 + 1/2d$ with the same multiplicity. 
        In the first case $\mathcal{P}'$ is isoclinic of slope $1/2d$.
        In the second case it has the two slopes $0, 1/d$ with the same
        multiplicity.
        \end{altenumerate}\qed
      \end{corollary}
    \begin{proof}
      We still have to prove the claimed equivalences of categories.
      We begin with the case where $R$ is reduced. 

Let us consider first the case where   $r$ is special and $K/F$  unramified.
It is enough to invert the construction of the functor $\mathfrak{C}'_{r,R}$.
For any $\psi$,  $P_{\psi}/I(R)P_{\psi}$ 
is locally on $\Spec R$ a free $(O_K/pO_K) \otimes_{\kappa_K, \psi} R$-module
of rank 2, cf. Lemma \ref{Displaykristall1l}. Let $\mathcal{P}' = (P', Q', F', \dot{F}')$ be an object of
$\mathfrak{d}\mathfrak{P}'_{r,R}$. We define as follows an object
$\mathcal{P} = (P, Q, F, \dot{F})$ of $\mathfrak{d}\mathfrak{P}_{r,R}$ such
that $\mathcal{P}'$ is the image of $\mathcal{P}$ by the functor
$\mathfrak{C}'_{r,R}$.   
We set $P_{\psi} = P'_{\psi}$ for $\psi \in \Psi$, and for
$\psi \notin  \{\psi_0, \bar{\psi}_0 \}$ we set 
    \begin{equation}\label{KottwitzC27e}
      \quad Q_{\psi} = \pi^{a_{\psi}}P_{\psi} + I(R)P_{\psi}. 
    \end{equation}
    Since $P_{\psi} = Q'_{\psi}$, we have
    $\dot{F}'(I(R)P_{\psi}) \subset W(R)F'P_{\psi} \subset pP_{\psi \sigma}$.
    By (\ref{KottwitzC27e}) we find
    $\dot{F}'(Q_{\psi}) \subset \pi^{a_{\psi}} P_{\psi \sigma}$. Since $p$ is not
    a zero divisor in $W(R)$, the element $\pi \in O_K$ acts injectively on
    $P_{\psi_0}$. Therefore we may define
    \begin{displaymath}
      \dot{F} = \pi^{-a_{\psi}} \dot{F}': Q_{\psi} \longrightarrow P_{\psi \sigma}, \quad
      F = \pi^{-a_{\psi}} F': Q_{\psi} \longrightarrow P_{\psi \sigma}. 
      \end{displaymath}
    If $\psi \in \{\psi_0, \bar{\psi}_0 \}$ we consider the split homomorphism
    of $R$-modules
    \begin{equation}\label{KottwitzC28e}
\pi^{a_{\psi}}: P_{\psi}/I(R)P_{\psi} \longrightarrow P_{\psi}/I(R)P_{\psi}. 
      \end{equation}
It is split because $P_{\psi}/I(R)P_{\psi}$ is a free
$(O_K/pO_K) \otimes_{\kappa_K, \psi} R$-module. We set
\begin{displaymath}
Q_{\psi} = \pi^{a_{\psi}} Q'_{\psi} + I(R) P_{\psi}. 
\end{displaymath}
If we apply $\dot{F}'$ to the last equation we obtain that
$\dot{F}'(Q_{\psi}) \subset \pi^{a_{\psi}} P_{\psi \sigma}$. Indeed, because the
action of $O_F$ on $\mathcal{P}'$ is strict $\pi$ annihilates
$P_{\psi}/Q'_{\psi}$. We conclude that 
$F'(\pi P_{\psi})\subset F'(Q'_{\psi})=p\dot{F}'(Q_{\psi})\subset pP_{\psi\sigma}$
and therefore
$F'(P_{\psi}) \subset \pi^{e-1}P_{\psi} \subset \pi^{a_{\psi}} P_{\psi}$. 
This justifies the following definition: 
\begin{displaymath}
  \dot{F} := \pi^{-a_{\psi}} \dot{F}': Q_{\psi} \longrightarrow P_{\psi \sigma}, \quad
  F := \pi^{-a_{\psi}} F': P_{\psi} \longrightarrow P_{\psi \sigma}.   
\end{displaymath} 
It is obvious that we obtain a display $\mathcal{P} = (P, Q, F, \dot{F})$.
We need to verify that the condition $({\rm EC}_r)$ is satisfied. We check the conditions (2) and (3) of Proposition \ref{KottwitzC8p}. 
By definition of $\mathfrak{d}\mathfrak{P}'_{r,R}$, the $R$-module
$P'_{\psi}/Q'_{\psi}$ is annihilated by $\pi$. The kernel of (\ref{KottwitzC28e})
is $\pi^{e - a_{\psi}}P_\psi$ and therefore contained in $Q'_{\psi}/I(R)P_{\psi}$. The
image of the last module by (\ref{KottwitzC28e}) is therefore a direct summand
of $P_{\psi}/I(R)P_{\psi}$. This image is $Q_{\psi}/I(R)P_{\psi}$. Therefore
condition (2) holds. Moreover, we obtain an isomorphism
\begin{displaymath}
P_{\psi}/Q'_{\psi} \isoarrow \pi^{a_{\psi}} P_{\psi} + I(R)P_{\psi}/Q_{\psi}. 
\end{displaymath}
In particular, the last module is locally free of rank 1 and the action
of $\pi$ on this module coincides with multiplication by $\varphi_0(\pi)$
if $\psi = \psi_0$, resp., by $\bar{\varphi}_0(\pi)$ if
$\psi = \bar{\psi}_0$. Hence condition (3) holds. 

In the split case the same arguments hold but we need only the easy part
because $\psi_0$ and $\bar{\psi}_0$ don't exist. 

Next we consider the case where $r$ is special and $K/F$  ramified. Again we
reverse the construction of the functor
$\mathfrak{C}'_{r,R}$. Let $(P',Q',F',\dot{F}')$ be an object of
$\mathfrak{d}\mathfrak{P}'_{r,R}$. We associate to it as follows an object
$(P, Q, F, \dot{F}) \in \mathfrak{d}\mathfrak{P}_{r,R}$. We set
$P_{\psi} = P'_{\psi}$ for all $\psi \in \Psi$. Assume that
$\psi \neq \psi_0$. We have $Q'_{\psi} = P'_{\psi}$ because the action of
$O_{F^t}$ is strict.  We set
      \begin{displaymath}
        Q_{\psi} = \Pi^{e}P'_{\psi} + I(R)P'_{\psi}. 
      \end{displaymath}
It follows from Lemma \ref{Displaykristall1l} that $P_{\psi}/I(R)P_{\psi}$ 
is locally on $\Spec R$ a free $(O_K/pO_K) \otimes_{\kappa_K, \psi} R$-module.
Therefore $P_{\psi}/Q_{\psi}$ is a locally free $R$-module. From $F'P_{\psi} = p \dot{F}' P_{\psi}$, we find that
$\dot{F}'Q_{\psi} \subset \Pi^e P'_{\psi\sigma}$.
 Since $R$ is reduced, the ring $W(R)$ has no $p$-torsion. It
  follows that the
map $\Pi^e: P_{\psi \sigma} \longrightarrow \Pi^{e}P_{\psi \sigma}$ is bijective. 
Therefore we may define 
      \begin{displaymath}
        \dot{F} := \Pi^{-e}\dot{F}'\colon Q_{\psi} \longrightarrow P_{\psi \sigma}.
      \end{displaymath} 
It is clear that this map is a Frobenius-linear epimorphism. 
Next, we set 
      \begin{displaymath}
        P_{\psi_0} = P'_{\psi_0}, \quad Q_{\psi_0} = \Pi^{e-1}Q'_{\psi_0} +
        I(R)P'_{\psi_0}. 
      \end{displaymath}
      Since the action of $O_F$ on $\mathcal{P}'$ is strict, we find 
      \begin{equation}\label{KottwitzC24e}
        \Pi^{e+1}P'_{\psi_0} \subset \Pi^2P'_{\psi_0} \subset Q'_{\psi_0}.
      \end{equation}
We consider the split homomorphism of $R$-modules,
\begin{equation}\label{KottwitzC25e}
\Pi^{e-1}\colon P'_{\psi_0}/I(R) P'_{\psi_0} \longrightarrow P'_{\psi_0}/I(R) P'_{\psi_0}.
  \end{equation}
The kernel of this map is the image of $\Pi^{e+1}$. Therefore the kernel is
contained in $Q'_{\psi_0}/I(R) P'_{\psi_0}$. This implies that the image of
$Q'_{\psi_0}/I(R) P'_{\psi_0}$under \eqref{KottwitzC25e} is a direct summand of
$P'_{\psi_0}/I(R) P'_{\psi_0}$. Hence the cokernel 
$P_{\psi_0}/Q_{\psi_0}$ is a locally free $R$-module. We apply $F'$ to
(\ref{KottwitzC24e}) and obtain 
      \begin{displaymath}
        \Pi^2F'P'_{\psi_0} \subset F'Q'_{\psi_0} = p  \dot{F}'Q'_{\psi_0} \subset
        pP'_{\psi_0 \sigma}.
      \end{displaymath}
      Using this, we get
      \begin{displaymath}
        \dot{F}'Q_{\psi_0} = \dot{F}'\big(\Pi^{e-1}Q'_{\psi_0} +
        I(R)P'_{\psi_0}\big) \subset \Pi^{e-1}P'_{\psi_0 \sigma} +
        F'P'_{\psi_0} \subset \Pi^{e-1}P'_{\psi_0 \sigma}. 
      \end{displaymath}
It follows that the following definitions of maps
$Q_{\psi_0} \longrightarrow P_{\psi_0}$, resp. $P_{\psi_0} \longrightarrow P_{\psi_0}$,
make sense:
      \begin{displaymath}
        \dot{F} = (1/\Pi^{e-1}) \dot{F}' ,  \quad {F} = (1/\Pi^{e-1}) {F}'.
      \end{displaymath}
Therefore we have defined $\mathcal{P} = (P,Q,F,\dot{F})$. It is clear
that we obtain a display. We have to verify the condition $({\rm EC}_r)$.
Only $({\rm EC}_{\psi_0,r})$ is not completely obvious. We prove the
conditions of Proposition \ref{KottwitzC6p}. By the $R$-module homomorphism
(\ref{KottwitzC25e}), $P'_{\psi_0}/I(R)P'_{\psi_0}$ is mapped to the direct summand
$(\Pi^{e-1}P'_{\psi_0} + I(R)P'_{\psi_0})/I(R)P'_{\psi_0}$, and
$Q'_{\psi_0}/I(R) P'_{\psi_0}$ is mapped to the direct summand
$Q_{\psi_0}/I(R) P'_{\psi_0}$.
We obtain an isomorphism
\begin{displaymath}
P'_{\psi_0}/Q'_{\psi_0} \isoarrow (\Pi^{e-1}P_{\psi_0} + I(R)P_{\psi_0})/Q_{\psi_0}. 
  \end{displaymath}
Therefore by the strictness of the $O_F$-action, the right hand side is a
locally free $R$-module of rank $2$ and $\iota(\pi)$ acts on the right hand
side as $\varphi_0(\pi)$. These are exactly the conditions of Proposition
\ref{KottwitzC6p}. The rank condition is now obvious for
$(P, Q, F, \dot{F})$.

Finally, in the case where $r$ is banal, including the split case  (and $R$ is
reduced), we can reverse the functor
$\mathfrak{C}'_{r, R}$ using the arguments for banal $\psi$ given above. 

Now we consider  assertion (iii)  of Theorem \ref{KottwitzC7p}  when $R$ is
not reduced. It follows from Corollary \ref{C'onslopes1c} $(3)$ that
$\mathcal{P}$ is isoclinic of slope $1/2$ because $\mathcal{P}'$ is
nilpotent. Therefore we may apply Grothendieck-Messing for displays
Corollary \ref{GM1c}. We consider a surjective homomorphism
$S \longrightarrow R$ of $O_{E'}$-algebras and assume that
the kernel $\mathfrak{a}$ is endowed with a divided power structure. 

We define the category $\mathfrak{d}\mathfrak{P}_{r, S/R}$
\index[NO]{PCJ@$\mathfrak{d}\mathfrak{P}_{r, S/R}$} as the full
subcategory of the category of pairs $(\mathcal{P}_1, \iota_1)$ where
$\mathcal{P}_1$ is a $\mathcal{W}({S/R})$-display, 
cf. Example \ref{ex:crystdisp}, and where 
    \begin{displaymath}
      \iota_1: O_K \longrightarrow \End \mathcal{P}_1 
    \end{displaymath}
is an action such that the base change $(\mathcal{P}, \iota)$ of such a pair
by the morphism of frames $\mathcal{W}({S/R}) \longrightarrow \mathcal{W}(R)$
lies in the category $\mathfrak{d}\mathfrak{P}_{r, R}$. We also say that
$(\mathcal{P}_1, \iota_1)$ is a lift of $(\mathcal{P}, \iota)$ to a relative
display. By Theorem \ref{DspKristall1t}, the lift $(\mathcal{P}_1, \iota_1)$
is uniquely determined by $(\mathcal{P}, \iota)$ if $\mathcal{P}$ satisfies
the nilpotence condition.

In the same way we define the category $\mathfrak{d}\mathfrak{P}'_{r, S/R}$, cf.
Definition \ref{Kontrakt1d}. Then  the functor $\mathfrak{C}'_{r,R}$ of Theorem
\ref{KottwitzC7p} extends to a functor
    \begin{equation}\label{contract3e}
\mathfrak{C}'_{r, S/R}: \mathfrak{d}\mathfrak{P}_{r, S/R} \longrightarrow
      \mathfrak{d}\mathfrak{P}'_{r, S/R}. 
    \end{equation}
Indeed, the definition of $\mathfrak{C}'_{r, S/R}$ is essentially the same as
that of $\mathfrak{C}'_{r, R}$. We indicate it  in the case where $K/F$ is
unramified or split. By the $O_K$-action, we have for the relative display
$\mathcal{P}_1$ a decomposition,
    \begin{displaymath}
      P_1 = \oplus_{\psi} P_{1, \psi}, \quad Q_1 = \oplus_{\psi} Q_{1, \psi}. 
    \end{displaymath}
We are going to define a $\mathcal{W}({S/R})$-display
$\CP'_1=(P'_1, Q'_1, F'_1, \dot{F}'_1)$. We set $P'_{1, \psi} = P_{1, \psi}$ for all
$\psi \in \Psi$. Since $\mathcal{P}_1$ is a lifting of $\mathcal{P}$, we have
a natural isomorphism
\begin{equation}\label{KottwitzC30e}  
      P_{1,\psi}/Q_{1, \psi} = P_{\psi}/Q_{\psi}.
\end{equation}

If $\psi \notin \{\psi_0, \bar{\psi}_0 \}$,  we set $Q'_{1, \psi} = P'_{1, \psi}$.
By the condition $({\rm EC}_r)$ for $\mathcal{P}$ we conclude that
$\tilde{\mathbf{E}}_{A_{\psi}} P_{1,\psi} \subset Q_{1, \psi}$. Therefore we can
define
    \begin{equation}\label{KottwitzC31e}
      \begin{aligned}
               F'_1\colon& P'_{1,\psi} \longrightarrow P_{1,\psi \sigma},\,\,  & F'_1(x) =
        F_1(\tilde{\mathbf{E}}_{A_{\psi}} x), &\quad x\in P'_{1, \psi}\\
         \dot{F}'_1\colon& Q'_{1, \psi} \longrightarrow P'_{1,\psi \sigma},\,\, & \dot{F}'_1(x) =
        \dot{F}_1(\tilde{\mathbf{E}}_{A_{\psi}} x), &\quad  x \in Q'_{1,\psi} .
      \end{aligned}
    \end{equation}
In the split case this decribes $\mathcal{P}'_1$ already completely.  
Now we consider the case $\psi \in \{\psi_0, \bar{\psi}_0 \}$. Then we
define $Q'_{1,\psi}$ as the kernel of the map
  \begin{displaymath}
  P_{1, \psi} \longrightarrow P_{\psi}/Q_{\psi}
  \xra{\mathbf{E}_{A_{\psi}}} P_{\psi}/Q_{\psi}. 
  \end{displaymath}
This implies $\tilde{\mathbf{E}}_{A_{\psi}}Q'_{1,\psi} \subset Q_{1,\psi}$.
Therefore we can define
    \begin{equation}\label{KottwitzC32e}
      \begin{aligned}
        F'_1\colon & P'_{1,\psi} \longrightarrow P'_{1, \psi \sigma}, \,\,&
        F'_1(x) = F_1(\tilde{\mathbf{E}}_{A_{\psi}} x), &\quad x \in P'_{1,\psi},\\ 
        \dot{F}'_1\colon & Q'_{1,\psi} \longrightarrow P'_{1,\psi_0 \sigma},\,\, &
        \dot{F}'_1(y) = \dot{F}_1(\tilde{\mathbf{E}}_{A_{\psi}} y), &
        \quad y \in Q'_{1,\psi}. 
      \end{aligned}
    \end{equation}
 We then define $\mathcal{P}'_1 = (P'_1, Q'_1, F'_1, \dot{F}'_1)$, where
 $P'_1 = \oplus P'_{1,\psi}$ and $Q'_1 = \oplus Q'_{1,\psi}$. In the ramified case the same definition holds with slight modifications.
 
 The functor $\mathfrak{C}'_{r, S/R}$ defines
a natural isomorphism
\begin{equation}\label{KottwitzC34e}
\mathbb{D}_{\mathcal{P}}(S) \cong \mathbb{D}_{\mathcal{P}'}(S).  
\end{equation}
This relates deformations of $\mathcal{P}$ and deformations of its image
$\mathcal{P}'$ under $\mathfrak{C}_{r, R}$ since $\mathcal{P}$ and
$\mathcal{P}'$ are nilpotent.

Let $(\mathcal{P}, \iota) \in \mathfrak{d}\mathfrak{P}^{\rm ss}_{r, R}$.
It has a unique lift $\mathcal{P}_1 \in \mathfrak{d}\mathfrak{P}_{r, S/R}$. 
The image $\mathcal{P}'_1$ by the functor $\mathfrak{C}'_{r, S/R}$ is the
unique lift of $\mathcal{P}'$ to an object of
$\mathfrak{d}\mathfrak{P}'_{r, S/R}$, cf. Theorem \ref{DspKristall1t}. 

 Let us fix $\mathcal{P}$. Let $\CM$ be the set of all
 isomorphism classes of deformations of $(\mathcal{P}, \iota)$ to an object
 in $\mathfrak{d}\mathfrak{P}_{r, S}$. Let $\CM'$ be the set of
 isomorphism classes of deformations of $(\mathcal{P}', \iota')$ to an
 object of $\mathfrak{d}\mathfrak{P}'_{r, S/R}$. 
 We claim that the functor $\mathfrak{C}'_{r, S}$  defines a bijection,
 \begin{equation}\label{KottwitzC35e}
\mathfrak{C}'_{r, S}: \CM \longrightarrow \CM' .
 \end{equation}
 We indicate this  when $K/F$ is unramified.
 Let $\bar{Q} \subset \mathbb{D}_{\mathcal{P}}(R) = P/I(R)P$ be 
 the image of $Q$, i.e., the Hodge filtration. The set $\CM$ is identified
 with the set of liftings of $\bar{Q}$ to a direct summand
 $\bar{Q}_1 \subset \mathbb{D}_{\mathcal{P}}(S) = P_1/I(S)P_1$ which is a
 $O_K \otimes_{\mathbb{Z}_p} S$-submodule and such that the factor module
 satisfies the Eisenstein condition. The $O_K$-action gives a decomposition
 $\bar{Q}_1 = \oplus \bar{Q}_{1, \psi}$. For $\psi$ banal, we must have by
 Proposition \ref{KottwitzC8p} that
 \begin{equation}\label{KottwitzC36e}
\tilde{\mathbf{E}}_{A_{\psi}}\mathbb{D}_{\mathcal{P}}(S)_{\psi} = \bar{Q}_{1,\psi} .
   \end{equation}
 We note that the left hand side is a direct summand of
 $\mathbb{D}_{\mathcal{P}}(S)_{\psi}$ as an $S$-module. This follows from the fact
 that $P_{1, \psi}/I(S)P_{1, \psi}$ is a free module over
 $S[T]/\mathbf{E}_\psi S[T]$. Therefore, there is exactly one possibility
 to lift the $\psi$-component of the Hodge filtration. We consider now
 liftings of $\bar{Q}_{\psi}$ when $\psi$ is not banal. In this case the
 Eisenstein condition implies that
 \begin{displaymath}
   \mathbf{S}_{\psi} \tilde{\mathbf{E}}_{A_{\psi}} \mathbb{D}_{\mathcal{P}}(S)_{\psi} \subset
   \bar{Q}_{1,\psi} \subset \tilde{\mathbf{E}}_{A_{\psi}} \mathbb{D}_{\mathcal{P}}(S)_{\psi}. 
   \end{displaymath}
 By the freeness of $\mathbb{D}_{\mathcal{P}}(S)_{\psi}$ just mentioned, the
 multiplication by $\tilde{\mathbf{E}}_{A_{\psi}}$ gives an isomorphism
 \begin{displaymath}
  \tilde{ \mathbf{E}}_{A_{\psi}}: \mathbb{D}_{\mathcal{P}}(S)_{\psi}/ \mathbf{S}_{\psi}
   \mathbb{D}_{\mathcal{P}}(S)_{\psi} \cong
   \tilde{\mathbf{E}}_{A_{\psi}} \mathbb{D}_{\mathcal{P}}(S)_{\psi}/
   \mathbf{S}_{\psi}\tilde{\mathbf{E}}_{A_{\psi}} \mathbb{D}_{\mathcal{P}}(S)_{\psi}. 
 \end{displaymath}
 This shows that it is the same thing to lift $\bar{Q}_{\psi}$ to a direct
 summand $\bar{Q}_{1,\psi} \subset \mathbb{D}_{\mathcal{P}}(S)_{\psi}$ such that
 the Eisenstein condition is satisfied or to lift
 $\tilde{\mathbf{E}}_{A_{\psi}}^{-1} \bar{Q}_{\psi}$ to a direct summand $\bar{Q}'_{1,\psi}$
 such that $\mathbb{D}_{\mathcal{P}}(S)_{\psi}/\bar{Q}'_{1,\psi}$ is annihilated by
 $\mathbf{S}_{\psi}$. The last condition  means that the action of $O_F$ is strict with
 respect to $\varphi_0$, resp., $\bar{\varphi}_0$. In other words,
 \begin{displaymath}
   \bar{Q}'_{1} = \bar{Q}_{1,\psi_0} \oplus \bar{Q}_{1,\bar{\psi}_0} \oplus
   (\oplus_{\psi \neq \psi_0 \bar{\psi}_0} \mathbb{D}_{\mathcal{P}}(S)_{\psi}) 
   \end{displaymath}
 is a lift of the Hodge filtration
 $\bar{Q}' \subset \mathbb{D}_{\mathcal{P}'}(R) = P/I(R)P$ to
 a Hodge filtration $\bar{Q}'_{1} \subset P_1/I(S)P_1$ such that the action
 of $O_F$ is strict, i.e., the Hodge filtration $\bar{Q}'_{1}$ defines a
 point of $\CM'$. This shows that (\ref{KottwitzC35e}) is bijective because the
 functor $\mathfrak{C}'_{r, S}$ maps the Hodge filtration $\bar{Q}_{1,\psi}$ to
 $\mathbf{E}_{A_{\psi}}^{-1} \bar{Q}_{1,\psi}$ when $\psi$ is special by the
 definition (\ref{KottwitzC19e}). 
 We leave the ramified case to the reader. 

 Finally we prove assertion (i) of Theorem \ref{KottwitzC7p}, i.e., we assume
 that $r$ is banal. We begin with the case where $K/F$ is a field extension.
 Then $\mathcal{P}$ is by Corollary \ref{C'onslopes1c} $(1)$ of slope $1/2$.
 By (\ref{KottwitzC36e}) there is a unique way to lift the Hodge filtration
 and therefore the Grothendieck-Messing criterion implies that there is
 a unique way to lift $\mathcal{P}$ to an object
 $\mathcal{P}_1 \in \mathfrak{d}\mathfrak{P}_{r, S/R}$. On the other hand
 $\mathcal{P}'$ is \'etale. Therefore it lifts obviously uniquely, and (i) follows. 
 In the case where $K/F$ is split the same argument applies if $\mathcal{P}$
 is local. If not, we consider the decomposition
 $\mathcal{P} = \mathcal{P}_{\alpha} \oplus \mathcal{P}_{\beta}$ induced by
 $O_K = O_F \times O_F$. By Corollary \ref{C'onslopes1c} $(2)$, in each
 geometric point of $\Spec R$ one of the factors of this decomposition is
 isoclinic of slope $0$ and the other is isoclinic of slope $1$. 
 That $P_{\alpha}$ is \'etale means that the locally free module
 $P_{\alpha}/Q_{\alpha}$ is zero. This is true on an open and closed subset of
 $\Spec R$. Therefore we may assume without loss of generality that $P_{\alpha}$
 is \'etale. Then $\mathcal{P}_{\alpha}$ has a unique lift and
 $\mathcal{P}_{\beta}$ has a unique lift by Grothendieck-Messing. Since
 $\mathcal{P}'$ is \'etale it has also a unique lift. This completes the
 proof in the split case. 
    \end{proof} 
 
    \begin{corollary}\label{KottwitzC7c}
      Let $R \in \Nilp_{O_{E'}}$ be such  that the ideal of nilpotent elements in $R$ is nilpotent. We denote by $\mathfrak{P}^{\rm ss}_{r,R}$
      \index[NO]{PCK@$\mathfrak{P}^{\rm ss}_{r,R}$} the full
      subcategory of objects of $\mathfrak{P}_{r,R}$ whose displays lie in $\mathfrak{d}\mathfrak{P}^{\rm ss}_{r, R}$. Let
      $\mathfrak{P}^{' {\rm form}}_{r,R}$
      \index[NO]{PCM@$\mathfrak{P}^{' {\rm form}}_{r,R}$} be the full subcategory of
      $\mathfrak{P}'_{r,R}$ whose objects are formal $p$-divisible groups.
      Then $\mathfrak{C}'_{r, R}$  induces an equivalence of categories
     \begin{displaymath}
  \mathfrak{C}'_{r, R}: \mathfrak{P}^{\rm ss}_{r,R} \longrightarrow
  \mathfrak{P}^{' {\rm form}}_{r,R} .
       \end{displaymath}\qed
      \end{corollary}
\subsection{The contracting functor  in the case of a special CM-type}\label{ss:contrpol}
In this subsection, $r$ will denote a special CM-type. In this case, we will  compose the functor  $\mathfrak{C}'_{r, R}$  with the Ahsendorf functor.
    \begin{definition}\label{Kontrakt2d}
      Let $r$ be special.   Let
      $R \in \Nilp_{O_F}$. We denote by
      $\mathfrak{d}\mathfrak{R}_{R}$\footnote{The symbol $\mathfrak{R}$ is to remind us that this is a category of \emph{relative displays}.}
      \index[NO]{RCA@$\mathfrak{d}\mathfrak{R}_{R}$}the category of
      $\mathcal{W}_{O_F}(R)$-displays $\mathcal{P}_{\rm{c}}$ endowed with a
      homomorphism of $O_F$-algebras 
      \begin{displaymath}
\iota_{\rm{c}}: O_K \longrightarrow \End \mathcal{P}_{\rm{c}}, 
        \end{displaymath}
      such that $\mathcal{P}_{\rm{c}}$ is of height $4$ and dimension $2$.
      In the case where $K/F$ is unramified, we require moreover that
      $\Lie \mathcal{P}_{\rm{c}}$ is locally on $\Spec R$ a free
      $O_K \otimes_{O_F} R$-module of rank $1$. 
    \end{definition}
We note that in the ramified case, the $O_K \otimes_{O_F} R$-module
$\Lie \mathcal{P}_{\rm{c}}$ is in general not locally  free on $\Spec R$. 
    
\begin{definition}\label{def:contr}
Let $r$ be special. Let $R$ be a $O_{E'}$-algebra.  We regard $R$ as a $O_F$-algebra via
$  \varphi_{0,R}: O_F \xra{\varphi_0} O_{E} \longrightarrow R. 
$
The contracting functor \index{contracting functor}
\index[NO]{CCB@$\mathfrak{C}_{r,R}$}
\begin{equation*}\label{Kontrakt1e}  
        \mathfrak{C}_{r,R}: \mathfrak{d}\mathfrak{P}_{r,R} \longrightarrow
        \mathfrak{d}\mathfrak{R}_{R} 
\end{equation*}
is the composition of $\mathfrak{C}'_{r,R}$ with the Ahsendorf functor 
$\mathfrak{A}_{O_F/\mathbb{Z}_p,R}$. 
\end{definition}
\begin{theorem}\label{Kontrakt1p}
  Let $ r$ be special. Let $R \in \Nilp_{O_{E'}}$ be such that the ideal of
  nilpotent elements of $R$ is nilpotent. Then the functor $\mathfrak{C}_{r,R}$
  induces an equivalence of categories 
  \begin{displaymath}
        \mathfrak{C}_{r, R}: \mathfrak{d}\mathfrak{P}^{\rm ss}_{r,R} \longrightarrow
        \mathfrak{d}\mathfrak{R}^{\rm nilp}_{R}.
    \end{displaymath}
      {\it Here $ \mathfrak{d}\mathfrak{R}^{\rm nilp}_{R}$ denotes the full subcategory of nilpotent displays in $ \mathfrak{d}\mathfrak{R}_{R}$.}
      \index[NO]{RCB@$\mathfrak{d}\mathfrak{R}^{\rm nilp}_{R}$}
      \end{theorem}
\begin{proof}
This follows from Proposition \ref{KottwitzC7p} and Theorem \ref{Adorf1t}. 
  \end{proof}

\begin{remark}\label{CRF1r}
  Let $R = k$ be a perfect field with an $O_{E'}$-algebra structure. Then the
  construction of the functor $\mathfrak{C}_{r, k}$ 
  simplifies.

  We begin with the unramified case. Let
  $\mathcal{P} = (P,F,V) \in \mathfrak{d}\mathfrak{P}_{r,k}$ viewed as a
  Dieudonn\'e module. The display 
  $\mathcal{P}' = (P, F',V')$ is described after (\ref{KottwitzC22e}). Applying the Ahsendorf functor to it, 
  we obtain the image $\mathcal{P}_{\rm{c}} = (P_{\rm{c}}, F_{\rm{c}}, V_{\rm{c}})$ of $\mathcal{P}$ by
  the functor $\mathfrak{C}_{r, k}$. 
  The $P_0$ of (\ref{Adorf23e}) is in our case
  $P_{\rm{c}} = P_{\psi_0} \oplus P_{\bar{\psi}_0}$ and $V_{\rm{c}} = (V')^f$ is the $V_{\rm a}$
  of (\ref{Adorf23e}). We know that the restriction of $V'$ to $P_{\psi \sigma}$
  is
  \begin{displaymath}
V' = \pi^{-a_{\psi}} V: P_{\psi \sigma} \longrightarrow P_{\psi}  .
    \end{displaymath}
  We conclude that $(V')^f: P_{\psi_0} \longrightarrow P_{\bar{\psi}_0}$ is equal to
  $\pi^{-g_{\bar{\psi}_0}}V^f$ where
  \begin{equation}
    \begin{aligned}
g_{\bar{\psi}_0} & = a_{\bar{\psi}_0} + a_{\bar{\psi}_0 \sigma} + \ldots + 
a_{\bar{\psi}_0 \sigma^{f-1}}\\
 & = a_{\psi_0 \sigma^{-f}} + a_{\psi_0 \sigma^{-(f-1)}} + \ldots + a_{\psi_0 \sigma^{-1}}.
      \end{aligned}
  \end{equation}
In the same way $(V')^f: P_{\bar{\psi}_0} \longrightarrow P_{\psi_0}$ is equal to
$\pi^{-g_{\psi_0}}V^f$ where
\begin{equation}
g_{\psi_0} = a_{\psi_0} + a_{\psi_0 \sigma} + \ldots + a_{\psi_0 \sigma^{f-1}}. 
\end{equation}
From (\ref{furtherDr2unram1e}) we obtain $g_{\psi_0} + g_{\bar{\psi}_0} = ef-1$.
In summary, $P_{\rm{c}} = P_{\psi_0} \oplus P_{\bar{\psi}_0}$ as a
$W_{O_F}(k) = O_F \otimes_{O_{F^t}, \tilde{\psi}_0} W(k)$-module,  and $V_{\rm{c}}$ 
is given by the matrix
  \begin{equation}\label{P'Pol14e}
    \left(
    \begin{array}{cc}
      0 & \pi^{-g_{\psi_0}} V^f\\
      \pi^{-g_{\bar{\psi}_0}} V^f & 0
    \end{array}
    \right). 
  \end{equation}
Finally $F_{\rm{c}}$ is determined by the equation $F_{\rm{c}} V_{\rm{c}} = \pi$.
For instance,  $F_{\rm{c}}: P_{\bar{\psi}_0} \longrightarrow P_{\psi_0}$ is equal to
$(\pi^{g_{\bar{\psi}_0 }+1}/p^f)F^f$. We obtain a Dieudonn\'e module $(P_{\rm{c}}, F_{\rm{c}}, V_{\rm{c}})$
with respect to the perfect frame
\begin{equation}\label{P'Pol22e}
  \mathcal{W}_{O_F}(k) = (O_F \otimes_{O_{F^t}, \tilde{\psi}_0} W(k),
  \pi O_F \otimes_{O_{F^t}, \tilde{\psi}_0} W(k), k, F^f, F^f \pi^{-1}).
  \end{equation}

In the ramified case we have $P_{\rm{c}} = P_{\psi_0}$ as a module over 
$W_{O_F}(k) = O_F \otimes_{O_{F^t}, \tilde{\psi}_0} W(k)$. If we apply the Ahsendorf
functor to $\mathcal{P}'$, we obtain by (\ref{defP'1e})
\begin{equation}\label{CRF1e}
V_{\rm{c}} = \Pi^{-ef+1} V^f: P_{\rm{c}} \longrightarrow P_{\rm{c}}.  
\end{equation}
$F_{\rm{c}}$ is determined by the equation $F_{\rm{c}} V_{\rm{c}} = \pi$, i.e.,
$F_{\rm{c}} = -(\Pi^{ef +1}/p^f)F^f$. We obtain a Dieudonn\'e module $(P_{\rm{c}}, F_{\rm{c}}, V_{\rm{c}})$
for the frame (\ref{P'Pol22e})
  \end{remark}

\label{page:theta}We next add polarizations to the picture. 
We set $\mathbf{t}(a) = \Trace_{F/\mathbb{Q}_p} \vartheta^{-1} a$
\index[NO]{TAA@$\mathbf{t}$} where
$\vartheta$ \index[NO]{ZZHB@$\vartheta$}is the different of $F/\mathbb{Q}_p$.
\begin{proposition}\label{P'Pol1p}
  Let $ r$ be special. 
  Let $R \in \Nilp_{O_{E'}}$. Let $(\mathcal{P}_1, \iota_1)$ and
  $(\mathcal{P}_2, \iota_2)$ be objects of $\mathfrak{d}\mathfrak{P}_{r,R}$.
  Let 
  $(\mathcal{P}'_1, \iota'_1)$ and $(\mathcal{P}'_2, \iota'_2)$ be their images by the functor $\mathfrak{C}'_{r, R}$. 
  Assume given a bilinear form of displays 
  \begin{displaymath}
\beta: \mathcal{P}_1 \times \mathcal{P}_2 \longrightarrow \mathcal{P}_{m,R},  
  \end{displaymath}
where $\mathcal{P}_m$ is the multiplicative display of $\mathcal{W}(R)$. 
Assume that $\beta$ is anti-linear for the $O_K$-actions $\iota_1$, resp. $\iota_2$, i.e,  
\begin{equation}\label{P'Pol1e} 
  \beta(\iota_1(a) x_1, x_2) = \beta(x_1, \iota_2(\bar{a}) x_2), \quad
  x_1 \in P_1, \; x_2 \in P_2, \,  a \in O_K .
  \end{equation}
Define
\begin{displaymath}
\tilde{\beta}: P_1 \times P_2 \longrightarrow O_F \otimes_{\mathbb{Z}_p} W(R) 
  \end{displaymath}
by the equation
\begin{displaymath}
  \mathbf{t}(\xi\tilde{\beta}(x_1, x_2)) = \beta(\xi x_1, x_2),  \quad x_1 \in P_1, \; x_2 \in P_2,\;
  \xi \in O_F \otimes_{\mathbb{Z}_p} W(R). 
\end{displaymath}
Then $\tilde{\beta}$ is a $O_F$-bilinear form of displays, 
\begin{displaymath}
  \tilde{\beta}: \mathcal{P}'_1 \times \mathcal{P}'_2 \longrightarrow
  \mathcal{L}_R,   
\end{displaymath}
where $\mathcal{L}_R$ is the Lubin-Tate display associated to the local field
$F$ and the algebra structure $\varphi_0: O_F \longrightarrow O_{E'} \longrightarrow R$,
cf. Definition \ref{LTD1d}. Furthermore, $\tilde{\beta}$ is anti-linear for the $O_K$-actions $\iota_1$, resp. $\iota_2$. 
\end{proposition}
\begin{proof}  To avoid a conflict with the present notations, we adapt some of the notation of section \ref{ss:tltd} to
  our  situation. What was $K$ in section \ref{ss:tltd} is now $F$. We set 
   $L^t = \psi_0(F^t) \subset E'$. We write the polynomials of that 
  section as follows:

\begin{displaymath}
  \tilde{\mathbf{E}}_{F,\psi}(Z) =
  \prod_{\chi:F \longrightarrow \tilde{E},\; \chi_{|F^t} = \psi}
  (Z - [\chi(\pi)]) \in W(O_{L^t})[Z]. 
  \end{displaymath}
We stress that here $\psi$ denotes an embedding of $F^t$ into $E'$, not as elsewhere in this section an embedding of $K^t$ into $E'$. For $\psi = \psi_0$, we consider the decomposition 
$\tilde{\mathbf{E}}_{F,\psi_0}(Z) = (Z- [\varphi_0(\pi)])\cdot
\tilde{\mathbf{E}}_{F,0}(Z)$ 
in $W(\varphi_0(F))[Z]$.   
In particular all of these polynomials lie in $W(O_{E'})[Z]$.

Let $M$ be an $O_F \otimes_{O_{F^t},\psi} W(R)$-module. Then we write by our
convention 
\begin{displaymath}
\tilde{\mathbf{E}}_{F,\psi} m = \tilde{\mathbf{E}}_{F,\psi}(\pi \otimes 1) m, 
\end{displaymath}
where
$\tilde{\mathbf{E}}_{F,\psi}(\pi \otimes 1) \in O_F \otimes_{O_{F^t},\psi} W(R)$ 
is the evaluation at $\pi \otimes 1$ in this $W(O_{L^t})$-algebra.  

We  first consider the assertion of Proposition \ref{P'Pol1p} in the ramified case. We have the decomposition
$P_{i} = \oplus_{\psi} P_{i,\psi}$. By (\ref{P'Pol1e}), we find
$\beta(P_{1,\psi}, P_{2,\psi'}) = 0$ for $\psi \neq \psi'$. 
We consider the restrictions of our
bilinear forms 
  \begin{equation*}
  \begin{aligned}
    \beta_{\psi}:& P_{1,\psi} \times P_{2,\psi} \longrightarrow W(R) \\
    \tilde{\beta}_{\psi}:& P_{1\psi} \times P_{2,\psi} \longrightarrow
    O_F \otimes_{O_{F^t}, \psi} W(R) .
        \end{aligned}
  \end{equation*}
\begin{lemma}\label{P'Pol1l}
  Let $K/F$ be ramified. Then:
  \begin{equation*}
  \begin{aligned}
    \tilde{\beta}_{\psi}(\tilde{\mathbf{E}}_{A_{\psi}} x_1,
    \tilde{\mathbf{E}}_{A_{\psi}} x_2) &= \tilde{\mathbf{E}}_{F,\psi}
    \tilde{\beta}_{\psi}(x_1, x_2) , \quad
    &x_1 \in P_{1,\psi}, \; x_2 \in P_{2,\psi},\; &\psi \neq \psi_0\\
    \tilde{\beta}_{\psi_0}(\tilde{\mathbf{E}}_{A_{\psi_0}} x_1,
    \tilde{\mathbf{E}}_{A_{\psi_0}} x_2) &= \tilde{\mathbf{E}}_{F,0}
    \tilde{\beta}_{\psi}(x_1, x_2) , \quad &x_1 \in P_{1,\psi_0}, \; x_2 \in P_{2,\psi_0}&.
    \end{aligned}
  \end{equation*}
\end{lemma}
\begin{proof}
  We can restrict ourselves to the case where $R$ is a $O_{\tilde{E}}$-algebra.
  Then we obtain 
  \begin{equation*}
    \begin{aligned} 
    \tilde{\beta}_{\psi}\big((\Pi \otimes 1 - 1 \otimes [\varphi(\Pi)])x_1, x_2\big) &= 
    \tilde{\beta}_{\psi}\big(x_1,
    (-\Pi \otimes 1 - 1 \otimes [\varphi(\Pi)])x_2\big)\\
    &= -
    \tilde{\beta}_{\psi}\big(x_1, (\Pi\otimes 1 - 1 \otimes
          [\bar{\varphi}(\Pi)])x_2\big). 
      \end{aligned}
  \end{equation*}
  In the case were $\psi \neq \psi_0$ we deduce
  \begin{displaymath}
    \tilde{\beta}_{\psi}(\tilde{\mathbf{E}}_{A_{\psi}}x_1, x_2) = 
    (-1)^{e}\tilde{\beta}_{\psi}(x_1, \tilde{\mathbf{E}}_{B_{\psi}}x_2) .
  \end{displaymath}
  We find
   \begin{equation*}
    \begin{aligned}
    \tilde{\mathbf{E}}_{B_{\psi}}(\Pi \otimes 1) \cdot
    \tilde{\mathbf{E}}_{A_{\psi}}(\Pi \otimes 1) &=
    \prod_{\varphi \in A_{\psi}} (\Pi \otimes 1 - 1 \otimes [\varphi(\Pi)])
    (\Pi \otimes 1 - 1 \otimes [\bar{\varphi}(\Pi)]) \\ 
   & = \prod_{\varphi \in A_{\psi}} (\Pi \otimes 1 - 1 \otimes [\varphi(\Pi)])
    (\Pi \otimes 1 + 1 \otimes [\varphi(\Pi)]) \\
    &=
    \prod_{\varphi \in A_{\psi}} (-\pi \otimes 1 + 1 \otimes [\varphi(\pi)])   
    = (-1)^{e} \tilde{\mathbf{E}}_{F,\psi}(\pi \otimes 1). 
      \end{aligned}
  \end{equation*}
  Therefore we obtain
  \begin{equation*}
   \begin{aligned}
    \tilde{\beta}_{\psi}(\tilde{\mathbf{E}}_{A_{\psi}}x_1,
    \tilde{\mathbf{E}}_{A_{\psi}}x_2) &=
    (-1)^e\tilde{\beta}_{\psi}(x_1,
    \tilde{\mathbf{E}}_{B_{\psi}}\tilde{\mathbf{E}}_{A_{\psi}}x_2) 
    = \tilde{\beta}_{\psi}(x_1, \tilde{\mathbf{E}}_{F,\psi} x_2) \\
    &= \tilde{\mathbf{E}}_{F,\psi} \tilde{\beta}_{\psi}(x_1, x_2), 
    \end{aligned}
  \end{equation*}
  which finishes the proof for $\psi\neq\psi_0$. 
  
  We turn now to the case $\psi_0$. The polynomials
  $\tilde{\mathbf{E}}_{A_{\psi_0}}$,   $\tilde{\mathbf{E}}_{B_{\psi_0}}$, and
  $\tilde{\mathbf{E}}_{F,0}$ are of degree $e-1$. The same computations
  yield for $x_1 \in P_{1,\psi_0}$ and $ x_2 \in P_{2,\psi_0}$,
  \begin{equation*}
    \begin{aligned}
      \tilde{\beta}_{\psi_0}(\tilde{\mathbf{E}}_{A_{\psi_0}}x_1, x_2) &= 
      (-1)^{e-1}\tilde{\beta}_{\psi_0}(x_1, \tilde{\mathbf{E}}_{B_{\psi_0}}x_2),
      \\
      \tilde{\mathbf{E}}_{B_{\psi_0}}(\Pi \otimes 1) \cdot
    \tilde{\mathbf{E}}_{A_{\psi_0}}(\Pi \otimes 1) &= 
     (-1)^{e-1} \tilde{\mathbf{E}}_{F,0}(\pi \otimes 1). 
      \end{aligned}
  \end{equation*}
  The  assertion for $\psi_0$  follows as before. 
\end{proof}
We continue with the proof of  Proposition \ref{P'Pol1p} in the ramified case.
We begin by showing that
  \begin{equation}\label{P'Pol3e}
\tilde{\beta}_{\psi}(Q'_{1,\psi}, Q'_{2,\psi}) \subset Q_{\mathcal{L}, \psi}.
  \end{equation}
  This is trivial for $\psi \neq \psi_0$. Let $y_1 \in Q'_{1,\psi_0}$ and
  $y_2 \in Q'_{2,\psi_0}$.
  By Lemma \ref{LTD5l}, the inclusion (\ref{P'Pol3e}) is equivalent to 
  \begin{displaymath}
    \tilde{\mathbf{E}}_{F,0}\tilde{\beta}_{\psi_0}(y_1, y_2) \in
    O_F \otimes_{O_{F^t},\psi_0} I(R)
  \end{displaymath}
  By  Lemma \ref{P'Pol1l} we find
  \begin{displaymath}
    \tilde{\mathbf{E}}_{F,0}\tilde{\beta}_{\psi_0}(y_1, y_2)  =
    \tilde{\beta}_{\psi_0}(\tilde{\mathbf{E}}_{A_{\psi_0}}y_1, 
    \tilde{\mathbf{E}}_{A_{\psi_0}}y_2) .
  \end{displaymath}
  The elements $u_1 = \tilde{\mathbf{E}}_{A_{\psi_0}}y_1$, resp., 
  $u_2 = \tilde{\mathbf{E}}_{A_{\psi_0}}y_2$, lie, by the definition of
  $Q'_{1,\psi_0}$, in $Q_{1,\psi_0}$, resp., by the definition of
  $Q'_{2,\psi_0}$, in $Q_{2,\psi_0}$.
  But for arbitrary elements
  $u_1 \in Q_{1,\psi_0}$ and $u_2 \in Q_{2,\psi_0}$, we have
  $\beta(u_1, u_2) \in I(R)$.
  By the definition of $\tilde{\beta}_{\psi_0}$,  we find
  \begin{displaymath}
    \Trace_{F/F^t} \vartheta^{-1} \xi \tilde{\beta}_{\psi_0}(u_1,u_2) =
    \beta(\xi u_1, u_2) \in I(R), 
    \end{displaymath}
  for all $\xi \in O_F \otimes_{O_{F^t},\psi_0} W(R)$. But this implies
  $\tilde{\beta}_{\psi_0}(u_1, u_2) \in O_F \otimes_{O_{F^t},\psi_0} I(R)$, as
  desired.

  Finally we have to check for $y_1 \in Q'_{1,\psi}$ and $y_2 \in Q'_{2,\psi}$
  that 
  \begin{displaymath}
    \tilde{\beta}_{\psi}(\dot{F}'y_1, \dot{F}'y_2) =
    \dot{F}_{\mathcal{L}}\tilde{\beta}_{\psi}(y_1, y_2). 
    \end{displaymath}
  If $\psi \neq \psi_0$ we find for the left hand side
    \begin{displaymath} 
      \tilde{\beta}_{\psi}\big(\dot{F}(\tilde{\mathbf{E}}_{A_{\psi}})y_1,
      \dot{F}(\tilde{\mathbf{E}}_{A_{\psi}})y_2\big) =
      ~^{\dot{F}}\big(\tilde{\mathbf{E}}_{F,\psi}
      \tilde{\beta}_{\psi}(y_1,y_2)\big) =
      \dot{F}_{\mathcal{L}}\big(\tilde{\beta}_{\psi}(y_1,y_2)\big).  
    \end{displaymath}
For $\psi_0$ we obtain 
    \begin{displaymath} 
      \tilde{\beta}_{\psi_0}\big(\dot{F}(\tilde{\mathbf{E}}_{A_{\psi_0}})y_1,
      \dot{F}(\tilde{\mathbf{E}}_{A_{\psi_0}})y_2\big) =
      ~^{\dot{F}}\big(\tilde{\mathbf{E}}_{F,0} 
      \tilde{\beta}_{\psi_0}\big(y_1,y_2)\big) =
      \dot{F}_{\mathcal{L}}\big(\tilde{\beta}_{\psi_0}(y_1,y_2)\big). 
    \end{displaymath}
This ends the proof of Proposition \ref{P'Pol1p} in  the ramified case. 

Now we consider the unramified case. We consider the decomposition
(\ref{OWdecomp1e}). Let us denote by $e_{\psi}$ the idempotents corresponding
to this decomposition. The conjugation of $K/F$ maps $e_{\psi}$ to
$e_{\bar{\psi}}$.
We consider the corresponding decompositions $P_i = \oplus P_{i, \psi}$, for
$i = 1,2$. We obtain that $\beta(P_{1, \psi_1}, P_{2, \psi_2}) = 0$ for
$\psi_2 \neq \bar{\psi_1}$, and
\begin{displaymath}
  \tilde{\beta}(P_{1, \psi}, P_{2, \bar{\psi}}) \subset
  O_F \otimes_{O_{F^t}, \tilde{\psi}} W(R). 
\end{displaymath}
We note that there are natural isomorphisms
\begin{equation}\label{P'Pol4e} 
O_K \otimes_{O_{K^t}, \tilde{\psi}} W(R) \cong O_F \otimes_{O_{F^t}, \tilde{\psi}} W(R) 
\cong O_K \otimes_{O_{K^t}, \widetilde{\bar{\psi}}} W(R).
  \end{equation}
Therefore $\tilde{\beta}$ induces an
$O_F \otimes_{O_{F^t}, \tilde{\psi}} W(R)$-bilinear form
\begin{displaymath}
  \tilde{\beta}_{\psi}: P_{1, \psi} \times P_{2, \bar{\psi}} \longrightarrow
  O_F \otimes_{O_{F^t}, \tilde{\psi}} W(R). 
\end{displaymath}
With the identification (\ref{P'Pol4e}),  we have
$\tilde{\mathbf{E}}_{A_{\bar{\psi}}}(\pi \otimes 1) =
\tilde{\mathbf{E}}_{B_{\psi}}(\pi \otimes 1) \in
O_F \otimes_{O_{F^t}, \tilde{\psi}} W(R)$.
The analogue of Lemma \ref{P'Pol1l} is 
  \begin{equation}\label{P'Pol5e} 
  \begin{aligned}
  \tilde{\beta}_{\psi}(\tilde{\mathbf{E}}_{A_{\psi}}x_1,
  \tilde{\mathbf{E}}_{A_{\bar{\psi}}}x_2) &= \tilde{\mathbf{E}}_{F,\psi}
  \tilde{\beta}_{\psi}(x_1,x_2), &\quad x_1 \in P_{1, \psi}, \; x_2 \in P_{2, \bar{\psi}},
  \; \psi \neq \psi_0, \bar{\psi}_0,\\
  \tilde{\beta}_{\psi}(\tilde{\mathbf{E}}_{A_{\psi}}x_1,
  \tilde{\mathbf{E}}_{A_{\bar{\psi}}}x_2) &= \tilde{\mathbf{E}}_{F,0}
  \tilde{\beta}_{\psi}(x_1,x_2), &\quad x_1 \in P_{1, \psi}, \; x_2 \in P_{2, \bar{\psi}},
  \; \psi = \psi_0, \bar{\psi}_0. 
  \end{aligned}
  \end{equation} 
  Here we recall again the notation introduced in the beginning of the proof: to be very precise, the expression $\tilde{\mathbf{E}}_{F,\psi}$ should be written as $\tilde{\mathbf{E}}_{F,\psi_{|{F^t}}}$.  These identities follow from  the identities
  \begin{displaymath}
    \tilde{\mathbf{E}}_{A_{\psi}} \tilde{\mathbf{E}}_{B_{\psi}} =
    \begin{cases}
      \tilde{\mathbf{E}}_{F,\psi}, \quad \psi \neq \psi_0, \bar{\psi}_0,\\ 
      \tilde{\mathbf{E}}_{F,0}, \quad \psi = \psi_0, \bar{\psi}_0.
      \end{cases}
  \end{displaymath}
  We need to check
  \begin{equation}\label{unrincl}
\tilde{\beta}_{\psi}(Q'_{1,\psi}, Q'_{2,\bar{\psi}}) \subset Q_{\mathcal{L}, \psi}.
    \end{equation}
It suffices to consider the case $\psi = \psi_0$. 
 By Lemma \ref{LTD5l}, the inclusion (\ref{unrincl}) is equivalent to
  \begin{displaymath}
    \tilde{\mathbf{E}}_{F,0}\tilde{\beta}_{\psi_0}(y_1, y_2) \in
    O_F \otimes_{O_{F^t},\psi_0} I(R), \quad  y_1 \in Q'_{1,\psi_0}, \;
    y_2 \in Q'_{2,\bar{\psi_0}}. 
  \end{displaymath}
  But, as in the ramified case, this is an immediate consequence of
  (\ref{P'Pol5e}).
  Finally we have to check that for $y_1 \in Q'_{1,\psi}$ and
  $y_2 \in Q'_{2,\bar{\psi}}$ 
  \begin{displaymath}
    \tilde{\beta}_{\psi}(\dot{F}'y_1, \dot{F}'y_2) =
    \dot{F}_{\mathcal{L}}\tilde{\beta}_{\psi}(y_1, y_2). 
    \end{displaymath}
  For this we can repeat the last two formulas in the proof of the ramified
  case.
\end{proof}
Let $R \in \Nilp_{O_F}$. Let $(\mathcal{P}_1, \iota_1)$ be an object of
$\mathfrak{d}\mathfrak{P}'_{r,R}$. We denote by 
$(\mathcal{P}^{\Delta}_1, \iota^{\Delta}_1)$ the conjugate Faltings dual.\index{conjugate Faltings dual} \index[NO]{PBF@$(\mathcal{P}^{\Delta}_1, \iota^{\Delta}_1)$}
It is defined from the Faltings dual exactly as the conjugate dual from
the dual. 
 
\begin{corollary}\label{P'Pol1c}
  Let $ r$ be special.  Let $R \in \Nilp_{O_{E'}}$ be such that the ideal of nilpotent
  elements of $R$ is nilpotent. We regard $R$ as an
  $O_F$-algebra via $\varphi_0$. Let $(\mathcal{P}, \iota)$ be an object of $\mathfrak{d}\mathfrak{P}_{r,R}$
  and let $(\mathcal{P}', \iota') \in \mathfrak{d}\mathfrak{P}'_{r,R}$ be its
  image under the pre-contracting  functor $\mathfrak{C}'_{r, R}$. 
  Then the image of the conjugate dual $(\mathcal{P}^{\wedge}, \iota^{\wedge})$
  under $\mathfrak{C}'_{r, R}$ is the  conjugate Faltings  dual
  $((\mathcal{P}')^{\Delta}, (\iota')^{\Delta})$, cf. Proposition \ref{LTD7p}. 

  With the notation of Proposition \ref{P'Pol1p}, assume that
  $\mathcal{P}_1^{\wedge}$ and $\mathcal{P}_2$ are in
  $\mathfrak{d}\mathfrak{P}^{\rm ss}_{r,R}$. Then the canonical map
  \begin{displaymath}
    \Bil_{O_K\text{\rm -anti-linear}}(\mathcal{P}_1 \times \mathcal{P}_2,
    \mathcal{P}_{m,R}) \longrightarrow
    \Bil_{O_K\text{\rm-anti-linear}}(\mathcal{P}'_1 \times \mathcal{P}'_2, \mathcal{L}_R)
  \end{displaymath}
  is bijective. Here these sets of bilinear forms $\Bil$ are meant
  as in Proposition \ref{P'Pol1p}. 
\end{corollary}
\begin{proof}
 We apply Proposition \ref{P'Pol1p} to the canonical
bilinear form
$\beta_{\rm can}: \mathcal{P} \times \mathcal{P}^{\wedge}
\longrightarrow \mathcal{P}_{m,R}$ and obtain
\begin{displaymath}
  \tilde{\beta}_{\rm can}: \mathcal{P}' \times (\mathcal{P}^{\wedge})' \longrightarrow
  \mathcal{L}_R.
\end{displaymath}
By Proposition \ref{LTD9p},  we obtain a morphism of displays
\begin{equation}\label{P'Pol6e}
\varkappa: (\mathcal{P}^{\wedge})' \longrightarrow (\mathcal{P}')^{\Delta}.  
\end{equation}
By definition, $\tilde{\beta}_{\rm can}$ is given by a perfect
$O_F \otimes W(R)$-bilinear form
\begin{displaymath}
P \times P^{\ast} \longrightarrow O_F \otimes W(R).  
  \end{displaymath}
(Recall that $P^{\ast} = \Hom_{W(R)}(P, W(R))$.) We obtain an isomorphism
\begin{displaymath}
P^{\ast} \isoarrow \Hom_{O_F \otimes W(R)}(P, O_F \otimes W(R)). 
\end{displaymath}
But this says exactly that the map which $\varkappa$ induces on the
''$P$-components'' of the displays (\ref{P'Pol6e}) is an isomorphism.
It is elementary to see that a morphism of displays
$\varkappa: \mathcal{P}_1 \longrightarrow \mathcal{P}_2$ which induces a
$W(R)$-module isomorphism $P_1 \longrightarrow P_2$ is an isomorphism of displays.

Finally we prove the bijectivity of the last map in the corollary. The left
hand side is, by (\ref{dual4e}),
\begin{displaymath}
\Hom_{\mathfrak{d}\mathfrak{P}_{r,R}} (\mathcal{P}_1, (\mathcal{P}_2)^{\wedge}). 
\end{displaymath}
This group is, by (iii) of Theorem \ref{KottwitzC7p}, equal to
\begin{displaymath}
  \Hom_{\mathfrak{d}\mathfrak{P}'_{r,R}}(\mathcal{P}'_1, (\mathcal{P}_2^{\wedge})') 
  \cong
  \Hom_{\mathfrak{d}\mathfrak{P}'_{r,R}}(\mathcal{P}'_1, (\mathcal{P}'_2)^{\Delta})
  \cong
  \Bil_{O_K\text{\rm -anti-linear}}(\mathcal{P}'_1 \times \mathcal{P}'_2, \mathcal{L}_R).
  \end{displaymath}
\end{proof}
We now combine the last corollary with Theorem \ref{LTD10p}. 
\begin{theorem}\label{P'Pol2t}
Let $ r$ be special. Let $R \in \Nilp_{O_{E'}}$ be such that the ideal of nilpotent elements of $R$ is nilpotent. 
Let $(\mathcal{P}_1, \iota_1)$ and $(\mathcal{P}_2, \iota_2)$ be objects of
$\mathfrak{d}\mathfrak{P}_{r,R}$, with  images  $(\mathcal{P}'_1, \iota'_1)$ and
$(\mathcal{P}'_2, \iota'_2)$ under the pre-contracting functor
$\mathfrak{C}'_{r,R}$, cf. Proposition \ref{KottwitzC7p}. Since the actions
$\iota'_i$ restricted to $O_F$ are strict with respect to
$\varphi_0: O_F \longrightarrow O_{E'} \longrightarrow R$,  the
Ahsendorf functor $\mathfrak{A}_{O_F/\mathbb{Z}_p, R}$ may be applied to them. For $i=1, 2$, let $\mathcal{P}_{i,c}=\mathfrak{A}_{O_F/\mathbb{Z}_p, R}(\CP'_i)$, $i=1, 2$, with its $O_F$-algebra homomorphism
\begin{displaymath}
\iota_{i, c}: O_K \longrightarrow \End_{\mathcal{W}_{O_F}(R)} \mathcal{P}_{i,c}.  
  \end{displaymath}
If $\mathcal{P}_1^{\wedge}$ and $\mathcal{P}_2$ are in
$\mathfrak{d}\mathfrak{P}^{\rm ss}_{r,R}$, then the  natural homomorphism 
\begin{displaymath}
    \Bil_{O_K\text{\rm-anti-linear}}(\mathcal{P}_1 \times \mathcal{P}_2,
    \mathcal{P}_{m,R}) \longrightarrow
 \Bil_{O_K\text{\rm -anti-linear}}(\mathcal{P}_{1, c} \times \mathcal{P}_{2, c},
  \mathcal{P}_{m, \mathcal{W}_{O_F}(R)}(\pi^{ef}/p^f)). 
\end{displaymath}
is a bijection. 

{\rm The twist $\mathcal{P}_{m, \mathcal{W}_{O_F}(R)}(\pi^{ef}/p^f)$ of the multiplicative
display is defined in Example \ref{def:twistdisplay}. More precisely, this is the twist by the image of $(\pi^{ef}/p^f)$ under  the canonical map
$O_F \longrightarrow W_{O_F}(R)$.  }
\end{theorem}
\begin{proof}
This follows from Corollary \ref{P'Pol1c} and Theorem  \ref{LTD10p}. 
\end{proof}

\begin{remark}\label{r:breve} Let $\breve{E} \subset \hat{\bar{\mathbb{Q}}}_p$ be the completion of the
maximal unramified extension of the reflex field $E$ of $r$. We extend
$\varphi_0: O_F \longrightarrow O_{\breve{E}}$ to an embedding 
$\breve\varphi_0: O_{\breve{F}} \longrightarrow O_{\breve{E}}$.
\index[NO]{ZZSC@$\breve\varphi_0$}We denote by
$\tau \in \Gal(\breve{F}/F)$ the Frobenius automorphism. We apply the
definition of $\eta_0$ after Definition \ref{LTD3d},
\begin{equation}\label{P'Pol8e}
\tau(\eta_0) \eta_0^{-1} = \pi^{e}/p, \quad \eta_0 \in O_{\breve{F}}^{\times}. 
  \end{equation}

Let $R \in \Nilp_{O_{\breve{E}}}$. Via $\breve\varphi_0$ we consider $R$ as an
$O_{\breve{F}}$-algebra. Therefore $\eta_{0,R}$ is defined, and  multiplication
by $\eta_{0,R}^f$ defines an isomorphism
\begin{equation}\label{P'Pol7e1}
  \mathcal{P}_{m, \mathcal{W}_{O_F}(R)}(\pi^{ef}/p^f) \isoarrow
  \mathcal{P}_{m, \mathcal{W}_{O_F}(R)}, 
  \end{equation}
cf. (\ref{P'Pol7e}). Therefore,  if $R \in \Nilp_{O_{\breve{E}}}$, we can ignore the twist by $(\pi^{e}/p)$ in  Theorem \ref{P'Pol2t}. 
\end{remark}

We recall the definition of polarized CM-pairs  
   $\mathfrak P^{\rm pol}_{r,S}$, cf. Definition \ref{KatCMtriple1d}. We also introduce the analogous category of polarized objects of
$\mathfrak{d}{\mathfrak R}_{R}$, as follows. 
\begin{definition}\label{P'Pol3d} 
  Let $R \in \Nilp_{O_F}$. We denote by $\mathfrak{d}\mathfrak{R}^{\rm pol}_{R}$
  \index[NO]{RCC@$\mathfrak{d}\mathfrak{R}^{\rm pol}_{R}$} 
  the category of triples $(\mathcal{P}_{\rm{c}}, \iota_{\rm{c}}, \beta_{\rm{c}})$ where
  $(\mathcal{P}_{\rm{c}}, \iota_{\rm{c}}) \in \mathfrak{d}{\mathfrak R}_{R}$ (cf.
  Definition \ref{Kontrakt2d}) and where 
  \begin{displaymath}
    \beta_{\rm{c}}: \mathcal{P}_{\rm{c}} \times \mathcal{P}_{\rm{c}} \longrightarrow
    \mathcal{P}_{m, \mathcal{W}_{O_F}(R)} 
  \end{displaymath}
  is a polarization which is anti-linear for the $O_K$-action $\iota_{\rm{c}}$. 
\end{definition}

Let $r$ be a special local CM-type with reflex field $E$. We regard an
algebra $R \in \Nilp_{O_{\breve{E}}}$ as an $O_{\breve{F}}$-algebra via $\breve\varphi_0$.
We now define the contracting  functor for polarized CM-pairs, \index{contracting  functor for polarized CM-pairs}\index[NO]{CCC@$\mathfrak{C}^{\rm pol}_{r, R}$}
\begin{equation}\label{P'Pol12e} 
  \mathfrak{C}^{\rm pol}_{r, R}: \mathfrak{d}\mathfrak{P}^{\rm pol}_{r,R}
  \longrightarrow \mathfrak{d}{\mathfrak R}^{\rm pol}_{R}.
\end{equation}
Let $(\mathcal{P}, \iota, \beta) \in \mathfrak{d}\mathfrak{P}^{\rm pol}_{r,R}$.
We apply the contracting functor $\mathfrak{C}_{r, R}$  to
$(\mathcal{P}, \iota)$ and obtain
$(\mathcal{P}_{\rm{c}}, \iota_{\rm{c}}) \in \mathfrak{d}{\mathfrak R}_{R}$, cf. Definition \ref{def:contr}. By Theorem
\ref{P'Pol2t}, the polarization
$\beta: \mathcal{P} \times \mathcal{P} \longrightarrow \mathcal{P}_{m,R}$ induces
an alternating bilinear form
\begin{equation}\label{P'Pol13e}
  \tilde{\beta}_{\rm{c}}: \mathcal{P}_{\rm{c}} \times \mathcal{P}_{\rm{c}} \longrightarrow
  \mathcal{P}_{m, \mathcal{W}_{O_F}(R)}(\pi^{ef}/p^f). 
\end{equation}
If we combine this with the chosen isomorphism (\ref{P'Pol7e1}), we obtain a
polarization of the $\mathcal{W}_{O_F}(R)$-display $\mathcal{P}_{\rm{c}}$,
\begin{displaymath}
  \beta_{\rm{c}}: \mathcal{P}_{\rm{c}} \times \mathcal{P}_{\rm{c}} \longrightarrow
  \mathcal{P}_{m, \mathcal{W}_{O_F}(R)}. 
\end{displaymath}
Then $(\mathcal{P}_{\rm{c}}, \iota_{\rm{c}}, \beta_{\rm{c}})$ is  defined to be the image of $(\mathcal{P}, \iota, \beta)$ by the functor
$\mathfrak{C}^{\rm pol}_{r, R}$. 

\begin{theorem}\label{P'Pol2c}
  Let $R \in \Nilp_{O_{\breve{E}}}$ be such that the ideal of nilpotent elements is nilpotent. The contracting functor $\mathfrak{C}^{\rm pol}_{r,R}$ induces an
  equivalence of categories
  \index[NO]{RCD@$\mathfrak{d}\mathfrak{R}^{\rm nilp, pol}_{R}$}
  \index[NO]{PCN@$\mathfrak{d}\mathfrak{P}^{\rm ss, pol}_{r,R}$}
  \begin{displaymath}
    \mathfrak{C}^{\rm pol}_{r,R}: \mathfrak{d}\mathfrak{P}^{\rm ss, pol}_{r,R}
    \longrightarrow \mathfrak{d}{\mathfrak R}^{\rm nilp, pol}_{R}.
  \end{displaymath}
  Let $(\mathcal{P}_{\rm{c}}, \iota_{\rm{c}}, \beta_{\rm{c}})$ the image of
  $(\mathcal{P}, \iota, \beta)$ under the functor $\mathfrak{C}^{\rm pol}_{r, R}$. Then
  \begin{displaymath}
\height_{O_F} \beta_{\rm{c}} = \frac{1}{f} \height \beta ,
    \end{displaymath}
cf. Definition \ref{def:pol}. 

{\rm Here the index ''$\text{\rm ss}$'' indicates the full subcategory of
  supersingular displays and the index ''$\text{\rm nilp}$'' the full
  subcategory of nilpotent displays.}
  \end{theorem}
\begin{proof}
  We use the notation of Proposition \ref{P'Pol1p}. We have a
  commutative diagram
  \begin{displaymath}
  \xymatrix{
    P \ar[rrd]_{\alpha} \ar[rr] & \overset{\tilde{\alpha}}{-} &
    {\Hom_{O_F \otimes_{\mathbb{Z}_p} W(R)}(P, O_F \otimes_{\mathbb{Z}_p} W(R))}
    \ar[d]^{\mathbf{t}_{\ast}}\\
 & & {\Hom_{W(R)}(P, W(R))}
  }
  \end{displaymath}
  Here the map $\tilde{\alpha}$ is induced by $\tilde{\beta}$ and the map
  $\alpha$ is induced by $\beta$. The vertical map is defined by
  $\mathbf{t}_{\ast}(\ell) = \mathbf{t} \circ \ell$ and is an isomorphism. The
  map $\alpha$ induces the isogeny $\mathcal{P} \longrightarrow \mathcal{P}^{\vee}$
  associated to $\beta$ and the map $\tilde{\alpha}$ induces the isogeny 
  $\mathcal{P}' \longrightarrow (\mathcal{P}')^{\Delta}$. Therefore these isogenies
  have the same height. If we apply the Ahsendorf functor to the last isogeny
  we obtain the map
  $\mathcal{P}_{\rm{c}} \longrightarrow \mathcal{P}_{\rm{c}}^{\vee}(\pi^{ef}/p^f)$ which is
  associated to $\tilde{\beta}_{\rm{c}}$. By Proposition \ref{Adorfslopes1p}
  we obtain
  \begin{displaymath}
    \height_{O_F} \beta_{\rm{c}}  = \height_{O_F} \tilde{\beta}_{\rm{c}} =
    \frac{1}{f} \height \beta. 
    \end{displaymath}
  \end{proof}
\begin{remark}\label{P'Pol2r}
  Let us explain how the bijection between bilinear forms of Theorem
  \ref{P'Pol2t} simplifies when $R = k$ is a perfect field in $\Nilp_{O_{E'}}$. We take $\CP_1=\CP_2$. 
  
  We consider the Dieudonn\'e module $(P, F, V)$ of $\mathcal{P}$. We consider $\beta: P \times P \longrightarrow W(k)$ as
  a bilinear form of Dieudonn\'e modules. Here we mean by $W(k)$ the
  Dieudonn\'e module $(W(k), F, V)$, cf. (\ref{bilin2e}). We define
  \begin{equation}
\tilde{\beta}: P \times P \longrightarrow O_F \otimes_{\mathbb{Z}_p} W(k)  
    \end{equation}
  as in Proposition \ref{P'Pol1p}. We know that $\tilde{\beta}$ induces a
  bilinear form of displays
  $\mathcal{P} \times \mathcal{P} \longrightarrow \mathcal{L}_k$. In terms of
  Dieudonn\'e modules, this means that the following equation holds,
  \begin{equation}\label{P'Pol23e}
\tilde{\beta}(V' x_1, V' x_2) = V_{\mathcal{L}} \tilde{\beta}(x_1, x_2). 
    \end{equation}
In terms of the decomposition (\ref{OWdecomp2e}), the
operator $V'$ is given by (\ref{KottwitzC39e}). 

By (\ref{LTD19e}), the Ahsendorf functor applied to $\mathcal{L}_k$ gives the
$\mathcal{W}_{O_F}(k)$-Dieudonn\'e module
\begin{equation}\label{P'Pol9e}
(O_F \otimes_{O_{F^t}, \tilde{\psi}_0} W(k), \frac{\pi^{ef}}{p^f} F^f,
    \frac{p^f}{\pi^{ef-1}} F^{-f}). 
  \end{equation}
  The bilinear form $\tilde{\beta}$ gives by restriction  to $P_{\rm{c}}=P\otimes_{O_{F^t}, \tilde{\psi}_0} W(k)\subset P$ the 
$O_F \otimes_{O_{F^t}, \tilde{\psi}_0} W(k)$-bilinear form
\begin{equation}\label{P'Pol17e} 
  \tilde{\beta}_{\rm{c}}: P_{{\rm{c}}} \times P_{{\rm{c}}} \longrightarrow
  O_F \otimes_{O_{F^t}, \tilde{\psi}_0} W(k). 
  \end{equation}
  Because this is obtained by applying the Ahsendorf functor to (\ref{P'Pol23e}), 
$\tilde{\beta}_{\rm{c}}$ is a bilinear form of $\mathcal{W}_{O_F}(k)$-Dieudonn\'e
modules if we equip  the right hand side with the
$\mathcal{W}_{O_F}(k)$-Dieudonn\'e module structure (\ref{P'Pol9e}).
 Therefore we obtain  
\begin{equation}\label{P'Pol10e}
  \tilde{\beta}_{\rm{c}}(V_{{\rm{c}}} x_{1}, V_{{\rm{c}}} x_{2}) = 
  \frac{p^f}{\pi^{ef-1}} ~^{F^{-f}}\tilde{\beta}_{\rm{c}}(x_{1}, x_{2}), \quad
  x_1, x_2 \in P_{{\rm{c}}}. 
  \end{equation}
In the case where $K/F$ is ramified, we have $P_{{\rm{c}}} = P_{ \psi_0}$, and 
\begin{displaymath}
V_{{\rm{c}}} = \Pi^{-ef+1}V^f: P_{ \psi_0} \longrightarrow P_{ \psi_0} , 
  \end{displaymath}
cf.  (\ref{CRF1e}). Note that  (\ref{P'Pol10e}) can be checked easily from these expressions.

In the case where $K/F$ is unramified, we have  $P_{{\rm{c}}} = P_{\psi_0} \oplus P_{\bar{\psi}_0}$, and 
 $V_{{\rm{c}}}$ is the endomorphism of
  $P_{{\rm{c}}} = P_{ \psi_0} \oplus P_{ \bar{\psi}_0}$ given by the matrix
  \begin{equation*}
    \left(
    \begin{array}{cc}
      0 & \pi^{-g_{\psi_0}} V^f\\
      \pi^{-g_{\bar{\psi}_0}} V^f & 0
    \end{array}
    \right),  
  \end{equation*}
  cf. \eqref{P'Pol14e}. 
  Before (\ref{P'Pol4e}) we already remarked that
  $\tilde{\beta}(P_{\psi_0}, P_{\psi_0}) = 0  = \tilde{\beta}(P_{\bar{\psi}_0}, P_{ \bar{\psi}_0})$.   
  Again (\ref{P'Pol10e}) can be checked directly on these  descriptions of $V_{{\rm{c}}}$. 

  Now we assume moreover that $k$ is a 
  $O_{\breve{E}}$-algebra. We have the map (\ref{P'Pol8e}),
  \begin{displaymath}
O_{\breve{F}} \longrightarrow O_F \otimes_{O_{F^t}, \tilde{\psi}_0} W(k) = W_{O_F}(k).  
  \end{displaymath}
  We consider the image $\eta_{0,k} \in O_F \otimes_{O_{F^t}, \psi_0} W(k)$ of
  $\eta_0$. We set
  \begin{equation}\label{P'Pol15e}
\beta_{\rm{c}} = \eta_{0,k}^f \tilde{\beta}_{\rm{c}}: P_{{\rm{c}}} \times P_{{\rm{c}}} \longrightarrow
  O_F \otimes_{O_{F^t}, \tilde{\psi}_0} W(k).
  \end{equation}
  Then we find
  \begin{displaymath}
    \begin{aligned}
    \beta_{\rm{c}}(V_{{\rm{c}}} x_{1}, V_{{\rm{c}}} x_{2}) &=
    \eta_{0,k}^f \tilde{\beta}_{\rm{c}}(V_{{\rm{c}}} x_{1}, V_{{\rm{c}}} x_{2}) =  \eta_{0,k}^f
    \frac{p^f}{\pi^{ef-1}} ~^{F^{-f}}\tilde{\beta}_{\rm{c}}(x_{1}, x_{2})\\
    &= \eta_{0,k}^f \frac{p^f}{\pi^{ef-1}}
    ~^{F^{-f}}(\eta_{0,k}^{-f}\beta_{\rm{c}}(x_{1}, x_{2})) 
    = \pi ~^{F^{-f}}\beta_{\rm{c}}(x_{1}, x_{2}), 
      \end{aligned}
  \end{displaymath}
  since $\eta_{0,k}^f ~^{F^{-f}}\!(\eta_{0,k}^{-f}) = \pi^{ef}/p^f$. Indeed, the
  left hand side of the last identity is the image of
  $(\eta_0 \tau^{-1}(\eta_0^{-1}))^f = \pi^{ef}/p^f$.

  This shows that $\beta_{\rm{c}}$ is a bilinear form of
  $\mathcal{W}_{O_F}(k)$-Dieudonn\'e modules, if we consider on
  $O_F \otimes_{O_{F^t}, \tilde{\psi}_0} W(k)$ the $\mathcal{W}_{O_F}(k)$-Dieudonn\'e
  module structure which corresponds to $\mathcal{P}_{m, \mathcal{W}_{O_F}(R)}$, namely 
  \begin{displaymath}
(O_F \otimes_{O_{F^t}, \tilde{\psi}_0} W(k), F^f, \pi F^{-f}). 
  \end{displaymath}
\end{remark}
 \begin{remark}
 Let us discuss the height identity in Theorem \ref{P'Pol2c} in a more direct way. We may
  assume that $R$ is a perfect field. We may write the equation in the form
  \begin{equation}
\ord_p \textstyle{\det_{W(k)}} \beta = f \ord_{\pi} \det_{W_{O_F}(k)} \beta_{\rm{c}}.
  \end{equation}
  On the left hand side the determinant is taken with respect to an arbitrary
  basis of the $W(k)$-module $P$. After we take $\ord_p$, the result is independent
  of the choice of the basis. 
    The right hand side of this equation does not change if we replace $\beta_{\rm{c}}$
    by the form $\tilde{\beta}_{\rm{c}}$ of (\ref{P'Pol17e}). 
    We begin with the ramified case. The decomposition $P = \oplus P_{\psi}$
    is orthogonal with respect to $\beta$. Let $\beta_{\psi}$ be the restriction
    to $P_{\psi}$. Let $\psi$ be banal. The map
    $\Pi^{-e}V: P_{\psi \sigma} \longrightarrow P_{\psi}$ is a $F^{-1}$-linear
    isomorphism. From the equation 
    \begin{displaymath}
      \beta_{\psi}(\Pi^{-e}V x, \Pi^{-e}V y) = \beta_{\psi}(V(\pi^{-e}x), Vy) =
      ~^{F^{-1}}\beta_{\psi \sigma}(\pi^{-e} px, y) 
      \end{displaymath}
    we conclude that
    $\ord_p \det_{W(k)} \beta_{\psi} = \ord_p \det_{W(k)} \beta_{\psi \sigma}$.
    Therefore this value is independent of $\psi$. In particular we obtain 
    \begin{displaymath}
      \ord_p \textstyle{\det_{W(k)}} \beta = f \ord_p \det_{W(k)} \beta_{\psi_0} =
      f \ord_{\pi} \det_{W_{O_F}(k)} \tilde{\beta}_{{\rm{c}}}.  
      \end{displaymath}
    The last equation follows because
    \begin{displaymath}
      \tilde{\beta}_{{\rm{c}}}: P_{\psi_0} \times P_{\psi_0} \longrightarrow
      W_{O_F}(k) = O_F \otimes_{O_{F^t}, \tilde{\psi_0}} W(k)
      \end{displaymath}
    is obtained from $\beta_{\psi_0}$ by the equation
    \begin{displaymath}
      \Trace_{W_{O_F}(k)/W(k)} (\vartheta^{-1} a \tilde{\beta}_{{\rm{c}}}(x,y)) =
      \beta_{\psi_0}(ax, y), \quad x, y \in P_{\psi_0}, \;
      a \in O_F \otimes_{O_{F^t}, \tilde{\psi_0}} W(k), 
    \end{displaymath}
    and since the pairing
    \begin{displaymath}
      \Trace_{W_{O_F}(k)/W(k)}(\vartheta^{-1} a_1 a_2):
      (O_F \otimes_{O_{F^t}, \tilde{\psi_0}} W(k)) \times
      (O_F \otimes_{O_{F^t}, \tilde{\psi_0}} W(k)) \longrightarrow W(k) 
    \end{displaymath}
    is perfect.

    In the unramified case we write $\psi \sigma^f = \bar{\psi}$. The
    modules $P_{\psi_1}$ and $P_{\psi_2}$ are orthogonal for
    $\psi_1 \neq \bar{\psi_2}$. We denote by $\beta_{\psi}$ the restriction
    of $\beta$ to $P_{\psi} \times P_{\bar{\psi}}$.
    We define $\ord_p \det \beta_{\psi}$ by taking an arbitrary basis of
    $P_{\psi}$ and an arbitrary basis of $P_{\bar{\psi}}$. Assume that $\psi$
    is banal; then $\bar{\psi}$ is also banal. We obtain two $F^{-1}$-linear
    isomorphisms
    \begin{displaymath}
      \pi^{-a_{\psi}} V: P_{\psi \sigma} \longrightarrow P_{\psi}, \quad
      \pi^{-a_{\bar{\psi}}} V: P_{\bar{\psi} \sigma} \longrightarrow  P_{\bar{\psi}}.  
    \end{displaymath}
    We have
    \begin{displaymath}
      \beta_{\psi} (\pi^{-a_{\psi}} V x, \pi^{-a_{\bar{\psi}}} V y) =
      \beta_{\psi} (\pi^{-e} V x, V y) = ~^{F^{-1}} \beta_{\psi \sigma}(p\pi^{-e}x,y). 
    \end{displaymath}
    We conclude that
    $\ord_{p} \det_{W} \beta_{\psi} = \ord_{p} \det_{W} \beta_{\psi \sigma}$.
    Because $\beta$ is alternating, $\beta_{\psi_0}$ and $\beta_{\bar{\psi}_0}$
    have the same order of determinant. We conclude that
    $h := \ord_{p} \det_{W} \beta_{\psi}$ is independent of $\psi \in \Psi$.
    We find $\ord_{p} \det_{W} \beta = 2fh$.
    The form $\tilde{\beta}_{\rm{c}}$ is obtained from the restriction
    $\beta_{\psi_0}$ by the equation
    \begin{displaymath}
      \Trace_{O_F \otimes_{O_{F^t}, \tilde{\psi_0}} W(k)/W(k)}
      (a\vartheta^{-1}\tilde{\beta}_{\rm{c}}(x,y)) =
      \beta_{\psi_0}(ax,y), \quad x \in P_{\psi_0}, \; y \in P_{\tilde{\psi}_0},\; a\in O_F. 
    \end{displaymath}
    Therefore we obtain 
    \begin{displaymath}
      2h = 2 \ord_p \textstyle{\det_W} \beta_{\psi_0} =
      \ord_{\pi} \det_{O_F \otimes_{O_{F^t}, \tilde{\psi_0}} W(k)} \tilde{\beta}_{\rm{c}}. 
      \end{displaymath}
\end{remark}

\subsection{The contracting functor in the case of a banal CM-type}
In the banal case we will associate to an object of the category
$\mathfrak{d}\mathfrak{P}'_{r,R}$ (cf. Definition \ref{Kontrakt1d}) an \'etale
sheaf on $\Spec R$. The construction does not use  the Ahsendorf functor $\mathfrak{A}_{O_F/\BZ_p}$, which is not useful here. 

\begin{definition}\label{Cbanal1d} 
Let $R$ be a ring. 
An \emph{\'etale Frobenius module} \index{\'etale Frobenius module}is a pair $(M, \Theta)$,
where $M$ is a finitely generated $W(R)$-module which is locally
on $\Spec R$  free and where $\Theta: M \longrightarrow M$ is a Frobenius linear isomorphism, i.e.,
$\Theta:\sigma^*(M)\to M$ is an isomorphism. 
\end{definition}
The following proposition is a variant of a result of Drinfeld, comp. \cite[Prop. 2.1]{DrinF}.  It can also be proved using the theory of displays. When $R$ is an algebraically closed field, the proposition
is a theorem of Dieudonn\'e.

\begin{proposition}\label{Cbanal1p} 
Let $R$ be a ring such that $p$ is nilpotent in $R$. There is a functor $\mathfrak{A}$  
from the category of \'etale Frobenius modules over $R$ to the category 
of locally constant $p$-adic \'etale sheaves which are finitely generated and
flat over $\mathbb{Z}_p$.   The functor $\mathfrak{A}$ is an equivalence of categories which commutes with arbitrary
base change. It is compatible with the tensor product of \'etale Frobenius modules, resp., of $p$-adic \'etale sheaves. 
\end{proposition}
\begin{proof}
  We give a sketch of the proof which shows how this equivalence is
  constructed. By \cite[Lem. 42]{Zi}  it follows that the category of
  \'etale Frobenius modules over $R$ and $R/pR$ are equivalent. Indeed,
  a Frobenius module lifts locally and by loc.~cit. two liftings are
  canonically isomorphic. Therefore we may assume  
  $pR = 0$. Let $(M, \Theta)$ be an \'etale Frobenius module. We set
$M_n = W_n(R) \otimes_{W(R)} M$.  Because $pR = 0$, the Frobenius $F$ on $W(R)$ induces a
Frobenius $F: W_n(R) \longrightarrow W_n(R)$. By base change
we obtain a $F$-linear map $\Theta_n: M_n \longrightarrow M_n$. We define a
functor $\mathfrak{A}_n$ on the category of $R$-algebras,
\begin{displaymath}
  \mathfrak{A}_n(S) = \{ x \in W_n(S) \otimes_{W_n(R)} M_n \; | \; \Theta_{n,S} (x) = x \} .
\end{displaymath}
One can show that $\mathfrak{A}_n$ is representable by a finite \'etale scheme over
$\Spec R$. Clearly $\mathbb{Z}/p^n \mathbb{Z} = W_n(\mathbb{F}_p)$ acts on
$\mathfrak{A}_n$. We define the associated $p$-adic sheaf 
\begin{displaymath}
  \mathfrak{A}_{(M, \Theta)} =
  \underset{\leftarrow}{\lim}{}_{{}_n} \; \mathfrak{A}_n.
\end{displaymath}
Let $W$ be the \'etale sheaf of Witt vectors. We have  $W \otimes_{\mathbb{Z}_p} \mathfrak{A}_{(M, \Theta)} = M$ in the sense of
\'etale sheaves, where the action
of $\Theta$ corresponds on the left hand side to the action of
$F \otimes \id$.

Finally, we show the compatibility
with tensor products. If $(M', \Theta')$ is a second \'etale Frobenius
module, we set $(N,\Xi) = (M \otimes_{W(R)} M', \Theta \otimes \Theta')$.  
We obtain a  natural homomorphism
\begin{equation}
\mathfrak{A}_{(M, \Theta)} \otimes_{\mathbb{Z}_p} \mathfrak{A}_{(M', \Theta')} \longrightarrow \mathfrak{A}_{(N, \Xi)} . 
\end{equation}
 To prove that this is an isomorphism, we may reduce by base change to the case
where $R$ is an algebraically closed field. Then the assertion is clear by
the theorem of Dieudonn\'e. 
\end{proof}
\begin{definition}\label{Cfunctor1d}
  Let $ r$ be  banal. Let $R\in \Nilp_{O_{E'}}$. Let ${\rm Et}(O_K)_R$
  \index[NO]{EAG@${\rm Et}(O_K)_R$} be 
 the category of locally
constant $p$-adic \'etale sheaves $G$ over $\Spec R$ which are $\mathbb{Z}_p$-flat
with $\rank_{\mathbb{Z}_p} G = 4d$ and which are equipped with an action
\begin{displaymath}
\iota: O_K \longrightarrow \End_{\mathbb{Z}_p} G. 
\end{displaymath}
The contracting functor  is the functor
$$\mathfrak{C}_{r, R}\colon\mathfrak{d}\mathfrak{P}_{r,R}\to {\rm Et}(O_K)_R, $$  
which is the composite of the pre-contracting functor $\mathfrak{C}'_{r, R}$ of Theorem \ref{KottwitzC7p} and the functor $\mathfrak{A}$
\index[NO]{ACC@$\mathfrak{A}$} of
  Proposition \ref{Cbanal1p}, applied to the \'etale Frobenius module $(P', \dot{F}')$. 
The functor commutes with arbitrary base change $R \longrightarrow R'$.  
\end{definition}

\begin{theorem}\label{Cbanal2p}
  Let $ r$ be  banal.
  Let $R \in \Nilp_{O_{E'}}$ be such that the ideal of nilpotent elements of $R$
  is nilpotent. Then the contracting functor\index{contracting functor} is an equivalence of categories,
  \begin{displaymath}
  \mathfrak{C}_{r, R}\colon \mathfrak{d}\mathfrak{P}_{r,R}\to {\rm Et}(O_K)_R .
    \end{displaymath}
  \end{theorem}
\begin{proof}
  Since the objects in $\mathfrak{d}\mathfrak{P}'_{r,R}$ are \'etale, this is
  simply a combination of Proposition \ref{KottwitzC7p} and Proposition
  \ref{Cbanal1p}
\end{proof}

\begin{remark}\label{Cbanal1r}
  In the banal case there is a functor
  \begin{equation}\label{BT-Dsp1e} 
  \mathfrak{P}_{r,R} \longrightarrow \mathfrak{d}\mathfrak{P}_{r,R}
    \end{equation}
  from $p$-divisible groups to displays which is an equivalence of categories.
  Indeed, in the case when $K/F$ is a field extension, the displays of objects in
  $\mathfrak{d}\mathfrak{P}_{r,R}$ are by Corollary \ref{C'onslopes1c}
  isoclinic of constant slope $1/2$ and therefore nilpotent. Therefore
  they are displays of formal $p$-divisble groups, cf. Theorem \ref{maindisp}.
  In the split case $O_K = O_F \times O_F$ we have a corresponding
  decomposition of a display $\mathcal{P} \in \mathfrak{d}\mathfrak{P}_{r,R}$:
  $\mathcal{P} = \mathcal{P}_{1} \oplus \mathcal{P}_{2}$. In the case
  where $\mathcal{P}$ is nilpotent we can argue as before. If not, one of the
  summands is \'etale and the other is isoclinic of slope $1$, cf.
  Corollary \ref{C'onslopes1c}. But we have an equivalence between \'etale
  $p$-divisible groups over $R$ and \'etale displays over $R$, which is
  easily defined by the $\mathfrak{A}$-functor. Therefore we conclude also
  in this case that the equivalence \eqref{BT-Dsp1e} exists.

  For more information about the resulting functor
  $\mathfrak{P}_{r,R} \longrightarrow {\rm Et}(O_K)_R$ cf. section
  \ref{ss:deeper}.  
  \end{remark} 

We now add polarizations to the picture.

\begin{lemma}\label{Cbanal1l}  
  Let $r$ be banal.  Let $R \in \Nilp_{O_{E'}}$.  Let $(\mathcal{P}_1, \iota_1)$ and
$(\mathcal{P}_2, \iota_2)$ be in $\mathfrak{d}\mathfrak{P}_{r, R}$.  
Let
  \begin{displaymath}
\bar{\beta}: P_1/I(R) P_1 \times P_2/I(R) P_2 \longrightarrow R 
  \end{displaymath}
  be an $R$-bilinear form such that
  \begin{displaymath}
    \bar{\beta}(\iota_1(a)x_1, x_2) = \bar{\beta}(x_1, \iota_2(\bar{a}) x_2),
    \quad \; a \in O_K.
    \end{displaymath}
  Then the restriction of $\bar{\beta}$ to $Q_1/I(R)P_1 \times Q_2/I(R)P_2$
  is zero. 
\end{lemma}
\begin{proof}
  We first consider  the case where $K/F$ is ramified. 
  We consider $P_{i,\psi}/I(R)P_{i,\psi}$ as an $O_K \otimes_{O_{F^t}, \psi} R$-module
  for $i =1,2$. Because of the isomorphism \eqref{KottwitzC6e}, it suffices to show that 
  \begin{displaymath}
    \bar{\beta}(\mathbf{E}_{A_{\psi}} x_1, \mathbf{E}_{A_{\psi}} x_2)
    = 0, \quad \; x_1 \in (P_1/I(R)P_1)_{\psi}, \;   x_2 \in
    (P_2/I(R)P_2)_{\psi}. 
    \end{displaymath}
  We consider
  $\mathbf{E}_{A_{\psi}}(\Pi \otimes 1) \in O_K \otimes_{O_{F^t}, \psi} O_{E'}$.
  The image of this element by the conjugation of $K/F$ is
  $(-1)^e \mathbf{E}_{B_{\psi}}(\Pi \otimes 1)$, cf. the proof of Lemma
  \ref{P'Pol1l}. Therefore we find
  \begin{displaymath}
    \bar{\beta}(\mathbf{E}_{A_{\psi}} x_1, \mathbf{E}_{A_{\psi}} x_2) =
    (-1)^e \bar{\beta}(x_1, \mathbf{E}_{B_{\psi}} \mathbf{E}_{A_{\psi}} x_2) =
    (-1)^e \bar{\beta}(x_1, \mathbf{E}_{\psi}  x_2) = 0. 
        \end{displaymath}
  Now we assume that $K/F$ is unramified. Then the condition on $\bar{\beta}$
  implies that $P_{1,\psi_1}/I(R) P_{1,\psi_1}$ and $P_{2,\psi_2}/I(R) P_{2,\psi_2}$ are
  orthogonal with respect to $\bar{\beta}$ if $\psi_1 \neq \bar{\psi}_2$.
Again by the isomorphism \eqref{KottwitzC6e}, it suffices to show that  
  \begin{displaymath}
    \bar{\beta}(\mathbf{E}_{A_{\psi}} x_1, \mathbf{E}_{A_{\bar{\psi}}} x_2)
    = 0, \quad \; x_1 \in (P_1/I(R)P_1)_{\psi}, \; 
    x_2 \in (P_2/I(R)P_2)_{\bar{\psi}}. 
    \end{displaymath}
  In this case the conjugation of $K/F$ maps $\mathbf{E}_{A_{\psi}}$ to
  $\mathbf{E}_{B_{\bar{\psi}}}$. Therefore the last equation follows from
  \begin{displaymath}
    \mathbf{E}_{B_{\bar{\psi}}}(\pi \otimes 1)
    \mathbf{E}_{A_{\bar{\psi}}}(\pi \otimes 1) =
    \mathbf{E}_{\bar{\psi}}(\pi \otimes 1) = 0.
  \end{displaymath}
  Exactly the same argument applies to the split case. 
\end{proof}
\begin{lemma}\label{Cbanal2l}
  In the situation of the last lemma, assume that $R$ is
  a reduced ring. Let $\beta\colon P_1 \times P_2 \longrightarrow W(R)$ be a
  $W(R)$-bilinear form such that $\beta$ is anti-linear for the $O_K$-actions $\iota_1$, resp. $\iota_2$, and such that 
  \begin{displaymath}
\beta(F_1x_1, F_2x_2) = p~^{F}\beta(x_1, x_2), \quad x_1 \in P_1, 
      x_2 \in P_2.
  \end{displaymath}
  Then $\beta$ induces a bilinear form of displays
  \begin{displaymath}
\beta: \mathcal{P}_1 \times \mathcal{P}_2 \longrightarrow \mathcal{P}_{m,R}. 
    \end{displaymath}
\end{lemma} 
\begin{proof}
  We must verify that $\beta(Q_1, Q_2) \subset I(R)$ and that
  \begin{displaymath}
    \beta(\dot{F}_1 y_1, \dot{F}_2 y_2) = ~^{\dot{F}}\beta(y_1,y_2), \quad
     \; y_1 \in Q_1, \; y_2 \in Q_2.     
    \end{displaymath}
  The inclusion is a consequence of Lemma \ref{Cbanal1l}. To verify the last 
  equation, we may multiply it by $p^2$ because $p$ is not a zero divisor
  in $W(R)$. But then it follows from the assumptions on $\beta$. 
  \end{proof}

\begin{definition}
  Let $\rho \in O_F \otimes_{\mathbb{Z}_p} W(R)$ be a unit. We define
  $O_F(\rho)$\index[NO]{OAA@$O_F(\rho)$} as the $p$-adic \'etale sheaf associated by Proposition
  \ref{Cbanal1p} to the \'etale Frobenius module
  $(O_F \otimes_{\mathbb{Z}_p} W(R), \Theta_{\rho})$, where
\begin{equation}\label{Cbanal6e}
  \Theta_{\rho} (a \otimes \xi) = \rho\cdot (a \otimes ~^{F}\xi),
  \quad a \in O_F,\; \xi \in W(R).  
\end{equation}
\end{definition}
When $\rho = 1$ we obtain the constant $p$-adic \'etale sheaf
$O_F = O_F(1)$. 

Let $\rho=\pi^e/p$.
Let $R\in \Nilp_{O_{E'}}$ and let
$(\mathcal{P}_i, \iota_i) \in \mathfrak{P}_{r, R}$ for $i = 1,2$. We will
associate to a bilinear form of displays
\begin{equation}\label{Cbanal3e}
  \beta: \mathcal{P}_1 \times \mathcal{P}_2 \longrightarrow \mathcal{P}_{m, R} 
    \end{equation}
which is anti-linear for the $O_K$-actions $\iota_1$, resp. $\iota_2$, 
a  bilinear form of $p$-adic \'etale sheaves which is anti-linear for the $O_K$-actions on ${C}_{\mathcal{P}_1}={\mathfrak C}_{r, R}(\CP_1)$, resp. ${C}_{\mathcal{P}_2}={\mathfrak C}_{r, R}(\CP_2)$,
\begin{equation}\label{Cbanal4e}
\phi: {C}_{\mathcal{P}_1} \times {C}_{\mathcal{P}_2} \longrightarrow O_F(\rho). 
\end{equation}
For the construction we may assume that $R$ is a $\kappa_{E'}$-algebra because
\'etale sheaves are insensitive to nilpotent elements. 

Let first $K/F$ be ramified. Then we find for $x_1 \in P_1$
and $x_2 \in P_2$ that
\begin{equation}\label{Cbanal7e}
  \begin{aligned}
    \beta(\dot{F}'_1 x_1, \dot{F}'_2 x_2) =
    \beta(\dot{F}_1 \Pi^{e}x_1, \dot{F}_2 \Pi^{e}x_2) &
    = ~^{\dot{F}}\beta(\Pi^e x_1, \Pi^e x_2)\\ 
    & = ~^{\dot{F}}\beta(\pi^e x_1, x_2) =
    ~^{F}\beta(\frac{\pi^e}{p} x_1, x_2).
  \end{aligned}
\end{equation}
We used the equation
\begin{equation}\label{Cbanal12e}
\dot{F}'_1 x_1 = \dot{F}_1 \Pi^{e}x_1  ,
  \end{equation}
which follows from (\ref{KottwitzC33e}) and (\ref{KottwitzC17e}).
Recall the function
$\mathbf{t}(a) = \Trace_{O_F/\mathbb{Z}_p} \vartheta^{-1}a$, for $a \in O_{F}$, 
where $\vartheta \in O_F$ is the different of $F/\mathbb{Q}_p$, cf. p. \pageref{page:theta}.
We define
$\tilde{\beta}: P_1 \times P_2 \longrightarrow O_F \otimes_{\mathbb{Z}_p} W(R)$
by the equation
\begin{equation}\label{Cbanal-t}
  \mathbf{t}(\xi \tilde{\beta}(x_1, x_2)) = \beta(\xi x_1, x_2), \quad
  \xi \in O_F \otimes_{\mathbb{Z}_p} W(R). 
  \end{equation}
Then $\tilde{\beta}$ is a bilinear form of $O_F \otimes_{\mathbb{Z}_p} W(R)$-modules.
We conclude from (\ref{Cbanal7e}) that

\begin{equation}\label{Cbanal8e}
  \tilde{\beta}(\dot{F}'x_1, \dot{F}'x_2) =
  \frac{\pi^e}{p} \cdot ~^{F}\tilde{\beta}( x_1, x_2).
\end{equation}
Hence $\tilde \beta$  is a bilinear form of Frobenius modules. Since the functor $\mathfrak{A}$ of
Proposition \ref{Cbanal1p} commutes with tensor products, it induces a bilinear
form (\ref{Cbanal4e}). 

Now we consider the case where $K/F$ is an unramified field extension.
For each $O_{E'}$-algebra $R$ we have the decomposition

\begin{equation}\label{Cbanal9e}
  O_K \otimes_{\mathbb{Z}_p} W(R) = \prod_{\psi \in \Psi}
  O_K \otimes_{O_{K^t}, \tilde{\psi}} W(R),
  \end{equation}
which is induced by (\ref{OWdecomp1e}). The conjugation of $K/F$ acts on  
$O_K \otimes_{\mathbb{Z}_p} W(R)$ via the first factor. We denote this by
$\eta \mapsto \bar{\eta}$. We denote
by $~^{F}\eta$ the action of the Frobenius via the second factor.
On the right hand side of (\ref{Cbanal9e}) these actions become 
\begin{equation}\label{Cbanal23e}
  \begin{aligned}
    O_K \otimes_{O_{K^t}, \tilde{\psi}} W(O_{E'}) & \longrightarrow 
    O_K \otimes_{O_{K^t}, \widetilde{\psi} \sigma^f} W(O_{E'})\\
    a \otimes \xi & \mapsto  \quad\bar{a} \otimes \xi, \\[2mm] 
  O_K \otimes_{O_{K^t}, \tilde{\psi}} W(R) & \longrightarrow 
    O_K \otimes_{O_{K^t}, \sigma \tilde{\psi}} W(R)\\  
a \otimes \xi & \mapsto \quad  a \otimes ~^{F}\xi. 
  \end{aligned}
\end{equation}
Here $\sigma$ denotes the Frobenius automorphism of $\Gal(K^t/\mathbb{Q}_p)$.
Looking at the right hand side of (\ref{Cbanal9e}), we define
\begin{equation}\label{pirunr}
  \pi_r := (\pi^{a_{\psi}} \otimes 1)_{\psi \in \Psi} \in
  O_K \otimes_{\mathbb{Z}_p} W(O_{E'}).
  \end{equation}
It follows that
\begin{displaymath}
\pi_r \bar{\pi}_r = \pi^e \otimes 1.
\end{displaymath}
Let $(\mathcal{P}, \iota) \in \mathfrak{P}_{r, R}$. We note that, since $R$ is a
$\kappa_{E'}$-algebra, the definition of
$(\mathcal{P}', \iota')=\mathfrak{C}'_{r, R}(\CP, \iota)$ in
\eqref{KottwitzC17e} takes the form
\begin{equation}\label{comFF'unr}
\dot{F'}(x) = \dot{F}(\pi_r x), \quad F'(x) = F(\pi_r x).
\end{equation}
Now let us start with  a bilinear form 
\begin{displaymath}
  \beta: \mathcal{P}_1 \times \mathcal{P}_2 \longrightarrow \mathcal{P}_m 
\end{displaymath}
which is anti-linear for the $O_K$-actions $\iota_1$, resp. $\iota_2$. We find 
\begin{equation}\label{Cbanal10e}
  \begin{aligned}
    \beta(\dot{F}'_1 x_1, \dot{F}'_2 x_2) =
    \beta(\dot{F}_1\pi_r x_1, \dot{F}_2\pi_r x_2)
    &  = ~^{\dot{F}}\beta(\pi_r x_1, \pi_r x_2)\\
    & = ~^{\dot{F}}\beta(\pi^e x_1, x_2) =
    ~^{F}\beta(\frac{\pi^e}{p} x_1, x_2).
  \end{aligned}
\end{equation}
As before, $\beta$ defines the $O_F \otimes_{\mathbb{Z}_p} W(R)$-bilinear form 
\begin{displaymath}
\tilde{\beta}: P_1 \times P_2 \longrightarrow  O_F \otimes_{\mathbb{Z}_p} W(R), 
  \end{displaymath}
which by (\ref{Cbanal10e}) satisfies

\begin{equation}\label{Cbanal11e}
  \tilde{\beta}(\dot{F}'_1 x_1, \dot{F}'_2 x_2) = \frac{\pi^e}{p}\cdot 
 ~^{F}\tilde{\beta}( x_1, x_2). 
\end{equation}
Applying, as before, the $\mathfrak{A}$-functor to $\CP'$, we obtain the desired bilinear form \eqref{Cbanal4e}.

Finally we consider the split case. In this case we consider in the
decomposition (\ref{zerlegt1e}) the element
\begin{equation}\label{zerlegt9e} 
  \pi_r = \pi_{r, 1} \times \pi_{r, 2} = \big((\pi^{a_{\theta_1}} \otimes 1)_{\theta \in \Theta}\big)
  \times \big((\pi^{a_{\theta_2}} \otimes 1)_{\theta \in \Theta}\big) 
\end{equation}
of (\ref{zerlegt1e}).
The conjugation acts on $O_K \otimes_{\mathbb{Z}_p} W(O_{E'})$ via the first factor.
On the right hand side of (\ref{zerlegt1e}), the conjugation just interchanges
the two factors in parentheses. This shows that 
$\pi_r \bar{\pi}_r = \pi^e \otimes 1$. Now starting with a bilinear form\footnote{One should not confuse
the notation $\mathcal{P}_1$ and $\mathcal{P}_2$ with the decomposition
(\ref{zerlegt8e}) which continues to  exist, e.g., 
$\mathcal{P}_1 = \mathcal{P}_{1, 1} \oplus \mathcal{P}_{1, 2}$.}  (\ref{Cbanal3e}), the formulas (\ref{Cbanal10e}), (\ref{Cbanal11e}) from the unramified case continue to hold, and this finishes the construction in the split case. 

\begin{proposition}\label{Cbanal3p}
  Let $ r$ be banal. Let $R \in \Nilp_{O_{E'}}$ be such that the
  ideal of nilpotent elements of $R$ is nilpotent. Let
  $(\mathcal{P}_1, \iota_1)$ and $(\mathcal{P}_2, \iota_2)$ be objects of
  $\mathfrak{d}\mathfrak{P}_{r, R}$. The construction above, which associates
  to a bilinear form of displays (\ref{Cbanal3e}) a bilinear form of
  $p$-adic \'etale sheaves (\ref{Cbanal4e}) is a bijection,
  \begin{displaymath}
    \Bil_{O_K\text{\rm-anti-linear}}(\mathcal{P}_1 \times \mathcal{P}_2,
    \mathcal{P}_{m,R}) \longrightarrow
 \Bil_{O_K\text{\rm -anti-linear}}(C_{\mathcal{P}_{1}} \times C_{\mathcal{P}_{2}},
  O_F(\rho)). 
\end{displaymath}
  \end{proposition}
\begin{proof}
  We reduce the question to the case where $R$ is reduced. Indeed, let
  $S \longrightarrow R$ be a pd-thickening in the category $\Nilp_{O_{E'}}$.
  Assume that $(\mathcal{P}_i, \iota_i) \in \mathfrak{d}\mathfrak{P}_{r, S}$ for $i=1,2$.
  It follows from Proposition \ref{GM2p}
  and Lemma \ref{Cbanal1l} that any bilinear form
  \begin{displaymath}
    \bar{\beta}: \mathcal{P}_{1,R} \times \mathcal{P}_{2,R} \longrightarrow
    \mathcal{P}_{m,R}  
  \end{displaymath}
  with the properties of (\ref{Cbanal3e}) lifts uniquely to a bilinear form
  \begin{displaymath}
    \beta : \mathcal{P}_{1} \times \mathcal{P}_{2} \longrightarrow
    \mathcal{P}_{m,S}.  
  \end{displaymath}
  Since bilinear forms of \'etale sheaves have the same property, we can assume
  that $R$ is reduced. 

  We begin with the ramified case.  For $i =1,2$,  let $(\mathcal{P}_i, \iota_i) \in \mathfrak{d}\mathfrak{P}_{r, R}$,  which correspond to $(P'_i, \dot{F}'_i)$ under the pre-contracting functor $\mathfrak{C}'_{r, R}$, cf.  Theorem 
  \ref{KottwitzC7p}, and to ${C}_{\mathcal{P}_i}$ under the contraction functor $\mathfrak{C}_{r, R}$. We start with a bilinear form of
  $p$-adic sheaves 
  \begin{displaymath}
    \tilde{\phi}: {C}_{\mathcal{P}_1} \times {C}_{\mathcal{P}_2} \longrightarrow
    O_F(\frac{\pi^e}{p}) 
  \end{displaymath}
  with the properties of (\ref{Cbanal4e}). We have to construct a bilinear
  form of displays (\ref{Cbanal3e}) which induces $\phi$. 
  By Proposition \ref{Cbanal1p}, $\phi$ comes from  a bilinear form of \'etale
  Frobenius modules
  \begin{displaymath}
\tilde{\beta}: P'_1 \times P'_2 \longrightarrow O_F \otimes_{\mathbb{Z}_p} W(R)
  \end{displaymath}
  which satisfies
  \begin{equation}\label{Cbanal15e}
    \tilde{\beta}(\dot{F}'_1 x_1, \dot{F}'_2 x_2) =
    \frac{\pi^e}{p}~^{F} \tilde{\beta}( x_1, x_2).
  \end{equation}
  After applying  $\mathbf{t}$ we obtain a bilinear form $\beta$ which satisfies
  \begin{equation}\label{Cbanal14e}
\beta(\dot{F}'_1 x_1, \dot{F}'_2 x_2) =  ~^{F}\beta(\frac{\pi^e}{p} x_1, x_2).
  \end{equation}
 By (\ref{Cbanal12e}) we may write
  \begin{displaymath}
  F_i = p\dot{F}_i = \dot{F}'_i \frac{p}{\Pi^e} ,
    \end{displaymath}
  because multiplication by $p$ is injective on $P_i$. We deduce from
  (\ref{Cbanal14e})
  \begin{displaymath}
    \beta(F_1 x_1, F_2 x_2) =
    \beta(\dot{F}'_1 \frac{p}{\Pi^e}x_1, \dot{F}'_2 \frac{p}{\Pi^e}x_2) =
     ~^{F}\beta(\frac{\pi^e}{p} \frac{p}{\Pi^e}x_1, \frac{p}{\Pi^e}x_2) =
     p~^{F}\beta(x_1, x_2).
  \end{displaymath}
 By Lemma \ref{Cbanal2l}, it follows that $\beta$ is a bilinear form of displays
  $\beta: \mathcal{P}_1 \times \mathcal{P}_2 \longrightarrow \mathcal{P}_m$. This proves the ramified case.

  Now let $K/F$ be an unramified field extension. We begin with a bilinear form of
  $p$-adic \'etale sheaves (\ref{Cbanal4e}) as before. This induces a bilinear
  form of \'etale Frobenius modules
  $\tilde{\beta}: P'_1 \times P'_2 \longrightarrow O_F \otimes_{\mathbb{Z}_p} W(R)$
  which satisfies (\ref{Cbanal15e}). Using \eqref{comFF'unr}, we  rewrite this as
  \begin{displaymath}
    \tilde{\beta}(\dot{F}_1 \pi_r x_1, \dot{F}_2 \pi_r x_2) = 
    \frac{\pi^e}{p} ~^{F}\tilde{\beta}( x_1, x_2).
  \end{displaymath}
  We multiply this equation with $p^2$ and find for the left hand side 
  \begin{displaymath}
    \tilde{\beta}(F_1 \pi_r x_1, F_2 \pi_r x_2) =
    \tilde{\beta}(~^{F}\pi_r F_1 x_1, ~^{F}\pi_r F_2 x_2) =
    \tilde{\beta}(~^{F}(\bar{\pi}_r \pi_r) F_1 x_1, F_2 x_2) =
    \pi^{e}\tilde{\beta}(F_1 x_1, F_2 x_2).
  \end{displaymath}
  If we compare this to the right hand side multiplied with $p^2$, we obtain
  \begin{displaymath}
\tilde{\beta}(F_1 x_1, F_2 x_2) = p ~^{F}\tilde{\beta}( x_1, x_2).
  \end{displaymath}
  Setting now $\beta = \mathbf{t} \circ \tilde{\beta}$,  the assumptions of
  Lemma \ref{Cbanal2l} are satisfied. Therefore $\beta$ induces a bilinear
  form of displays
  $\beta: \mathcal{P}_1 \times \mathcal{P}_2 \longrightarrow \mathcal{P}_m$. 
  
In the split case the argument is the same using the $\pi_r$ which appeared in this context. 
\end{proof}

On $\Spec \bar{\kappa}_{E}$ we can choose a trivialization of the twisted constant \'etale sheaf,  
\begin{equation}\label{KneunBa6e}
O_F(\pi^e/p) \isoarrow O_F , 
\end{equation}
as follows. Choose $\eta \in O_F \otimes_{\mathbb{Z}_p} W(\bar{\kappa}_{E})$
such that
$~^{F}\eta \eta^{-1} = \pi^e/p$ (this is equivalent to the choice of $\eta_0$
in (\ref{P'Pol8e})). Then the multiplication by $\eta$ 
\begin{displaymath}
  \eta: (O_F \otimes_{\mathbb{Z}_p} W(\bar{\kappa}_{E}), (\pi^e/p) \otimes F)
  \longrightarrow (O_F \otimes_{\mathbb{Z}_p} W(\bar{\kappa}_{E}), 1 \otimes F) 
\end{displaymath}
is an isomorphism of \'etale Frobenius modules, which induces
(\ref{KneunBa6e}) under the $\mathfrak{A}$-functor into ${\rm Et}(O_K)_{\bar{\kappa}_{E}}$, cf. Proposition \ref{Cbanal1p}. 
\begin{definition}\label{def:polet}
  Let $R \in \Nilp_{O_{\breve{E}}}$. Let ${\rm Et}(O_K)_R^{\rm pol}$
  \index[NO]{EAH@${\rm Et}(O_K)_R^{\rm pol}$} be the category of 
$p$-adic \'etale sheaves $(G, \iota)\in {\rm Et}(O_K)_R$, equipped with a $O_F$-linear alternating form 
\begin{equation}\label{KneunBa9e}
\phi:G \times G \longrightarrow O_F, 
  \end{equation}
which is anti-linear for the $O_K$-action. 
\end{definition}

Using the trivialization \eqref{KneunBa6e}, and applying   Proposition \ref{Cbanal3p}, we now obtain the contracting functor with polarizations   which is a functor from   $ \mathfrak{d}\mathfrak{P}^{\rm pol}_{r,R}$ to 
    ${\rm Et}(O_K)_R^{\rm pol}$. 
  \begin{theorem}\label{KneunBa4p}
  Let $r$ be banal. Let $R\in \Nilp{O_{\breve E}}$ be such that 
   the ideal of nilpotent elements in $R$ is nilpotent. Then the contracting functor
   $\mathfrak{C}_{r, R}^{\rm pol}$\index[NO]{CCC@$\mathfrak{C}^{\rm pol}_{r, R}$}
   is an equivalence of categories,
  \begin{displaymath}
    \mathfrak{C}_{r, R}^{\rm pol}\colon \mathfrak{d}\mathfrak{P}^{\rm pol}_{r,R}\to
    {\rm Et}(O_K)_R^{\rm pol}.
    \end{displaymath}\qed
  \end{theorem}
  \begin{remark}\label{etformsplit} In the split case, let
    $C_{\mathcal{P}}=\mathfrak{C}_{r, R}^{\rm pol}(\CP)$. Let $\CP=\CP_1\times \CP_2$
    be the decomposition induced by $O_K = O_F\times O_F$. This induces a
    decomposition $C_{\mathcal{P}} = C_{\mathcal{P},1} \times C_{\mathcal{P},2}$, 
where $C_{\mathcal{P},i}$ is the \'etale sheaf associated to the Frobenius module
$(P_i, \dot{F}\pi_{r,i})$, $i =1,2$. Here the elements $\pi_{r, i}$ are defined in \eqref{zerlegt9e}. The subsheaves $C_{\mathcal{P},i}$ of $C_\CP$ are isotropic
with respect to $\phi$ as in \eqref{KneunBa9e}, and hence $\phi$ corresponds to an $O_F$-bilinear form
\begin{equation*}
\phi: C_{\mathcal{P}, 1} \times C_{\mathcal{P}, 2} \longrightarrow O_F . 
\end{equation*}
\end{remark}

\begin{remark}\label{CPexplban}
  Let $k \in \Nilp_{O_{E'}}$ be an algebraically closed field. Let
  $(\mathcal{P}, \beta) \in \mathfrak{d}\mathfrak{P}^{\rm pol}_{r,k}$. We will
  give a description of
  $\mathfrak{C}_{r, R}^{\rm pol}(\mathcal{P}, \beta) = (C_{\mathcal{P}}, \phi)$. 
We write  $\mathcal{P} = (P ,F ,V)$ as a
Dieudonn\'e module. The image of $\CP$ under the contracting functor,
a sheaf $C_{\mathcal{P}}$ on $\Spec k$,  is simply an $O_K$-module.

Assume that $K/F$ is ramified. 
From the definition of the pre-contracting functor (cf. Theorem
\ref{KottwitzC7p}) and the $\mathfrak A$-functor we have 
\begin{displaymath}
  C_{\mathcal{P}} =
  \{x \in P \; | \; V^{-1} \Pi^{e} x = x \} .
  \end{displaymath}
To describe this further, with its bilinear form of displays, we extend the bilinear form $\beta$ to
\begin{displaymath}
\tilde{\beta}: P \times P \longrightarrow O_F \otimes W(k) , 
  \end{displaymath}
cf. (\ref{Cbanal-t}). The decomposition $P = \oplus_{\psi \in \Psi} P_{\psi}$
is orthogonal
with respect to $\tilde{\beta}$ and, by restriction, we obtain for every $\psi$
\begin{displaymath}
  \tilde{\beta}_{\psi}: P_{\psi} \times P_{\psi} \longrightarrow
  O_F \otimes_{O_{F^t}, \tilde{\psi}} W(k) \subset O_F \otimes W(k).  
\end{displaymath}
Let $x_{\psi}, x'_{\psi} \in P_{\psi}$. Since $\beta$ is a polarization, we obtain 
\begin{equation}\label{Cbanal20e}
  \tilde{\beta}_{\psi \sigma} (V^{-1} \Pi^{e} x_{\psi}, V^{-1} \Pi^{e} x'_{\psi}) = 
  \frac{\pi^e}{p}  ~^{F}\tilde{\beta}_{\psi} (x_{\psi}, x'_{\psi}).  
  \end{equation}
The action of $F$ on the right hand side is defined by (\ref{Cbanal23e}). 
Fix $\psi_a \in \Psi$.  The projection $x \mapsto x_{\psi_a}$
is an isomorphism
\begin{equation}\label{Cbanal21e}
  C_{\mathcal{P}} \cong 
  \{x_{\psi_a} \in P_{\psi_a} \; | \; V^{-f} \Pi^{ef} x_{\psi_a} = x_{\psi_a} \}.  
\end{equation}
In particular, we see that $C_{\mathcal{P}}$ is indeed a free $O_K$-module of
rank $2$.
For $x, x' \in C_{\mathcal{P}}$ we obtain from (\ref{Cbanal20e})
\begin{displaymath}
   \tilde{\beta}_{\psi \sigma} (x_{\psi \sigma}, x'_{\psi \sigma}) =
  \frac{\pi^e}{p}  ~^{F}\tilde{\beta}_{\psi} (x_{\psi}, x'_{\psi}).  
\end{displaymath}
Since $\psi_a = \psi_a \sigma^f$, we obtain
\begin{displaymath}
  \tilde{\beta}_{\psi_a}(x_{\psi_a}, x'_{\psi_a}) = \left(\frac{\pi^e}{p}\right)^{f}
  ~^{F^f} \tilde{\beta}_{\psi_a} (x_{\psi_a}, x'_{\psi_a}).  
\end{displaymath}
In the same way we may interpret the sheaf $O_F(\pi^e/p)$: the projection
\begin{displaymath}
O_F \otimes_{\mathbb{Z}_p} W(k) \longrightarrow  
  O_F \otimes_{O_{F^t}, \tilde{\psi}_a} W(k)
\end{displaymath}
defines an isomorphism 
\begin{equation}\label{Cbanal24e}
  O_F(\frac{\pi^e}{p}) \cong \{a_{\psi_a} \in O_F \otimes_{O_{F^t}, \tilde{\psi}_a} W(k) \; | \;
  a_{\psi_a} = \big(\frac{\pi^e}{p}\big)^f ~^{F^f} a_{\psi_a}\}.
    \end{equation}
For the last equation we may write
$\eta a_{\psi_a} =  ~^{F^f} (\eta a_{\psi_a})$ (cf. \eqref{KneunBa6e} for $\eta$) or,
equivalently, $\eta a_{\psi_a} \in O_F$. 
Therefore, using the expression (\ref{Cbanal21e}) for $C_{\mathcal{P}}$,
the restriction of $\tilde{\beta}$ to $C_{\mathcal{P}}$ multiplied by $\eta$
gives the desired bilinear form
\begin{equation}
  \begin{aligned}
   \phi\colon C_{\mathcal{P}} \times C_{\mathcal{P}}&  \longrightarrow O_F\\
    (x_{\psi_a} , x'_{\psi_a})&  \mapsto  \eta \tilde{\beta}_{\psi_a}(x_{\psi_a}, x'_{\psi_a}) .
    \end{aligned}
  \end{equation}

Now let $K/F$ be unramified. In this case, in the decomposition $P=\oplus_{\psi\in\Psi} P_\psi$, the summands $P_{\psi_1}$ and $P_{\psi_2}$ are orthogonal, unless $\psi_1=\bar\psi_2$. The $O_K$-module $C_\CP$ is, in this case, given by
\begin{equation}
C_\CP=\{ x=(x_\psi)\in P\mid V^{-1}\pi_r x_\psi=x_{\psi\sigma} \},
\end{equation}
where we recall the element $\pi_r$ from \eqref{pirunr}. After fixing $\psi_a$, we can write 
\begin{equation}\label{banalunra}
C_\CP=\{(x_{\psi_a}, x_{\bar{\psi}_a})\in P_{\psi_a}\oplus P_{\bar{\psi_a}}\mid V^{-f}\pi^{ g} x_{\psi_a}=x_{\bar{\psi}_a}, V^{-f}\pi^{\bar g} x_{\bar{\psi}_a}=x_{{\psi}_a} \} ,
\end{equation}
where $g=a_{\psi_a}+a_{\psi_a\sigma}+\cdots +a_{\psi_a\sigma^{f-1}} $ and
$\bar g=a_{\bar{\psi}_a}+a_{\bar{\psi}_a\sigma}+\cdots +a_{\bar{\psi}_a\sigma^{f-1}} $.
Using the expression \eqref{banalunra} for $C_\CP$, we may write 
\begin{equation}\label{altfunra}
  \begin{aligned}
  \phi: \quad C_{\mathcal{P}} \times C_{\mathcal{P}} & \longrightarrow  O_F(\rho)\\
  (x_{\psi_a}+x_{\bar{\psi}_a} , y_{\psi_a}+y_{\bar{\psi}_a} ) & \mapsto 
  \eta \tilde{\beta}_{\psi_a}(x_{\psi_a}, y_{\bar{\psi}_a})
  + \eta \tilde{\beta}_{\bar{\psi}_a}(x_{\bar{\psi}_a}, y_{{\psi}_a}) .
    \end{aligned}
\end{equation}
We have for arbitrary elements $x_{\psi} \in P_{\psi}$ and
$y_{\bar{\psi}} \in P_{\bar{\psi}}$ that
\begin{displaymath}
  \tilde{\beta}_{\bar{\psi}}(V^{-f} \pi^{g} x_{\psi},
  V^{-f} \pi^{\bar{g}} y_{\bar{\psi}}) =
  \frac{\pi^{ef}}{p^f} ~^{F^f}\tilde{\beta}_{\psi}(x_{\psi}, y_{\bar{\psi}}). 
  \end{displaymath}
If $x = x_{\psi_a}+x_{\bar{\psi}_a}$ and $ y = y_{\psi_a}+y_{\bar{\psi}_a}$ in $C_{\mathcal{P}}$
the last formula becomes
\begin{displaymath}
  \tilde{\beta}_{\bar{\psi}}(x_{\bar{\psi}}, y_{\psi}) = 
  \frac{\pi^{ef}}{p^f} ~^{F^f}\tilde{\beta}_{\psi}(x_{\psi}, y_{\bar{\psi}}). 
\end{displaymath}
By the formula (\ref{altfunra}) for $\phi$ we obtain
\begin{displaymath}
\phi(x,y) =  ~^{F^f} \phi(x,y).
  \end{displaymath}
This shows again that $\phi(x, y) \in O_F$, cf. (\ref{Cbanal24e}).  

Finally let $K = F \times F$. We use the notation of the last Remark. 
We obtain 
\begin{displaymath}
  C_{\mathcal{P}_i} = \{x \in P_i \; | Vx = \pi_{r,i} x \; \}, 
\end{displaymath}
where we recall the element $\pi_{r,i}$ from \eqref{zerlegt9e}
We have the decomposition
\begin{displaymath}
C_{\mathcal{P}} = C_{\mathcal{P}_1} \oplus C_{\mathcal{P}_2}
\end{displaymath}
If we fix $\theta_0 \in \Theta = \Hom_{\text{$\BQ_p$-Alg}}(F^t, \bar{\mathbb{Q}}_p)$, the
natural projection $P_i \rightarrow P_{i, \theta_0}$ defines an isomorphism
\begin{equation}\label{banalsplit1e}
  C_{\mathcal{P}_i} = \{x \in P_{i,\theta_0} \; | V^f x = \pi^{a_i} x \; \},
    \quad \text{for}\; i = 1,2, 
\end{equation}
where $a_i = \sum_{\theta} a_{\theta_i}$, cf. \eqref{sumsplit}.
The bilinear form $\tilde{\beta}$ induces by restriction 
\begin{displaymath}
  \tilde{\beta}_{\theta_0}: P_{1,\theta_0} \times P_{2,\theta_0} \rightarrow
  O_F \otimes_{O_{F^t}, \tilde{\theta}_0} W(k). 
\end{displaymath}
In the notation of (\ref{banalsplit1e}) we obtain
\begin{displaymath}
  \begin{array}{lcc} 
  \phi: C_{\mathcal{P}_1} \times C_{\mathcal{P}_2} & \rightarrow & O_F,\\
  (x_1, x_2) & \mapsto & \eta \tilde{\beta}_{\theta_0}(x_1, x_2).  
    \end{array}
  \end{displaymath}
This determines $\phi$ on $C_{\mathcal{P}}$ since the subspaces
$C_{\mathcal{P}_i}$ for $i = 1,2$ are isotropic, cf. Remark \ref{etformsplit}.
  \end{remark}

\begin{proposition}\label{CFbanal5p}
  Let $r$ be banal. Let $k \in \Nilp_{O_{E'}}$ be an algebraically closed
  field. Let
  $(\mathcal{P}, \iota, \beta)$ and $ (\mathcal{P}^{+}, \iota^{+}, \beta^{+})$ be two objects in  $\mathfrak{d}\mathfrak{P}_{r,k}^{\rm pol}$.
  \begin{altenumerate}

 \item If $K/F$ is split, then there exists a quasi-isogeny
  \begin{equation}\label{CFbanal2e}
    (\mathcal{P}, \iota, \beta) \rightarrow
    (\mathcal{P}^{+}, \iota^{+}, \beta^{+}) . 
    \end{equation}

\item Let $K/F$ be a field extension. Then there exists a quasi-isogeny    
(\ref{CFbanal2e}) iff
$\inv(\mathcal{P}, \iota, \beta)=\inv(\mathcal{P}^{+}, \iota^{+}, \beta^{+})$,  
cf. (\ref{Kneun3e}).

\item Let $K/F$ be an unramified field extension. If $\beta$ is a polarization
of height $2fh$ with $h \in \{0 , 1\}$ then
$\inv^r (\mathcal{P}, \iota, \beta) = (-1)^h$. For a given $h$, there 
exists $(\mathcal{P}, \iota, \beta)$ with these properties. 
\end{altenumerate}
\end{proposition} 
\begin{proof}
  To prove the first assertion, we may apply the polarized contraction functor
  $\mathfrak{C}_{r, k}^{\rm pol}$ of Theorem \ref{KneunBa4p}. We choose
  an arbitrary isomorphism $\alpha_1$ of the $F$-vector spaces
  $C_{\mathcal{P}_1} \otimes \mathbb{Q}$ and
  $C_{\mathcal{P}^{+}_1} \otimes \mathbb{Q}$. Since $\phi$, resp. $\phi^{+}$,
  define dualities of these spaces with $C_{\mathcal{P}_2} \otimes \mathbb{Q}$,
  resp. $C_{\mathcal{P}^{+}_2} \otimes \mathbb{Q}$, we can extend $\alpha_1$
  to an isomorphism
  $\alpha: (C_{\mathcal{P}}, \phi)\otimes\BQ \rightarrow (C_{\mathcal{P}^{+}}, \phi^{+})\otimes\BQ$.

  If $K/F$ is a field extension, we conclude by Proposition \ref{KneunBa1l} 
  that the equality $\inv(\mathcal{P}, \iota, \beta)=\inv(\mathcal{P}^{+}, \iota^{+}, \beta^{+})$ is equivalent to the equality $\inv((C_{\mathcal{P}}, \iota, \phi) \otimes \mathbb{Q})=\inv((C_{\mathcal{P}^{+}}, \iota^{+}, \phi^{+}) \otimes \mathbb{Q})$. 
  Therefore, by Definition \ref{invdet},  these  anti-hermitian $K$-vector spaces  are isomorphic, which proves our assertion.

  Finally we prove the last assertion. We consider the bilinear form
  $\beta_{\psi}: P_{\psi} \times P_{\bar{\psi}} \rightarrow W(k)$. If we choose 
  a $W(k)$-basis of $P_{\psi}$ and $P_{\bar{\psi}}$, we can speak of
  $\ord_p \det_{W(k)} \beta_{\psi}$. This number is independent of $\psi$ and
  equals $h$. Let $\beta(\psi)$ be the restriction of $\beta$ to
  $P_{\psi} \oplus P_{\bar{\psi}}$. We obtain $\ord_p \det_{W(k)} \beta(\psi) = 2h$.
  Recall $\tilde\beta$, cf. \eqref{Cbanal-t}. Let $\tilde{\beta}(\psi)$ the restriction of $\tilde{\beta}$, 
  \begin{displaymath}
    \tilde{\beta}(\psi): (P_{\psi} \oplus P_{\bar{\psi}}) \times
    (P_{\psi} \oplus P_{\bar{\psi}}) \rightarrow
    O_F \otimes_{O_{F^t}, \tilde{\psi}} W(k) .
  \end{displaymath}
  Then we have
  \begin{equation}\label{banalsplit4e}
    \ord_p \textstyle{\det_{W(k)}} \beta(\psi) = \ord_{\pi}
    \det_{O_F \otimes_{O_{F^t}, \tilde{\psi}} W(k)} \tilde{\beta}(\psi) .
    \end{equation}
 Indeed, the function $\mathbf{t}(a) = \Trace_{O_F/O_{F^t}} \vartheta^{-1}a$,
  for $a \in O_{F}$, where $\vartheta \in O_F$ is the different of
  $F/\mathbb{Q}_p$, defines for an arbitrary
  $O_F \otimes_{O_{F^t}, \tilde{\psi}} W(k)$-module $U$ an isomorphism of
  $O_F \otimes_{O_{F^t}, \tilde{\psi}} W(k)$-modules,
  \begin{equation}\label{banalsplit3e}
    \Hom_{O_F \otimes_{O_{F^t}, \tilde{\psi}} W(k)} (U, O_F \otimes_{O_{F^t}, \tilde{\psi}} W(k))
    \isoarrow \Hom_{W(k)}(U, W(k)), \quad \tilde{\alpha} \mapsto
    \tilde{\alpha} \circ \mathbf{t}.
  \end{equation}
  We apply this to $U = P_{\psi} \oplus P_{\bar{\psi}}$.
  If we regard $\tilde{\beta}$ as a homomorphism of $U$ to the left hand side
  of (\ref{banalsplit3e}) and $\beta$ as a homomorphism from $U$ to the
  right hand side, they correspond to each other. Therefore the cokernels
  of these two homomorphisms are isomorphic and have the same length.
  This shows (\ref{banalsplit4e}).
  By  Remark \ref{CPexplban}, we have for each $\psi$ an isomorphism
  $C_{\mathcal{P}} \otimes_{O_F} (O_F \otimes_{O_{F^t}, \tilde{\psi}} W(k) )\cong \mathcal{P}$.
  Since $\phi$ coincides with the restriction of $\tilde{\beta}(\psi)$ up to
  a unit,  we conclude that $\ord_{\pi} \det_{O_F} \phi = 2h$. 
  By Lemma \ref{Kneun0l}, we have 
  $\inv (C_{\mathcal{P}}, \iota, \phi) = (-1)^h$.
  By Proposition \ref{KneunBa1l} we are done. 
\end{proof}

    \section{The alternative moduli problem revisited}\label{s:altmodprob}
    In this section we give another  proof of the main result of
    \cite{KRalt} which gives an alternative interpretation of the Drinfeld
    moduli space of special formal $O_D$-modules in the case of a
    quaternion division algebra $D$ over a $p$-adic local field $F$. We also prove a refinement concerning  descent data.  The original proof was already simplified by Kirch \cite{Ki}, but the argument here is different and is based on the theory of displays.

    \subsection{Special formal $O_D$-modules}\label{ss:sfmod} 
We fix the finite extension $F$ of $\BQ_p$, with
uniformizer $\pi$ and residue field $\kappa_{F}$. Let $R$ be an
$O_F$-algebra. Let $(X, \iota)$ be a $p$-divisible group
over $R$ with a strict action $\iota: O_F \longrightarrow \End X$.
A relative polarization of $X$ is a relative
polarization of the display of $X$. Here, by a relative polarization we mean one with respect to $O_F$, cf.
Definition \ref{LTD3d}. 

If $R = k$ is a perfect field, we may work with the associated
$\mathcal{W}_{O_F}(k)$-Dieudonn\'e module $(M, F, V)$ of $X$. It is obtained
from the display of $X$ by the Ahsendorf functor
$\mathfrak{A}_{O_F/\mathbb{Z}_p,R}$, cf. Remark \ref{Adorf1r}. In this language, 
a relative polarization is a $\mathcal{W}_{O_F}(k)$-alternating pairing
\begin{displaymath}
  \psi: M \times M \longrightarrow \mathcal{W}_{O_F}(k), 
\end{displaymath}
such that
\begin{equation}\label{Pol15e}
  \psi(Fx, Fy) = \pi ~^{F}\psi(x,y), \quad x, y \in M, 
\end{equation}
cf. (\ref{bilin2e}).     
If the bilinear form $\psi$ is perfect,  we will say that the polarization
is principal.

We denote by $D$ the quaternion division algebra with center $F$. Let
$F' \subset D$ be a quadratic unramified extension of $F$. Let
$O_D \subset D$ be the ring of integers. Recall that a \emph{special formal
$O_D$-module} $X$ \index{special formal
$O_D$-module} over $R$ is a $p$-divisible group $X$ over $R$ of
height $[D:\mathbb{Q}_p]$ with an action $\iota: O_D \longrightarrow \End X$
such that the restriction of $\iota$ to $O_F$ is strict and such that
$\Lie X$ is locally on $\Spec R$ a free $O_{F'} \otimes_{O_F} R$-algebra, cf. \cite{Dr}.
One can check that this condition is independent of the choice of $F'$.       
    
    \begin{proposition}\label{Pol3p} Let $k$ be an algebraically closed
      field of characteristic $p$ with an
      $O_F$-algebra structure $O_F \longrightarrow k$.  Let $F'$ be  an unramified
      quadratic extension $F'$ of $F$. We denote  by $F$ and $V$ the Frobenius and the
      Verschiebung acting on $W_{O_F}(k)$.   

      Let $(M, F, V)$ be a $\mathcal{W}_{O_F}(k)$-Dieudonn\'e module (see
      Definition \ref{Adorf1ex}) of height $4$ and dimension $2$ which is
      endowed with an $O_F$-algebra homomorphism  
      \begin{displaymath}
        \iota:  O_{F'} \longrightarrow \End (M, F, V). 
      \end{displaymath}
      Assume that $\iota$ makes  $M/VM$ into a  free module  of rank $1$ over  $\kappa_{F'} \otimes_{\kappa_{F}} k$. Then there exists a principal relative polarization $\psi$ on $(M, F, V)$
      such that  
      \begin{equation}\label{Pol20e}
        \psi(\iota(u)x, y) = \psi(x, \iota(u)y), \quad \text{for} \; u \in
        O_{F'}, \; x, y \in M
      \end{equation}
Any other relative polarization $\phi$ of $(X, \iota)$ with the property
      (\ref{Pol20e}) (with $\psi$ replaced by $\phi$) is of the form  
      \begin{displaymath}
        \phi(x, y) = \psi(\iota(c)x, y)
      \end{displaymath}
      for some element $c \in O_{F'}$. 
    \end{proposition}
    \begin{proof}
      We choose an embedding $O_{F'} \longrightarrow W_{O_F}(k)$. We set, for
      $i\in \mathbb{Z}/2\mathbb{Z}$, 
      \begin{equation*}
        M_i = \{x \in M \; | \; \iota(u) x = ~^{F^{i}}\!u x, \; \text{for}
        \; u \in O_{F'} \}.  
      \end{equation*}
      We have the decomposition
      \begin{equation}\label{decunr} 
        M = M_0 \oplus M_1. 
      \end{equation}
      The operators $F$ and $V$ are of degree $1$. The $k$-vector spaces
      $M_0/VM_1$ and $M_1/VM_0$ are by assumption both of rank $1$. 

      If $\psi$ is a bilinear form with the properties (\ref{Pol20e}), then
      the decomposition (\ref{decunr}) is orthogonal. We choose alternating
      perfect forms $\tilde{\psi}_0$ resp. $\tilde{\psi}_1$ on the free
      $W_{O_F}(k)$-modules $M_0$ resp. $M_1$ of rank 2. These forms are unique
      up to a unit in $W_{O_F}(k)$. By this uniqueness we find an
      equation of the form
      \begin{equation}\label{Pol22e}
        ~^{F^{-2}}\tilde{\psi}_0(F^{2}x_0, F^{2}x'_0) = \xi \pi^{2}
        \tilde{\psi}_0 (x_0, x'_0) , \quad  \; \xi \in
        W_{O_F}(k), \;\text{ for all } x_0, x'_0 \in M_0. 
      \end{equation}
      By assumption we  have $\ord_{\pi} \det (F^2 |M_0) =2$. Comparing the
      determinants on both sides of (\ref{Pol22e}), we conclude that $\xi$ is
      a unit. Since $k$ is algebraically closed we may write
      \begin{displaymath}
        \xi = ~^{F^{-2}}\eta \eta^{-1}.
      \end{displaymath} 
      Replacing $\tilde{\psi}_0$ by $\psi_0 := \eta \tilde{\psi}_0$ we may
      assume that we have $\xi = 1$ in equation (\ref{Pol22e}). 

      With the same argument as before we find an equation
      \begin{displaymath}
        ~^{F^{-1}} \psi_0(Fx_1, Fx'_1) = \xi_1 \pi \tilde{\psi}_1(x_1,
        x'_1), \quad \xi_1 \in W_{O_F}(k), \; \text{ for all }x_1, x'_1 \in M_1. 
      \end{displaymath} 
      Comparing the determinants we see that $\xi_1 \in W_{O_F}(k)$ is a
      unit. We set $\psi_1 = \xi_1 \tilde{\psi}_1$ and
      $\psi = \psi_0 \oplus \psi_1$ (orthogonal sum). Then $\psi $ satisfies
      (\ref{Pol15e}).  To prove (\ref{Pol20e}) it suffices to show that 
      \begin{displaymath}
        \psi_i(\iota(u) x_i, x'_i) = \psi_i(x_i, \iota(u) x'_i) \quad \text{for}
        \; i = 0, 1.
      \end{displaymath} 
      This is trivial from the definition of $(M_i, \psi_i)$. 

      If we have a second $\phi$ satisfying \eqref{Pol20e}, we find
      $c \in W_{O_F}(k)$ such that
      \begin{displaymath}
        \phi_0 = c \psi_0.
      \end{displaymath}
      Since both sides of this equation satisfy (\ref{Pol22e}) with
      $\xi = 1$ we obtain $~^{F^{2}}c = c$. Therefore we have
      $c \in O_{F'} \subset W_{O_F}(k)$. We obtain that
      $\phi(x, y) = \psi(\iota(c) x, y)$.   
    \end{proof} 

    \begin{corollary}
      Let $(N, \iota)$ be a second $\mathcal{W}_{O_F}(k)$-Dieudonn\'e module
      of height $4$ and dimension $2$ with an action of $O_{F'}$ such that
      $N/VN$ is a free $\kappa_{F'} \otimes_{\kappa_F} k$-module of rank $1$.  Let $\rho: N \otimes \mathbb{Q} \longrightarrow M \otimes \mathbb{Q}$ be a
      quasi-isogeny of height $0$. Then
      \begin{equation}\label{Pol20e1} 
        \psi(\rho (z), \rho (w)), \quad z,w \in N
      \end{equation}
      is a perfect bilinear form on $N$.
    \end{corollary}
    \begin{proof}
      Let $\psi_N$ be a perfect alternating form on $N$ given by 
      Proposition \ref{Pol3p}, and let $\psi_{N_0}$ be its restriction to $N_0$. This
      form differs from the form (\ref{Pol20e1}) restricted to $N_0$ by a
      factor in $\zeta \in F'$. Since $\rho$ has height $0$ we conclude that
      $\zeta$ is a unit.  
    \end{proof}

    Let $K/F$ be a ramified quadratic extension of $F$ generated by  a prime
    element $\Pi \in K$ such that $\Pi^2 = -\pi$.
    Let $\tau \in \Gal(F'/F)$ be the Frobenius automorphism. Let 
    \begin{displaymath}
      O_D = O_{F'}[\Pi],
    \end{displaymath}
    such that the following relations hold: 
    \begin{displaymath}
      \Pi u = \tau(u) \Pi, \quad \Pi^{2} = -\pi, \quad u \in O_{F'}.
    \end{displaymath}
    Then $O_D$ is the maximal order in the quaternion division algebra over $F$. 

    We have $O_K = O_F[\Pi] \subset O_D$.  We consider on $O_D$ the involution:
    \begin{equation}\label{Pol1e}
      d = u + v\Pi \mapsto d' = u - \Pi v, \quad u,v \in O_{F'}. 
    \end{equation}
    It is trivial on $O_{F'}$ and induces the conjugation of $O_K$ over $O_F$.

    \begin{proposition}\label{Pol1p} 
      Let $k$ be an algebraically closed field of characteristic $p$ which
      is endowed with an algebra structure $O_F \longrightarrow k$. Let
      $X$ be a special formal $O_D$-module over $k$. Let $(M,F,V)$ be the
      $\mathcal{W}_{O_F}(k)$-Dieudonn\'e module of $X$.   Then there exists a principal relative polarization 
      \begin{displaymath}
        \psi:  M \times M \longrightarrow W_{O_F}(k),
      \end{displaymath} 
      on $X$ such that 
      \begin{equation}\label{Pol2e}
        \psi(\iota(d)x_1, x_2) = \psi(x_1, \iota(d')x_2).
      \end{equation}
  Any other polarization with the property (\ref{Pol2e}) is of the form
      $u\psi$, with $u \in O_F$.
    \end{proposition}
    \begin{proof}
      We take $\psi$ as in Proposition
      \ref{Pol3p}. Then we consider the alternating bilinear form 
      \begin{displaymath}
        \psi_1(x,y) = \psi(\iota(\Pi)x, \iota(\Pi)y), \quad x,y \in M.
      \end{displaymath} 
      Then $\psi_1$ is by the uniqueness part of  Proposition \ref{Pol3p} of the form 
      \begin{displaymath}
        \psi_1(x, y) = \psi (\iota(c)x, y), \quad c \in O_{F'}.
      \end{displaymath}
      If we apply the last equations to $\pi^{2}\psi(x, y) =
      \psi(\iota(\Pi)x, \iota(\Pi)y)$, we obtain $c \tau(c) =
      \pi^2$. Therefore $c$ is divisible by $\pi$. We write $c = a \pi$ for
      some unit $a \in O_{F'}$ with $a \tau(a) = 1$. We write $a = u \tau(u)^{-1}$ by Hilbert 90
      and consider the form
      \begin{displaymath}
        \psi_2(x, y) = \psi(\iota(u)x, y).
      \end{displaymath}
      Then we have
      \begin{equation*}
        \begin{aligned}
          \psi_2(\iota(\Pi)x, \iota(\Pi)y) &= \psi(\iota(u)\iota(\Pi)x,
          \iota(\Pi)y) = \psi(\iota(\Pi)\iota(\tau(u))x, \iota(\Pi)y)= \\
          &=\psi(\iota(c) \iota(\tau(u))x, y) = \psi(\iota(\pi u)x, y) = \pi
          \psi_2(x,y). 
        \end{aligned}
      \end{equation*}
      Therefore $\psi_2$ satisfies the requirements \eqref{Pol2e}. The
      uniqueness assertion is proved as before.  
    \end{proof}

Let $(X, \iota)$ be a special formal $O_D$-module over a ring
$R \in \Nilp_{O_F}$. We denote by $\mathcal{P}$ the corresponding
$\mathcal{W}_{O_F}(R)$-display. Let $\mathcal{P}^{\vee}$ be the dual 
$\mathcal{W}_{O_F}(R)$-display. Then $\iota$ induces a homomorphism
$\iota^{\vee}: O_D^{\mathrm{opp}} \rightarrow \End \mathcal{P}^{\vee}$. Let
$\iota': O_D \rightarrow \End \mathcal{P}^{\vee}$ be the composite of 
$\iota^{\vee}$ with the involution
$O_D \rightarrow O_D^{\mathrm{opp}}$, cf.  (\ref{Pol1e}). By Theorem \ref{maindisp}, 
$(\mathcal{P}^{\vee}, \iota')$ corresponds to
a special formal $O_D$-module $(X',\iota')$. 

Let $\psi$ be a relative polarization of $X$, i.e.
a $\mathcal{W}_{O_F}(R)$-polarization  
$\psi: \mathcal{P} \times \mathcal{P} \rightarrow \mathcal{P}_m$.
We assume that $\psi$ induces on $O_D$ the involution (\ref{Pol1e}).
In other words, (\ref{Pol2e}) is satisfied in this context, i.e., for 
$x_1 ,x_2 \in P$ and $d \in O_D$. To give such a polarization $\psi$ is by
(\ref{dual4e}) the same thing as to give an isogeny of special formal
$O_D$-modules $\lambda_{\psi}: (X,\iota) \rightarrow (X',\iota')$ which is
anti-symmetric with respect to the duality $X \mapsto X'$. 

\begin{definition}\label{Pol3d}
Let $(X, \iota)$ be a special formal $O_D$-module over $R \in \Nilp_{O_F}$.  
A \emph{Drinfeld polarization} \index{Drinfeld polarization} is a principal relative polarization on $X$
which induces on $O_D$ the involution (\ref{Pol1e}). Alternatively, a Drinfeld
polarization is given by an anti-symmetric isomorphism of special formal
$O_D$-modules $\lambda: (X,\iota) \rightarrow (X',\iota')$. 
\end{definition}

    \begin{proposition}\label{Pol2p}
      Let $\mathbb{X}$ be a  special formal $O_D$-module over $\bar{\kappa}_F$. 
   Let $S$ be a connected scheme over $\Spf O_{\breve{F}}$.
   In other words, $p$ is locally nilpotent on $S$. We set
   $\bar{S} = S \times_{\Spf O_{\breve{F}}} \Spec \bar{\kappa}_F$. Let $(X, \iota)$
   be a special formal $O_D$-module over $S$ such that there exists a
   quasi-isogeny of special formal $O_D$-modules 
   \begin{equation}\label{Pol9e}
     \mathbb{X} \times_{\Spec \bar{\kappa}_F} \bar{S} \longrightarrow X \times_S
     \bar{S}. 
   \end{equation}
   
 Then there is a Drinfeld polarization $\psi$ on $X$. Any other relative
 polarization on $X$ which induces the involution $d \mapsto d'$ is of the
 form $f\psi$ for some $f\in O_F$.   

 In particular, $\psi$ is, up to a factor in $F^{\times}$, compatible with
 a Drinfeld polarization on $\mathbb{X}$ by the quasi-isogeny (\ref{Pol9e})
    \end{proposition}
    \begin{proof}
We begin with a $\bar{\kappa}_F$-scheme $S$ which is not necessarily
connected. We fix a Drinfeld polarization $\lambda_{\mathbb{X}}$ on
$\mathbb{X}$. We will show that for each point of $S$ there is an open
neighbourhood $U$ and an integer $c$ such that the quasi-polarization induced
by $\pi^c \lambda_{\mathbb{X}}$ via (\ref{Pol9e}) on $X$ is a Drinfeld
polarization over $U$. For this we can assume that $S = \Spec R$. Because $X$
is the quotient of $\mathbb{X} \times_{\Spec \bar{\kappa}_F} \times S$ by a finite
locally free subgroup scheme, the quasi-isogeny (\ref{Pol9e}) is defined
over a subalgebra $R_0 \subset R$ which is finitely generated over
$\bar{\kappa}_F$. Therefore we may assume that $R = R_0$ and, in particular,
that $R$ is noetherian.

Once we know the existence of $c$, it follows immediately that the
quasi-polarization induced by $\pi^c \lambda_{\mathbb{X}}$ via (\ref{Pol9e}) 
is a polarization on an open and closed subset $U \subset S$ which contains
the point we started with. This will prove the Proposition in the case
where $S$ is an $\bar{\kappa}_F$-scheme. 

To prove the existence of $c$, we recall some generalities from \cite{Dr} which
are formulated there for Cartier modules. Let $R \in \Nilp_{O_{\breve{F}}}$.  We fix
an embedding $O_{F'} \longrightarrow O_{\breve{F}}$. From this we obtain
homomorphisms $O_{F'} \longrightarrow R$ and
$\lambda: O_{F'} \longrightarrow W_{O_{F}}(O_{F'}) \longrightarrow W_{O_F}(R)$.
Let $\bar{\lambda}$ be the composite of $\lambda$ with the conjugation of
$F'/F$.
      
Let $\mathcal{P}$ be the $\mathcal{W}_{O_F}(R)$-display of $X$. The
action of $O_D$ on $\mathcal{P}$ is also denoted by $\iota$.
We have the decompositions 
\begin{equation}\label{decunr2e}
        P = P_0 \oplus P_1, \quad Q = Q_0 \oplus Q_1
\end{equation}
such that for $a \in O_{F'} \subset O_D$ the action of $\iota(a)$ on
$P_0$ is multiplication by $\lambda(a)$ and the action on $P_1$ is
multiplication by $\bar{\lambda}(a)$, and $Q_i = Q \cap P_i$.
We regard \eqref{decunr2e} as a $ \mathbb{Z}/2 \mathbb{Z}$-grading. Then
$F, \Pi, \dot{F}$ are all homogeneous of degree $1$, 
\begin{equation*}
  F: P_i \longrightarrow P_{i+1}, \quad  \Pi: P_i \longrightarrow P_{i+1},
  \quad  \dot{F}: Q_i \longrightarrow P_{i+1}.  
\end{equation*}
Let $\psi: \mathcal{P} \times \mathcal{P} \rightarrow \mathcal{P}_m$ be a
polarization which induces on $O_D$ the involution (\ref{Pol1e}). Because
this involution is trivial on $O_{F'}$ we obtain that $P_0$ is orthogonal to
$P_1$ with respect to $\psi$. Therefore $\psi$ is given by two
alternating $W_{O_F}(R)$-bilinear forms
\begin{displaymath}
  \psi_0: P_0 \times P_0 \rightarrow W_{O_F}(R), \quad
  \psi_1: P_1 \times P_1 \rightarrow W_{O_F}(R) .
\end{displaymath}
In our case $Q_0/I_{O_F}(R) P_0 \subset P_0/I_{O_F}(R) P_0$ is a direct summmand
of rank $1$ and therefore an arbitrary alternating $R$-bilinear form on
$P_0/I_{O_F}(R) P_0$ is zero on this direct summand. This implies that for an
arbitrary alternating form $\psi_0$ the inclusion
$\psi_0(Q_0, Q_0) \subset I_{O_F}(R)$ holds. The same remark applies to $\psi_1$. 
If $\psi$ is a polarization the equation 
\begin{equation}\label{decunr3e}
\psi_1(\dot{F} y, \dot{F} y') = ~^{\dot{F}}\psi_0(y,y')
\end{equation}
holds for $y, y' \in Q_0$. 
Since $\dot{F}: Q_0 \rightarrow P_1$ is a Frobenius-linear epimorphism, we see
that $\psi_1$ is uniquely determined by $\psi_0$. In fact, we can construct
$\psi_1$ from $\psi_0$ as follows. We take a normal decomposition
$P_0  = L_0 \oplus T_0$, $Q_0  = L_0 \oplus I_{O_F}(R) T_0$. The linearizations 
of $F$ and $\dot{F}$ define an isomorphism
\begin{displaymath}
  \dot{F}^{\sharp} \oplus F^{\sharp}: W_{O_F}(R) \otimes_{F,W_{O_F}(R)} L_0 \; \oplus 
  \; W_{O_F}(R) \otimes_{F,W_{O_F}(R)} T_0 \overset{\sim}{\longrightarrow} P_1.   
  \end{displaymath}
Therefore we can define a bilinear form $\psi_1$ on $P_1$ by the equations
\begin{displaymath}
  \begin{array}{cl}
    \psi_1(\xi\dot{F} l_0, \eta\dot{F} l'_0) = \xi \eta
    ^{\dot{F}}\!\!\psi_0(l_0, l'_0), & \xi, \eta \in W_{O_F}(R),
    \quad l_0, l'_0 \in L_0, \\[2mm] 
     \psi_1(\xi\dot{F} l_0, \eta F t_0) = \xi \eta
     ^{F}\!\!\psi_0(l_0, t_0), & t_0 \in T_0, \\[2mm]
     \psi_1(\eta Ft_0, \xi \dot{F} l_0) = \eta \xi 
     ^{F}\!\!\psi_0(t_0, l_0), &  \\[2mm]
     \psi_1(\xi Ft_0, \eta F t'_0) = \xi \eta \pi 
    ^{F}\!\!\psi_0(t_0, t'_0), & t'_0 \in T_0. \\[2mm]
    \end{array}
\end{displaymath}
One checks that with this definition of $\psi_1$ the identity (\ref{decunr3e}) holds. 
Therefore it makes sense to ask whether an alternating form $\psi_0$ on $P_0$
is a polarization.  

To show the existence of a principal $\psi$, we begin with the case where
$S = \Spec R$ and where $\mathcal{P}$ has a critical index
$i \in \mathbb{Z}/2\mathbb{Z}$. Assume that $i=0$ is critical,
i.e., the homomorphism 
\begin{equation}\label{Pol12e}
  \Pi: P_0/Q_0 \longrightarrow P_1/Q_1
\end{equation}
is zero. This implies that $\pi = -\Pi^2$ is zero in $R$. We consider the
composite
$\Phi: P_0\overset{\Pi}{\longrightarrow} Q_1 \overset{\dot{F}}{\longrightarrow} P_0$. 
We claim that $\Phi$ is a Frobenius-linear isomorphism. It is enough to show
that $\det \Phi \in W_{O_F}(R)$ is a unit. By base change, we
may assume that $R = k$ is a perfect field. Since $i = 0$ is critical,
we find $\Pi P_0  \subset Q_1 = VP_0$. Since $P_1/VP_0$ and
$P_1/\Pi P_0$ are $k$-vector spaces of dimension $1$ we obtain
$\Pi P_0 = V P_0$. Therefore $V^{-1} \Pi = \Phi$ is bijective and
therefore a Frobenius-linear isomorphism.
      
  For each $n \in \mathbb{N}$ we consider on the category of affine schemes
  $\Spec A \rightarrow S = \Spec R$ the functor 
  \begin{equation}\label{Pol10e}
    \underline{U}_0(n): \Spec A\mapsto (P_0\otimes_{W_{O_F}(R)} W_{O_F, n}(A))^{\Phi},
  \end{equation}
  where the RHS denotes invariants of the Frobenius-linearly extended operator $\Phi$.
  This functor is representable by a scheme which is finite and \'etale over
  $S$ of degree $\sharp (O_F/\pi^n O_F)^2$, cf. Proposition \ref{Cbanal1p}. 
  Moreover, the existence of the quasi-isogeny (\ref{Pol9e})
  implies that this scheme is a constant finite scheme. This means that
  $\underline{U}_0(n) = S \times U_0(n)$ where $U_0(n)$ is the set of sections
  of $\underline{U}_0(n) \rightarrow S$. We note that the category of finite
  constant sheaves on $S$ is equivalent to the category of finite sets because
  $S$ is connected.
  We conclude that $U_0(n)$ is an $O_F$-module isomorphic to
  $(O_F/\pi^nO_F)^2$. We set $U_0 = \projlim U_0(n)$. It is a free $O_F$-module
  of rank $2$. We obtain a canonical isomorphism $W_{O_F}(R)$-modules 
  \begin{displaymath}
    P_0 \cong W_{O_F}(R) \otimes_{O_F} U_0. 
  \end{displaymath}
Let $\phi$ be a relative polarization on $X$ which induces the
involution $d \mapsto d'$ on $O_D$. Then we obtain for $x, x' \in P_0$ 
\begin{displaymath}
  \phi_0(\dot{F} \Pi x, \dot{F} \Pi x') = ^{\dot{F}}\!\! \phi_1(\Pi x, \Pi x') 
  = ^{\dot{F}}\!\! \phi_0(x, -\Pi^2 x') = \pi\; ^{\dot{F}}\!\! \phi_0( x, x') =
  ^{F}\!\! \phi_0( x, x').   
  \end{displaymath}
Therefore the restriction of $\phi_0$ induces an alternating $O_F$-bilinear
form $\bar{\phi}_0: U_0 \times U_0 \rightarrow O_F$. The form $\bar{\phi}_0$ 
determines $\phi_0$ and then $\phi_1$, as we have seen above. Therefore
$\bar{\phi}_0$ determines $\phi$ uniquely. We conclude that any 
relative polarization $\phi$ of $\mathcal{P}$ that induces the given
involution on $O_D$ is of the form $\phi = f \psi$ for some $f \in F^{\times}$.  

Conversely, we start with a perfect alternating pairing
$\alpha: U_0 \times U_0  \longrightarrow O_F$  and extend it by base change
to a perfect bilinear pairing 
\begin{displaymath}
  \psi_0: P_0 \times P_0 \longrightarrow W_{O_F}(R). 
\end{displaymath}
As explained above, $\psi_0$ extends to an alternating bilinear form 
\begin{displaymath}
  \psi: P  \times P \longrightarrow W_{O_F}(R) ,
\end{displaymath}
such that $P_0$ and $P_1$ are orthogonal with respect to $ \psi$ and
(\ref{decunr3e}) holds. 

To prove that $\psi$ is a polarization we begin with the case where, moreover,
$R$ is a reduced ring. In this case $P \subset P \otimes \mathbb{Q}$. 
We may assume without loss of generality that both indices $i=0$ and $i=1$
are critical for $\mathbb{X}$. We choose a Drinfeld polarization 
$\psi_{\mathbb{X}}$ on $\mathbb{X}$. Let $\alpha_{\mathbb{X}}$ be the alternating
$O_F$-bilinear form on $U_{0, \mathbb{X}}$. 
The quasi-isogeny (\ref{Pol9e}) induces an isomorphism
\begin{displaymath}
U_{0,\mathbb{X}} \otimes \mathbb{Q} \rightarrow U_0 \otimes \mathbb{Q}.  
\end{displaymath}
Since both sides are two-dimensional $F$-vector spaces, the bilinear forms
$\alpha$ and $\alpha_{\mathbb{X}}$ differ by a factor in $F^{\times}$ under the 
isomorphism. Let $\phi$ be the relative quasi-polarization induced  by $\psi_{\mathbb{X}}$ on $X$ 
via (\ref{Pol9e}). We have seen that $\phi_0$ and $\psi_0$
differ by a factor in $F^{\times}$. Because for the bilinear forms
$\phi: P \otimes\mathbb{Q} \times P\otimes\mathbb{Q} \rightarrow W_{O_F}(R) \otimes \mathbb{Q}$
and $\psi$ the identity  (\ref{decunr3e}) holds, these bilinear forms also differ by the
same factor in $F^{\times}$. In particular $\psi$ inherits from $\phi$ the 
identities
\begin{displaymath}
    \psi(\dot{F} x, \dot{F} z) = ^{\dot{F}}\!\!\psi(x,y), \quad
    \psi(\Pi x, z) = \psi(x, \Pi z), \quad x,z \in P \otimes\mathbb{Q}.  
\end{displaymath}
This proves that $\psi: P \times P \rightarrow W_{O_F}(R)$ is a polarization.
To show that $\psi$ is perfect it is enough to show that
$\det \psi \in W_{O_F}(R)$ is a unit. This may be reduced to the case where
$R$ is a perfect field. In this case we know it by Proposition \ref{Pol1p}. 
This proves that $\psi$ a Drinfeld polarization on $X$. Since the restrictions
of $\psi$ and $\phi$ to $U_0 \otimes \mathbb{Q}$ differ by a factor in
$F^{\times}$, we conclude that there is $c \in \mathbb{Z}$ such that the
quasi-polarization $\pi^c \psi_{\mathbb{X}}$ induces via (\ref{Pol9e}) a Drinfeld
polarization on $X$.

We make a general remark on liftings before we continue.  
Let $\tilde{R} \rightarrow R$ be a surjection in $\Nilp_{O_{\breve{F}}}$ with
nilpotent kernel. Let $\tilde{X}$ be a special formal $O_D$-module over
$\tilde{R}$ and let $X$ be its base change to $R$.
Then any relative polarization $\psi$ of $X$
which induces on $O_D$ the given involution (\ref{Pol1e}) lifts uniquely to a
relative polarization on $\tilde{X}$ which also induces the given involution.
Indeed, to see this we may assume that the map $\tilde{R} \rightarrow R$ is an
$O_F$-$pd$-thickening, because any surjection with nilpotent kernel breaks up
into such thickenings. In this situation we apply Proposition \ref{GM2p}. Let
$\mathcal{P}$ resp. $\tilde{\mathcal{P}}$ be the displays of $X$ resp.
$\tilde{X}$. Then
$\psi: \mathcal{P} \times \mathcal{P} \rightarrow \mathcal{P}_m$ lifts to
$\tilde{\mathcal{P}}$ iff  
\begin{displaymath}
  \psi_{\mathrm{crys}} (\tilde{Q}/I_{O_F}(\tilde{R})\tilde{P},
  \tilde{Q}/I_{O_F}(\tilde{R})\tilde{P}) = 0, 
\end{displaymath}
and this lifting is unique. As we remarked above, this condition follows because 
$\tilde{Q}_0/I_{O_F}(\tilde{R})\tilde{P}_0$ and 
$\tilde{Q}_1/I_{O_F}(\tilde{R})\tilde{P}_1$ are locally free of rank $1$. 

Now we return to the case of a reduced ring $R$ such that $S = \Spec R$ is
connected, but we do not assume that $X$ has a critical index. We consider the
closed subscheme $S_0 = \Spec R/\mathfrak{a}_0$ where $i=0$ is critical, i.e.
where the homomorphism (\ref{Pol12e}) is zero. Let
$S_1 = \Spec R/\mathfrak{a}_1$ be the closed subscheme where $i=1$ is
critical.  We obtain a fiber product diagram of rings
\begin{displaymath}
   \xymatrix{
     R/(\mathfrak{a}_0 \cap \mathfrak{a}_1) \ar[r] \ar[d] &
     R/\mathfrak{a}_0\ar[d]\\
     R/\mathfrak{a}_1 \ar[r] & R/(\mathfrak{a}_0 + \mathfrak{a}_1) .\\ 
   }
  \end{displaymath} 
Let $G$, $H$ be $p$-divisible groups over
$R/(\mathfrak{a}_0\cap\mathfrak{a}_1)$. Let $G_i$, $H_i$ be the restrictions of
these $p$-divisible groups to $R/\mathfrak{a}_i$ for $i = 1, 2$. 
It is easy to see that two homomorphisms $\gamma_i: G_i \rightarrow H_i$
which agree on $R/(\mathfrak{a}_0 + \mathfrak{a}_1)$ come from a unique
homomorphism $G \rightarrow H$. This implies that two Drinfeld polarizations
of $X$ over $R/\mathfrak{a}_1$ and $R/\mathfrak{a}_2$ which agree on
$R/(\mathfrak{a}_0 + \mathfrak{a}_1)$ are the restriction of a single Drinfeld
polarization of $X$ over $R/(\mathfrak{a}_0 \cap \mathfrak{a}_1)$.
Finally we can lift  a Drinfeld polarization from
$R/(\mathfrak{a}_0 \cap \mathfrak{a}_1)$ to $R$ because the ideal
$\mathfrak{a}_0 \cap \mathfrak{a}_1$ is nilpotent. Note that
$\mathfrak{a}_0 \mathfrak{a}_1 = 0$ because $\Pi^2 = -\pi = 0$. 

Let $\{U_r\}_{r \in \mathbf{r}}$ be the connected components of $S_0$ and let
$\{V_t\}_{t \in \mathbf{t}}$ be the connected components of $S_1$. These are
closed subschemes of $\Spec R$. We find
integers $c_r$ and $d_t$ such that $\pi^{c_r} \psi_{\mathbb{X}}$ induces a 
Drinfeld polarization of $X$ over $U_r$ and $\pi^{d_t} \psi_{\mathbb{X}}$ induces
a Drinfeld polarization of $X$ over $V_t$. If $U_r \cap V_t \neq \emptyset$
we obtain $c_r = d_t$. Since $S = \Spec R$ is connected we conclude that
$c_r = d_t = c$ is independent of $r$ and $t$. This shows that
$\pi^{c} \psi_{\mathbb{X}}$ defines a Drinfeld polarization of $X$ over $S$.  
If $\phi$ and $\psi$ are two relative polarizations of $X$ which induce the
given involution (\ref{Pol1e}), we have already seen that their restrictions
to $U_r$ resp. $V_t$ differ by a factor in $F^{\times}$. The same argument as
before shows that $\phi = f \psi$ holds over $S$  for some $f \in F^{\times}$.

Next we consider the case that $\Spec R$ is noetherian and connected. Since
the kernel $R \rightarrow R_{\red}$ has nilpotent kernel, the relative
polarizations of $X$ over $R_{\red}$ inducing (\ref{Pol1e}) on $O_D$ lift
uniquely to $R$. This shows that two polarizations over $R$ of this type
differ by a factor in $F^{\times}$ and that there exist Drinfeld polarizations
over $R$.

As we said at the beginning of the proof, this implies that for each connected
$\bar{\kappa}_F$-scheme $S$ the Proposition holds. The general case follows
because the kernel of $\mathcal{O}_S \rightarrow \mathcal{O}_{\bar{S}}$ is
nilpotent.
\end{proof} 
 
\begin{corollary}\label{Pol3c}
  With the assumptions of Proposition \ref{Pol2p}, let $\psi_{\mathbb{X}}$ be a Drinfeld
  polarization on $\mathbb{X}$. Moreover, let 
  \begin{equation*}
    \rho: \mathbb{X} \times_{\Spec \bar{\kappa}_F} \bar{S} \longrightarrow
    X \times_S \bar{S}
  \end{equation*}
  be a quasi-isogeny of height $0$. Then the relative quasi-polarization
  $\bar{\psi}$ on $X \times_S  \bar{S}$ induced by $\psi_{\mathbb{X}}$ is a
  Drinfeld polarization of $X \times_S \bar{S}$ that lifts to a Drinfeld
  polarization $\psi$ on $X$. 
    \end{corollary}
\begin{proof}
  In the notation introduced before Definition \ref{Pol3d}, $\psi_{\mathbb{X}}$
  induces an isomorphism
  $\lambda_{\mathbb{X}}: \mathbb{X} \rightarrow \mathbb{X}'$. By the definition
  of $\bar{\psi}$ we obtain a commutative diagram 
\begin{displaymath}
   \xymatrix{
     \mathbb{X} \times_{\Spec \bar{\kappa}_F} \bar{S}
     \ar[r]^{\rho} \ar[d]_{\lambda_{\mathbb{X}}} &
     X \times_S \bar{S} \ar[d]^{\lambda_{\bar{\psi}}}\\
     \mathbb{X}' \times_{\Spec \bar{\kappa}_F} \bar{S}  & \ar[l]_{\quad \rho'}
     X' \times_S  \bar{S} .\\ 
   }
\end{displaymath}
Since $\rho$ and its dual $\rho'$ have height zero, we conclude that
$\lambda_{\bar{\psi}}$ is a quasi-isogeny of height zero. On the other hand,
there exists by  Proposition \ref{Pol2p} a Drinfeld polarization $\phi$ on $X$.
Moreover, there is $f \in F^{\times}$ such that
\begin{displaymath}
\lambda_{\bar{\psi}} =  \lambda_{\bar{\phi}} \circ \iota(f). 
\end{displaymath}
Hence $\iota(f): X \times_S  \bar{S} \rightarrow X \times_S  \bar{S}$ is
a quasi-isogeny of height zero. Since $4 \ord_{\pi} f = \height_{O_F} \iota(f)$,
we conclude that $f$ is a unit in $O_F$. Therefore $\lambda_{\bar{\psi}}$  is
an isomorphism, and $\bar{\psi}$ is a Drinfeld polarization which lifts to
the Drinfeld polarization $\lambda_{\phi} \circ \iota(f)$ on $X$.  
\end{proof}  

Let us recall the Drinfeld moduli functor $\mathcal{M}_{\rm Dr}$ 
\index[NO]{MBD@$\mathcal{M}_{\rm Dr}$} on the
    category of schemes  $S$ over $ \Spf O_{\breve{F}}$. We will use the
    notation $\ov S=S\otimes_{\Spf O_{\breve F}} \Spec \ov\kappa_F$. 
    We fix a special formal
    $O_D$-module $(\mathbb{Y}, \iota_{\mathbb{Y}})$
    \index[NO]{YDA@$(\mathbb{Y}, \iota_{\mathbb{Y}})$} 
    over the $O_{\breve{F}}$-algebra $\ov\kappa_F$. We call $\BY$ a
    {\it framing object. } \index{framing object} By \cite{Dr} there is a quasi-isogeny of
    height $0$ between any two choices. For a scheme
    $S \longrightarrow \Spf O_{\breve{F}}$, a point of
    $\mathcal{M}_{\rm Dr}(S)$ consists of the following data up to isomorphism:
    \begin{enumerate} 
    \item[(1)] A special formal $O_D$-module $(Y, \iota)$ over $S$.
    \item[(2)] A quasi-isogeny of $O_D$-modules of height $0$ 
      \begin{equation}\label{Pol26e}
        \rho: Y \times_{S} \ov S  \longrightarrow \mathbb{Y}
        \times_{\Spec\bar{\kappa}_F} \ov S.
      \end{equation}
    \end{enumerate} 
    The functor is representable by  the $p$-adic formal
    $O_{\breve{F}}$-scheme
    $\wh{\Omega}_{F} \times_{\Spf O_F} \Spf O_{\breve F}$.  

    We define the functor $\mathcal{M}_{\text{Dr}}(i)$ by replacing in (2) height
    \index[NO]{MBE@$\tilde{\mathcal{M}}_{\text{Dr}}$} 
 $0$ by the condition $\height_{O_F} \rho = 2i$. We set 
    \begin{displaymath}
      \tilde{\mathcal{M}}_{\text{Dr}} =
      \coprod_{i \in \mathbb{Z}} \mathcal{M}_{\text{Dr}}(i).  
    \end{displaymath}
    
Let $(Y, \iota)$ be a special formal $O_D$-module. Let $u \in D^{\times}$.
Then we define a new special formal $O_D$-module $(Y^{u}, \iota^{u})$ by setting
    \begin{displaymath}
      Y^{u} = Y, \quad \iota^{u}(d) = \iota(u^{-1}du), \quad \text{for}\;
      d \in O_D. 
    \end{displaymath}
The multiplication $\iota(u): (Y^{u}, \iota^{u}) \longrightarrow (Y, \iota)$
is a quasi-isogeny of special formal $O_D$-modules. 
We obtain for each $i\in \mathbb{Z}$ an isomorphism of functors
    \begin{equation}\label{FuncMDr1e}
      u:  \mathcal{M}_{\text{Dr}}(i) \isoarrow \mathcal{M}_{\text{Dr}}(i+\ord_D u),
      \quad (Y, \rho) \mapsto (Y^{u}, \iota_{\mathbb{Y}}(u) \rho^{u}),  
    \end{equation}
This defines an action of $D^{\times}$ on $\tilde{\mathcal{M}}_{\text{Dr}}$. If
$u \in O_D^{\times}$, the multiplication by $\iota(u)$ defines an isomorphism
$\iota(u): (Y^{u}, \iota_{\mathbb{Y}}(u) \rho^{u}) \rightarrow  (Y, \rho)$.
Therefore the action of $D^{\times}$ factors through
$\ord_D: D^{\times} \rightarrow \mathbb{Z}$. We will call this action the
\emph{translation}.  \index{translation action}
 
    We endow $\tilde{\mathcal{M}}_{\text{Dr}}$ with a Weil descent datum
    relative to
    $O_{\breve{F}}/O_{F}$. Let $\tau \in \Gal(\breve{F}/F)$ be the Frobenius
    automorphism. Let $\varepsilon: O_{\breve{F}} \longrightarrow R$ be an algebra in
    $\Nilp_{O_{\breve{F}}}$. We denote by $R_{[\tau]}$\index[NO]{RAC@$R_{[\tau]}$}
    the ring $R$ with the new
    $O_{\breve{F}}$-algebra structure $\varepsilon\circ \tau$.
    The Frobenius $\tau$ induces
    $\bar{\tau}: \bar{\kappa}_F \longrightarrow \bar{\kappa}_F$.
    We have the Frobenius morphism
    \begin{equation}\label{WD1e}
      F_{\mathbb{Y},\tau}: \mathbb{Y} \longrightarrow \tau_{\ast}\mathbb{Y}.  
    \end{equation}
For a $\bar{\kappa}_F$-algebra
$\varepsilon: \bar{\kappa}_F \longrightarrow \bar{R}$, we set
$\phi(r) = r^{p^f}$ for $r \in \bar{R}$. This defines a $\bar{\kappa}_F$-algebra
homomorphism $\bar{R} \longrightarrow \bar{R}_{[\tau]}$.
    If we apply the functor $\mathbb{Y}$ we obtain (\ref{WD1e}). 
    We will define a morphism  
    \begin{equation}\label{WD2e}
      \omega_{\mathcal{M}_{\text{Dr}}}: \mathcal{M}_{\text{Dr}}(i)(R) \longrightarrow
      \mathcal{M}_{\text{Dr}}(i+1)(R_{[\tau]}).
    \end{equation}
    Let $(Y,\rho) \in \mathcal{M}_{\text{Dr}}(i)(R)$. We define $\rho'$ as the
    composite
    \begin{displaymath}
      Y_{R \otimes_{O_{\breve{F}}} \bar{\kappa}_F} \overset{\rho}{\longrightarrow}
      \varepsilon_{\ast} \mathbb{Y}
      \overset{\varepsilon_{\ast} F_{\mathbb{Y}, \tau}}{\longrightarrow}
      \varepsilon_{\ast} \tau_{\ast} \mathbb{Y}.
    \end{displaymath}
The image of $(Y, \rho)$ under (\ref{WD2e}) is by definition $(Y, \rho')$. 
Since $\height_{O_F} F_{\mathbb{Y}, \tau} = 2$, we obtain that
$\height_{O_F} \rho' = 2i + 2$. 
From (\ref{WD2e}) we obtain a Weil descent datum
\index[NO]{ZZWD@$\omega_{\mathcal{M}_{\text{Dr}}}$}
\begin{equation}\label{WD3e}
  \omega_{\mathcal{M}_{\text{Dr}}}: \tilde{\mathcal{M}}_{\text{Dr}}(R) \longrightarrow
  \tilde{\mathcal{M}}_{\text{Dr}}(R_{[\tau]}) 
\end{equation}
on the functor $\tilde{\mathcal{M}}_{\text{Dr}}$ (compare \cite{RZ}).
We introduce the notation
\begin{equation}\label{WD8e}
  \tilde{\mathcal{M}}_{\text{Dr}}^{(\tau)} = \tilde{\mathcal{M}}_{\text{Dr}}
  \times_{\Spf O_{\breve{F}}, \Spf \tau} \Spf O_{\breve{F}}.
  \end{equation}
Then we have
$\tilde{\mathcal{M}}_{\text{Dr}}^{(\tau)}(R) = \tilde{\mathcal{M}}_{\text{Dr}}(R_{[\tau]})$. We write (\ref{WD3e}) in the form
\begin{equation}
  \omega_{\mathcal{M}_{\text{Dr}}}: \tilde{\mathcal{M}}_{\text{Dr}} \longrightarrow
  \tilde{\mathcal{M}}_{\text{Dr}}^{(\tau)}. 
\end{equation}

The translation
$\Pi: \mathcal{M}_{\text{Dr}}(i) \rightarrow \mathcal{M}_{\text{Dr}}(i+1)$ is
an isomorphism. We use it to identify these functors. By Drinfeld's theorem
we obtain an isomorphism 

\begin{equation}\label{FuncMDr2e}
  \tilde{\mathcal{M}}_{\text{Dr}} \cong
  (\wh{\Omega}_{F} \times_{\Spf O_F} \Spf O_{\breve F} ) \times \mathbb{Z}.
\end{equation}
We denote by $\omega_{\tau}$\index[NO]{ZZWE@$\omega_{\tau}$} the action of $\tau$ via the second
factor on $\wh{\Omega}_{F} \times_{\Spf O_F} \Spf O_{\breve F}$.
\begin{proposition}
The  Weil descent datum $\omega_{\mathcal{M}_{\text{Dr}}}$ induces
on the right hand side of (\ref{FuncMDr2e}) the Weil descent datum
\begin{equation}\label{WD13e} 
   \omega_{\mathcal{M}_{\text{Dr}}}: \; (\xi, i) \mapsto (\omega_{\tau}(\xi), i+1). 
 \end{equation}
The translation functor is on the right hand side
$(\xi, i) \mapsto (\xi, i+1)$. 
  \end{proposition}
\begin{proof}
  Let $(Y, \rho) \in \mathcal{M}_{\text{Dr}}(R)$. Composing
  $\omega_{\mathcal{M}_{\text{Dr}}}$ with the translation we obtain a
  Weil-descent datum on $\mathcal{M}_{\text{Dr}}$,
  \begin{displaymath}
\alpha: \mathcal{M}_{\text{Dr}}(R) \rightarrow \mathcal{M}_{\text{Dr}}(R_{[\tau]}).  
  \end{displaymath}
  It associates to $(Y, \rho)$ the point $(Y^{\Pi^{-1}}, \rho_1)$, where
  $\rho_1$ is the composite
  \begin{equation}\label{Pol28e}
      Y^{\Pi^{-1}}_{R \otimes_{O_{\breve{F}}} \bar{\kappa}_F} \overset{\rho}{\longrightarrow}
      \varepsilon_{\ast} \mathbb{Y}^{\Pi^{-1}}
      \overset{\varepsilon_{\ast} F_{\mathbb{Y}, \tau}}{\longrightarrow}
      \varepsilon_{\ast} \tau_{\ast} \mathbb{Y}^{\Pi^{-1}}
      \overset{\iota(\Pi^{-1})}{\longrightarrow}
\varepsilon_{\ast} \tau_{\ast} \mathbb{Y}. 
    \end{equation}
  Our assertion says that Drinfeld's morphism
  $\mathcal{M}_{\text{Dr}} \rightarrow \wh{\Omega}_{F}$ fits into a
  commutative diagram
    \begin{displaymath}
    \xymatrix{
      \mathcal{M}_{\text{Dr}}(R) \ar[rr]^{\alpha} \ar[rd] & &
      \mathcal{M}_{\text{Dr}}(R_{[\tau]}) \ar[ld]\\ 
      & \wh{\Omega}_{F}(R) & \\
    }   
    \end{displaymath}
 This is stated as an exercise in the proof of \cite[Prop. 3.77]{RZ},  
 but we give the verification. We have to go back to Drinfeld's proof and
 therefore we use his notation. A point of $\wh{\Omega}_{F}(R)$ is given
 by data $(\eta, T, u, r)$ (\cite{Dr}, \S 2, Thm.). Drinfeld constructs
 the data $(\eta, T, u)$ entirely from a graded Cartier module
  $M = \oplus M_i$. The Cartier modules
  $M$ and $M^{\Pi^{-1}}_{[\tau]}$ are the same and the gradings are also the
  same because $M \mapsto M^{\Pi^{-1}}$ shifts the grading by $1$ and
  $M \mapsto M_{[\tau]}$ shifts the grading by $1$ in the opposite
  direction. Finally,  we have to see that the rigidification $r$ is not changed
  by the application of $\alpha_{\breve{G}}$. This can be checked on the geometric
  points of $\Spec R$.  But
  over an algebraically closed field $L$, the rigidification is obtained
  as follows. We take the morphism of rational Dieudonn\'e modules
  \begin{displaymath}
    N \rightarrow \mathbf{N}  
  \end{displaymath}
  induced by (\ref{Pol26e}) for $S = \ov S = \Spec L$. Then $r$ is obtained by
  taking the
  invariants by $V^{-1}\Pi$ on both sides. We see from the definition 
  (\ref{Pol28e}) that $\alpha_{\breve{G}}$ does not change $r$. 
  \end{proof}
Let $\Aut_D^o(\mathbb{Y})$\index[NO]{AAE@$\Aut_D^o(\mathbb{Y})$}  be the group of quasi-isogenies of $\mathbb{Y}$
which commute with the action of $\iota_{\mathbb{Y}}$. With the notation of
(\ref{decunr2e}),  let $N = N_0 \oplus N_1 = P \otimes \mathbb{Q}$ be the
rational Dieudonn\'e module of $\mathbb{Y}$. The natural map
\begin{displaymath}
 \Aut_D^o(\mathbb{Y}) \rightarrow \GL_F (N_0^{V^{-1}\Pi}) \cong \GL_2(F)
\end{displaymath}
is an isomorphism. This group acts on the functor
$\tilde{\mathcal{M}}_{\text{Dr}}$ as follows. For $g \in \Aut_D^o(\mathbb{Y})$
we define
\begin{displaymath}
g: \mathcal{M}_{\text{Dr}}(i)\rightarrow \mathcal{M}_{\text{Dr}}(i + \ord \det g),  
\quad (Y, \rho) \mapsto (Y, g\rho). 
\end{displaymath}
The action commutes with the translation. Let $J_{\text{Dr}}$
\index[NO]{JAC@$J_{\text{Dr}}$}  be the cokernel 
\begin{displaymath}
  \begin{array}{cccl}
      \mathbb{Z} & \rightarrow & \Aut_D^o(\mathbb{Y}) \times \mathbb{Z} &
      \rightarrow J_{\text{Dr}} \rightarrow 0 .\\
  i & \mapsto & (\pi^i, -2i) & \\
    \end{array}
\end{displaymath}
The second group acts on $\tilde{\mathcal{M}}_{\text{Dr}}$ such that the factor
$\mathbb{Z}$ acts by translation. We obtain an action of $J_{\text{Dr}}$ on
$\tilde{\mathcal{M}}_{\text{Dr}}$. We introduce the groups $J^{\ast {\rm r}}$
\index[NO]{JAD@$J^{\ast {\rm r}}$} and $J^{\ast{\rm ur}}$
\index[NO]{JAE@$J^{\ast{\rm ur}}$}  as cokernels 
\begin{equation}\label{Pol29J}
  \begin{array}{cccl}
      F^{\times} & \rightarrow & \Aut_D^o(\mathbb{Y}) \times K^{\times} &
      \rightarrow J^{\ast{\rm r}} \rightarrow 0,\\
      f & \mapsto & (f, f^{-1}) & \\[3mm]
     F^{\times} & \rightarrow & \Aut_D^o(\mathbb{Y}) \times (F')^{\times} &
      \rightarrow J^{\ast{\rm ur}} \rightarrow 0.\\
  f & \mapsto & (f, f^{-1}) & \\
\end{array}
  \end{equation}
The homomorphisms $\ord_K: K^{\times} \rightarrow \mathbb{Z}$, resp.
$2\ord_{F'}: (F')^{\times} \rightarrow \mathbb{Z}$, induce homomorphisms
$J^{\ast{\rm r}} \rightarrow J_{\text{Dr}}$, resp. $J^{\ast{\rm ur}} \rightarrow J_{\text{Dr}}$. 
Therefore these groups act on $\tilde{\mathcal{M}}_{\text{Dr}}$.

    \subsection{The alternative theorem in the ramified case}\label{ss:altram}

Let $K$ be a ramified quadratic extension of $F$. We also assume that
$p\neq 2$. We choose prime elements $\Pi \in O_K$ and $\pi \in O_F$ such
that $\Pi^2 = - \pi$ as in section \ref{s:almostprinc}. With the notation
before Proposition \ref{Pol1p}, we regard $O_K$ as a subring of $O_D$. 

For each $i \in \mathbb{Z}$ we define the functor $\CN(i) =\CN_{K/F,r}(i)$
on the category of schemes $S \longrightarrow \Spf O_{\breve{F}}$.
We fix a special formal $O_D$-module $\mathbb{Y}$ over 
$\bar{\kappa}_F$ and we fix a Drinfeld polarization $\psi_{\mathbb{Y}}$.
We denote by $\lambda_{\mathbb{Y}}: \mathbb{Y} \longrightarrow \mathbb{Y}^{\Delta}$
the isomorphism associated to $\psi_{\mathbb{Y}}$, cf. Definition \ref{Pol3d}. 
We will consider $p$-divisible groups $X$ on $S$ with an action
$\iota: O_K \longrightarrow \End X$ such that the restriction of this action to
$O_F$ is strict. By duality we obtain an action of $O_K$ on the Faltings dual
$X^{\nabla}$. If we compose this action with the conjugation of $K/F$ we
obtain $\iota^{\Delta}: O_K \longrightarrow \End X^{\nabla}$. We write
$X^{\Delta} = (X^{\nabla}, \iota^{\Delta})$
\index[NO]{XAB@$X^{\Delta} = (X^{\nabla}, \iota^{\Delta})$} and call this the \emph{Faltings conjugate
dual} of $(X, \iota)$. \index{Faltings conjugate
dual}

\index[NO]{NBC@$\mathcal{N}(i)$}
\begin{definition}\label{Ni1d} 
A point of $\mathcal{N}(i)(S)$ consists of the following data: 
    \begin{enumerate}
    \item A formal $p$-divisible group $X$ over $S$ with an action
      \begin{displaymath}
        \iota: O_K \longrightarrow \End X,
      \end{displaymath}
      such that the restriction of $\iota$ to $O_F$ is a strict action. 
    \item An isomorphism of $O_K$-modules $\lambda: X \longrightarrow X^{\Delta}$
      which induces a relative polarization on $X$, cf. Corollary
      \ref{LTD9c}. 
      \item A quasi-isogeny of $O_K$-modules 
      \begin{displaymath}
        \rho: X \times_{S} \ov S \longrightarrow   \mathbb{Y}
        \times_{\Spec \bar{\kappa}_F} \ov S. 
      \end{displaymath}
    \end{enumerate}

    We require that the following conditions are satisfied.
    \begin{enumerate}
    \item[a)] $\rho$  respects the $O_K$-actions.
      There is an element $u \in O_{F}^{\times}$ such that the following diagram of
      quasi-isogenies is commutative
    \begin{equation}\label{Mi4e} 
      \xymatrix{
        X \times_{S} \bar{S} \ar[r]^{\rho} \ar[d]_{u \pi^i \lambda} & 
        \mathbb{Y} \times_{\Spec \bar{\kappa}_F} \bar{S}\ar[d]^{\lambda_{\mathbb{Y}}}\\
        X^{\Delta} \times_{S} \bar{S}  & \ar[l]^{\rho^{\Delta}}
        \mathbb{Y}^{\Delta} \times_{\Spec \bar{\kappa}_F} \bar{S}.\\
      } 
    \end{equation}
    \item[b)] 
      \begin{equation}\label{TraceCond1e}
        \Trace( \iota(\Pi) \; | \; \Lie X) = 0.
      \end{equation}
    \end{enumerate}
    Two such data $(X_1, \iota_1, \lambda_1, \rho_1)$ and
    $(X_2, \iota_2, \lambda_2, \rho_2)$ define the same point of
    $\mathcal{N}(i)(S)$ iff there is an  isomorphism
    $\alpha: (X_1, \iota_1) \longrightarrow (X_2, \iota_2)$ which respects the
    polarizations up to a factor in $O_{F}^{\times}$ and such that $\alpha$
    commutes with $\rho_1$ and $\rho_2$.
\end{definition}
We note that changing $\lambda$ by a factor in $O_F^{\times}$ does not alter
the points of $\mathcal{N}(S)$. The existence of $\rho$ implies that
$\dim X = 2$ and that the $O_F$-height of $X$ is $4$. The condition b) 
implies the following Kottwitz condition for the characteristic polynomial,  
    \begin{equation}\label{kottforrel}
      {\rm char} (\iota(a) \; | \; \Lie X) = (T- a)(T - \bar{a}), \quad a \in
      O_K. 
    \end{equation}
Clearly the functor $\mathcal{N}(i)$ does not depend on the choice of
the Drinfeld polarization $\lambda_{\mathbb{Y}}$.

It follows from \cite{RZ} that $\mathcal{N}(i)$ is representable by a
formal scheme which is locally formally of finite type over $\Spf O_{\breve{F}}$. 

Let $S = \Spec R$, $R \in \Nilp_{O_{\breve{F}}}$. Let $\mathcal{P}_X$ be the
$\mathcal{W}_{O_F}(R)$-display associated to the $p$-divisible group $X$.
The conjugate dual $\mathcal{W}_{O_F}(R)$-display $\mathcal{P}_X^{\Delta}$ is
nilpotent. It corresponds to $X^{\Delta}$. We denote by
$\psi:\mathcal{P}_X \times \mathcal{P}_X \longrightarrow \mathcal{P}_{m,W_{O_F}(R)}$
the bilinear form of displays which corresponds to $\lambda$.
We may reformulate the commutativity of the diagram
(\ref{Mi4e}) as follows:  the quasi-polarization $\rho^{\ast} \psi_{\mathbb{Y}}$
coincides with $\pi^i \psi_{R/\pi R}$ of $(\mathcal{P}_X)_{R/\pi R}$ up to a
factor in $O_F^{\times}$. 

We obtain from (\ref{Mi4e}) that 
    \begin{displaymath}
      4i = 2\height_{O_F} \rho. 
    \end{displaymath}

    As for the functors $\mathcal{M}_{\text{Dr}}(i)$ we have functor isomorphisms
    \begin{equation}\label{FuncN1e}
      \Pi: \mathcal{N}(i) \isoarrow \mathcal{N}(i+1), \quad (X, \rho)
      \mapsto (X, \iota(\Pi) \rho) ,
    \end{equation}
which we call the translations. Let $\tau \in \Gal(\breve{F}/F)$ be the
Frobenius automorphism. Using the Frobenius 
$F_{\mathbb{Y}, \tau}: \mathbb{Y} \longrightarrow \tau_{\ast} \mathbb{Y}$ we define
\index[NO]{ZZWF@$\omega_{\mathcal{N}}$}
    \begin{equation}\label{WD4e}
      \omega_{\mathcal{N}}: \mathcal{N}(R)(i) \longrightarrow
      \mathcal{N}(i+1)(R_{[\tau]}).
    \end{equation}
    exactly as $\omega_{\mathcal{M}_{\text{Dr}}}$ in (\ref{WD2e}). This defines a Weil
    \index[NO]{NBA@$\tilde{\mathcal{N}}$} 
descent datum $\omega_{\mathcal{N}}$  relative to $O_{\breve{F}}/O_{F}$ on
the functor
    \begin{displaymath}
    \tilde{\mathcal{N}} = \coprod_{i \in \mathbb{Z}} \mathcal{N}(i).  
      \end{displaymath}
    \begin{lemma}\label{J^v1l} 
      The action of the group $J^{\ast{\rm r}}$ from (\ref{Pol29J}) on the $O_K$-module
      \index[NO]{JAF@$J^{\cdot}$} $\mathbb{Y}$ gives an isomorphism 
\begin{equation*}
J^{\ast{\rm r}}\isoarrow J^{\cdot} ,
\end{equation*}
where  
\begin{displaymath}
  J^{\cdot} = \{\alpha \in \Aut^o_{K}(\mathbb{Y}) \; | \;
  \psi_{\mathbb{Y}}(\alpha(x), \alpha(y)) = \mu(\alpha) \psi_{\mathbb{Y}}(x,y),
  \; \text{for some} \; \mu(\alpha) \in F^{\times}, \; 
      x,y \in P_{\mathbb{Y}} \otimes \mathbb{Q}\}.
      \end{displaymath}\qed
    \end{lemma}The group $J^{\ast{\rm r}}$ acts on the functor $\tilde{\mathcal{N}}$ by
    \begin{displaymath}
      (Y, \iota, \rho) \mapsto (Y, \iota, g \rho), \quad \text{for}\;
      g \in J^{\ast{\rm r}}. 
      \end{displaymath}
    We have a natural morphism of functors on $\Nilp_{O_{\breve{F}}}$
\begin{equation}\label{KRfunctor1e}
  \CM_{\rm Dr}(i)
  \longrightarrow \mathcal{N}(i).  
\end{equation}
This is defined as follows. Let
$(Y, \iota, \rho) \in \mathcal{M}_{\rm Dr}(i)(S)$ be a point. Let $\psi$ be a
Drinfeld polarization on $Y$ which is compatible with the quasi-isogeny $\rho$, 
cf. Proposition \ref{Pol2p}. It is uniquely determined up to a factor in
$O_F^{\times}$. Locally
on $\bar{S}$ we have $\rho^{\ast} \psi_{\mathbb{Y}} = f \psi$ for $f \in F$.
Since $\height_{O_F} \rho = 2i$ we obtain $\ord_{\pi} f = i$. 
Therefore $(Y, \iota_{| O_K}, \psi, \rho) \in \mathcal{N}(i)(S)$.

The main result of \cite{KRalt} may now be formulated as follows. Note that in loc.~cit. Weil descent data were not considered. 
\begin{theorem}[\cite{KRalt}]\label{altDrinf}
  Assume that $p\neq 2$. 
    The functor morphisms \eqref{KRfunctor1e}  define  a functor isomorphism 
    \begin{displaymath}
         \tilde{\mathcal{M}}_{\rm Dr}
  \isoarrow \tilde{\mathcal{N}}
      \end{displaymath}
which commutes with the Weil  descent data and the action of the group
$J^{\cdot} = J^{\ast{\rm r}}$ on both sides. In particular it commutes with the
translations.   
\end{theorem}
It is clear that the morphism of functors 
$
  \tilde{\mathcal{M}}_{\rm Dr}
  \longrightarrow \tilde{\mathcal{N}},  
$
is compatible with the  Weil descent data
$\omega_{\mathcal{M}_{\text{Dr}}}$ and $\omega_{\mathcal{N}}$ relative to  
$O_{\breve{F}}/O_{F}$ and  with the
translations (\ref{FuncMDr1e}) and (\ref{FuncN1e}). From this we see that
it commutes also with the actions of $J^{\ast{\rm r}}$. 
We need to prove that it is an isomorphism. For the proof we need some preparations. 

Let $k$ be an algebraically closed field which is an
$O_{\breve{F}}$-algebra.  We consider a $W_{O_F}(k)$-Dieudonn\'e
module $M$ of height $4$ and dimension $2$. We assume that an
$O_K$-action $ \iota: O_K \longrightarrow \End M$ on $M$ is given such
that the restriction to $O_F$ is via $O_F \longrightarrow W_{O_F}(k)$.

    Let
    \begin{equation}
      \psi : M \times M \longrightarrow W_{O_F}(k)
    \end{equation}
    be a relative polarization, i.e., an alternating
    $W_{O_F}(k)$-bilinear form such that
    \begin{displaymath}
      \psi(Fx_1, Fx_2) = \pi ~^F\psi (x_1, x_2).
    \end{displaymath} 
    We require that
    \begin{displaymath}
      \psi(\iota(a)x, y) = \psi (x, \iota(\bar{a})y), \quad a \in O_K.
    \end{displaymath}

    \begin{proposition}\label{Pol4p}
      Let $\mathbb{M}$ be the $W_{O_F}(k)$-Dieudonn\'e module of a
      special formal $O_D$-module with a Drinfeld polarization
      $\psi_{\mathbb{M}}$, cf. Definition  \ref{Pol3d}.
Let $(M, \iota, \psi)$ be as above and such that $\psi$ is perfect. We
assume that there exists an isomorphism of rational
$W_{O_F}(k)$-Dieudonn\'e modules
$\rho: \mathbb{M} \otimes \mathbb{Q} \longrightarrow  M \otimes
\mathbb{Q}$, such that $\rho$ is a 
homomorphism of $O_K$-modules and such that $\rho$ respects the
polarizations $\psi_{\mathbb{M}}$ and $\psi$ up to a factor in $F^{\times}$. 

Then there exists a unique $O_D$-module structure on $M$ such that $M$  
becomes the Dieudonn\'e module of a special formal $O_D$-module and such
that $\rho$ is a quasi-isogeny of $O_D$-modules.  
    \end{proposition}
    \begin{proof}
We will write $ax := \iota(a)x$, for $a \in O_K$ and $x \in
M$. For $a \in O_F$ this coincides by definition with the action
via $O_F \longrightarrow W_{O_F}(k)$.

We define $\tilde{W} = O_K \otimes_{O_F} W_{O_F}(k)$. We extend
the conjugation of $K$ over $F$ by linearity to $\tilde{W}$ over
$W_{O_F}(k)$. We denote the traces of $K/F$ and of
$\tilde{W}/W_{O_F}(k)$ both by $\Trace$. The Frobenius
endomorphism of $W_{O_F}(k)$ extends $O_K$-linearly to
$\tilde{W}$ and is denoted by $F$. It will be impossible to
confuse this with the field $F$.

We define a hermitian form
      \begin{displaymath}
        h: M \times M \longrightarrow \tilde{W},
      \end{displaymath} 
      by requiring that
      \begin{displaymath}
        \Trace \xi{\Pi}^{-1} h(x, y) =  \psi( x, \xi y) 
        \quad \xi \in \tilde{W}, \; x,y \in M.
      \end{displaymath}
      Then $h$ is $\tilde{W}$-linear in the second variable and hermitian,
      \begin{displaymath}
        h(x, y) = \overline{h(y, x)}.
      \end{displaymath} 
      The pairing $h$ is perfect and  satisfies the equation
      \begin{displaymath}
        h(Fx, Fy) = \pi ~^{F} h(x, y).
      \end{displaymath}
      Since $N := M \otimes \mathbb{Q}$ is the rational Dieudonn\'e module
      of a special formal $O_D$-module with its Drinfeld polarization, we have a decomposition 
      \begin{equation}\label{Pol14e}
        N = N_0 \oplus N_1, 
      \end{equation}
      which is orthogonal with respect to $\psi$ (see
      (\ref{decunr})). One should note that $N_0$ and $N_1$ are not
      $\tilde{W}$-modules. 

      We note that for $n_0, n'_0 \in N_0$, $n_1, n'_1 \in N_1$ we have

      \begin{equation}\label{Pol13e}
        \begin{array}{c}
          h(n_0, n'_0) = \frac{1}{2} \Pi \psi(n_0, n'_0), \quad    h(n_1, n'_1)
          = \frac{1}{2} \Pi \psi(n_1, n'_1),\\[2mm]  
          h(n_0, n_1) =  \frac{1}{2} \psi(n_0, \Pi n_1). 
        \end{array}
      \end{equation}
      Indeed, the equation 
      \begin{displaymath}
        \Trace \frac{\Pi}{\Pi} h(n_0,  n'_0) = \psi( n_0, \Pi n'_0) = 0   
      \end{displaymath}
      implies that $\Pi^{-1}h(n_0, n'_0) \in W_{O_F}(k) \otimes
      \mathbb{Q}$. We obtain the first equation of (\ref{Pol13e}):
      \begin{displaymath}
        {2}{\Pi}^{-1} h(n_0, n'_0) = \Trace {\Pi}^{-1} h(n_0, n'_0) =
        \psi (n_0, n'_0). 
      \end{displaymath} 
      The proof of the next equation is the same. We have
      \begin{displaymath}
        \Trace {\Pi}^{-1}h(n_0, n_1) = \psi(n_0, n_1) = 0.
      \end{displaymath}
      This implies $h(n_0, n_1) \in W_{O_F}(k) \otimes \mathbb{Q}$. We
      obtain the last equation of (\ref{Pol13e}):
      \begin{displaymath}
        2h(n_0, n_1) = \Trace h(n_0, n_1) = \psi(n_0, \Pi n_1).
      \end{displaymath}
      In particular we see from (\ref{Pol13e}) that an element $n_0 \in N_0$
      is isotropic for $h$. 

We call an element $x \in M$ primitive if it is not in $\Pi M$.
We find an element $x \in M \cap N_0$ such that $x \notin \pi M$. Assume that
$x = \Pi y$ for some $y \in M$. Then $y \in M \cap N_1$. Then it is clear that
$y$ is a primitive element in $M$. Interchanging the role of the indices $0$
and $1$, we may assume that $x \in M \cap N_0$ is primitive. 

Since the pairing $h$ is perfect and $x$ is isotropic
for $h$, we find an element $y' \in M$, such that $h(x, y') = 1$. We can even
choose $y'$ to be isotropic for $h$. Indeed, we set $y = y' + \lambda x$ for
some $\lambda \in \tilde{W}$. Then $h(x, y) = 1$. We compute:
      \begin{displaymath}
        h(y, y) = h(y', y') + h(y', \lambda x) + h(\lambda x, y') = h(y',
        y') + \lambda + \bar{\lambda}. 
      \end{displaymath}
      We choose $\lambda = -(1/2)h(y', y')$ (which is legitimate, as $p\neq 2$)
      and obtain $h(y, y) = 0$. According to (\ref{Pol14e}) we write
      \begin{displaymath}
        y = y_0 + y_1, \quad y_0 \in N_0, \; y_1 \in N_1. 
      \end{displaymath}
      We write
      \begin{displaymath}
        1 = h(x, y) = h(x, y_0) + h(x, y_1). 
      \end{displaymath}
      We have $h(x,y_0) \in \Pi W_{O_F}(k)
      \otimes \mathbb{Q}$ and $h(x, y_1) \in W_{O_F}(k) \otimes
      \mathbb{Q}$ by the formulas (\ref{Pol13e}). 
      This implies $h(x, y_0) = 0$ and $h(x, y_1) = 1$. 
      On the other hand, we find by (\ref{Pol14e}) 
      \begin{displaymath}
        0 = h(y, y) = h(y_0 +y_1, y_0 + y_1) = h(y_0, y_1) + h(y_1, y_0).
      \end{displaymath}
      Since $h(y_0, y_1) \in W_{O_F}(k) \otimes \mathbb{Q}$, this implies
      $h(y_0, y_1) = 0$ and then $h(y, y_0) = 0$. The elements $x$ and $y$
      generate $N$ as a $\tilde{W} \otimes \mathbb{Q}$-vector space. Because
      we already proved that $h(x, y_0) = 0$, we conclude $y_0$ =
      0. Therefore $y = y_1 \in M \cap N_1$. We obtain
      \begin{displaymath}
        M = \tilde{W}x + \tilde{W}y,
      \end{displaymath}
      because $h$ is unimodular on the right hand side. Then the elements $x,
      \Pi x, y, \Pi y$  are a basis of the $W_{O_F}(k)$-module $M$. We have
      $x, \Pi y \in M \cap N_0$ and $y, \Pi x \in M \cap N_1$ and therefore 
      \begin{displaymath}
        M = (M \cap N_0) \oplus (M \cap N_1).
      \end{displaymath}
      This shows that the $O_D$-module structure on $N$ induces an
      $O_D$-module structure on $M$. 
    \end{proof}
    \begin{proof}[Proof of Theorem \ref{altDrinf}]  
      We consider the morphism (\ref{KRfunctor1e}) for $i = 0$ and
      denote it by
      \begin{equation}\label{KRfunctor3e}
  \CM_{\rm Dr}  \longrightarrow \mathcal{N}. 
        \end{equation}
Clearly it is enough to show that this is an isomorphism. Proposition
\ref{Pol4p} shows that for any algebraically closed field $k$  which is an
$O_{\breve{F}}$-algebra,  the induced map $\CM_{\rm Dr}(k) \longrightarrow \mathcal{N}(k)$
is bijective.

We note that the morphism (\ref{KRfunctor3e}) is formally
unramified. Indeed, let $S \longrightarrow R$ be a surjective morphism in
$\Nilp_{O_{\breve{F}}}$ with nilpotent kernel. Let $X$ be a
$p$-divisible group over $S$ with base change $X_{R}$ over $R$. Then
an $O_D$-module structure on $X_R$ lifts by rigidity in at most one
way to $X$. We consider  the underlying topological spaces in
(\ref{KRfunctor1e}) with their induced structure of reduced
schemes. Then we obtain a formally unramified morphism of
$\bar{\kappa}_F$-schemes 
\begin{equation}\label{KR-Bw1e}
  \CM_{{\rm Dr},{\rm red}} \longrightarrow \mathcal{N}_{\rm red}
\end{equation}
These schemes are locally of finite type over $\bar{\kappa}_F$ and have
irreducible components which are proper over $\bar{\kappa}_F$, cf. \cite[Prop. 2.32]{RZ}. Moreover, the morphism is bijective on geometric points. Then
the irreducible components of both schemes correspond  bijectively to
each other. We consider a point $x \in \CM_{\rm Dr}(\bar{\kappa}_F)$. Let $X$
be the union of all irreducible components which pass through $x$ with the
reduced scheme structure. Let $y \in \mathcal{N}(\bar{\kappa}_F)$ be the image
of $x$ and define $Y \subset \mathcal{N}$ in the same way as $X$. Then
$X \longrightarrow Y$ is a finite morphism. If we remove all points in $X$
resp. $Y$ which belong to components not passing through $x$, resp. $y$,
we obtain a finite morphism of open neighbourhoods $U \longrightarrow V$
of $x \in \CM_{{\rm Dr},{\rm red}}$ and $y \in \mathcal{N}_{r\rm ed}$. Therefore
(\ref{KR-Bw1e}) is a finite morphism of
schemes locally of finite type over the algebraically closed field
$\bar{\kappa}_F$. Since this morphism is unramified and bijective on geometric points, it is an isomorphism. 

      \begin{lemma}\label{KR-Bw1l}
        Let $S \longrightarrow \bar{\kappa}_F$ be a surjective morphism in
        $\Nilp_{O_{\breve{F}}}$ such that the kernel is nilpotent and endowed
        with divided powers. Then the map
        \begin{displaymath}
          \CM_{\rm Dr}(S) \longrightarrow \mathcal{N}(S) 
        \end{displaymath}
        is bijective. 
      \end{lemma}

Let us assume that the lemma is proved. Then we consider points $x$
and $y$ as above. We consider an open affine neighbourhood $U$ of
$x$. By the isomorphism (\ref{KR-Bw1e}) we regard $U$ also as a
neighbourhood of $y$. Let $n \in \mathbb{N}$.  For a suitable ideal
sheaf of definition $\mathcal{J}$ of $\CN$, we have a homomorphism
      \begin{displaymath}
        (\mathcal{O}_{\mathcal{N}}/\mathcal{J})(U) \longrightarrow
        \mathcal{O}_{\CM_{\rm Dr}}(U)/\pi^{n} \mathcal{O}_{\CM_{\rm
            Dr}}(U).
      \end{displaymath}
This map is surjective modulo $\pi$ by (\ref{KR-Bw1e}) and is
therefore surjective. We note that by EGA$0_{I},$ Prop.~7.2.4 the ring
$(\mathcal{O}_{\mathcal{N}}/\mathcal{J})(U)$ is $\pi$-adic. It follows
that
      \begin{displaymath}
        \mathcal{O}_{\mathcal{N}} (U) \longrightarrow
        \mathcal{O}_{\CM_{\rm Dr}}(U)
      \end{displaymath}
is surjective. Taking the inductive limit over $U$ we obtain an
epimorphism of local rings
      \begin{equation}\label{KR-Bw2e}
        \mathcal{O}_y \longrightarrow \mathcal{O}_{x}.
      \end{equation}
By \cite[Thm.~2.16]{RZ} this is a homomorphism of noetherian
adic rings (comp. EGA I, Prop.~10.1.6). The ring $\mathcal{O}_x$
is, as a local ring of the scheme $\CM_{\rm Dr}$, regular of
dimension $2$. Let $\mathfrak{m}_y$ and $\mathfrak{m}_x$ be the
maximal ideals of the local rings. We remark that the squares of
the ideals are open because the topologies are adic.

We apply  Lemma \ref{KR-Bw1l} to $S = \mathcal{O}_y/\mathfrak{m}_y^2$. Then
we obtain an oblique arrow which makes the following diagram
commutative, 

      \begin{displaymath}
        \xymatrix{
          \mathcal{O}_y  \ar[r] \ar[d] & \mathcal{O}_x \ar[ld]\\
          \mathcal{O}_y/\mathfrak{m}_y^2 \, .
        } 
      \end{displaymath}
      It follows that there is a surjective homomorphism
      $\mathfrak{m}_x/\mathfrak{m}_x^2 \longrightarrow
      \mathfrak{m}_y/\mathfrak{m}_y^2$.
      The epimorphism of local rings also gives  a surjection in the other
      direction. We conclude
      \begin{displaymath}
        \dim_{\bar{\kappa}_F} \mathfrak{m}_y/\mathfrak{m}_y^2 = 2.
      \end{displaymath}
      Therefore $\mathcal{O}_y$ is a regular local ring of dimension
      $2$, and the map (\ref{KR-Bw2e}) is an isomorphism. It follows
      that the map of sheaves
      \begin{displaymath}
        \mathcal{O}_{\mathcal{N}}  \longrightarrow
        \mathcal{O}_{\CM_{\rm Dr}}
      \end{displaymath}
      is an isomorphism. 
      Finally let $J$ be the maximal ideal sheaf of definition of
      $\mathcal{O}_{\mathcal{N}}$. By the isomorphism (\ref{KR-Bw2e}) we obtain an
      isomorphism 
      \begin{displaymath}
        \mathcal{O}_y/J \mathcal{O}_y \longrightarrow \mathcal{O}_x/\pi
        \mathcal{O}_x.  
      \end{displaymath}
      Therefore $J \mathcal{O}_y = \pi \mathcal{O}_y$.  Therefore
      $J = \pi \mathcal{O}_{\mathcal{N}}$ is an ideal sheaf of
      definition. We obtain that (\ref{KRfunctor3e}) is an isomorphism of
      formal schemes. 

      It remains to prove  Lemma \ref{KR-Bw1l}. We denote  by $\mathfrak{m}$  the kernel of
      $S \longrightarrow \bar{\kappa}_F$. Let
      $\xi: \Spec S \longrightarrow \mathcal{N}$ be a morphism. We show that it lifts
      uniquely to $\Spec S \longrightarrow \CM_{\rm Dr}$. We denote by
      $y \in \mathcal{N}(\bar{\kappa}_F)$ the point induced by $\xi$. Let
      $x \in \CM_{\rm Dr}(\bar{\kappa}_F)$ be the unique point over
      $y$.       

      We denote by $\mathcal{P}$ the $O_F$-display of the special formal
      $O_D$-module over $\bar{\kappa}_F$ which corresponds to $x$.  We
      denote by $\tilde{\mathcal{P}}$ the unique
      $\mathcal{W}_{O_F}(S/\bar{\kappa}_F)$-display which lifts
      $\mathcal{P}$, cf. Theorem \ref{DspKristall1t}.  We write
      $\tilde{\CP} = (\tilde{P}, \hat{Q}, F, \dot{F})$. The $O_D$-action on
      $\mathcal{P}$ extends to $\tilde{\mathcal{P}}$. Therefore we have the
      decompositions 
      \begin{displaymath}
        \tilde{P} = \tilde{P}_0 \oplus \tilde{P}_1, \quad \hat{Q} =
        \hat{Q}_0 \oplus \hat{Q}_1. 
      \end{displaymath}
      We consider only the most interesting case where $\Pi$ acts trivially
      on $\Lie \mathcal{P}$, i.e.,
      $\Spec \bar{\kappa}_F \longrightarrow \CM_{\rm Dr}$ is a singular
      point of the special fibre, cf. \cite{Dr}. In this case we obtain Frobenius-linear isomorphisms
      \begin{displaymath}
        \dot F \circ \Pi: \tilde{P}_0 \longrightarrow \hat{Q}_1 \longrightarrow
        \tilde{P}_0, \quad   \dot F \circ \Pi: \tilde{P}_1 \longrightarrow
        \hat{Q}_0 \longrightarrow \tilde{P}_1. 
      \end{displaymath}
      We set $\tilde{U}_i = \{ x \in \tilde{\mathcal{P}}_i \; | \; \dot F
      \circ \Pi (x) = x \}$. Then the canonical morphism $W_{O_F}(S)
      \otimes_{O_F} \tilde{U}_i \longrightarrow \tilde P_i$ is an isomorphism. 

      We can make the same construction with the display $\mathcal{P}$. Then
      we obtain $U_i \subset P_i$ such that the canonical $O_F$-module
      homomorphism $\tilde{U}_i \longrightarrow U_i$ is an
      isomorphism. Using our knowledge about $\mathcal{P}$ we find elements
      $\tilde{e}_i \in U_i$, for $i = 0, 1$ such that 
      \begin{displaymath}
        \tilde{e}_0, \Pi \tilde{e}_1 \in \tilde{P}_0, \quad   \tilde{e}_1,
        \Pi \tilde{e}_0 \in \tilde{P}_1, 
      \end{displaymath}
      are a basis of the $W_{O_F}(S)$-module $\tilde{P}$. The natural
      polarization $\psi$ on $\mathcal{P}$ extends to a polarization
      $\tilde{\psi}$ on $\tilde{\mathcal{P}}$ which is given by the conditions 
      \begin{displaymath}
        \tilde{\psi}(\tilde{e}_0, \Pi \tilde{e}_1) = 1 =
        \tilde{\psi}(\tilde{e}_1, \Pi \tilde{e}_0),  
      \end{displaymath}  
      and such that the decomposition
      $\tilde{P} = \tilde{P}_0 \oplus \tilde{P}_1$ is orthogonal with
      respect to $\tilde{\psi}$.   

      We classify now the liftings of
      $\Spec \bar{\kappa}_F \longrightarrow \mathcal{N}$ to a
      point $\Spec S \longrightarrow \mathcal{N}$. We consider the Hodge
      filtration $L = Q/I_{O_F}(k)P \subset P/I_{O_F}(k)P$. Since we 
      compute now all the time modulo the augmentation ideal
      $I_{O_F}(k) \subset W_{O_F}(k)$, resp.,
      $I_{O_F}(S) \subset W_{O_F}(S)$, we continue to simply  write $\tilde{e}_0$ when we
      mean the residue class in $\tilde{P}/I(S) \tilde{P}$. The $k$-vector
      space $L$ has the basis $\Pi e_0, \Pi e_1$. Therefore a lifting of $L$
      to a direct summand $\tilde{L} \subset \tilde{P}/I(S)\tilde{P}$ has a
      unique basis of the form 
      \begin{displaymath}
        f_0 = \Pi \tilde{e}_1 + \gamma \tilde{e}_0 + \delta \tilde{e}_1,
        \quad f_1 = \Pi \tilde{e}_0 + \alpha \tilde{e}_0 + \beta \tilde{e}_1,
      \end{displaymath} 
      because it is complementary to the module generated by $\tilde{e}_0,
      \tilde{e}_1$. Since we want a lifting of $L$ we have $\alpha, \beta,
      \gamma, \delta \in \mathfrak{m}$. 
      The lifting $\tilde{L}$ determines a lifting of the
      display to $S$. The form $\psi$ lifts to a polarization of this
      display if and only if $\tilde{L}$ is isotropic under
      $\psi$. Therefore we must have
      \begin{displaymath}
        0 = \tilde{\psi}(\Pi \tilde{e}_0 + \alpha \tilde{e}_0 + \beta
        \tilde{e}_1, \Pi \tilde{e}_1 + \gamma \tilde{e}_0 + \delta \tilde{e}_1)
      \end{displaymath}
      One obtains easily that the right hand side is 
      \begin{displaymath}
        \tilde{\psi}(\Pi \tilde{e}_0, \delta \tilde{e}_1) +
        \tilde{\psi}(\alpha \tilde{e}_0, \Pi \tilde{e}_1) = - \delta + \alpha. 
      \end{displaymath}
      Since the lifting $\tilde{L}$ should define a point of $\mathcal{N}$, 
      the condition 2) in the definition of points of $\mathcal{N}$ implies
      \begin{displaymath}
        0 = \Trace(\Pi \; | \; \tilde{P}/\tilde{L}) = \alpha + \delta. 
      \end{displaymath}
      Because $p \neq 2$ we obtain $\alpha = \delta = 0$. This implies that
      $\tilde{L} = (\tilde{L} \cap \tilde{P}_0) \oplus (\tilde{L} \cap
      \tilde{P}_1)$. 
      This shows that the display over $S$ defined by $\tilde{L}$ is the
      display of a special formal $O_D$-module. Therefore the liftings of
      $\Spec \bar{\kappa}_F \longrightarrow \mathcal{N}$ to $\mathcal{N}(S)$
      correspond via (\ref{KRfunctor1e}) bijectively to the liftings of
      $\Spec \bar{\kappa}_F \longrightarrow \CM_{\rm Dr}$ to a point of
      $\mathcal{M}(S)$. This proves  Lemma \ref{KR-Bw1l} and Theorem \ref{altDrinf}.

    \end{proof}

    The properties of Drinfeld's moduli scheme $\CM_{\rm Dr}$ imply the following corollary, cf., e.g.,  \cite{BC}. 
    \begin{corollary}\label{spfibred}
      The formal scheme $\CN$ is $\pi$-adic and has semi-stable reduction. The special fiber $\CN\otimes_{O_{\breve F}}\ov\kappa_{F}$ of $\CN$ is a reduced scheme. \qed
    \end{corollary}\label{cor:structureN}

Finally we prove the uniqueness of the framing object, cf.  (i) of subsection \ref{ss:uniquefram}.  We begin
with this question in the category $\mathfrak{d}{\mathfrak R}^{\rm pol}_{R}$,
cf. Definition \ref{P'Pol3d}. 
\begin{proposition}\label{ramFrame1p}
  Let $r$ be special  and let  $K/F$ be ramified. Let $k \in \Nilp_{O_F}$ be an
  algebraically closed field. Let $(\mathcal{P}_{{\rm{c}},1}, \iota_{{\rm{c}},1}, \beta_{{\rm{c}},1})$
  and $(\mathcal{P}_{{\rm{c}},2}, \iota_{{\rm{c}},2}, \beta_{{\rm{c}},2})$ be two objects in 
  $\mathfrak{d}{\mathfrak R}^{\rm pol}_{k}$. Assume  that
  $\inv (\mathcal{P}_{{\rm{c}},i}, \iota_{{\rm{c}},i}, \beta_{{\rm{c}},i})) = -1$ for $i = 1,2$. Then there exists a quasi-isogeny
  $\alpha: \mathcal{P}_{{\rm{c}},1} \longrightarrow \mathcal{P}_{{\rm{c}},2}$ which respects
  $\iota_{{\rm{c}},i}$ and $\beta_{{\rm{c}},i}$. 

 If the forms $\beta_{{\rm{c}},i}$ are perfect,  then the
  actions $\iota_{{\rm{c}}, i}$ extend to actions
  $\tilde{\iota}_{{\rm{c}},i}: O_D \longrightarrow \End_{O_F} \mathcal{P}_{{\rm{c}},i}$ such that
$\mathcal{P}_{{\rm{c}},i}$ becomes a special formal $O_D$-module with Drinfeld
polarization $\beta_i$ and such that $\alpha$ becomes a homomorphism
of $O_D$-modules. 
    \end{proposition}
For the proof we need some preparations.
\begin{lemma}\label{ramFrame1l}
  Let $K/F, r, k$ as in the last Proposition and let
  $(\mathcal{P}_{\rm{c}}, \iota_{\rm{c}}, \beta_{\rm{c}}) \in \mathfrak{d}{\mathfrak R}^{\rm pol}_{k}$.
  Assume that $\inv (\mathcal{P}_{\rm{c}}, \iota_{\rm{c}}, \beta_{\rm{c}}) = -1$, cf. Definition
  \ref{invcontr1d}.  Then the $\mathcal{W}_{O_F}(k)$-display $\mathcal{P}_{\rm{c}}$ is isoclinic of
  slope $1/2$. 
\end{lemma}
\begin{proof} The $K \otimes_{O_F} W_{O_F}(k)$-vector
  space
$N = P_{\rm{c}} \otimes \mathbb{Q}$ has  dimension $2$. The isoclinic decomposition of the
  $\mathcal{W}_{O_F}(k)$-isocrystal $N$ is invariant under the action
  of $O_K$ and has therefore at most two summands. We have to show that
  there is only one summand. If not, we have $N = N_0 \oplus N_1$, where
  $N_0$ is \'etale and $N_1$ is dual to $N_0$. The dimension of each
  $N_i$ as a $K \otimes_{O_F} W_{O_F}(k)$-vector space is one. Therefore we find
  a generator $e_0 \in N_0$ such the $V_{\rm{c}} e_0 = e_0$. We use the notation
  of before Definition \ref{invcontr1d}. Let $e_1 \in N_1$ be the
  generator such that $\varkappa_{\rm{c}} (e_0, e_1) = 1$. Let $\tau$ be the
  Frobenius acting via the second factor on $K \otimes_{O_F} W_{O_F}(k)$. 
  From the equation
  \begin{displaymath}
    \varkappa_{\rm{c}}(V_{\rm{c}} e_0, V_{\rm{c}} e_1) = \pi \varkappa_{\rm{c}}(e_0, e_1)^{\tau^{-1}} = \pi,
  \end{displaymath}
  we conclude that $V_{\rm{c}} e_1 = \pi e_1$. Therefore
  $V_{\rm{c}}(e_0 \wedge e_1) = \pi e_0 \wedge e_1$. This implies that the
  invariant of  $(\mathcal{P}_{\rm{c}}, \iota_{\rm{c}}, \beta_{\rm{c}})$ is $1$, which  contradicts
  the assumption $(\mathcal{P}_{\rm{c}}, \iota_{\rm{c}}, \beta_{\rm{c}})=-1$.
  \end{proof}

Let $(\mathcal{P}_{\rm{c}}, \iota_{\rm{c}}) \in \mathfrak{d}{\mathfrak R}_{k}$ be isoclinic
of slope $1/2$. Then there is an $W_{O_F}(k)$-lattice 
$\Lambda \subset P_{\rm{c}} \otimes \mathbb{Q}$ which is invariant by $\pi^{-1}V_{\rm{c}}^2$.
Then there is also a lattice invariant by the ''square root''
$\Pi^{-1}V_{\rm{c}}$. One deduces that the invariants $C$ of $\Pi^{-1}V_{\rm{c}}$ acting on 
$P_{\rm{c}} \otimes \mathbb{Q}$ form  a $K$-vector space of dimension $2$ and
\begin{equation}\label{Kneun23e}
  P_{\rm{c}} \otimes \mathbb{Q} = C \otimes_{O_F} W_{O_F}(k). 
\end{equation}
The anti-hermitian form $\varkappa_{\rm{c}}$ associated to $\beta_{\rm{c}}$ by
\eqref{betacontr} induces the anti-hermitian form on the $K$-vector space $C$  
\begin{equation}\label{Kneun22e}
 \varkappa_{\rm{c}}: C \times C \longrightarrow K.
\end{equation}
Indeed, for $x, y \in C$ we find
\begin{displaymath}
  \pi \varkappa_{\rm{c}}(x, y) = \varkappa_{\rm{c}}(\Pi x, \Pi y) = \varkappa_{\rm{c}}(V_{\rm{c}} x, V_{\rm{c}} y)
  = \pi ~^{F^{-1}} \varkappa_{\rm{c}}(x,y).
\end{displaymath}
This shows that $\varkappa_{\rm{c}}(x, y) \in K \otimes_{O_F} W_{O_F}(k)$ is invariant
by the Frobenius $F$ acting on $W_{O_F}(k)$, and therefore this element is in
$K$. The same argument shows that $\beta_{\rm{c}}(x,y) \in F$. 
The form $\varkappa_{\rm{c}}$ restricted to $C$ is obtained from the restriction
of $\beta_{\rm{c}}$ to $C$ by the formula
\begin{displaymath}
  \Trace_{K/F} (a \varkappa_{\rm{c}}(x, y)) = \beta_{\rm{c}}(a x, y), \quad x, y \in C, \;
  a \in K. 
\end{displaymath}
\begin{lemma}\label{ramFrame2l}
  Let $(\mathcal{P}_{{\rm{c}},1}, \iota_{{\rm{c}},1}, \beta_{{\rm{c}},1})$ and
  $(\mathcal{P}_{{\rm{c}},2}, \iota_{{\rm{c}},2}, \beta_{{\rm{c}},2})$ be objects of
  $\mathfrak{d}{\mathfrak R}^{\rm pol}_{k}$ such that $\mathcal{P}_{{\rm{c}},1}$
  and $\mathcal{P}_{{\rm{c}},2}$ are isoclinic of slope $1/2$. Then the canonical map
  \begin{displaymath}
    \Hom \big((\mathcal{P}_{{\rm{c}},1}, \iota_{{\rm{c}},1}, \beta_{{\rm{c}},1}),
    (\mathcal{P}_{{\rm{c}},2}, \iota_{{\rm{c}},2}, \beta_{{\rm{c}},2})\big) \otimes \mathbb{Q} \longrightarrow
    \Hom_K\big((C_1, \beta_{{\rm{c}},1}), (C_2, \beta_{{\rm{c}},2})\big)
  \end{displaymath}
  is an isomorphism. 
\end{lemma}
\begin{proof}
  This is an immediate consequence of the isomorphism (\ref{Kneun23e}) because
  the $K$-action, $\beta_{\rm{c}}$, and $V_{\rm{c}}$ on $P_{\rm{c}} \otimes \mathbb{Q}$ can be
  recovered from the right hand side of the isomorphism. The map $V_{\rm{c}}$ is
  induced from $\Pi \otimes F^{-1}$ on the right hand side. 
  \end{proof}
\begin{lemma}\label{Kneun2l}
There is the following relation between the invariants defined in Definition \ref{invcontr1d} and in Definition  \ref{invdet}, 
  \begin{equation}\label{Kneun20e}
    \inv(\mathcal{P}_{\rm{c}}, \iota, \beta_{\rm{c}}) = - \inv (C, \beta_{\rm{c}}). 
  \end{equation}
\end{lemma}
We remark here that $(C, \beta_{\rm{c}})$ determines
$(\mathcal{P}_{\rm{c}}, \iota, \beta_{\rm{c}})$ up to isogeny.
\begin{proof}
  Let $x_1, x_2$ be a basis of the $K$-vector space $C$. Then the right hand
  side  of (\ref{Kneun20e}) is given by the $2 \times 2$-determinant 
  \begin{displaymath}
\det(\varkappa_{\rm{c}}(x_i,x_j)). 
  \end{displaymath}
  By definition of $C$ we have $V_{\rm{c}} x_i = \Pi x_i$. We conclude that
  $V_{\rm{c}}(x_1 \wedge x_2) = -\pi (x_1 \wedge x_2)$ in $\wedge^{2}_K C$.
  From Lemma \ref{Kneuninv1l} we obtain that the determinant above gives
  $-\inv (\mathcal{P}_{\rm{c}}, \iota, \beta_{\rm{c}})$. 
\end{proof}

\begin{lemma}\label{ramFrame4l}
  Let $(\mathcal{P}_{\rm sp}, \iota_{\rm sp})$ the $\mathcal{W}_{O_F}(k)$-display of
  a special formal $O_D$-module. We denote by $\psi$ a Drinfeld polarisation.
  Let $\iota'_{\rm sp}$ be the restriction of $\iota_{\rm sp}$ to $O_K \subset O_D$.
  Then
  $(\mathcal{P}_{\rm sp},\iota'_{\rm sp},\psi)\in \mathfrak{d}{\mathfrak R}^{\rm pol}_{k}$, and 
  \begin{displaymath}
\inv (\mathcal{P}_{\rm sp},\iota'_{\rm sp},\psi) = -1. 
    \end{displaymath}
\end{lemma}
\begin{proof}
  We write $M = P_{\rm sp}$ and consider it as a $\mathcal{W}_{O_F}(k)$-Dieudonn\'e
  module. Let $N = M \otimes \mathbb{Q}$. Then $\psi$ is a relative
  polarization that satisfies (\ref{Pol20e}). By the decomposition
  (\ref{decunr}) (or (\ref{decunr2e})) we obtain a decomposition
  \begin{displaymath}
N = N_0 \oplus N_1 ,
  \end{displaymath}
  which is orthogonal with respect to $\psi$. As in the proof of 
  Lemma \ref{ramFrame2l}, we consider the invariants $C_{\rm sp} = N^{V^{-1}\Pi}$. Because $V^{-1}\Pi$ is
  homogenous of degree zero, the decomposition of $N$ induces
  $C_{\rm sp} = C_0 \oplus C_1$. Each $C_i$ is a $F$-vector space of dimension $2$.
  The restriction of $\psi$ is a nondegenerate alternating pairing
  \begin{displaymath}
\psi: C_{\rm sp} \times C_{\rm sp} \longrightarrow F
    \end{displaymath}
  Let $\varkappa: C_{\rm sp} \times C_{\rm sp} \longrightarrow K$ be the anti-hermitian form
  associated to $\psi$ as before  Lemma \ref{ramFrame2l}. 
  We choose a basis $e_0, f_0$ of the $F$-vector space $C_0$ such that
  $\psi(e_0, f_0) = 2$.  We claim that
  \begin{equation}\label{Kneun21e} 
\varkappa(e_0, e_0) =  \varkappa(f_0, f_0) = 0, \quad \varkappa(e_0, f_0) = 1 .
    \end{equation}
  Indeed, we write $\varkappa(e_0, e_0) = a + \Pi b$, $a,b \in F$.
  By definition of $\varkappa$ we find 
  \begin{displaymath}
    \Trace_{K/F} (\varkappa(e_0, e_0)) =  \psi(e_0, e_0) = 0, \quad
    \Trace_{K/F} (\Pi\varkappa(e_0, e_0)) =  \psi(\Pi e_0, e_0) = 0. 
  \end{displaymath}
  The last equation follows because $C_0$ and $C_1$ are orthogonal. This
  implies $a = b = 0$. Clearly it is enough to verify the last equation
  of (\ref{Kneun21e}). Again we write $\varkappa(e_0, f_0) = a + \Pi b$,
  $a,b \in F$. Then we find
    \begin{displaymath}
    \Trace_{K/F} (\varkappa(e_0, f_0)) =  \psi(e_0, f_0) = 2, \quad
    \Trace_{K/F} (\Pi\varkappa(e_0, f_0)) =  \psi(\Pi e_0, f_0) = 0,  
    \end{displaymath}
    and therefore $a=1$ and $b=0$.
    Since $e_0, f_0$ is a basis of the $K$-vector space $C_{\rm sp}$, the determinant
    \begin{displaymath}
      \det \left(
      \begin{array}{cc} 
        \varkappa(e_0, e_0) & \varkappa(e_0,f_0)\\
        \varkappa(f_0, e_0) & \varkappa(f_0, f_0) 
      \end{array}
      \right)
      = 1
      \end{displaymath}
    gives the invariant
    $1 = \inv(C_{\rm sp}, \psi) = - \inv (\mathcal{P}_{\rm sp},\iota'_{\rm sp},\psi)$
    by the last Lemma. 
\end{proof}

\begin{proof} (of Proposition \ref{ramFrame1p})
  By Lemma \ref{ramFrame1l} we know that $\mathcal{P}_{{\rm{c}},i}$ is isoclinic of
  slope $1/2$ for $i = 1,2$. Therefore Lemma \ref{ramFrame2l} is applicable.
  By Lemma \ref{Kneun2l}, the associated $K$-vector spaces $(C_i, \beta_{_{{\rm{c}},i}})$
  have the same  invariant $1$ and are therefore isomorphic.
  Therefore we find the quasi-isogeny  $\alpha$ by Lemma \ref{ramFrame2l}.

  We use the notations of Lemma \ref{ramFrame4l}. By what we just proved we
  find a quasi-isogeny
  $(\mathcal{P}_{\rm sp},\iota'_{\rm sp},\psi) \longrightarrow (\mathcal{P}_{{\rm{c}},1}, \iota_{{\rm{c}},1}, \beta_{{\rm{c}},1})$.
  If $\beta_{{\rm{c}},1}$ is perfect, this quasi-isogeny extends by Proposition
  \ref{Pol4p} to a quasi-isogeny of special formal $O_D$-modules and
  so does $\alpha$. 
\end{proof}

We can now prove the uniqueness of the framing object.  
\begin{proposition}\label{ramFrame2p}
  Let $r$ be special and $K/F$  ramified. Let $k$ be an algebraically
  closed field in $\Nilp_{O_E}$. 
  Let $(\mathcal{P}, \iota, \beta) \in \mathfrak P^{\rm pol}_{r,k}$ be an object such that $\beta$ is
  perfect, cf.
  Definition \ref{KatCMtriple1d}. Assume that $\inv^r (\mathcal{P}, \iota, \beta) = -1$. Then $\mathcal{P}$ is isoclinic of slope $1/2$.  
  
  If moreover
  $(\mathcal{P}_1, \iota_1, \beta_1)$ is a second triple
  with the same properties, then there is a quasi-isogeny of height zero
  \begin{displaymath}
    \rho: (\mathcal{P}, \iota, \beta) \longrightarrow
    (\mathcal{P}_1, \iota_1, \beta_1) ,
  \end{displaymath}
  such that there is an $f \in O_F^{\times}$ with
  \begin{displaymath}
\beta_1(\rho(x), \rho(y)) = \beta(fx, y), \quad x,y \in P. 
    \end{displaymath}
\end{proposition}
\begin{proof}
  We apply the functor $\mathfrak{C}^{\rm pol}_{r, k}$ to
  $(\mathcal{P}, \iota, \beta)$ and obtain $(\mathcal{P}_{\rm{c}}, \iota_{\rm{c}}, \beta_{\rm{c}})$, cf. (\ref{P'Pol12e}).
  By the definition of this functor, $\beta_{\rm{c}}$ is perfect. We conclude from
  Proposition \ref{Kneuninv1p} that
  $\inv (\mathcal{P}_{\rm{c}}, \iota_{\rm{c}}, \beta_{\rm{c}}) = -1$. By Lemma \ref{ramFrame1l},
  $\mathcal{P}_{\rm{c}}$ is isoclinic of slope $(1/2)$. By Corollary
  \ref{C'onslopes1c} and  Proposition \ref{Adorfslopes1p},
  $\mathcal{P}$ is isoclinic of slope $1/2$.  
  
  By Proposition \ref{ramFrame1p}, we find a quasi-isogeny
  $\alpha: (\mathcal{P}_{\rm{c}}, \iota_{\rm{c}}, \beta_{\rm{c}}) \longrightarrow                          (\mathcal{P}_{{\rm{c}},1}, \iota_{{\rm{c}},1}, \beta_{{\rm{c}},1})$
  which we can make into a quasi-isogeny of special formal $O_D$-modules. 
  The height of $\alpha$ is then a multiple of $2$. Composing $\alpha$
  with an endomorphism of the special formal $O_D$-module
  $(\mathcal{P}_{\rm{c}}, \iota_{\rm{c}})$, we can obtain a quasi-isogeny of height $0$ of $O_D$-modules
  $\rho_{\rm{c}}: (\mathcal{P}_{\rm{c}}, \iota_{\rm{c}}) \longrightarrow (\mathcal{P}_{{\rm{c}},1}, \iota_{{\rm{c}},1})$. 
  Then $\rho_{\rm{c}}$ respects the Drinfeld polarizations $\beta_{\rm{c}}$ and $\beta_{{\rm{c}},1}$
  up to a constant in $O_F^{\times}$. By Theorem \ref{P'Pol2c}, we obtain a
  quasi-isogeny of height zero as claimed in the proposition. 
  \end{proof}
 
    \begin{remark}
      We chose here the framing object for $\CN$ as coming from the Drinfeld moduli problem. It can also be characterized in terms of the moduli problem $\CN$, cf. \cite{Ki}: it is a triple $(X, \iota, \lambda)$ consisting of a $p$-divisible strict formal $O_F$-module $X$ over $\ov \kappa_F$, with an action $\iota$ of $O_K$ satisfying the Kottwitz condition \eqref{kottforrel}, and a perfect relative polarization   $\lambda$ such that the special  automorphism group is isomorphic to $\SL_2(F)$, comp. \cite[Prop. 3.2]{Ki}.
    \end{remark}

\subsection{The alternative theorem in the unramified case}\label{ss:altunram} 

In this subsection  $K$ denotes an  {\it unramified} quadratic
extension of $F$.  
Let $k$ be an algebraically closed field of characteristic $p$ which
is endowed  with an $O_{F}$-algebra structure.
We will sometimes write $F' = K$ if we refer to subsection \ref{ss:sfmod}.
Let $\tau$ be the Frobenius of $F'/F$. We write $\tau(a) = \bar{a}$ for
$a \in O_{F'}$. 
    
    Let $M$ be the $\mathcal{W}_{O_F}(k)$-Dieudonn\'e module of a special
    formal $O_D$-module over $k$, as in Proposition \ref{Pol1p}.
    In
    addition to the Drinfeld polarization, we use another type of polarization
    of $M$,    
    \begin{displaymath}
      \theta: M \times M \longrightarrow W_{O_F}(k).
    \end{displaymath} 
    This is an alternating bilinear form of $M$ which satisfies
    \begin{equation}\label{Pol21e}
      \begin{aligned}
        \theta(Fx_1, Fx_2) &= \pi ~^F\theta(x_1, x_2), & x_1, x_2 \in M\\
        \theta(\iota(a)x_1, x_2) &= \theta(x_1, \iota(\bar{a})x_2), &
        a \in O_{F'},\\
        \theta(\iota(\Pi)x_1, x_2) &= \theta(x_1, \iota(\Pi)x_2),\\
        \ord_{\pi} \det \theta &= 2. 
      \end{aligned}
    \end{equation}   
    The polarization $\theta$ is unique up to a constant in $O_F^{\times}$. It
    is constructed as follows: We choose an element $\delta \in
    O_{F'}^{\times}$, such that $\delta + \tau(\delta) = 0$.
    We set $\Pi_1 = \delta \Pi$. Then $\Pi_1$ is invariant under the
    involution (\ref{Pol1e}) and therefore we have
    \begin{displaymath}
      {\psi}(\iota(\Pi_1) x, y) = {\psi}(x, \iota(\Pi_1)y). 
    \end{displaymath}
    We define
    \begin{equation}\label{def:theta} 
      \theta(x, y) = {\psi}(\iota(\Pi_1) x, y). 
    \end{equation}
    We see that $\theta$ is alternating. It induces on $D$ the involution
    given by
    \begin{equation}\label{Pol31e}
      \Pi^{\dagger} = \Pi, \quad u^{\dagger} = \tau(u), \quad \text{for}
      \; u \in F'.  
    \end{equation}
    Conversely, assume that $\theta$ is a polarization with the properties \eqref{Pol21e}. Let ${\psi}_1(x, y) = \theta(\Pi_1x, y)$. Using
    $\Pi_1 \Pi = - \Pi \Pi_1$ we see that ${\psi}_1$ satisfies
    the properties 
    (\ref{Pol2e}). By Proposition \ref{Pol1p} this shows the uniqueness
    of $\theta$ with the properties above.

    Let $(\mathbb{Y}, \iota_{\mathbb{Y}})$ be a special formal
    $O_{D}$-module over the $O_F$-algebra $\ov\kappa_{F}$ 
    and such that $\iota(\Pi)$ acts as zero on $\Lie \mathbb{Y}$. We
    endow $\mathbb{Y}$ with the polarization $\theta_{\mathbb{Y}}$
    defined above, cf. (\ref{def:theta}). 
    
    \begin{definition}\label{N-unram1d} 
    We define for each $i \in \mathbb{Z}$ the functor $\CN(i) =\CN_{K/F}(i)$
    on the category $({\rm Sch}/\Spf O_{\breve{F}})$. A point of
    $\mathcal{N}(i)(S)$ consists of the following data: 
    \begin{enumerate}
    \item A formal $p$-divisible group $X$ over $S$ with an action
      \begin{displaymath}
        \iota: O_K \longrightarrow \End X,
      \end{displaymath}
      such that the restriction of $\iota$ to $O_F$ is a strict action. 
    \item
      A  relative polarization $\theta$ on $X$ such that the
      determinant of $\theta$ is $\pi^{2}$ up to a unit and such that
      $\theta$ induces on $O_K$ the conjugation over $O_F$. 
    \item A quasi-isogeny of $O_K$-modules 
        \begin{displaymath}
        \rho: X \times_{S} \ov S \longrightarrow   \mathbb{Y}
        \times_{\Spec \bar{\kappa}_F} \ov S. 
      \end{displaymath}
    \end{enumerate}
    Here, if $S=\Spec R$,  the condition in (2) means that the
    polarization of the corresponding $O_F$-display $\CP$ of $Y$ has
    determinant $\pi^2$, up to a unit in $W_{O_F}(R)$.  
    We require that the following conditions are satisfied.
    \begin{enumerate}
    \item[a)] $\rho$ respects $O_K$-actions. The relative quasi-polarization
      $\rho^{\ast} \theta_{\mathbb{Y}}$ differs from $\pi^i\theta$ by a factor
      in $O_F^{\times}$.  
          \item[b)] $\Lie X$ is locally on $S$ a free
      $O_K \otimes_{O_F} \mathcal{O}_S$-module of rank $1$. 
         \end{enumerate}

    We note that, as in the ramified case, the $O_F$-height of $X$ is $4$ and
    the dimension $2$. The condition b) implies the
    following Kottwitz condition for the characteristic polynomial: 
    \begin{displaymath}
      {\rm char} (\iota(a) \; | \; \Lie X) = (T- a)(T - \bar{a}), \quad a \in
      O_K. 
    \end{displaymath}

    Two data $(X_1, \iota_1, \theta_1, \rho_1)$ and
    $(X_2, \iota_2, \theta_2, \rho_2)$ define the same point of
    $\mathcal{N}(i)(S)$ iff there is an  isomorphism
    $\alpha: (X_1, \iota_1) \longrightarrow (X_2, \iota_2)$ which respects the
    polarizations up to a factor in $O_{F}^{\times}$ and such that $\alpha$
    commutes with $\rho_1$ and $\rho_2$.   
      \end{definition}

    It follows from \cite{RZ} that $\mathcal{N}(i)$ is representable by a
    formal scheme which is locally formally of finite type over
    $\Spec O_{\breve{F}}$. The functor $\mathcal{N}(0)$ will be also denoted
    by $\CN$.

    We have a natural functor morphism
    \begin{equation}\label{Pol27e}
      \CM_{\rm Dr}(i) \longrightarrow \mathcal{N}(i) .
    \end{equation}
    Indeed, let $(Y, \rho) \in   {\CM}_{\rm Dr}(i) (R)$.  Then we have the
    Drinfeld polarization ${\psi}$ of $Y$ and we define $\theta_{Y}$ by
    the formula (\ref{def:theta}). This gives a point of $\CN(R)(i)$. 

    The diagram similiar to (\ref{Mi4e}) shows that
    \begin{displaymath}
\height_{O_F} \rho = 2i. 
    \end{displaymath}
    We will define a translation functor isomorphism 
    \begin{equation}\label{translation3e} 
\Pi: \mathcal{N}(i) \isoarrow \mathcal{N}(i+1).  
      \end{equation}
    Let $(Y, \iota)$ be a special formal $O_D$-module over $R \in \Nilp_{O_F}$.
    We fix a Drinfeld polarization $\psi$. This is also a Drinfeld
    polarization for $(Y^{\Pi}, \iota^{\Pi})$. For the polarizations $\theta$
    and $\theta^{\Pi}$ derived by (\ref{def:theta}), we obtain
    $\theta^{\Pi} = - \theta$. We consider the morphism
    $\iota^{\Pi}: Y^{\Pi} \longrightarrow Y$. One easily checks that
    \begin{displaymath}
\theta(\iota(\Pi) x, \iota(\Pi) y) = \pi \theta^{\Pi}(x, y).
      \end{displaymath}
    This is an identity of bilinear forms on the $\mathcal{W}_{O_F}(R)$-display
    of $Y$.

    If $(X,\iota)$ is a $p$-divisible $O_K$-module, we define the conjugate\index{conjugate
    $p$-divisible $O_K$-module}\index[NO]{XAC@$(X^{\rm c}, \iota^{\rm c})$}
    $p$-divisible $O_K$-module $(X^{\rm c}, \iota^{\rm c})$ by setting $X^{\rm c} = X$ and
    $\iota^{\rm c}(a) = \iota(\bar{a})$ for $a \in O_K$. For the special formal
    $O_D$-module $Y$ we have
    \begin{displaymath}
(Y^{\rm c}, (\iota_{|O_K})^{\rm c}) = (Y^{\Pi}, \iota^{\Pi}_{|O_K}).
      \end{displaymath}

    Let $R \in \Nilp_{O_{\breve{F}}}$ and let 
   $(X, \iota, \theta, \rho) \in \mathcal{N}(i)(R)$. We define
    \begin{displaymath}
      \rho^{\rm c}: X^{\rm c}_{\bar{R}} \; \overset{\rho}{\longrightarrow}\;
      \mathbb{Y}^{\rm c}_{\bar{R}} = \mathbb{Y}^{\Pi}_{\bar{R}} \;
      \overset{\iota(\Pi)}{\longrightarrow} \; \mathbb{Y}_{\bar{R}}. 
      \end{displaymath}
    We set $\theta^{\rm c} = -\theta$. Then
    $(X^{\rm c}, \iota^{\rm c}, \theta^{\rm c}, \rho^{\rm c}) \in \mathcal{N}(i+1)(R)$. This defines
    the translation functor morphism (\ref{translation3e}). It is clearly an
    isomorphism. With this definition, the functor morphism
    (\ref{Pol27e}) commutes with the translations on source and target. 

Let $\tau \in \Gal(\breve{F}/F)$ be the Frobenius automorphism. Using the
Frobenius $F_{\mathbb{Y}, \tau}:\mathbb{Y}\longrightarrow\tau_{\ast}\mathbb{Y}$,
we obtain a morphism
\begin{displaymath}
\omega_{\mathcal{N}}: \mathcal{N}(i) \longrightarrow \mathcal{N}(i+1)^{(\tau)} 
  \end{displaymath}
with the same definition as (\ref{WD2e}). This induces a Weil  descent datum
$\omega_{\mathcal{N}}$\index[NO]{ZZWF@$\omega_{\mathcal{N}}$}  on
\index[NO]{NBB@$\tilde{\mathcal{N}}$}
\begin{displaymath}
    \tilde{\mathcal{N}} = \coprod_{i \in \mathbb{Z}} \mathcal{N}(i).  
\end{displaymath}
\begin{lemma}\label{J^uv1l}
  The action of the group $J^{\ast{\rm ur}}$ (cf. (\ref{Pol29J}) for $F'$, which is
  now denoted by $K$) on the $O_K$-module $\mathbb{Y}$ gives an isomorphism
 \begin{equation*}
 J^{\ast{\rm ur}}\isoarrow J^{\cdot} ,
 \end{equation*} 
 where 
\begin{displaymath}
  J^{\cdot} = \{\alpha \in \Aut^o_{K}(\mathbb{Y}) \; | \;
  \theta_{\mathbb{Y}}(\alpha(x), \alpha(y)) =
  \mu(\alpha) \theta_{\mathbb{Y}}(x,y), \; \text{for some} \;
  \mu(\alpha) \in F^{\times}, \;  
      x,y \in P_{\mathbb{Y}} \otimes \mathbb{Q}\}
      \end{displaymath}\qed
\end{lemma}
The group $J^{\cdot}$ acts via the rigidification $\rho$ on the functor 
$\tilde{\mathcal{N}}$. 

    \begin{theorem} [\cite{KRalt}]\label{thm:altdrinunr} 
      The morphisms of functors (\ref{Pol27e})
       for varying $i$ extend to a functor isomorphism 
      \begin{displaymath} 
         \tilde{\mathcal{M}}_{\rm Dr}
  \isoarrow \tilde{\mathcal{N}}
        \end{displaymath}
      which commutes with the Weil descent data, the actions of
      $J^{\ast{\rm ur}} = J^{\cdot}$, and the translations on both sides. 
    \end{theorem}
    \begin{proof} We already checked that (\ref{Pol27e}) extends to a functor morphism which  respects translations and
      Weil  descent data on both sides. Therefore it suffices to see that
      (\ref{Pol27e}) is an isomorphism for $i = 0$,
      \begin{equation}\label{Pol29e}
\mathcal{M}_{\mathrm{Dr}} \isoarrow \mathcal{N}.  
        \end{equation}
          We begin with the case where $R = k$ is an algebraically closed
      field. Let $Y \in \CN(k)$. 
      Let $M$ be the  $O_F$-Dieudonn\'e module of $Y$ and let
      $\mathbb{M}$ be the $O_F$-Dieudonn\'e module of $\mathbb{Y}$.  The
      quasi-isogeny $\rho$ induces an isomorphism
      $M \otimes \mathbb{Q} \cong \mathbb{M} \otimes \mathbb{Q}$. The
      polarization $\psi_{\mathbb{M}}$ induces a polarization
      ${\psi}_1$ on $M \otimes \mathbb{Q}$. Since $\rho$ is of height zero
      and $\ord_{\pi} \det \psi_{\mathbb{M}} = 0$ we conclude that
      $\ord_{\pi} \det {\psi}_1 = 0$. On the other hand, we have by
      Proposition \ref{Pol3p} a perfect 
      pairing ${\psi}$ on $M$ which differs from ${\psi}_1$ by an element $f \in
      F'$. This shows that ${\psi}_1$ is perfect on $M$. Then we define the
      action $\iota(\Pi_1)$ by the equation
      \begin{equation}\label{ThetaPol1e}
        \theta(x, y) = {\psi}(\iota(\Pi_1) x, y),  \quad x,y \in M.
      \end{equation}
      Therefore the morphism (\ref{Pol29e}) evaluated at $k$ is bijective. 

      Since both functors of (\ref{Pol29e}) are
      representable by formal schemes locally of finite type, it suffices
      now to check the following statement. Let $S \longrightarrow R$ be a
      surjective $O_{\breve{F}}$-algebra homomorphism
      such that $S$ and $R$ are artinian local rings with
      algebraically closed residue class field.  Assume that (\ref{Pol29e}) is
      bijective when evaluated at $R$. Then  it is
      bijective when evaluated at $S$. We may assume that the kernel of
      $S\to R$ is endowed with divided powers. 

      We consider a point $\wt{Y} \in
      \CN(S)$ and we denote by $Y \in
      \CN(R)$ its reduction. By our assumption,  $Y$
      carries the structure of an $O_D$-module compatible with
      $\rho$. Therefore ${\psi}_{\mathbb{Y}}$ induces a perfect 
      polarization on the $O_F$-display $\mathcal{P}$ of $Y$. The $O_F$-display $\wt{\mathcal{P}}$ of $\wt Y$ is a lifting of
      $\mathcal{P}$. By the crystalline property of displays \cite{ACZ}, cf. end of subsection \ref{ss:reldisp}, we
      obtain a perfect pairing
      \begin{displaymath}
        \wt{{\psi}}: \wt{P} \times \wt{P} \longrightarrow W_{O_F}(S). 
      \end{displaymath}
      The involution induced by $\wt{\psi}$ on $O_{F'}$ is trivial. It
      follows that the decomposition 
      \begin{displaymath}
        \wt{P} = \wt{P}_0 \oplus \wt{P}_1
      \end{displaymath}
      according to the two $O_F$-algebra embeddings $O_{F'} \longrightarrow
      O_{\breve{F}}$ is orthogonal with respect to $\wt{\psi}$. The Hodge
      filtration
      \begin{displaymath}
        \tilde{Q}_i/I_{O_F}\tilde{P}_i   \subset
        \tilde{P}_i/I_{O_F}\tilde{P}_i, \quad  i = 0,1 
      \end{displaymath} 
      is isotropic with respect to $\wt{\psi}$ because these direct
      summands are of rank $1$.  Therefore $\wt{\psi}$ is a 
      polarization of the $O_F$-display $\wt{\CP}$. Using the given polarization
      $\wt{\theta}$ on $\wt{\CP}$, we can define the endomorphism
      $\iota(\Pi_1)=\iota(\delta\Pi)$ of $\wt{P}$ by
      \begin{displaymath}
        \wt{\theta}(x, y) = \wt{{\psi}}(\iota(\Pi_1) x, y).
      \end{displaymath}
      This gives the desired $O_D$-module structure on $\wt{\CP}$ and
      therefore on $\wt Y$.
    \end{proof}
  
    The analogue of Corollary \ref{spfibred} follows as before from
    the properties of the Drinfeld moduli space. 
     \begin{corollary}\label{spfibred2}
      The formal scheme $\CN$ is $\pi$-adic and has semi-stable reduction. The special fiber $\CN\otimes_{O_{\breve F}}\ov\kappa_{F}$ of $\CN$ is a reduced scheme. \qed
    \end{corollary}
We next prove the uniqueness of the framing object, cf. (i) of subsection \ref{ss:uniquefram}. We start with the following statement. 

    \begin{proposition}\label{unrFrame1p}
Let $k$ be an algebraically closed field which contains $\kappa_F$.       
Let $M$ be a $\mathcal{W}_{O_F}(k)$-Dieudonn\'e module of height $4$ and
dimension $2$. Let $\iota$ be a homomorphisms of $O_F$-algebras
\begin{displaymath}
\iota: O_{F'} \longrightarrow \End M.  
\end{displaymath}
Assume that $M/VM$ is a free $\kappa_{F'} \otimes_{\kappa_F} k$-module of
rank 1. Let $\theta$ be a relative polarization on $M$ which satisfies  
\begin{displaymath}
      \begin{aligned}
        \theta(\iota(a)x_1, x_2) &= \theta(x_1, \iota(\bar{a})x_2), &
        a \in O_{F'},\\
  \ord_{\pi} \det \theta &= 2. 
      \end{aligned}
  \end{displaymath}
Then the action $\iota$ extends to an action $\iota: O_D \longrightarrow \End M$
 such that $\theta$ satisfies (\ref{Pol21e}). In particular $M$ is isoclinic
 of slope $1/2$. Furthermore,  $\inv(M, \iota, \theta) = -1$ {\rm (see Definition \ref{invcontr1d} for this invariant)}.

 If $(M',\iota', \theta')$ is a second triple with the same properties, then
 there exists a quasi-isogeny $(M, \iota) \longrightarrow (M', \iota')$ of height
 $0$ which respects the polarizations $\theta$ and $\theta'$ up to a factor in 
 $O_F^{\times}$. 
    \end{proposition}
    \begin{proof}
 Let $\psi$ be the principal relative polarization on $M$ which exists by
 Proposition \ref{Pol3p}. We define an endomorphism
 $\rho: M \longrightarrow M$ by the equation
 \begin{displaymath}
\theta(x,y) = \psi(x, \rho(y)), \quad x,y \in M. 
 \end{displaymath}
 One checks that $\rho$ is an endomorphism of the Dieudonn\'e module $M$
 such that
 \begin{equation}\label{uRahmen1e}
\rho(\iota(a)x) = \iota(\bar{a}) \rho(x), \quad a \in O_{F'}. 
 \end{equation}
 As in the proof of Proposition \ref{Pol3p}, we choose an embedding
 $\lambda: O_{F'} \longrightarrow W_{O_F}(k)$ and obtain a decomposition
 $M = M_0 \oplus M_1$. We note that
 \begin{displaymath}
M_1 = \{ x \in M \; | \; \iota(a)x = \lambda(\bar{a})x \}. 
 \end{displaymath}
 It follows from (\ref{uRahmen1e}) that $\rho(M_0) \subset M_1$ and
 $\rho(M_1) \subset M_0$. We obtain a commutative diagram
 \begin{displaymath}
   \xymatrix{
M_0 \ar[d]_{V} \ar[r]^{\rho} & M_1 \ar[d]^{V}\\
M_1 \ar[r]^{\rho} & M_0.\\
}
 \end{displaymath}
 By our assumption on $M/VM$, the cokernels of both vertical maps have
 $W_{O_F}(k)$-length $1$. Therefore the cokernels of the horizontal maps
 have also the same length. This length must be $1$ because
 $\ord_{\pi} \det (\rho | M) = \ord_{\pi} \det \theta = 2$.

 We have 
 \begin{displaymath}
\theta(\rho(x), y) = \theta(x, \rho(y)), 
 \end{displaymath}
 because both sides are equal to $\psi(\rho(x), \rho(y))$.
 We consider the form
 \begin{displaymath}
\psi_1(x, y) := \psi(\rho(x), \rho(y)). 
 \end{displaymath}
 This relative polarization satisfies the assumptions of the  last part of
 Proposition \ref{Pol3p}. Therefore there exists $c \in O_{F'}$ such that
 $\psi_1(x,y) = \psi(\iota(c)x, y)$. We find
 \begin{displaymath}
   \psi(\iota(c)x, y) = \theta(\rho(x), y) = - \theta(y, \rho(x)) =
   - \psi(y, \rho^2(x)) = \psi(\rho^2(x), y).
   \end{displaymath}
 This shows that
 \begin{displaymath}
\rho^2 = \iota(c). 
 \end{displaymath}
 Since $\rho$ commutes with the left hand side it commutes with $\iota(c)$.
 Comparing this with (\ref{uRahmen1e}), we obtain $c \in O_F$. Since $\rho^2$
 has height $4$ we obtain $\ord_\pi c = 1$.

 For $\alpha \in O_{F'}^{\times}$ we consider the endomorphism
 $\rho_{\alpha}(x) = \iota(\alpha) \rho(x)$ of $M$. We obtain
 $\rho_{\alpha}^2 = \iota(\alpha \bar{\alpha}) \rho^2$. Since each unit of
 $F$ is a norm in the unramified extension $F'/F$ we can arrange that
 $\rho_{\alpha}^2 = -\pi$. We set $\Pi = \rho_{\alpha}$. Then we obtain an
 action of $O_D = O_{F'}[\Pi]$ on the Dieudonn\'e module $M$. Since
 $M_0/\Pi M_1$ and $M_1/\Pi M_0$ have length $1$, we have obtained a special
 formal $O_D$-module. The equations (\ref{Pol21e}) are satisfied for
 $\theta$. Therefore $\theta$ is up to a factor in $O_F^{\times}$ uniquely
 determined by the $O_D$-action. This implies the first and the last
 assertion of the Proposition. 

 Because the invariant depends only on the isogeny class, it is enough to
 compute it for a special formal $O_D$-module with two critical indices
 and the canonical form $\theta$ from (\ref{Pol21e}).

 We use the isomorphism
 \begin{equation}\label{uRahmen2e}
O_{F'} \otimes_{O_F} W_{O_F}(k) \longrightarrow W_{O_F}(k) \times W_{O_F}(k), 
   \end{equation}
 which maps $a \otimes \xi$ to $(\lambda(a)\xi, \lambda(\bar{a}) \xi)$. 
 Let $\sigma = F$ be the Frobenius automorphism of $W_{O_F}(k)$. It acts on the
 left hand side of (\ref{uRahmen2e}) via the second factor. This induces on
 the right hand side the action
 $\sigma: (\xi_1, \xi_2) \mapsto (\sigma (\xi_2), \sigma(\xi_1))$.

 We set $N_i = M_i \otimes \mathbb{Q}$. This is a $W_{O_F}(k)_{\mathbb{Q}}$-vector
 space of dimension $2$. In the decomposition
 $N = M \otimes \mathbb{Q} = N_0 \otimes N_1$, the summands are isotropic with respect to
 $\theta$.  We consider the invariants
 \begin{displaymath}
U_i = N_i^{V^{-1} \iota(\Pi)}. 
 \end{displaymath}
 (In the notation of the proof of Proposition \ref{Pol2p}, this is
 $U_i \otimes \mathbb{Q}$.) The $U_i$ are $F$-vector spaces of dimension 2.
 Let
 \begin{displaymath}
   \wedge^2 \theta: \wedge^2 N_0 \times \wedge^2 N_1 \longrightarrow
   W_{O_F}(k)_{\mathbb{Q}} 
   \end{displaymath}
 be the bilinear form defined by
 \begin{displaymath}
   \theta (n_0 \wedge n'_0, n_1 \wedge n'_1) = 
  \det \left(
      \begin{array}{cc} 
        \theta(n_0, n_1) & \theta(n_0, n'_1)\\
        \theta(n'_0, n_1) & \theta(n'_0, n'_1) 
      \end{array}
      \right), 
 \end{displaymath}
for $n_0, n'_0 \in N_0$, $n_1, n'_1 \in N_1$. In the same way we can define
$\wedge^2\theta:\wedge^2 N_1 \times\wedge^2 N_0\longrightarrow W_{O_F}(k)_{\mathbb{Q}}$.
Then we obtain $\wedge^2 \theta(x_0, x_1) = \wedge^2 \theta(x_1, x_0)$ for
$x_0 \in \wedge^2 N_0, x_1 \in \wedge^2 N_1$. From $\theta$ we pass to
$\varkappa$, cf. (\ref{Kneunpsi1e}) (there our $F'$ is called $K$), 
\begin{displaymath}
  \varkappa: N \times N \longrightarrow F' \otimes_{F} W_{O_F}(k)_{\mathbb{Q}} \cong
  W_{O_F}(k)_{\mathbb{Q}} \times W_{O_F}(k)_{\mathbb{Q}}. 
\end{displaymath}
Explicitly we have
\begin{displaymath}
  \varkappa(n_0 + n_1, n'_0 + n'_1) = (\theta(n_0,n'_1), \theta(n_1,n'_0)) \in
  W_{O_F}(k)_{\mathbb{Q}} \times W_{O_F}(k)_{\mathbb{Q}}. 
\end{displaymath}
We take $\wedge^2 \varkappa$ on the $F' \otimes_{F} W_{O_F}(k)_{\mathbb{Q}}$-module
\begin{displaymath}
  \bigwedge_{F' \otimes_{F} W_{O_F}(k)_{\mathbb{Q}}}^2 N \cong \wedge^2 N_0 \oplus
  \wedge^2 N_1. 
  \end{displaymath}
From the expression for $\varkappa$ we obtain
\begin{equation}\label{uRahmen4e}
  \wedge^2 \varkappa(x_0 + x_1, x'_0 + x'_1) =
  (\wedge^2 \theta (x_0, x'_1), \wedge^2 \theta (x_1, x'_0)) \in
  W_{O_F}(k)_{\mathbb{Q}} \times W_{O_F}(k)_{\mathbb{Q}}
\end{equation}
The restriction of $\theta$ to $U_0 \times U_1$ induces a nondegenerate
$F$-bilinear form
\begin{equation}\label{uRahmen3e}
  \theta: U_0 \times U_1 \longrightarrow \delta F \subset F' \subset
  W_{O_F}(k)_{\mathbb{Q}}, 
\end{equation}
where $\delta$ was defined after (\ref{Pol21e}). 
Indeed, for $u_0 \in U_0$ and $u_1 \in U_1$ we have by definition
\begin{displaymath}
  \theta(V u_0, V u_1) = \theta(\iota(\Pi)u_0, \iota(\Pi)u_1) =
  \theta(\iota(\Pi)^2 u_0, u_1) = -\pi \theta(u_0, u_1). 
  \end{displaymath}
Because $\theta$ is a polarization, we have on the other hand 
\begin{displaymath}
\theta(V u_0, V u_1) = \pi \sigma^{-1} \big(\theta(u_0, u_1)\big). 
\end{displaymath}
Therefore $\theta(u_0, u_1)$ is anti-invariant by $\sigma$ and \eqref{uRahmen3e}
is proved.

We choose a nonzero element $u_0 \in U_0$. Then we find $u_1 \in U_1$ such that
$\theta(u_0, u_1) = \delta$. We remark that $\theta(\iota(\Pi) n, n) = 0$, for an arbitrary $n \in N$. This is clear from the third equation of
$(\ref{Pol21e})$ because $\theta$ is alternating.
Since $\theta(u_0, \iota(\Pi) u_0) = 0$, the vectors
$u_1, \iota(\Pi) u_0 \in U_1$ are linearly independent. We set 
\begin{displaymath}
x = x_0 + x_1 := u_0 \wedge \iota(\Pi) u_1 + u_1 \wedge \iota(\Pi) u_0 \in
\wedge_{F' \otimes_{F} W_{O_F}(k)_{\mathbb{Q}}}^2 N. 
\end{displaymath}
It satisfies $\wedge^2 V x = \pi x$. Indeed,
\begin{displaymath}
  \wedge^2 V (u_0 \wedge \iota(\Pi) u_1) = (Vu_0 \wedge \iota(\Pi) Vu_1) =
  (\iota(\Pi) u_0 \wedge \iota(\Pi)^{2} u_1) =
  (\iota(\Pi) u_0 \wedge (- \pi) u_1) = \pi (u_1 \wedge \iota(\Pi) u_0).   
\end{displaymath}
The similiar equation holds for the second summand in the definition of $x$.
By Definition \ref{invcontr1d} we obtain
\begin{displaymath}
\inv(M, \iota, \theta) = (-1)^{\ord_{\pi} \wedge^2 \varkappa (x,x)} .
  \end{displaymath}
By (\ref{uRahmen4e}) we have
$\ord_{\pi} \wedge^2 \varkappa (x,x) = \ord_{\pi} \wedge^2 \theta (x_0,x_1)$.
We compute
\begin{displaymath}
  \wedge^2 \theta (x_0, x_1)
  = \det \left(
      \begin{array}{cc} 
        \theta(u_0, u_1) & \theta(u_0, \iota(\Pi)u_0)\\
        \theta(\iota(\Pi) u_1, u_1) &
        \theta(\iota(\Pi)u_1, \iota(\Pi)u_0) 
      \end{array}
      \right)
    = \det \left(
      \begin{array}{cc} 
        \delta & 0\\
        0 & \delta \pi 
      \end{array}
      \right)
      = \pi \delta^2.
\end{displaymath}
This shows $\ord_{\pi} \wedge^2 \theta (x_0,x_1) = 1$ and therefore
$\inv(M, \iota, \theta) = -1$. 
      \end{proof}

We can now prove the uniqueness of the framing object. 
\begin{proposition}\label{unrFrame2p}
  Let $r$ be special and $K/F$  unramified. Let $k$ be an algebraically
  closed field in $\Nilp_{O_E}$. 
  Let $(\mathcal{P}, \iota, \beta) \in \mathfrak P^{\rm pol}_{r,k}$ be 
  such that $\ord_p \det_{W(k)} \beta = 2f$, cf.
  Definition \ref{KatCMtriple1d}.
Then $\mathcal{P}$ is isoclinic of slope $1/2$ and
  $\inv^{r} (\mathcal{P}, \iota,\beta) = -1$. 
  
  If
  $(\mathcal{P}_1, \iota_1, \beta_1)$ is a second triple
  with the same properties, then there exists a quasi-isogeny of height zero,
  \begin{displaymath}
    \rho: (\mathcal{P}, \iota, \beta) \longrightarrow
    (\mathcal{P}_1, \iota_1, \beta_1) , 
  \end{displaymath}
  such that there is an $f \in O_F^{\times}$ with
  \begin{displaymath}
\beta_1(\rho(x), \rho(y)) = \beta(fx, y), \quad x,y \in P. 
    \end{displaymath}
\end{proposition}
\begin{proof}
  We apply the functor $\mathfrak{C}^{\rm pol}_{r, k}$ to
  $(\mathcal{P}, \iota, \beta)$ and obtain $(\mathcal{P}_{\rm{c}}, \iota_{\rm{c}}, \beta_{\rm{c}})$, cf. (\ref{P'Pol12e}).
  By Theorem \ref{P'Pol2c} we find $\ord_{\pi} \det_{W_{O_F}(k)} \beta_{\rm{c}} = 2$.
  Therefore we can apply Proposition \ref{unrFrame1p} to
  $(\mathcal{P}_{\rm{c}}, \iota_{\rm{c}}, \beta_{\rm{c}})$. We obtain that $\mathcal{P}_{{\rm{c}}}$ is
  isoclinic of slope $1/2$ and $\inv(\mathcal{P}_{\rm{c}}, \iota_{\rm{c}}, \beta_{\rm{c}})=-1$. By Corollary \ref{C'onslopes1c} and
  Proposition \ref{Adorfslopes1p}, we find that
  $\mathcal{P}$ is isoclinic of slope $1/2$, and by Proposition
  \ref{Kneuninv1p} we obtain $\inv^r (\mathcal{P}, \iota, \beta) = -1$.  

  By Proposition \ref{unrFrame1p}, there is a quasi-isogeny of height $0$
  between $(\mathcal{P}_{\rm{c}}, \iota_{\rm{c}}, \beta_{\rm{c}})$ and
  $(\mathcal{P}_{1, c}, \iota_{1, c}, \beta_{1, c})$. It induces by Theorem
  \ref{P'Pol2c} a quasi-isogeny of height zero as claimed in the Proposition. 
  \end{proof}

We end this section by justifying the footnote in Definition \ref{CRFd2}. Let $S \in \Nilp_{O_{\breve{F}}}$ and let
$(Y, \iota, \theta, \rho) \in \mathcal{N}(S)$. Since there is an
$O_D$-module structure on $Y$ such that $\theta$ is of the form
\eqref{ThetaPol1e}, the kernel of $\theta: Y \longrightarrow Y^{\wedge}$,
considered as morphism to the dual relative to $O_F$, is annihilated
by $\pi$. More generally we prove:

    \begin{proposition}\label{SFB701p}
Let $K/F$ be an unramified quadratic field extension. 
Let $R$ be an $O_K$-algebra. 
Let $\mathcal{P}$ and $\mathcal{P}'$ be
$\mathcal{W}_{O_F}(R)$-displays of height $4$ with an action 
$\iota: O_K \longrightarrow \End \mathcal{P}$, resp. $\iota': O_K
\longrightarrow \End \mathcal{P}'$.  Assume that $\Lie \mathcal{P}$,
resp. $\Lie \mathcal{P}'$, is locally on $\Spec R$ a free $O_K
\otimes_{O_F} R$-module of rank $2$.  
Let $\alpha: \mathcal{P} \longrightarrow \mathcal{P}'$ be an isogeny of
$O_F$-height $2$. 
Then there exists locally on $\Spec R$ an isogeny $\beta: \mathcal{P}'
\longrightarrow \mathcal{P}$ such that
\begin{displaymath}
   \beta \circ \alpha = \pi \id_{\mathcal{P}}, \quad \alpha \circ
   \beta = \pi \id_{\mathcal{P}'}.
\end{displaymath}

Let $\mathcal{P}$ and $\mathcal{P}'$ be the displays of formal
$p$-divisible groups $X$ and $X'$ with an $O_K$-action. Then the
kernel of any isogeny $\alpha: X \longrightarrow X'$ of height 2 is
annihilated by $\pi$. 

\end{proposition} 
    \begin{proof}
The proof is a variant of the proof of Proposition 1.6.4 in \cite{Zi3}.
We will use notation from that proof.       
The $O_K$-algebra structure on $R$ induces a natural homomorphism
$O_K \longrightarrow W_{O_F}(R)$ which is equivariant with respect
to the Frobenius $\tau \in \Gal(K/F)$ and the Frobenius on
$W_{O_F}(R)$. The composition with $\tau$ gives a second homomorphism
$O_K \longrightarrow W_{O_F}(R)$. We denote by $\Psi$ the set of these two
homomorphism.  We write $\bar{\psi} = \psi \circ \tau$ for $\psi \in \Psi$. 

The $O_K$-action gives the usual decompositions,
\begin{displaymath}
   P = \oplus_{\psi \in \Psi} P_{\psi}, \quad P' = \oplus_{\psi \in \Psi} P'_{\psi} .
\end{displaymath}
 We have the same kind of decompositions  for $Q \subset P$ and $Q'
 \subset P'$. We choose normal decompositions
\begin{displaymath}
   P_{\psi} = T_{\psi} \oplus L_{\psi}, \qquad P'_{\psi} = T'_{\psi} \oplus L'_{\psi} .
\end{displaymath}
The $T$ and $L$ on the right hand sides are by assumption locally free
of rank $1$. Using these decompositions, we write $\alpha_{\psi}:
P_{\bar{\psi}} \longrightarrow P'_{\psi}$ in matrix form
\begin{equation}\label{Isg13e}
M_{\psi} = \left(
\begin{array}{rr} 
X_{\psi} & ~^V Y_{\psi}\\
U_{\psi} & Z_{\psi} 
\end{array}
\right).
\end{equation}
The maps $\dot{F}_{\psi} : I_{O_F}(R) T_{\psi} \oplus L_{\psi}
\longrightarrow T_{\psi \tau} \oplus L_{\psi \tau}$, resp. $\dot{F}'_{\psi}
: I_{O_F}(R) T'_{\psi} \oplus L'_{\psi} \longrightarrow T'_{\psi \tau}
\oplus L'_{\psi \tau}$, are given by invertible matrices  
\begin{displaymath}
   \Phi_{\psi} = \left(
\begin{array}{cc}
A_{\psi} & B_{\psi}\\
C_{\psi} & D_{\psi}
\end{array}
\right), \qquad
   \Phi'_{\psi} = \left(
\begin{array}{cc}
A'_{\psi} & B'_{\psi}\\
C'_{\psi} & D'_{\psi}
\end{array}
\right).
\end{displaymath}
The matrices (\ref{Isg13e}) define a morphism $\alpha$ of displays iff
\begin{equation}\label{Isg14e}
   M_{\psi \tau} \Phi_{\psi} = \Phi'_{\psi}\; ~^s M_{\psi}, \quad
   \text{for} \; \psi \in \Psi.
\end{equation}
We will argue as in \cite{Zi3}. The meaning of the upper left index $~^s$ is the same as there. Taking  determinants we obtain
\begin{equation}\label{Isg15e}
   \det M_{\psi \tau} \; \det \Phi_{\psi} = \det \Phi'_{\psi}\; \det
   ~^s M_{\psi}.  
\end{equation}
In particular $\det M_{\psi}$ and
$~^{F^2}\det M_{\psi}$ differ by a unit in $W_{O_F}(R)$. As in
\cite{Zi3}  we obtain that 
\begin{displaymath}
   \det M_{\psi} = \pi^{\tt h} \epsilon_{\psi},
\end{displaymath} 
for some units $\epsilon_{\psi} \in W_{O_F}(R)$. By (\ref{Isg15e}), $\tt h$ is
independent of $\psi$.  Since $\alpha$ is an isogeny of height $2$ we
conclude that $\tt h = 1$. Now we pass to the adjucate matrices 
\begin{displaymath}
   ~^{\ad}\Phi_{\psi} \; ~^{\ad}M_{\psi \tau} = ~^{\ad}(~^s M_{\psi})\;
   ~^{\ad}\Phi'_{\psi}.
\end{displaymath}
Since the matrices $\Phi$ are invertible, we conclude
\begin{equation}\label{Isg16e}
   (\det \Phi_{\psi})\; ~^{\ad} M_{\psi \tau}\; \Phi'_{\psi} = (\det
   \Phi'_{\psi}) \Phi_{\psi} ~^{\ad}(~^s M_{\psi}).
\end{equation}
We consider first the case where $R$ is reduced. Then $\pi$ is not a
zero divisor in $W_{O_F}(R)$. In this case, we conclude from
(\ref{Isg15e}) 
\begin{displaymath}
   \det \Phi_{\psi}\; \epsilon_{\psi \tau} = \det \Phi'_{\psi} \;
   ~^F\epsilon_{\psi}. 
\end{displaymath} 
Thus we may rewrite equation (\ref{Isg16e}) as
\begin{displaymath}
   \epsilon_{\psi \tau}^{-1} ~^{\ad} M_{\psi \tau}\; \Phi'_{\psi} = 
   \Phi_{\psi}  ~^F\epsilon_{\psi}^{-1} ~^{\ad}(~^s M_{\psi}).
\end{displaymath}
This shows that the matrices $\epsilon_{\psi}^{-1} ~^{\ad}
M_{\psi}$ define the desired morphism $\beta: \mathcal{P}' \longrightarrow
\mathcal{P}$. In the case where $R$ is not reduced we may argue as in
the proof of Proposition 1.6.4 in \cite{Zi3}.
\end{proof}

    \section{Moduli spaces of formal local CM-triples}\label{s:modformal}

In this section we prove Theorems \ref{mainUR} and \ref{main-banal}. Recall that $d=[F:\BQ_p]$, and that we write $d=ef$.

\subsection{The case $r$ special and $K/F$ ramified} 

Let $(\mathbb{Y}, \iota_{\mathbb{Y}}, \lambda_{\mathbb{Y}})$ be a special formal
$O_D$-module over $\bar{\kappa}_F$ with a Drinfeld polarization
$\lambda_{\mathbb{Y}}$, cf. Definition \ref{Pol3d} and section \ref{ss:altram}.
We denote by $\breve{F}$ resp. $\breve{E}$,  the completions of the maximal
unramified extension of $F$, resp. of the reflex field $E$. Their residue class
fields are identified with $\bar{\kappa}_F$ resp. $\bar{\kappa}_E$.
We note that in the ramified case $E = E'$. We extend the embedding
$\varphi_0: O_F \longrightarrow O_E$ to an embedding
$\breve \varphi_0: O_{\breve{F}} \longrightarrow O_{\breve{E}}$. The base change
$(\mathbb{Y}, \iota_{\mathbb{Y}}, \lambda_{\mathbb{Y}})_{\bar{\kappa}_E}$ is
the base change via $\breve\varphi_0$.
 The relative display of this base change is an object of
 $\mathfrak{d}{\mathfrak R}^{\mathrm{nilp, pol}}_{\bar{\kappa}_E}$
if we restrict the action $\iota_{\mathbb{Y}}$ to $O_K$, cf. Definition
\ref{P'Pol3d} and Theorem \ref{Kontrakt1p}. By  Theorem \ref{P'Pol2c},
this relative display is the image of an object in
$\mathfrak{d}\mathfrak{P}^{\mathrm{ss, pol}}_{r,\bar{\kappa}_E}$ 
by the contracting functor. To the latter object corresponds a polarized
$p$-divisible formal group $(\mathbb{X}, \iota_{\mathbb{X}}, \lambda_{\mathbb{X}})$
with an action $\iota_{\mathbb{X}}: O_K \rightarrow \End \mathbb{X}$. 
The polarization $\lambda_{\mathbb{X}}$ is again principal. Using the
equivalence between formal $p$-divisible groups and nilpotent displays
one can say that
$(\mathbb{X}, \iota_{\mathbb{X}}, \lambda_{\mathbb{X}}) \in \mathfrak{P}^{\mathrm{pol}}_{r,\bar{\kappa}_E}$
is the object which is mapped by the contracting functor
$\mathfrak{C}^{\mathrm{pol}}_{r,\bar{\kappa}_E}$ to the special formal $O_D$-module
$(\mathbb{Y}, \iota_{\mathbb{Y}}, \lambda_{\mathbb{Y}})$. We note that we use
for the definition of
$(\mathbb{X}, \iota_{\mathbb{X}}, \lambda_{\mathbb{X}}) \in \mathfrak{P}^{\mathrm{pol}}_{r,\bar{\kappa}_E}$
only classical Dieudonn\'e theory over a perfect field.  

We consider the functor $\CM_r = \mathcal{M}_{K/F, r}$
\index[NO]{MBA@$\mathcal{M}_{K/F, r}$} of Definition \ref{CRFd2} in the
ramified case (where ${\tt h} = 0$). By Proposition \ref{KottwitzC6p}, this
Definition may be reformulated as follows.
\begin{proposition}\label{FuncMr1p}
Let $S$ be a scheme over $\Spf O_{\breve{E}}$. A point of
  $\mathcal{M}_{r}(S)$ consists of
  \begin{altenumerate}
  \item[(1) ]\, a local CM-pair $(X, \iota)$ of CM-type $r$ over $S$
    which satisfies the Eisenstein conditions
    $({\rm EC}_r)$ relative to a fixed uniformizer $\Pi$ of $K$ and such
    that
    \begin{equation}\label{FuncMr1e}
\Trace(\iota(\Pi) \; \mid \; \mathbf{E}_{A_{\psi_0}} \Lie_{\psi_0} X) = 0.  
      \end{equation}
    \item[(2) ] an isomorphism of $p$-divisible $O_K$-modules
    \begin{displaymath}
      \lambda: X \isoarrow X^{\wedge}, 
    \end{displaymath}
    which is a polarization of $X$. 
     \item[(3) ] a quasi-isogeny of height zero
    of $p$-divisible $O_K$-modules
    \begin{displaymath}
      \rho: X \times_{S}  \bar{S} \longrightarrow \mathbb{X}
      \times_{\Spec \bar{\kappa}_E}  \bar{S}, 
    \end{displaymath}
    such that the pullback quasi-isogeny $\rho^*(\lambda_\BX)$ differs
    from $\lambda_{|X \times_S \ov{S}}$ by a scalar in $O_F^{\times}$, locally
    on $\ov S$.  Here
    $\ov{S} = S \times_{\Spec O_{\breve E}} \Spec \bar{\kappa}_E$.  
  \end{altenumerate}
  Two data $(X_1, \iota_1, \lambda_1, \rho_1)$ and
  $(X_2, \iota_2, \lambda_2, \rho_2)$
define the same point iff there is an isomorphism of $O_K$-modules
$\alpha: X_1 \rightarrow X_2$ such that
$\rho_2 \circ \alpha_{\bar{S}} = \rho_1$.  (This implies that $\alpha$ respects
the polarizations up to a factor in $O_F^{\times}$).  
  \end{proposition}
In Definition \ref{CRFd2} we required that the scalar is in $F^{\times}$
but because the polarizations $\lambda$ and $\lambda_{\mathbb{X}}$ are
principal and $\rho$ has height $0$ the scalar is automatically in
$O_F^{\times}$. We remark that the condition (\ref{FuncMr1e}) depends only on the restriction
of the structure morphism $S \longrightarrow \Spf O_E$.  

We define for $i \in \mathbb{Z}$ the functor $\mathcal{M}_r(i)$ on the
category of schemes $S \longrightarrow \Spf O_{\breve{E}}$  by replacing $(3)$ in
Proposition \ref{FuncMr1p} by
\begin{enumerate}
\item[($3^\prime$)] {\it  a quasi-isogeny of $p$-divisible $O_K$-modules}
    \begin{displaymath}
      \rho: X \times_{S}  \bar{S} \longrightarrow \mathbb{X}
      \times_{\Spec \bar{\kappa}_E}  \bar{S}, 
    \end{displaymath}
        {\it  such that the pullback quasi-isogeny $\rho^*(\lambda_{\mathbb{X}})$
          differs from $p^i\lambda_{|X \times_S \bar{S}}$ by a scalar in
          $O_F^{\times}$, locally on $\bar{S}$}.
  \end{enumerate}
 It follows from the last condition that
  \begin{equation}\label{iM1e}
2\height \rho = \height (p^i \; | \; X) = 4di.
  \end{equation} 
  We have an isomorphism of functors
  \begin{displaymath}
\mathcal{M}_r \longrightarrow \mathcal{M}_r(i) ,
    \end{displaymath}
  which associates to a point $(X, \iota, \lambda, \rho) \in \mathcal{M}_r(S)$ 
  the point $(X, \iota, \lambda, \Pi^{ei}\rho) \in \mathcal{M}_r(i)(S)$. 
  We set \index[NO]{MBF@$\tilde{\mathcal{M}}_r$}
  \begin{displaymath}
 \tilde{\mathcal{M}}_r = \coprod_{i \in \mathbb{Z}} \mathcal{M}_r(i).
  \end{displaymath}
    We define a Weil descent datum on $\tilde{\mathcal{M}}_r$ relative to
  $O_{\breve{E}}/O_E$. Let $\tau_E$ be the Frobenius of $\breve E/E$.
  It is enough to consider the functor $\tilde{\mathcal{M}}_r$
  for affine schemes $S = \Spec R$. We write
  $\varepsilon: O_{\breve{E}} \longrightarrow R$ for the given algebra structure.
  We write  $R_{[\tau_E]}$  for the ring $R$ with the new algebra structure
  $\varepsilon \circ\tau_E$.  
  By base change to $\bar{\kappa}_{E}$, we obtain
  \begin{displaymath}
    \bar{\varepsilon}: \bar{\kappa}_{E} \longrightarrow \bar{R} :=
    R \otimes_{O_{\breve{E}}} \bar{\kappa}_{E}. 
    \end{displaymath}
  We consider a point $(X, \iota, \lambda, \rho) \in \mathcal{M}_r(i)(R)$,
  where $\rho$ is a quasi-isogeny
  \begin{displaymath}
\rho: X_{\bar{R}} \longrightarrow \bar{\varepsilon}_{\ast} \mathbb{X}. 
    \end{displaymath}
  Since the notion of a CM-triple depends only on the induced
  $O_{E}$-algebra structure on $R$, we may regard $(X, \iota, \lambda)$
  as a CM-triple on $R_{[\tau_E]}$. We set
  \begin{displaymath}
    \tilde{\rho}: X_{\bar{R}} \overset{\rho}{\longrightarrow}
    \bar{\varepsilon}_{\ast} \mathbb{X} 
    \overset{F_{\mathbb{X},\tau_E}}{\longrightarrow}
    \bar{\varepsilon}_{\ast} (\tau_E)_{\ast} \mathbb{X}. 
  \end{displaymath}
  The assignment $(X, \iota, \lambda, \rho) \mapsto (X, \iota, \lambda, \tilde{\rho})$
  defines a morphism\index[NO]{ZZWG@$\omega_{\mathcal{M}_r}$}
  \begin{equation}\label{WD5e} 
    \omega_{\mathcal{M}_r}: \mathcal{M}_r(i)(R) \longrightarrow
    \mathcal{M}_r(i + f_E)(R_{[\tau_E]}) 
    \end{equation}
  where $f_E = [\kappa_E : \mathbb{F}_p]$. Here we note that the inverse image
  of the polarization $(\tau_E)_{\ast} \lambda_{\mathbb{X}}$ on
  $(\tau_E)_{\ast} \mathbb{X}$ by
  $F_{\mathbb{X},\tau_E}: \mathbb{X} \longrightarrow (\tau_E)_{\ast} \mathbb{X}$
  is $p^{f_E}\lambda_{\mathbb{X}}$.
  From (\ref{WD5e}) we obtain the Weil  descent datum
  \begin{equation}\label{WD7e}
    \omega_{\mathcal{M}_r}: \tilde{\mathcal{M}}_r \longrightarrow
    \tilde{\mathcal{M}}_r^{(\tau_E)} ,
  \end{equation}
  where the upper index $(\tau_E)$ denotes the base change via
  $\Spec \tau_E : \Spf O_{\breve{E}} \longrightarrow \Spf O_{\breve{E}}$.

Let $\mathcal{N}(i)$ be the functor of Definition \ref{Ni1d}. Note that we
took $\BY$ for the framing object appearing in the definition of
$\mathcal{N}(i)$.  We consider a point 
$(X, \iota,\lambda, \rho) \in \mathcal{M}_{r}(i)(R)$, where
$R\in \Nilp_{O_{\breve E}}$.
 Applying the contracting functor $\mathfrak{C}^{\rm pol}_{r, R}$
 of Theorem \ref{P'Pol2c} to its display, 
we obtain a quadruple $(X_{\rm{c}},\iota_{\rm{c}},\lambda_{X_{\rm{c}}}, \rho_{\rm{c}})$.
It follows from the isomorphism (\ref{KottwitzC20e}) (which also holds in the
ramified case, cf. a few lines below \eqref{KottwitzC20e}) that the condition (\ref{FuncMr1e}) implies 
\begin{displaymath}
\Trace(\iota_{\rm{c}}(\Pi) \; \mid \; \Lie X_{\rm{c}}) = 0.  
  \end{displaymath}
By functoriality, the polarizations $\rho_{\rm{c}}^{\ast} (\lambda_{\mathbb{Y}})$ and
$p^{i} \lambda_{(X_{\rm{c}})_{\bar{R}}}$ differ by a unit in $O_F$. Hence 
$(X_{\rm{c}}, \iota_{\rm{c}},\lambda_{X_{\rm{c}}}, \rho_{\rm{c}})$ defines a point of $\mathcal{N}(i)$.
Therefore we obtain  from Theorem \ref{P'Pol2c} an isomorphism of
functors, 
\begin{equation}\label{CRFram11e}
    \mathcal{M}_{r}(i) \isoarrow \mathcal{N}(ei)\times_{\Spf O_{\breve F}}\Spf{O_{\breve{E}}}.  
\end{equation}
The base change on the right hand side is via
$\breve\varphi_0: O_{\breve{F}} \longrightarrow O_{\breve{E}}$

We set \index[NO]{NBD@$\tilde{\mathcal{N}}[e]$}
\begin{equation}\label{N[e]}
\tilde\CN[e]=\coprod_{i\in\BZ}\CN(ei) .
\end{equation}
We endow $\tilde\CN[e]$ with a Weil descent datum relative to $O_{\breve E}/O_E$. Let $R\in\Nilp_{O_{\breve E}}$. We consider the map
\begin{equation*}
\Pi^{(d-1)f_E/f }\omega_\CN^{f_E/f}\colon \CN(ei)(R)\to\CN(e(i+f_E)(R_{[\tau_E]}) .
\end{equation*}
Here on the left hand side appears the iterate of the Weil descent datum   $\omega_{\mathcal{N}}: \mathcal{N}(i) \longrightarrow
      \mathcal{N}(i+1)^{({\tau})}$  of $\tilde\CN$ relative to $O_{\breve F}/O_F$ from \eqref{WD4e} and the translation functor $ \Pi: \mathcal{N}(i) \longrightarrow \mathcal{N}(i+1)$, cf. \eqref{FuncN1e}. This defines a Weil descent datum relative to $O_{\breve E}/O_E$,
      \begin{equation}\label{weildescN[e]}
 \Pi^{(d-1)f_E/f }\omega_\CN^{f_E/f}\colon \tilde\CN[e]\to\tilde\CN[e]^{(\tau_E)} .
      \end{equation}

  We define \index[NO]{JAG@$J'$}
  \begin{equation}\label{J'ram1e}
    \begin{array}{rl} 
    J' = & \{\alpha \in \Aut^{o}_K \mathbb{X} \; | \;
    \alpha^{\ast}( \lambda_{\mathbb{X}}) = u \lambda_{\mathbb{X}}, \; \text{for}\;
    u \in p^{\mathbb{Z}} O_F^{\times} 
    \}\\[2mm] = & 
    \{\alpha \in \Aut^{o}_K \mathbb{Y} \; | \;
    \alpha^{\ast} (\lambda_{\mathbb{Y}}) = u \lambda_{\mathbb{Y}}, \; \text{for}\;
    u \in p^{\mathbb{Z}} O_F^{\times} 
    \} .
      \end{array}
  \end{equation}
  The last equation holds because of the contraction functor. 
  This group acts via the rigidifications $\rho$ on the functors 
  $\tilde{\mathcal{M}}_{r}$ and $\tilde{\mathcal{N}}[e]$. By the last
  equation of (\ref{J'ram1e}), we may regard $J'$ as a subgroup
  of $J^{\ast {\rm r}}$ of Lemma \ref{J^v1l}. 
  
 \begin{proposition}\label{CRFram5p}  
 There is an isomorphism of formal schemes over $\Spf O_{\breve{E}}$ 
  \begin{displaymath}
    \tilde{\mathcal{M}}_{r} \longrightarrow
    \tilde{\mathcal{N}}[e]\times_{\Spf O_{\breve F}}\Spf{O_{\breve{E}}} ,  
  \end{displaymath}
  where the right hand side denotes the base change via 
  $\breve\varphi_0: O_{\breve{F}} \longrightarrow O_{\breve{E}}$. This isomorphism
  is compatible with the action of $J'$ on both sides. 
  
  Under the isomorphism the Weil descent datum \eqref{WD7e} on the left hand
  side corresponds to the Weil descent datum 
   $$
     \Pi^{(d-1)f_E/f} \omega_{\mathcal{N}}^{f_E/f}: \tilde{\mathcal{N}}[e]
  \longrightarrow \tilde{\mathcal{N}}[e]^{(\tau_E)} 
   $$
   on the right hand side.  More explicitly,  for any $i\in\BZ$, there is a commutative diagram 
\begin{equation}\label{WD6e}
\begin{aligned}
  \xymatrix{
   {\mathcal{M}}_r(i) \ar[r]^-\sim \ar[d]_{\omega_{\mathcal{M}_r}} &
    {\mathcal{N}}(ei)_{O_{\breve{E}}} \ar[d]^{\Pi^{(d-1)f_E/f} \omega_{\mathcal{N}}^{f_E/f}}\\
    {\mathcal{M}}_r(i+f_E)^{(\tau_E)}  \ar[r]^-\sim &
    {\mathcal{N}}(e(i+f_E))_{O_{\breve{E}}}^{(\tau_E)}.\\
  }
  \end{aligned}
  \end{equation} 
\end{proposition}
\begin{proof} The isomorphism of formal schemes over $\Spf O_{\breve E}$ comes from \eqref{CRFram11e}. 
  It remains to show that the diagram is commutative. Let
  $R \in \Nilp_{O_{\breve{E}}}$ with structure morphism
  $\varepsilon: O_{\breve{E}} \longrightarrow R$. Let
  $\bar{\varepsilon}: \bar{\kappa}_E \longrightarrow \bar{R} = R \otimes_{O_{\breve{E}}}  \bar{\kappa}_E$
  be the induced morphism.  We start with a point
  $(X, \iota, \lambda, \rho) \in \mathcal{M}_r(i)(R)$.
  If we apply $\omega_{\mathcal{M}_r}$, we change $\rho$ to
  \begin{displaymath}
   \tilde{\rho}: X_{\bar{R}} \longrightarrow \bar{\varepsilon}_{\ast} \mathbb{X}
   \overset{F_{\mathbb{X}, \tau_E}}{\longrightarrow} 
   \bar{\varepsilon}_{\ast} (\tau_E)_{\ast} \mathbb{X} .
  \end{displaymath}
  
  The lower horizontal arrow in (\ref{WD6e}) applies to $\tilde{\rho}$ the
  contracting  functor $\mathfrak{C}_{r, R}$, cf. Definition \ref{Kontrakt1e}. We have $\mathbb{X}_{\rm{c}} = \mathbb{Y}$. Let
  $P_{\mathbb{X}}$ be the $\mathcal{W}(\kappa_{\bar{E}})$-Dieudonn\'e module of
  $\mathbb{X}$. Then the $\mathcal{W}_{O_F}(\kappa_{\bar{E}})$-Dieudonn\'e module of
  $\mathbb{Y}$ is by definition the degree-zero component of $P'_{\mathbb{X}}$ defined by (\ref{defP'1e}), comp. Remark \ref{CRF1r}.   From this definition we obtain
  \begin{displaymath}
(V')^f = \Pi^{1-d}V^f.
    \end{displaymath}
  In terms of Dieudonn\'e modules,  $F_{\mathbb{X}, \tau_E}$ is given by
  \begin{displaymath}
    P_{\mathbb{X}} \longrightarrow W(\kappa_{\bar{E}}) \otimes_{W(\tau_E), W(\kappa_{\bar{E}})}
    P_{\mathbb{X}}, \quad x \mapsto 1 \otimes V^{f_E} x. 
  \end{displaymath}
  In terms of the relative Dieudonn\'e module, $F_{\mathbb{Y}, \tau_E}$ is given by $(V')^{f_E}$. Therefore the contracting 
  functor applied to $F_{\mathbb{X}, \tau_E}$ gives 
  \begin{displaymath}
    \iota_{\mathbb{Y}}(\Pi^{(d-1)f_E/f}) F_{\mathbb{Y}, \tau_E}: \mathbb{Y}
    \longrightarrow (\tau_E)_{\ast} \mathbb{Y}, 
  \end{displaymath}
  since $(V')^{f_E} = \Pi^{(1-d)f_E/f} V^{f_E}$. 
  On the other hand $\omega_{\mathcal{N}}$ just multiplies $\rho_{\rm{c}}$ by
  $F_{\mathbb{Y}, \tau_E}$. Therefore we obtain the commutativity of the diagram. 
  \end{proof}

\begin{corollary}\label{CRFram1c}
Let $\omega_{\tau_E}$ denote the action of $\tau_E$  on the formal scheme 
scheme  $\wh{\Omega}_{F} \times_{\Spf O_F} \Spf O_{\breve{E}}$ via the second factor. 
   There exists an isomorphism of formal schemes over $\Spf O_{\breve E}$
  $$
       \tilde{\mathcal{M}}_{r} \isoarrow   (\wh{\Omega}_{F} \times_{\Spf O_F} \Spf O_{\breve{E}})  \times \mathbb{Z} ,
$$
such that the Weil  descent datum $\omega_{\mathcal{M}_r}$ induces on the right
hand side the Weil  descent datum
\begin{displaymath}
   (\xi, i) \mapsto (\omega_{\tau_E}(\xi), i+f_E).
\end{displaymath}
In particular the formal scheme $\mathcal{M}_{K/F, r}$ over $ \Spf O_{\breve{E}}$
is a $p$-adic formal scheme which has semi-stable reduction; hence it is also flat, with reduced special fiber. 
\end{corollary}
\begin{proof}
  We consider the isomorphism 
$\tilde{\mathcal{M}}_{\mathrm{Dr}} \longrightarrow \tilde{\mathcal{N}}$ 
  of formal schemes over $\Spf O_{\breve{F}}$ from Theorem \ref{altDrinf}. It is compatible with 
  translations and the Weil descent data on both sides. 
 Combining this with the Drinfeld isomorphism  \eqref{FuncMDr2e},  we obtain an isomorphism 
  \begin{equation}\label{CRFram3e}
  (\wh{\Omega}_{F} \times_{\Spf O_F} \Spf O_{\breve F} ) \times \mathbb{Z}e  
\overset{\sim}{\longrightarrow} \tilde{\mathcal{N}}[e] .
  \end{equation}
  We consider on the right hand side the Weil  descent datum
  $\Pi^{(d-1)f_E/f} \omega_{\mathcal{N}}^{f_E/f}$ which is a composite of
  an iterate of the translation functor and an iterate of the Weil
   descent datum $\omega_{\mathcal{N}}$. By (\ref{WD13e}) we see that
  $\Pi^{(d-1)f_E/f} \omega_{\mathcal{N}}^{f_E/f}$ induces on the left hand side
  of (\ref{CRFram3e}) the Weil  descent datum
  \begin{displaymath}
(\xi, ei) \mapsto (\omega_{\tau_E}, e(i+f_E)).
  \end{displaymath}
  The assertion about descent data follows by forgetting $e$. The last assertion follows from Corollary \ref{spfibred}.
  \end{proof}

\subsection{The case $r$ special and $K/F$ unramified} 

Let $\varphi_0, \bar{\varphi}_0 \in \Phi$ be the special embeddings. Their
restrictions to $F$ are the same. We extend the resulting embedding
$O_F \longrightarrow O_E$ to an embedding
$O_{\breve{F}} \longrightarrow O_{\breve{E}}$. The two embeddings $\varphi_0, \bar{\varphi}_0$ then factor over the two
$O_F$-algebra
homomorphisms $\breve\varphi_0, \bar{\breve\varphi}_0: O_K \longrightarrow O_{\breve{F}}$.  

Let $(\mathbb{Y}, \iota_{\mathbb{Y}})$ be a special formal $O_D$-module
over $\bar{\kappa}_F$. We endow the $\mathcal{W}_{O_F}(\bar{\kappa}_F)$-Dieudonn\'e module $P_{\mathbb{Y}}$ of $\mathbb{Y}$ with the polarization
$\theta$, cf. (\ref{Pol21e}). (One should note that we call now $K$ what was
$F'$ in that section.) Then $\theta$ defines an isogeny to the Faltings dual,
\begin{displaymath}
\lambda_{\mathbb{Y}}: \mathbb{Y} \longrightarrow \mathbb{Y}^{\nabla}.  
\end{displaymath}
If we compose the action $\iota_{\mathbb{Y}}^{\nabla}$ of $(O_D)^{\rm opp}$ on 
$\mathbb{Y}^{\nabla}$ with the anti-involution $\dagger$ of (\ref{Pol31e}),
the isogeny  $\lambda_{\mathbb{Y}}$ becomes an isogeny of special formal $O_D$-modules.
We indicate this by rewriting the isogeny as
\begin{equation}
\lambda_{\mathbb{Y}}: \mathbb{Y} \longrightarrow \mathbb{Y}^{\Delta}.   
  \end{equation}

In section \ref{ss:altunram} we used $(\mathcal{P}_{\mathbb{Y}}, \theta)$ to
define the functor $\mathcal{N}(i)$. Together with the restriction to $O_K$ of
the action of $O_D$ on $\mathcal{P}_{\mathbb{Y}}$, we obtain an object of
$\mathfrak{d}{\mathfrak R}^{\rm pol}_{\bar{\kappa}_F}$, cf. Definition \ref{P'Pol3d}.
By Theorem \ref{P'Pol2c}, this object is the image of an object in
$\mathfrak{d}\mathfrak{P}^{\rm ss, pol}_{r,\bar{\kappa}_E}$. The latter is the
display of an object
$(\mathbb{X}, \iota_{\mathbb{X}}, \lambda_{\mathbb{X}}) \in
\mathfrak{P}^{\rm pol}_{r, \bar{\kappa}_E}$, cf. Definition \ref{KatCMtriple1d}.
The height of the polarization $\lambda_{\mathbb{X}}$ is $2f$, and the associated Rosati involution induces on $O_K$
the conjugation over $O_F$.  

We consider the functor $\mathcal{M}_r=\CM_{K/F, r}$ of Definition \ref{CRFd2}
with the framing object $(\mathbb{X}, \iota_{\mathbb{X}}, \lambda_{\mathbb{X}})$ as
defined above. We can give an alternative description of that functor. 

\begin{proposition}\label{FuncMr2p}
 Let
  $S$ be a scheme over  $ \Spf O_{\breve{E}}$.  A point of
  $\mathcal{M}_{r}(S)$ consists of
  \begin{enumerate}
  \item[(1) ] a local CM-pair $(X, \iota)$ of CM-type $r$ over $S$
    which satisfies the Eisenstein conditions $({\rm EC}_r)$ relative to
    the fixed uniformizer $\pi$ of $F$. 
    \item[(2) ] an isogeny of height $2f$ of $p$-divisible $O_K$-modules
    \begin{displaymath}
      \lambda: X \longrightarrow  X^{\wedge}, 
    \end{displaymath}
    which is a polarization of $X$. 
     \item[(3) ] a quasi-isogeny of $p$-divisible $O_K$-modules
    \begin{displaymath}
      \rho: X \times_{\Spec R}  \Spec \bar{R} \longrightarrow \mathbb{X}
      \times_{\Spec \bar{\kappa}_E}  \Spec \bar{R}, 
    \end{displaymath}
    such that the pullback quasi-isogeny $\rho^*(\lambda_{\mathbb{X}})$ differs
    from $\lambda_{|X \times_{S} \Spec \bar{S}}$ by a scalar in $O_F^{\times}$,
    locally on $\bar{S}$, where
    $\bar{S} = S \otimes_{\Spf O_{\breve E}} \Spec \bar{\kappa}_E$. 
  \end{enumerate}
  Two data $(X_1, \iota_1, \lambda_1, \rho_1)$ and
  $(X_2, \iota_2, \lambda_2, \rho_2)$
define the same point iff there is an isomorphism of $O_K$-modules
$\alpha: X_1 \rightarrow X_2$ such that
$\rho_2 \circ \alpha_{\bar{S}} = \rho_1$.  (This implies that $\alpha$ respects
the polarizations up to a factor in $O_F^{\times}$).  
\end{proposition}
\begin{proof}
  We may assume that $S = \Spec R$. 
  Let $(X_0, X_1, \lambda, \rho_X) \in \mathcal{M}_r(R)$ be a point as in
  Definition \ref{CRFd2}. We obtain a point $(X, \lambda, \rho_X)$ as in the
  Proposition if we set
  $X := X_0$, keep $\rho_X$ and redefine $\lambda$ to be the composite
  \begin{displaymath}
X_0 \longrightarrow X_1 \overset{\lambda}{\longrightarrow} X_0^{\wedge}.   
  \end{displaymath}
  We note that $\rho$ is automatically of height zero because by the last
  condition $\rho^*(\lambda_{\mathbb{X}})$ and $\lambda_{\mathbb{X}}$ have
  the same height $2f$.
  
  Conversely, assume $(X, \lambda, \rho_X)$ is as in the Proposition.
  Then we set $X_0 = X$ and $X_1 = X^{\wedge}$. By Corollary \ref{dualEisen2c},
  $X_1$ satisfies the condition $({\rm EC}_r)$ and, by Proposition
  \ref{KottwitzC5p}, $X_0$ and $X_1$ satisfy $({\rm KC}_r)$ . 
  The polarization $\lambda$ defines
  an isogeny $a: X_0 \longrightarrow X_1$ of $p$-divisible $O_K$-modules which
  has  height $2f$. To the morphism induced by $a$ on the displays  
  we apply  the contracting functor of Definition \ref{Kontrakt1e}. We obtain an isogeny
  of $\mathcal{W}_{O_F}(R)$-displays
  $\alpha: \mathcal{P}_0 \longrightarrow \mathcal{P}_1$ of height $2$. To this
  isogeny we may apply Proposition \ref{SFB701p}. We find an isogeny
  $\beta: \mathcal{P}_1 \longrightarrow \mathcal{P}_0$ such that
  $\beta \circ \alpha = \pi \id_{\mathcal{P}_0}$. 

The existence of $\rho$ guarantees that $X_0$, $X_1$ and $a$ are defined over
an $O_{\breve{E}}$-subalgebra $R' \subset R$ which
is of finite type over $O_{\breve{E}}$. Therefore we may apply Theorem 
\ref{Kontrakt1p}. It gives us a morphism $b: X_1 \longrightarrow X_0$ such that
$b \circ a = \pi \id_{X_0}$. We see that $X_0, X_1, \rho$ together with the
defining isomorphism $X_1 \overset{\sim}{\longrightarrow} X_0$ defines a point
of the functor $\mathcal{M}_r$ of Definition \ref{CRFd2}.
\end{proof}

\begin{definition}\label{iMr2d} 
  We define the functor $\mathcal{M}_r(i)$ in the same way as in Definition
  \ref{FuncMr2p}, but we replace the data $(3)$ by the following
  \begin{enumerate}
  \item[($3^\prime$)] A quasi-isogeny of $p$-divisible $O_K$-modules
    \begin{displaymath}
      \rho: X \times_{S} \bar{S} \longrightarrow \mathbb{X}
      \times_{\Spec \bar{\kappa}_E} \bar{S}, 
    \end{displaymath}
    such that the pullback quasi-isogeny $\rho^*(\lambda_{\mathbb{X}})$ differs
    from $p^i\lambda_{|X \times_{S} \bar{S}}$ by a scalar in
    $O_F^{\times}$, locally on $\bar{S}$. 
    \end{enumerate}
  \end{definition}
As in the ramified case (\ref{iM1e}) we conclude that $\height \rho = 2di$. 
We set\index[NO]{MBG@$\tilde{\mathcal{M}}_r$}
  \begin{displaymath}
 \tilde{\mathcal{M}}_r = \coprod_{i \in \mathbb{Z}} \mathcal{M}_r(i).
  \end{displaymath}
  Let $R \in \Nilp_{O_{\breve{E}}}$.  
  Exactly as in the ramified case we obtain a morphism 
  \begin{equation}\label{WDu5e} 
    \omega_{\mathcal{M}_r}: \mathcal{M}_r(i)(R) \longrightarrow
    \mathcal{M}_r(i + f_E)(R_{[\tau_E]}) ,
    \end{equation}
  where, as in the ramified case, $f_E = [\kappa_E : \mathbb{F}_p]$.
  With the notation used in the ramified case, it changes the
  datum $\rho$ in point ($4^\prime$) of Definition \ref{iMr2d} to
   \begin{displaymath}
    \tilde\rho: (X_0)_{\bar{R}} \overset{\rho}{\longrightarrow}
    \bar{\varepsilon}_{\ast} \mathbb{X} 
    \overset{F_{\mathbb{X},\tau_E}}{\longrightarrow}
    \bar{\varepsilon}_{\ast} (\tau_E)_{\ast} \mathbb{X}. 
  \end{displaymath}
   From (\ref{WDu5e}) we obtain the Weil  descent datum,
  \begin{equation}\label{WDu7e}
    \omega_{\mathcal{M}_r}: \tilde{\mathcal{M}}_r \longrightarrow
    \tilde{\mathcal{M}}_r^{(\tau_E)} .
  \end{equation}
We define an isomorphism of functors on $\Nilp_{O_{\breve E}}$,
  \begin{equation}\label{MNmorph1e}
\mathcal{M}_r(i) \isoarrow \mathcal{N}(ei). 
    \end{equation}
  Let $(X, \iota, \lambda, \rho) \in \mathcal{M}_r(i)(R)$. Applying the
contracting functor $\mathfrak{C}^{\rm pol}_{r, R}$, we obtain a quadruple $(X_{\rm{c}}, \iota_{\rm{c}},\lambda_{X_{\rm{c}}}, \rho_{\rm{c}})$, where
  $\rho_{\rm{c}}: X_{{\rm{c}}}\otimes_R,\bar{R} \longrightarrow \mathbb{Y} \otimes_{\bar{\kappa}_E} \bar{R}$.
This gives a point of
$\mathcal{N}(i)(R)$. The functor $\mathfrak{C}^{\rm pol}_{r, R}$ is
  an equivalence of categories if the ideal of nilpotent elements of $R$ is
  nilpotent, cf. Theorem~\ref{P'Pol2c}.Therefore we may reverse the
  construction of (\ref{MNmorph1e}).
  Therefore $\mathcal{M}_r(i)(R) \longrightarrow \mathcal{N}(ei)(R)$ is
  bijective for those $R$. For a general $R$ we obtain the bijectivity as
  in the proof of Proposition \ref{FuncMr2p}. With the notation (\ref{N[e]})
  we have a bijection $\tilde{\CM}_{r}(R) \isoarrow \tilde{\CN}[e](R)$.  
  
  We define\index[NO]{JAG@$J'$} 
  \begin{equation}\label{J'unram1e}
    \begin{array}{rl} 
    J' = & \{\alpha \in \Aut^{o}_K \mathbb{X} \; | \;
    \alpha^{\ast} (\lambda_{\mathbb{X}}) = u \lambda_{\mathbb{X}}, \; \text{for}\;
    u \in p^{\mathbb{Z}} O_F^{\times} 
    \}\\[2mm] = & 
    \{\alpha \in \Aut^{o}_K \mathbb{Y} \; | \;
    \alpha^{\ast} (\lambda_{\mathbb{Y}}) = u \lambda_{\mathbb{Y}}, \; \text{for}\;
    u \in p^{\mathbb{Z}} O_F^{\times} 
    \} .
      \end{array}
  \end{equation}
  The last equation holds because of the contraction functor. 
  This group acts via the rigidifications $\rho$ on the functors 
  $\tilde{\mathcal{M}}_{r}$ and $\tilde{\mathcal{N}}[e]$. By the last
  equation of (\ref{J'unram1e}) we may regard $J'$ as a subgroup
  of $J^{\ast{\rm ur}}$ of Lemma \ref{J^uv1l}.

  \begin{proposition}\label{MNmorph2p}
    There exists an isomorphism of formal schemes
    $\Spf O_{\breve{E}}$ 
    $$
    \tilde{\CM}_{r}\isoarrow \tilde{\CN}[e] \times_{\Spf O_{\breve{F}}}
    \Spf O_{\breve{E}} , 
    $$
  where the right hand side denotes the base change via 
  $\breve\varphi_0: O_{\breve{F}} \longrightarrow O_{\breve{E}}$. This isomorphism
  is compatible with the action of $J'$ on both sides, and the  
  Weil  descent datum \eqref{WD7e} on the left hand side corresponds to the
  Weil descent datum 
  \begin{displaymath}
    \pi^g\omega_{\mathcal{N}}^{f_E/f}: \tilde{\mathcal{N}}[e]
  \longrightarrow \tilde{\mathcal{N}}[e]^{(\tau_E)}.  
    \end{displaymath}
   on the right hand side.  Here 
    $g = f_E (d-1)/2f$. More explicitly,  there is a commutative diagram
\begin{equation}\label{WD6ee}
\begin{aligned}
  \xymatrix{
    {\mathcal{M}}_r(i)  \ar[r]^-\sim \ar[d]_{\omega_{\mathcal{M}_r}} &
    {\mathcal{N}}_{O_{\breve{E}}}(ei) \ar[d]^{\pi^g\omega_{\mathcal{N}}^{f_E/f}}\\
    {\mathcal{M}}_r(i+f_E)^{(\tau_E)}   \ar[r]^-\sim &
    {\mathcal{N}}_{O_{\breve{E}}}(e(i+f_E))^{(\tau_E)}.\\
  }
\end{aligned}
  \end{equation} 
The multiplication by $\pi$ is the morphism 
$\mathcal{N}(j) \rightarrow \mathcal{N}(j + 2)$ which is obtained by
multiplying $\rho$ in Definition \ref{N-unram1d} by $\pi$. Equivalently one
can apply two times the translation functor (\ref{translation3e}). 
  \end{proposition} 
  \begin{proof}
    We have already proved the isomorphism of formal schemes over  $\Spf O_{\breve{E}}$. The compatibility
    with the Weil descent datum follows from the following lemma. 
    \end{proof}

  \begin{lemma}\label{Frobstrich1l}
    Let $\mathbb{X}$ be the framing object over $\bar{\kappa}_E$, with corresponding Frobenius morphism $F_{\mathbb{X}, \kappa_E}: \mathbb{X} \longrightarrow (\tau_{E})_{\ast} \mathbb{X}$. Let $\BY$ be its image under the contracting functor $\mathfrak{C}_{r, \bar{\kappa}_E}^{\rm pol}$. Then the image of $F_{\mathbb{X}, \kappa_E}$ under $\mathfrak{C}_{r, \bar{\kappa}_E}^{\rm pol}$ is given as 
    \begin{displaymath}
      \pi^g F_{\mathbb{Y}, \kappa_E}: \mathbb{Y}\longrightarrow
      (\tau_{E})_{\ast} \mathbb{Y}. 
    \end{displaymath}
    \end{lemma}
  \begin{proof}
   Let
  $M=P_{\mathbb{X}}$ be the $\mathcal{W}(\bar{\kappa}_{E})$-Dieudonn\'e module of
  $\mathbb{X}$. Then the $\mathcal{W}_{O_F}(\kappa_{\bar{E}})$-Dieudonn\'e module $M'$ of
  $\mathbb{Y}$ is by definition the $\psi_0$-component of $P'_{\mathbb{X}}$ defined by (\ref{moreeq}).
  In terms of Dieudonn\'e modules,  $F_{\mathbb{X}, \tau_E}$ is given by
  \begin{displaymath}
  M \longrightarrow W(\kappa_{\bar{E}}) \otimes_{\tau_E, W(\kappa_{\bar{E}})}
    M, \quad x \mapsto 1 \otimes V^{f_E} x. 
  \end{displaymath}
 This induces
    \begin{equation}\label{Frobstrich1e}
      V^{f_E}: M' \longrightarrow W(\bar{\kappa}_E) \otimes_{\tau_E, W(\bar{\kappa}_E)} M'.
    \end{equation}
    We consider the decomposition
    \begin{displaymath}
M = \oplus M_{\psi}.
      \end{displaymath}
    We denote by $\sigma$ the Frobenius on the Witt ring $W(\bar{\kappa}_E)$
    and also the Frobenius of $K^t/\mathbb{Q}_p$. We note that
    \begin{displaymath}
      W(\bar{\kappa}_E) \otimes_{\tau_E, W(\bar{\kappa}_E)} M_{\psi} =
      (W(\bar{\kappa}_E) \otimes_{\tau_E, W(\bar{\kappa}_E)} M)_{\tau_E \psi} =
      (W(\bar{\kappa}_E) \otimes_{\tau_E, W(\bar{\kappa}_E)} M)_{\psi \sigma^{f_E}}. 
      \end{displaymath}
     In terms of the relative Dieudonn\'e module, $F_{\mathbb{Y}, \tau_E}$ is given by $(V')^{f_E}$. Our problem is to express (\ref{Frobstrich1e}) in terms of $V'$. By
    definition we have
    \begin{displaymath}
V' = \pi^{-a_{\psi}} V : M_{\psi \sigma} \longrightarrow M_{\psi}.  
    \end{displaymath}
    We consider
    \begin{displaymath}
      V^{f_E} : M_{\psi} \overset{V}{\longrightarrow} M_{\psi \sigma^{-1}}
      \overset{V}{\longrightarrow} \ldots \overset{V}{\longrightarrow}
      M_{\psi \sigma^{-f_E}} .
      \end{displaymath}
    Therefore we obtain for this map
    \begin{displaymath}
V^{f_E} = \pi^{a_{\psi \sigma^{-1}}} \cdot \ldots \cdot \pi^{a_{\psi \sigma^{-f_E}}} (V')^{f_E}.  
      \end{displaymath}
    By definition of the reflex field $E$, we have $a_{\psi} = a_{\psi \sigma^{-f_E}}$. This
    implies that the number
    \begin{displaymath}
g = {a_{\psi \sigma^{-1}}} + \ldots + {a_{\psi \sigma^{-f_E}}}
    \end{displaymath}
    is independent of $\psi$. We add to the last equation
        \begin{displaymath}
g = {a_{\bar{\psi} \sigma^{-1}}} + \ldots + {a_{\bar{\psi} \sigma^{-f_E}}} .
        \end{displaymath}
        Since  $a_{\psi} + a_{\bar{\psi}}$ is $e$ or $e-1$ we obtain
        \begin{displaymath}
2g = \frac{f_E}{f} \cdot (d-1). 
        \end{displaymath}
        This proves the Lemma. 
  \end{proof}

\begin{corollary}\label{MNmorph2c}
Let $\omega_{\tau_E}$ denote the action of $\tau_E$  on the formal scheme 
scheme  $\wh{\Omega}_{F} \times_{\Spf O_F} \Spf O_{\breve{E}}$ via the second factor, i.e., in the notation of (\ref{WD13e}),
$\omega_{\tau_E} = \omega_{\tau}^{f_E/f}$. 
   There exists an isomorphism of formal schemes over $\Spf O_{\breve E}$
  $$
       \tilde{\mathcal{M}}_{r} \isoarrow
  (\wh{\Omega}_{F} \times_{\Spf O_F} \Spf O_{\breve{E}})  \times \mathbb{Z} ,
$$
such that the Weil  descent datum $\omega_{\mathcal{M}_r}$ induces on the right
hand side the Weil  descent datum
\begin{displaymath}
   (\xi, i) \mapsto (\omega_{\tau_E}(\xi), i+f_E).
\end{displaymath}
In particular the formal scheme $\mathcal{M}_{K/F, r}$ over $ \Spf O_{\breve{E}}$
is a $p$-adic formal scheme which has semi-stable reduction; hence it is also flat, with reduced special fiber. 
\end{corollary}
\begin{proof}
  We consider the isomorphisms of functors
  \begin{equation}\label{MNmorph6e}
   \tilde{\mathcal{M}}_r \isoarrow \tilde{\mathcal{N}}_{O_{\breve E}}[e]
\stackrel{\sim}{\longleftarrow} (\tilde{\mathcal{M}}_{\text{Dr}}[e])_{O_{\breve E}} \isoarrow  
  (\wh{\Omega}_{F} \times_{\Spf O_F} \Spf O_{\breve E} ) \times \mathbb{Z}e.
  \end{equation}
  The last arrow is the isomorphism (\ref{FuncMDr2e}) and the left arrow in
  the middle is the isomorphism of Theorem \ref{thm:altdrinunr}. We must see
  what the Weil descent data $\pi^g\omega_{\mathcal{N}}^{f_E/f}$ on
  $\tilde{\mathcal{N}}_{O_{\breve E}}[e]$ does on the last functor. By Theorem \ref{thm:altdrinunr},
  it induces on $(\tilde{\mathcal{M}}_{\text{Dr}}[e])_{O_{\breve E}}$ the
  Weil descent datum $\omega_{\mathcal{M}_{\text{Dr}}}^{f_E/f}$ multiplied $2g$-times
  with the translation (cf. last sentence of Proposition \ref{MNmorph2p}).
  By (\ref{WD13e}), the induced Weil descent datum on the last functor of
  (\ref{MNmorph6e}) is
  \begin{displaymath}
(\xi, ie) \mapsto (\omega_{\tau_E}, ie + (f_E/f) + 2g) .
  \end{displaymath}
  But we have
  \begin{displaymath}
ie + (f_E/f) + 2g = ie + (f_E/f) + (d-1)(f_E/f) = ie + f_E e.  
  \end{displaymath}
  This proves the Corollary. 
\end{proof}

\subsection{The case $r$ banal and $K/F$ ramified}\label{ss:tcrab} 

Let $\varepsilon \in \{\pm 1\}$. There is up to isomorphism a unique
anti-hermitian $K$-vector space $(V, \psi)$ of dimension $2$ such that
$\inv (V, \psi) = \varepsilon$, cf. Definition \ref{invdet}. Let $\Lambda \subset V$ be an $O_K$-lattice
such that $\psi$ induces a perfect pairing on $\Lambda$, cf. Lemma
\ref{Kneun01l}. By Theorem \ref{KneunBa4p}, $(\Lambda, \psi)$ corresponds to
a display
$(\mathcal{P}, \iota, \beta) \in \mathfrak{d}\mathfrak{P}^{\rm pol}_{r, \bar{\kappa}_E}$
over the residue class field of $O_{\breve{E}}$ which is unique up to
isomorphism. Let $(\mathbb{X}, \iota_{\mathbb{X}}, \beta_{\mathbb{X}})$ be the
corresponding polarized $p$-divisible $O_K$-module of type $r$. It is uniquely
determined by the conditions that $\beta$ is principal and that
$\inv^{r} (\mathbb{X}, \iota_{\mathbb{X}}, \beta_{\mathbb{X}}) = \varepsilon$. Then $(C_\BX\otimes\BQ, \phi)\simeq(V, \psi)$.

\begin{definition}\label{KneunBa1d}
  We define a functor $\mathcal{M}_{r, \varepsilon}(i)$ on the category
  $\Nilp_{O_{\breve{E}}}$. For $R \in \Nilp_{O_{\breve{E}}}$, a point of
  $\mathcal{M}_{r, \varepsilon}(i)(R)$ is given by the following data: 
\begin{enumerate}
\item[(1) ] a local CM-pair $(X, \iota)$ of CM-type $r$ over $R$
  which satisfies the Eisenstein conditions
  $({\rm EC}_r)$ relative to a fixed uniformizer $\Pi \in K$. 
\item[(2) ] an isomorphism of $p$-divisible $O_K$-modules
  \begin{displaymath}
    \lambda: X \longrightarrow X^{\wedge}, 
  \end{displaymath}
  which is a polarization of $X$. 
     \item[(3) ] a quasi-isogeny of $p$-divisible $O_K$-modules
    \begin{displaymath}
      \rho: X_{\bar{R}} \longrightarrow \mathbb{X}
      \times_{\Spec \bar{\kappa}_E}  \Spec \bar{R}, 
    \end{displaymath}
    such that the pullback quasi-isogeny $\rho^*(\lambda_{\mathbb{X}})$ differs
    from $p^i\lambda_{|X \times_S \bar{S}}$ by a scalar in $O_F^{\times}$, locally
    on $\Spec \bar{R}$.  Here $\bar{R} = R \otimes_{O_{\breve E}} \bar{\kappa}_E$.  
\end{enumerate}
Two data $(X, \iota, \lambda, \rho)$ and $(X_1, \iota_1, \lambda_1, \rho_1)$
define the same point iff there is an isomorphism of $O_K$-modules
$\alpha: X \rightarrow X_1$ such that
$\rho_1 \circ \alpha_{\bar{R}} = \rho$.  (This implies that $\alpha$ respects
the polarizations up to a factor in $O_F^{\times}$). 
\end{definition}

By Proposition \ref{KottwitzC2p}, 
$\mathcal{M}_{r, \varepsilon}(0)$ is the functor $  \mathcal{M}_{K/F, r, \varepsilon}$ used in Theorem \ref{main-banal}.
We will also consider the functor
\index[NO]{MBH@$\tilde{\mathcal{M}}_{r, \varepsilon}$}
\begin{displaymath}
 \tilde{\mathcal{M}}_{r, \varepsilon} = \coprod_{i \in \mathbb{Z}} \mathcal{M}_{r, \varepsilon}(i).
  \end{displaymath}

Let $i \in \mathbb{Z}$. We consider the following functor $\mathcal{G}_\varepsilon(i)$
on the category $\Nilp_{O_{\breve{E}}}$. A point of $\mathcal{G}_\varepsilon(i)(R)$ is given
by the following data:
\begin{enumerate}
\item[(1)] a locally
  constant $p$-adic \'etale sheaf $C$ on $\Spec R$ which is
  $\mathbb{Z}_p$-flat with $\rank_{\mathbb{Z}_p} C = 4d$ and with an action
\begin{displaymath}
\iota: O_K \longrightarrow \End_{\mathbb{Z}_p} C. 
\end{displaymath}
\item[(2)] a perfect alternating $O_F$-bilinear pairing
  \begin{displaymath}
  \phi: C \times C \longrightarrow O_F, 
  \end{displaymath}
  such that $\phi(\iota(a) c_1, c_2) = \phi (c_1, \iota(\bar{a})c_2)$ for
  $c_1, c_2 \in C$ and $a \in O_K$. 
\item[(3)] a quasi-isogeny of $O_K$-module sheaves on $\Spec R$ 
  \begin{equation}\label{KneunBa8e}
    \rho: (C, \iota) \longrightarrow (\underline{C}_{\mathbb{X}}, \iota) 
  \end{equation}
  such that locally on $\Spec R$ there is an $f \in O_F^{\times}$ with
  \begin{displaymath}
    p^if\phi(c_1, c_2) = \phi_{\mathbb{X}}(\rho(c_1), \rho(c_2)). 
    \end{displaymath}
 \end{enumerate}

Another set of data $(C', \iota', \lambda', \rho')$ defines the same point
iff there is an isomorphism $\alpha: C \isoarrow C'$ such that
$\rho' \circ \alpha = \rho$. Then $\alpha$ respects $\phi$ and $\phi'$ up
to a factor in $O_F^{\times}$.

We remark that in (\ref{KneunBa8e}) we regard
$\underline{C}_{\mathbb{X}}$ as the constant
sheaf on $\Spec R$. The existence of the quasi-isogeny implies that $C$ is locally
constant for the Zariski topology. Therefore locally on $\Spec R$ the sheaf
$C$ is the constant sheaf associated to an $O_K$-submodule
$C \subset C_{\mathbb{X}} \otimes_{\mathbb{Z}_p} \mathbb{Q}_p$ and $\rho$ is
given by the last inclusion. 

The polarized contraction functor $\mathfrak{C}_{r}^{\rm pol}$ defines a morphism of functors
\begin{equation}\label{KneunBa13e} 
  \mathcal{M}_{r, \varepsilon}(i) \longrightarrow \mathcal{G}_\varepsilon(i).
\end{equation}

To describe the functor $\mathcal{G}_\varepsilon(i)$, we may restrict to the case
where the sheaf $C$ is given by an $O_K$-submodule of
$C_{\mathbb{X}} \otimes_{\mathbb{Z}_p} \mathbb{Q}_p$. 
Then $C$ defines a point of $\mathcal{G}_\varepsilon(i)(R)$ iff
$(1/p^i) \phi_{\mathbb{X}}$ is a perfect alternating pairing on $C$.
We define an algebraic group over $\mathbb{Z}_p$, and its $\BZ_p$-rational points,
\begin{displaymath}
  J'(\mathbb{Z}_p) = \{g \in \GL_{O_K}( C_{\mathbb{X}} )\; | \;
  \phi_{\mathbb{X}}(gc_1, gc_2) = f\cdot \phi_{\mathbb{X}}(c_1, c_2) \;
  \text{for some} \; f \in O_F^{\times}\}. 
  \end{displaymath}
By Lemma \ref{Kneun01l}, there is an isomorphism
$g: (C_{\mathbb{X}}, \phi_{\mathbb{X}}) \longrightarrow (C, \frac{1}{p^i} \phi_{\mathbb{X}})$.
This means that $gC_{\mathbb{X}} = C$ and
\begin{equation}\label{KneunBa12e}
  \phi_{\mathbb{X}}(gc_1, gc_2) = p^i \phi_{\mathbb{X}}(c_1, c_2) , \;
  \quad  c_1, c_2 \in C_{\mathbb{X}} \otimes_{\mathbb{Z}_p} \mathbb{Q}_p .
  \end{equation}
We define 
\begin{displaymath}
  J'(i) = \{ g \in \GL_{K}(C_{\mathbb{X}} \otimes_{\mathbb{Z}_p} \mathbb{Q}_p)\; | \;
  \phi_{\mathbb{X}}(gc_1, gc_2) = p^if \phi_{\mathbb{X}}(c_1, c_2), \;
  \text{for some} \; f \in O_F^{\times}\}.
  \end{displaymath}
This construction gives us a functor isomorphism
\begin{displaymath}
\mathcal{G}_\varepsilon(i) \isoarrow J'(i)/J'(\mathbb{Z}_p) ,
  \end{displaymath}
where the right hand side is considered as the restriction of the constant
sheaf to $\Nilp_{O_{\breve{E}}}$.
Let $J'\subset \GL(C_\BX)$ be the union of the $J'(i)$. Using the contraction
functor we may write
\begin{equation}\label{J'banalram1e}
  J' = \{\alpha \in \Aut^{o}_{K} \mathbb{X} \; |\;
  \alpha^{\ast} (\lambda_{\mathbb{X}}) = \mu(\alpha) \lambda_{\mathbb{X}} \;
   \text{for some} \; \mu(\alpha) \in p^{\mathbb{Z}} O_F^{\times}\}. 
  \end{equation}
Therefore $J'$ acts via $\rho$ on the functor
$\tilde{\mathcal{M}}_{r, \varepsilon}$. 
\begin{proposition}\label{KneunBa5p} 
  The morphism of functors on $\Nilp_{O_{\breve{E}}}$ obtained from 
  (\ref{KneunBa13e}) is a $J'$-equivariant isomorphism,
  \begin{equation}\label{KneunBa10e}
\tilde{\mathcal{M}}_{r, \varepsilon} \isoarrow J'/J'(\mathbb{Z}_p). 
  \end{equation}
  The left hand side is endowed with the Weil  descent datum
  $\omega_{\mathcal{M}_r}$ relative to $O_{\breve{E}}/O_E$ which is defined exactly
  in the same way as (\ref{WD7e}). This Weil  descent datum 
  corresponds on the right hand side  to the Weil  descent
  datum given by multiplication with $\Pi^{ef_E}$. Here we regard $\Pi^{ef_E}$
  as an automorphism of the $K$-vector space
  $C_{\mathbb{X}} \otimes_{\mathbb{Z}_p} \mathbb{Q}_p$.
\end{proposition}
\begin{proof}
  We show that (\ref{KneunBa13e}) is an isomorphism of functors.
  Let $R' \rightarrow R$ be an epimorphism of $O_{\breve{E}}$-algebras with
  nilpotent kernel. We claim that the induced map 
  $\mathcal{M}_{r,\varepsilon}(i)(R') \rightarrow \mathcal{M}_{r,\varepsilon}(i)(R)$
  is bijective. We may assume that the kernel is endowed with divided powers.
  Let $(X, \iota, \lambda, \rho)$ be as in Definition \ref{KneunBa1d}. To see
  that $(X,\iota)$ lifts uniquely to $R'$ we apply Corollary \ref{GM1c}. Let
  $\mathcal{P}$ be the display of $(X, \iota)$ over $R$ and let $\mathbb{D}$
  be the associated crystal (\ref{displaycrystal1e}). By Proposition
  \ref{KottwitzC2p}, the Hodge filtration of $\mathcal{P}$ is
  $\mathbf{E}_{A_{\psi}} \mathbb{D}_{\psi}(R) \subset \mathbb{D}_{\psi}(R)$. The
  only possibility to lift this filtration to $R'$ such that the condition
  $({\rm EC}_r)$ continues to hold is to take
  $\mathbf{E}_{A_{\psi}} \mathbb{D}_{\psi}(R') \subset \mathbb{D}_{\psi}(R')$. Hence
  we obtain a unique lifting $(X', \iota')$ over $R'$. That $\lambda$ lifts
  to a principal polarization $\lambda'$ of $X'$ follows from Proposition
  \ref{GM2p} and Lemma \ref{Cbanal1l}. Since the quasi-isogeny $\rho$ lifts
  automatically, we obtain the claimed bijectivity. The functor on the right
  hand side of (\ref{KneunBa13e}) also induces a bijection when applied to
  $R' \rightarrow R$ because it is defined in terms of \'etale sheaves.  

  Using these bijections we see that it is enough to prove the Proposition for
  the restrictions of the functors to the  category of
  $\bar{\kappa}_E$-algebras.  
  Both functors commute with inductive limits of rings, i.e., they are locally
  of finite presentation. To see this, we consider the special fiber
  $\mathcal{M}_{{r, \varepsilon},\bar{\kappa}_E}(0)$.
  Let $m \in \mathbb{N}$. We consider the subfunctor $U_m$ which consists of
  points such that $p^m\rho^{-1}$ is an isogeny. Then $X$ is the quotient of
  $\mathbb{X} \times_{\Spec \bar{\kappa}_E}  \Spec R$ by a finite locally free
  subgroup scheme of $\mathbb{X}(4dm) \times_{\Spec \bar{\kappa}_E}  \Spec R$.
  (We have denoted by $\mathbb{X}(4dm)$ the kernel of the multiplication by
  $p^{4dm}$.) This shows that $U_m$ is a scheme of finite presentation over
  $\Spec \bar{\kappa}_E$. Therefore $U_m$ is locally of finite representation
  as a functor, cf. EGA IV, Thm.~(8.8.2). One easily deduces from this that
  $\mathcal{M}_{{r, \varepsilon},\bar{\kappa}_E}(0)$ is locally of finite presentation as a
  functor. In the same way we see that $\mathcal{G}_\varepsilon(i)$ is locally of finite
  presentation. 
  To show that (\ref{KneunBa13e}) is an isomorphism of functors we can
  therefore restrict to $\bar{\kappa}_E$-algebras $R$ which are of finite
  type over $\bar{\kappa}_E$. For such $R$, (\ref{KneunBa13e}) induces a
  bijection by Theorem \ref{KneunBa4p}.

  It remains to compare the Weil descent data on both sides of
  (\ref{KneunBa10e}).  
  It is enough to make the comparison
  for the restriction of (\ref{KneunBa10e}) to the category of
  $\bar{\kappa}_E$-algebras. Let $\varepsilon: \bar{\kappa}_E \longrightarrow R$ be a $\bar{\kappa}_E$-algebra.  The restriction of the functor
  $\mathcal{M}_{r, \varepsilon}$ to $\bar{\kappa}_E$-algebras has a Weil descent datum $\omega_{\mathcal{M}_{r, \varepsilon}, \mathbb{F}_p}: \mathcal{M}_{r, \varepsilon}(R) \longrightarrow \mathcal{M}_{r, \varepsilon}(R_{[\sigma]})$
  over $\mathbb{F}_p$, given by
  \begin{equation}\label{KneunWAb3e}
    \omega_{\mathcal{M}_{r, \varepsilon}, \mathbb{F}_p} ((X, \iota, \lambda, \rho)) =
    (X, \iota, \lambda, \rho_{[\sigma]}),
    \end{equation}
  where $\rho_{[\sigma]}$ is the composite
  \begin{displaymath}
    X \overset{\rho}{\longrightarrow} \varepsilon_{\ast}\mathbb{X}
    \overset{\varepsilon_{\ast} F_{\mathbb{X}}}{\longrightarrow}
    \varepsilon_{\ast} \sigma_{\ast}\mathbb{X}.
    \end{displaymath}
  Here $\sigma$ denotes the Frobenius automorphism of
  $\bar{\kappa}_E$ over $\mathbb{F}_p$. To see that this makes sense, we have to check that all $p$-divisible
  $O_K$-modules above satisfy the rank condition $({\rm RC}_r)$ and the Eisenstein conditions $({\rm EC}_r)$.
 The condition $({\rm EC}_r)$   says that  $\Pi^{e}\Lie X = 0$.
  This remains true if we regard $X$ as a $p$-divisible $O_K$-module
  over $R_{[\sigma]}$. For the condition $({\rm RC}_r)$, the claim is obvious.

  Therefore it suffices to show that $\omega_{\mathcal{M}_r, \mathbb{F}_p}$
  induces on the
  right hand side of (\ref{KneunBa10e}) the Weil descent datum
  $g \mapsto \Pi^eg$. This follows from the following Lemma.
\end{proof}

\begin{lemma}
  There is an identification 
  $C_{\mathbb{X}} = C_{\sigma_{\ast} \mathbb{X}}$. The functor $\mathfrak{C}_{r, \bar{\kappa}_E}^{\rm pol}$applied to the Frobenius morphism 
  $F_{\mathbb{X}}: \mathbb{X} \longrightarrow \sigma_{\ast} \mathbb{X}$ 
  yields $\Pi^{e}: C_{\mathbb{X}} \longrightarrow C_{\mathbb{X}}$. 
  \end{lemma}
\begin{proof} The first assertion follows because  the functor $\mathfrak{C}_{r, \bar{\kappa}_E}^{\rm pol}$ commutes with base change. 
  Consider the Dieudonn\'e module $P_{\mathbb{X}}$ of $\mathbb{X}$. We have
  \begin{displaymath}
C_{\mathbb{X}} = \{c \in P_{\mathbb{X}} \; | \; Vc = \Pi^e c \}, 
  \end{displaymath}
  cf. Remark \ref{CPexplban}.
  The map
  \begin{displaymath}
    P_{\mathbb{X}} \longrightarrow W(\bar{\kappa}_E)
    \otimes_{\sigma, W(\bar{\kappa}_E)} P_{\mathbb{X}}, \quad c \mapsto 1 \otimes c 
  \end{displaymath}
  defines the identification $C_{\mathbb{X}} = C_{\sigma_{\ast} \mathbb{X}}$. The
  Frobenius $F_{\mathbb{X}}$ induces on the Dieudonn\'e modules
  \begin{displaymath}
    P_{\mathbb{X}} \overset{V^{\sharp}}{\longrightarrow} W(\bar{\kappa}_E)
    \otimes_{\sigma, W(\bar{\kappa}_E)} P_{\mathbb{X}}, \quad x \mapsto 1 \otimes Vx. 
  \end{displaymath}
  For $c \in C_{\mathbb{X}}$ we obtain
  $V^{\sharp}c = 1 \otimes Vc = 1 \otimes \Pi^ec$.
  \end{proof}

We can reformulate a part of Proposition \ref{KneunBa5p} as follows.
We consider the algebraic group over $\mathbb{Z}_p$ such that 
\index[NO]{JAH@$J(\mathbb{Z}_p)$}
\begin{displaymath}
  J(\mathbb{Z}_p) = \{g \in \GL_{O_K} (C_{\mathbb{X}}) \; | \;
  \phi_{\mathbb{X}}(gc_1, gc_2) = u\cdot \phi_{\mathbb{X}}(c_1, c_2) \;
  \text{for some} \; u \in \mathbb{Z}_p^{\times}\}. 
  \end{displaymath}
We define 
\begin{displaymath}
  J(i) = \{g \in \GL_{K} (C_{\mathbb{X}} \otimes \mathbb{Q}) \; | \;
  \phi_{\mathbb{X}}(gc_1, gc_2) = up^i\cdot \phi_{\mathbb{X}}(c_1, c_2) \;
  \text{for some} \; u \in \mathbb{Z}_p^{\times}\}. 
\end{displaymath}
The union of the $J(i)$ is the group $J(\mathbb{Q}_p)$
\index[NO]{JAG@$J(\mathbb{Q}_p)$} of unitary similitudes with similitude factor in $\BQ_p^\times$. 
\begin{corollary}\label{KneunBa4c}
  Let $J(\BQ_p)$ be the group of unitary similitudes of $C_\BX\otimes_{\BZ_p}\BQ_p$ with similitude factor in $\BQ_p^\times$, and let $J(\BZ_p)$ be its subgroup stabilizing the lattice $C_\BX$.  There are isomorphisms of functors on $\Nilp_{O_{\breve E}}$, 
  \begin{displaymath}
\tilde{\mathcal{M}}_{r, \varepsilon} \isoarrow J(\mathbb{Q}_p)/J(\mathbb{Z}_p), \, \, \quad \CM_{K/F, r, \varepsilon} \isoarrow J(\mathbb{Q}_p)^o/J(\mathbb{Z}_p) .
    \end{displaymath}
  Here  $J(\mathbb{Q}_p)^o$\index[NO]{JAI@$J(\mathbb{Q}_p)^o$} denotes
  the group of unitary similitudes with similitude factor in $\BZ_p^\times$.
    \end{corollary}
\begin{proof}
  It is enough to show that $\mathcal{G}_\varepsilon(i)(\bar{\kappa}_E)$ is in bijection
  with $J(i)/J'(\mathbb{Z}_p)$. For this is enough to show that
  for each $C \in \mathcal{G}_\varepsilon(i)(\bar{\kappa}_E)$ there exists $g \in J(i)$
  such that $g C_{\mathbb{X}} = C$. This we have already shown before
  (\ref{KneunBa12e}). 
\end{proof}
In this reformulation it is less obvious what the Weil  descent datum is.

\subsection{The case  $r$ banal and $K/F$ unramified}\label{sec6unr}

Let $\varepsilon \in \{ \pm 1\}$. We consider a  CM-triple
$(\mathbb{X}, \iota, \lambda_{\mathbb{X}})$ over $\bar{\kappa}_E$ such that
$\lambda_{\mathbb{X}}$ is principal if $\varepsilon =1$ and is almost principal
if $\varepsilon = -1$. By Proposition \ref{CFbanal5p} such a CM-triple exists
and $\inv^r (\mathbb{X}, \iota, \lambda_{\mathbb{X}}) = \varepsilon$. In fact, by Lemma \ref{Kneun0l} and Theorem \ref{KneunBa4p}, $(\mathbb{X}, \iota, \lambda_{\mathbb{X}})$ is unique up to isomorphism. 

We recall the functor $\CM_{K/F, {r, \varepsilon}}$  from  section \ref{ss:defofform}.   
\begin{definition}\label{ubF1d}
  Let $(\mathbb{X}, \iota, \lambda_{\mathbb{X}})$ be a framing object with an
  almost principal polarization. We define a functor $\mathcal{M}_{r^-}(i)$ on the category
$\Nilp_{O_{\breve{E}}}$. For $R \in \Nilp_{O_{\breve{E}}}$, a point of
$\mathcal{M}_{r^-}(i)(R)$ is given by the following data: 
  \begin{enumerate}
  \item[(1) ] a local CM-triple $(X, \iota, \lambda)$ of type $r$ over $R$
    which satisfies the Eisenstein conditions $({\rm EC}_r)$ relative to
    the fixed uniformizer $\pi \in K$. 
  \item[(2) ] the polarization $\lambda$ is almost principal. 
  \item[(3) ] a quasi-isogeny of $p$-divisible $O_K$-modules
    \begin{displaymath}
      \rho: X_{\bar{R}} \longrightarrow \mathbb{X}
      \times_{\Spec \bar{\kappa}_E}  \Spec \bar{R}, 
    \end{displaymath}
    such that the pullback quasi-isogeny $\rho^*(\lambda_{\mathbb{X}})$ differs
    from $p^i\lambda_{|X_{\bar{R}}}$ by a scalar in $O_F^{\times}$, locally
    on $\Spec \bar{R}$.  Here $\bar{R} = R \otimes_{O_{\breve E}} \bar{\kappa}_E$.  
  \end{enumerate}
  \end{definition}
  Now let $(\mathbb{X}, \iota, \lambda_{\mathbb{X}})$ be a framing object with a 
principal polarization. Then we have $\tt{h}(\lambda_{\mathbb{X}}) = 0$
and $\inv^{r} (\mathbb{X}, \iota, \lambda_{\mathbb{X}}) = 1$.
We define the functor $\mathcal{M}_{r^+}(i)$ by exactly the same data
but we replace the condition (2) above by the condition 
\begin{enumerate}
  \item[$(2^\prime)$ ] the polarization $\lambda$ is principal.
  \end{enumerate}

In the almost principal case there exists an isogeny
$X^{\wedge} \longrightarrow X$ such that the composite
\begin{displaymath}
  X \overset{\lambda}{\longrightarrow} X^{\wedge} \longrightarrow X   
\end{displaymath}
is $\pi\id_{X}$. This follows from the following analogue of Proposition \ref{SFB701p}. 
\begin{proposition}\label{ubF1c} 
  Let $\alpha: \mathcal{P}_1 \longrightarrow \mathcal{P}_2$ be an isogeny of
  CM-pairs of type $r$ over $R\in\Nilp_{O{\breve E}}$ which both satisfy the Eisenstein condition $({\rm EC}_r)$. Let
  $\alpha_C: C_1 \longrightarrow C_2$ be the morphism in ${\rm Et}(O_K)_R$ associated by the contracting functor $\mathfrak{C}_{r, R}$.
  Then
  \begin{displaymath}
    \height \alpha = 2f\cdot \mathrm{length}_{O_K} \Coker \alpha_C.
        \end{displaymath}
  If $\height \alpha = 2f$ then there exists  a unique morphism
  $\beta: \mathcal{P}_2 \longrightarrow \mathcal{P}_1$ such that
  \begin{displaymath}
    \alpha \circ \beta = \pi \id_{\mathcal{P}_2}, \quad
    \beta \circ \alpha = \pi \id_{\mathcal{P}_1}, 
    \end{displaymath}
  \end{proposition} 
\begin{proof}
  To prove the first assertion we can assume that $R = k$ is an algebraically
  closed field. Then we can use that
  $\Coker \alpha = \Coker \alpha_C \otimes_{\mathbb{Z}_p} W(k)$.
If $\mathrm{length}_{O_K} \Coker \alpha_C = 1$, then $\beta_C$
  clearly exists.   
\end{proof}

In the case where $i = 0$, the quasi-isogeny $\rho$ is of height zero because
the polarizations $\lambda_{\mathbb{X}}$ and $\lambda$ have the same degree.
We set $X_0 =X$ and $X_1 = X^{\wedge}$. Since $X$ and
$X^{\wedge}$ satisfy the Eisenstein condition $({\rm EC}_r)$ and, by Proposition
\ref{KottwitzC2p},  also the Kottwitz condition $({\rm KC}_r)$, we obtain a point of the functor
$\mathcal{M}_{K/F, r, -1}$ of Definition \ref{CRFd2}.  We conclude that 
$\mathcal{M}_{r^-}(0) = \mathcal{M}_{K/F, r, 1}$. The index $r^-$ on the left hand
side indicates that we are in the case where the adjusted invariant of the
framing object is $-1$. Similarly, $\mathcal{M}_{r^+}(0) = \mathcal{M}_{K/F, r, 1}$. The index $r^+$ on the left hand
side indicates that we are in the case where the adjusted invariant of the
framing object is $1$.  

We will describe the formal scheme which represents the functor
\index[NO]{MBI@$\tilde{\mathcal{M}}_{r^{\pm}}$}
$\mathcal{M}_{K/F, {r, \varepsilon}}$. 
More precisely, consider the functors on $\Nilp_{O_{\breve{E}}}$,
  \begin{displaymath}
 \tilde{\mathcal{M}}_{r^\pm} = \coprod_{i \in \mathbb{Z}} \mathcal{M}_{r^\pm}(i).
  \end{displaymath} 
  These functors are endowed with a Weil  descent datum
  \index[NO]{ZZWH@$\omega_{\mathcal{M}_{r^\pm}}$}
  $\omega_{\mathcal{M}_{r^\pm}}: \mathcal{M}_{r^\pm}(i) \longrightarrow \mathcal{M}_{r^\pm}(i+f_E)^{\tau_E}$ 
  relative to $\breve{E}/E$ using exactly the same definition as in (\ref{WD7e}).  Recall $(C_{\mathbb{X}}, \phi_{\mathbb{X}})=(C_{\mathbb{X}},\iota_\BX, \phi_{\mathbb{X}})$. We define\index[NO]{JAJ@$J'(\mathbb{Z}_p)$} 
  \begin{displaymath}
  \begin{aligned} 
  J'(\mathbb{Z}_p) &=  \{g \in \Aut_{O_K} C_{\mathbb{X}} \; | \;
  \phi_{\mathbb{X}}(gc_1, gc_2) = f\cdot \phi_{\mathbb{X}}(c_1, c_2) \;
  \text{for some} \; f \in O_F^{\times}\},\\ 
  J'(i) &=  \{ g \in \Aut C_{\mathbb{X}} \otimes_{\mathbb{Z}_p} \mathbb{Q}_p\; | \;
  \phi_{\mathbb{X}}(gc_1, gc_2) = p^if \phi_{\mathbb{X}}(c_1, c_2), \;
  \text{for some} \; f \in O_F^{\times}\}.
    \end{aligned}
\end{displaymath}

As in the ramified case,
we see that a point of $\mathcal{M}_{r^\pm}(i)(R)$ is locally on $\Spec R$ given by
a lattice $C \subset C_{\mathbb{X}} \otimes_{\mathbb{Z}_p} \mathbb{Q}_p$ such
that the restriction of $(1/p^i) \phi_{\mathbb{X}}$ to $C$ induces a
$O_F$-bilinear form
\begin{displaymath}
  \frac{1}{p^i} \phi_{\mathbb{X}}: C \times C \longrightarrow O_F 
  \end{displaymath}
which is perfect in the case of $r^+$ and such that $\ord_{\pi} \frac{1}{p^i} \phi_{\mathbb{X}} = 1$ in the case of $r^-$. For the case of $r^-$, we are using here Proposition \ref{ubF1c}. 

We deduce that there is an isometry
$g: (C_{\mathbb{X}}, \phi_{\mathbb{X}}) \longrightarrow (C, (1/p^i) \phi_{\mathbb{X}})$, cf. Lemma \ref{Kneun0l}.
Consequently we have $g \in J'(i)$. Conversely, if $g \in J'(i)$, the sublattice
$C = gC_{\mathbb{X}}$ with the bilinear form $(1/p^i) \phi_{\mathbb{X}}$ gives rise
to a point of $\mathcal{M}_{r^\pm}(i)(R)$.  

We will denote by $J' \subset \Aut C_{\mathbb{X}}$ the union \index[NO]{JAG@$J'$}
of the $J'(i)$.  
We can identify $J'$ with a subgroup of $\Aut^{o}_K \mathbb{X}$ exactly
as in the ramified case, cf. \eqref{J'banalram1e}. It acts via $\rho$ on  the
functor $\tilde{\mathcal{M}}_{r^{\pm}}$. 

\begin{proposition}\label{ubF2p}
  There is a $J'$-equivariant isomorphism of functors on $\Nilp_{O_{\breve{E}}}$,
  \begin{equation}\label{ubF9e}
\tilde{\mathcal{M}}_{r^\pm} \isoarrow J'/J'(\mathbb{Z}_p).  
    \end{equation}
Here the right hand side is the constant sheaf on $\Nilp_{O_{\breve{E}}}$. 
The Weil  descent datum $\omega_{\mathcal{M}_{r,^\pm}}$ on the left hand side
corresponds on the right hand side  to the Weil descent datum given by
multiplication with $\pi^{ef_E/2}$. Here we view
$\pi^{ef_E/2}$ as an automorphism of the $K$-vector space
$C_{\mathbb{X}} \otimes_{\mathbb{Z}_p} \mathbb{Q}_p$ by multiplication. 
\end{proposition}
\begin{proof}
  That (\ref{ubF9e}) is an isomorphism of functors follows from Theorem
  \ref{KneunBa4p} in the same way as in the proof of Proposition
  \ref{KneunBa5p}.
 
  Let us recall the definition of the Weil  descent relative to
  $O_{\breve{E}}/O_E$ on the functor
  $\tilde{\mathcal{M}}_{r^\pm}$.  We write
  \begin{equation}\label{frobunr}
F_{\mathbb{X}, \tau_E}: \mathbb{X} \longrightarrow (\tau_{E})_{\ast} \mathbb{X}  
  \end{equation} 
  for the Frobenius relative to $\kappa_E$. Let $\varepsilon: O_{\breve{E}} \longrightarrow R$ be an object of
  $\Nilp_{O_{\breve{E}}}$.   We write
  $\bar{\varepsilon} = \varepsilon \otimes \bar{\kappa}_E: \bar{\kappa}_E \longrightarrow \bar{R}$.
Let 
  $(X, \iota, \lambda, \rho) \in \mathcal{M}_{r^\pm}(i)(R)$ be a point. We view $(X, \iota, \lambda)$ as a CM-triple on $R_{[\tau_E]}$ and we endow it 
  with the framing 
  \begin{displaymath}
    \tilde\rho: X_{\bar{R}} \overset{\rho}{\longrightarrow}
    \bar{\varepsilon}_{\ast}\mathbb{X} \; 
    \overset{\bar{\varepsilon}_{\ast} F_{\mathbb{X}, \tau_E}}{\longrightarrow}
    \bar{\varepsilon}_{\ast} (\tau_E)_{\ast}\mathbb{X} .
  \end{displaymath}
  Then $(X, \iota, \lambda, \tilde\rho)$ defines a point of
  $\mathcal{M}_{r^\pm}(i+f_E)(R_{[\tau_E]})$.
  Varying $i \in \mathbb{Z}$, we obtain the Weil  descent datum
  \begin{displaymath}
\omega_{\mathcal{M}_{r^\pm}}: \tilde{\mathcal{M}}_{r^\pm} \longrightarrow \tilde{\mathcal{M}}_{r^\pm}^{(\tau_E)}.
    \end{displaymath}
    We note that the inverse image of $(\tau_E)_{\ast} \lambda$ by (\ref{frobunr}) is
  $p^{f_E} \lambda$. 
The compatibility of the Weil  descent data follows as in the proof of
Proposition \ref{KneunBa5p} from the following Lemma. 
 \end{proof}
\begin{lemma}
The contracting functor applied to the Frobenius morphism 
  $F_{\mathbb{X}, \tau_E}: \mathbb{X} \longrightarrow (\tau_E)_{\ast}\mathbb{X}$  
yields the multiplication by
$\pi^{ef_E/2}: C_{\mathbb{X}} \longrightarrow C_{\mathbb{X}}$. 
  \end{lemma}
\begin{proof}
  We use the Dieudonn\'e module $P$ of $\mathbb{X}$ over $\bar{\kappa}_E$.
  The map $F_{\mathbb{X}, \tau_E}$ induces on the Dieudonn\'e modules the map
  \begin{displaymath}
    V^{f_E, \sharp}: P \longrightarrow W(\bar{\kappa}_E) \otimes_{F^{f_E}, \bar{\kappa}_E} P,
    \quad  x \mapsto 1 \otimes V^{f_E} x. 
  \end{displaymath}
  By definition we have
  \begin{displaymath}
C:= C_{\mathbb{X}} = \{c \in P \; | \; Vc = \pi_r c \} ,
    \end{displaymath}
    where we recall $\pi_r$ from \eqref{pirunr}, cf. Remark \ref{CPexplban}. 
  With the identification $C_{\mathbb{X}} = C_{{(\tau_{E})_\ast}\mathbb{X}}$, the
  restriction of $V^{f_E, \sharp}$ to $C$ gives
  \begin{equation}\label{ubF14e}
    ~^{F^{-f_E+1}}\pi_r\cdot \ldots \cdot ~^{F^{-1}}\pi_r \cdot \pi_r:
    C \longrightarrow C.  
  \end{equation}
  The right hand side is a module over 
  \begin{displaymath}
    O_K \otimes_{\mathbb{Z}_p} W(\bar{\kappa}_E) = \prod_{\psi \in \Psi}
    O_K \otimes_{O_{K^t}, \tilde{\psi}} W(\bar{\kappa}_E).
  \end{displaymath}
  On the right hand side, $F^{-1}$ is given by the map
  \begin{displaymath}
    O_K \otimes_{O_{K^t}, \tilde{\psi \sigma}} W(\bar{\kappa}_E) \longrightarrow
    O_K \otimes_{O_{K^t}, \tilde{\psi}} W(\bar{\kappa}_E), \quad a \otimes \xi
    \mapsto a \otimes ~^{F^{-1}} \xi. 
  \end{displaymath}
  Therefore the components of the element on the left hand side of
  (\ref{ubF14e}) are
  \begin{displaymath}
\pi^{a_{\psi \sigma^{(f_E-1)}}} \cdot \ldots \cdot \pi^{a_{\psi \sigma}} \cdot \pi^{a_{\psi}}.
  \end{displaymath}
  Since $\sigma^{f_E}$ fixes $\kappa_E$ we have $a_{\psi \sigma^{f_E}} = a_{\psi}$.
  It follows that the numbers
  \begin{displaymath}
g_{\psi} := a_{\psi \sigma^{(f_E-1)}} +  \ldots a_{\psi\sigma} + a_{\psi}
    \end{displaymath}
  are independent of $\psi$. We call this number $g$. We find:
  \begin{displaymath}
2g = g_{\psi} + g_{\bar{\psi}} = ef_E
  \end{displaymath}
  because $a_{\psi} + a_{\bar{\psi}} = e$.
  We conclude that (\ref{ubF14e}) is the multiplication by $\pi^{ef_E/2}$. 
\end{proof}
\begin{corollary}\label{J'unr}
 Let $J(\BQ_p)$ be the group of unitary similitudes of $C_\BX\otimes_{\BZ_p}\BQ_p$ with similitude factor in $\BQ_p^\times$, and let $J(\BZ_p)$ be its subgroup stabilizing the lattice $C_\BX$. 
  There are isomorphisms of functors on $\Nilp_{O_{\breve E}}$, 
  \begin{displaymath}
\tilde{\mathcal{M}}_{r, \varepsilon} \isoarrow J(\mathbb{Q}_p)/J(\mathbb{Z}_p), \, \, \quad \CM_{K/F, r, \varepsilon} \isoarrow J(\mathbb{Q}_p)^o/J(\mathbb{Z}_p) .
    \end{displaymath}
    Here $J(\mathbb{Q}_p)^o$ denotes the group of unitary similitudes with similitude factor in $\BZ_p^\times$.\qed
    \end{corollary}

\subsection{The banal split case}\label{ss:tbsc} 

We assume that $K = F \times F$. 
Let $R = k$ be an algebraically closed field. There is up to isomorphism a
unique anti-hermitian $O_K$-module $(C, \phi)$ of rank $2$ with $\phi$ perfect. Hence 
there is a unique CM-triple
\begin{displaymath}
  (\mathbb{X}, \iota_{\mathbb{X}}, \lambda_{\mathbb{X}}) \in
  \mathfrak{P}^{\rm pol}_{r,\bar{\kappa}_E}
   \end{displaymath}
with principal $\lambda_{\mathbb{X}}$. We take this as  framing object. 

  We define  functors $\mathcal{M}_r(i)$ on the
  category $\Nilp_{O_{\breve{E}}}$.  For $R \in \Nilp_{O_{\breve{E}}}$, a point of
  $\mathcal{M}_{r}(i)(R)$ is given by the following data: 
  \begin{enumerate}
  \item[(1) ] a local CM-triple $(X, \iota, \lambda)$ of type $r$ over $R$
    which satisfies the Eisenstein conditions $({\rm EC}_r)$ relative to
    the fixed uniformizer $\pi \in F$. 
  \item[(2) ] the polarization $\lambda\colon X\to X^\wedge$ is principal. 
  \item[(3) ] a quasi-isogeny of $p$-divisible $O_K$-modules
    \begin{displaymath}
      \rho: X_{\bar{R}} \longrightarrow \mathbb{X}
      \times_{\Spec \bar{\kappa}_E}  \Spec \bar{R}, 
    \end{displaymath}
    such that the pullback quasi-isogeny $\rho^*(\lambda_{\mathbb{X}})$ differs
    from $p^i\lambda_{|X_{\bar{R}}}$ by a scalar in $O_F^{\times}$, locally
    on $\Spec \bar{R}$.  Here $\bar{R} = R \otimes_{O_{\breve E}} \bar{\kappa}_E$.  
  \end{enumerate}  

  Note that $\mathcal{M}_{r}(0)=\mathcal{M}_{r, 1}$ of  before Theorem \ref{main-banal}. Consider a point $(X, \iota, \lambda, \rho)$ as above. Let $(C, \phi)$ be
  the $p$-adic \'etale sheaf associated to $(X, \iota, \lambda)_{\bar{R}}$ on
  $(\Spec \bar{R})_{\text{\'et}} = (\Spec R)_{\text{\'et}}$.  Let $\underline{C}_{\mathbb{X}}$ be
  the constant sheaf on $(\Spec R)_{\text{\'et}}$
  of the $O_K$-module $C_{\mathbb{X}}$. The existence of $\rho$ implies that $C$ is locally constant for the Zariski topology. Therefore, locally on $\Spec R$, we may regard $C$
  as a submodule of $C_{\mathbb{X}} \otimes \mathbb{Q}$. By the definition of
  $\mathcal{M}_r(i)$, we have
  \begin{displaymath}
    fp^i \phi(x, y) = \phi_{\mathbb{X}}(x, y), \quad x, y \in
    C_{\mathbb{X}} \otimes \mathbb{Q}, 
    \end{displaymath}
  for some $f \in O_F^{\times}$. We see by Theorem \ref{KneunBa4p} that a
  point of $\mathcal{M}_r(i)(R)$ is the same as a $O_K$-sublattice
  $C \subset (C_{\mathbb{X}})_R \otimes \mathbb{Q}$ such that the restriction of
  $(1/p^i)\phi_{\mathbb{X}}$ to $C$ induces a perfect pairing
  \begin{displaymath}
C \times C \longrightarrow O_F.  
    \end{displaymath}
  This is directly clear if the ideal of nilpotent elements of $R$ is
  nilpotent, and follows from the argument in the proof of Proposition
  \ref{KneunBa5p} in the general case. 
  
  Again we set
  \begin{displaymath}
  \begin{aligned}
  J'(\mathbb{Z}_p) &=  \{g \in \GL_{O_K} (C_{\mathbb{X}}) \; | \;
  \phi_{\mathbb{X}}(gc_1, gc_2) = f\cdot \phi_{\mathbb{X}}(c_1, c_2), \;
  \text{for some} \; f \in O_F^{\times}\},\\ 
  J'(i) &=  \{ g \in \GL( C_{\mathbb{X}} \otimes_{\mathbb{Z}_p} \mathbb{Q}_p)\; | \;
  \phi_{\mathbb{X}}(gc_1, gc_2) = p^if \phi_{\mathbb{X}}(c_1, c_2), \;
  \text{for some} \; f \in O_F^{\times}\}.
    \end{aligned}
  \end{displaymath}
  There is an isometry up to a constant in $O^{\times}_F$,
  \begin{displaymath}
    g: (C_{\mathbb{X}}, \phi_{\mathbb{X}}) \longrightarrow
    (C, \frac{1}{p^i}\phi_{\mathbb{X}}).
    \end{displaymath}
  Then $g \in J'(i)$ and $gC_{\mathbb{X}} = C$. Any other isometry of this type
  is of the form $gh$, with $h \in J'(\mathbb{Z}_p)$. Therefore we have
  associated to the point $(X, \iota, \lambda, \rho)$ a section of the
  constant sheaf $J'(i)/J'(\mathbb{Z}_p)$. 

  We set
  \begin{displaymath}
    \tilde{\mathcal{M}}_r = \coprod_{i \in \mathbb{Z}} \mathcal{M}_r(i), \quad
    J' = \cup_{i \in \mathbb{Z}} J'(i). 
  \end{displaymath}
  As in the ramified case,  the group $J'$ acts via $\rho$
  on the functor $\tilde{\mathcal{M}}_r$, cf.  \eqref{J'banalram1e}. Moreover,  this functor 
  is endowed with the Weil descent datum
$\omega_{\mathcal{M}_r}: \mathcal{M}_r(i) \longrightarrow \mathcal{M}_r(i+f_E)^{(\tau_E)}$ 
  relative to $O_{\breve{E}}/O_{E}$.

  Let $\sigma \in \Gal(F^t/F)$ be the Frobenius. We use the notation introduced below \eqref{zerlegt1e}. We fix $\theta \in \Theta$, with $\theta_1, \theta_2\in \Psi$. 
 Set
  \begin{equation}\label{zerlegt6e}
a_{1,E} = a_{\theta_1 \sigma^{f_E-1}} + \ldots + a_{\theta_1 \sigma} + a_{\theta_1}
  \end{equation}
  This number is independent of the choice of $\theta$ because, by the
  definition of the reflex field $E$,
  \begin{displaymath}
a_{\theta_1 \sigma^{f_E}} = a_{\theta_1}.  
  \end{displaymath}
  We already defined $a_1 = \sum_{\theta \in \Theta} a_{\theta_1}$,
  cf. (\ref{sumsplit}). If we sum the right hand side of (\ref{zerlegt6e})
  over all $\theta \in \Theta$ we obtain
  \begin{displaymath}
f a_{1,E} = f_E a_1. 
    \end{displaymath}
  In the same way we define $a_{2,E}$ by using $\theta_2$. Then we obtain
  $f a_{2,E} = f_E a_2$. Since $a_1 + a_2 = d$ we find 
  \begin{displaymath}
a_{1,E} + a_{2,E} = ef_E. 
  \end{displaymath}
 The endomorphism 
  \begin{equation}\label{banaldiag}
    \pi^{a_E} := \pi^{a_{1,E}} \oplus \pi^{a_{2,E}}:
    C_{\mathbb{X},1} \oplus C_{\mathbb{X},2} \longrightarrow  
    C_{\mathbb{X},1} \oplus C_{\mathbb{X},2} 
  \end{equation}
  is an element of $J'(f_E)$.
  \begin{proposition}\label{zerlegt1p} 
  The polarized contraction functor defines an isomorphism of functors on $\Nilp_{O_{\breve{E}}}$,
  \begin{equation}\label{zerlegt4e}
\tilde{\mathcal{M}}_{r} \isoarrow J'/J'(\mathbb{Z}_p).  
    \end{equation}
The Weil  descent datum $\omega_{\mathcal{M}_{r}}$ relative to $O_{\breve{E}}/O_E$
  corresponds on the right hand
  side  to the Weil  descent datum given by multiplication
  with $\pi^{a_E} \in J'(f_E)$.  
\end{proposition}
  We note that $J'/J'(\mathbb{Z}_p) = J'/J'(\mathbb{Z}_p)^{(\tau_E)}$ because this
  is true for any constant sheaf. 
   Proposition \ref{zerlegt1p} is the consequence of the definition
  of $\omega_{\mathcal{M}_{r}}$ and the following Lemma. 

  \begin{lemma}
    The Frobenius
    $F_{\mathbb{X}, \tau_E}: \mathbb{X} \longrightarrow (\tau_E)_{\ast} \mathbb{X}$
    induces on $C_{\mathbb{X}}$ the multiplication by $\pi^{a_E}$. 
  \end{lemma}
  \begin{proof}
    The statement needs an explanation. 
    Because the functor $\mathfrak{C}_{r, \bar{\kappa}_E}^{\rm pol}$ commutes with base change, we have a canonical
    isomorphism $C_{\mathbb{X}} = C_{(\tau_E)_{\ast} \mathbb{X}}$. Indeed, the inverse
    image of the constant sheaf $C_{\mathbb{X}}$ by $\Spec \tau_E$ is the
    constant sheaf $C_{\mathbb{X}}$.

    Let $M = P_{\mathbb{X}}$ be the Dieudonn\'e module of $\mathbb{X}$. The
    Frobenius $F_{\mathbb{X}, \tau_E}$ is induced by the Verschiebung
    \begin{displaymath}
V^{f_E}: M \longrightarrow M.  
    \end{displaymath}
    We write in this proof $C := C_{\mathbb{X}}$. By definition we have
    \begin{displaymath}
 C = M^{\pi_r V^{-1}} =C_1 \oplus C_2,
      \end{displaymath}
       cf. Remark \ref{CPexplban}.
    Therefore the action of $V^{f_E}$ on $C$ coincides with the action of
\begin{equation}\label{zerlegt5e} 
    ~^{F^{-f_E+1}}\pi_r\cdot \ldots \cdot ~^{F^{-1}}\pi_r \cdot \pi_r:
    C \longrightarrow C.  
  \end{equation}
We look at the components of the element on the left hand side in 
(\ref{zerlegt1e}). Let us consider the components of the first set of
factors of (\ref{zerlegt1e}) which act on $C_1$. The component of
\eqref{zerlegt5e} in the factor $O_F \otimes_{O_{F^t}, \tilde{\theta}} W(O_{E'})$
is
\begin{displaymath}
  \pi^{a_{\theta_1\sigma^{f_E-1}}} \cdot \ldots \cdot \pi^{a_{\theta_1\sigma}} \cdot
  \pi^{a_{\theta_1}} \otimes 1 = \pi^{a_{1,E}} \otimes 1. 
  \end{displaymath}
Therefore $V^{f_E}$ induces on $C_1$ the multiplication by $\pi^{a_{1,E}}$. 
By the same argument it induces on $C_2$ the multiplication by $\pi^{a_{2,E}}$. 
  \end{proof}
\begin{corollary}\label{J'split} Let $J(\BQ_p)=\GL_K(C_{\BX, 1}\otimes_{\BZ_p}\BQ_p)$  and  $J(\BZ_p)=\GL_{O_K}(C_{\BX, 1})$.
  There are isomorphisms of functors on $\Nilp_{O_{\breve E}}$, 
  \begin{displaymath}
\tilde{\mathcal{M}}_r \isoarrow J(\mathbb{Q}_p)/J(\mathbb{Z}_p), \, \, \quad \CM_{K/F, r, 1} \isoarrow J(\mathbb{Q}_p)^o/J(\mathbb{Z}_p) .
    \end{displaymath}
    Here $J(\mathbb{Q}_p)^o$ denotes the subgroup of elements with determinant in $O_K^\times$. \qed
    \end{corollary}

        \section{Application to  $p$-adic uniformization}\label{s:attpu}

        In this section, we reap the global fruits from our local work in the preceding sections. This section is modelled on the case of {\it $p$-adic uniformization of the first kind}  of the previous paper \cite[section 4]{KRnew}.

\subsection{The Shimura variety and its $p$-integral model}\label{ss:tsv} 

In this section, $K$ and $F$ will be number fields.         
Let $K/F$ be a CM-field. We fix an archimedean place $w_0$ of $F$. We denote by $a \mapsto \bar{a}$ the
conjugation acting on $K$.

Let $V$ be a $K$-vector space of dimension $2$. Let\index[NO]{ZZQA@$\varsigma$}
\begin{displaymath}
  \varsigma (\; , \;): V \times V \longrightarrow \mathbb{Q} 
\end{displaymath}
be a non-degenerate alternating $\mathbb{Q}$-bilinear form such that
\begin{displaymath}
  \varsigma (a v_1, v_2) = \varsigma (v_1, \bar{a}v_2), \quad a \in K, \;
  v_1, v_2 \in V. 
\end{displaymath}

We define three algebraic groups over $\mathbb{Q}$. For a $\mathbb{Q}$-algebra
$R$, the $R$-valued points are:
\begin{equation}\label{Unif3e}
  \begin{aligned}
  U(V, \varsigma)(R) &= \{g \in \GL_{K \otimes R} (V \otimes R) 
  \mid \varsigma_R(gx_1, gx_2) = \varsigma_R(x_1, x_2)\}\\
    G(V, \varsigma)(R) &= \{ g \in \GL_{K \otimes R}(V \otimes R) \mid
    \varsigma_R (gv_1, gv_2) = \mu(g) \varsigma_R(v_1, v_2),
    \; \mu(g) \in R^{\times}\},\\
    \dot{G}(V, \varsigma)(R)& = \{ g \in \GL_{K \otimes R}(V \otimes R) \mid
    \varsigma_R(gv_1, gv_2) = \varsigma_R (\mu(g)v_1, v_2), \;
    \mu(g) \in (F \otimes R)^{\times}\}.
  \end{aligned}
\end{equation}
If $(V, \varsigma)$ is fixed, we write $U, G, \dot{G}$. 

Equivalently, we can replace the form $\varsigma$ by the anti-hermitian
form \index[NO]{ZZJA@$\varkappa$}
\begin{displaymath}
  \varkappa: V \; \times \; V \longrightarrow K 
\end{displaymath}
on the $K$-vector space $V$ which is 
defined by the equation
\begin{equation}\label{KneunM2e}
  \Trace_{K/\mathbb{Q}} a \varkappa(v_1, v_2) = \varsigma(av_1, v_2),
  \quad a \in K.  
\end{equation}
Then $\varkappa$ is linear in the first argument and anti-linear in the
second.

For each place $w$ of $F$ we obtain an anti-hermitian pairing
\begin{equation}\label{KneunM1e} 
  \varkappa_w: V \otimes_{F} F_w \; \times \; V \otimes_{F} F_w
  \longrightarrow K \otimes_{F} F_w. 
\end{equation}
Let $w: F \longrightarrow \mathbb{R}$ be an archimedean place. We choose an
extension of $w$ to $\varphi: K \longrightarrow \mathbb{C}$. This defines an
isomorphism $K \otimes_{F} F_w \cong \mathbb{C}$ and  
$\varkappa_{w}$ becomes an anti-hermitian pairing
\begin{displaymath}
  \varkappa_{\varphi}: V \otimes_{K, \varphi} \mathbb{C} \; \times \; 
  V \otimes_{K, \varphi} \mathbb{C} \longrightarrow \mathbb{C}.  
\end{displaymath}

Note that the space $V$ is determined  up to isomorphism by the signature at
the archimedean places $w$ of $F$ and the local invariants
$\inv_v(V) := \inv (V \otimes_F F_{v}, \varkappa \otimes_F F_{v})$ at the
non-archimedean places $v$ of $F$, cf. Definition \ref{invdet}. If $v$ splits
in $K$, we set $\inv_v(V) = 1$.  
We will impose the following signature condition on $V$. Let
\index[NO]{ZZSA@$\Phi$} $\Phi = \Hom_{\text{$\mathbb{Q}$-Alg}}(K, \mathbb{C})$
and let \index{special CM-type} 
$r$ be a special CM-type of rank 2 wrt. $w_0$, i.e., a function 
\begin{equation}
r: \Phi\lra \BZ_{\geqslant 0}, \qquad \,\, \varphi \mapsto r_\varphi,
\end{equation}
such that $r_{\varphi}+r_{\bar{\varphi}} =2$ for all $\varphi\in \Phi$ and such
$r_{\varphi} = 1$ iff $\varphi$ restricts to $w_0: F \rightarrow \mathbb{R}$.
We write ${\varphi_0}, {\bar{\varphi}_0} \in \Phi$ for the two extensions of
$w_0$. 

We require that $ \varkappa_{\varphi}$ is isomorphic to the anti-hermitian form on
$\mathbb{C}^2$ given by the matrix
\begin{equation}\label{KneunM3e}
  \left(
  \begin{array}{cc} 
    \mathbf{i} E_{r_{\varphi}} & \mathbf{0}\\
    \mathbf{0} & -\mathbf{i} E_{r_{\bar{\varphi}}}
  \end{array}
  \right),
\end{equation}
for every ${\varphi} \Phi$. Here $E_{r_{\varphi}}$ is denotes the unit matrix of
size $r_\varphi$, and $\mathbf{i}$ the imaginary unit. We note that the last
requirement is independent of the choice of $\varphi$ above $w$. 
We endow
\begin{displaymath}
  V \otimes_{\mathbb{Q}} \mathbb{R} = \prod_{w: F \longrightarrow
    \mathbb{R}} V \otimes_{F, w} \mathbb{R},
\end{displaymath}
with the complex structure $\mathcal{J}$ such that
$\varkappa_{w}(v_1, \mathcal{J}v_2)$ is  hermitian and positive-definite for
all $w$. This defines a Shimura datum $(G,h)$, resp., $(\dot{G}, h)$,
cf. \cite[4.9, 4.13]{De} and an associated Shimura variety ${\rm Sh}(G, h)$
with canonical model over  the reflex field $E$ of $r$.

Let $p$ be a prime number. We also impose a condition on $(V, \varsigma)$ at a
$p$-adic place of $F$. The condition
is with respect to a chosen embedding $E \rightarrow \bar{\mathbb{Q}}_{p}$. 
Let $\nu$ be the place of $E$ defined by this embedding. We denote by $v$
\index[NO]{VAB@$v$}\index[NO]{ZZMA@$\nu$} 
the $p$-adic place of $F$ defined by 
\begin{displaymath}
  F \overset{\varphi_0}{\longrightarrow} E\longrightarrow \bar{\mathbb{Q}}_p.   
\end{displaymath}
We require that
\begin{displaymath}
\inv_v(V)=-1. 
\end{displaymath}
In particular, this implies that the place $v$ does not split in $K$. 
We will denote by $\mathfrak{p}_{\nu}$ the prime ideal of $O_E$ which corresponds
to $\nu$, and by $\mathfrak{p}_v$ the prime ideal of $O_F$ which corresponds
to $v$.  

We have an isomorphism
\begin{equation}\label{KneunM4e} 
  V \otimes \mathbb{Q}_p \cong \oplus_{\mathfrak{p}|p} V \otimes_{F} F_{\mathfrak{p}}.
\end{equation}
Here $\mathfrak{p}$ runs over all prime ideals of $F$ which divide $p$. We will write
\begin{displaymath}
  V_{\mathfrak{p}} := V \otimes_{F} F_{\mathfrak{p}}, \quad
  K_{\mathfrak{p}} = K \otimes_{F} F_{\mathfrak{p}}. 
\end{displaymath}
The decomposition (\ref{KneunM4e}) is orthogonal with respect to $\varsigma$ 
and, for
each prime ideal $\mathfrak{p}$ of $F$ over $p$, we obtain a bilinear form
\begin{displaymath}
  \varsigma_{\mathfrak{p}}: V_{\mathfrak{p}} \times V_{\mathfrak{p}}
  \longrightarrow \mathbb{Q}_p.     
\end{displaymath}
It is related to $\varkappa_{\mathfrak{p}}$ by
\begin{displaymath}
  \Trace_{K \otimes_{F} F_{\mathfrak{p}}/\mathbb{Q}_p} a \varkappa_{\mathfrak{p}}(x_1, x_2) =
  \varsigma_{\mathfrak{p}}(ax_1, x_2),
  \quad a \in K \otimes_{F}  F_{\mathfrak{p}}, \;
  x_1, x_2 \in V_{\mathfrak{p}}. 
\end{displaymath}
One defines algebraic groups over $\mathbb{Q}_p$ as in (\ref{Unif3e}) above:
\index[NO]{UAA@$U_{\mathfrak{p}}$} \index[NO]{GAA@$G_{\mathfrak{p}}$}
\index[NO]{GAB@$\dot{G}_{\mathfrak{p}}$} 
\begin{displaymath}
  U_{\mathfrak{p}} = U(V_{\mathfrak{p}}, \varsigma_{\mathfrak{p}}), \quad 
   G_{\mathfrak{p}} = G(V_{\mathfrak{p}}, \varsigma_{\mathfrak{p}}), \quad 
  \dot{G}_{\mathfrak{p}} = \dot{G}(V_{\mathfrak{p}}, \varsigma_{\mathfrak{p}}). 
  \end{displaymath}

Let $\mathfrak{p}|p$ be  such that $\mathfrak{p}$ is
unramified (and hence non-split) in $K/F$ and such that
$\inv(V_{\mathfrak{p}}, \varkappa_{\mathfrak{p}}) = -1$, cf. Definition \ref{invdet}. Then\index[NO]{ZZKB@$\Lambda_{\mathfrak{p}}$}
there is a $O_{K_{\mathfrak{p}}}$-lattice
\begin{displaymath}
  \Lambda_{\mathfrak{p}} \subset V_{\mathfrak{p}}
\end{displaymath}
such that $\varsigma_{\mathfrak{p}}$ induces a bilinear form
\begin{equation}\label{Unif2e}
  \varsigma_{\mathfrak{p}}:  \Lambda_{\mathfrak{p}} \times
  \Lambda_{\mathfrak{p}} \longrightarrow \mathbb{Z}_p  
\end{equation}
such that $\Lambda_{\mathfrak p}$ is almost self-dual, i.e., ${\tt h}(\Lambda_{\mathfrak{p}}, \varsigma_{\mathfrak{p}}) = 1$ (compare
(\ref{Kneunpsi3e})). Any other lattice with these properties has the
form $g \Lambda_{\mathfrak{p}}$ where
$g \in \U(V_{\mathfrak{p}}, \varsigma_{\mathfrak{p}})(\mathbb{Q}_p)$. 
This follows from Lemma \ref{Kneun0l}.

In all other cases, we apply the following lemma.

\begin{lemma}
  Let $\mathfrak{p}$ be a $p$-adic place of $F$. Assume that 
 $\inv(V_{\mathfrak{p}}, \varkappa_{\mathfrak{p}}) = 1$ if $\mathfrak{p}$
  is unramified in $K/F$.   
There is an $O_{K_{\mathfrak{p}}}$-lattice
$\Lambda_{\mathfrak{p}} \subset V \otimes_F F_{\mathfrak{p}}$ such that
$\varsigma_{\mathfrak{p}}$ induces a perfect pairing
\begin{displaymath}
  \varsigma_{\mathfrak{p}}:  \Lambda_{\mathfrak{p}} \times
  \Lambda_{\mathfrak{p}} \longrightarrow \mathbb{Z}_p.  
\end{displaymath}
Any other such lattice  is of the form 
$g \Lambda_{\mathfrak{p}}$ where
$g \in U_{\mathfrak{p}}(\BQ_p)$. 
\end{lemma}

\begin{proof} Indeed, because of  Lemmas \ref{Kneun0l} and \ref{Kneun01l},
we need only a justification in the case where $\mathfrak{p}$ is split.
In this case we have
$K_{\mathfrak{p}} = F_{\mathfrak{p}} \times F_{\mathfrak{p}}$ and,
accordingly, a decomposition
$V \otimes_{F} F_{\mathfrak{p}} = U_1 \oplus U_2$. The vector spaces
$U_1$ and $U_2$ are isotropic with respect to $\varsigma_{\mathfrak{p}}$
and therefore $\varsigma_{\mathfrak{p}}$ induces an isomorphism
$U_2 = \Hom_{\mathbb{Q}_p}(U_1,\mathbb{Q}_p)$. The form
$\varsigma_{\mathfrak{p}}$ becomes 
\begin{displaymath}
  \varsigma_{\mathfrak{p}} (x + x^{\ast}, y + y^{\ast}) = x^{\ast}(y) -
  y^{\ast}(x), \quad x,y \in U_1, \; x^{\ast}, y^{\ast} \in U_2.  
\end{displaymath}
The existence and uniqueness of $\Lambda_{\mathfrak{p}}$
follows easily. 
\end{proof}
To pass to a $p$-integral model  over
$O_{E, ({\mathfrak{p}_{\nu}})}$ of ${\rm Sh}(G, h)$, we restrict the choice of the level structure. To do this, we  choose for each $\mathfrak{p}|p$ a
$O_{K \otimes_F F_{\mathfrak{p}}}$-lattice
$\Lambda_{\mathfrak{p}} \subset V_{\mathfrak{p}}$ as above.
\index[NO]{KAE@$\mathbf{K}_{\mathfrak{p}}$}\index[NO]{KAF@$\mathbf{K}_{p}$}
We define 
\begin{equation}\label{boldK1e}
  \begin{aligned}
    \bK_{\mathfrak{p}} &= \{g \in G_{\mathfrak{p}}  \; | \; g\Lambda_{\mathfrak{p}}
    = \Lambda_{\mathfrak{p}} \}\\
    {\bK}_{p} &= \{g \in G(\mathbb{Q}_{p})  \; | \; g\Lambda_{\mathfrak{p}}
    = \Lambda_{\mathfrak{p}} \; \text{for all} \; \mathfrak{p} | p\}. 
  \end{aligned}
\end{equation}
We choose an open compact subgroup ${\bK}^p \in G(\mathbb{A}_f^{p})$ and set
${\bK} = {\bK}_p \cdot {\bK}^p \subset G(\mathbb{A}_f)$.

We extend the embedding $\nu\colon E \longrightarrow \bar{\mathbb{Q}}_p$ to an embedding
$\bar{\mathbb{Q}} \longrightarrow \bar{\mathbb{Q}}_p$. We obtain a decomposition
\begin{equation}\label{decPhi}
  \Phi = \Hom_{\text{$\mathbb{Q}$-Alg}}(K, \bar{\mathbb{Q}}) =
  \coprod_{\mathfrak{p}} \Hom_{\text{$\mathbb{Q}_p$-Alg}}(K_{\mathfrak{p}}, \bar{\mathbb{Q}}_p) =
  \coprod_{\mathfrak{p}} \Phi_{\mathfrak{p}}. 
  \end{equation}
The restriction of the function $r$ to $\Phi_{\mathfrak{p}}$ will be denoted by
$r_{\mathfrak{p}}$. \index[NO]{RAD@$r_{\mathfrak{p}}$}
The group
$\Gal(\bar{\mathbb{Q}}_p/E_{\nu})$ acts on $\Phi$ via the restriction
\begin{displaymath}
  \Gal(\bar{\mathbb{Q}}_p/E_{\nu}) \longrightarrow \Gal(\bar{\mathbb{Q}}/E).
  \end{displaymath}
Therefore $r_{\tau \varphi} = r_{\varphi}$ for $\varphi \in \Phi$ and
$\tau \in \Gal(\bar{\mathbb{Q}}_p/E_{\nu})$. This implies that the local
reflex fields $E(K_{\mathfrak{p}}/F_{\mathfrak{p}}, r_{\mathfrak{p}})$ are all contained
in $E_{\nu}$.

Let $R$ be an $O_{E_{\nu}}$-algebra. Let $\CL$ be an $R$-module with an
$O_K$-action. We have decompositions
\begin{displaymath}
  O_K \otimes R = \prod_{\mathfrak{p}} O_{K_{\mathfrak{p}}} \otimes_{\mathbb{Z}_p} R, \quad
  \CL = \oplus_{\mathfrak{p}} \CL_{\mathfrak{p}}.
\end{displaymath}
We will say that $\CL$ satisfies the Eisenstein condition $({\rm EC}_r)$
if each $\CL_{\mathfrak{p}}$ satisfies the Eisenstein condition
    $({\rm EC}_{r_{\mathfrak{p}}})$, cf. \eqref{Eisen}. We use a similar terminology for the Kottwitz condition $({\rm KC}_r)$.
\begin{definition}\label{KneunA1d}
  We define the groupoid $\mathcal{A}_{\bK}$\index[NO]{ABA@$\mathcal{A}_{\bK}$}
  on the category of $O_{E, ({\mathfrak{p}_\nu})}$-algebras.
  A point of $\mathcal{A}_{\bK}(R)$ consists  of the  following
  data:

\begin{enumerate}
\item[(1)] An abelian scheme $A$ over $R$, up to isogeny of degree
  prime to $p$, with an algebra homomorphism
  \begin{displaymath}
    \iota: O_K \longrightarrow \End A \otimes \mathbb{Z}_{(p)}. 
  \end{displaymath}
  such that $\Lie A\otimes_{O_{E, ({\mathfrak{p}_\nu})}}O_{E_{\nu}}$
  satisfies the conditions
  $({\rm KC}_r)$ and  $({\rm EC}_r)$.
\item[(2)] A $\mathbb{Q}$-homogeneous polarization $\bar{\lambda}$ of $A$
  such that the Rosati involution induces the conjugation of $K/F$.
\item[(3)] A class of $O_K$-linear isomorphisms
  \begin{displaymath}
    \bar{\eta}^p: V \otimes \mathbb{A}_f^{p} \longrightarrow V^p(A)  \;
    \text{mod} \; {{\bK}}^{p}  
  \end{displaymath}
 which respect the forms on both sides up to a constant in
  $\mathbb{A}_f^p(1)^{\times}$. 
  \end{enumerate}
  We impose the following two conditions.
  \begin{enumerate}
\item[(i)] There exists a polarization
$\lambda \in \bar{\lambda}$ such that the induced map to the dual variety
$\lambda: A \rightarrow A^{\wedge}$ has the following property. Let $(\ker \lambda)_p$ be
the $p$-primary part of the kernel of $\lambda$. It has the decomposition
$(\ker \lambda)_p = \prod_{\mathfrak{p}|p} (\ker \lambda)_{\mathfrak{p}}$. We require that $(\ker \lambda)_{\mathfrak{p}} $ is trivial, 
unless $\mathfrak{p}$ is unramified in $K$ and
$\inv(V_{\mathfrak{p}}, \varsigma_{\mathfrak{p}}) = -1$. In the latter case the height of $(\ker \lambda)_{\mathfrak{p}}$ is $2f_{\mathfrak{p}}$. 

\item[(ii)] For each geometric point $\omega: R \rightarrow k$  There is an identity of invariants, 
\begin{equation*}\label{invcond}
  \inv_{\mathfrak p}^r(A_{\omega}, \iota, \lambda) 
  = \inv_{\mathfrak p}(V_{\mathfrak p}, \varsigma_{\mathfrak{p}}),
  \quad \text{ for all }\, \mathfrak p|p,  
\end{equation*}
cf. the explanation after this Definition. 

\end{enumerate}
An isomorphism of such data $(A, \iota, \bar{\lambda}, \bar{\eta}^{p})\to (A', \iota', \bar{\lambda'}, \bar{\eta'}^{p})$  is given by a $O_K$-linear quasi-isogeny $\phi\colon A\to A'$ of degree prime to $p$ compatible with the $\BQ$-homogeneous polarizations and the level structures. 
\end{definition}
We briefly explain the requirement (ii), with references to the
Appendix. Assume that the characteristic of $\omega$ is $p$. Then the invariant
on the left hand side is defined in terms of the Dieudonn\'e module $M$ of 
$A_{\omega}$. The definitions of section \ref{ss:inv} apply to the Dieudonn\'e
module $M_{\mathfrak{p}} := M \otimes_{F \otimes \mathbb{Q}_p} F_{\mathfrak{p}}$ with the
action of $O_{K_{\mathfrak{p}}}$. The adjusted invariant (\ref{Kneun4e}) of the latter
Dieudonn\'e module is by definition
$\inv_{\mathfrak p}^r(A_{\omega}, \iota, \lambda)$. In the case where the
characteristic of $\omega$ is zero, $\inv_{\mathfrak p}^r(A_{\omega}, \iota, \lambda)$
is the invariant of the $K_{\mathfrak{p}}$-vector space 
$(V_{\mathfrak{p}}(A_{\omega}), \mathcal{E}_{\mathfrak{p}})$ with the Riemann form
induced by $\lambda$, cf. Definition \ref{invdet}. It is easily seen that
the left hand side of (ii) depends only on the image $s \in \Spec R$ of
$\omega$. By  Proposition \ref{contcorr}
we regard $\inv_{\mathfrak p}^r(A, \iota, \lambda)$ as a locally constant function
on $\Spec R$ and write the condition (ii) in the form
\begin{displaymath}
  \inv_{\mathfrak p}^r(A, \iota, \lambda)
  = \inv_{\mathfrak p}(V_{\mathfrak p}, \varsigma_{\mathfrak{p}}). 
  \end{displaymath}  
 
We will denote a point of $\mathcal{A}_{\bK}(R)$ simply by
$(A, \iota, \bar{\lambda}, \bar{\eta}^{p})$.
\begin{remark} It is equivalent to consider in (2) a $\mathbb{Z}_{(p)}$-homogeneous polarization $\bar{\lambda}$ of $A$
  such that the elements $\lambda\in\bar\lambda$ satisfy the condition (i) on the $p$-primary part of the kernel of $\lambda$. 
  \end{remark}

\begin{remark}\label{lattforlambda}
Let $(A, \iota, \bar{\lambda}, \bar{\eta}^{p})\in \CA_\bK(R)$. Let  $\lambda \in \bar{\lambda}$ be  as in condition (i) of Definition \ref{KneunA1d}.   For each
  geometric point $\omega$ of characteristic $0$ of $\Spec R$ the pairing
  induced by $\lambda$   on the $\mathfrak p$-adic Tate module, 
  \begin{equation}\label{KneunM6e2}
    \mathcal{E}_{\mathfrak{p}}:  T_{\mathfrak{p}}(A_{\omega}) \times
    T_{\mathfrak{p}}(A_{\omega}) \longrightarrow \mathbb{Z}_p(1) 
  \end{equation}
  has the following properties. 
  If $\mathfrak{p}$ is ramified in $K/F$, the pairing is perfect and
  $\inv(V_{\mathfrak{p}}(A_{\omega}), \mathcal{E}_{\mathfrak{p}}) =
  \inv(V_{\mathfrak{p}}, \varsigma_{\mathfrak{p}})$. 
  If $\mathfrak{p}$ is unramified, then 
  (\ref{KneunM6e2}) is perfect if
  $\inv(V_{\mathfrak{p}}, \varsigma_{\mathfrak{p}}) = 1$ and is almost perfect if
  $\inv(V_{\mathfrak{p}}, \varsigma_{\mathfrak{p}}) = -1$.
  If $\mathfrak{p}$ is split in $K/F$, then (\ref{KneunM6e2}) is perfect. 

  For each geometric point $\omega$ of $R$ of characteristic $p$, the polarization 
  $\lambda$ induces a pairing on the Dieudonn\'e module $M$ of $A_{\omega}$
  and therefore for each prime $\mathfrak{p} | p$ of $F$ a pairing
  \begin{equation}\label{KneunM8e}
    \mathcal{E}_{\mathfrak{p}}:  M_{\mathfrak{p}}(A_{\omega}) \times
    M_{\mathfrak{p}}(A_{\omega}) \longrightarrow W(\kappa(\omega)) ,  
  \end{equation}
  with  the following properties. 
  If $\mathfrak{p}$ is ramified in $K/F$, the pairing is perfect and
  $\inv^{r}(M_{\mathfrak{p}}(A_{\omega}), \mathcal{E}_{\mathfrak{p}}) = \inv(V_{\mathfrak{p}}, \varsigma_{\mathfrak{p}})$. 
  If $\mathfrak{p}$ is unramified, then 
  (\ref{KneunM8e}) is perfect if
  $\inv(V_{\mathfrak{p}}, \varsigma_{\mathfrak{p}}) = 1$ and is almost perfect if
  $\inv(V_{\mathfrak{p}}, \varsigma_{\mathfrak{p}}) = -1$.
  If $\mathfrak{p}$ is split in $K/F$, then  (\ref{KneunM8e}) is perfect.   
  \end{remark}

\begin{proposition}\label{proprepr}
Assume that ${\bK}^p$ is small enough. 
Then the functor $\mathcal{A}_{\bK}$ is representable by a projective scheme over  $\Spec O_{E, ({\mathfrak{p}_\nu})}$ whose generic
fiber is the Shimura variety $\mathrm{Sh}_{\bK}$ associated to the Shimura datum
$(G, h)$. For general ${\bK}^p$, $\mathcal{A}_{\bK}$ is a DM-stack proper over $\Spec O_{E, ({\mathfrak{p}_\nu})}$ whose generic fiber is the Shimura variety  $\mathrm{Sh}_{\bK}$ considered as the classifying stack of a group action. 
\end{proposition}
\begin{proof}
  Let $(A, \iota, \bar{\lambda}, \bar{\eta}^{p})\in \CA_\bK(R)$, and fix $ \eta^{p} \in \bar{\eta}^{p}$. Let $\Lambda \subset V$ be a $O_K$-lattice on which $\varsigma$ is integral.  We find an abelian variety $A_1$ in
the class $A$ up to isogeny prime to $p$ such that for each $\ell \neq p$
\begin{displaymath}
  \eta^{p} (T_{\ell}(A_1)) = \Lambda \otimes \mathbb{Z}_{\ell}.
\end{displaymath}
In this way we obtain also a polarization on $A_1$ whose degree is bounded
in terms of $\varsigma$ and $\Lambda$. If ${\bK}^p$ is small enough we obtain
a level structure on the $m$-division points for some $m \geq 3$.

The fact that the moduli problem of abelian varieties with
  a polarization of given degree and a $m$-level structure for $m \geq 3$
  is a quasi-projective scheme implies that the functor of $(A, \iota, \bar\lambda, \bar\eta^p)$ as in (1)--(3) of Definition \ref{KneunA1d}, is  representable by a
  quasi-projective scheme. Now the conditions (i) and (ii) define open and closed subschemes (this is easy for condition (i), and follows from  Proposition \ref{contcorr} for condition (ii)). The representability by a Deligne-Mumford stack for general ${\bK}^p$ follows. 
  
  To compare the generic fiber of $\CA_\bK$ with $\mathrm{Sh}_{\bK}$, recall from \cite[\S 8]{K} that $\mathrm{Sh}_{\bK}$ represents the following functor $\CA_{\bK, E}$ on the category of $E$-algebras, comp. section \ref{ss:globalres}. A point of $\CA_{\bK, E}(R)$ consists of the following data.
  
\begin{enumerate}
\item[(1)] An abelian scheme $A$ over $R$, up to isogeny, with an algebra homomorphism
  \begin{displaymath}
    \iota: O_K \longrightarrow \End A \otimes \mathbb{Q} 
  \end{displaymath}
  such that the Kottwitz condition  $({\rm KC}_r)$ is satisfied.
\item[(2)] A $\mathbb{Q}$-homogeneous polarization $\bar{\lambda}$ of $A$
  such that the Rosati involution induces the conjugation of $K/F$.
\item[(3)] A class of $K$-linear isomorphisms
  \begin{displaymath}
    \bar{\eta}: V \otimes \mathbb{A}_f \longrightarrow \hat V(A) \;
    \text{mod} \; {{\bK}}  
  \end{displaymath}
 which respect the forms on both sides up to a constant in
  $\mathbb{A}_f(1)^{\times}$. 
  \end{enumerate}
Here we are implicitly using the fact that $G$ satisfies the Hasse principle, cf. \cite[\S 7]{K}. We define a map $\CA_{\bK, E}(R)\to\CA_\bK(R)$. We fix a $O_K$-lattice $\Lambda$ in $V$ whose localizations at $\mathfrak{p}|p$ are the given lattices $\Lambda_{\mathfrak p}$ above. Let $(A, \iota, \bar{\lambda}, \bar{\eta})\in \CA_{\bK, E}(R)$, and fix $\eta\in\bar{\eta}$. We find an abelian variety $A_1$ in
the isogeny class $A$ such that for each $\ell $
\begin{displaymath}
  \eta (T_{\ell}(A_1)) = \Lambda \otimes \mathbb{Z}_{\ell}.
\end{displaymath}
 Then we obtain $\iota_1\colon O_K\to \End(A_1)\otimes\BZ_{(p)}$. The Eisenstein condition $({\rm EC}_r)$ is automatically satisfied, cf. Proposition \ref{unram}. We also find a polarization $\lambda_1\in\bar\lambda$ which satisfies the condition (i) in Definition \ref{KneunA1d}.  The existence of $\eta$ implies Condition (ii). By forgetting the $p$-component of $\bar\eta$, we have associated to $(A, \iota, \bar\lambda, \bar\eta)$ a well-defined object $(A_1, \iota_1, \bar\lambda_1, \bar\eta^p)$ of $\CA_\bK(R)$. By the uniqueness property of the lattices $\Lambda_{\mathfrak p}$ mentioned above, this map is bijective.
  
  The properness of
  $\mathcal{A}_{\bK} \longrightarrow \Spec O_{E, ({\mathfrak{p}_\nu})}$ is a consequence
  of Proposition \ref{proprop} below.
\end{proof}
\begin{remark}\label{gf1r}
If $\mathfrak{p}_v$ is the only prime ideal of $F$ over $p$, then it follows from the product formula that condition (ii) of Definition \ref{KneunA1d} is 
automatically satisfied. Indeed, condition (ii) defines an open and closed subscheme which in this case has the same generic fiber. 
\end{remark}

\begin{proposition}\label{proprop}
  The morphism $\mathcal{A}_{\bK} \longrightarrow \Spec O_{E, ({\mathfrak{p}_\nu})}$ is
  proper. 
  \end{proposition}
For the proof of Proposition \ref{proprop}, we need two lemmas. 
\begin{lemma}
 Let $K/F$ be a CM-field and let $r$ be a generalized CM-type of
 rank $2$. Let $E \subset \bar{\mathbb{Q}}$ be the reflex field.
 Let $R$ be a complete discrete valuation ring with an 
 $O_E$-algebra structure. Let $L$ be the field of fractions of $R$.
 We assume that the residue characteristic
 of $R$ is $p > 0$. Let $w$ be a finite place of $F$ of residue
  characteristic $\ell$, such that  $K_{w}/F_{w}$ is a field extension.
  We assume that $L$ is of characteristic zero or that $\ell \neq p$. 
  
  Let $(A, \iota, \lambda)$ be a CM-triple of
  type $r$ over $L$. The polarization induces on the rational
  Tate module $V_{w}(A)$ an alternating pairing 
  \begin{equation}\label{K9red1e}
    \psi_{w}: V_{w}(A) \times V_{w}(A) \longrightarrow \mathbb{Q}_{\ell}(1). 
  \end{equation}
If $\inv(V_{w}(A), \psi_{w}) = -1$, then  the abelian variety $A$ has
  potentially good reduction.
\end{lemma}
\begin{proof}
  We consider only the case $\ell = p$. We may assume that $A$ has
  semistable reduction. We choose an isomorphism
  $\mathbb{Q}_{p} \cong \mathbb{Q}_{p}(1)$ over $\bar{L}$. We obtain from
  $\psi_w$ the anti-hermitian form
  \begin{displaymath}
    \varkappa_{w}: V_{w}(A) \times V_{w}(A) \longrightarrow K_{w}, 
  \end{displaymath}
  cf. (\ref{Kneunpsi1e}). 
    Let $T$ be the toric part of the special fiber of the
  N\'eron model of $A$. Then $O_K$ acts on the character group
  $X_{\ast}(T)$. If  $T$ is non-trivial, we obtain that
  \begin{displaymath}
    \dim T = [K:\mathbb{Q}] = \dim A.
  \end{displaymath}
  This implies that the toric part $V^t_{w}(A) \subset V_{w}(A)$ is
  a $K_{w}$-vector subspace of dimension $1$. By the orthogonality
  theorem [SGA7, Exp IX, Thm. 5.2], the anti-hermitian form
  $\varkappa_{w}$ is zero on this subspace. Let $u_1, u_2$ be a basis
  of $V_{w}(A)$ such that $u_1$ is a basis of $V^t_{w}(A)$. Then we
  obtain, in the notation of \eqref{KneunDisc1e}, 
  \begin{displaymath}
    \mathfrak{d}_{K/F} (V_{w}(A), \varkappa_{w}) = - \varkappa_{w}(u_1, u_2)
    \varkappa_{w}(u_2, u_1) = \varkappa_{w}(u_1, u_2)
    \overline{\varkappa_{w}(u_1, u_2)} \equiv 1 
  \end{displaymath}
  modulo $\Nm_{K_{w}/F_{w}} K_{w}^{\times}$. This contradicts the assumption
  $\inv(V_{w}(A), \varkappa_{w}) = -1$. 
\end{proof}

With the notation of the last lemma,  we consider the case where $\ell = p$
and where the characteristic of $L$ is also $p$. The $O_E$-algebra structure
on $R$ factors
\begin{displaymath}
O_E \longrightarrow \kappa_{\nu} \longrightarrow R,  
  \end{displaymath}
where $\kappa_{\nu}$ is the residue field of $E_\nu$.  We fix
a commutative diagram
\begin{displaymath}
   \xymatrix{
E_\nu \ar[r]  & \bar{\mathbb{Q}}_p \\
E \ar[r]\ar[u] & \bar{\mathbb{Q}} .\ar[u] 
}
  \end{displaymath}
Let $w$ be a $p$-adic place of $F$. By the last diagram we can restrict $r$
to a local CM-type $r_w$ of $K_w/F_w$. Then $E_\nu$ is the composite of the subfields $E(K_w/F_w, r_w)$, for $w$ running over all places of $F$ over $p$. 

Let $(A, \iota, \lambda)$ be a CM-triple of type $r$ over $L$. The action
of $O_F \otimes \mathbb{Z}_p = \prod_{w} O_{F_w}$, where $w$ runs over all
$p$-adic places of $F$, induces a decomposition of the $p$-divisible group
of $A$:
\begin{displaymath}
X = \prod_{w} X_w. 
  \end{displaymath}
Then $(X_w, \iota_w, \lambda_w)$ is a local CM-triple with respect to
$K_w/F_w, r_w$ over the field $L$. 

\begin{lemma}\label{Rofcharp}
  Let $R$ be a discrete valuation ring of equal characteristic
  $p > 0$,  and let $L$ be the field of fractions. Let $R$ be an
  $O_E$-algebra. Let $\nu$ be the $p$-adic place of $E$ induced by this
  algebra structure.

  Let $(A, \iota, \lambda)$ be a CM-triple of type $r$ over $L$.
  We assume that there is a $p$-adic place $w$ of $F$ such that one
  of the following conditions is satisfied.
  \begin{enumerate}
  \item $K_w/F_w$ is a ramified field extension. The local CM-triple
    $(X_w, \iota_w, \lambda_w)$ satisfies the Eisenstein condition
    $({\rm EC}_{r_w})$,  and
    $\inv^r (X_w, \iota_w, \lambda_w) = -1$.
    \item $K_w/F_w$ is an unramified field extension. The local CM-triple
    $(X_w, \iota_w, \lambda_w)$ satisfies the Eisenstein condition
    $({\rm EC}_{r_w})$, and $\lambda_w$ is almost principal.
    \end{enumerate}

Then the abelian variety $A$ has potentially good reduction over $R$. 
\end{lemma}

\begin{proof}
  We may assume that $A$ has semistable reduction over $R$. Let $\tilde{A}$
  be the N\'eron model over $R$, and let $B$ be the identity component of
  the special fibre of $\tilde{A}$. Let us assume that the torus part
  $T \subset B$ is nontrivial. Since $O_K$ acts on $T$, we obtain that
  $X_{\ast}(T)_{\mathbb{Q}}$ is a $K$-vector space of dimension one. Let
  $Y$ be the $p$-divisible group of $T$. We obtain a decomposition
  \begin{displaymath}
Y = \prod_u Y_u
  \end{displaymath}
  where $u$ runs over all places of $K$ over $p$ and $Y_u$ is an
  $O_{K_u}$-module which is of height $[K_u:\mathbb{Q}_p]$ and of multiplicative type.
  We pass from $R$ to the completion $\hat{R}$. Let $\hat{X} = X_{\hat{K}}$
  be the $p$-divisible group of $A_{\hat{K}}$. By [SGA7, Exp IX, \S 5], the
  multiplicative group $Y_w$ lifts to a multiplicative group
  $\tilde{Y}_w \subset \hat{X}^f_w$ of the finite part of $\hat{X}_w$ over
  $\hat{R}$. If we pass to the general fibre of the last inclusion we obtain
  a nontrivial multiplicative subgroup
  $(\tilde{Y}_w)_{\hat{L}} \subset \hat{X}_w$. But our assumption implies, by
  Lemma  \ref{ramFrame1l} in the ramified case,  and by Proposition
  \ref{unrFrame1p} in the unramified case, that $\hat{X}_w$ is isoclinic of
  slope $1/2$. This contradicts the existence of a nontrivial multiplicative
  part and therefore the assumption that the torus part of $B$ is nontrivial.
  \end{proof}

\begin{proof}(of Proposition \ref{proprop})
  We check the valuative criterion. Let $R$ a discrete valuation ring with a $O_{E_{\nu}}$-algebra structure.
  Let $L$ be the field of fractions of $R$. Let
  $\alpha: \Spec L \longrightarrow \mathcal{A}_{\bK}$ be a  $O_{E_{\nu}}$-morphism.
  We have to show that $\alpha$ extends to $\Spec R \longrightarrow \mathcal{A}_{\bK}$.
  It is enough to show that  for a discrete valuation ring $R'$ which dominates $R$,
  the morphism $\Spec L' \longrightarrow \mathcal{A}_{\bK}$ induced by $\alpha$
  extends $\Spec R' \longrightarrow \mathcal{A}_{\bK} $. 
  The map $\alpha$ gives a point
  $(A, \iota, \bar{\lambda}, \bar{\eta}^{p}) \in \mathcal{A}_{\bK}(L)$.
  Since we may replace $L$ by $L'$ we may assume that $A$ has semistable
  reduction.  Let $\omega$ be a geometric point concentrated in the generic point of $\Spec R$. We are assuming that
$\inv^{r}(T_{\mathfrak{p_{v}}}(A_{\omega}), \mathcal{E}_{\mathfrak{p_{v}}}) =
  \inv(V_{\mathfrak{p_{v}}}, \varsigma_{\mathfrak{p}_{v}}) = -1$ when ${\rm char}\, L=0$, resp.  $\inv^{r}(M_{\mathfrak{p_{v}}}(A_{\omega}), \mathcal{E}_{\mathfrak{p_{v}}}) =
  \inv(V_{\mathfrak{p_{v}}}, \varsigma_{\mathfrak{p}_{v}}) = -1$, when ${\rm char}\, L=p$. 
  Hence we conclude by the
  last two lemmas that $A$ has good reduction. Let $\tilde{A}/R$ denote  the abelian
  scheme which extends $A$. Then $\iota$ extends to an action $\tilde{\iota}$
  of $O_K$ on $\tilde{A}$. The Kottwitz condition and the Eisenstein
  condition are closed conditions and hold therefore for $\tilde{A}$.
  The polarization $\lambda$ extends to
  $\tilde{\lambda}: \tilde{A} \longrightarrow \tilde{A}^{\wedge}$. The condition (i) from Definition \ref{KneunA1d} extends from $A$ to $\tilde A$. 
  For a geometric point $\omega_0$ concentrated in the closed point of
  $\Spec R$ we find 
  \begin{displaymath}
    \inv^r_{\mathfrak{p}} (\tilde{A}_{\omega_0}, \tilde{\iota}_{\omega_0},
    \lambda_{\omega_0}) = 
    \inv^{r}(M_{\mathfrak{p}}(\tilde{A}_{\omega_0}), \mathcal{E}_{\mathfrak{p}}) =
    \inv(V_{\mathfrak{p}}, \varsigma_{\mathfrak{p}_{v}})
  \end{displaymath}
  because the left hand is by Proposition \ref{contcorr} equal to
 $\inv^r_{\mathfrak{p}} (\tilde{A}_{{L}}, \iota_{{L}}, \lambda_{{L}})$. Hence condition (ii) from Definition \ref{KneunA1d} also extends from $A$ to $\tilde A$. 
  From this we obtain an extension of
  $(A, \iota, \bar{\lambda}, \bar{\eta}^{p})$ to a point 
of $ \mathcal{A}_{\bK}(R)$.
\end{proof}
\begin{remark}
The scheme $\CA_\bK$ turns out to be flat over $\Spec O_{E, { (\mathfrak p_{v})}}$, cf. Theorem  \ref{maint}, (i).  
Hence its  generic fiber is dense.  It follows  that it is enough to check  the valuative criterion on discrete valuation rings $R$ with fraction field $L$
of characteristic zero. Hence Lemma \ref{Rofcharp}  is not needed. 
\end{remark}
The following proposition shows that there is only one isogeny class in the special fiber of $\mathcal{A}_{\bK}$. This is the underlying reason why there is $p$-adic uniformization.

\begin{proposition}\label{KneunUni3p}
  Let $\kappa_{\nu}$ be the residue class field of $E_{\nu}$. Let
  $(A_1, \iota_1, \bar{\lambda}_1, \bar{\eta}_1^{p})$ and
  $(A_2, \iota_2, \bar{\lambda}_2, \bar{\eta}_2^{p})$
  be two points of $\mathcal{A}_{\bK}(\bar{\kappa}_{\nu})$. Then there exists a quasi-isogeny
  \begin{displaymath}
    (A_1, \iota_1, \bar{\lambda}_1)  \longrightarrow
    (A_2, \iota_2, \bar{\lambda}_2) ,
  \end{displaymath}
  i.e.,  a quasi-isogeny which respects the actions $\iota_i$ and the
  $\mathbb{Q}$-homogeneous polarizations $\bar{\lambda}_i$. In fact, there exists such a quasi-isogeny of degree prime to $p$. 
\end{proposition} 
\begin{proof} Let $X_i$ be the $p$-divisible group of $A_i$, with its decomposition $X_i=\prod\nolimits_{\mathfrak{p}|p}X_{i, \mathfrak p}$. 
  It follows from
Proposition \ref{ramFrame1p} (jointly with Lemma \ref{ramFrame1l}) and Proposition \ref{unrFrame2p} that 
$X_{i,\mathfrak{p}_{v}}$ is isoclinic. In the banal
cases $\mathfrak{p} \neq \mathfrak{p}_{v}$, the same  follows  from Lemma \ref{C'onslopes1c} for $X_{i,\mathfrak{p}}$.  By \cite[Cor. 6.29]{RZ} we find a quasi-isogeny
  \begin{displaymath}
    a: (A_1, \iota_1) \longrightarrow (A_2, \iota_2).
  \end{displaymath}
  We choose $\lambda_i \in \bar{\lambda}_i$. We set
  $\lambda = a^{\ast} (\lambda_2)$. We find an endomorphism $u\in \End^o(A_1)$ 
   such that
  \begin{displaymath}
    \lambda = \lambda_1 u.
  \end{displaymath}
  Since $\lambda_1$ and $\lambda$ induce the conjugation on $K$, we
  conclude that $u \in \End^o_{K} A_1$. Moreover $u$  is
  fixed by the Rosati involution $\ast$ induced by $\lambda_1$ on
  $D := \End^{o}_{K} A_1$. It is enough to find an element $d \in
D^{\times}$ such that
\begin{equation}\label{KneunU1e}
  u = f d^{\ast} d
\end{equation}
for some element $f \in \mathbb{Q}^{\times}$. The solutions of these equations
form a torsor under the algebraic group $J$ over $\mathbb{Q}$
\index[NO]{JAA@$J$} such that 
\begin{equation}\label{KneunU2e}
  J(\mathbb{Q}) = \{ e \in D^{\times} \; | \; e^{\ast}e \in \mathbb{Q}^{\times} \}. 
\end{equation}
 By \cite[\S7]{K}, this group satisfies the Hasse principle.
Therefore it is enough to find a solution of the equation
(\ref{KneunU1e}) in $D \otimes \mathbb{Q}_w^{\times}$ for all places
$w$ of $\mathbb{Q}$. If $w$ is a finite place $w \neq p$ we have, by
\cite[Cor. 6.29]{RZ}, that
\begin{displaymath}
  D \otimes \mathbb{Q}_w = \End_{K \otimes \mathbb{Q}_w} V_w(A_1)
\end{displaymath}
such that the Riemann form $\mathcal{E}^{\lambda_1}_w$ induces the
involution $\ast$. A solution of (\ref{KneunU1e}) exists iff the
symplectic $K \otimes \mathbb{Q}_w$-modules
\begin{displaymath}
  (V_w(A_1), \mathcal{E}^{\lambda_1}_{w}), \quad  \quad
  (V_w(A_2), \mathcal{E}^{\lambda_2}_{w})
\end{displaymath}
are similar up to a factor in $\mathbb{Q}_w^{\times}$. But this follows from the
existence of $\bar{\eta}^p_1$ and $\bar{\eta}^p_2$.

In the case $w = p$ we can use Dieudonn\'e modules. In this case we
 know, by condition (ii) in Definition \ref{KneunA1d}, that the rational Dieudonn\'e modules of $A_1$ and $A_2$ together
with their polarizations are isomorphic. 

If $w$ is the infinite place, one can deduce the assertion from the fact
that $u$ in (\ref{KneunU1e}) is totally positive. 

Now let us prove the second assertion. 
  We consider the Dieudonn\'e modules $M_1$, resp. $M_2$, of $A_1$, resp. $A_2$.
  We choose the polarizations $\lambda_1 \in \bar{\lambda}_1$, resp.
  $\lambda_2 \in \bar{\lambda}_2$, as in condition (i)  of Definition
  \ref{KneunA1d}. Using the contracting functor,
  it is clear that there is a quasi-isogeny of height zero
  $\alpha: (M_1, \lambda_1) \longrightarrow (M_2, \lambda_2)$.
Let $\rho\colon (M_1, \bar{\lambda}_1) \longrightarrow
    (M_2, \bar{\lambda}_2)$ be an arbitrary quasi-isogeny. Consider the morphism
  \begin{displaymath}
    \alpha \circ \rho^{-1}: (M_2, \bar{\lambda}_2) \longrightarrow
    (M_2, \bar{\lambda}_2).  
  \end{displaymath}
  We consider the group $J$ for $(A_2, \iota_2, \bar{\lambda}_2)$
  (compare (\ref{KneunU2e})). Then $\alpha \circ \rho^{-1}$ is an
  element of $J(\mathbb{Q}_p)$ by Tate's theorem (\cite[Cor. 6.29]{RZ}).
  We approximate it by an element $\alpha_1 \in J(\mathbb{Q})$. Then
  \begin{displaymath}
    \rho \circ \alpha_1: (A_1, \iota_1, \bar{\lambda}_1)  \longrightarrow
    (A_2, \iota_2, \bar{\lambda}_2)
  \end{displaymath}
  is the desired quasi-isogeny of order prime to $p$. 
\end{proof}

\subsection{The RZ-space $\tilde\CM_r$}\label{RZtildeM}
We fix a point
$(A_0, \iota_0, \bar{\lambda}_0, \bar{\eta}^{p}_0)$ of $ \mathcal{A}_{\bK}(\bar{\kappa}_{{\nu}})$. 
 We also fix a polarization
$\lambda_0 \in \bar{\lambda}_0$ which satisfies the
condition (i) of Definition \ref{KneunA1d}. We denote by $\mathbb{X}$ the $p$-divisible group of $A_0$. The action $\iota_0$
induces an action $\iota_{\mathbb{X}}$ on $\mathbb{X}$ and $\lambda_0$
induces a polarization $\lambda_{\mathbb{X}}$ on $\mathbb{X}$. We denote by
$q_{\nu}$ the number of elements in $\kappa_\nu=\kappa_{E_{\nu}}$.

Let $R \in \Nilp_{O_{E_{\nu}}}$ and let $(X, \iota)$ be a $p$-divisible
group over $\Spec R$ with an action
\begin{displaymath}
  \iota: O_K \otimes \mathbb{Z}_p \longrightarrow \End X.  
\end{displaymath}
The notion of a \emph{semi-local CM-triple $(X, \iota, \lambda)$ \index{semi-local CM-triple} relative to
$K \otimes \mathbb{Q}_p/ F \otimes \mathbb{Q}_p$ and $r$} should be obvious but 
we explain it more precisely: the decomposition
\begin{displaymath}
  O_F \otimes \mathbb{Z}_p = \prod_{\mathfrak{p}} O_{F_{\mathfrak{p}}}
\end{displaymath}
induces the decomposition
\begin{displaymath}
  X = \prod X_{\mathfrak{p}}.
\end{displaymath}
Let $\lambda$ be a polarization of $X$ which induces the conjugation
on $K/F$. Then the decomposition extends to 
\begin{equation}\label{Unif1e}
  (X, \iota, \lambda) =
  \prod (X_{\mathfrak{p}}, \iota_{\mathfrak{p}}, \lambda_{\mathfrak{p}}). 
\end{equation}
We call $(X, \iota, \lambda)$ a semi-local CM-triple of type $r$ if each 
$(X_{\mathfrak{p}}, \iota_{\mathfrak{p}}, \lambda_{\mathfrak{p}})$ is a local
CM-triple of type $r_{\mathfrak p}$ with respect to $K_{\mathfrak{p}}/F_{\mathfrak{p}}$. This makes sense
because $E(K_{\mathfrak{p}}/F_{\mathfrak{p}}, r_{\mathfrak{p}}) \subset E_{\nu}$. 

\begin{definition}\label{Kneun5d}
  A semi-local CM-triple $(X, \iota, \lambda)$ of type
  $(K \otimes \mathbb{Q}_p / F \otimes \mathbb{Q}_p, r)$ over an
  algebraically closed field with a $\kappa_{E_{\nu}}$-algebra structure
  is said to be \emph{compatible with} $(V,\varsigma)$ if, for each $\mathfrak p | p$,
  \begin{displaymath}
    \inv^{r}(X_{\mathfrak{p}}, \iota_{\mathfrak{p}}, \lambda_{\mathfrak{p}}) =
    \inv(V_{\mathfrak{p}}, \varsigma_{\mathfrak{p}}), 
  \end{displaymath}
  and if $\lambda_{\mathfrak{p}}$ is principal, except in the case where\index{compatible semi-local CM-triple}
  $K_{\mathfrak{p}}/F_{\mathfrak{p}}$ is unramified and
  $\inv(V_{\mathfrak{p}}, \varsigma_{\mathfrak{p}}) = -1$. 
  In the latter case $\lambda_{\mathfrak{p}}$ is almost principal.   
\end{definition}
We note that the  CM-triple
$(\mathbb{X}, \iota_{\mathbb{X}}, \lambda_{\mathbb{X}})$ over $\bar\kappa_{\nu}$
is compatible with $(V,\varsigma)$ and satisfies the conditions $({\rm KC}_r)$
and $({\rm EC}_r)$, in the sense explained before Definition \ref{Kneun5d} 

\begin{definition}\label{KneunU2d} 
  Let $i \in \mathbb{Z}$.  
  Let $\mathcal{M}_{r}(i)$ be the following functor on the category
  $\Nilp_{O_{\breve{E}_{\nu}}}$. For an object $R \in \Nilp_{O_{\breve{E}_{\nu}}}$,
   write $\bar{R} = R \otimes_{O_{\breve{E}_{\nu}}} \bar{\kappa}_{E_{\nu}}$. 
  A point of $\mathcal{M}_{r}(i)(R)$ is given
  by the following data: 

  \begin{enumerate}
  \item[(1)] A CM-triple $(X, \iota, \lambda)$ of type
    $(K \otimes \mathbb{Q}_p / F \otimes \mathbb{Q}_p, r)$ over $\Spec R$ which satisfies
    the conditions $({\rm KC}_r)$ and $({\rm EC}_r)$ and is compatible with
    $(V,\varsigma)$. 
  \item[(2)] A $O_K\otimes\BZ_p$-linear quasi-isogeny
    \begin{displaymath}
  \rho: \bar{X} := X \times_{\Spec R} \Spec \bar{R} \longrightarrow \mathbb{X}
      \times_{\Spec \kappa_{\breve{E}_{\nu}}} \Spec \bar{R}  
    \end{displaymath}
     such that $\rho$ respects the polarization $p^i \lambda$ on $X$ and
    $\lambda_{\mathbb{X}}$ up to a factor in $(O_F \otimes \mathbb{Z}_p)^{\times}$. 
      \end{enumerate}
  We denote these data by $(X, \iota, \lambda, \rho)$. Two data
  $(X, \iota, \lambda, \rho)$ and $(X', \iota', \lambda', \rho')$  define the same point of
  $\mathcal{M}_{r}(i)$ iff there is an isomorphism
  $\alpha: (X, \iota)\longrightarrow (X', \iota')$, such that
  $\rho' \circ \alpha_{\bar{R}} = \rho$. In particular, $\alpha$
  respects the polarizations $\lambda$ and $\lambda'$ up to a factor
  in $(O_F \otimes \mathbb{Z}_p)^{\times}$. 
\end{definition}

\begin{remark}\label{KneunU2dRm}
  In (2) we could replace the last condition on $\rho$ by
    \begin{enumerate}
\item[$(2')$] {\it The quasi-isogeny $\rho$ respects the polarizations as follows,}
    \begin{displaymath}
p^i \lambda = \rho^{\ast} (\lambda_{\mathbb{X}}). 
    \end{displaymath}
  \end{enumerate}
Then we obtain a functor which is naturally isomorphic to $\mathcal{M}_r(i)$.   
This follows because for $a \in (O_F \otimes \mathbb{Z}_p)^{\times}$ the
points $(X, \iota, a\lambda, \rho)$ and $(X, \iota, \lambda, \rho)$
of $\mathcal{M}_r(R)$ are isomorphic.
We could also require $u p^i \lambda = \rho^{\ast} (\lambda_{\mathbb{X}})$ for some
$u \in \mathbb{Z}_p^{\times}$ without changing the functor.
We use different descriptions of the functor $\mathcal{M}_r(i)$ in order to
describe better different group actions. 
  \end{remark}

Let $\tau_{E_{\nu}} \in \Gal(\breve{E}_{\nu}/E_{\nu})$ be the Frobenius
automorphism and let $f_{E_{\nu}}$ be the inertia index of
$E_{\nu}/\mathbb{Q}_p$, i.e., $q_\nu=p^{f_{E_{\nu}}}$. 
As earlier, the Frobenius $F_{\mathbb{X},\tau_{E_{\nu}}}$
defines a Weil descent datum on these functors,
\index[NO]{ZZWG@$\omega_{\mathcal{M}_r}$}
\begin{equation}\label{KneunU14e}
\omega_{\mathcal{M}_r}: \mathcal{M}_r(i)(R) \longrightarrow
    \mathcal{M}_r(i + f_{E_{\nu}})(R_{[\tau_{E_{\nu}}]}) , 
  \end{equation}
  cf.  (\ref{WD5e}), (\ref{WDu5e}).
Since the degrees of the polarizations $\lambda$ and $\lambda_{\mathbb{X}}$
are the same, it follows that
\begin{displaymath}
2 \height \rho = \height (p^i \; | \; X) = 4[F : \mathbb{Q}]i.
\end{displaymath}
More precisely, $\rho = \prod_{\mathfrak{p}} \rho_{\mathfrak{p}}$ where
$\mathfrak{p}$ runs over the prime ideals of $F$ over $p$. For each
$\mathfrak{p}$ we have
\begin{displaymath}
  2 \height \rho_{\mathfrak{p}} = \height (p^i \; | \; X_{\mathfrak{p}}) =
  4[F_{\mathfrak{p}} : \mathbb{Q}_p]i.
  \end{displaymath}
We define\index[NO]{MBG@$\tilde{\mathcal{M}}_r$}
\begin{equation}\label{KneunU8e}
\tilde{\mathcal{M}}_r = \coprod_{i \in \mathbb{Z}} \mathcal{M}_r(i).
  \end{equation}
We describe the functor $\tilde{\mathcal{M}}_r$ with its Weil descent datum. Let \index[NO]{JAG@$J(\mathbb{Q}_p)$}
\begin{equation}\label{Jpadic}
J(\mathbb{Q}_p) = \{\alpha \in \End^o_{K \otimes \mathbb{Q}_p} \mathbb{X}
\; | \; \alpha^{\ast} (\lambda_{\mathbb{X}}) = c \lambda_{\mathbb{X}}, \;
\text{for some} \; c \in \mathbb{Q}_p^{\times}  \} .
  \end{equation}
This group acts naturally on $\tilde{\mathcal{M}}_r$ via the rigidification
$\rho$.
We consider the decomposition (\ref{Unif1e}) for $\mathbb{X}$. We set 
\index[NO]{JAK@$J_{\mathfrak{p}}$}
\begin{equation}\label{Endo9e} 
  J_{\mathfrak{p}} = \{\alpha \in \End^o_{K_{\mathfrak{p}}}
  \mathbb{X}_{\mathfrak{p}} \; |
  \; \alpha^{\ast} (\lambda_{\mathfrak{p}}) = c \lambda_{\mathfrak{p}},
  \; \text{for some} \; c \in \mathbb{Q}_p^{\times} \}
\end{equation}
For all $\mathfrak{p}$ the groups $J_{\mathfrak{p}}$ are subgroups of $J'$
as introduced in section  \ref{s:modformal} in the local cases and they agree with $J$
introduced in the banal cases.

We will give an explicit description of these groups. For this, it is
convenient to replace the bilinear form $\varsigma_{\mathfrak{p}}$ by the
$F_{\mathfrak{p}}$-bilinear form
\begin{displaymath}
  \tilde{\varsigma}_{\mathfrak{p}}: V_{\mathfrak{p}} \times V_{\mathfrak{p}} \rightarrow
  F_{\mathfrak{p}}, 
\end{displaymath}
which is defined by
\begin{displaymath}
  \mathbf{t}(a\tilde{\varsigma}_{\mathfrak{p}}(x_1, x_2)) =
  \varsigma_{\mathfrak{p}}(a x_1, x_2), \quad a \in F_{\mathfrak{p}}, 
  \end{displaymath}
for $\mathbf{t}(a) = \Trace_{F_{\mathfrak{p}}/\mathbb{Q}_p} \vartheta^{-1}a$, where as usual $\vartheta \in O_F$ is the different of
  $F/\mathbb{Q}_p$. 
The restriction to the lattices $\Lambda_{\mathfrak{p}}$ gives
\begin{displaymath}
  \tilde{\varsigma}_{\mathfrak{p}}: \Lambda_{\mathfrak{p}} \times \Lambda_{\mathfrak{p}}
  \rightarrow O_{F_{\mathfrak{p}}}.  
  \end{displaymath}

Let us consider the prime $\mathfrak{p} = \mathfrak{p}_{v}$. We denote by 
$D_{v}$  the quaternion division algebra over $F_{{v}}$.
We choose a two-dimensional
$K_{{v}}$-vector space with an anti-hermitian form
\begin{displaymath}
  \bar{\varsigma}_{\mathfrak{p}_{v}}: \bar{V}_{\mathfrak{p}_{v}} \times
  \bar{V}_{\mathfrak{p}_{v}} \rightarrow F_{\mathfrak{p}_{v}}  
\end{displaymath}
of invariant $+1$. The contraction functor associates to
$(\mathbb{X}_{\mathfrak{p}_{v}}, \iota_{\mathbb{X}_{\mathfrak{p}_{v}}},
\lambda_{\mathbb{X}_{\mathfrak{p}_{v}}})$
a special formal $O_{D_{v}}$-module $\mathbb{Y}$ with the relative
polarization $\lambda_{v} = \psi_{v}$, resp. $\lambda_{v} = \theta_{v}$, as
in section 5.2 resp. 5.3.
Since the endomorphism ring is not changed by the contraction functor,
it follows from  Lemmas \ref{J^v1l} and \ref{J^uv1l} 
that there is an isomorphism 
\begin{equation}\label{Endo12e}
  J_{\mathfrak{p}_{v}} = G(\bar{V}_{\mathfrak{p}_{v}},
  \bar{\varsigma}_{\mathfrak{p}_{v}}). 
\end{equation}
Indeed, to see that the two groups agree we can assume that
$\bar{V}_{\mathfrak{p}_{v}} = K_{\mathfrak{p}_{v}}^2$ and that the antihermitian form
$\varkappa$ associated to  $\bar{\varsigma}_{\mathfrak{p}_{v}}$, cf. section
\ref{ss:binform}, is given by the
matrix
\begin{displaymath}
\left( 
  \begin{array}{rr}
    0 & 1\\
    -1 & 0
  \end{array}
  \right). 
  \end{displaymath}
Then an elementary computation yields that the right hand side of
(\ref{Endo12e}) coincides with the groups defined by the exact sequences
(\ref{Pol29J}). 

For a banal prime $\mathfrak{p} |p$ of $F$, we consider
  the image $(C_{\mathbb{X}_{\mathfrak{p}}}, \phi_{\mathbb{X}_{\mathfrak{p}}})$ by the
polarized contraction functor $\mathfrak{C}_{r, k}^{\rm pol}$ of
Theorem \ref{KneunBa4p}. By Proposition \ref{KneunBa1l}, it follows from
Condition (ii) in Definition \ref{KneunA1d} that there is an isomorphism
\begin{equation}\label{Endo10e} 
  (C_{\mathbb{X}_{\mathfrak{p}}}, \phi_{\mathbb{X}_{\mathfrak{p}}}) \cong
  (\Lambda_{\mathfrak{p}}, \tilde{\varsigma}_{\mathfrak{p}}).
  \end{equation}
More precisely, Condition (ii) implies that the corresponding vector
spaces are isomorphic; the integral isomorphism follows from Lemmas
 \ref{Kneun0l} and \ref{Kneun01l}.
Therefore we obtain
\begin{equation}\label{KneunU17e} 
  J_{\mathfrak{p}} = G(V_{\mathfrak{p}}, \tilde{\varsigma}_{\mathfrak{p}}) =
  G_{\mathfrak{p}},
  \quad \text{for} \; \mathfrak{p} \neq \mathfrak{p}_{v}. 
\end{equation}
 Since we want a uniform notation, we set
$(\bar{V}_{\mathfrak{p}}, \bar{\varsigma}_{\mathfrak{p}}) = (V_{\mathfrak{p}}, \tilde{\varsigma}_{\mathfrak{p}})$ 
for $\mathfrak{p} \neq \mathfrak{p}_{v}$.
We set\index[NO]{GAB@$\bar{G}_{\mathfrak{p}}$}
\begin{displaymath}
\bar{G}_{\mathfrak{p}} =  G(\bar{V}_{\mathfrak{p}}, \bar{\varsigma}_{\mathfrak{p}}). 
  \end{displaymath}
We now have fixed an isomorphism $J_{\mathfrak{p}} \cong \bar{G}_{\mathfrak{p}}$
for all $\mathfrak p|p$. For $\mathfrak{p}$ banal, we have
$ G_{\mathfrak{p}} = \bar{G}_{\mathfrak{p}}$. 

Let
\begin{displaymath}
\bar{V}_p = \oplus_{\mathfrak{p}|p} \bar{V}_{\mathfrak{p}} .
\end{displaymath}
This is an $K \otimes \mathbb{Q}_p$-module. Let 
\begin{displaymath}
  \bar{\varsigma}_p: \bar{V}_{p} \times \bar{V}_{p}
  \rightarrow F \otimes \mathbb{Q}_p 
  \end{displaymath}
be the orthogonal sum of the forms $\bar{\varsigma}_{\mathfrak{p}}$. We define
\index[NO]{GAC@$\bar{G}(\mathbb{Q}_p)$}
\begin{equation}\label{defbarG}
  \bar{G}(\mathbb{Q}_p) := G(\bar{V}_p, \bar{\varsigma}_p) := 
  \{g \in \Aut_{K \otimes \mathbb{Q}_p} (\bar{V}_p) \; | \;
  \bar{\varsigma}_p (gx, gy) = c \bar{\varsigma}_p (x, y),
  \; \text{for some} \; c \in \mathbb{Q}_p^{\times}\}. 
  \end{equation}
We have shown that $\bar{G}(\mathbb{Q}_p) = J(\mathbb{Q}_p)$. 
In the description of the descent data, the following slightly larger group will be needed. We  define the group $\bar{G}'_p \supset \bar{G}(\mathbb{Q}_p)$
\index[NO]{GZD@$\bar{G}'_p$} via
\begin{displaymath}
  \bar{G}'_p = \{g \in \Aut_{K \otimes \mathbb{Q}_p} (\bar{V}_p) \; |
  \; \bar{\varsigma}_p(gx, gy) = \mu_p(g) \bar{\varsigma}_p(x,y), \; \text{for}
  \; \mu_p(g) \in p^{\mathbb{Z}} (O_F \otimes \mathbb{Z}_p)^{\times}  \}, 
  \end{displaymath}
and
\begin{displaymath}
  \bar{G}'_{\mathfrak{p}} = \{g \in \Aut_{K_{\mathfrak{p}}} (\bar{V}_{\mathfrak{p}}) \; |
  \; \bar{\varsigma}_{\mathfrak{p}}(gx, gy) =
  \mu_{\mathfrak{p}}(g) \bar{\varsigma}_{\mathfrak{p}}(x,y), \; \text{for}
  \; \mu_{\mathfrak{p}}(g) \in p^{\mathbb{Z}} O_{F_{\mathfrak{p}}}^{\times} \} . 
\end{displaymath}
The groups $\bar{G}'_{\mathfrak{p}}$ are isomorphic to the groups
$J'_{\mathfrak{p}} =J'$ introduced in section 6 in the local cases. We fix these
isomorphisms which are associated to the framing objects. Therefore the groups
$\bar{G}'_{\mathfrak{p}}$ act on the local moduli spaces $\mathcal{M}$ of section
6 and the subgroup $\bar{G}'_p \subset \prod_{\mathfrak{p}|p} \bar{G}'_{\mathfrak{p}}$
acts on $\tilde{\mathcal{M}}_r$, cf. (\ref{KneunU8e}).

We define the group $\hat{G}'(\mathbb{Q}_p)$
\index[NO]{GAE@$\hat{G}'(\mathbb{Q}_p)$} as the union of the following
sets for $i \in \mathbb{Z}$,
  \begin{equation}\label{defhatG'}
  \hat{G}'(i) = \{ (c, g_{\mathfrak{p}}) \in
  p^{i} O_{F_{\mathfrak{p}_{v}}}^{\times} \times \prod_{\mathfrak{p}\; \text{banal}}
  \bar{G}'_{\mathfrak{p}} \; | \;
  \mu_{\mathfrak{p}}(g_{\mathfrak{p}}) \in p^{i} O_{F_{\mathfrak{p}}}^{\times}, \;
  \text{for all} \; \mathfrak{p}\} .
\end{equation}
  Let $\hat{G}'(\mathbb{Z}_p) \subset \hat{G}'(\mathbb{Q}_p)$
  \index[NO]{GAF@$\hat{G}'(\mathbb{Z}_p)$} be the subgroup
  of elements
 $(c, g_{\mathfrak{p}})$ such that $c \in O_{F_{\mathfrak{p}_{v}}}^{\times}$ and
 $g_{\mathfrak{p}}(\Lambda_{\mathfrak{p}}) = \Lambda_{\mathfrak{p}}$.
 The multiplicator
 $\mu_{\mathfrak{p}_{v}}: \bar{G}'_{\mathfrak{p}_{v}} \rightarrow p^{\mathbb{Z}} O_{F_{\mathfrak{p}_{v}}}^{\times}$
 induces  homomorphisms
 \begin{equation}\label{defmaps}
 \bar{G}'_p \rightarrow \hat{G}'(\mathbb{Q}_p)  \quad \text{ and } \quad G(\mathbb{Q}_p) \rightarrow \hat{G}'(\mathbb{Q}_p) .
\end{equation}  For the second map we used the identification 
 $\bar{G}_{\mathfrak{p}} = G_{\mathfrak{p}}$ for $\mathfrak{p}$ banal. 

 \begin{definition}\label{Unif3d}
  We consider the following element \index[NO]{WAA@$w'_r$}
  $w'_r = (c,w_{\mathfrak{p}}) \in \hat{G}'(\mathbb{Q}_p)$.
  \begin{enumerate}
  \item[(1)] $c = p^{f_{E_{\nu}}}$. 
  \item[(2)] If $K_{\mathfrak{p}}/F_{\mathfrak{p}}$ is ramified and hence 
    $\lambda_{\mathbb{X}_\mathfrak{p}}$ is principal,   $w_{\mathfrak{p}}$
    is the multiplication
    \begin{displaymath}
      \Pi_{\mathfrak{p}}^{e_{\mathfrak{p}}f_{E_{\nu}}}: \bar V_{\mathfrak{p}}
      \longrightarrow  \bar V_{\mathfrak{p}} ,
    \end{displaymath}
    see Proposition \ref{KneunBa5p}. 
  \item[(3)] If $K_{\mathfrak{p}}/F_{\mathfrak{p}}$ is unramified,  then both principal and almost principal
    $\lambda_{\mathbb{X}_{\mathfrak{p}}}$ are allowed.
    In both cases we define $w_{\mathfrak{p}}$ as the multiplication
    \begin{displaymath}
      \pi_{\mathfrak{p}}^{e_{\mathfrak{p}}f_{E_{\nu}}/2}: \bar V_{\mathfrak{p}}
      \longrightarrow \bar V_{\mathfrak{p}} ,
    \end{displaymath}
    see Proposition \ref{ubF2p}. 
  \item[(4)]
    In the case where 
    $K_{\mathfrak{p}} = F_{\mathfrak{p}} \times F_{\mathfrak{p}}$
    is split and hence  $\lambda_{\mathbb{X}_{\mathfrak{p}}}$ is principal, we have the
    decompositions
    $\mathbb{X}_{\mathfrak{p}} = \mathbb{X}_{\mathfrak{p},1} \times \mathbb{X}_{\mathfrak{p},2}$ and 
   $\bar{V}_{\mathfrak{p}} = \bar{V}_{\mathfrak{p}, 1}\oplus\bar{V}_{\mathfrak{p}, 2}$.
    We set
    $a_{\mathfrak{p},i,E_{\nu}} = a_{\mathfrak{p},i} \frac{f_{E_{\nu}}}{f_{\mathfrak{p}}}$
    for $i = 1,2$,  where $2a_{\mathfrak{p},i} = \dim \mathbb{X}_{\mathfrak{p}_i}$,
    cf. (\ref{sumsplit})  and the dicussion before
    Proposition \ref{zerlegt1p}.
    We define $w_{\mathfrak{p}}$ to be the multiplication by
    $\pi_{\mathfrak{p}}^{a_{\mathfrak{p},1,E_{\nu}}}$ on $\bar V_{\mathfrak{p}, 1}$ and the
    multiplication by $\pi_{\mathfrak{p}}^{a_{\mathfrak{p},2,E_{\nu}}}$ on
    $\bar V_{\mathfrak{p}, 2}$. 
  \end{enumerate}  
\end{definition}

\begin{proposition}\label{compfo}
    There exists an isomorphism
         \begin{displaymath}
    \tilde{\mathcal{M}}_{r} \isoarrow  
    (\wh{\Omega}_{F_{v}} \times_{\Spf O_{F_{v}}} \Spf O_{\breve{E}_{\nu}}) \times \hat{G}'(\mathbb{Q}_p)/\hat{G}'(\mathbb{Z}_p) 
        \end{displaymath}
which is equivariant with respect to the action of $\bar{G}'_p$
on both sides. This extends the action of
$J(\mathbb{Q}_p) = \bar{G}(\mathbb{Q}_p) \subset \bar{G}'_p$.  
  
The Weil  descent datum $\omega_{\mathcal{M}_r}$ relative to
$O_{\breve{E}_{\nu}}/O_{E_{\nu}}$ on the left hand side \eqref{KneunU14e} corresponds on the
right hand side to
      \begin{displaymath}
        (\xi, g) \mapsto (\omega_{\tau_{E_{\nu}}}(\xi), w'_rg), \quad
        g \in \hat{G}'(\mathbb{Q}_p). 
      \end{displaymath} 
\end{proposition}
\begin{proof} We use the decomposition
\begin{displaymath}
\mathcal{M}_r(i) = \prod_{\mathfrak{p}|p} \mathcal{M}_{r_{\mathfrak{p}}}(i) , 
  \end{displaymath}
which follows immediately from (\ref{Unif1e}). 
Then we conclude by the results of section \ref{s:modformal}, in particular  
Propositions \ref{KneunBa5p}, \ref{ubF2p}, \ref{zerlegt1p}. 
    \end{proof}
  
  \begin{remark}\label{compfoRm}  
    We may multiply each $w_{\mathfrak{p}}$ by a unit in $K_{\mathfrak{p}}$ in the
    Definition \ref{Unif3d} of $w'_r$. This does not change the assertion of
    the last Proposition. 
    \end{remark}

  We introduce the group
\begin{equation}\label{KneunU18e} 
 \hat{G}(\mathbb{Q}_p) = \{ (c, g_{\mathfrak{p}}) \in
 \mathbb{Q}_{p}^{\times} \times \prod_{\mathfrak{p}\; \text{banal}} G_{\mathfrak{p}}
 \; | \;
 \mu_{\mathfrak{p}}(g_{\mathfrak{p}}) = c, \; \text{for all} \; \mathfrak{p}\;
 \text{banal}\}. 
\end{equation}
There are natural homomorphisms 
\begin{equation}
G(\mathbb{Q}_p) \rightarrow \hat{G}(\mathbb{Q}_p) \quad\text{ and }\quad \bar{G}(\mathbb{Q}_p) \rightarrow \hat{G}(\mathbb{Q}_p) .
\end{equation}  
For the second map, we used that  in the definition of $\hat{G}(\mathbb{Q}_p)$ we can replace $G_{\mathfrak{p}}$ by
$\bar{G}_{\mathfrak{p}}$. In particular
the groups $G(\mathbb{Q}_p)$ and 
$J(\mathbb{Q}_p)$   act on $\hat{G}(\mathbb{Q}_p)$. 
We denote by
$\hat{G}(\mathbb{Z}_p) \subset \hat{G}(\mathbb{Q}_p)$ the subgroup of all 
$(c, g_{\mathfrak{p}})$ such that $c \in \mathbb{Z}_p^{\times}$ and 
$g_{\mathfrak{p}} \Lambda_{\mathfrak{p}} = \Lambda_{\mathfrak{p}}$. By  Corollaries \ref{KneunBa4c},
\ref{J'unr}, \ref{J'split}, we obtain a bijection 
\begin{equation}\label{G^G'1e}  
  \hat{G}(\mathbb{Q}_p)/\hat{G}(\mathbb{Z}_p) \isoarrow
  \hat{G}'(\mathbb{Q}_p)/\hat{G}'(\mathbb{Z}_p).
\end{equation}
\begin{corollary}\label{KneunU1c} 
 There exists an isomorphism
         \begin{displaymath}
    \tilde{\mathcal{M}}_{r} \isoarrow  
    (\wh{\Omega}_{F_{v}} \times_{\Spf O_{F_{v}}} \Spf O_{\breve{E}_{\nu}}) \times \hat{G}(\mathbb{Q}_p)/\hat{G}(\mathbb{Z}_p) 
        \end{displaymath}
which is equivariant with respect to the action of $J(\BQ_p)$
on both sides. \qed
\end{corollary}
Note that in this version of Proposition \ref{compfo} we loose control of the
descent data.
\subsection{The $p$-adic uniformization}
We will now define a uniformization morphism in the sense of \cite{RZ}. 
We fix a point
$(A_0, \iota_0, \bar{\lambda}_0, \bar{\eta}^{p}_0)$ of
$ \mathcal{A}_{\bK}(\bar{\kappa}_{{\nu}})$.  
The uniformization morphism will depend on the choice of
$\eta_0^{p} \in \bar{\eta}^{p}_0$. This choice defines a point of the
proscheme $\projlim\nolimits_{\bK^p} \mathcal{A}_{\bK}$ for all congruence
subgroups ${\bK} = {\bK}_p{\bK}^{p}$ as above. We also fix a polarization
$\lambda_0 \in \bar{\lambda}_0$ which satisfies the
condition (i) of Definition \ref{KneunA1d}. Let $(\BX, \iota_\BX, \lambda_\BX)$
be the $p$-divisible group  corresponding to $(A_0, \lambda_0)$. 

We denote by
$\hat{\mathcal{A}}_{\bK}$ the restriction of $\mathcal{A}_{\bK}$ to the category
$\Nilp_{O_{E_{\nu}}}$. The uniformization morphism \index{uniformization morphism}
\index[NO]{ZZHA@$\Theta$}
\begin{equation}\label{KneunU5e}
  \Theta: 
  \tilde{\mathcal{M}}_{r} \times G(\mathbb{A}_f^{p})/{\bK}^p \longrightarrow
  \hat{\mathcal{A}}_{\bK} \times_{\Spf O_{E_{\nu}}} \Spf O_{\breve{E}_{\nu}} 
\end{equation}
is defined as follows. Let
$(X, \iota, \lambda, \rho) \in \mathcal{M}_{r}(i)(R)$ and let
$g \in G(\mathbb{A}_f^{p})$. Recall the notation
$\bar{R}=R\otimes_{O_{\breve{E}_{\nu}}}\bar{\kappa}_\nu$.
  There exists an abelian scheme $\bar{A}$ over $\Spec \bar{R}$ endowed with an
  isomorphism of the $p$-divisible group of $\bar{A}$ with $\bar{X}$ and 
  with a quasi-isogeny of abelian schemes of order a power of $p$
\begin{equation}\label{KneunU4e}
  \rho: \bar{A} \longrightarrow
  A_0 \times_{\Spec \bar{\kappa}_{{\nu}}} \Spec \bar{R}.
\end{equation}
 which induces on the $p$-divisible groups the given map  
 $\rho: \bar{X}\rightarrow \mathbb{X} \times_{\Spec \bar{\kappa}_{{\nu}}}\Spec\bar{R}$.
 The pair $(\bar{A}, \rho)$ is unique up to canonical isomorphism. 
Because $O_K$ acts on $\bar{X}=X\otimes_R \bar R$, we obtain a map
$O_K\to \End(\bar A)\otimes\BZ_{(p)}$. Moreover the polarization
$\lambda_0: A_0 \longrightarrow A_0^{\wedge}$
induces on $\bar{A}$ a quasi-polarization
$\lambda'_{\bar{A}}: \bar{A} \longrightarrow \bar{A}^{\wedge}$ and $\bar{\eta}_0^{p}$
induces
\begin{displaymath}
  \bar{\eta}^p_{\bar{A}} = V^p(\rho^{-1}) \circ \eta^p_0:
  V \otimes \mathbb{A}_f^p
  \isoarrow V^{p}(\bar{A})   \; \mod {\bK}^p. 
\end{displaymath}
On the $p$-divisible groups, $\lambda'_{\bar{A}}$ differs from $p^i\lambda$ by a factor from $(O_F \otimes \mathbb{Z}_p)^{\times}$ and therefore
$\lambda_{\bar{A}} := p^{-i} \lambda'_{\bar{A}}$ satisfies the condition  (i)
in the Definition \ref{KneunA1d} of the functor $\mathcal{A}_{\bK}$. 

We associate to the pair $(X, g)$ from the left hand side of (\ref{KneunU5e})
the point
\begin{equation}\label{KneunU6e}
  (\bar{A}, \iota_{\bar{A}}, \lambda_{\bar{A}}, \bar{\eta}^p_{\bar{A}} g)
  \in \mathcal{A}_{\bK}(\bar{R}). 
\end{equation}
The CM-triple $(X, \iota, \lambda)$ over $R$ defines by the Serre-Tate theorem
a lifting of (\ref{KneunU6e}) to a point of $\mathcal{A}_{\bK}(R)$. This
finishes the definition of the uniformization morphism $\Theta$ in
(\ref{KneunU5e}). 

\begin{lemma}
  The uniformization morphism is compatible with the
  Weil descent data $\omega_{\mathcal{M}_r}$  acting on the
  first factor on the left hand side and the natural Weil descent data on
  $\hat{\mathcal{A}}_{\bK} \times_{\Spf O_{E_{\nu}}} \Spf O_{\breve{E}_{\nu}}$. 
\end{lemma}
\begin{proof}
  This is essentially \cite[Thm. 6.21]{RZ}  but we repeat the simple argument
  in our context. By definition of the Weil descent data repeated below, it is
  enough to consider both sides of (\ref{KneunU6e}) on the category of
  $\bar{\kappa}_{{\nu}}$-algebras $R$. We will denote by
  $\varepsilon: \bar{\kappa}_{{\nu}} \longrightarrow R$ the algebra structure. 
  Consider a point $(X, \iota, \lambda, \rho) \in \tilde{\mathcal{M}}_r(R)$.
  The Weil  descent datum $\omega_{\mathcal{M}_r}$ is obtained by changing
  $\rho$ to $\rho'$: 
\begin{displaymath}
  \rho': X \overset{\rho}{\longrightarrow}
  \varepsilon_{\ast} \mathbb{X}
  \overset{\varepsilon_{\ast}F_{\mathbb{X},\tau_{E_{\nu}}}}{\longrightarrow}
  \varepsilon_{\ast}(\tau_{E_{\nu}})_{\ast} \mathbb{X}. 
\end{displaymath}
This gives a point
$(X, \iota, \lambda, \rho') \in \tilde{\mathcal{M}}_r(R_{[\tau_{E_{\nu}}]})$.
The point $(X, \iota, \lambda, \rho)$ defines a quasi-isogeny of abelian varieties
\begin{displaymath}
\rho: A \longrightarrow \varepsilon_{\ast }A_0 ,
\end{displaymath}
as explained in the definition of $\Theta$. The point
$(X, \iota, \lambda, \rho')$ defines in the same way the quasi-isogeny of
abelian varieties over $R_{[\tau_{E_{\nu}}]}$,
\begin{displaymath}
  A \longrightarrow \varepsilon_{\ast} A_0
  \overset{\varepsilon_{\ast} F_{A_0, \tau_{E_{\nu}}}}{\longrightarrow}
  (\varepsilon \tau_{E_{\nu}})_{\ast} A_0.
  \end{displaymath}
Here $A$ with its additional structure is regarded as a point of
$\hat{\mathcal{A}}_{{\bK}^p}(R_{[\tau_{E_{\nu}}]})$. This makes sense because to be a
point of $\hat{\mathcal{A}}_{{\bK}^p}(R)$ depends only on the
$\kappa_{{\nu}}$-algebra structure on $R$. In other words 
\begin{equation}
\hat{\mathcal{A}}_{{\bK}^p}(R) = \hat{\mathcal{A}}_{{\bK}^p}(R_{[\tau_{E_{\nu}}]}).
\end{equation}
But this equation is the Weil descent datum on the right hand side of
(\ref{KneunU5e}).
  \end{proof}

We define the  group\index[NO]{JAL@$J(\mathbb{Q})$}
\begin{equation}\label{KneunU2e-alt}
  J(\mathbb{Q}) = \{\gamma \in \End^o_K A_0 \; | \; \gamma^{\ast} (\lambda_0) =
  u \lambda_0, \; \text{ for some } u \in \mathbb{Q}^\times \}, 
  \end{equation}
cf. \eqref{KneunU2e}. Regarded as an algebraic group over $\mathbb{Q}$, the
group $J$ is an inner form of $G$. In the proof of Proposition \ref{KneunUni3p}
we saw that the $\BQ_p$-valued points of $J$ coincide with the group $J(\BQ_p)$
of \eqref{Jpadic}. We proved in section \ref{RZtildeM} that
$\bar{G}(\mathbb{Q}_p) = J(\mathbb{Q}_p)$. Let $\gamma \in J(\mathbb{Q})$. 
With the chosen $\eta_0^p$, we define $\omega(\gamma) \in G(\mathbb{A}_f^{p})$
by the equation
\begin{equation}\label{KneunU10e}
   V^p(\gamma) \circ\eta_0^p  =  \eta_0^p \omega(\gamma).
\end{equation}
This defines a homomorphism
\begin{displaymath}
  \omega: J(\mathbb{Q}) \longrightarrow G(\mathbb{A}_f^{p}), 
\end{displaymath}
and an isomorphism $J(\mathbb{A}_f^{p}) \cong G(\mathbb{A}_f^{p})$.
Therefore $J$ and $G$ are isomorphic over the finite places $w \neq p$ of
$\mathbb{Q}$. At the infinite place $J$ is anisotropic because the Rosati
involution is positive. 

The group $J(\mathbb{Q}_p)$ acts on $\tilde{\mathcal{M}}_{r}$, 
\begin{displaymath}
  (X, \iota, \lambda, \rho) \mapsto (X, \iota, \lambda, \gamma \rho), \quad
   \; \gamma \in J(\mathbb{Q}_p). 
\end{displaymath}
Let $((X, \iota, \lambda, \rho), g)$, with $g \in G(\mathbb{A}_f^{p})$ be a
point from the left hand side of (\ref{KneunU5e}) and let
$(A, \iota_A, \lambda_A, \eta^p_A g)$ be its image by $\Theta$, cf.
(\ref{KneunU6e}).
If $\gamma \in J(\mathbb{Q})$, the quasi-isogeny $\gamma \rho$ extends
to the quasi-isogeny of abelian schemes
\begin{displaymath}
  \bar{A} \overset{\rho}{\longrightarrow} (A_0)_{\bar{R}} 
  \overset{\gamma}{\longrightarrow} (A_0)_{\bar{R}}. 
\end{displaymath}
In follows from (\ref{KneunU10e}) that the image of
$((X, \iota, \lambda, \gamma\rho), g)$ by the morphism $\Theta$ is
\begin{displaymath}
(A, \iota_A, \lambda_A, \eta^p_A \omega(\gamma^{-1}) g)
  \end{displaymath}

We define an action of $J(\mathbb{Q})$ on the left hand side of
(\ref{KneunU5e}) by
\begin{displaymath}
  ((X, \iota, \lambda, \rho), g) \mapsto
  ((X, \iota, \lambda, \gamma \rho), \omega(\gamma) g). 
  \end{displaymath} 
\begin{proposition}\label{unifth}
The
uniformization morphism (\ref{KneunU5e}) factors through an isomorphism  
\begin{equation*}\label{KneunU7e}
  \Theta: 
  J(\mathbb{Q})\backslash
  (\tilde{\mathcal{M}}_{r} \times G(\mathbb{A}_f^{p})/{\bK}^p)
  \isoarrow \hat{\mathcal{A}}_{\bK} \times_{\Spf O_{E_{\nu}}} \Spf O_{\breve{E}_{\nu}}. 
\end{equation*}
This isomorphism is compatible with the Weil descent data  relative to
$O_{\breve{E}_{\nu}}/O_{{E_\nu}}$. Here  the Weil descent datum on the left is
induced from $\omega_{\mathcal{M}_r}$, cf. Proposition \ref{KneunU14e}. 
  \end{proposition}
\begin{proof}
  We have just proved that the morphism is well-defined. The bijectivity follows
  from the Proposition \ref{KneunUni3p} and \cite[Thm. 6.30]{RZ}. 
  \end{proof}
By inserting Proposition \ref{compfo} in this result, we  obtain our main theorem about uniformization. 
  \begin{theorem}\label{maint} 
  Assume that $\bK^p$ is sufficiently small. In particular, $\CA_\bK$ is representable, cf. Proposition \ref{proprepr}.
 \smallskip
 
\noindent (i)  The $O_{E_{\nu}}$-scheme $\CA_\bK$ is a  projective and flat relative curve, which is
stable in the sense of Deligne-Mumford \cite{DM}. 

\smallskip

\noindent (ii) Let $\hat\CA_\bK$
\index[NO]{ABB@$\hat{\mathcal{A}}_{\mathbf{K}}$}
be the completion of $\CA_\bK$ along its special fiber, which is a formal scheme over $\Spf O_{E_\nu}$. There exists an isomorphism of formal schemes over $\Spf O_{\breve E_\nu}$,
\begin{equation}\label{maint1e} 
  J(\BQ)\bs\big[\big(\widehat{\Omega}_{F_v} \times_{{\rm Spf}\,O_{F_v}}{\rm Spf}\,
    O_{\breve E_\nu}\big) \times \hat{G}'(\mathbb{Q}_p)/\hat{G}'(\mathbb{Z}_p)
    \times G(\BA_f^p)/{\bK}^p\big]
  \isoarrow  \hat\CA_\bK\times_{{\rm Spf}\,O_{E_{\nu}}}{\rm Spf}\, O_{\breve E_\nu}  \, .
  \end{equation}
For varying $\bK^p$, this isomorphism is compatible with the action of
$G(\BA_f^p)$ through Hecke correspondences on both sides.

Let $w'_r$ the element in the center of $\hat{G}'(\mathbb{Q}_p)$ of Definition \ref{Unif3d}.
We endow
the left hand with the Weil descent datum
\begin{displaymath}
  (\xi, h, g) \mapsto (\omega_{\tau_{E_{\nu}}}(\xi), w'_rh, g), \quad 
  h \in \hat{G}'(\mathbb{Q}_p), \; g \in G(\BA_f^p).
\end{displaymath}
Then the isomorphism (\ref{maint1e}) is compatible with the Weil descent
data on both sides. \qed
\end{theorem} 
 \begin{proof}
 The only assertion that remains to be proved is that  the geometric special fiber of $\CA_\bK$ is a stable  curve. This geometric special fiber  is a finite disjoint sum of schemes of the form 
 \begin{equation}\label{specOm}
 \bar\Gamma\backslash(\widehat{\Omega}_{F_v} \times_{{\rm Spf}\,O_{F_v}}{\rm Spec}\, \bar\kappa_\nu) ,
 \end{equation}
where $\bar\Gamma\subset \PGL_2(F_v)$ is a discrete group, comp. \cite[Rmks. 5.4]{BC} or  \cite[Cor. 6.8]{BZ1}. Consider the action of $\bar\Gamma$ on the Bruhat-Tits tree $\sB$ of $ \PGL_2(F_v)$. By making $\bK^p$ sufficiently small,  we can make sure that no non-trivial element of $\bar\Gamma$ takes a vertex
of $\sB$  into itself or to an adjacent vertex. Therefore no non-trivial element of $\bar\Gamma$
takes an irreducible component of $\widehat{\Omega}_{F_v} \times_{{\rm Spf}\,O_{F_v}}{\rm Spec}\, \bar\kappa_\nu$ into itself or into a second component meeting the first one. From the structure of $\widehat{\Omega}_{F_v} \times_{{\rm Spf}\,O_{F_v}}{\rm Spec}\, \bar\kappa_\nu$, it follows that \eqref{specOm} is reduced and that each irreducible component of \eqref{specOm} is a projective line meeting $q+1$ other irreducible components. Here $q$ is the number of elements in the residue field of $O_{F_v}$. Hence each irreducible component meets at least three other components, which proves the stability of  \eqref{specOm}. 
 \end{proof} 
 Under additional assumptions there is a much simpler statement of 
 Theorem \ref{maint} which uses $\hat{G}$ instead of $\hat{G}'$. We can replace $\hat{G}'$
 by $\hat{G}$ in (\ref{maint1e}) using (\ref{G^G'1e}), where we recall $\hat{G}(\mathbb{Q}_p)$ from \eqref{KneunU18e}. We define
 $\hat{G}(\mathbb{A}_f) = \hat{G}(\mathbb{Q}_p) \times G(\mathbb{A}_f^p)$ and
  $\hat{\mathbf{K}} = \hat{G}(\mathbb{Z}_p) \times \mathbf{K}^p$, cf. 
 \eqref{KneunU18e}. 
 
   \begin{corollary}\label{cortomain}
    There is a natural isomorphism of formal schemes
    \begin{equation}\label{maint2e} 
J(\BQ)\bs\big[\big(\widehat{\Omega}_{F_v} \times_{{\rm Spf}\,O_{F_v}}{\rm Spf}\,
    O_{\breve E_\nu}\big) 
    \times \hat{G}(\BA_f)/\hat{\mathbf{K}}\big]
  \isoarrow  \hat\CA_\bK\times_{{\rm Spf}\,O_{E_{\nu}}}{\rm Spf}\, O_{\breve E_\nu}  \, .
    \end{equation}
    Assume that the inertia index $f_{E_{\nu}}$ is even. Assume moreover
    for  prime ideals $\mathfrak{p}|p$ of $F$ which split in $K$ that  
$a_{\mathfrak{p},1} = a_{\mathfrak{p},2} = [F_{\mathfrak{p}}:\mathbb{Q}_p]/2$
in the notation of Definition \ref{Unif3d}. 
The multiplication by
    $p$ on  $V \otimes \mathbb{Q}_p$ defines an element of $G(\mathbb{Q}_p)$.
    Let $\hat{p}$ be the image in $\hat{G}(\mathbb{Q}_p)$. We also denote
    by $\hat{p}$ the element
    \begin{displaymath}
      (\hat{p}, 1) \in \hat{G}(\mathbb{Q}_p) \times G(\mathbb{A}_f^p) =
      \hat{G}(\mathbb{A}_f).
    \end{displaymath}
    If we endow the left hand side of (\ref{maint2e}) with the Weil descent
    datum
    \begin{displaymath}
  (\xi, g) \mapsto (\omega_{\tau_{E_{\nu}}}(\xi), \hat{p}^{f_{E_{\nu}}/2} g), \quad 
  \; g \in \hat{G}(\BA_f) ,
    \end{displaymath}
   then  the morphism (\ref{maint2e}) is compatible with the Weil descent data. 
  \end{corollary}
  \begin{proof}
    We use the notations of the Theorem. By Remark (\ref{compfoRm}), we may change the components $w_{\mathfrak{p}}$ for banal $\mathfrak{p}$
    in $w'_r = (p^{f_{E_{\nu}}}, w_{\mathfrak{p}})$ by units in $K_{\mathfrak{p}}$. It follows easily from our assumptions
    that $w_{\mathfrak{p}}$ and $p^{f_{E_{\nu}}/2}$ differ by an element in
    $O_{F_{\mathfrak{p}}}^{\times}$. The Corollary follows.       
  \end{proof}
  
  Note that $  {G}(\mathbb{Q}_p)/{G}(\mathbb{Z}_p)\isoarrow \hat{G}(\mathbb{Q}_p)/\hat{G}(\mathbb{Z}_p)$, as follows from
  $\bK_{{\mathfrak p}_v}=\ker(\mu_{{{\mathfrak p}_v}}\colon G_{{\mathfrak p}_v}(\BQ_p)\to \BQ_p^\times/\BZ_p^\times)$.
    Let us indicate briefly the last identity. We use the exact sequence
  \[
  1\to \SU(V_{{{\mathfrak p}_v}})\to G_{{{\mathfrak p}_v}}(\BQ_p)\to \BQ_p^\times\times F_{{{\mathfrak p}_v}}^1 \to 1.
  \]
 (This exact sequence is induced by an exact sequence of algebraic groups over $\BQ_p$.) The right map is given by $g\mapsto (\mu_{{{\mathfrak p}_v}}(g), \frac{\det(g)}{\mu_{{{\mathfrak p}_v}}(g)})$. Since $V_{{{\mathfrak p}_v}}$ is anisotropic, $\SU(V_{{{\mathfrak p}_v}})$ is compact. Furthermore, $\BZ_p^\times\times F_{{{\mathfrak p}_v}}^1$ is the  unique maximal compact subgroup of the target group. Hence $\bK'_{{{\mathfrak p}_v}}:=\ker(\mu_{{{\mathfrak p}_v}}\colon G_{{\mathfrak p}_v}(\BQ_p)\to \BQ_p^\times/\BZ_p^\times)$ is the unique maximal compact subgroup of $G_{{{\mathfrak p}_v}}(\BQ_p)$.  On the other hand, $\bK'_{{{\mathfrak p}_v}}$ stabilizes $\Lambda_{{{\mathfrak p}_v}}$. Hence $\bK'_{{{\mathfrak p}_v}}=\bK_{{{\mathfrak p}_v}}$ by the maximality of $\bK'_{{{\mathfrak p}_v}}$.

  Hence 
  \begin{equation}
 {G}(\BA_f)/{\mathbf{K}}\simeq  \hat{G}(\BA_f)/\hat{\mathbf{K}} .
\end{equation}
Using these facts, Theorem \ref{maint} and Corollary \ref{cortomain} imply Theorem \ref{ThmB} in the
Introduction.

\subsection{The uniformization for deeper level structures at $p$}\label{ss:deeper}
We now pass to deeper level structures. For each prime ideal $\mathfrak{p}$
of $O_F$ with $\mathfrak{p}|p$ we have the group
\begin{displaymath}
  G_{\mathfrak{p}} = \{g \in \GL_{O_{K_{\mathfrak{p}}}} (V_{\mathfrak{p}}) \; | \;
  \varsigma_{\mathfrak{p}}(g x_1, g x_2) = \mu(g)
  \varsigma_{\mathfrak{p}}(x_1, x_2), \; \text{for some}\; \mu_{\mathfrak{p}}(g) \in
  \mathbb{Q}_p^{\times}\}, 
\end{displaymath}
and the open compact subgroup $\mathbf{K}_{\mathfrak{p}} \subset G_{\mathfrak{p}}$
cf. (\ref{boldK1e}).  
We will assume that there exist prime ideals $\mathfrak{p}$ which are banal
since our deeper level structures exist only in this case. 
For each  banal $\mathfrak{p}$, we choose an open subgroup of  
$\mathbf{K}_{\mathfrak{p}}^{\ast} \subset \mathbf{K}_{\mathfrak{p}}$.
For a natural number $M$ we consider the subgroup
$\mathbf{K}_{\mathfrak{p}}(p^M) \subset \mathbf{K}_{\mathfrak{p}}$ which consists of 
the elements that act trivially on $\Lambda_{\mathfrak{p}}/p^M\Lambda_{\mathfrak{p}}$.
We will assume that for some $M$
\begin{equation}\label{banK2e}
\mathbf{K}_{\mathfrak{p}}(p^M) \subset \mathbf{K}_{\mathfrak{p}}^{\ast} .
  \end{equation}
For the special prime $\mathfrak{p}_{v}$ we set
$\mathbf{K}_{\mathfrak{p}_v}^{\star} = \mathbf{K}_{\mathfrak{p}_{v}}$.
We set
\begin{displaymath}
  \mathbf{K}_{p}^{\star} = \{g = (g_{\mathfrak{p}}) \in G(\mathbb{Q}_p) \; | \;
  g_{\mathfrak{p}} \in \mathbf{K}_{\mathfrak{p}}^{\star} \}. 
\end{displaymath}
This says that $\mu_{\mathfrak{p}}(g_{\mathfrak{p}})$ is independent of
$\mathfrak{p}$. We also introduce\index[NO]{KAG@$\mathbf{K}_{p}^{\star,{\rm ba}}$}
\begin{equation}\label{banK}
  \mathbf{K}_{p}^{\star, {\rm ba}} = \{(g_{\mathfrak{p}}) \in \prod_{\mathfrak{p}, \;
    \text{banal}} \mathbf{K}^{\star}_{\mathfrak{p}}\; | \;
  \mu_{\mathfrak{p}}(g_{\mathfrak{p}}) = c \in \mathbb{Z}_p^{\times}, \;
  \text{independent of} \; \mathfrak{p}\} .
\end{equation}
This is a subgroup of\index[NO]{GAG@$G^{{\rm ba}}(\mathbb{Q}_p)$}
\begin{displaymath}
  G^{{\rm ba}}(\mathbb{Q}_p) = \{(g_{\mathfrak{p}}) \in 
  \prod_{\mathfrak{p}, \; \text{banal}} G_{\mathfrak{p}}\; | \;
  \mu_{\mathfrak{p}}(g_{\mathfrak{p}}) = c \in \mathbb{Q}_p^{\times}, \;
  \text{independent of} \; \mathfrak{p}\}. 
  \end{displaymath}
Also, let $O_K^{{\rm ba}} = \prod_{\mathfrak{p}, \; \text{banal}} O_{K_{\mathfrak{p}}}$.  Since the multiplier $\mu_{\mathfrak{p}_{v}}\colon G_{\mathfrak{p}_{v}}\to\BZ_p^\times$ is surjective, the groups $\mathbf{K}_{p}^{\star, {\rm ba}}$ and $\mathbf{K}_{p}^{\star}$ determine each other. 

We need  some generalities on $p$-divisible groups suited for our special
case. Let $X$ and $Y$ be $p$-divisible groups on a scheme $S$.
We consider the category of \'etale morphisms $U \rightarrow S$ with
the \'etale topology. 
\begin{definition}\label{LS1d}
  Let $n \in \mathbb{N}$. We define a sub-presheaf
  \begin{displaymath}
G^{\mathrm{p}}_n \subset \underline{\Hom}^{\rm et}(X(n), Y(n)) ,
    \end{displaymath}
  where the right hand side denotes the Hom in the category of \'etale sheaves.
  A homomorphism $\alpha: X(n)_{U} \rightarrow Y(n)_{U}$ belongs to
  $G^{\mathrm{p}}_n(U)$ if there is a profinite \'etale covering
$\tilde{U} \rightarrow U$ and a homomorphism of $p$-divisible groups
  $\tilde{\alpha}: X_{\tilde{U}} \rightarrow Y_{\tilde{U}}$ such that
  the restriction of $\tilde{\alpha}$ to $X(n)_{\tilde{U}}$ is $\alpha_{\tilde{U}}$.

  We denote the sheafification of $G^{\mathrm{p}}_n$ by $G_n$. We define
  the prosheaf
  \begin{displaymath}
\underline{\Hom}^{\rm et}(X, Y) = \; ''\underset{\leftarrow}{\lim}'' G_n.  
  \end{displaymath}
  The limit is taken with respect to the natural restriction maps
  $G_n \rightarrow G_m$ for $n > m$. 
  \end{definition}
We note that a homomorphism of $p$-divisible groups
$\tilde{\alpha}: X_{\tilde{U}} \rightarrow Y_{\tilde{U}}$ defines a homomorphism
$\alpha: X(n)_{U} \rightarrow Y(n)_{U}$ iff  
\begin{equation}\label{LS21e}
    \mathrm{pr}^{\star}_1 \tilde{\alpha} - \mathrm{pr}^{\star}_2 \tilde{\alpha} \in
    p^n \Hom (X_{\tilde{U} \times_{U} \tilde{U}}, Y_{\tilde{U} \times_{U} \tilde{U}}). 
    \end{equation}

We consider now a banal local CM-type
$(K_{\mathfrak{p}}/F_{\mathfrak{p}}, r_{\mathfrak{p}})$. Let $E_{\mathfrak p}=E(r_{\mathfrak{p}})$ be the corresponding 
reflex field. Let  $(X, \iota_X)$ and $(Y, \iota_Y)$ be local CM-pairs over $S/ \Spf O_{E_{\mathfrak{p}}}$, which satisfy the Eisenstein condition.
As above, we define
$\underline{\Hom}^{\rm et}_{O_{K_{\mathfrak{p}}}}(X, Y)$
by replacing throughout
homomorphisms by homomorphisms of $O_{K_{\mathfrak{p}}}$-modules. The presheaf
$G^{\rm p}_n$ is now meant in this sense. 
The contracting functor (cf. Definition \ref{Cfunctor1d}) associates
$p$-adic \'etale sheaves $C_X$ and $C_Y$ with an $O_{K_{\mathfrak{p}}}$-module
structure. By Theorem \ref{Cbanal2p}, 
\begin{displaymath}
\Hom_{O_{K_{\mathfrak{p}}}} (X, Y) \cong \Hom_{O_{K_{\mathfrak{p}}}} (C_X, C_Y). 
\end{displaymath}
We set $C_{n, X} = C_X/ p^n C_X$. One checks easily by the remark after
Definition \ref{LS1d} that
\begin{equation}\label{LS13e}
  G^{\rm p}_n(U) = \Hom_{O_{K_{\mathfrak{p}}}} (C_X, C_{n, Y}) = 
  \Hom_{O_{K_{\mathfrak{p}}}} (C_{n,X}, C_{n, Y}). 
\end{equation}
In particular $G^p_n = G_n$. We conclude that, for a scheme
$S/ \Spf O_{E_{\mathfrak{p}}}$, the pro-sheaf
$\underline{\Hom}^{\rm et}_{O_{K_{\mathfrak{p}}}}(X, Y)$
is a $p$-adic \'etale sheaf. Let $\iota: \omega \rightarrow S$ be a geometric
point. Then we find for the fiber
\begin{displaymath}
  \underline{\Hom}^{\rm et}_{O_{K_{\mathfrak{p}}}}(X, Y)_{\omega} = 
  \Hom_{O_{K_{\mathfrak{p}}}}(X_{\omega}, Y_{\omega}), 
  \end{displaymath}
where the right hand side is the Hom in the cateory of $p$-divisible
$O_{K_{\mathfrak{p}}}$-modules. 

Let us assume that $S$ is a scheme over $\Spf O_{\breve{E}_{\mathfrak{p}}}$.
The contracting functor of Theorem \ref{KneunBa4p} 
associates to a CM-triple $(X, \iota_{X}, \lambda_{X})$ which satisfies the
Eisenstein condition a $p$-adic \'etale sheaf $C_{X} \in \mathrm{Et}(O_K)_S$
with an alternating form
\begin{equation}\label{LS22e}
\phi_X: C_{X} \times C_{X} \rightarrow O_{F_{\mathfrak{p}}}.
\end{equation}
We set
$\xi_X = \Trace_{F_{\mathfrak{p}}/\mathbb{Q}_p} \vartheta^{-1} \phi_X$.
In particular there is a CM-triple 
$(\mathbb{X}_{\mathfrak{p}}, \iota_{\mathbb{X}_{\mathfrak{p}}}, \lambda_{\mathbb{X}_{\mathfrak{p}}})$
over $\bar{\kappa}_{\breve{E}_{\mathfrak{p}}}$ such that
\begin{equation}\label{LS18e}
  (C_{\mathbb{X}_{\mathfrak{p}}}, \xi_{\lambda_{\mathbb{X}_{\mathfrak{p}}}}) \cong
  (\Lambda_{\mathfrak{p}}, \varsigma_{\mathfrak{p}}),  
  \end{equation}
cf. (\ref{Endo10e}).

The group $\mathbf{K}_{\mathfrak{p}}$ acts on the right hand side by similitudes.
Therefore we obtain a homomorphism
$\mathbf{K}_{\mathfrak{p}} \rightarrow \Aut_{O_{K_{\mathfrak{p}}}} \mathbb{X}_{\mathfrak{p}}$
such that the automorphisms in the image respect the polarization
$\lambda_{\mathbb{X}_{\mathfrak{p}}}$ up to a factor in $\mathbb{Z}_p^{\times}$. 

\begin{definition}\label{banK1l}
  Let $\mathfrak{p}$ be banal and let $(X, \iota_{X}, \lambda_{X})$ be a
  CM-triple on $S$ which satisfies the Eisenstein condition as above. 
  A CL-level structure\index{CL-level structure} on $(X, \iota_X, \lambda_X)$ as a class of
isomorphisms of $p$-adic \'etale sheaves 
\begin{equation}\label{LS2e}
   (\Lambda_{\mathfrak{p}}, \varsigma_{\mathfrak{p}}) \isoarrow (C_X, \xi_{\lambda_X})
  \; \mod \mathbf{K}_{\mathfrak{p}}^{\star} ,
  \end{equation}
which respect the bilinear forms on both sides up to a factor in 
$\mathbb{Z}_p^{\times}$.
We will write
\begin{equation}\label{LS2ee}
   \mathbb{X}_{\mathfrak{p}} \isoarrow X \; \mod \mathbf{K}_{\mathfrak{p}}^{\star}. 
  \end{equation}
   for a CL-structure.
\end{definition}
More precisely this means the following. Let $M \geq 1$ such that (\ref{banK2e}) holds. Then
a CL-level structure is a right 
$\mathbf{K}_{\mathfrak{p}}^{\star}/(\mathbf{K}_{\mathfrak{p}}(p^M))$-torsor
\begin{displaymath}
  T \subset \underline{\Isom}_{O_{K_{\mathfrak{p}}}}(\Lambda_{\mathfrak{p}} \otimes
  \mathbb{Z}/p^M\mathbb{Z}, C_X \otimes \mathbb{Z}/p^M\mathbb{Z})
  \end{displaymath}
such that the inclusion is equivariant with respect to the right actions
of $\mathbf{K}_{\mathfrak{p}}^{\star}/(\mathbf{K}_{\mathfrak{p}}(p^M))$ on both sides
and such that the local sections of $T$ respect the bilinear forms on
$\Lambda_{\mathfrak{p}}$ and $C_X$ up to a factor in
$(\mathbb{Z}/p^M\mathbb{Z})^{\times}$. If $S$ is connected and
$\omega \rightarrow S$ is a geometric point a CL-structure is given by a
$\mathbf{K}^{\star}_{\mathfrak{p}}$-orbit of an isomorphism
$\Lambda_{\mathfrak{p}} \rightarrow (C_X)_{\omega}$ which respects the bilinear
forms on both sides by a factor in $\mathbb{Z}_p^{\times}$ and such that the
orbit is preserved by the action of $\pi_1(S, \omega)$. This explains the
notation (\ref{LS2ee}).

Let $(X, \iota, \lambda)$ be a semi-local CM-triple relative to
$(K \otimes \mathbb{Q}_p/F \otimes \mathbb{Q}_p, r)$ over a scheme
$S\in ({\rm Sch}/\Spf O_{\breve{E}_{\nu}})$, cf. the beginning of section
\ref{RZtildeM}. We set  \index[NO]{XAD@$X^{{\rm ba}}$}
$X^{{\rm ba}} = \prod_{\mathfrak{p}, \; \text{banal}} X_{\mathfrak{p}}$.
We choose $(\mathbb{X}, \iota_{\mathbb{X}}, \lambda_{\mathbb{X}})$ as in
section \ref{RZtildeM}. Then (\ref{LS18e}) holds. 
From this we obtain an action of $\mathbf{K}_{p}^{\star, {\rm ba}}$ on
$\mathbb{X}^{{\rm ba}}$
which respects the polarization
$\prod_{\mathfrak{p}, \, \text{banal}} \lambda_{\mathbb{X}_{\mathfrak{p}}}$ up to a 
factor in $\mathbb{Z}_p^{\times}$.

We define a CL-level structure on $\mathbb{X}^{{\rm ba}}$ modulo
$\mathbf{K}_{p}^{\star, {\rm ba}}$ as a CL-level stuctures
$\bar{\eta}_{\mathfrak{p}}:\mathbb{X}_{\mathfrak{p}}\isoarrow X_{\mathfrak{p}} \; \mod \mathbf{K}_{\mathfrak{p}}^{\star}$,
for each banal $\mathfrak{p}$ which respect the bilinear forms up to a factor
in $\mathbb{Z}_p^{\times}$ that is independent of $\mathfrak{p}$. 
\begin{definition}\label{banK2d}
  With the notations of Definition \ref{KneunU2d},  let $i \in \mathbb{Z}$.  
  Let $\mathcal{M}_{\mathbf{K}_{p}^{\star}}(i)$ be the following functor on the
  category of schemes $S$ over $\Spf O_{\breve{E}_{\nu}}$. We will write
  $\bar{S} = S \times_{\Spf O_{\breve{E}_{\nu}}} \Spec \bar{\kappa}_{E_{\nu}}$. 
  A point of $\mathcal{M}_{\mathbf{K}_{p}^{\star}}(i)(S)$ is given
  by the following data: 

  \begin{enumerate}
  \item[(1)] A CM-triple $(X, \iota, \lambda)$ of type
    $(K \otimes \mathbb{Q}_p / F \otimes \mathbb{Q}_p, r)$ over $S$
    which satisfies the conditions $({\rm KC}_r)$ and $({\rm EC}_r)$ and is
    compatible with $(V,\varsigma)$. 
  \item[(2)] A $O_K\otimes \mathbb{Z}_p$-linear quasi-isogeny
    \begin{displaymath}
  \rho: \bar{X} := X \times_{S} \bar{S} \longrightarrow \mathbb{X}
      \times_{\Spec \kappa_{\breve{E}_{\nu}}} \bar{S}  
    \end{displaymath}
     such that $\rho$ respects the polarization $p^i \lambda$ on $X$ and
    $\lambda_{\mathbb{X}}$ up to a factor in $\mathbb{Z}_p^{\times}$. 
\item[(3)] Let $X^{{\rm ba}} = \prod_{\mathfrak{p}\; \text{banal}} X_{\mathfrak{p}}$. A
  CL-level structure 
  \begin{displaymath}
    \bar{\eta}: \mathbb{X}^{{\rm ba}} \isoarrow X^{{\rm ba}} \;
    \mathrm{mod}\; \mathbf{K}_{p}^{\star, {\rm ba}}. 
  \end{displaymath}
      \end{enumerate}
  We  set\index[NO]{MBJ@$\tilde{\mathcal{M}}_{\mathbf{K}_{p}^{\star}}$}
  \begin{displaymath}
    \tilde{\mathcal{M}}_{\mathbf{K}_{p}^{\star}} = \coprod_{i \in \mathbb{Z}}
    \mathcal{M}_{\mathbf{K}_{p}^{\star}}(i).
    \end{displaymath}
Two data $(X_1, \iota_1, \lambda_1, \rho_1, \bar{\eta}_1)$ and
  $(X_2, \iota_2, \lambda_2, \rho_2, \bar{\eta}_2)$
define the same point iff there is an isomorphism of $O_K$-modules
$\alpha: X_1 \rightarrow X_2$ such that
$\rho_2 \circ \alpha_{\bar{S}} = \rho_1$ and such that the level structures are
respected (in particular, $\alpha$ respects
the polarizations up to a factor in $\mathbb{Z}_p^{\times}$).  
\end{definition}
We note that in the case without level structure we used a different version
of the functor $\mathcal{M}_r(i)$, cf. Remark \ref{KneunU2dRm}. We prefer here
to consider $\mathbb{Z}_p^{\times}$-homogeneous polarizations. We prove directly
a version of Corollary \ref{KneunU1c} in this setting. Since we assume that
banal places exist,  we do not need the group $\hat{G}$.     
\begin{proposition}\label{LS1p}
 There exists an isomorphism
         \begin{displaymath}
    \tilde{\mathcal{M}}_{\mathbf{K}_{p}^{\star}} \overset{\sim}{\longrightarrow}
    (\hat{\Omega}_{F_{v}} \times_{\Spf O_{F_{v}}} \Spf O_{\breve{E}_{\nu}}) \times
    G^{{\rm ba}}(\mathbb{Q}_p)/\mathbf{K}_{p}^{\star, {\rm ba}} ,
        \end{displaymath}
which is equivariant with respect to the action of $J(\mathbb{Q}_p)$
on both sides. 
\end{proposition}
\begin{proof}
  At the banal places we may use lisse $p$-adic
  \'etale sheaves to describe a point of
$\tilde{\mathcal{M}}_{\mathbf{K}_{p}^{\star}}(S)$. A point consists of
a CM-triple $(X_{\mathfrak{p}_{v}},\iota_{\mathfrak{p}_{v}},\lambda_{\mathfrak{p}_{v}})$
and a quasi-isogeny
$\rho_{\mathfrak{p}_{v}}: X_{\mathfrak{p}_{v},\bar{S}}\rightarrow\mathbb{X}_{\mathfrak{p}_{v}}$ 
such that
$\rho^{\star}_{\mathfrak{p}_v} (\lambda_{\mathbb{X}_{\mathfrak{p}_v}}) = u p^i \lambda_{X_{\mathfrak{p}_v}}$ 
for some $u \in \mathbb{Z}_p^{\times}$, $i \in \mathbb{Z}$ and an  isomorphism 
of lisse $p$-adic \'etale sheaves on $\bar{S}$, 
\begin{equation}\label{LS19e}
  (C_{\mathbb{X}^{{\rm ba}}}, \xi_{\mathbb{X}^{{\rm ba}}}) \overset{\eta}{\longrightarrow} (C, \xi)
  \subset (C_{\mathbb{X}^{{\rm ba}}}, \xi_{\mathbb{X}^{{\rm ba}}}) \otimes \mathbb{Q}, 
  \end{equation}
where $\eta$ respects the alternating forms up to a factor in
$\mathbb{Z}_p^{\times}$ and such that the restriction of $\xi_{\mathbb{X}^{{\rm ba}}}$
with respect to the last inclusion is equal to $u p^i \xi$ with the same $u$
and $i$ as above.
By (\ref{Endo10e}) we have
\begin{equation}\label{trivCban}
  (C_{\mathbb{X}^{{\rm ba}}}, \xi_{\mathbb{X}^{\rm ba}}) \cong
  (\Lambda^{{\rm ba}}, \varsigma^{{\rm ba}}),
\end{equation}
where the right hand side is the orthogonal direct sum over all
$(\Lambda_{\mathfrak{p}}, \varsigma_{\mathfrak{p}})$ for $\mathfrak{p}$ banal. 

We denote by $\tilde{\mathcal{M}}^{{\rm ba}}_{\mathbf{K}_{p}^{\star}}$ 
\index[NO]{MBK@$\tilde{\mathcal{M}}^{{\rm ba}}_{\mathbf{K}_{p}^{\star}}$} the moduli functor
described by the data (\ref{LS19e}).
We claim that there is a natural isomorphism
\begin{equation}\label{LS20e}
  \tilde{\mathcal{M}}^{{\rm ba}}_{\mathbf{K}_{p}^{\star}} \cong
  G^{{\rm ba}}(\mathbb{Q}_p)/\mathbf{K}_{p}^{\star, {\rm ba}}. 
  \end{equation}
Indeed, the group $G^{{\rm ba}}(\mathbb{Q}_p)$ acts
naturally on this functor: Let $g \in G^{{\rm ba}}(\mathbb{Q}_p)$ such that  
\begin{displaymath}
  \phi_{\mathbb{X}}(g x_1, g x_2) = u' p^{j} \xi_{\mathbb{X}}(x_1, x_2). 
\end{displaymath}
Then $g$ maps (\ref{LS19e}) to 
\begin{displaymath}
  (C_{\mathbb{X}^{{\rm ba}}}, \xi_{\mathbb{X}^{{\rm ba}}}) \overset{g \eta}{\longrightarrow}
  (g C, (1/u' p^{j}) \xi) 
    \subset C_{\mathbb{X}^{{\rm ba}}} \otimes \mathbb{Q}.  
  \end{displaymath}
If we have an arbitrary point (\ref{LS19e}), then the composite of the arrow
with the inclusion is an element $g \in G^{{\rm ba}}(\mathbb{Q}_p)$ and therefore
(\ref{LS19e}) is isomorphic to
\begin{displaymath}
  C_{\mathbb{X}^{{\rm ba}}} \overset{g}{\longrightarrow}
  (g C, (1/u' p^{j}) \xi_{\mathbb{X}^{{\rm ba}}}) 
    \subset C_{\mathbb{X}^{{\rm ba}}} \otimes \mathbb{Q}.
  \end{displaymath}
We see that the action is transitive and that the stabilizer of the base point
\begin{displaymath}
  (C_{\mathbb{X}^{{\rm ba}}}, \xi_{\mathbb{X}^{{\rm ba}}}) \overset{\id}{\longrightarrow} 
  (C_{\mathbb{X}^{{\rm ba}}}, \xi_{\mathbb{X}^{{\rm ba}}}) \subset C_{\mathbb{X}^{{\rm ba}}}
  \otimes \mathbb{Q}  
\end{displaymath}
is $\mathbf{K}_{p}^{\star, {\rm ba}}$. This shows (\ref{LS20e}).

Now we fix $i \in \mathbb{Z}$. We denote by $\mathcal{M}_{r_{\mathfrak{p}_v}}$ the
functor of section 6 associated to the special local CM-type
$(K_{\mathfrak{p}_v}/F_{\mathfrak{p}_v}, r_{\mathfrak{p}_v})$. There is the natural
injection of functors
\begin{displaymath}
  \mathcal{M}_{\mathbf{K}_{p}^{\star}}(i) \rightarrow \mathcal{M}_{r_{\mathfrak{p}_v}}(i)
  \times \tilde{\mathcal{M}}^{{\rm ba}}_{\mathbf{K}_{p}^{\star}}(i).
\end{displaymath}
We claim that this map is surjective. Indeed, assume we are given a point
$(X_{\mathfrak{p}_{v}},\iota_{\mathfrak{p}_{v}},\lambda_{\mathfrak{p}_{v}}, \rho_{\mathfrak{p}_v})$, where 
$\rho^{\star}_{\mathfrak{p}_v} (\lambda_{\mathbb{X}_{\mathfrak{p}_v}}) = u_1 p^i \lambda_{\mathfrak{p}_v}$, with $u_1 \in \mathbb{Z}_p^{\times}$,  from the first factor on the right
hand side, and a point $(C, \xi) \subset C_{\mathbb{X}^{{\rm ba}}}$ (endowed with $\eta$), 
where $u_2 p^i \xi = \xi_{\mathbb{X}^{{\rm ba}}}$ with $u_2 \in \mathbb{Z}_p^{\times}$, from the second factor.
These two points form a point of $\mathcal{M}_{\mathbf{K}_{p}^{\star}}(i)$ iff
$u_1 = u_2$ because only then the condition (2) of  Definition \ref{banK2d}
is fulfilled for the resulting polarization on
$X = X_{\mathfrak{p}_{v}} \times X^{\mathrm{ba}}$.  
But in the point from the first factor we can replace
$\lambda_{\mathfrak{p}_{v}}$ by $(u_1/u_2)\lambda_{\mathfrak{p}_{v}}$ without changing
the isomorphism class of this point.

Therefore the surjectivity holds and
the Proposition follows as Proposition \ref{compfo}.
\end{proof}

Let us fix an open and compact subgroup
$\mathbf{K}^p \subset G(\mathbb{A}_f^p)$. We set 
$\mathbf{K}^{\star} = \mathbf{K}^{\star}_p \mathbf{K}^p$ and 
$\mathbf{K} = \mathbf{K}_p \mathbf{K}^p$ as after (\ref{boldK1e}).
We choose $(\mathbb{X}, \iota_{\mathbb{X}}, \lambda_{\mathbb{X}})$ as above. 

\begin{definition}\label{balevel2d}
  We define a functor $\hat{\mathcal{A}}^{\star}_{\mathbf{K}^{\star}}$ on the
  category of schemes $S$ over $\Spf O_{\breve{E}_{\nu}}$.
  A point of $\hat{\mathcal{A}}^{\star}_{\mathbf{K}^{\star}}(S)$ consists of
  the following data:
\begin{enumerate}
\item[(1)] a point $(A, \iota, \bar{\lambda}, \bar{\eta}^p)$ of
  $\mathcal{A}_{\mathbf{K}}(S)$,
\item[(2)]   a CL-level structure 
  \begin{displaymath}
\bar{\eta}_p: \mathbb{X}^{{\rm ba}} \rightarrow A[p^\infty]^{{\rm ba}} \;
    \mathrm{mod}\; \mathbf{K}_{p}^{\star, {\rm ba}}. 
    \end{displaymath}
\end{enumerate}
\end{definition}
We denote here by $A[p^\infty]^{{\rm ba}}$\index[NO]{AAF@$A[p^\infty]^{{\rm ba}}$}
the banal part of the $p$-divisible group of $A$ with its structure of a
semi-local CM-triple. The morphism
$\hat{\mathcal{A}}^{\star}_{\mathbf{K}^{\star}} \rightarrow \hat{\mathcal{A}}_{\mathbf{K}}$
is a finite \'etale covering of formal schemes. Since we assume that
$\mathbf{K}^p$ is small enough, $\mathcal{A}_{\mathbf{K}}$ is a proper scheme
over $\Spec O_{\breve{E}_{\nu}}$. By the algebraization theorem, there is a unique
finite \'etale morphism of schemes over $\Spec O_{\breve{E}_{\nu}}$
\begin{equation}\label{balevel2e}
\mathcal{A}^{\star}_{\mathbf{K}^{\star}} \rightarrow \mathcal{A}_{\mathbf{K}}\times_{\Spec\, O_{E, ({\mathfrak p}_\nu)}}\Spec\, O_{\breve E_\nu} , 
\end{equation}
such that the $p$-adic completion of $\mathcal{A}^{\star}_{\mathbf{K}^{\star}}$ is
$\hat{\mathcal{A}}^{\star}_{\mathbf{K}^{\star}}$.
\begin{remark}
We note that the scheme  $\mathcal{A}^{\star}_{\mathbf{K}^{\star}}$ is defined over the ring $O_{\breve E_\nu}$, which is not of finite type over $\BZ_p$. As mentioned after Theorem \ref{ThmC}, this scheme is closely related to an integral model of a Shimura variety which is a central twist of ${\rm Sh}_{\bK^\star}(G, \{h\})$. 
\end{remark}

Recall the projective scheme $\CA_{\bK^\star, E}$ over $E$ from section \ref{ss:globalres} (the canonical model of ${\rm Sh}_{\bK^\star}(G, \{h\})$), comp. the proof of Proposition \ref{proprepr}. We will now relate $\CA_{\bK^\star, E}$ with the general fiber
$\mathcal{A}^{\star}_{\mathbf{K}^{\star}}\times_{\Spec O_{\breve{E}_{\nu}}}\Spec\breve{E}_{\nu}$. 
 We start with a reformulation of the
level structure $\bar{\eta}_p$ in Definition \ref{balevel2d}.

We assume that $S$ is a scheme over $\Spf O_{\breve{E}_{\mathfrak{p}}}$, i.e., we pass
to the completion of the maximal unramified extension of $E_{\mathfrak{p}}$. 
We consider now a  polarized local CM-pair $(X, \iota_X, \lambda_X)$ over $S$
of CM-type $(K_{\mathfrak{p}}/F_{\mathfrak{p}}, r_{\mathfrak{p}})$, cf. Definition
\ref{KatCMtriple1d}. We will always assume that the Eisenstein conditions are
satisfied. By  Theorem \ref{KneunBa4p}, $\lambda_X$ is described by a
$O_{F_{\mathfrak{p}}}$-bilinear form $\phi_{\lambda_X}$, or also 
\begin{displaymath}
\xi_{\lambda_X}: C_X \times C_X \rightarrow \mathbb{Z}_p ,
\end{displaymath}
as defined after (\ref{LS22e}).
Equivalently, we can consider the $O_{K_{\mathfrak{p}}}$-anti-hermitian form
\begin{equation}\label{LS23e}
\varkappa_{\lambda_X}: C_X \times C_X \rightarrow K_{\mathfrak{p}}    ,
  \end{equation}
which is defined by
\begin{displaymath}
  \Trace_{K_{\mathfrak{p}}/F_{\mathfrak{p}}} a \varkappa_{\lambda_X}(c_1, c_2) =
  \phi_{\lambda_X}(ac_1, c_2),
  \quad a \in O_{K_{\mathfrak{p}}}, \; c_1, c_2 \in C_X. 
\end{displaymath}
Then $\varkappa_{\lambda_X}$ is $O_{K_{\mathfrak{p}}}$-linear in the first variable and
$O_{K_{\mathfrak{p}}}$-anti-linear in the second variable.

 If we define
$(\mathbb{X}, \iota_{\mathbb{X}}, \lambda_{\mathbb{X}})$ by
$(C_{\mathbb{X}}, \xi_{\lambda_{\mathbb{X}}}) \cong (\Lambda_{\mathfrak{p}}, \varsigma_{\mathfrak{p}})$ 
(cf. (\ref{LS18e})), we can reformulate (\ref{LS2e}):  A CL-structure is a
class of isomorphisms
\begin{displaymath}
  (\mathbb{X}, \lambda_{\mathbb{X}}) \rightarrow (X, \lambda_X)
  \; \mod \mathbf{K}_{\mathfrak{p}}^{\star} 
\end{displaymath}
which respect the polarizations up to a factor in $\mathbb{Z}_p^{\times}$. 
This agrees with Definition \ref{balevel2d}.

Let $\bar{\kappa}_{E_\mathfrak{p}}$ be the residue class field of
$\breve{E}_{\mathfrak{p}}$. We will consider CM-pairs $(Z, \iota_Z)$ of CM-type
$r_{\mathfrak{p}}/2$ over a scheme $S/ \Spf O_{\breve{E}_{\mathfrak{p}}}$.
Then $Z$ is a $p$-divisible group of height $2d_{\mathfrak{p}}$
and dimension $d_{\mathfrak{p}}$, where
$d_{\mathfrak{p}} = [K_{\mathfrak{p}} : F_{\mathfrak{p}}]$.
We will always assume that the Eisenstein condition is fulfilled. 
Proposition \ref{KottwitzC2p} continues to hold with the same polynomials
$\mathbf{E}_{A_{\psi}}$. The functor $C_Z$ (cf. Definition \ref{Cfunctor1d})
exists for
local CM-pairs of type $(K_{\mathfrak{p}}/F_{\mathfrak{p}}, r_{\mathfrak{p}}/2)$. 

We will reformulate CL-level structures as suggested by \cite{RSZ}.
There is up to isomorphism a unique CM-pair
$(\bar{X}_0, \bar{\iota}_0)$ of CM-type $r_{\mathfrak p}/2$ over $\bar\kappa_{E_{\mathfrak p}}$. It lifts
uniquely to a CM-pair $(X_0, \iota_0)$ over
$O_{\breve{E}_{\mathfrak{p}}}$, and 
\begin{equation}\label{LS4e}
C_{X_0} \cong O_{K_{\mathfrak{p}}} 
  \end{equation}
is the constant $p$-adic sheaf. We consider biextensions 
\begin{displaymath}
\beta: X_{0} \times X_{0} \rightarrow \hat{\mathbb{G}}_m   
\end{displaymath}
or, equivalently, bilinear forms of  displays as in Proposition
\ref{Cbanal3p}. They are in bijection with bilinear forms
\begin{equation}\label{LS5e}
\phi: C_{X_0} \times C_{X_0} \rightarrow O_{F_{\mathfrak{p}}}.  
\end{equation}
Equivalently we use $\xi=\xi_\phi$ or $\varkappa=\varkappa_\phi$ as before (\ref{LS23e}). 

We define
$\varkappa_0: C_0 \times C_0 \rightarrow O_{K_{\mathfrak{p}}}$ using
(\ref{LS4e}), by
\begin{equation}\label{LS6e}
  \varkappa_0(x,y) = x \bar{y}. \quad x, y \in O_{K_{\mathfrak{p}}}. 
\end{equation}
 We denote by
$\mathfrak{s}_0: X_0 \rightarrow X_0^{\wedge}$ the homomorphism associated to
$\phi_0$. This homomorphism is symmetric. 

We note that there is a principal polarization $\lambda$ on $X_0$. It is a generator of the free $O_{K_{\mathfrak{p}}}$-module of rank one $\Hom_{O_{K_{\mathfrak{p}}}}(X_0, X_0^\wedge)$. In
the case where $K_{\mathfrak{p}}/F_{\mathfrak{p}}$ is ramified, the corresponding
form under the bijection  (\ref{LS5e}) is
\begin{displaymath}
\phi_{\lambda}(x, y) = \Trace_{K_{\mathfrak{p}}/F_{\mathfrak{p}}} \Pi^{-1} x \bar{y}.
\end{displaymath}
Then $\mathfrak{s}_0=\lambda \Pi$. In the case where $K_{\mathfrak{p}}/F_{\mathfrak{p}}$ is unramified,
we choose a unit $\varepsilon \in O_{K_{\mathfrak{p}}}$ such that
$\varepsilon + \bar{\varepsilon} = 0$. Then the corresponding
form under the bijection  (\ref{LS5e}) is
\begin{displaymath}
\phi_{\lambda}(x, y) = \Trace_{K_{\mathfrak{p}}/F_{\mathfrak{p}}} \varepsilon^{-1} x \bar{y}.
\end{displaymath}
Then $\mathfrak{s}_0=\lambda \varepsilon$.

  Let $S$ be a $p$-adic formal scheme over $\Spf O_{\breve{E}_{\mathfrak{p}}}$. Let
  $(X, \iota_X, \lambda_X)$ be a polarized CM-pair  of type
  $(K_{\mathfrak{p}}/F_{\mathfrak{p}}, r_{\mathfrak{p}})$ which satisfies the Eisenstein
  condition as always required. We endow
$\underline{\Hom}_{O_{K_{\mathfrak{p}}}}^{\rm et}(X_0, X)$
 with an
 $O_{K_{\mathfrak{p}}}$-anti-hermitian form with values in $K_{\mathfrak{p}}$.
 Let $u_1, u_2 \in C_n(U)$. They are given by homomorphisms
 $\tilde{u}_1, \tilde{u}_2: X_0 \rightarrow X$ which are defined over a
 profinite \'etale covering $\tilde{U} \rightarrow U$. We consider the
 homomorphism 
  \begin{displaymath}
    \tilde{u}_2^{\wedge} \lambda_X \tilde{u}_1: X_0 \rightarrow X
    \rightarrow X^{\wedge} \rightarrow X_0^{\wedge}.     
  \end{displaymath}
This element of $\Hom_{O_{K_{\mathfrak{p}}}}((X_0)_{\tilde{U}}, (X_0^{\wedge})_{\tilde{U}})$
may be written as
  \begin{equation}\label{LS7e}
    \tilde{u}_2^{\wedge} \lambda_X \tilde{u}_1 =
    \tilde{\delta}(\tilde{u}_1, \tilde{u}_2) \mathfrak{s}_0 , 
  \end{equation}
  with some constant
  $\tilde{\delta}(\tilde{u}_1, \tilde{u}_2) \in K_{\mathfrak{p}}$. 
  In the ramified case,
  \begin{displaymath}
    \delta(u_1, u_2) := \tilde{\delta}(\tilde{u}_1, \tilde{u}_2)
    \mod p^n \Pi^{-1} O_{O_{K_{\mathfrak{p}}}} 
  \end{displaymath}
  is well defined. In the unramified case,  the element $\tilde{\delta}(\tilde{u}_1, \tilde{u}_2)$ is well-defined
  modulo $ p^n O_{K_{\mathfrak{p}}}$. 
  Varying $n$, we therefore  obtain a bilinear form 
  \begin{displaymath}
 \delta: \underline{\Hom}_{O_{K_{\mathfrak{p}}}}^{\rm et}(X_0, X) \times
    \underline{\Hom}_{O_{K_{\mathfrak{p}}}}^{\rm et}(X_0, X) \rightarrow K_{\mathfrak{p}}.  
    \end{displaymath}
  This is a $O_{K_{\mathfrak{p}}}$-anti-hermitian form. We set
  \begin{displaymath}
    \mathfrak{e} = \Trace_{F_{\mathfrak{p}}/\mathbb{Q}_p} \Trace_{K_{\mathfrak{p}}/F_{\mathfrak{p}}}
    \vartheta_{F_{\mathfrak{p}}/\mathbb{Q}_p}^{-1} \delta.
    \end{displaymath}
  Then $\mathfrak{e}$ is an alternating form
  \begin{displaymath}
\mathfrak{e}: \underline{\Hom}_{O_{K_{\mathfrak{p}}}}^{\rm et}(X_0, X) \times
    \underline{\Hom}_{O_{K_{\mathfrak{p}}}}^{\rm et}(X_0, X) \rightarrow \mathbb{Z}_p. 
  \end{displaymath}
  which satisfies $\mathfrak{e}(a u_1, u_2) = \mathfrak{e}(u_1, \bar{a} u_2)$, $a\in O_{K_{\mathfrak p}}$. 
  \begin{proposition}\label{LS3p}
    A CL-level structure on a polarized CM-pair $(X, \iota_X, \lambda_X)$
     of type $(K_{\mathfrak{p}}/F_{\mathfrak{p}}, r_{\mathfrak{p}})$ over the
    $p$-adic formal scheme $S$ can equivalently be given as a class of
    isomorphisms of $p$-adic \'etale sheaves 
    \begin{equation}\label{LS24e}
      \eta: (\Lambda_{\mathfrak{p}}, \varsigma_{\mathfrak{p}}) \isoarrow
      (\underline{\Hom}_{O_{K_{\mathfrak{p}}}}^{\rm et}(X_0, X), \mathfrak{e}) \;       
      \mod \mathbf{K}_{\mathfrak{p}}^{\star},  
    \end{equation}
    which respect the bilinear forms on both sides up to a constant in
    $\mathbb{Z}_p^{\times}$. 
  \end{proposition}
  \begin{proof}
    Indeed, we apply the contracting functor to the right hand side of the
    isomorphism \eqref{LS24e}. We view
    $\tilde{u}_1, \tilde{u}_2$ from (\ref{LS7e}) as homomorphisms
    \begin{displaymath}
\tilde{u}_i: O_{K_{\mathfrak{p}}} = C_{X_0} \rightarrow C_X.  
      \end{displaymath}
Let $\varkappa_{\lambda_X}: C_X \times C_X \rightarrow K_{\mathfrak{p}}$ be the
anti-hermitian form induced by $\lambda_X$. The definition (\ref{LS7e}) of the
sesqui-linear form $\delta$ which gives rise to $\mathfrak{e}$,
reads in terms of the contracting functor  as
defined by (\ref{LS7e}) 
  \begin{equation}\label{LS8e}
    \varkappa_{\lambda_X}(\tilde{u}_1(x), \tilde{u}_2(y)) =
    \tilde{\delta}(\tilde{u}_1, \tilde{u}_2) x \bar{y}. 
  \end{equation}
  If we identify
  \begin{displaymath}
    \underline{\Hom}_{O_{K_{\mathfrak{p}}}}^{\rm et}(X_0, X) =
    \underline{\Hom}_{O_{K_{\mathfrak{p}}}}^{\rm et}(O_{K_{\mathfrak{p}}}, C_X) = C_X  
  \end{displaymath}
  by sending $\tilde{u}$ to $\tilde{u}(1)$, the form $\tilde{\delta}$ is mapped
  to the form $\varkappa_{\lambda_X}$. This is immediate by setting $x = y = 1$  in
  (\ref{LS8e}).
 Therefore we have
  identified the right hand side of (\ref{LS24e}) with $(C_X, \xi_{\lambda_X})$.
  This proves the assertion. 
\end{proof}

   \begin{proposition}\label{LS5p}
Let $S$ be a flat proper scheme over $\Spec O_{\breve{E}_{\mathfrak{p}}}$. Let
    $(X, \iota_{X})$ and $(Y, \iota_Y)$ be CM-pairs of type
    $(K_{\mathfrak{p}}/F_{\mathfrak{p}}, r_{\mathfrak{p}})$ or
$(K_{\mathfrak{p}}/F_{\mathfrak{p}}, r_{\mathfrak{p}}/2)$ over $S$.
Let $U \rightarrow S$ be a finite \'etale covering, and 
$\hat{U} \rightarrow \Spf O_{\breve{E}_{\mathfrak{p}}}$ its formal completion
along the special fibre. Set
$U_{\eta} = U \times_{\Spec O_{\breve{E}_{\mathfrak{p}}}} \Spec \breve{E}_{\mathfrak{p}}$. 
Then there is a natural
bijective homomorphism
\begin{equation}\label{LS15e}
  G_n(\hat{U}) \rightarrow
  \Hom_{O_{K_{\mathfrak{p}}}} (X(n)_{U_{\eta}}, Y(n)_{U_{\eta}}).  
  \end{equation}
In particular, the $p$-adic \'etale sheaf
$\underline{\Hom}_{O_{K_{\mathfrak{p}}}}(T_p(X_{S_{\eta}}), T_p(Y_{S_{\eta}}))$
is unramified along the special fibre of $S$. 
  \end{proposition}
  \begin{proof}
We consider the natural map
$G_n(\hat{U}) \rightarrow  \Hom_{O_{K_{\mathfrak{p}}}}(X(n)_{\hat{U}}, Y(n)_{\hat{U}})$.  
By Grothendieck's existence theorem (EGA III, Thm. 5.1.4), the target
of this arrow coincides with $\Hom_{O_{K_{\mathfrak{p}}}}(X(n)_U, Y(n)_U)$.
If we restrict the last set of homomorphisms to the generic fibre
we obtain the map (\ref{LS15e}).

The injectivity of (\ref{LS15e}) follows from the definition of $G^n$. 
To prove  surjectivity, we can assume that $U$ is
connected. By Grothendieck's existence theorem we find a finite connected 
\'etale covering $U_1 \rightarrow U$ such that the sheaves
$C_{n,X_{\hat{S}}}$ and $C_{n,Y_{\hat{S}}}$ become trivial over $\hat{U}_1$.

We write the proof only in the case where $X$ and $Y$ are of CM-type
$r_{\mathfrak{p}}$. The cases where $r_{\mathfrak{p}}/2$ appears will be obvious.   
By the choice of $U_1$, we deduce the isomorphism 
    \begin{displaymath}
G_n(\hat{U}_1) \cong \Hom_{O_{K_{\mathfrak{p}}}} ((O_{K_{\mathfrak{p}}}/p^n O_{K_{\mathfrak{p}}})^2,
(O_{K_{\mathfrak{p}}}/p^n O_{K_{\mathfrak{p}}})^2). 
    \end{displaymath}
    We choose a geometric point $\omega$ of $(U_1)_{\eta}$. Then we obtain
    injective homomorphisms
    \begin{displaymath}
      \begin{array}{c} 
G_n(\hat{U}_1) \rightarrow
\Hom_{O_{K_{\mathfrak{p}}}} (X(n)_{U_{1,\eta}}, Y(n)_{U_{1,\eta}}) = \\[2mm]
\underline{\Hom}_{O_{K_{\mathfrak{p}}}} (T_p(X_{\eta}) \otimes \mathbb{Z}/(p^n),
T_p(Y_{\eta}) \otimes \mathbb{Z}/(p^n))(U_{1, \eta}) \longrightarrow \\[2mm]  
\Hom_{O_{K_{\mathfrak{p}}}} (T_p(X_{\omega}), T_p(Y_{\omega})) \otimes \mathbb{Z}/(p^n)
\cong \Hom_{O_{K_{\mathfrak{p}}}} ((O_{K_{\mathfrak{p}}})^2, (O_{K_{\mathfrak{p}}})^2) \otimes
\mathbb{Z}/(p^n).  
        \end{array}
      \end{displaymath}
Since we have the same number of elements on both sides,  the arrows are
bijective. In particular this shows that the \'etale sheaf
$\underline{\Hom}_{O_{K_{\mathfrak{p}}}} (T_p(X_{\eta}), T_p(Y_{\eta}))
\otimes \mathbb{Z}/(p^n)$ 
becomes trivial over the finite \'etale covering
$U_{1,\eta} \rightarrow S_\eta$. Therefore it is unramified along the special
fibre of $S$.

Finally, we obtain the bijectivity of (\ref{LS15e}) by
exploiting the sheaf property with respect to the covering
 \begin{displaymath}
      U_1 \times_{U} U_1 
      \begin{array}{c} 
        \rightarrow\\[-2mm] 
        \rightarrow 
        \end{array}
      U_1 \rightarrow U. 
    \end{displaymath}
  \end{proof}
 \begin{corollary}\label{padicetale}
    With the assumptions of the last Proposition, there is a $p$-adic \'etale
    sheaf $\underline{\Hom}_{O_{K_{\mathfrak{p}}}}(X, Y)$ on $S$ whose restriction
    to the special fibre 
    $S \times_{\Spec O_{\breve{E}_{\mathfrak{p}}}} \Spec \bar\kappa_{{E}_{\mathfrak{p}}}$ is
    $\underline{\Hom}^{\rm et}_{O_{K_{\mathfrak{p}}}}(X_{\bar\kappa_{{E}_{\mathfrak{p}}}}, Y_{\bar\kappa_{{E}_{\mathfrak{p}}}})$
    and whose restriction to the general  fibre
$S \times_{\Spec O_{\breve{E}_{\mathfrak{p}}}} \Spec \breve{E}_{\mathfrak{p}}$  is
$\underline{\Hom}_{O_{K_{\mathfrak{p}}}} (T_p(X_{\breve{E}_{\mathfrak{p}}}), T_p(Y_{\breve{E}_{\mathfrak{p}}}))$.   
  \end{corollary}
  \begin{proof}
    The sheaves $G_n$ over $\hat{S}$ are representable by  finite \'etale
    morphisms of formal schemes. They come therefore from finite \'etale
    morphisms $G^{\rm al}_n \rightarrow S$. We have to compare the general fibre of
    $G^{\rm al}_n$ with
$\underline{\Hom}_{O_{K_{\mathfrak{p}}}} (T_p(X_{\eta}), T_p(Y_{\eta})) \otimes \mathbb{Z}/(p^n)$.   

We have shown that both sheaves are trivialized by a finite \'etale covering
$S_1 \rightarrow S$. The homomorphism (\ref{LS15e}) gives a canonical
isomorphism between these sheaves with constant \'etale sheaves on
$S_1 \times_{\Spec O_{\breve{E}_{\mathfrak{p}}}} \Spec \kappa_{\breve{E}_{\mathfrak{p}}}$.
Finally, we consider descent for the general fiber of the covering 
    \begin{displaymath}
      S_1 \times_{S} S_1 
      \begin{array}{c} 
       \overset{p_1}{\rightarrow}\\[-2mm] 
        \underset{p_2}{\rightarrow} 
        \end{array}
      S_1 \rightarrow S 
    \end{displaymath}
We see that the descent data for the two sheaves agree since they are induced
from the descent datum on the \'etale sheaf $\underline{\Hom}(X(n), Y(n))$. 
  \end{proof}
  We now go back to the Definition \ref{balevel2d}. Let
  $\Lambda^{\rm ba}=\prod_{\mathfrak{p}, \text{banal}}\Lambda_{\mathfrak{p}}$.
  \index[NO]{ZZKC@$\Lambda^{\rm ba}$}   We choose for each banal
  $\mathfrak{p}$ a CM-pair $(X_{\mathfrak{p},0}, \iota_{\mathfrak{p},0})$ of local CM-type
  $(K_{\mathfrak{p}}/F_{\mathfrak{p}}, r_{\mathfrak{p}}/2)$ over $\Spf O_{\breve{E}_{\nu}}$.
  We may assume that $C_{X_{\mathfrak{p},0}} = O_{K_{\mathfrak{p}}}$. We endow $C_{X_{\mathfrak{p},0}}$  with
  the hermitian form (\ref{LS6e}) which corresponds to  the symmetric homomorphism
$\mathfrak{s}_{\mathfrak{p},0}: X_{\mathfrak{p},0}\rightarrow X_{\mathfrak{p},0}^{\wedge}$. 
  We define $X_0^{{\rm ba}} = \prod_{\mathfrak{p}, \; \text{banal}} X_{\mathfrak{p},0}$ and we endow it
  with $\mathfrak{s}^{{\rm ba}}_0 = \prod \mathfrak{s}_{\mathfrak{p},0}$. Then by Proposition \ref{LS3p} 
  we may replace $(2)$ in Definition \ref{balevel2d} by
  \begin{enumerate}\label{altCLlevel}
 \item[$(2^\prime)$]  {\it A class $\bar{\eta}_p$ of isomorphisms of $p$-adic \'etale sheaves,
  \begin{displaymath}
    \eta_p: \Lambda^{{\rm ba}} \rightarrow
    \underline{\Hom}_{O^{{\rm ba}}_K} (X_0^{{\rm ba}}, A[p^\infty]^{{\rm ba}}) \;
    \mathrm{mod}\; \mathbf{K}_{p}^{\star, {\rm ba}} ,
    \end{displaymath}
  which respect the forms on both sides up to a constant in
  $\mathbb{Z}_p^{\times}$.}
  \end{enumerate}
  The lisse $p$-adic sheaf on $\hat{\mathcal{A}}_{\mathbf{K}}$ given by the right hand
  side of $(2')$  is the algebraization
  of a lisse $p$-adic sheaf on $\mathcal{A}_{\mathbf{K}}$ which exists because
  this scheme is  proper over $\Spec O_{\breve{E}_{\nu}}$. We denote this sheaf
  by the same symbol. Then the scheme $\mathcal{A}^{\star}_{\mathbf{K}^{\star}}$
  \index[NO]{ABD@$\mathcal{A}^{\star}_{\mathbf{K}^{\star}}$} is
  given by the following functor on the category of schemes $S$ over
  $\Spec O_{\breve{E}_{\nu}}$: A point of $\mathcal{A}^{\star}_{\mathbf{K}^{\star}}(S)$
  consists of a point $(A, \iota, \bar{\lambda}, \bar{\eta}^p)$ of
  $\mathcal{A}_{\mathbf{K}}(S)$ and a class $\bar{\eta}_p$ as in $(2^\prime)$. We deduce the following description of  $\mathcal{A}^{\star}_{\mathbf{K}^{\star}}$. 
   
  \begin{proposition}\label{defA*}
The scheme
$\mathcal{A}^{\star}_{\mathbf{K}^{\star}}\times_{\Spec O_{\breve{E}_{\nu}}}\Spec \breve{E}_{\nu}$
represents the following functor on the category of $\breve{E}_{\nu}$-schemes. 
A $T$-valued point is a point $(A, \iota, \bar{\lambda}, \bar{\eta}^p)$ of 
$\mathcal{A}_{\mathbf{K}}(T)$ and a class $\bar{\eta}_p$ of
  isomorphisms of $p$-adic \'etale sheaves
  \begin{equation}\label{LS26e} 
    \eta_p: \Lambda^{{\rm ba}} \isoarrow
    \underline{\Hom}_{O^{{\rm ba}}_K} (T_p((X_0^{{\rm ba}})_{\breve{E}_{\nu}}),
    T_p(A)^{{\rm ba}}) \;
    \mathrm{mod}\; \mathbf{K}_{p}^{\star, {\rm ba}} ,
    \end{equation}
  which respect the forms on both sides up to a constant in
  $\mathbb{Z}_p^{\times}$.
  
  The scheme $\mathcal{A}^{\star}_{\mathbf{K}^{\star}}$ is the normalization of $\mathcal{A}_{\mathbf{K}}\times_{\Spec O_{E, ({\mathfrak p}_\nu)}}\Spec O_{\breve{E}_{\nu}}$ in $\mathcal{A}^{\star}_{\mathbf{K}^{\star}}\times_{\Spec O_{\breve{E}_{\nu}}}\Spec \breve{E}_{\nu}$ and is finite and \'etale over  $\mathcal{A}_{\mathbf{K}}\times_{\Spec O_{E, ({\mathfrak p}_\nu)}}\Spec O_{\breve{E}_{\nu}}$. \qed
  \end{proposition}

  \begin{theorem}\label{isovoerab}
Let $\breve{E}^{\rm ab}_{\nu}$ be the maximal abelian extension of $\breve{E}_{\nu}$.
Then there is an isomorphism
\begin{equation}\label{LS28e}
  \CA_{\bK^\star, E} \times_{\Spec E} \Spec \breve{E}^{\rm ab}_{\nu} \cong
  \mathcal{A}^{\star}_{\mathbf{K}^{\star}} \times_{\Spec O_{\breve{E}_{\nu}}} \Spec \breve{E}^{\rm ab}_{\nu}
\end{equation}
which is natural in $\mathbf{K}^{\star}$. 
  \end{theorem}
  \begin{proof}
    We make explicit what a level structure (\ref{LS26e}) means after base
    change to $\breve{E}^{\rm ab}_{\nu}$. Over $\breve{E}^{\rm ab}_{\nu}$ we may choose an
    isomorphism $\mathbb{Z}_p(1) \cong \mathbb{Z}_p$ and therefore we do not
    need to worry about Tate twists. The Tate module $T_p(X_{\mathfrak{p},0})$ of
    $X_{\mathfrak{p},0}$ over an algebraic closure  of $\breve E_\nu$ is an $O_{K_{\mathfrak{p}}}$-module which
    is free of rank $1$. Therefore the Galois group of $\breve E_\nu$  acts on the
    Tate-module via its maximal abelian quotient. We choose an
    isomorphism
    \begin{equation}\label{LS27e}
T_p(X_{\mathfrak{p},0}) \cong O_{K_{\mathfrak{p}}}
    \end{equation}
such that the action of the Galois group of $\breve E^{\rm ab}_\nu$  on both sides is
trivial. The symmetric map
$\mathfrak{s}_{\mathfrak{p}}: X_{\mathfrak{p},0}\rightarrow X_{\mathfrak{p},0}^{\wedge}$
induces a hermitian form $\varkappa_{\mathfrak{p},0}$ on the Tate-module
(\ref{LS27e}). We find
\begin{displaymath}
\varkappa_{\mathfrak{p},0}(x, y) = c_{\mathfrak{p},0} x \bar{y}, \quad x, y \in O_K  
  \end{displaymath}
for some constant $c_{\mathfrak{p},0} \in O_{F_{\mathfrak{p}}}^{\times}$. Note that in
the ramified case two isomorphism classes are possible for
$\varkappa_{\mathfrak{p},0}$. 

We consider a $T$-valued point
$(A, \iota, \bar{\lambda}, \bar{\eta}^p, \bar{\eta}_p)$ from the right hand
side of (\ref{LS28e}). Let $X = \prod X_{\mathfrak{p}}$ be the $p$-divisible
group of $A$. A polarization from $\bar{\lambda}$ induces an anti-hermitian
pairing $\varkappa_{\mathfrak{p}}$ on $T_p(X_{\mathfrak{p}})$. The anti-hermitian form
$\delta_{\mathfrak{p}}$ on
$\underline{\Hom}(T_p(X_{\mathfrak{p},0}), T_p(X_{\mathfrak{p}}))$ is given by
\begin{equation}\label{LS30e}
  \varkappa_{\mathfrak{p}}(u_1(x), u_2(y)) = \delta_{\mathfrak{p}}(u_1, u_2)
  c_{\mathfrak{p},0} x \bar{y}, \quad x, y \in O_{K_{\mathfrak{p}}},
\end{equation}
where $u_1, u_2 \in \underline{\Hom}(T_p(X_{\mathfrak{p},0}), T_p(X_{\mathfrak{p}}))$
are sections.

For an  $O_{K_{\mathfrak{p}}}$-lattice $(\Gamma, \varkappa_{\Gamma})$ with an
anti-hermitian form
$\varkappa_{\Gamma}: \Gamma \times \Gamma \rightarrow K_{\mathfrak{p}}$,  we write
$\Gamma[c] = (\Gamma, c \varkappa_{\Gamma})$. The equation (\ref{LS30e}) gives an isomorphism
\begin{displaymath}
  (\underline{\Hom}(T_p(X_{\mathfrak{p},0}), T_p(X_{\mathfrak{p}}),
  c_{\mathfrak{p},0} \delta_{\mathfrak{p}}) \cong
  (T_p(X_{\mathfrak{p}}), \varkappa_{\mathfrak{p}}).
\end{displaymath}
We see that a level structure (\ref{LS26e}) at the banal prime $\mathfrak{p}$
is given by an isomorphism
\begin{displaymath}
  \Lambda_{\mathfrak{p}}[c_{\mathfrak{p},0}] \rightarrow
  (T_p(X_{\mathfrak{p}}), \varkappa_{\mathfrak{p}}) \; \mathrm{mod}\;
  \mathbf{K}_{\mathfrak{p}}^{\star}. 
\end{displaymath}
Choosing a fixed isomorphism
$\Lambda_{\mathfrak{p}}[c_{\mathfrak{p},0}] \cong \Lambda_{\mathfrak{p}}$,  we see that
such a level structure at $\mathfrak{p}$ is the same as a class of isomorphisms 
\begin{displaymath}
\bar{\eta}_{\mathfrak{p}}: \Lambda_{\mathfrak{p}} \rightarrow
  (T_p(X_{\mathfrak{p}}), \varkappa_{\mathfrak{p}}) \; \mathrm{mod}\;
  \mathbf{K}_{\mathfrak{p}}^{\star}. 
\end{displaymath}
Since we want a level structure for $T_p(A)^{{\rm ba}}$, we require that the
$\eta_{\mathfrak{p}}$ must respect the bilinear forms on both sides by the
same factor $u \in \mathbb{Z}_p^{\times}$. For the special prime $\mathfrak{p}_v$
we take an arbitrary isomorphism 
\begin{displaymath}
\eta_{\mathfrak{p}_v}: \Lambda_{\mathfrak{p}_v} \rightarrow T_p(X_{\mathfrak{p}_v})  
\end{displaymath}
which respects the bilinear forms on both sides up to the same factor
$u \in \mathbb{Z}_p^{\times}$. This is possible because by (ii) of Definition
\ref{KneunA1d},  the $O_{K_{\mathfrak{p}_v}}$-lattices of both sides are
isomorphic and since,  by  Lemmas \ref{Kneun0l} and \ref{Kneun01l}, there exist 
isomorphisms with an arbitrary multiplicator $u \in \mathbb{Z}_p^{\times}$.
We set
\begin{displaymath}
  \tilde{\eta}_p = \eta_{\mathfrak{p}_v} \eta_p: \Lambda \otimes \mathbb{Z}_p
  \rightarrow T_p(A).  
  \end{displaymath}
Finally we set
\begin{displaymath}
\eta = \tilde{\eta}_p \eta^p: V \otimes \mathbb{A}_f \rightarrow  \hat{V}(A).  
\end{displaymath}
Let $\bar{\eta}$ be the class of this isomorphism modulo $\mathbf{K}^{\star}$. 
Then $(A, \iota, \bar{\lambda}, \bar{\eta})$ is a $T$-valued point of
$\CA_{\bK^\star, E}$. Since the last construction can be reversed, we obtain the
isomorphism of the theorem.
  \end{proof}
We will formulate a more precise version of the last Theorem. Let $E^c_{\nu}$
be the algebraic closure of $\breve{E}_{\nu}$. The action of the Galois group
$\Gal(E^c_{\nu}/\breve{E}_{\nu})$ on $T_p(X_{\mathfrak{p},0})$ is given by a character
\index[NO]{ZZTA@$\chi_{\mathfrak{p},0}$}
  \begin{equation}\label{LS37e}
    \chi_{\mathfrak{p},0}: \Gal(E^c_{\nu}/\breve{E}_{\nu}) \rightarrow
    O_{K_{\mathfrak{p}}}^{\times},   
  \end{equation}
such that $\sigma (t) = \chi_{\mathfrak{p},0}(\sigma) t$ for
$t \in T_p(X_{\mathfrak{p},0})$ and $\sigma \in \Gal(E^c_{\nu}/\breve{E}_{\nu})$.  
Since the polarization of $X_{\mathfrak{p}, 0}$ is defined over $\breve{E}_{\nu}$,
we obtain that $\Nm_{K_{\mathfrak{p}}/F_{\mathfrak{p}}} \chi_{\mathfrak{p},0}(\sigma) = 1$. 
We define
\begin{displaymath}
  \chi^{\mathrm{ba}}_{0}(\sigma) =
  \prod_{\mathfrak{p}, \text{banal}} \chi_{\mathfrak{p},0}(\sigma)
  \in G^{\mathrm{ba}}(\mathbb{Q}_p).  
\end{displaymath}
Finally we define
$\chi_0: \Gal(E^c_{\nu}/\breve{E}_{\nu}) \rightarrow G(\mathbb{Q}_p)$ by setting
\begin{displaymath}
  \chi_0(\sigma) = 1 \times \chi^{\mathrm{ba}}_{0}(\sigma) \in
  G_{\mathfrak{p}_v} \times G^{\mathrm{ba}}(\mathbb{Q}_p). 
\end{displaymath}
We note that this element is in the center of the group $G(\mathbb{Q}_p)$. 
By definition of the functor $\mathcal{A}_{ \mathbf{K}^{\star}, E}$ before Remark
\ref{gf1r},  $\chi_0(\sigma)$ acts on $\mathcal{A}_{ \mathbf{K}^{\star}, \breve E_{\nu}}=\mathcal{A}_{ \mathbf{K}^{\star},  E}\times_{\Spec E}\Spec \breve E_\nu$ via the
datum $(3)$, i.e. it acts by Hecke operators. We obtain the homomorphism
\index[NO]{ZZTB@$\chi^{\rm h}_0$}
\begin{equation}\label{LS33e} 
  \chi^{\rm h}_0: \Gal(E^c_{\nu}/\breve{E}_{\nu}) \rightarrow
  \Aut^{\rm opp} \mathcal{A}_{ \mathbf{K}^{\star}, \breve E_{\nu}}.
\end{equation}
(We write here the opposite group because the Hecke operators act by definition
from the right.) 
\begin{corollary}\label{isovoerab1c}
Let $\sigma \in \Gal(E^c_{\nu}/\breve{E}_{\nu})$. Then the action of
$\id_{\mathcal{A}^{\star}_{\mathbf{K}^{\star}}} \times \Spec \sigma$ on the right hand side of
(\ref{LS28e}) induces on the left hand side the automorphism
$\chi_0^{\rm h}(\sigma) \times \Spec \sigma$. 
\end{corollary}
\begin{remark}\label{gf2r}
In general, let $X$ a quasi-projective scheme over $\breve{E}_{\nu}$. Let
$\chi: \Gal(E^c_{\nu}/\breve{E}_{\nu}) \rightarrow  \Aut^{\rm opp} X$ be a continuous
homomorphism. Then descent says that there is a unique quasi-projective scheme
$X(\chi)$ over $\breve{E}_{\nu}$ and an isomorphism
\begin{displaymath}
  X \times_{\Spec \breve{E}_{\nu}} \Spec E^c_{\nu} \rightarrow
  X(\chi) \times_{\Spec \breve{E}_{\nu}} \Spec E^c_{\nu}
  \end{displaymath}
such that, for all $\sigma\in \Gal(E^c_{\nu}/\breve{E}_{\nu})$, the action of $\id_{X(\chi)} \times \Spec \sigma$ on the right hand
side induces on the left hand side the action $\chi(\sigma) \times \Spec \sigma$.  
We will call $X(\chi)$ the \emph{Galois twist} of $X$ by $\chi$. \index{Galois twist of a scheme}
  \end{remark}
\begin{proof}(of Corollary \ref{isovoerab1c}) 
  We take (\ref{LS28e}) over the algebraic closure $E^c_{\nu}$. For
  $\sigma \in \Gal(E^c_{\nu}/\breve{E}_{\nu})$, we write
  $\hat{\sigma} := \Spec \sigma$. We consider the non-commutative diagram
  \begin{equation}\label{LS36e}
  \begin{aligned}
\xymatrix{
\mathcal{A}_{\mathbf{K}^{\star},\breve{E}_{\nu}}\times_{\Spec \breve{E}_{\nu}}\Spec E^c_{\nu}
\ar[r] \ar[d]_{\id_{\mathcal{A}} \times \hat{\sigma}} &
\mathcal{A}^{\star}_{\mathbf{K}^{\star},\breve{E}_{\nu}}\times_{\Spec\breve{E}_{\nu}}\Spec E^c_{\nu}
\ar[d]^{\id_{\mathcal{A}^{\star}} \times \hat{\sigma}}\\
\mathcal{A}_{\mathbf{K}^{\star},\breve{E}_{\nu}}\times_{\Spec \breve{E}_{\nu}}\Spec E^c_{\nu}
\ar[r]  &
\mathcal{A}^{\star}_{\mathbf{K}^{\star},\breve{E}_{\nu}}\times_{\Spec\breve{E}_{\nu}}\Spec E^c_{\nu} . 
  }
  \end{aligned}
  \end{equation}
  To understand how this does not commute we consider more generally a scheme
  $S$ of finite type over $\breve E_{\nu}$ and write
  $S_{E^c_{\nu}} = S \times_{\Spec \breve E_{\nu}} \Spec E^c_{\nu}$. The morphism
  $\hat{\sigma}_S:=\id_S \times \hat{\sigma}: S_{E^c_{\nu}}\rightarrow S_{E^c_{\nu}}$
  induces maps
  \begin{displaymath}
    \sigma_{\mathcal{A}}: \mathcal{A}_{\mathbf{K}^{\star},\breve{E}_{\nu}}(S_{E^c_{\nu}})
    \rightarrow \mathcal{A}_{\mathbf{K}^{\star}, \breve{E}_{\nu}}(S_{E^c_{\nu}}), \quad
    \sigma_{\mathcal{A}^{\star}}:
    \mathcal{A}^{\star}_{\mathbf{K}^{\star},\breve{E}_{\nu}}(S_{E^c_{\nu}})
    \rightarrow \mathcal{A}^{\star}_{\mathbf{K}^{\star}, \breve{E}_{\nu}}(S_{E^c_{\nu}}). 
  \end{displaymath}
  Our task is to compare the effect of these maps on an element
  $\xi \in \mathcal{A}_{\mathbf{K}^{\star},\breve{E}_{\nu}}(S_{E^c_{\nu}}) = \mathcal{A}^{\star}_{\mathbf{K}^{\star},\breve{E}_{\nu}}(S_{E^c_{\nu}})$. 
The moduli interpretation describes $\xi$ as a point of
$\mathcal{A}_{\mathbf{K}^{\star}, \breve{E}_{\nu}}$ by a point
$(A,\iota,\lambda,\bar{\eta}^p)\in\mathcal{A}_{\mathbf{K},\breve{E}_{\nu}}(S_{E^c_{\nu}})$
and a rigidification
$\bar{\eta}_p: \Lambda^{\mathrm{ba}} \rightarrow T_p(A)^{\mathrm{ba}}  \; \text{mod}  \;  \mathbf{K}^{\star, \mathrm{ba}}$.
To make this more precise, we choose a geometric point
$\omega: \Spec E^c_{\nu} \rightarrow S$ which extends naturally to a point
$\omega: \Spec E^c_{\nu} \rightarrow S_{E^c}$. We define $\omega'$ by the
commutative diagram
\begin{displaymath}
  \xymatrix{
S_{E^c_{\nu}} \ar[r]^{\hat{\sigma}_S}  &  S_{E^c_{\nu}} \\
\Spec E^c_{\nu} \ar[u]_{\omega'} \ar[r]^{\hat{\sigma}}&
\Spec E^c_{\nu} \ar[u]_{\omega} .\\
}
  \end{displaymath}
The rigidification is  given by a homomorphism 
\begin{displaymath}
\eta_p: \Lambda^{\mathrm{ba}} \rightarrow T_p(A_{\omega})^{\mathrm{ba}}.
\end{displaymath}
There is an isomorphism
\begin{displaymath}
T_p((\hat{\sigma}_S^{\ast} A))_{\omega'} = \hat{\sigma}^{\ast} (T_p(A_{\omega})). 
  \end{displaymath}
By the moduli interpretation, the point $\sigma_{\mathcal{A}}(\xi)$ is given by
$(\hat{\sigma}_S^{\ast}A,  \hat{\sigma}_S^{\ast}\iota, \hat{\sigma}_S^{\ast}\lambda,   \hat{\sigma}_S^{\ast}\bar{\eta}^p)$
and the rigidification is given by 
\begin{equation}\label{LS35e}
  \Lambda^{\mathrm{ba}} = \hat{\sigma}^{\ast} (\Lambda^{\mathrm{ba}}) \;
  \overset{\hat{\sigma}^{\ast}(\eta_p)}{\longrightarrow} \; 
  \hat{\sigma}^{\ast} (T_p(A_{\omega})).  
\end{equation}
Now we consider $\sigma_{\mathcal{A}^{\star}}(\xi)$. We can give the sheaf
$T_p((X_0^{{\rm ba}})_{\breve{E}_{\nu}})$ in (\ref{LS26e}) equivalently by
the $\Gal(E^c_{\nu}/E_{\nu})$-module $\Lambda^{\mathrm{ba}}(\chi^{\mathrm{ba}}_0)$, 
where we indicate that the Galois group acts via the character
$\chi^{\mathrm{ba}}_0$. Then  $\bar{\eta}_p$ of (\ref{LS26e}) can be considered
as a class of maps
\begin{equation}\label{LS34e} 
  \Lambda^{\mathrm{ba}}(\chi^{\mathrm{ba}}_0)_{E^c_{\nu}} \rightarrow
  T_p(A_{\omega})^{\mathrm{ba}} \; \text{mod}  \;  \mathbf{K}^{\star, \mathrm{ba}}.
  \end{equation}
Since we are over $E^c_{\nu}$, the action via $\chi^{\mathrm{ba}}_0$ is trivial
and therefore $(A,\iota,\lambda,\bar{\eta}^p, \eta_p)$ describes also a point
of $\mathcal{A}^{\star}_{\mathbf{K}^{\star}, \breve{E}_{\nu}}(S_{E^c_{\nu}})$. But if we want
to identify the inverse image of this point by $\hat{\sigma}_S$ we must take
into account the twist $\chi^{\mathrm{ba}}_0$. This inverse image is again given
by $(\hat{\sigma}_S^{\star}A,  \hat{\sigma}_S^{\star}\iota, \hat{\sigma}_S^{\star}\lambda,   \hat{\sigma}_S^{\ast}\bar{\eta}^p)$
as before, but the new rigidification at the banal places is
\begin{equation}\label{LS35f}
  \Lambda^{\mathrm{ba}}(\chi^{\mathrm{ba}}_0)_{E^c_{\nu}} \cong 
  \hat{\sigma}^{\ast}( \Lambda^{\mathrm{ba}})(\chi^{\mathrm{ba}}_0)_{E^c_{\nu}} \;
  \overset{\hat{\sigma}^{\ast}(\eta_p)}{\longrightarrow} \; 
  \hat{\sigma}^{\ast} (T_p(A_{\omega})).  
\end{equation}
The first isomorphism comes from the fact that both sides are the inverse image
of $\Lambda^{\mathrm{ba}}(\chi^{\mathrm{ba}}_0)$ considered as a sheaf on
$\Spec \breve E_{\nu}$. Therefore this isomorphism is the descent datum on the constant
sheaf, which  is the multiplication by $\chi^{\mathrm{ba}}_0(\sigma)$. We obtain
\begin{displaymath}
  \Lambda^{\mathrm{ba}}(\chi^{\mathrm{ba}}_0) \;
  \overset{\chi^{\mathrm{ba}}_0(\sigma)}{\longrightarrow} \; 
  \hat{\sigma}^{\ast}( \Lambda^{\mathrm{ba}})(\chi^{\mathrm{ba}}_0) \;
  \overset{\hat{\sigma}^{\ast}(\eta_p)}{\longrightarrow} \; 
  \hat{\sigma}^{\ast} (T_p(A_{\omega})) .
\end{displaymath}
This proves that
$\sigma_{\mathcal{A}^{\ast}} = \chi^{\mathrm{h}}_0(\sigma) \sigma_{\mathcal{A}}$. If
we apply this to the diagram (\ref{LS36e}), we obtain
\begin{displaymath}
  \id_{\mathcal{A}^{\ast}} \times \hat{\sigma} =
  \chi^{\mathrm{h}}_0(\sigma) (\id_{\mathcal{A}} \times \hat{\sigma}). 
  \end{displaymath}
\end{proof}
  
 We now drop the assumption on $\bK_p^\star$ that it be contained in $ \bK_p$ and come from a product of $\bK^\star_{\mathfrak p}$. More precisely, let $\bK_p^\star\subset G(\BQ_p)$ be of the form 
 \begin{equation}\label{intKP}
 \bK_p^\star=G(\BQ_p)\cap \bK_{\mathfrak{p}_v}\bK_p^{\star, {\rm ba}} ,  
 \end{equation}
 where  $\bK_p^{\star, {\rm ba}}$ is an arbitrary open compact subgroup of $G^{\rm ba}(\BQ_p)$.  
 Since $\bK_{\mathfrak{p}_v}$ is a normal subgroup of $G_{\mathfrak{p}_v}$ and $G(\BQ_p)$, this class of subgroups is stable under conjugation by elements of $G(\BQ_p)$. Therefore, using the naturality of the construction in Proposition \ref{defA*}, we can extend the definition of $\mathcal{A}^{\star}_{\mathbf{K}^{\star}}$ to all such $\bK^\star=\bK^\star_p\bK^p$ by first passing to a small enough normal subgroup of finite index and then dividing out by the factor group.  
 
 We make this extension process more explicit by defining the functor
 $\hat{\mathcal{A}}^{\star}_{\mathbf{K}^{\star}}$ without using the choice of $\Lambda_p$.
 Thereby the action of the Hecke operators becomes more obvious.

  Let $\theta: X \rightarrow Y$ be an isogeny of $p$-divisible
 $O_F \otimes \mathbb{Z}_p$-modules. Let $\mathfrak{p}_{v}$ be a prime of
  $O_F$  over $p$. We say that $\theta$ is an isogeny of order prime to
  $\mathfrak{p}_{v}$ if $\theta_{\mathfrak{p}_{v}}$ is an isomorphism. We use a similar terminology for abelian varieties with action by $O_F$.
  \index[NO]{ABC@$\hat{\mathcal{A}}^{\star}_{\mathbf{K}^{\star}}$}

 We consider a scheme $S$ over $\Spf O_{\breve{E}_{\nu}}$. We consider abelian
 schemes $\bar{A}$ over $S$ up to isogeny of order prime to
 $\mathfrak{p}_{v}$ which are endowed with an action
 $\iota: O_K \rightarrow \End \bar{A}$ and with a $\mathbb{Q}$-homogeneous
 polarization $\bar{\lambda}$ such that the Rosati involution induces the
 conjugation on $O_K$. Moreover, we assume that there is a triple
 $(A, \iota, \lambda)$ as in Definition \ref{KneunA1d} which represents
 $(\bar{A}, \iota, \bar{\lambda})$ such that $(A, \iota)$ satisfies the
 conditions $({\rm KC}_r)$ and  $({\rm EC}_r)$ and such that
 $(A, \iota, \lambda)$ satisfies the conditions (i) and (ii) of
 Definition \ref{KneunA1d}. Then call $(\bar{A}, \iota, \bar{\lambda})$
 an \emph{admissible prime-to-$\mathfrak{p}_{v}$-isogeny class}.
 Let $X = \prod_{\mathfrak{p}} X_{\mathfrak{p}}$ be the $p$-divisible group of $A$.
 Then $(C_{X^{\mathrm{ba}}}, \xi_{\lambda})$ makes sense and
 $(C_{X^{\mathrm{ba}}} \otimes \mathbb{Q}, \xi_{\lambda})$ depends only
 on $(\bar{A}, \iota, \bar{\lambda})$.
 Let $(\mathbb{X}, \iota_{\mathbb{X}}, \lambda_{\mathbb{X}})$ as in Definition
 \ref{balevel2d}. 
 \begin{definition}\label{balevel3d}
  We define a functor $\hat{\mathcal{A}}^{\star}_{\mathbf{K}^{\star}}$ on the
  category of schemes $S$ over $\Spf O_{\breve{E}_{\nu}}$.
  A point of $\hat{\mathcal{A}}^{\star}_{\mathbf{K}^{\star}}(S)$ consists of
  the following data:
\begin{enumerate}
\item[(1)] an admissible prime-to-$\mathfrak{p}_v$-isogeny class
  $(\bar{A}, \iota, \bar{\lambda})$ over $S$. 
\item[(2)]   a class of isomorphisms  
  \begin{displaymath}
    \bar{\eta}_p: C_{\mathbb{X}^{{\rm ba}}} \otimes \mathbb{Q} \rightarrow
    C_{\bar{A}[p^\infty]^{{\rm ba}}} \otimes \mathbb{Q} \;
    \mathrm{mod}\; \mathbf{K}_{p}^{\star, {\rm ba}},
  \end{displaymath}
  which respects the bilinear forms on both sides up to a factor in
  $\mathbb{Q}_p^{\times}$. 
\end{enumerate}
   \end{definition}
 We explain in more detail what is meant by $(2)$. We assume that $S$ is
 connected and we choose a geometric point $\omega$ of $S$. Then the meaning of $(2)$ is 
 that we have a class of isomorphisms
 \begin{displaymath}
  \eta_p: (C_{\mathbb{X}^{{\rm ba}}})_{\omega} \otimes \mathbb{Q} \rightarrow
   (C_{\bar{A}[p^\infty]^{{\rm ba}}})_{\omega} \otimes \mathbb{Q} \; \mathrm{mod} \;
   \mathbf{K}_{p}^{\star, {\rm ba}} 
 \end{displaymath}
 which respects the bilinear
forms on both sides up to a factor in $\mathbb{Q}_p^{\times}$ and such that the
class is preserved by the action of $\pi_1(S, \omega)$.

Let $\mathbf{K}_{p}^{\star}$ as in Definition \ref{balevel2d}. Then the functors
of the Definitions \ref{balevel2d} and \ref{balevel3d} coincide. Indeed, let
us start with a point of Definition \ref{balevel3d}. Also, fix a triple $(A, \iota, \lambda)$ which represents $(\bar A, \iota, \bar \lambda)$, as before Definition \ref{balevel3d}. The sublattice 
$\Lambda_p^{\mathrm{ba}} \subset C_{\mathbb{X}^{{\rm ba}}})_{\omega}$ is fixed by
$\mathbf{K}_{p}^{\star, {\rm ba}}$. Therefore the image $C$ of
$\Lambda_p^{\mathrm{ba}}$ by $\eta_p$ depends only on the class $\bar{\eta}_p$
and is invariant by $\pi_1(S, \omega)$. Therefore $C$ defines a $p$-adic
\'etale sheaf on $S$. We endow it with the polarization induced by
$\varsigma^{\mathrm{ba}}$, cf. (\ref{LS18e}). Therefore, using the contracting
functor $C$ defines a
$p$-divisible $O_{K^{\mathrm{ba}}}$-module $Y^{\mathrm{ba}}$ with a polarization.  
Then 
$Y := X_{\mathfrak{p}_v} \times Y^{\mathrm{ba}}$ is isogenous to the $p$-divisible
group $X$ of $A$. The polarization on $Y$ differs from the polarization
induced from $\lambda$ on $A$ by a factor in $\mathbb{Z}_p^{\times}$, as we see
by comparing the degrees of the polarizations. Therefore we obtain a
point $(A_1, \iota_1, \bar{\lambda}_1)$ of
$\hat{\mathcal{A}}^{\star}_{\mathbf{K}^{\star}}(S)$ which is isogenous to
$(\bar{A}, \iota,\bar\lambda)$. This proves that the point we started with comes
from a point of the functor in Definition \ref{balevel2d}. It is clear that
we have a bijection. 

We have an action of $G^{\mathrm{ba}}(\mathbb{Q}_p)$ on the tower
$\hat{\mathcal{A}}^{\star}_{\mathbf{K}^{\star}}$ for varying $\mathbf{K}_{p}^{\star}$.
This action extends to the algebraization $\mathcal{A}^{\star}_{\mathbf{K}^{\star}}$ 
and coincides via Theorem \ref{isovoerab} with the Hecke operators on the
tower $\mathcal{A}_{\mathbf{K}^{\star}, E}$. 

\begin{corollary}\label{genKPschemes}
  For every $\bK^\star=\bK^\star_p\bK^p$ with \eqref{intKP} , there 
  exists a normal scheme $\mathcal{A}^{\star}_{\mathbf{K}^{\star}}$ over $\Spec O_{\breve E_\nu}$ 
  such that for the $p$-adic completion of this scheme there is an isomorphism 
  \begin{displaymath}
\hat{\mathcal{A}}^{\star}_{\mathbf{K}^{\star}} \simeq
J(\mathbb{Q}) \backslash [(\hat{\Omega}_{F_{v}} \times_{\Spf O_{F_{v}}}
\Spf O_{\breve{E}_{\nu}})\times  G^{{\rm ba}}(\mathbb{Q}_p)/\mathbf{K}_{p}^{\star, {\rm ba}}
\times G(\mathbb{A}_f^p)/\mathbf{K}^p].
  \end{displaymath}
 For varying $\mathbf{K}^{\star}$, these schemes form a tower with an action of the group
$G(\mathbb{Q}_p) \times G(\mathbb{A}_f^{p})$, where the action of $G(\BQ_p)$ factors through  $G(\mathbb{Q}_p) \rightarrow G^{{\rm ba}}(\mathbb{Q}_p)$. The isomorphism of formal schemes 
is compatible with these actions. 
  
The general fiber of $\mathcal{A}^{\star}_{\mathbf{K}^{\star}}$ is a Galois twist of
 $\mathcal{A}_{ \mathbf{K}^{\star}, E}\times_{\Spec E}\Spec \breve E_\nu$ by the character
$\chi_0^{\rm h}$, cf. (\ref{LS33e}) and Remark \ref{gf2r}. The Galois twist
respects the Hecke operators (cf. section \ref{ss:TTm} for an explicit
description of $\chi_0^{\rm h}$).  
  \end{corollary}
\begin{proof}
  This is a consequence of Proposition \ref{LS1p} and the general pattern of
  $p$-adic uniformization, cf. (\ref{KneunU5e}). The last assertion follows
  because $\chi_0: \Gal(E^{\rm ab}_{\nu}/E_{\nu}) \rightarrow G(\mathbb{Q}_p)$
  factors through the center. 
  \end{proof}

\subsection{The rigid-analytic uniformization}
Let $\CA_\bK^{\rm rig}$ denote the rigid-analytic space over $\Sp\, E_\nu$ associated to $\CA_{\bK, E}$. Then Theorem \ref{maint} implies the following corollary concerning generic fibers.
\begin{corollary}\label{rigmain7a} Let $\bK=\bK_p\bK^p$ as in \eqref{boldK1e}. 
There exists an isomorphism of rigid-analytic spaces over $\Sp\,  {\breve E_\nu}$,
\begin{equation*}\label{rigmainequa} 
 \CA_\bK^{\rm rig}\times_{{\rm Sp}\,{E_{\nu}}}{\rm Sp}\, {\breve E_\nu} \simeq J(\BQ)\bs\big[\big({\Omega}_{F_v} \times_{{\rm Sp}\,{F_v}}{\rm Sp}\,
    {\breve E_\nu}\big) \times \hat{G}'(\mathbb{Q}_p)/\hat{G}'(\mathbb{Z}_p)
    \times G(\BA_f^p)/{\bK}^p\big]  \, .
  \end{equation*}
For varying $\bK^p$, this isomorphism is compatible with the action of
$G(\BA_f^p)$ through Hecke correspondences on both sides.\qed

{\rm Here ${\Omega}_{F_v}=\BP^1_{F_v}\setminus \BP^1(F_v)$ is Drinfeld's $p$-adic halfspace corresponding to the $p$-adic field $F_v$.  }
\end{corollary}
Similarly, Corollary \ref{genKPschemes} implies the following corollary concerning  generic fibers for deeper level structures.
\begin{corollary}\label{rigmain7} Assume that there are banal primes. Let $\bK^\star=\bK^\star_p\bK^p$ with \eqref{intKP}. Let $\CA_{\bK^\star}^{\rm rig}$ denote the rigid-analytic space over $\Sp\, E_\nu$ associated to $\CA_{\bK^\star, E}$.
There exists an isomorphism of rigid-analytic spaces over $\Sp\,  {\breve E^{\rm ab}_\nu}$, 
    \begin{displaymath}
{\mathcal{A}}^{\rm rig}_{\mathbf{K}^{\ast}}\times_{\Sp\, E_\nu}\Sp\, \breve E^{\rm ab}_\nu \simeq
J(\mathbb{Q}) \backslash [({\Omega}_{F_{v}} \times_{\Sp\, {F_{v}}}
\Sp\, {\breve E^{\rm ab}_\nu})\times  G^{{\rm ba}}(\mathbb{Q}_p)/\mathbf{K}_{p}^{\ast, {\rm ba}}
\times G(\mathbb{A}_f^p)/\mathbf{K}^p].
  \end{displaymath}
 For variable $\bK^\star$, this isomorphism is
  compatible with the Hecke correspondences by $G(\BQ_p)\times G(\BA_f^p)$. \qed
\end{corollary}

\subsection{Determination of the character $\chi_0^{\mathrm{h}}$}\label{ss:TTm}  
In this section we give an explicit description of the character
$\chi_0^{\mathrm{h}}$ (\ref{LS33e}) which is used in Corollary
\ref{genKPschemes}. 
In the case where $\mathfrak{p}_{v}$ is ramified in $K/F$, we only obtain
the restriction of $\chi_0^{\mathrm{h}}$ to the Galois group of a
quadratic extension of $\breve{E}_{\nu}$. It is enough to describe
$\chi_{\mathfrak{p},0}$ (\ref{LS37e}) for each banal prime $\mathfrak{p}$.
This is done by Proposition \ref{KM2p} below.

Let $K/F$ be a CM-field. Let $\Xi \subset \Hom_{\text{$\mathbb{Q}$-Alg}}(K, \mathbb{C})$\index[NO]{ZZNA@$\Xi$}
be a CM-type. We denote the reflex field by $H$. We define an algebraic
torus over $\mathbb{Q}$, with $\BQ$-valued points 
\begin{displaymath}
T(\mathbb{Q}) = \{a \in K^{\times} \; | \; a \bar{a} \in \mathbb{Q}^{\times} \}. 
\end{displaymath}
We use the notation $V = K$ for $K$ regarded as a $K$-vector space. 

We recall the reciprocity law. We define the homomorphism
\begin{displaymath}
\mu: \mathbb{C}^{\times} \rightarrow (K \otimes_{\mathbb{Q}} \mathbb{C})^{\times} \cong
  \prod_{\varphi: K \rightarrow \mathbb{C}} \mathbb{C}^{\times}. 
\end{displaymath}
The element $\mu(z)$, for $z \in \mathbb{C}$, has component $z$ for
$\varphi \in \Xi$ and has component $1$ for $\varphi \notin \Xi$ on the right hand
side. We find
$\mu \bar{\mu} = 1 \otimes z \in (K \otimes_{\mathbb{Q}} \mathbb{C})^{\times}$.
We obtain a homomorphism of algebraic tori
\begin{displaymath}
\mu: \mathbb{G}_{m, \mathbb{C}} \rightarrow T_{\mathbb{C}}.  
  \end{displaymath}
This homomorphism is defined over $H$, 
\begin{displaymath}
\mu: \mathbb{G}_{m,H} \rightarrow T_H. 
\end{displaymath}
From this we deduce the reciprocity map\index[NO]{RAE@$\mathfrak{r}$}
\begin{equation}\label{KM1e}
  \mathfrak{r}: \mathrm{Res}_{H/\mathbb{Q}}( \mathbb{G}_{m,H}) 
  \overset{\mu}{\longrightarrow} \mathrm{Res}_{H/\mathbb{Q}} (T_H )
  \overset{\Nm_{H/\mathbb{Q}}}{\longrightarrow} T.  
  \end{equation}
We consider over the algebraic closure $\bar{H} = \bar{\mathbb{Q}}$ the set of
tuples $(A, \iota, \bar{\lambda}, \kappa)$, where $(A, \iota)$ is an abelian
variety over $\bar{H}$ of CM-type $\Xi$, endowed with a
$\mathbb{Q}$-homogeneous polarization $\bar{\lambda}$ which induces on $K$
the conjugation over $F$ and an isomorphism
$\kappa: \hat{V}(A) \rightarrow V \otimes \mathbb{A}_f$ of
$K \otimes \mathbb{A}_f$-modules. We call a second tuple
$(A', \iota', \bar{\lambda}', \kappa')$ equivalent to
$(A, \iota, \bar{\lambda}, \kappa)$ if there is  a quasi-isogeny
\begin{equation}\label{KM3e}
    \alpha: (A, \iota, \bar{\lambda}) \rightarrow
    (A', \iota', \bar{\lambda}')
    \end{equation}
such that the following diagram 
\begin{displaymath}
\xymatrix{
     \hat{V}(A) \ar[rd]_{\kappa} \ar[rr]^{\alpha} & &
     \hat{V}(A') \ar[ld]^{\kappa'}\\
     & V \otimes \mathbb{A}_f & \\
} 
  \end{displaymath}
commutes. We also say that $(A, \iota, \bar{\lambda}, \kappa)$ is quasi-isogenous to $(A', \iota', \bar{\lambda'}, \kappa')$.

Let $\mathcal{C}_{\Xi}$ be the set of tuples $(A, \iota,\bar{\lambda},\kappa)$
up to equivalence.
Let $\sigma \in \Gal(\bar{H}/H)$. Taking the inverse image of
$(A, \iota, \bar{\lambda}, \kappa)$ by
$\hat{\sigma} := \Spec \sigma: \Spec \bar{H} \rightarrow  \Spec \bar{H}$ 
gives a left action of $\Gal(\bar{H}/H)$ on $\mathcal{C}_{\Xi}$. We denote
the inverse image by $\sigma (A, \iota, \bar{\lambda}, \kappa)$. 

We formulate the main theorem of complex multiplication of Shimura and
Taniyama. 
\begin{theorem}\label{ShimTan}(\cite[Thm. 4.19]{De}) 
  The Galois group $\Gal(\bar{H}/H)$ acts on $\mathcal{C}_{\Xi}$ via its maximal
  abelian quotient $\Gal({H}^{\rm ab}/H)$. Let
  $e \in (H \otimes \mathbb{A})^{\times}$ and let
  $\mathrm{rec}(e) \in \Gal({H}^{\rm ab}/H)$ be the automorphism given by the
  reciprocity law of class field theory. The following tuples are equivalent:
  \begin{displaymath}
    \mathrm{rec}(e) (A, \iota, \bar{\lambda}, \kappa) \equiv 
    (A, \iota, \bar{\lambda}, \mathfrak{r}(e_f) \kappa)  ,
    \end{displaymath}
    where $e_f$ is the finite part of the id\`ele $e$.  
\end{theorem}
\begin{remark}
Let $(H^{\times})^{\wedge} \subset (H \otimes \mathbb{A}_f)^{\times}$ be the closure
of $H^{\times}$. We deduce a homomorphism 
\begin{equation}\label{KM5e}
  \Gal(\bar{H}/H) \rightarrow (H \otimes \mathbb{A}_f)^{\times}/(H^{\times})^{\wedge}
  \overset{\mathfrak{r}}{\longrightarrow}
  T(\mathbb{A}_f)/T(\mathbb{Q}) ,
\end{equation}
where the first arrow is deduced from class field reciprocity and the second
arrow exists because $T(\mathbb{Q}) = T(\mathbb{Q})^{\wedge}$. To see this last fact, we note that 
the group of units in $T(\BQ)$  is finite. Indeed, the units are elements of
$K^{\times}$ with all absolute values equal to $1$ at all places including the
infinite ones. Therefore
$T(\mathbb{Q}) = T(\mathbb{Q})^{\wedge}$ by Chevalley's theorem.

 Theorem \ref{ShimTan} says
that the action of $\Gal(\bar{H}/H)$ on $\mathcal{C}_{\Xi}$ is via (\ref{KM5e}).
One can consider the Shimura variety $\mathrm{Sh}_{T}$.
We may choose as usual a bijection
 \begin{displaymath}
  \mathrm{Sh}_{T}(\bar{H}) = \mathrm{Sh}_{T}(\mathbb{C}) =
  T(\mathbb{A}_f)/T(\mathbb{Q})^{\wedge}. 
 \end{displaymath}
 Then the theorem may be regarded as a consequence of Langlands' description of
 the reduction of this Shimura variety at good places \cite{K}. 
  \end{remark}
We  fix an embedding $\bar\BQ\to\bar\BQ_p$. 
The $p$-adic place which is induced on a subfield of $\bar{\mathbb{Q}}$ will
be denoted by $\nu$. 
\begin{proposition}\label{KM1p}
  Let $L\subset \bar\BQ$ be a number field such that $H \subset L$. Let
  $(A_0, \iota_0, \lambda_0)$ be an abelian variety over $L$ with an action
  $\iota_0: O_K \rightarrow \End A_0$ which is of CM-type $\Xi$. We assume
  that $A_0$ has good reduction at
  $\nu$. The group $\Gal(\bar{L}_{\nu}/L_{\nu})$ acts on the Tate module
  $T_p(A_0)$ via its maximal abelian
  quotient $\Gal(L^{\rm ab}_{\nu}/L_{\nu})$. Let
  $I_{\nu} \subset \Gal(L^{\rm ab}_{\nu}/L_{\nu})$ be the inertia group.
  The action of $I_{\nu}$ on the Tate module can be described as follows.

  The inverse of the map (\ref{KM1e}) induces a homomorphism
  \begin{displaymath}
   \rho: L_{\nu}^{\times} \overset{\Nm_{L_{\nu}/H_{\nu}}}{\longrightarrow} H_{\nu}^{\times}
      \subset (H \otimes \mathbb{Q}_p)^{\times}
      \overset{\mathfrak{r}^{-1}}{\longrightarrow} (K \otimes \mathbb{Q}_p)^{\times}. 
  \end{displaymath}
 Composing $\rho$ with the reciprocity law of local class field theory yields 
  \begin{displaymath}
    I_{\nu} \overset{\mathrm{rec}}{\simeq} O_{L_{\nu}}^{\times}
    \overset{\rho}{\longrightarrow} (O_K \otimes {\mathbb{Z}_p})^{\times}. 
  \end{displaymath}
  The action of an element $\sigma \in I_{\nu}$ on the Tate-module is the
  multiplication by the image in the right hand side. 
\end{proposition}
\begin{proof}
  We set $A = A_0\otimes_H{\bar{H}}$ with the $O_K$-action and the induced
  polarization. We set $\Lambda = O_K \subset V$ and
  $\hat{\Lambda} = O_K \otimes \hat{\mathbb{Z}}$. We choose a rigidification
  $\kappa: \hat{T}(A) \overset{\sim}{\rightarrow} \hat{\Lambda}$.
  We consider the tuple $(A, \iota, \bar{\lambda}, \kappa)$. Let
  \begin{displaymath}
    \sigma \in \tilde{I}_{\nu} \subset \Gal(\bar{H}_{\nu}/H_{\nu}) \subset
    \Gal(\bar{H}/H) 
  \end{displaymath}
  be an element of the inertia group at $\nu$. The image in
  $\Gal({H}^{\rm ab}/H)$ corresponds to an id\`ele in
  $(H \otimes \mathbb{A})^{\times}$ which has components
  $1$ outside $\nu$ and a component $e_{\nu} \in O_{H_{\nu}}^{\times}$ at the place
$\nu$. We denote the id\`ele also by $e_{\nu}$. By  Theorem \ref{ShimTan}, we have a
  quasi-isogeny
    \begin{displaymath}
    \hat{\sigma}^{\ast} (A, \iota, \bar{\lambda}, \kappa) \cong
    (A, \iota, \bar{\lambda}, \mathfrak{r}(e_{\nu})\kappa). 
  \end{displaymath}
Let us moreover assume that $\sigma$ fixes the elements of $L$. 
Since $(A, \iota, \bar{\lambda})$ is defined over $L$, it is not changed by
$\hat{\sigma}^{\ast}$. Now we consider the product of the Tate modules for all primes, 
\begin{displaymath}
  \hat{T}(A) = \hat{\sigma}^{\ast} (\hat{T}(A))
  \overset{\hat{\sigma}^{\ast}( \kappa)}{\longrightarrow} \hat{\Lambda}. 
\end{displaymath}
The first identification is due to the fact that $\hat{T}(A)$ is a projective
limit of \'etale sheaves on $\Spec L$.
\begin{lemma}\label{KM1l}
  Denote by $\hat{T}(\sigma)$ the action of $\sigma$ on $\hat{T}(A)$.
  Then 
  \begin{displaymath}
\hat{\sigma}^{\ast}(\kappa) \hat{T}(\sigma) = \kappa. 
    \end{displaymath}
\end{lemma}
We postpone the proof of the Lemma.  Because we have good
reduction, the element $\hat{T}(\sigma)$ acts trivially on the
Tate modules $T_{\ell}(A)$ for $\ell \neq p$. On $T_p(A)$ it acts by multiplication with an element
$u(\sigma) \in (O_K \otimes \mathbb{Z}_p)^{\times}$. Therefore we have a quasi-isogeny  
\begin{equation}\label{KM2e}
  (A, \iota, \bar{\lambda}, \mathfrak{r}(e_{\nu})\kappa) \cong
  (A, \iota, \bar{\lambda}, u(\sigma)^{-1}\kappa)
\end{equation}
The quasi-isogeny giving this equivalence must be trivial on the Tate modules
$V_{\ell}(A)$ for $\ell \neq p$. It is therefore the identity. The proposition follows therefore from Lemma \ref{KM1l}. 
\end{proof}

\begin{proof}(of Lemma \ref{KM1l}) We consider an \'etale sheaf $G$ over
  $\Spec L$ where $L$ is any field. Let $L^s$ be the separable closure of $L$.
  For $\sigma \in \Gal(L^s/L)$ we denote by
  $G(\sigma): G(L^s) \rightarrow G(L^s)$ the natural action. 
  Let $\underline{\Gamma}$ be a constant sheaf on $\Spec L$ associated to a
  set $\Gamma$. Let
  \begin{displaymath}
\kappa: G \rightarrow \underline{\Gamma} 
    \end{displaymath}
  be an isomorphism of sheaves on $(\Spec L^s)_{\text{\'et}}$. There are canonical
  isomorphisms $\hat{\sigma}^{\ast}( G )\cong G$ and
  $\hat{\sigma}^{\ast} (\Gamma) \cong \Gamma$ because both sheaves are defined
  over $L$. We must show that the map 
  \begin{displaymath}
  \kappa': G(L^{s}) \cong \hat{\sigma}^{\ast} (G)(L^{s})
    \overset{\hat{\sigma}^{\ast} (\kappa)}{\longrightarrow}
    \hat{\sigma}^{\ast} (\underline{\Gamma})(L^{s}) = \Gamma 
  \end{displaymath}
  coincides with $\kappa G(\sigma^{-1})$.
  
  Let $A$ be a finite \'etale algebra over $L^s$. By definition of
  the inverse image, we have 
  $\hat{\sigma}^{\ast}(G)(A) = G(A_{[\sigma]})$. Therefore the $L^s$-algebra
  isomorphism $\sigma: L^s \rightarrow L^s_{[\sigma]}$ induces a natural map
  $G[\sigma]\colon G(L^s) \rightarrow \hat{\sigma}^{\ast} (G)(L^s)$. Our assertion follows from
  the commutative diagram,   in which the composition of the two upper horizontal arrows is $G(\sigma)$, 
  \begin{displaymath}
   \xymatrix{
     G(L^s) \ar[r]^{G[\sigma]} \ar[d]_{\kappa} & \hat{\sigma}^{\ast}(G)(L^s)
     \ar[d]_{\hat{\sigma}^{\ast} (\kappa)} \ar[r]^{\sim} & G(L^s) \ar[ld]^{\kappa'} \\
\Gamma \ar[r]_{\Gamma[\sigma] = \id} & \Gamma \\ 
}
  \end{displaymath}
\end{proof}

Let $K/F, r, E \rightarrow \bar{\mathbb{Q}}_p, \nu$ be as in section \ref{ss:tsv}.  Let
$\mathfrak{p}$ be a banal prime of $K$. Let
$(X_{\mathfrak{p},0}, \iota_{\mathfrak{p},0})$ be the unique CM-pair of CM-type
$r_{\mathfrak{p}}/2$ over $\Spec O_{\breve{E}_{\nu}}$. We set
\begin{displaymath}
  \Xi_{\mathfrak{p}} = \{\varphi \in  \Phi_{\mathfrak{p}}    \mid r_{\varphi} =2  \},  
\end{displaymath}
where $\Phi_{\mathfrak{p}}=\Hom_{\text{$\BQ_p$-Alg}}(K_{\mathfrak{p}}, \bar{\mathbb{Q}}_p)$, as in \eqref{decPhi}. We consider the homomorphism
\begin{displaymath}
  \mu_{\mathfrak{p}}: \bar{\mathbb{Q}}_p^{\times} \rightarrow
  (K_{\mathfrak{p}} \otimes \bar{\mathbb{Q}}_p)^{\times} \overset{\sim}{\rightarrow}
  \prod_{\Phi_{\mathfrak{p}}} \bar{\mathbb{Q}}_p^{\times} 
\end{displaymath}
such that the component of $\mu_{\mathfrak{p}}(a)$, $a \in \bar{\mathbb{Q}}_p$,
is equal to $a$ for $\varphi \in \Xi_{\mathfrak{p}}$ and is $1$ for
$\varphi \notin \Xi_{\mathfrak{p}}$. This morphism is defined over $E_{\nu}$.
We define the local reciprocity law $\mathfrak{r}_{\mathfrak{p}}$ as 
\begin{equation}\label{KM4e}
  \mathfrak{r}_{\mathfrak{p}}: E_{\nu}^{\times}
  \overset{\mu_{\mathfrak{p}}}{\longrightarrow} 
  (K_{\mathfrak{p}} \otimes_{\mathbb{Q}} E_{\nu})^{\times}
  \overset{\Nm_{E_{\nu}/\mathbb{Q}_p}}{\longrightarrow} K_{\mathfrak{p}}^{\times}. 
\end{equation}
Let $I_{\nu} \subset \Gal(\breve{E}_{\nu}^{\rm ab}/\breve{E}_{\nu})$ be the
inertia group. As before \eqref{LS37e}, let $E^{c}_{\nu}$ be the algebraic closure of
$\breve{E}_{\nu}$ in the completion of $\bar{\mathbb{Q}}_p$. By the reciprocity
law of local class field theory, we define
\begin{displaymath}
  \rho_{\mathfrak{p}}: \Gal(E^c_{\nu}/\breve{E}_{\nu})
    \rightarrow \Gal(\breve{E}_{\nu}^{\rm ab}/\breve{E}_{\nu}) \overset{\mathrm{rec}}{\simeq} O_{E_{\nu}}^{\times}
    \overset{\mathfrak{r}_{\mathfrak{p}}^{-1}}{\longrightarrow}
    O_{K_{\mathfrak{p}}}^{\times}. 
\end{displaymath}
\begin{proposition}\label{KM2p}
  Let $\mathfrak{p}$ be a banal prime of $K$. Let
$\chi_{\mathfrak{p},0}:\Gal(E^c_{\nu}/\breve{E}_{\nu})\rightarrow O_{K_{\mathfrak{p}}}^{\times}$
be the character given by the action  on
$T_p(X_{\mathfrak{p},0})$, compare (\ref{LS37e}). 
Then the restriction of this character to the
subgroup
\begin{displaymath}
  \Gal(E^c_{\nu}/\breve{E}_{\nu} \varphi_0(K_{\mathfrak{p}_v})) \subset
  \Gal(E^c_{\nu}/\breve{E}_{\nu})
  \end{displaymath}
coincides with the restriction of $\rho_{\mathfrak{p}}$ to this subgroup. 
\end{proposition}
We remark that $\breve{E}_{\nu} \varphi_0(K_{\mathfrak{p}_v})$ equals
$\breve{E}_{\nu}$ if $\mathfrak{p}_v$ is unramified in $K/F$ and is a quadratic
extension of $\breve{E}_{\nu}$ if $\mathfrak{p}_v$ is ramified in $K/F$.
\begin{proof}
  It follows from the functoriality of $\mathrm{rec}$ that the proposition
  implies the same statement for a finite extension $E'_{\nu}$ of $E_{\nu}$. 

  We define a CM-type $\Xi \subset \Phi = \Hom_{\text{$\mathbb{Q}$-Alg}}(K, \mathbb{C})$
  by choosing $\varphi_0: K \rightarrow \mathbb{C}$ with $r_{\varphi_0} = 1$ and setting
  \begin{displaymath}
\Xi = \{\varphi \in \Phi \; | \; r_{\varphi} = 2 \} \cup \{ \varphi_0 \}. 
  \end{displaymath}
  We denote by $H$ the reflex field of $\Xi$. We find that
  $H \varphi_0(K) = E \varphi_0 (K)$. We claim that there exists  an extension of number
  fields $L/H$ which is unramified at $\nu$ and a tuple
  $(A, \iota, \bar{\lambda}, \kappa)$ which is defined over 
  $L$   and such that $A$ has
  good reduction $\tilde{A}$ over $O_{L_{\nu}}$. Let
  $Y$ be the $p$-divisible group of $\tilde{A}$, which we write as $Y = \prod_{\mathfrak{p}} Y_{\mathfrak{p}}$, where $\mathfrak{p}$ runs through the
  prime ideals of $O_K$ over $p$. Let
  $\mathfrak{p}$ be banal. Then $(Y_{\mathfrak{p}}, \iota)\otimes_{O_{L_\nu}}{O_{\breve{L}_{\nu}}}$ is a
    CM-pair of type $r_{\mathfrak{p}}/2$ which satisfies the Kottwitz condition
    and the Eisenstein condition. Therefore it is isomorphic to
    $(X_{\mathfrak{p},0}, \iota_{X_{\mathfrak{p},0}})$ which is defined over
    $O_{\breve{H}_{\nu}}$. Therefore the proposition follows from Proposition
    (\ref{KM1p}).

It remains to show the existence of $L$. We fix an open compact subgroup
$C \subset T(\mathbb{A}_f)$ which is maximal in $p$ and is small enough.
The Shimura variety $\mathrm{Sh}_{\Xi, C}$ which is associated to
$(T, \mu)$ and $C$ is representable by a moduli problem
$\mathcal{A}_{\Xi, C, H}$ which is finite and \'etale over $H$. Moreover
it has a model $\mathcal{A}_{\Xi, C}$ over $O_{H_{\nu}}$. It is defined
exactly in the same way as $\mathcal{A}_{\mathbf{K}}$. Since for the moduli
problem $\mathcal{A}_{\Xi, C}$ each prime $\mathfrak{p}$ of $O_K$ is banal,
it is representable by a finite \'etale scheme over $O_{H_{\nu}}$. We conclude
the $\mathcal{A}_{\Xi, C,H}= \coprod_{i=1}^{m} \Spec L_i$ for some finite field
extensions $L_i/H$ which are unramified over $\nu$. Restricting the universal
abelian scheme over $\mathcal{A}_{\Xi, C,H}$ to some $L = L_i$,  we obtain
a tuple as required. 
\end{proof}

\section{Appendix: Adjusted invariants}  \label{s:appA}
        In this appendix we first collect some facts about anti-hermitian forms. Then we give a correction to \cite[Prop. 3.2]{KRnew}, by introducing the $r$-adjusted invariant of a CM-triple. Finally, we relate the $r$-adjusted invariant to the contracting functor of section \ref{s:tcf}.

\subsection{Recollections on binary anti-hermitian forms over $p$-adic local fields}\label{ss:binform}

We first recall the invariant of an anti-hermitian form in the case relevant to us. A good reference for this material is \cite{J}. 

Let $K/F$ be a quadratic extension of fields of characteristic $0$. We
denote by $a \mapsto \bar{a}$ the non-trivial automorphism of $K$ over
$F$. Let $V$ be an 2-dimensional vector space over $K$. Let
\begin{displaymath}\label{Kneunkap1e} 
  \varkappa: V \times V \longrightarrow K, 
\end{displaymath}
be a sesquilinear form which is linear in the first argument
and anti-linear in the second. We assume that $\varkappa$ is
anti-hermitian: 
\begin{displaymath}
  \varkappa(x, y) =  - \overline{\varkappa(y,x)} .
\end{displaymath}
We choose a basis $\{ v_1, v_{2}\}$ of $V$. Then  $\det (\varkappa(v_i, v_j))_{i,j \in \{1,2\}} \in F^\times$.  We denote by
\index[NO]{DCA@$\mathfrak{d}_{K/F} (V, \varkappa)$}
\begin{equation}\label{KneunDisc1e}
  \mathfrak{d}_{K/F} (V, \varkappa) \in F^\times/\Nm_{K/F} K^{\times}
\end{equation}
the residue class of this element. It is independent of the choice of
the basis and is called the \emph{discriminant} of $(V, \varkappa)$. \index{discriminant of an anti-hermitian form}

\begin{definition}\label{invdet}
   Let $F$ be a $p$-adic local  field and $K/F$  a  quadratic field extension. Let $(V, \varkappa)$ be
   a $K$-vector space of dimension $2$ with an anti-hermitian form
   \index[NO]{IAC@$\inv(V,\varkappa)$}
  $\varkappa$ which is nondegenerate. We denote by $\inv (V, \varkappa) \in \{\pm 1\}$ the image of
  $\mathfrak{d}_{K/F} (V, \varkappa)$ under the canonical isomorphism
  $F^{\times}/\Nm_{K/F} K^{\times} \simeq \{\pm 1\}$. The invariant determines $(V, \varkappa)$ up to isomorphism, cf. \cite{J}.
\end{definition}

We note that an
anti-hermitian form $\varkappa$ can  equivalently be given by an alternating
non-degenerate $\mathbb{Q}_p$-bilinear form
\begin{equation}\label{Kneunpsi1e}
  \psi: V \times V \longrightarrow \mathbb{Q}_p 
\end{equation}
such that 
\begin{equation*}
  \psi(ax, y) = \psi(x, \bar{a}y), \quad x,y \in V, \; a \in K.
\end{equation*}
The anti-hermitian form $\varkappa$ is defined by the
equation
\begin{displaymath}
  \Trace_{K/\mathbb{Q}_p} a \varkappa(x, y) = \psi(ax, y). 
\end{displaymath}
In this case we set
\begin{displaymath}
  \inv(V,\psi) = \inv(V, \varkappa). 
\end{displaymath}
The invariant $\inv(V,\psi)$ determines $(V,\psi)$ up to isomorphism.

Let $\Lambda \subset V$ be an $O_K$-lattice such that $\psi$ induces a
pairing
\begin{equation}\label{Kneunpsi2e}
  \psi: \Lambda \times \Lambda \longrightarrow \mathbb{Z}_{p} ,
\end{equation}
i.e., $\psi$ is integral on $\Lambda$. 
We consider the map
\begin{equation}\label{Kneunpsi3e}
  \begin{aligned}
    \Lambda & \longrightarrow  \Hom_{\mathbb{Z}_p}(\Lambda, \mathbb{Z}_p)\\
    y  & \mapsto  \ell_y, \quad \text{where}\; \ell_y(x) = \psi(x,y)
  \end{aligned}
\end{equation}
This is an anti-linear map of $O_K$-modules. Therefore the image of
this map is an $O_K$-submodule. We denote the length of the cokernel
as an $O_K$-module by ${\tt {h}} {(\Lambda,\psi)}$.
\index[NO]{HDB@${\tt {h}}(\Lambda,\psi)$}

\begin{lemma}[\cite{J}, Thm. 7.1]\label{Kneun0l}
 Let  $F$ be a local $p$-adic field and $K/F$  an unramified field extension. Let $V$
  be a $2$-dimensional $K$-vector space. Let
  \begin{displaymath}
    \psi: V \times V \longrightarrow \mathbb{Q}_{p} 
  \end{displaymath}
  as in (\ref{Kneunpsi1e}). Then $\inv(V, \psi) = 1$ iff there exists an $O_K$-lattice
  $\Lambda \subset V$ such that $\psi$ is integral on $\Lambda$  and such that ${\tt h}(\Lambda, \psi) = 0$, i.e., such that  $ \psi|_{ \Lambda \times \Lambda}$ is a perfect pairing. 
  Moreover, $\Lambda$ is uniquely determined up to $\Aut(V, \psi)$.
  
Similarly, $\inv(V, \psi) = -1$ iff there exists an $O_K$-lattice
  $\Lambda \subset V$, such that $\psi$ is integral on $\Lambda$ 
  and such that ${\tt h} {(\Lambda, \psi)} = 1$.
  Moreover, $\Lambda$ is uniquely determined up to $\Aut(V, \psi)$. In this case, $ \psi|_{ \Lambda \times \Lambda}$ is called 
   \emph{almost perfect}.
\end{lemma}
\begin{proof} This  reduces to the analogous statement for the anti-hermitian form $\tilde\psi: V \times V \longrightarrow K$ defined by  
\begin{displaymath}
  \mathbf{t}(\xi\tilde{\psi}(x_1, x_2)) = \psi(\xi x_1, x_2),  \quad x_1, x_2 \in V, \, 
  \xi \in K. 
\end{displaymath}
where \index[NO]{TAA@$\mathbf{t}$} $\mathbf{t}: K\to \BQ_p$ is defined by $\mathbf{t}(a)=\tr_{K/\BQ_p}(\vartheta^{-1} a)$, where $\vartheta$
\index[NO]{ZZHB@$\vartheta$} denotes the different of $K/\BQ_p$.
Then it follows from loc.~cit.
\end{proof}
\begin{lemma}[\cite{J}, Prop. 8.1 a)]\label{Kneun01l} Let $p\neq 2$, and
  let $F$ be a local $p$-adic field and $K/F$  a ramified quadratic field extension. Let $V$
  be a $2$-dimensional $K$-vector space. Let
  \begin{displaymath}
    \psi: V \times V \longrightarrow \mathbb{Q}_{p}.
  \end{displaymath}
  as in (\ref{Kneunpsi1e}). Then there exists an $O_K$-lattice $\Lambda \subset V$ such that $\psi$
  induces a perfect form
  \begin{displaymath}
    \psi: \Lambda \times \Lambda \longrightarrow \mathbb{Z}_{p}. 
  \end{displaymath}
  Moreover $\Lambda$  is unique up to $\Aut(V, \psi)$.\qed
\end{lemma}
        
   \subsection{The $r$-adjusted invariant}\label{ss:radj} 
      
   Let $K$ be a CM-field, with totally real subfield $F$. We set\index{$r$-adjusted invariant}
   $\Phi = \Hom_{\text{$\mathbb{Q}$-Alg}}(K, \bar{\mathbb{Q}})$. Let $r$ be a generalized CM-type\index{generalized CM-type} 
   of  rank $n$, i.e., $r_{\varphi} + r_{\bar{\varphi}} = n$ for all $\varphi\in\Phi$. Throughout this subsection, we assume that $n$ is \emph{even}. Let $E=E_r$ be the reflex
   field, cf.  \cite[\S 2]{KRnew}. A  CM-\emph{triple}\index{CM-triple} over an $O_E$-algebra $R$ is
   a triple $(A, \iota, \lambda)$ where $A$ is an abelian scheme over $R$ with
   an action $\iota: O_K \rightarrow \End A$ with satisfies the Kottwitz
   condition $({\rm KC}_r)$ and a polarization $\lambda$ whose Rosati involution
   induces the conjugation of $K/F$. In the case $n=2$ this is a CM-triple
   with satisfies the Kottwitz condition, cf. section \ref{ss:loctriples}. Let $v$ be a place of $F$. We define
   an
   $r$-{\it adjusted} invariant $\inv^r_v(A, \iota, \lambda)$ attached to a triple  $(A, \iota, \lambda)$ of CM-type $r$, defined over a field $k$ that is at the same time an $O_E$-algebra.   When $v$ is non-archimedean split in $K$, then $\inv^r_v(A, \iota, \lambda)=\inv_v(A, \iota, \lambda)=1$. If $v$ is archimedean, or non-archimedean non-split in $K$, with residue characteristic of $v$  different from the characteristic of $k$, then $\inv^r_v(A, \iota, \lambda)=\inv_v(A, \iota, \lambda)$, i.e., the adjusted invariant coincides with the invariant of \cite[\S3]{KRnew}. Comp. section \ref{ss:inv} for the definition of the latter invariant  for  $n=2$. The case of general  even $n$ is substantially the same.

        Now let $v$ be non-split with residue characteristic equal to the characteristic $p$ of $k$. We may assume that $k$ is algebraically closed. Let us first assume that the $O_E$-algebra structure of $k$ is induced by a $O_{\ov\BQ}$-algebra structure. Let $\tilde\nu$ be the induced  $p$-adic place of $\ov\BQ$. Let
        \begin{equation}
          \Phi_v=\{\varphi: K\to \ov\BQ\mid \tilde\nu\circ\varphi\text{ induces $v$ }\} .
        \end{equation}
        Then
        $$
        \Phi_v=\Hom_{\BQ_p}(K_v, \ov\BQ_{\tilde\nu} ). 
        $$
        Also let 
        $$
        r_v=r_{|\Phi_v} . 
        $$
        Now define \index[NO]{IAD@$\inv^r_v(A, \iota, \lambda)$}
        \begin{equation}
          \inv^r_v(A, \iota, \lambda)=\inv_v(A, \iota, \lambda)\sgn(r_v) , 
        \end{equation}
        with\index[NO]{SAB@$\sgn(r_v)$}
        \begin{equation}\label{signr}
          \sgn(r_v)=(-1)^{( \frac{n}{2}d_v-\sum_{\varphi\in \Phi_v^+}r_\varphi )}=(-1)^{\frac{1}{2}\sum_{\varphi\in \Phi_v^+}(r_\varphi -r_{\bar\varphi})}.
        \end{equation}  
        Here $\Phi_v^+$ is a half-system of embeddings in $\Phi_v$, which has cardinality $d_v=[F_v:\BQ_p]$. Since $r_\varphi+r_{\ov\varphi}=n$ for all $\varphi\in\Phi_v$, and $n$ is supposed to be even, \eqref{signr} is independent of $\Phi_v^+$. 
        Note that $\sgn(r_v)$ only depends on the place $\nu$ of $E$ induced by $\wt\nu$. 

        The correct version of \cite[Prop. 3.2]{KRnew} is now as follows.
        \begin{proposition}\label{contcorr}
Let $S$ be an $O_E$-scheme. Let $(A, \iota, \lambda)$ be a CM-triple over $S$ 
which satisfies $({\rm KC}_r)$. Let $c \in \{\pm1\}$. Then for every place
$v$ of $F$, the set of points $s \in S$ such that
\begin{displaymath}
\inv^r_v(A_s, \iota_s, \lambda_s) = c 
  \end{displaymath}
is open and closed in $S$. 
  \end{proposition}
\begin{proof} Clearly we may assume that $S$ is an $O_E$-scheme of finite
  type. Further we can assume that $S$ is irreducible. Obviously the invariant
  is constant on the generic fiber of $S$. Also, we may assume that $v$ is non-archimedean non-split in $K$. 
  
  First we consider the case when $S$ is an irreducible scheme of finite type
over $\kappa_{E_\nu}$. Since each local ring of $S$ is dominated by a discrete
valuation ring $R$, it is enough to consider the case $S = \Spec R$. We may
replace $R$ by a discrete valuation ring that dominates $R$. Therefore
we can assume that $R$ is complete with algebraically closed residue class field,
i.e., $R \cong k[[t]]$ for an algebraically closed field $k$. According to the action of $F\otimes\BQ_p$, the $p$-divisible
group $X$ of $A$ is  isogenous to a product $\prod_{w\mid p} X_{w}$. We consider the factor $X_v$. Let 
 $\mathcal{P}$ be the display of $X_{v}$ over $R$, cf. (\ref{Laufunctor1e}).
We note that $P$ is the value of the crystal of $X_{v}$ at the
$pd$-thickening $W(R)/R$. 
 By Lemma \ref{exofx}
below,  there is an element
$x \in \wedge^n_{O_{K_{v}} \otimes_{\mathbb{Z}_p} W(R)} P$
such that $Fx = p^{n/2} x$. We define the anti-hermitian form 
\begin{displaymath}
  \varkappa: P_{\mathbb{Q}} \times P_{\mathbb{Q}} \rightarrow
  K_v \otimes_{\mathbb{Z}_p} W(R) 
\end{displaymath}
as in (\ref{Kneun03e}). We consider the hermitian form  
$h = \wedge^n_{O_{K_{v}} \otimes_{\mathbb{Z}_p} W(R)} \varkappa$ on
$\wedge^n_{O_{K_{v}} \otimes_{\mathbb{Z}_p} W(R)} P_{\mathbb{Q}}$. From the equation
\begin{displaymath}
  h(F y_1, F y_2) = p^n ~^{F} h(y_1, y_2), \quad y_1, y_2 \in
  \wedge^n_{O_{K_{v}} \otimes_{\mathbb{Z}_p} W(R)} P_{\mathbb{Q}} ,
\end{displaymath}
we obtain that $h(x,x)$ lies in the invariants
$(K_v \otimes_{\mathbb{Z}_p} W(R))^{F} = K_v$. Because $h$ is hermitian, we obtain
$h(x,x) \in F_v$. The element $x$ can be used to determine  the invariant of the Dieudonn\'e module $P\otimes_{W(R)} W(L)$ 
obtained for arbitrary base change $R \rightarrow L$  to a perfect field. 
Therefore $\inv_v(A_s, \iota_s, \lambda_s)=\inv_v(A_\eta, \iota_\eta, \lambda_\eta)$ and $\inv^r_v(A_s, \iota_s, \lambda_s)=\inv^r_v(A_\eta, \iota_\eta, \lambda_\eta)$, where $s$ and $\eta$ denote the special and the generic point of $\Spec R$. 
For the comparison with the definition of the invariant of a Dieudonn\'e
module we should remark that the equations $Fx = p^{n/2} x$ and
$Vx = p^{n/2} x$ are equivalent because $F V = p^n$ on
$\wedge^n_{O_{K_{v}} \otimes_{\mathbb{Z}_p} W(R)} P$.

  Now we consider the case when the function field of $S$ has characteristic $0$.  This case can be reduced to the
  case  when $S=\Spec O_L$, where $L$ is the completion of a
  subfield of $\ov\BQ_p$ which contains $E$ and such that its ring of integers
  $O=O_L$ is a discrete valuation ring with residue field $\ov\BF_p$.  We
  denote by $A_L$ the generic fiber of $A$, and by $A_k$ its special fiber.

          We decompose the rational $p$-adic Tate module of $A_L$, resp. the rational Dieudonn\'e module of $A_k$,  with respect   to the actions of $F\otimes\BQ_p$, 
          $$
          V_p(A_L)=\bigoplus_{w\mid p} V_w(A_L), \quad  M(A_k)_{\BQ}=\bigoplus_{w\mid p} M(A_k)_{\BQ, w}  .
          $$
          Here $V_w( A_L)$ is a free $K\otimes_F F_w$-module of rank $n$, and $M(A_k)_{\BQ, w}$ is a free $K\otimes_F F_w\otimes_{\BQ_p} W(k)_\BQ$-module of rank $n$. Set $\breve \BQ_p=W(k)_\BQ$.

          Let $ S_{v} = \bigwedge^n_{K_v} \, V_v( A_L) \, $ and $N_{\BQ, v}= \bigwedge^n_{K_v} \, M(A_k)_{\BQ, v}$. Both are equipped with  hermitian forms (for the first module, cf. \cite[section 3, case b)]{KRnew}; for the second module, cf. subsection \ref{ss:inv}). Also, we have  $N_{\BQ, v}={\bf 1}_v(\frac{n}{2})$, where ${\bf 1}_v$ is a multiple of the unit object in the category of Dieudonn\'e modules, comp. \eqref{Kneun02e}, or Lemma \ref{exofx}.  Let $U_v$ be the image under the Fontaine functor of $N_{\BQ,v}(-\frac{n}{2})$. We need to compare the two hermitian vector spaces $S_v(-\frac{n}{2})$ and $U_v$.

          Let $T$ be the torus over $\BQ_p$ which is the kernel of the  map defined by the norm of $K_v/F_v$,
          $$
          1\to T\to \Res_{K_v/\BQ_p}\BG_{m, K_v}\to \Res_{F_v/\BQ_p}\BG_{m, F_v}\to 1 .
          $$
          Then $H^1(\BQ_p, T)=F_v^\times/\Nm (K_v)$. We may regard the isomorphisms of hermitian vector spaces $\Isom (U_v, S_v(-\frac{n}{2}))$ as an etale sheaf on $\Spec F_v$. This is a $T$-torsor.  Its class $\cl(U_v,  S_v(-\frac{n}{2}))$ is calculated by \cite[Prop. 1.~20]{RZ}.

          To evaluate this formula, note that the first summand, $\kappa (b)$, in loc.~cit. is trivial. To evaluate the second summand, $\mu^\sharp$, we use the following description of the filtration on $N_{\BQ, v}\otimes_{\breve\BQ_p}\ov{\breve \BQ}_p$. For the filtration of $M(A_k)_{\BQ, v}\otimes_{\breve\BQ_p}\ov{\breve \BQ}_p=\oplus_{\varphi\in \Phi_v}M(A_k)_{\BQ, v, \varphi}$ we have that the jumps are  in degree $0$ and $ 1$, with 
          \begin{equation}
            (0)\subset {\rm Fil}^1_\varphi\subset^{r_\varphi} M(A_k)_{\BQ, v, \varphi} .
          \end{equation}
          The upper index means that the cokernel has dimension $r_\varphi$. For the filtration of the one-dimensional  vector space $N_{\BQ, v, \varphi}$, this means that the unique jump is in degree $n-r_\varphi$. We use the identification
          $$
          X_*(T)=\Ker \big(\Ind_{F_v}^{K_v}(\Ind_{\BQ_p}^{F_v}(\BZ))\to \Ind_{\BQ_p}^{F_v}(\BZ)\big) .
          $$
          Then the corresponding filtration on $N(A)_{\BQ, v, \varphi}(-\frac{n}{2})$ is  given by the cocharacter $\mu\in X_*(T)$ with 
          \begin{equation}
            \mu_\varphi=\frac{n}{2}-r_\varphi , \quad \varphi\in\Phi_v . 
          \end{equation}
          We have to determine the image $\mu^\sharp$ of $\mu$ in $X_*(T)_\Gamma$. Under the identification $X_*(T)_\Gamma=H^1(\BQ_p, T)=\BZ/2$, we obtain 
          \begin{equation}
            \cl(U_v,  S_v(-\frac{n}{2}))=\mu^\sharp=\sum_{\varphi\in\Phi_v^+} \mu_\varphi= \frac{n}{2}d_v-\sum_{\varphi\in \Phi_v^+}r_\varphi\, ,
          \end{equation}
          where we used the notation introduced for \eqref{signr}. We deduce $\inv^r_v(A_k, \iota_k, \lambda_k)=\inv_v(A_L, \iota_L, \lambda_L)$, as desired. 
          \end{proof}
In the proof of Proposition \ref{contcorr}, we used the following lemma. 
\begin{lemma}\label{exofx}
  Let $F/\mathbb{Q}_p$ be a finite field extension of degree $d$ and $K/F$
  be a quadratic field extension. Let $n$ be an even natural number. 
  Let $k$ be an algebraically closed field of characteristic $p$. Let
  $(X, \iota)$ be a $p$-divisible group over $k[[t]]$ of dimension $nd$ and
  height $2nd$ with an action $O_K \rightarrow \End X$. 
  Let $(\mathcal{P}, \iota)$ be the display of $X$, cf. (\ref{Laufunctor1e}). 
  Then there exists a non-zero element 
  $x \in \wedge^n_{O_K \otimes_{\mathbb{Z}_p} W(k)[[t]]} P$ such that
  \begin{displaymath}
\wedge^n F (x) = p^{n/2} x. 
    \end{displaymath}
    \end{lemma}
\begin{proof}
We consider the $\mathbb{Z}_p$-frame
$\mathcal{B}_k = (W(k)[[t]], pW(k)[[t]], k[[t]], \sigma, \dot{\sigma})$,
where $\sigma$ is the extension of the Frobenius  on $W(k)$ to the
power series ring given by $\sigma(t) = t^p$ and
where $\dot{\sigma} = (1/p) \sigma$. 
The evaluation $P_1$ of the crystal of $X$ at the $pd$-thickening
$W(k)[[t]]/k[[t]]$ has the structure of a $\mathcal{B}_k$-display.
The display $\mathcal{P}$ is obtained by base change with respect to a
morphism of frames 
$\mathcal{B}_k \rightarrow \mathcal{W}(k[[t]])$, cf. \cite{Zi2} and \cite{L3}. 
Therefore, it is enough to prove our assertion for the
$\mathcal{B}_k$-display of $X$ which we will now denote by $\mathcal{P}$.

We consider first the case when $K/F$ is ramified. 
  When writing $\det_{W(k)[[t]]} F$, we mean this with respect to an arbitrary
  $W(k)[[t]]$-basis of $P$. This determinant is well determined up to multiplication with
  a unit in $W(k)[[t]]$. We know that
  $\det_{W(k)[[t]]} F = p^{hd} u_1$ for some $u_1 \in (W(k)[[t]])^{\times}$.

  We consider the decomposition $P = \oplus P_{\psi}$ according to
  \begin{displaymath}
O_K \otimes W(k)[[t]] = \prod_{\psi} O_K \otimes_{O_{F^t}, \tilde{\psi}} W(k)[[t]]. 
    \end{displaymath}

The Frobenius is graded, $F: P_{\psi} \rightarrow P_{\psi \sigma}$. We conclude
that $\det_{W(k)[[t]]} (F^f | P_{\psi}) = p^{nd} u_2$ for some unit
$u_2 \in (W(k)[[t]])^{\times}$. We fix $\psi$. 
Up to a unit we have
\begin{equation}\label{Kneuninv31e}
  \Nm_{K/F^t} \textstyle{\det_{O_K \otimes_{O_{F^t}, \tilde{\psi}} W(k)[[t]]}} (F^f | P_{\psi}) 
 = \det_{W(k)[[t]]} (F^f | P_{\psi}). 
\end{equation}
We fix a normal extension $L$ of $W(k)_{\mathbb{Q}}$ which contains
$K \otimes_{O_{F^t}, \tilde{\psi}} W(k)$. The left hand side of (\ref{Kneuninv31e})
is the product of conjugates $c_1, \ldots, c_{2e} \in O_L[[t]]$ of
$\det_{O_K \otimes_{O_{F^t}, \tilde{\psi}} W(k)[[t]]} (F^f | P_{\psi})$. These elements
have the same order with respect to the prime element $\omega_L$ of $L$
which is a prime element in the regular local ring $O_L[[t]]$. We rewrite
(\ref{Kneuninv31e})
\begin{displaymath}
c_1 \cdot c_2 \cdot \ldots \cdot c_{2e} = p^{nd} \cdot u_3,
  \end{displaymath}
for some unit $u_3$. Since $O_L[[t]]$ is factorial, we find
$c_i = p^{fn/2} \mu_i$ for some units $\mu_i$. We conclude that
\begin{displaymath}
\textstyle{\det_{O_K \otimes_{O_{F^t}, \tilde{\psi}} W(k)[[t]]}} (F^f | P_{\psi}) = p^{fn/2} u_4
  \end{displaymath}
for some unit $u_4 \in (O_K \otimes_{O_{F^t}, \tilde{\psi}} W(k)[[t]])^{\times}$.
Since $P_{\psi}$ is a free $O_K \otimes_{O_{F^t}, \tilde{\psi}} W(k)[[t]]$-module
of rank $n$, we find that, for each element 
$y_{\psi} \in \wedge^n_{O_K \otimes_{O_{F^t}, \tilde{\psi}} W(k)[[t]]} P_{\psi}$, there is
an equation 
$\wedge^n F^f y_{\psi} = p^{fn/2} u(y_{\psi}) y_{\psi}$ for some unit
$u(y_{\psi}) \in O_K \otimes_{O_{F^t}, \tilde{\psi}} W(k)[[t]]$. On the last ring,
$\sigma^f$ acts via the second factor. There is a unit
$\zeta \in O_K \otimes_{O_{F^t}, \tilde{\psi}} W(k)[[t]]$ such that
\begin{displaymath}
\sigma^f (\zeta) \zeta^{-1} = u(y_{\psi}). 
  \end{displaymath}
Indeed, consider the image $\bar{u}$ of $u(y_{\psi})$ in
$O_K \otimes_{O_{F^t}, \tilde{\psi}} W(k)$ by setting $t=0$. It is well-known
that in this ring $\sigma^f (\bar{\zeta}) \bar{\zeta}^{-1} = \bar{u}$ is
solvable. One can lift $\bar{\zeta}$ successively modulo $t^n$ to a solution
$\zeta$. Then $x_{\psi} = \zeta y_{\psi}$ satisfies
\begin{displaymath}
\wedge^n F^f x_{\psi} = p^{fn/2} x_{\psi}. 
\end{displaymath}
We define $x_{\psi \sigma^i} \in P_{\mathbb{Q}}$ by
$\wedge^n F^i x_{\psi} = p^{in/2} x_{\psi \sigma^i}$ for
$i = 1, \ldots, f$. Then $x = (x_\psi) \in P \otimes \mathbb{Q}$ satisfies
$\wedge^h F (x) = p^{n/2} x$. Multiplying by a power of $p$ we can arrange that
$x \in P$. 

The proof in the unramified case is almost the same. We indicate the
differences. In this case $\Hom_{\text{$\mathbb{Q}_p$-Alg}}(F^t, W(k)_{\mathbb{Q}})$ has
$2f$ elements. Therefore we have the equation \eqref{Kneuninv31e} with $f$
replaced by $2f$,
\begin{displaymath}
\textstyle{\det_{O_K \otimes_{O_{F^t}, \tilde{\psi}} W(k)[[t]]}} (F^{2f} | P_{\psi}) = p^{fn} u_4 .
\end{displaymath}
We define $x_{\psi \sigma^i} \in P_{\mathbb{Q}}$ by
$\wedge^n F^i x_{\psi} = p^{in/2} x_{\psi \sigma^i}$ for $i = 1, \ldots, 2f$.
Then $x = (x_\psi) \in P \otimes \mathbb{Q}$ satisfies
$\wedge^h F (x) = p^{n/2} x$.
\end{proof}

        \begin{remarks}
          (i) The remarks and results  on a product formula at the end of \S 3 of \cite{KRnew} become correct when the invariants  $\inv_v(A, \iota, \lambda)$ are replaced by the adjusted invariants $\inv^r_v(A, \iota, \lambda)$. 

          (ii) In the definition of $\CM_{r, {\tt h}, V}$ in \cite[(4.3)]{KRnew}, the invariants  $\inv_v(A, \iota, \lambda)$ have to be replaced by the adjusted invariants $\inv^r_v(A, \iota, \lambda)$.

          (iii) One defines in the obvious way the $r$-adjusted invariant $\inv^r(X, \iota, \lambda)$ of a {\it local} CM-triple of type $r$, $(X, \iota, \lambda)$, over a field of characteristic $p$.
        \end{remarks}

\subsection{$r$-adjusted invariant and the contracting functor}\label{ss:invcont}
In this subsection, we return to the situation in section \ref{ss:specialbanal}. We assume that $K/F$ is a field extension. 
Let $k$ be an algebraically closed field of
characteristic $p$ with an $O_F$-algebra structure, i.e., $k \in \Nilp_{O_F}$.

We consider the case where $r$ is special. Consider an object
$(\mathcal{P}_{\rm{c}}, \iota_{\rm{c}}, \beta_{\rm{c}}) \in \mathfrak{d}{\mathfrak R}^{\rm pol}_{k}$
cf. Definition \ref{P'Pol3d}. We write $\mathcal{P}_{\rm{c}} = (P_{\rm{c}}, F_{\rm{c}}, V_{\rm{c}})$ for the corresponding 
$\mathcal{W}_{O_F}(k)$-Dieudonn\'e module. To avoid too many double notations
we denote the Frobenius automorphism on $W_{O_F}(k)$ by $\tau$. The
Verschiebung on $W_{O_F}(k)$  is then $\pi \tau^{-1}$. 
For our purposes it is more convenient to allow
quasi-polarizations, i.e., $\beta_{\rm{c}}$ is a $W_{O_F}(k)$-bilinear form
 \begin{displaymath}
      P_{\rm{c}} \otimes \mathbb{Q} \; \times \; P_{\rm{c}} \otimes \mathbb{Q} \longrightarrow
      W_{O_F}(k)_{\mathbb{Q}}, 
 \end{displaymath}
 such that 
$(\mathcal{P}_{\rm{c}}, \iota_{\rm{c}}, p^t \beta_{\rm{c}}) \in \mathfrak{d}{\mathfrak R}^{\rm pol}_{k}$ for large enough $t \in \mathbb{Z}$. 
 Then $\beta_{\rm{c}}$ is alternating and the following
 equations hold:
 \begin{displaymath}
   \begin{aligned}
     \beta_{\rm{c}}(F_{\rm{c}} u_1, F_{\rm{c}} u_2) &= \pi \tau(\beta_{\rm{c}}(u_1, u_2)),\quad
      u_1, u_2 \in P_{\rm{c}} \otimes \mathbb{Q},\\ 
\beta_{\rm{c}}(\iota_{\rm{c}}(a)u_1, u_2) &= \beta_{\rm{c}}(u_1, \iota_{\rm{c}}(\bar{a})u_2),  \quad a \in K.
     \end{aligned}
 \end{displaymath}\label{invcontr1e}
The polarization $\beta_{\rm{c}}$ defines an anti-hermitian form
\begin{equation}\label{betacontr}
  \varkappa_{\rm{c}}: P_{\rm{c}} \otimes \mathbb{Q} \times  P_{\rm{c}} \otimes \mathbb{Q}
  \longrightarrow K \otimes_{O_F} W_{O_F}(k),
\end{equation}
by the formula
\begin{displaymath}
  \Trace_{K/F} a \varkappa_{\rm{c}}(u_1, u_2) = \beta_{\rm{c}}(au_1, u_2), \; a \in K
  \otimes_{O_F} W_{O_F}(k), \; u_1, u_2 \in P_{\rm{c}} \otimes \mathbb{Q}. 
\end{displaymath}
We note that,  by Lemma \ref{Displaykristall1l}, $P_{\rm{c}} \otimes \mathbb{Q}$ is a free
$K \otimes_{O_F} W_{O_F}(k)$-module of rank two.

Since $\Lie P_{\rm{c}}$ has dimension 2, we have $\ord_{\pi} \det W_{O_F}(k)(V_{\rm{c}} | P_{\rm{c}}) = 2$.
We recall that, for an arbitrary $K \otimes_{O_F} W_{O_F}(k)$-linear map
$V_{\rm{c}}^{\sharp}: P_{\rm{c}} \otimes \mathbb{Q} \longrightarrow P_{\rm{c}} \otimes \mathbb{Q}$,
\begin{displaymath}
  \Nm_{K/F} \det\textstyle_{K \otimes_{O_F} W_{O_F}(k)}
  (V_{\rm{c}}^{\sharp} | P_{\rm{c}} \otimes \mathbb{Q}) =
  \det_{W_{O_F}(k)}(V_{\rm{c}}^{\sharp} | P_{\rm{c}} \otimes \mathbb{Q}). 
  \end{displaymath}
We conclude that 
\begin{equation}\label{Inv1e}
  \begin{aligned} 
    \ord_{\Pi} \textstyle{\det_{K \otimes_{O_F} W_{O_F}(k)}}(V_{\rm{c}} | P_{\rm{c}}) = 2, &
    \quad K/F \text{ ramified},\\
    \ord_{\pi} \textstyle{\det_{K \otimes_{O_F} W_{O_F}(k)}}(V_{\rm{c}} | P_{\rm{c}}) = 2, &
    \quad K/F \text{ unramified}. 
    \end{aligned}
\end{equation}
With our convention $\pi = \Pi$ in the unramified case, this is the same
formula. 

Let 
\begin{displaymath}
  H_{\rm{c}} = \bigwedge^2_{K \otimes_{O_F} W_{O_F}(k)} P_{\rm{c}} \otimes
  \mathbb{Q} .
\end{displaymath}
This is a free $K \otimes_{O_F} W_{O_F}(k)$-module of rank $1$. There
is an element $x_{\rm{c}} \in H_{\rm{c}}$ such that
\begin{equation}\label{Kneuninv5e}
  \wedge^2 V_{\rm{c}}\; x_{\rm{c}} = \pi x_{\rm{c}}.
\end{equation}
The existence of $x_{\rm{c}}$ follows from (\ref{Inv1e}) and the fact that the
$\tau$-conjugacy class of an element $\xi \in K \otimes_{O_F} W_{O_F}(k)$
is determined by its order, compare (\ref{Kneunord1e}). 

The anti-hermitian form $\varkappa_{\rm{c}}$ induces on $H_{\rm{c}}$ an hermitian
form 
\begin{displaymath}
  h_{\rm{c}} = \wedge^2 \varkappa_{\rm{c}}: H_{\rm{c}} \times H_{\rm{c}} \longrightarrow
  K \otimes_{O_F} W_{O_F}(k). 
\end{displaymath}
We find 

\begin{equation}
  h_{\rm{c}}(\wedge^2 V_{\rm{c}}\; x_1, \wedge^2 V_{\rm{c}}\; x_2) =
  \pi^2 \tau^{-1}(h_{\rm{c}}(x_1, x_2)) , 
\end{equation}
where $\tau$ acts on $K \otimes_{O_F} W_{O_F}(k)$ via the second factor. 
Using (\ref{Kneuninv5e}) this implies
\begin{displaymath}
  h_{\rm{c}} (x_{\rm{c}}, x_{\rm{c}}) \in F^{\times} \subset K \otimes_{O_F} W_{O_F}(k).
\end{displaymath}
 The following definition is analogous to  \eqref{Kneun3e}.
  \begin{definition}\label{invcontr1d}
 \label{def:invcontr}  The invariant  $\inv (\mathcal{P}_{\rm{c}}, \iota_{\rm{c}}, \beta_{\rm{c}}) \in \{\pm 1\}$ is defined as 
  the image of
  $h_{\rm{c}}(x_{\rm{c}}, x_{\rm{c}})$ by the canonical map
  \begin{displaymath}
    F^{\times} \longrightarrow F^{\times}/\Nm_{K/F} K^{\times}
    \isoarrow \{\pm 1\} ,  
  \end{displaymath}
  \end{definition}

The following proposition relates this invariant with the invariant \eqref{Kneun3e} under the contracting functor. 

\begin{proposition}\label{Kneuninv1p}
  Let $K/F$ be a field extension and let $r$ be special.
  Recall the reflex field $E$ associated to $r$.  Let
  $k \in Nilp_{O_{\breve{E}}}$ be an algebraically closed field. Let
  $(\mathcal{P}, \iota, \beta) \in \mathfrak{d}{\mathfrak P}^{\rm pol}_{r, k}$ and
  let
  $(\mathcal{P}_{\rm{c}}, \iota_{\rm{c}}, \beta_{\rm{c}}) \in \mathfrak{d}{\mathfrak R}^{\rm pol}_{k}$
  be its image by the contracting functor $\mathfrak{C}^{\rm pol}_{r,k}$, cf.
  (\ref{P'Pol12e}). Then 
      \begin{displaymath}
        \inv^{r} (\mathcal{P}, \iota, \beta)  = 
        \inv(\mathcal{P}_{\rm{c}}, \iota_{\rm{c}}, \beta_{\rm{c}}) . 
      \end{displaymath}
      Here the $r$-adjusted invariant is given by 
      \begin{equation*}
  \inv^{r} (\CP, \iota, \beta) =  (-1)^{d-1}\inv (\CP, \iota, \beta). 
\end{equation*} 
\end{proposition}
\begin{proof} The second assertion follows from the definition of $\sgn(r)$, cf. \eqref{signr}. Let us prove the first assertion. 

  We begin with the ramified case. We have the decomposition
  $P = \oplus_{\psi} P_{\psi}$, cf. (\ref{OWdecomp2e}). By the definition
  of the contracting functor for Dieudonn\'e modules, 
  we have
  \begin{displaymath}
P_{\rm{c}} = P_{\psi_0}, \quad V_{\rm{c}} = \Pi^{-ef+ 1} V, 
  \end{displaymath}
  cf. Remark \ref{P'Pol2r}.
  The bilinear form $\tilde{\beta}_{\rm{c}}$ on $P_{\rm{c}}$ is the restriction of
  $\tilde{\beta}$ of Proposition \ref{P'Pol1p}. Since we may change
  $\tilde{\beta}$ by a factor in $F^{\times}$ without changing the invariant,
  we may replace $\tilde{\beta}$ by $\vartheta^{-1} \tilde{\beta}$, i.e.,
  we may assume that $\Trace_{F/\mathbb{Q}_p} \tilde{\beta} = \beta$.
We define the anti-hermitian form
\begin{displaymath}
  {\varkappa}: P \otimes \mathbb{Q} \; \times \; P \otimes \mathbb{Q}
        \longrightarrow K \otimes_{\mathbb{Z}_p} W(k)
\end{displaymath}
by $\Trace_{K/F} {\varkappa} = \tilde{\beta}$. On
\begin{equation}\label{Inv2e}
  H = \bigwedge^2_{K \otimes W(k)} P \otimes \mathbb{Q}
\end{equation}
we obtain the hermitian form
${h} = \wedge^2 {\varkappa}$. We have the decomposition
\begin{displaymath}
  H = \bigoplus_{\psi} \bigwedge^2_{K \otimes_{O_{F^t}, \tilde{\psi}} W(k)} P_{\psi}
  \otimes \mathbb{Q} = \bigoplus_{\psi} H_{\psi}. 
\end{displaymath}
The hermitian form ${h}$ is the orthogonal sum of the induced forms
\begin{displaymath}
  {h}_{\psi}: H_{\psi} \times H_{\psi} \longrightarrow K \otimes_{O_{F^t},
    \tilde{\psi}} W(k). 
\end{displaymath}

To determine  $\inv(\mathcal{P}_{\rm{c}}, \iota_{\rm{c}}, \beta_{\rm{c}})$,  we consider $\varkappa_{\rm{c}}$
defined by $\Trace_{K/F} \varkappa_{\rm{c}} = \beta_{\rm{c}}$ and the hermitian form
$h_{\rm{c}} = \wedge^2 \varkappa_{\rm{c}}$ on $H_{\rm{c}} = H_{\psi_0}$.  The form ${h}_{\psi_0}$ coincides with the hermitian form deduced from the form $\tilde{\beta}_{\rm{c}}$ above.
By definition  $\beta_{\rm{c}} = \eta_{0.k}^f \tilde{\beta_{\rm{c}}}$, cf. 
(\ref{P'Pol15e}).  Hence we have 
\begin{displaymath}
  h_{\rm{c}} = \eta_{0,k}^{2f} {h}_{\psi_0}. 
\end{displaymath}
We choose an element $x \in H$ such that
\begin{displaymath}
  \wedge^{2}V ( x )= px.
\end{displaymath}
Let $x_{\psi_0}$ be the $\psi_0$-component of $x$. Then
$\inv (\mathcal{P}, \iota, \beta)$ is given by the element
\begin{displaymath}
  {h}_{\psi_0}(x_{\psi_0}, x_{\psi_0}) \in F^{\times}. 
\end{displaymath}
We set $z_{\psi_0} = \eta_{0,k}^{-f} x_{\psi_0}$. Then we find
\begin{displaymath}
  h_{\rm{c}}(z_{\psi_0}, z_{\psi_0}) = \eta_{0,k}^{2f} {h}_{\psi_0}(\eta_{0,k}^{-f}
  x_{\psi_0}, \eta_{0,k}^{-f} x_{\psi_0}) = {h}_{\psi_0}(x_{\psi_0},
  x_{\psi_0}). 
\end{displaymath}
From  $V_{\rm{c}} =  \Pi^{-d+1}V^f$ and $\wedge^{2}V^f x_{\psi_0} = p^f
x_{\psi_0}$, we obtain
\begin{equation*}
\begin{aligned}
  \wedge^2V_{\rm{c}} ( x_{\psi_0}) &= (-1)^{d-1} \pi (p/\pi^e)^{f} x_{\psi_0}\\
  \wedge^2V_{\rm{c}} ( z_{\psi_0}) &= \tau^{-1}(\eta_{0,k}^{-f}) (-1)^{d-1}
  \pi (p/\pi^e)^{f} \eta_{0,k}^{f} z_{\psi_0} =  (-1)^{d-1} \pi z_{\psi_0}.
\end{aligned}
\end{equation*}
By  Lemma \ref{Kneuninv1l} below, $h_{\rm{c}}(z_{\psi_0}, z_{\psi_0}) \in
F^{\times}$ defines $(-1)^{d-1}\inv (\mathcal{P}_{\rm{c}}, \iota_{\rm{c}}, \beta_{\rm{c}})$.

We consider now the unramified case. As before, we have $H$
with its hermitian form ${h}$, cf. (\ref{Inv2e}). We consider the decomposition
\begin{equation}\label{Kneuninv9e}
  H = \bigoplus_{\psi} \Big( \bigwedge^2_{K \otimes_{O_{K^t}, \tilde{\psi}} W(k)} P_{\psi} 
  \otimes \mathbb{Q} \Big)  = \bigoplus_{\psi} H_{\psi},  
\end{equation}
which has now $2f$ summands. Now $H_{\psi_1}$ and $H_{\psi_2}$ are orthogonal for
$\psi_1 \neq \bar{\psi}_2$. We denote by
\begin{displaymath}
  {h}_{\psi}: H_{\psi} \times H_{\bar{\psi}} \longrightarrow
  K \otimes_{O_K^t, \tilde{\psi}} W(k) 
  \end{displaymath}
the sesquilinear form induced by ${h}$. Let $x = (x_{\psi}) \in H$
such that $\wedge^2 V (x) = px$ or, equivalently,
$\wedge^2 V (x_{\psi}) = px_{\psi \sigma^{-1}}$ for all $\psi$. 
The invariant of $(\mathcal{P}, \iota, \beta)$
is the class in $F^{\times}/ \Nm_{K/F} K^{\times}$ of
${h}_{\psi}(x_{\psi}, x_{\bar{\psi}}) \in F\subset K \otimes_{O_{K^t}, \tilde{\psi}} W(k)$. This is independent of $\psi$. Equivalently,  we can consider
$\ord_{\pi} {h}_{\psi}(x_{\psi}, x_{\bar{\psi}}) \in \mathbb{Z}/2\mathbb{Z}$.
Note that $\ord_\pi$ makes sense for each element of
$K \otimes_{O_{K^t}, \tilde{\psi}} W(k)$.

The invariant of $(\mathcal{P}_{\rm{c}}, \iota_{\rm{c}}, \beta_{\rm{c}})$ is defined by
$H_{\rm{c}} = H_{\psi_0} \oplus H_{\bar{\psi}_0}$ and $h_{\rm{c}}$, via 
\begin{equation}\label{Inv3e}
\ord_{\pi} h_{\rm{c}}(x_{{\rm{c}},\psi_0}, x_{{\rm{c}}, \bar{\psi}_0}), 
\end{equation}
where $x_{\rm{c}} = (x_{{\rm{c}},\psi_0}, x_{{\rm{c}}, \bar{\psi}_0}) \in H_{\psi_0} \oplus H_{\bar{\psi}_0}$
is the element of (\ref{Kneuninv5e}). We note that we can change $h_{\rm{c}}$ and
the elements $x_{{\rm{c}},\psi_0}$, resp. $x_{{\rm{c}}, \bar{\psi}_0}$, by a unit in  
$K \otimes_{O_{K^t}, \tilde{\psi}} W(k)$ without changing (\ref{Inv3e}). In
particular, (\ref{Inv3e}) is equal to
$\ord_{\pi} h_{\psi_0}(x_{{\rm{c}},\psi_0}, x_{{\rm{c}}, \bar{\psi}_0})$.  

For an element $y = (y_{\psi_0}, y_{\bar{\psi}_0}) \in H_{\rm{c}}$ we obtain from
(\ref{P'Pol14e})
\begin{displaymath}
  \wedge^2 V_{\rm{c}} (y_{\psi_0}) = \pi^{-2g_{\bar{\psi}_0}} \wedge^2 V (y_{\psi_0}), \quad 
  \wedge^2 V_{\rm{c}} (y_{\bar{\psi}_0}) = \pi^{-2g_{\psi_0}} \wedge^2 V (y_{\bar{\psi}_0}).  
  \end{displaymath}
We set
\begin{displaymath}
z_{\psi_0} = \pi^{-g_{\psi_0}} x_{\psi_0}, \quad z_{\bar{\psi}_0} = \pi^{-g_{\bar{\psi}_0}}
x_{\bar{\psi}_0}. 
\end{displaymath}
We find
\begin{displaymath}
  \begin{aligned} 
  \wedge^2 V_{\rm{c}} (z_{\psi_0}) &= \pi^{-g_{{\psi}_0}} \pi^{-2g_{{\bar{\psi}_0}}} p^f x_{\bar{\psi}_0} =
  \pi^{-g_{\psi_0}-g_{\bar{\psi}_0}}p^f z_{\bar{\psi}_0}, \\
  \wedge^2 V_{\rm{c}} (z_{\bar{\psi}_0}) &= \pi^{-g_{{\bar{\psi}_0}} - g_{\psi_0}} p^f z_{\psi_0}. 
    \end{aligned}
  \end{displaymath}
We have $\ord_{\pi}(\pi^{-g_{\bar{\psi}_0} - g_{\psi_0}} p^f) = 1$. Therefore we obtain
an element $x_{\rm{c}}$ as in \eqref{Kneuninv5e} if we change $z_{\psi_0}$ and $z_{\bar{\psi_0}}$ by a unit, cf.
 Lemma \ref{lemminun} below. 
Therefore the invariant of $(\mathcal{P}_{\rm{c}}, \iota_{\rm{c}}, \beta_{\rm{c}})$ is
\begin{displaymath}
  \ord_{\pi} h_{\psi_0}(z_{\psi_0}, z_{{\bar{\psi}_0}}) =
  (-g_{\psi_0}-g_{\bar{\psi}_0}) +  
  \ord_{\pi} h_{\psi_0}(x_{\psi_0}, x_{{\bar{\psi}_0}}) =
(1-d) + \ord_{\pi} h_{\psi_0}(x_{\psi_0}, x_{{\bar{\psi}_0}}).
\end{displaymath}
This proves the unramified case. 
\end{proof}
In the previous proof, we used  two lemmas which we state as Lemmas \ref{Kneuninv1l} and \ref{lemminun}. 
    \begin{lemma}\label{Kneuninv1l}
      Let $K/F$ be  ramified and let $r$ be special. 
      Let $y_{\rm{c}} \in H_{\rm{c}}$ be an element such that
      \begin{displaymath}
        \wedge^2 V_{\rm{c}} (y_{\rm{c}}) = - \pi y_{\rm{c}}.
      \end{displaymath}
      Then $h_{\rm{c}}(y_{\rm{c}}, y_{\rm{c}}) \in F^{\times}$ and the image of this element in
      $\{\pm 1\}$ is $-\inv (\mathcal{P}_{\rm{c}}, \iota_{\rm{c}}, \beta_{\rm{c}})$. 
    \end{lemma}
\begin{proof}
We choose an element $\zeta \in W_{O_F}(k)^{\times}$ such that
\begin{displaymath}
  \tau^{-1}(\zeta) \zeta^{-1} = -1. 
\end{displaymath}
Then $\tau^{2}(\zeta) = \zeta$ and therefore
$\zeta \in O_{F'} \subset W_{O_F}(k)$ where $F'/F$ is the unramified
extension of degree $2$. More explicitly,  we take an element
$c \in \kappa_{F'}$ such that $\tau (c) = - c$ and define $\zeta =
[c]$ to be the Teichm\"uller representative.

We set $x_{\rm{c}} = \zeta y_{\rm{c}}$. Then the equation (\ref{Kneuninv5e}) is
satisfied. We find
\begin{equation}\label{Kneuninv7e}
  h_{\rm{c}}(x_{\rm{c}}, x_{\rm{c}}) = \zeta^2 h_{\rm{c}}(y_{\rm{c}},y_{\rm{c}}). 
\end{equation}
Since $\zeta^2\mod \pi = c^2 \in \kappa_F$ is not a square in this
field, we conclude that $\zeta^2$ is not in the image of
$\Nm_{K/F}: O_K^{\times} \longrightarrow O_F^{\times}$ since the norm is the
square on the residue fields. Therefore the image of the right hand
side of (\ref{Kneuninv7e}) in $\{\pm 1\}$ is different from the image of
$h_{\rm{c}}(y_{\rm{c}}, y_{\rm{c}})$.    
\end{proof}

The last lemma has the following variant which we need in the banal case. 
\begin{lemma}\label{Kneuninv3l}
  Let $K/F$ be ramified and let $r$ be arbitrary. 
  Let $(\mathcal{P}, \iota, \beta)$ be a CM-triple of type $r$ over an
  algebraically closed field $k$. Let
  $y \in \wedge^2_{K \otimes W(k)} P_{\mathbb{Q}}$ be an element such that
  \begin{displaymath}
    \wedge^{2} V (y) = -py.
  \end{displaymath}
  Set $h = \wedge^{2} \varkappa$. Then
  $h(y,y) \in F \subset F \otimes W(k)$ and 
  \begin{displaymath}
h(y,y) \equiv (-1)^f \inv (\mathcal{P}, \iota, \beta) \; \mod \Nm_{K/F} K^{\times}. 
    \end{displaymath}
\end{lemma}
\begin{proof}
  We consider the decomposition
  \begin{displaymath}
O_F \otimes W(k) = \prod_{\psi} O_F \otimes_{O_{F^t}, \tilde{\psi}} W(k). 
  \end{displaymath}
  We denote by $\sigma$ the Frobenius acting on $W(k)$. It induces via any of
  the embeddings $\tilde{\psi}$ the Frobenius
  $\sigma \in \Gal(F^t/\mathbb{Q}_p)$. The decomposition induces a
  decomposition $P = \oplus_{\psi} P_{\psi}$ and 
  \begin{displaymath}
    \wedge^2_{K \otimes W(k)} P_{\mathbb{Q}} = \oplus_{\psi} 
    \wedge^2_{K \otimes_{O_{F^t}, \tilde{\psi}} W(k)} P_{\psi, \mathbb{Q}} ,
    \end{displaymath}
  which is orthogonal with respect to $h$. By restriction of $h$, we obtain
  \begin{displaymath}
    h_{\psi}: \wedge^2_{K \otimes_{O_{F^t}, \tilde{\psi}} W(k)} P_{\psi, \mathbb{Q}} \times
    \wedge^2_{K \otimes_{O_{F^t}, \tilde{\psi}} W(k)} P_{\psi, \mathbb{Q}} \longrightarrow
    K \otimes_{O_{F^t}, \tilde{\psi}} W(k). 
    \end{displaymath}

  We find $\zeta \in O_F \otimes W(k)$
  such that $\sigma^{-1}(\zeta) \zeta^{-1} = -1$ or equivalently
  $\sigma(\zeta) = - \zeta$. We set $x = \zeta y$. Then we find
  \begin{displaymath}
    \wedge^{2} V (x) = \sigma^{-1}(\zeta) \wedge^{2} V (y) =
    - \sigma^{-1}(\zeta) py = - \sigma^{-1}(\zeta) \zeta^{-1}p x = px 
  \end{displaymath}
  Therefore  $\inv (\mathcal{P}, \iota, \beta)$ is the class of
  \begin{equation}\label{ramInv10e} 
h(\zeta y, \zeta y) = \zeta^{2} h(y, y) \; \mod \Nm_{K/F} K^{\times}.  
    \end{equation} 
  This shows in particular that
  $h(y,y) \in F^{\times}$ because $\zeta^2 \in F^{\times}$. We can replace in
  (\ref{ramInv10e}) the left hand side by $\zeta_{\psi} h_{\psi}(y,y)$ which
  gives for all $\psi$ the same element of $F$. The equation
  $\sigma(\zeta) = - \zeta$ may be written as
  $\sigma (\zeta_{\psi}) = - \zeta_{\psi \sigma}$. If we choose for a given $\psi$
  an element $\zeta_{\psi} \in O_F \otimes_{O_{F^t}} W(k)$ such that
  $\sigma^f (\zeta_{\psi}) = (-1)^f \zeta_{\psi}$, we obtain from this element
  a unique $\zeta$.

  In the case where $f$ is even, we can choose $\zeta_{\psi} = 1$, which proves
  the Lemma in this case. If $f$ is odd, we obtain that $\zeta \in F' \setminus F$. This
  implies as in the last Lemma that $\zeta^2 \notin \Nm_{K/F} K^{\times}$. This
  proves the case where $f$ is odd. 
\end{proof}

The following fact is well-known. 
\begin{lemma}\label{lemminun}
 Let $K/F$ be unramified. Let
  $u \in O_K \otimes_{O_{F^t}, \tilde{\psi}_0} W(k)$ be a unit.
Then there exists a unit $\zeta \in O_K \otimes_{O_{F^t}, \tilde{\psi}_0} W(k)$ such that
  \begin{displaymath}
\sigma^{-f}(\zeta) \cdot \zeta^{-1} = u. 
    \end{displaymath}\qed
  \end{lemma}

\begin{proposition}\label{KneunBa1l}
  Let $r$ be banal, and let $K/F$ be a field extension. Let $R = k$ be an
  algebraically closed field. Let $(C_{\mathcal{P}}, \iota, \phi)$  be the image
  of $(\mathcal{P}, \iota, \beta)\in \mathfrak{d}\mathfrak{P}^{\rm pol}_{r,k}$
  by the polarized contraction functor $\mathfrak{C}_{r, k}^{\rm pol}$, cf. Theorem \ref{KneunBa4p}. Then 
  \begin{displaymath}
    \inv(C_{\mathcal{P}}, \iota, \phi) = \inv^{r}(\mathcal{P}, \iota, \beta). 
  \end{displaymath}
       Here the $r$-adjusted invariant is given by 
      \begin{equation*}
  \inv^{r} (\CP, \iota, \beta) =  (-1)^{d}\inv (\CP, \iota, \beta). 
\end{equation*}
  \end{proposition} 
\begin{proof}
  We begin with the ramified case. We choose $\bar{\kappa}_E \subset k$. Let
  $~^{F} \eta \eta^{-1} = \pi^e/p=\rho$, $\eta \in O_F \otimes W(\bar{\kappa}_E)$
  as in (\ref{KneunBa6e}). We define
  $\tilde{\beta}: P \times P \longrightarrow O_F \otimes W(k)$ by (\ref{Cbanal-t})
  and the anti-hermitian form
  $\varkappa: P_{\mathbb{Q}} \times P_{\mathbb{Q}} \longrightarrow K \otimes W(k)$ by
  $\Trace \varkappa = \tilde{\beta}$. This $\varkappa$ differs from the
  $\varkappa$ of (\ref{Kneun03e}) by a constant in $F$. We can use it to
  compute $\inv^{r}(\mathcal{P}, \iota, \beta)$. We set 
  \begin{displaymath}
    \tilde{\beta}' = \eta \tilde{\beta}, \quad \varkappa' = \eta \varkappa.
  \end{displaymath}
We set $V' = \Pi^{-e} V$. Then we have 
$C_{\mathcal{P}} = \{y \in P \; | \; V'y = y \}$, cf. Remark \ref{CPexplban}.
  From this,  one deduces
  \begin{displaymath}
    \rho\cdot  ~^{F}\tilde{\beta}(y_1, y_2) =  \tilde{\beta}(y_1, y_2), \quad
    y_1, y_2 \in C_{\mathcal{P}}, 
  \end{displaymath}
  cf. (\ref{Cbanal8e}). This implies
  \begin{displaymath}
~^{F} \tilde{\beta}'(y_1, y_2) = \tilde{\beta}'(y_1, y_2) 
  \end{displaymath}
  The restriction of $\tilde{\beta}'$ to $C_{\mathcal{P}}$ is the form $\phi$, cf. Remark \ref{CPexplban}. 

We choose an element
$x \in \wedge^2 P_{\mathbb{Q}} := \wedge^2_{K \otimes W(k)} P_{\mathbb{Q}}$ such that
$\wedge^2 V (x) = (-1)^e p x$. By Lemma \ref{Kneuninv3l} the class of
$\wedge^2 \varkappa(x,x) \in F^{\times}/\Nm_{K/F} K^{\times}= \{ \pm1\}$ is
$(-1)^{fe}\inv(\mathcal{P}, \iota, \beta)= \inv^{r}(\mathcal{P}, \iota, \beta)$. 

We note that $\wedge^2 V' = (-1)^e \pi^{-e} \wedge^2 V$. We set
\begin{displaymath}
z = \eta^{-1} x \in \wedge^2 P_{\mathbb{Q}}.
  \end{displaymath}
Then we find
\begin{displaymath}
  \wedge^2 V' (z) = ~^{F^{-1}}(\eta^{-1}) (-1)^e \pi^{-e} \wedge^2 V (x) =
  ~^{F^{-1}}(\eta^{-1}) (-1)^e \pi^{-e} (-1)^e p x =
  ~^{F^{-1}}(\eta^{-1}) \eta \pi^{-e} p z = z. 
  \end{displaymath}
Therefore $z \in \wedge^2_{K} C_{\mathcal{P}} \otimes \mathbb{Q}$. The invariant
$\inv(C_{\mathcal{P}}, \iota, \phi)$ is given by $\wedge^2\varkappa'(z,z)$. 
Therefore the equality of invariants follows from
\begin{displaymath}
  \wedge^2 \varkappa'(z,z) = \eta^2 \wedge^2\varkappa (\eta^{-1}x, \eta^{-1}x) =
  \wedge^2 \varkappa(x,x). 
\end{displaymath}
Now we consider the case where $K/F$ is unramified. We use the notation
$H, h, H_{\psi}, h_{\psi}$ from (\ref{Kneuninv9e}). We have by
(\ref{comFF'unr}) that 
\begin{displaymath}
\wedge^2_{O_K} C_{\mathcal{P}} = \{z \in H \; | \; \wedge^2 V (z) = \pi^2_r z\}. 
\end{displaymath}
Using the decomposition (\ref{Kneuninv9e}),  the condition for
$z = (z_{\psi})$ becomes
\begin{displaymath}
\wedge^2 V (z_{\psi \sigma}) = \pi^{2a_{\psi}} z_{\psi}. 
\end{displaymath}
We choose $ z \neq 0$. Then
\begin{displaymath}
\inv(C_{\mathcal{P}}, \iota, \phi) = (-1)^{\ord_{\pi} h_{\psi}(z_{\psi}, z_{\bar{\psi}})},
  \end{displaymath}
 for any  $\psi$. We set $g_{\psi} = a_{\psi} + a_{\psi \sigma} + \ldots a_{\psi \sigma^{f-1}}$.
Then we obtain
\begin{displaymath}
  g_{\psi \sigma} - g_{\psi} = a_{\psi \sigma^{f}} - a_{\psi} = a_{\bar{\psi}} - a_{\psi}
  = e - 2a_{\psi}. 
  \end{displaymath}
We set $u_{\psi} = \pi^{g_{\psi}} z_{\psi}$. Then $u = (u_{\psi})$ satisfies
\begin{displaymath}
\wedge^2 V (u) = \pi^e u.
\end{displaymath}
Indeed, 
\begin{displaymath}
  \wedge^2 V (u_{\psi \sigma}) = \pi^{g_{\psi \sigma}} \wedge^2 V (z_{\psi \sigma}) =
  \pi^{g_{\psi \sigma}} \pi^{2a_{\psi}} z_{\psi} =
  \pi^{g_{\psi \sigma}} \pi^{2a_{\psi}} \pi^{-g_{\psi}} u_{\psi} = \pi^e u_{\psi}. 
\end{displaymath}
Let again $\eta \in O_F \otimes W(\bar{\kappa}_E)$ be the element defined after
(\ref{KneunBa6e}). It satisfies $\eta^{-1} ~^{F^{-1}}\eta = (p/\pi^e)$.
We set $x = \eta u$. Then $\wedge^2 V (x) = px$. Therefore
\begin{displaymath}
  \inv(\mathcal{P}, \iota, \beta) = (-1)^{\ord_{\pi} h_{\psi}(x_{\psi}, x_{\bar{\psi}})}, 
\end{displaymath}
for any $\psi$. Therefore the unramified case follows from
\begin{displaymath}
  \ord_{\pi} h_{\psi}(x_{\psi}, x_{\bar{\psi}}) =
  \ord_{\pi} h_{\psi}(\eta \pi^{g_{\psi}} z_{\psi},
  \eta \pi^{g_{\bar{\psi}}} z_{\bar{\psi}}) = (g_{\psi} + g_{\bar{\psi}}) +
  \ord_{\pi} h_{\psi}(z_{\psi}, z_{\bar{\psi}}) = ef +
  \ord_{\pi} h_{\psi}(z_{\psi}, z_{\bar{\psi}}). 
  \end{displaymath}
\end{proof}

\printindex[NO]

\printindex

\end{document}